%% file: Clean.IntroFlagVariety.tex
\documentclass[12pt]{amsart}
\usepackage{macros}

\title[Cohomology of the Flag Variety]{Introduction to the Cohomology of the Flag Variety}

\author{Sara C. Billey}
\address{SB: Department of Mathematics, University of Washington, Seattle, WA, US}
\email{billey@math.washington.edu}

\author{Yibo Gao}
\address{YG: Beijing International Center for Mathematical Research, Peking University. Beijing, China}
\email{gaoyibo@bicmr.pku.edu.cn}

\author{Brendan Pawlowski}
\address{BP: HRL Laboratories. Los Angeles, CA}
\email{br.pawlowski@gmail.com}

\date{\today} 

\begin{document}

\begin{abstract} 
One hundred years ago, Hilbert gave a list of important open problems
in mathematics in his day.  Hilbert's 15th problem asked for the
development of a rigorous calculus to explain the enumerative results
derived by Hermann Schubert in the later 1800s with regard to
intersecting varieties based on constraints imposed by rank conditions
on vector spaces.  Today by way of many contributions in algebraic
topology, geometry, and combinatorics, we consider this problem to be
fully solved.  Yet, deep questions remain in terms of the subtleties
of actually carrying out the process.  In this chapter, we hope to
summarize the contributions that have lead up to the rigorous
development of what has become known as Schubert calculus.  We will
survey some of the major developments in the past few decades with an
eye toward computation. We will discuss Grassmannian varieties and
flag varieties.  We will give a description of their cohomology rings
in algebraic terms inspired by the Chow ring of a smooth variety,
building on Monk's constructive approach.  From this description, we
will derive formulas for Schur polynomials and Schubert polynomials
which are the cohomology classes of Schubert varieties in these cases.
The Grassmannians, flag varieties, cohomology rings, Chow rings,
Schubert polynomials, Schur polynomials, etc. can be generalized in
many ways.  In this chapter, we will hint at the vast literature in
this area and point you to the other references in the Handbook for
more information.  Finally, we will identify some open problems that
remain a challenge even with all of the modern tools at our fingertips
in hopes of inspiring further contributions in this fascinating field.

\vspace{.5in}

\noindent Note to Readers: This article is intended to be the first
chapter of a book entitled ``Handbook of Combinatorial Algebraic
Geometry: Subvarieties of the Flag Variety'' that is a compendium of
topics in this area.  The book is being edited by Erik Insko, Martha
Precup, and Ed Richmond. In addition to this introductory chapter, the
other chapters, describing many of the other types of subvarieties of
the flag varieties and Grassmannians and some of their
generalizations, will each be written by experts in those areas. For
example, there are Kazhdan-Lusztig varieties (Chapter 2), generalized
smooth Schubert varieties (Chapter 3), Richardson varieties and
positroid varieties (Chapter 4), K-orbit closures (Chapter 5), torus
orbit closures (Chapter 6), spanning line configurations (Chapter 7),
different flavors of Hessenberg varieties (Chapters 8-11), and
generalizations to all Kac-Moody flag varieties (Chapter 12).  We hope
you will enjoy this chapter enough to seek out the others which are
either on the arXiv already or coming soon!  In the meantime, we hope
you will send us any comments or corrections you find as you read this
article.  We hope to have it polished up soon for publication in 2025.
Also note, the exercises presented here vary in difficulty.  The ones
prefaced by a citation are related to published work, so expect them
to be the more challenging problems.  You will benefit from attacking
them in the same way you would the easier problems, but you can always
consult the references for hints or more details.
\end{abstract}


\maketitle

\newpage

.\vspace{1in}
\begin{center}
 \includegraphics[height=10cm]{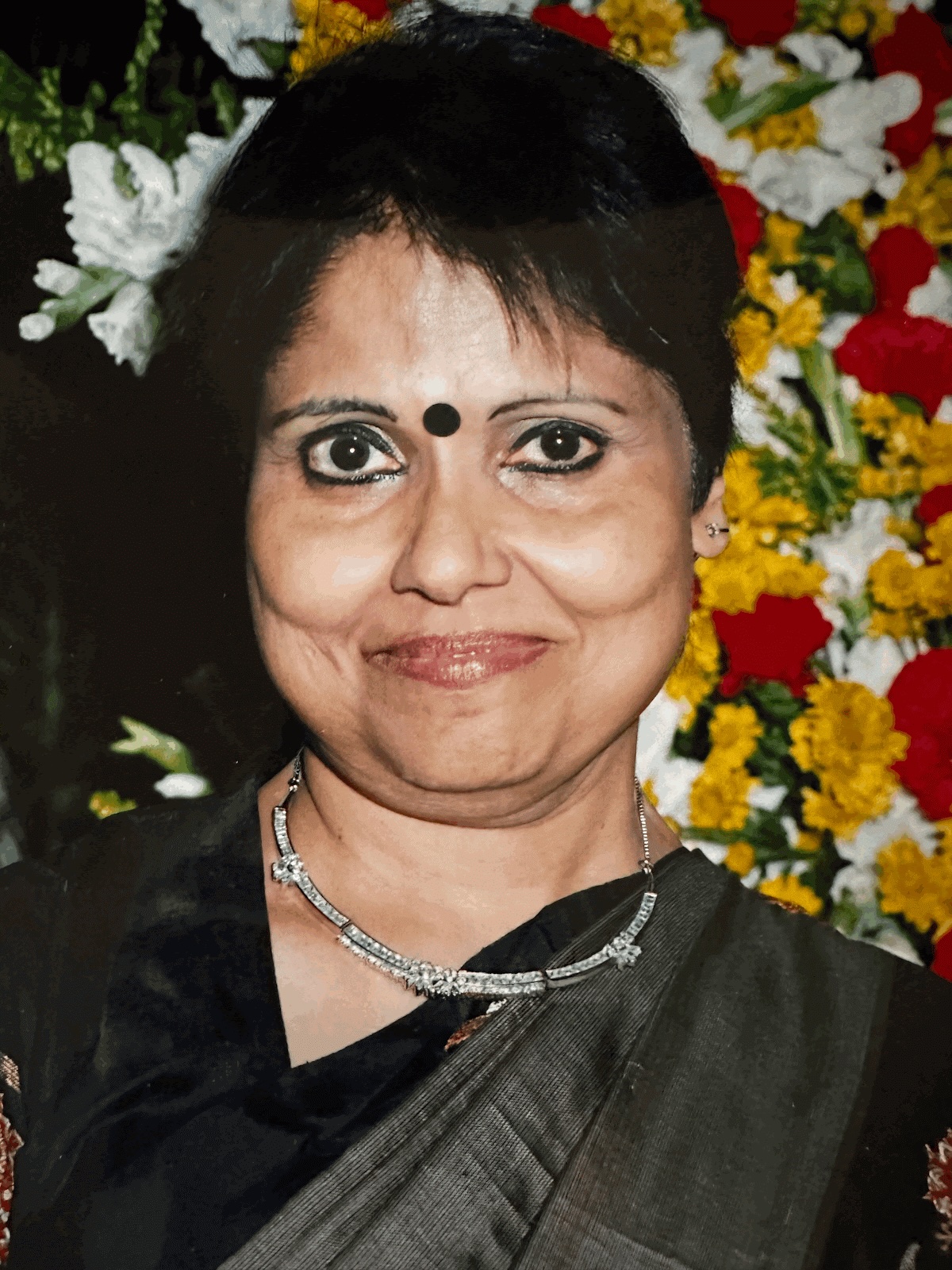}
\end{center}

\begin{center}
In Memoriam:  Lakshmibai 
%%%%%%%%%%%%%%%%%%%%%%%%%%%%%%%%%%%%%%%%%%%%%%%%%%%%%%%%%%%%
\end{center}

%%%%%%%%%%%%%%%%%%%%%%%%%%%%%%%%%%%%%%%%%%%%%%%%%%%%%%%%%%%%

\bigskip

This chapter is dedicated to the memory of Lakshmibai, who passed away
on December 2, 2023.  Her research contributions deeply enriched our
understanding of flag varieties and Grassmannian varieties.  She was a
dear friend and colleague, inspiring mathematicians around the
world. May her legacy endure in the annals of mathematical history.
Some of her work is highlighted in this chapter and beyond in later
chapters of the book.

\newpage

\vspace{1in}
\setcounter{tocdepth}{2}
\tableofcontents

\newpage
\input{section12.tex}

\input{section3.tex} %% formerly \input{flags.tex} \input{cohomology.tex}
 \input{section4.tex} %\input{CombinatoricsofSchubertPolys.tex}
 \input{section5.tex} %\input{CAGofSchubertVarieties.tex}  %%
		      %including {singular.tex}
\input{section6.tex} %% And Beyond
\newpage
\input{notation}

\newpage
\input{acknowledgments}

\newpage

\bibliographystyle{siam} \bibliography{flagrefs.bib}

\end{document}

%% file: section12.tex
\section{Introduction}\label{sec:intro}

What is ``Schubert calculus''?  What is a Schubert variety?  What is a
flag variety?  What is a cohomology ring?  What can I do with all of
that information if I were to learn about it?  Good questions!

The study of Schubert calculus and the associated structures in
algebraic geometry, algebraic topology, combinatorics, and
representation theory inspired some of the great advances in
mathematics in the twentieth century.  Inherently, the subject rests
firmly on the study of matrices, determinants, intersections of linear
spaces (meaning vector spaces and their translates), symmetry, and
computation.  If you know enough linear algebra, you know the basic
tools used in this chapter.  You can fill in missing details as you go
along using the references provided.  Our hope is that this chapter
will motivate you to learn more mathematics in the context of Schubert
calculus that also comes up in many applications in pure and applied
math.  The reason to learn this material now is because linear algebra
is revolutionizing our world.  Applications of linear algebra in web
search, linear optimization of business and government systems,
graphics, machine learning, DNA advances in biology, chemistry, and
physics can fill volumes.  Building a deep understanding of matrices
will serve you well no matter what hard problems you choose to attack
in the future.

\section{Enumerative Algebraic Geometry and Hilbert's 15th
Problem}\label{sec:enumerative.geom}

\subsection{Solving Schubert Problems in 1900}\label{sub:SchubertProblems1900}

The German and Italian schools of algebraic geometry were flourishing
in the mid 1800's. One area of intense study concerned questions in
enumerative geometry.

\begin{enumerate}
\item How many points are in the intersection of two lines in
  $\mathbb{R}^2$?  \\Answer: 0 or 1 or $\infty$.
\medskip
\item Given 2 circles in the plane, how many common tangents do they
have?  
\\Answer: 0, 1, 2, 3, 4, or $\infty$.  Draw pictures!
\medskip
\item Given 3 circles in the plane, how many circles are tangent to
all 3?  \\Answer: 0, 1, 2, 3, 4, 5, 6, 8, or $\infty$.  The generic solution has
8 circles known as the Circles of Apollonius, ca 200 BC.  The picture
below has colorfully rendered all 8 of the circles tangent to the
three black circles.

  \begin{center}
  \includegraphics[height=6cm]{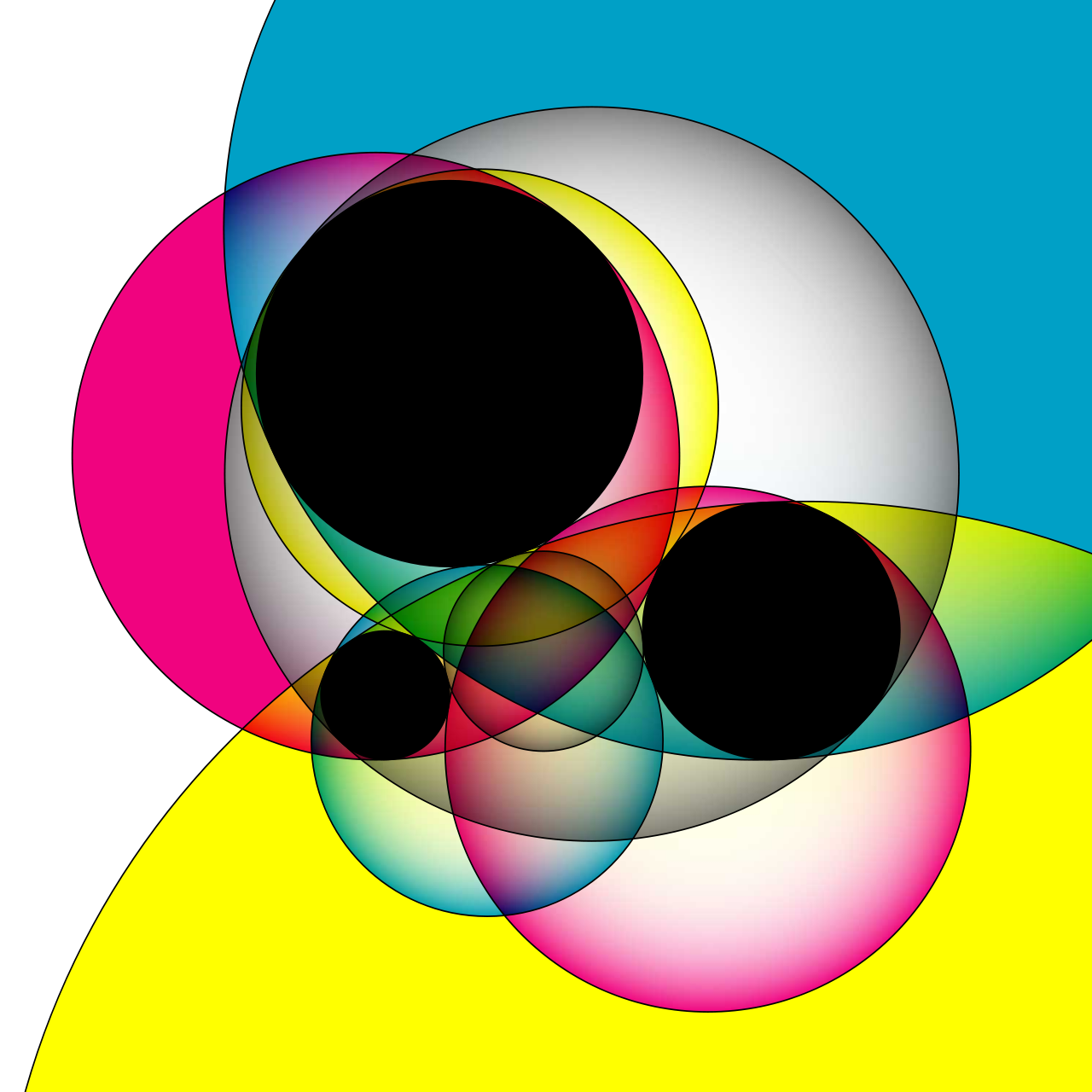}

  \begin{tiny}
\url{https://commons.wikimedia.org/wiki/File:Apollonius8ColorMultiplyV2.svg}
\end{tiny}
\end{center}
\end{enumerate}

Note, there are 33 different combinatorial types of circle
configurations to consider.  For example, the original circles could
all be the same or they could all be different but concentric.  It is
quite interesting that there is no way to get exactly 7 circles
tangent to 3 given circles.

What other enumerative questions might you ask about common objects
from geometry?  For example, how many intersection points are there in
the plane of two curves, one defined by a homogeneous polynomial of
degree $d_{1}$ and one of degree $d_{2}$?  A classical result known as
B\'ezout's Theorem from 1779 gives a precise answer: either there are
an infinite number of intersections or there are $d_{1}d_{2}$
intersection points, provided we count the intersection points with
appropriate multiplicity.  This type of product formula generalizes to
counting intersection points for any $n$ hypersurfaces defined by
homogeneous polynomials of degrees $d_{1}, \dots , d_{n}$. Another famous
fact from classical algebraic geometry is that a smooth cubic surface
contains exactly 27 lines; we will return to this result in
\Cref{ex:27.lines}.  Yet another impressive number along these lines is 3264,
which is the number of plane conics tangent to five general plane
conics \cite{eisenbud_harris_2016}.

Hermann C\"asar Hannibal Schubert (1848-1911) was a German
mathematician working in enumerative geometry.  He was particularly
interested in counting configurations of subspaces in given
arrangements.  The classic ``Schubert problem'' is:

\medskip
\begin{quote}
How many lines in $\mathbb{R}^{3}$ intersect 4 given lines?
\end{quote}
\medskip

\noindent The possible answers are 0, 1, 2, or infinity.  If three
distinct lines intersect a given family of 4 lines, then there must be
an infinite number of others because the position of the given family
of lines is not generic (cf. Exercise~\ref{exer:4-lines-3-lines}). For example, if the 4 given lines lie in a
plane, then an infinite number of other lines in that plane will
intersect all 4.  If you take 4 distinct parallel lines where three are in the
$z=0$ plane and one is the $z=1$ plane, then there are no lines in
$\mathbb{R}^{3}$ that go through all 4.  However, this is again not
the generic situation.  In the generic situation, that answer is
always 2.  While this question was originally inspired by pure
mathematical research, it also finds applications in modern day
computer graphics \cite{TH91}.   

More generally, Schubert was interested in testing when a family of
linear spaces in certain relative positions was intersected by only a
finite number of other linear spaces and determining the generic
number of solutions. Today, we call these types of questions
\textit{Schubert problems.} These types of relative positions give rise to
sets of matrices satisfying certain conditions on their determinantal
minors as we will explain in detail below.  His computations were all
on paper, without computers, modern tools in mathematics, or the
notation for matrices that we use today. How did he do it?  It remains
a mystery how Schubert accomplished as much as he did using the tools
at his disposal.  There was a sense among mathematicians of his day
that the theory Schubert was developing was not rigorous enough to
give formal algorithms.

One difficulty lies in properly interpreting the word
``generic''. Suppose we are trying to count the solutions to a
polynomial system whose coefficients are themselves polynomials in
some parameters, as in a Schubert problem. This number can, in
principle, be determined by checking polynomial equations and
inequalities in the parameters. For a simple example, the number of
solutions to a linear system of equations---or more generally, the dimension of
the solution space---is determined by the vanishing and nonvanishing of appropriate matrix minors. The subset of parameters where a polynomial equation holds will be strictly lower-dimensional, and so for almost every choice of parameters, {\bf none} of the equations will hold. The number of solutions in this case is the \emph{generic} number. For instance, as described above there will be exactly 2 lines intersecting 4 fixed lines in $\bR^3$ unless the configuration of 4 lines is degenerate in some way, and the degeneracy conditions can be expressed as polynomial equations in some parameters describing the lines.

Here is the trouble: Schubert's counting techniques relied on first
choosing a special kind of configuration which is easier to work
with. In our running example, if we choose two generic planes $P, Q
\subseteq \bR^3$ and generic lines $\ell_1, \ell_2 \subseteq P$ and
$\ell_3, \ell_4 \subseteq Q$, then it is easy to find two lines
intersecting $\ell_1, \ldots, \ell_4$, namely,  the line $P \cap Q$, and the
line through the points $\ell_1 \cap \ell_2$ and $\ell_3 \cap
\ell_4$. However, it is completely unclear if this special
configuration is actually generic in the required sense.

In 1900, David Hilbert famously summarized 10 major open problems in
Mathematics at the Paris conference of the International Congress of
Mathematicians (ICM) and published 23 problems in his follow-up
article in the Bulletin of the American Mathematical Society
\cite{Hilbert}.  Some information about the current status of these 23
turn of the 20th century problems is tracked on Wikipedia 
\cite{Hilbert.Problems.wiki}.  Among these influential problems was
Hilbert's 15th problem.

\bigskip

  \begin{center}
  \includegraphics[height=6cm]{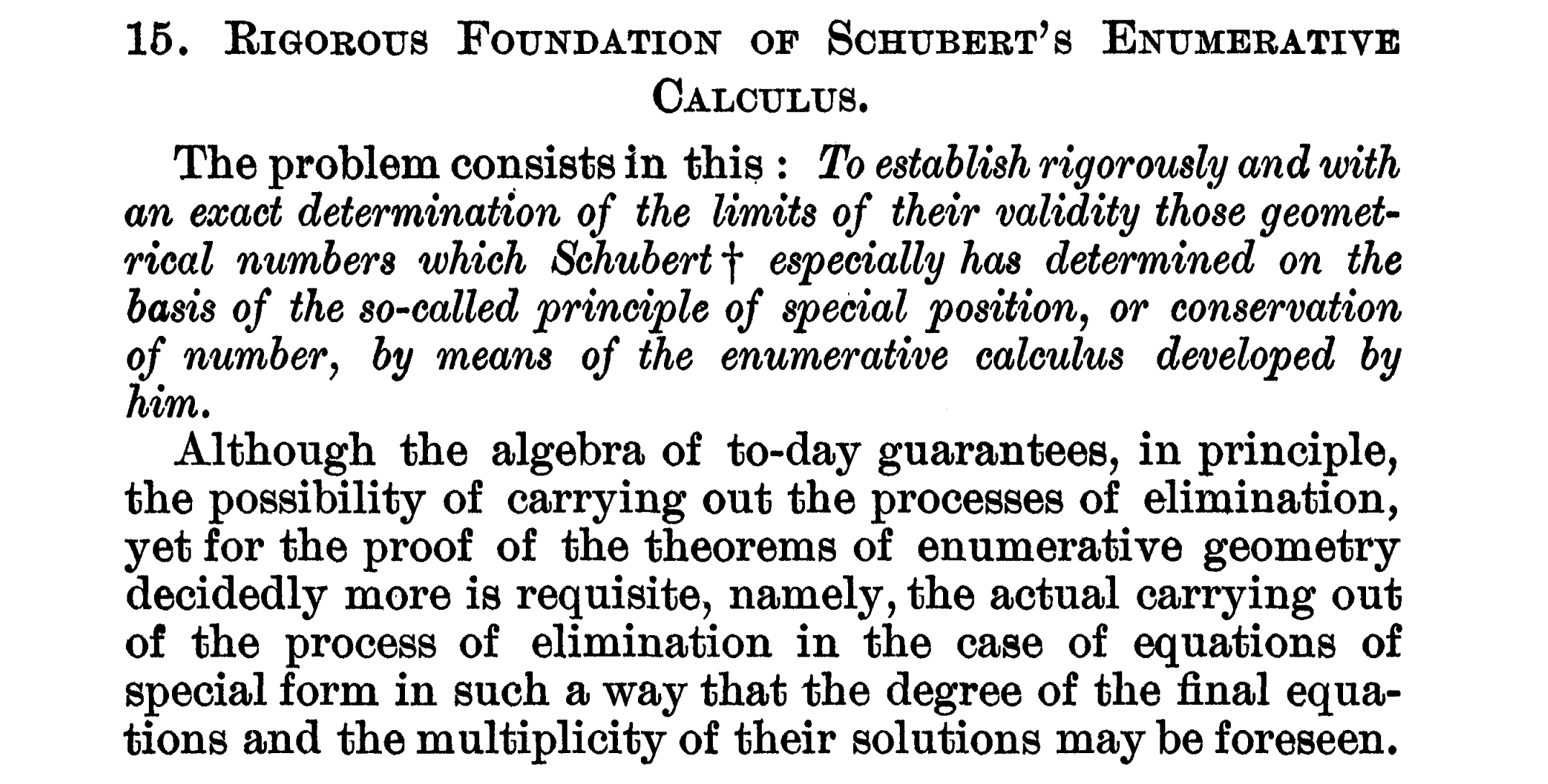}
  \end{center}

The problem of finding a ``rigorous foundation of Schubert's
enumerative calculus'' is not as specific as, say, the Riemann
Hypothesis (Problem 8).  One can still debate if Hilbert's 15th
problem is completely solved or not.  In our attempt to explain the
state of the art on solving Schubert-type problems, we will describe
an algorithm by which one could in principle solve any given problem
which is completely rigorous and has exact determination of the limits
of its validity.  Hence, we could say at this point that Hilbert's
15th problem is completely solved. That technology was a culmination
of research over the 20th century. However, these problems are bumping
into the limits of the computational hierarchy that has developed
around the Millennium Problem of determining if ``$P=NP?$'' or other
steps leading to the collapse of the polynomial hierarchy, which is
considered by many experts to be highly unlikely.  If Schubert was
able to do some of his more impressive calculations by hand, then we
must be missing something!

Schubert's calculus and Hilbert's 15th problem inspired many
developments in singular homology, cohomology, de Rham cohomology,
Chow rings, equivariant cohomology, quantum cohomology,
intersection theory, cobordism, combinatorics, representation theory
and beyond over the past 150 years.  Over the years, vocabulary and
methods have been developed that Schubert himself never used.  Around
1960, roughly a century after Hermann Schubert first started working
in this area, Bert Kostant named the corresponding sets of solutions
to intersection conditions on linear spaces \textit{Schubert
varieties}. The notion has been generalized far beyond the range of
problems Schubert himself was working on. So, a warning to the reader
is in order.  One might need to ask \dots
\medskip
\begin{quote}
Honey, where are my Schubert varieties? 
\end{quote}
\medskip
You will need to be mindful of the context for Schubert varieties.
Are they in a Grassmannian manifold or a flag manifold or do they
live in the setting of a classical Lie group, or an arbitrary Lie
group, or a Kac-Moody group, etc.  The context
will require a change in the definition of a Schubert variety.

\subsection{Solving Schubert Problems in 2000}\label{sub:SchubertProblems2000}

A modern \textit{Schubert problem} may be stated as
follows.
\begin{quote}
What is the expected number of points in the intersection of a given
family of Schubert varieties?
\end{quote}

\noindent To solve a Schubert problem in the 21st century, one
computes the product of two Schubert polynomials expanded in the basis
of Schubert polynomials and extracts a particular coefficient of a
basis element.  This approach leads to explicit algorithms using
linear algebra and avoids explicit computations of intersections and
questions of genericity.  However, these problems still become
intractable in high dimensions even with modern computers.  Narayanan
has shown that Schubert problems are at least as hard as the
computation class of $\#P$ problems \cite{Narayanan}, which includes
hard problems such as counting the number of Hamiltonian cycles in a
graph and finding the number of solutions to a binary linear
optimization problem.

The goal of this essay is a Revisionist History of Schubert problems,
how to solve them in the 21st century, and what it continues to
inspire.\footnote{Like Malcolm Gladwell's podcast ``Welcome to
Revisionist History'' where each episode he reinterprets something
from the past: an event, a person, an idea. Something
overlooked. Something misunderstood.}  We will follow a different
presentation path than what is considered the standard approach in the
literature.  The standard approach would begin with the work of
Bernstein-Gelfand-Gelfand \cite{BGG} and Demazure \cite{Dem}
introducing divided difference operators, the definition of Schubert
polynomials via the divided difference recurrence due to
Lascoux-Sch\"utzenberger \cite{LS1}, follow Macdonald's ``Notes on
Schubert Polynomials'' \cite{M2}, and arrive at the modern approach to
solving Schubert problems.  Instead, in \Cref{sec:flags} and
\Cref{sec:CombinatoricsofSchubertPolys}, we will follow Monk's earlier
approach to solving Schubert problems more directly as intersection
problems on Schubert varieties \cite{MONK}, which naturally could have
led to the Lascoux-Sch\"utzenberger \textit{transition equations} for
defining Schubert polynomials, and their visualizations using
\textit{pipe dreams} or the most recent tool \textit{bumpless pipe
dreams,} and then we will prove the divided difference recurrence as
an 
easy consequence of this approach in \Cref{sub:Fomin-Stanley}.  We
believe this nonstandard approach will be more intuitive and
constructive for the readers.  Along the way, we will define Schubert
varieties in a few different contexts.  We will survey some of the
beautiful known results on Schubert structure constants and related
combinatorics.  The last section returns to the geometry of Schubert
varieties, the exact equations defining them as projective varieties,
finding their singular loci, and some of the connections to
Kazhdan-Lusztig theory.

It is encouraged that readers do computations by hand and computer as
you read this chapter.  There are some computer packages available to
help you get started coding in Macaulay 2 \cite{MatrixSchubert.Mac2}
and Sage \cite{sage.combinat.schub}.

%% file: section3.tex
%%%%%%%%%%%%%%%%%%%%%%%%%%%%%%%%%%%%%%%%%%%%%%%%%%%%%%%%%
\section{Introduction to Flag Varieties and Schubert Varieties}\label{sec:flags}
%%%%%%%%%%%%%%%%%%%%%%%%%%%%%%%%%%%%%%%%%%%%%%%%%%%%%%%%%

This section lays out the basic definitions and notation that we will
need.  It builds on the ideas in Fulton's book on ``Young Tableaux: 
With Applications to Representation Theory and Geometry'' \cite[Part
3]{Fulton-book}. More information on these topics can also be found in
\cite{BLak,Gillespie.2019,kumar-book,Brown-Lakshmibai.Grassmannian,Brown-Lakshmibai.Flags,M2,manivel-book}.
See also a video version of this chapter given by the first author in 2021 at
an online conference at ICERM \cite{ICERM.videos}.

We do expect the readers to be familiar with some linear algebra,
combinatorics, and algebraic geometry.  Some excellent references
along these lines include the book by Cox, Little and
O'Shea entitled ``Ideals, Varieties, and Algorithms: An Introduction to
Computational Algebraic Geometry and Commutative Algebra,'' Richard Stanley's ``Enumerative
Combinatorics, Volumes 1 and 2,'' and Gilbert Strang's book ``Introduction to Linear
Algebra'' \cite{Cox-Little-OShea,ec1,ec2,Strang-LA}.

\subsection{Complete Flags}\label{sub:flags}

We begin with the basic definitions of flags in a complex
$n$-dimensional vector space, their matrix representations, and the
intuitive pictures we have in mind.  The field of complex numbers can
be replaced with other fields, but there are places where we are
assuming the field has characteristic zero and other places where we
need the field to be algebraically closed.  For this level of
exposition, we do not want to distract the reader with some of these
subtleties.

\begin{Definition}\label{def:flag}
Fix a positive integer $n$.  A \textit{complete flag} $F_{\ci}= (F_{1}
\subset \dots \subset F_{n})$ in $\mathbb{C}^{n}$ is a nested sequence
of subspaces such that $\mathrm{dim} (F_{i})=i$ for $1\leq i \leq n$.
A complete flag $F_\ci$ is determined by an ordered basis $(f_{1},
f_{2}, \dots, f_{n } )$ where $F_{i} = \mathrm{span}\langle
f_{1},\dots , f_{i} \rangle$.  Let $\Fl(n)$ denote the
\textit{complete flag variety}, which is the set of all complete flags
in $\mathbb{C}^{n}$.
\end{Definition}

More generally, given any subset ${\bf d} = \{d_1 < \cdots < d_m\}
\subseteq [n-1]$, a \emph{partial flag} with dimensions ${\bf d}$ is a
sequence of subspaces $F_1 \subset \cdots \subset F_m \subseteq
\C^n$ with $\dim F_i = d_i$.  Our main focus is on the complete flags,
so unless otherwise specified a ``flag'' means a complete flag.

\begin{Example}\label{example:n=4}
If $n=4$, let $(e_{1}, e_{2},e_{3},e_{4} )$ be the standard ordered
basis where $e_{i}$ has a 1 in the $i^{th}$ coordinate and 0's
elsewhere.  The corresponding \textit{standard flag} is denoted
$E_{\bullet }=(E_{1},E_{2},E_{3},E_{4})$ where each subspace
$E_{i}$ is the span of $e_{1},\dots ,e_{i}$.  The ordered basis
\[
(6 e_{1}+ 3e_{2},\hspace{.1in}4 e_{1} +2e_{3},\hspace{.1in}9e_{1}+ e_{3}+e_{4},\hspace{.1in}e_{2})
\]
determines another  flag $F_{\ci}=(F_{1} \subset F_{2}
\subset F_{3} \subset F_{4}) \in \Fl(4) $ where $F_{1}$ is the
1-dimensional subspace containing the origin and the point $6 e_{1}+
3e_{2}$,\ $F_{2}$ is the 2-dimensional subspace containing $F_{1}$ and
the point $4 e_{1} +2e_{3}$, \ $F_{3}$ is the 3-dimensional subspace
containing $F_{2}$ and the point $9e_{1}+ e_{3}+e_{4}$, and
$F_{4}=\mathbb{C}^{4}$. The reader should verify that these four
vectors are independent so each subspace $F_{i}$ has dimension $i$.
\end{Example}

\begin{figure}[h]
\begin{tikzpicture}[scale=0.6]
\draw[fill=black, fill opacity=0.1](0,0)--(4,0)--(4,3)--(0,3)--(0,0);
\draw(0,0)--(-1.5,-1.5)--(2.5,-1.5)--(4,0);
\draw(-1.5,-1.5)--(-1.5,1.5)--(0,3);
\draw(1,0)--(4,2);
\node at (2.2,0.8) {$\bullet$};
\end{tikzpicture}
\caption{Projective representation of a flag in $\mathbb{R}^{4}$ as a
point on a line in a plane, which is spanned by one wall of a
shoebox.  }\label{fig:one.blue.flag}
\end{figure}
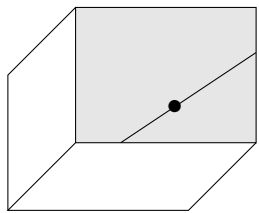

\begin{figure}[h]
\centering
\includegraphics[height=3cm]{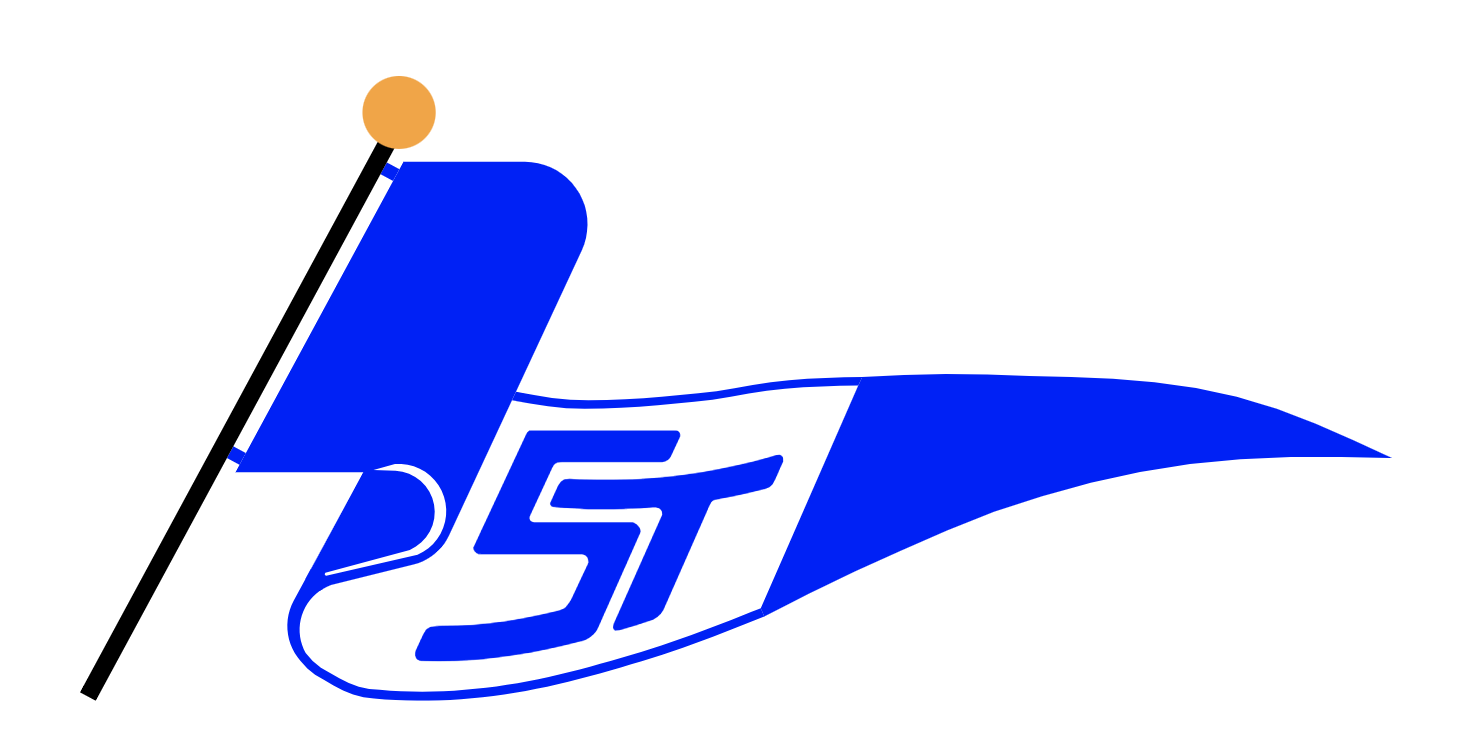}
\caption{A flag on a flag pole. Go Schubert Team!}
\label{fig:ST.flag}
\end{figure}

Every 1-dimensional subspace of $\mathbb{C}^{n}$ can be thought of as
a line through the origin.  Every 2-dimensional subspace can be
thought of as a plane containing the origin.  It is already hard to
visualize planes in say, $\mathbb{C}^{4}$, which is an 8-dimensional
space.  Intuitively, we prefer to think about $\mathbb{R}^{2}$ or
$\mathbb{R}^{3}$.  However, that will still limit our intuition to
small dimensions.  One additional trick that allows us to visualize
flags in $\mathbb{R}^{4}$ is to identify its linear subspaces by their
intersection with a fixed hyperplane at some distance from the
origin. The hyperplane is a flat 3-dimensional real object.  A typical
line through the origin meets the fixed hyperplane in exactly one
point.  We can ignore the case of a line being parallel to our chosen
hyperplane by perturbing it a little bit if necessary.  Every
2-dimensional subspace containing the origin in $\mathbb{R}^{4}$ meets
the hyperplane in a line.  Every 3-dimensional subspace of
$\mathbb{R}^{4}$ meets the hyperplane in a plane.  Therefore, drawn
projectively, a flag in $\Fl(4)$ is a point, on a line, in a plane,
contained in one side of a shoebox, which represents $\mathbb{C}^{4}$
projected on the fixed hyperplane.  See \Cref{fig:one.blue.flag}.  A
point, on a line, in a plane is reminiscent of a flag on a flag pole
such as we see in \Cref{fig:ST.flag}, hence the name
``flag''.\footnote{Image credit: Sloop Tavern Yacht Club.}

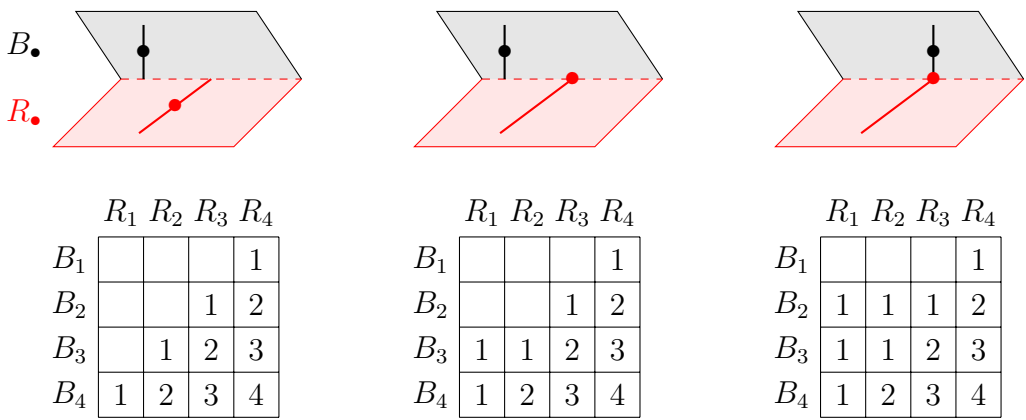
\begin{figure}[h]
\begin{tikzpicture}[scale=0.6]
\def\a{8};
\def\b{6};
\draw[step=1.0,thin] (0,0) grid (4,4);
\node at (0.5,0.5) {$1$};
\node at (1.5,0.5) {$2$};
\node at (2.5,0.5) {$3$};
\node at (3.5,0.5) {$4$};
\node at (0.5,1.5) {$ $};
\node at (1.5,1.5) {$1$};
\node at (2.5,1.5) {$2$};
\node at (3.5,1.5) {$3$};
\node at (0.5,2.5) {$ $};
\node at (1.5,2.5) {$ $};
\node at (2.5,2.5) {$1$};
\node at (3.5,2.5) {$2$};
\node at (0.5,3.5) {$ $};
\node at (1.5,3.5) {$ $};
\node at (2.5,3.5) {$ $};
\node at (3.5,3.5) {$1$};
\node[left] at (0,0.5) {$B_4$};
\node[left] at (0,1.5) {$B_3$};
\node[left] at (0,2.5) {$B_2$};
\node[left] at (0,3.5) {$B_1$};
\node[above] at (0.5,4) {$R_1$};
\node[above] at (1.5,4) {$R_2$};
\node[above] at (2.5,4) {$R_3$};
\node[above] at (3.5,4) {$R_4$};

\draw[step=1.0,thin] (\a,0) grid (\a+4,4);
\node at (\a+0.5,0.5) {$1$};
\node at (\a+1.5,0.5) {$2$};
\node at (\a+2.5,0.5) {$3$};
\node at (\a+3.5,0.5) {$4$};
\node at (\a+0.5,1.5) {$1$};
\node at (\a+1.5,1.5) {$1$};
\node at (\a+2.5,1.5) {$2$};
\node at (\a+3.5,1.5) {$3$};
\node at (\a+0.5,2.5) {$ $};
\node at (\a+1.5,2.5) {$ $};
\node at (\a+2.5,2.5) {$1$};
\node at (\a+3.5,2.5) {$2$};
\node at (\a+0.5,3.5) {$ $};
\node at (\a+1.5,3.5) {$ $};
\node at (\a+2.5,3.5) {$ $};
\node at (\a+3.5,3.5) {$1$};
\node[left] at (\a+0,0.5) {$B_4$};
\node[left] at (\a+0,1.5) {$B_3$};
\node[left] at (\a+0,2.5) {$B_2$};
\node[left] at (\a+0,3.5) {$B_1$};
\node[above] at (\a+0.5,4) {$R_1$};
\node[above] at (\a+1.5,4) {$R_2$};
\node[above] at (\a+2.5,4) {$R_3$};
\node[above] at (\a+3.5,4) {$R_4$};

\draw[step=1.0,thin] (2*\a,0) grid (2*\a+4,4);
\node at (2*\a+0.5,0.5) {$1$};
\node at (2*\a+1.5,0.5) {$2$};
\node at (2*\a+2.5,0.5) {$3$};
\node at (2*\a+3.5,0.5) {$4$};
\node at (2*\a+0.5,1.5) {$1$};
\node at (2*\a+1.5,1.5) {$1$};
\node at (2*\a+2.5,1.5) {$2$};
\node at (2*\a+3.5,1.5) {$3$};
\node at (2*\a+0.5,2.5) {$1$};
\node at (2*\a+1.5,2.5) {$1$};
\node at (2*\a+2.5,2.5) {$1$};
\node at (2*\a+3.5,2.5) {$2$};
\node at (2*\a+0.5,3.5) {$ $};
\node at (2*\a+1.5,3.5) {$ $};
\node at (2*\a+2.5,3.5) {$ $};
\node at (2*\a+3.5,3.5) {$1$};
\node[left] at (2*\a+0,0.5) {$B_4$};
\node[left] at (2*\a+0,1.5) {$B_3$};
\node[left] at (2*\a+0,2.5) {$B_2$};
\node[left] at (2*\a+0,3.5) {$B_1$};
\node[above] at (2*\a+0.5,4) {$R_1$};
\node[above] at (2*\a+1.5,4) {$R_2$};
\node[above] at (2*\a+2.5,4) {$R_3$};
\node[above] at (2*\a+3.5,4) {$R_4$};

\draw[draw=red, fill=red, fill opacity=0.1](1.5-1,\b+1.5)--(0-1,\b)--(4-1,\b)--(5.5-1,\b+1.5);
\draw[dashed,red](1.5-1,\b+1.5)--(5.5-1,\b+1.5);
\draw[fill=black, fill opacity=0.1](1.5-1,\b+1.5)--(0.5-1,\b+3)--(4.5-1,\b+3)--(5.5-1,\b+1.5);
\draw[thick](2-1,\b+1.5)--(2-1,\b+2.7);
\node at (2-1,\b+2.1) {$\bullet$};
\draw[thick,red](3.5-1,\b+1.5)--(1.9-1,\b+0.3);
\node[red] at (2.7-1,\b+0.9) {$\bullet$};
\node[left,red] at (0-1,\b+0.75) {$R_{\bullet}$};
\node[left] at (0-1,\b+2.25) {$B_{\bullet}$};

\draw[draw=red, fill=red, fill opacity=0.1](\a+1.5-1,\b+1.5)--(\a+0-1,\b)--(\a+4-1,\b)--(\a+5.5-1,\b+1.5);
\draw[dashed,red](\a+1.5-1,\b+1.5)--(\a+5.5-1,\b+1.5);
\draw[fill=black, fill opacity=0.1](\a+1.5-1,\b+1.5)--(\a+0.5-1,\b+3)--(\a+4.5-1,\b+3)--(\a+5.5-1,\b+1.5);
\draw[thick](\a+2-1,\b+1.5)--(\a+2-1,\b+2.7);
\node at (\a+2-1,\b+2.1) {$\bullet$};
\draw[thick,red](\a+3.5-1,\b+1.5)--(\a+1.9-1,\b+0.3);
\node[red] at (\a+3.5-1,\b+1.5) {$\bullet$};

\draw[draw=red, fill=red, fill opacity=0.1](2*\a+1.5-1,\b+1.5)--(2*\a+0-1,\b)--(2*\a+4-1,\b)--(2*\a+5.5-1,\b+1.5);
\draw[dashed,red](2*\a+1.5-1,\b+1.5)--(2*\a+5.5-1,\b+1.5);
\draw[fill=black, fill opacity=0.1](2*\a+1.5-1,\b+1.5)--(2*\a+0.5-1,\b+3)--(2*\a+4.5-1,\b+3)--(2*\a+5.5-1,\b+1.5);
\draw[thick](2*\a+3.5-1,\b+1.5)--(2*\a+3.5-1,\b+2.7);
\node at (2*\a+3.5-1,\b+2.1) {$\bullet$};
\draw[thick,red](2*\a+3.5-1,\b+1.5)--(2*\a+1.9-1,\b+0.3);
\node[red] at (2*\a+3.5-1,\b+1.5) {$\bullet$};
\end{tikzpicture}
\caption{Pairs of flags in 3 different relative positions along with
their intersection tables. The rows and columns are labeled the same
way in each case. The zero-dimensional intersections are
represented by empty cells in these tables for ease of reading.} \label{fig:rel.position}
\end{figure}

For each pair of flags, we can consider how their subspaces relate to
each other.  We can classify such pairs according to the
\textit{intersection table} of dimensions of the $i^{th}$ subspace in
the first flag intersected with the $j^{th}$ subspace of the second
flag.  For example, again in $n=4$, consider a flag $B_{\bullet}$ drawn
in black and $R_{\bullet}$ drawn in red as in \Cref{fig:rel.position}.
The tables of dimensions $\mathrm{dim}(B_{i} \cap R_{j})$ are shown
below each \emph{shoebox drawing}. Think of these drawings as planes
cutting through a box with a chosen special line and plane.\footnote{If you
have a shoebox, you might try to make a demo for yourself.}  The first
pair is in the most general position, and they get successively more
specialized in the sense that there are more incidences between the
components of the flags.

\bigskip

Recall that every complete flag $F_{\bullet}$ can be represented by an
ordered basis $(f_{1},f_{2},\dots f_{n })$.  In turn, every ordered
basis of $\mathbb{C}^{n}$ can be represented by an $n\times n$ complex
invertible matrix in $GL_{n}(\mathbb{C})$.  However, there are many
different ordered bases that represent the same flag, so we need to
spell out our conventions carefully.  Assume there is a fixed standard
ordered basis $(e_{1},e_{2},\ldots,e_n)$ of $\mathbb{C}^{n}$ and the
$f_{j}$'s are expressed in the standard basis as $f_{j}=
\sum_{i=1}^{n}m_{i,j}e_{i}.$  The column-wise matrix associated to an
ordered basis $( f_{1},f_{2},\dots, f_{n })$ is the matrix $M=M(
f_{1},f_{2},\dots, f_{n }) = (m_{i,j}) \in GL_{n}(\mathbb{C})$.  The
coefficients of the $f_{j}$'s expressed in terms of the $e_{i}$'s
become the columns of the matrix.  If we rescale a column of $M$ by a
nonzero complex number, the new matrix would represent the same
complete flag.  Furthermore, we could add any linear combination of
columns $1,2,\ldots, i$ to column $i+1$ and again the new matrix would
again represent the same flag $F_{\bullet}$. Thus, an infinite number
of different matrices in $GL_{n}(\mathbb{C})$ represent the same
complete flag.

Continuing with \Cref{example:n=4}, the flag
\[
F_{\ci} = (6 e_{1}+ 3e_{2},\hspace{.1in}4 e_{1}
+2e_{3},\hspace{.1in}9e_{1}+ e_{3}+e_{4},\hspace{.1in}e_{2} )
\]
can be represented by the matrices 

\begin{equation}\label{eq:example.options}
\left[\begin{array}{cccc}
6 &	4 &	 9 &	 0 \\
3 &	0 &	0 &	1\\
0 &    2 &     1 &	0\\
0 &    0 &      1 &      0
\end{array} \right]
  \sim
\left[\begin{array}{cccc}
2 &	2 &	 9 &	-2  \\
1 &	0 &	0 &	0\\
0 &   1 &    1 &	0\\
0 & 0 &	 1 &	0
\end{array} \right] \sim
\left[\begin{array}{cccc}
2 &	2 &	 7 &	 1 \\
1 &	0 &	0 &	0\\
0 &   1 &     0 &	0\\
0 & 0 &	 1 &	0
\end{array} \right].
\end{equation}
Hence, $F_{\ci}$ also can be determined by the ordered bases 
\[
(2 e_{1}+ e_{2},\hspace{.1in} 2 e_{1} +e_{3},\hspace{.1in}
9e_{1}+e_{3}+e_{4},\hspace{.1in} -2e_{1})
\]
or
\[
(2 e_{1}+ e_{2},\hspace{.1in} 2 e_{1} +e_{3},\hspace{.1in}
7e_{1}+e_{4},\hspace{.1in} e_{1}).
\]

Let $B$ be the set of upper triangular matrices in
$GL_{n}(\mathbb{C})$, and let $B_{-}$ be the set of lower triangular
matrices in $GL_{n}(\mathbb{C})$.  Let $w_{0}=[n,n-1,n-2,\dotsc , 1]$
be the permutation that reverses the normal order on $\{1,2,\ldots ,n
\}$.  We can represent $w_{0}$ as the matrix with 1's along the
antidiagonal and 0's elsewhere.  Then, $B=w_{0}B_{-}w_{0}$. Both subgroups are called
\textit{Borel subgroups}, hence the standard notation is to use the
letter $B$.  In general, a Borel subgroup of a linear algebraic group
like $GL_{n}$ is a maximal closed, connected, solvable algebraic
subgroup.  We won't need the full generality of Borel subgroups in
this chapter, but they do play a key role in the theory of flag
varieties more broadly.  Choosing a Borel subgroup is on par with
choosing an ordered basis for $\mathbb{C}^{n}$.

\begin{Exercise}\label{hw:flags.to.cosets}
Given a complete flag $F_{\ci}$ represented by the ordered basis
$(f_{1},f_{2},\dots, f_{n })$, let $M=M( f_{1},f_{2},\dots, f_{n })
\in GL_{n}(\mathbb{C})$.  Then, the set of all matrices in
$GL_{n}(\mathbb{C})$ that represent the same flag $F_{\ci}$ are
exactly the matrices of the form $Mb$ for $b \in B$.  
\end{Exercise}

By \Cref{hw:flags.to.cosets}, the set of all complete flags in
$\Fl(n)$ are in bijection with the cosets in
$GL_{n}(\mathbb{C})/B$. Therefore, some authors write $GL_{n}/B$
instead of $\Fl(n)$ to denote the complete flag variety.  More
generally, $G/B$ might refer to the complete flags or to the cosets of
another linear algebraic group $G$ such as $SL_{n}$, $Sp_{2n}$ or
$O_{n}$ modulo a Borel subgroup.  We will stick with $G=GL_{n}(\C)$ in
this chapter, and will simply write $\GL_n$ unless we explicitly want a different field or to emphasize the field $\C$.

It is always beneficial in any mathematical and/or computational
setting to choose a canonical representative of each object in
question.  This helps us determine when two objects are equal.  For
example, we could verify $6/14$ and $1047/2443$ are equal by reducing
both to their canonical representative $3/7$.  We use Gaussian
elimination as the guide to finding such a canonical representative for a flag,
but we must be mindful not to change the flag!

To perform \textit{column reduction} on a matrix $M$ that represents a
flag, one can do a sequence of elementary column operations of the
form:
\begin{itemize}
\item multiply column $j$ by a nonzero complex number, or
\item add a nonzero multiple of column $i$ to column $j$ for $i<j$.
\end{itemize}
The first type of elementary operation can be executed by multiplying
$M$ on the right by a diagonal matrix with 1's in every entry except
$(j,j)$ which is the nonzero complex number.  The second type of elementary
operation can be executed by multiplying $M$ on the right by an upper
triangular matrix with 1's along the diagonal, a nonzero value in
position $(i,j)$, and 0's everywhere else.  The group these two types
of matrices generate is exactly the Borel subgroup $B$ of invertible
upper triangular matrices.  Note, we cannot swap two columns of $M$
without changing the flag, so we may not be able to do the full
reduction to the column echelon form for $M$.  However, every matrix
representing a flag $F_{\bullet}$ will result in the same matrix after
column reduction.

\begin{Definition}\label{def:canonical.mat} 
A \textit{canonical matrix} in $GL_{n}$ is one that is column-reduced,
so it has exactly one pivot 1 in each row and column, 0's below each
pivot 1, and 0's to the right of each pivot 1.  Let $M(F_{\ci})$
denote the canonical matrix representing a flag $F_{\ci}$.
\end{Definition}

Continuing \Cref{example:n=4}, the canonical matrix representing 
$F_{\ci} $ is the simplest matrix in \eqref{eq:example.options}, namely
\begin{equation}\label{example:n=4.continued}
M( 6 e_{1}+ 3e_{2},\hspace{.1in}4 e_{1} +2e_{3},\hspace{.1in}9e_{1}+
e_{3}+e_{4},\hspace{.1in}e_{2})
=
\left[\begin{array}{cccc}
2 &	2 &	 7 &	 1 \\
1 &	0 &	0 &	0\\
0 &   1 &     0 &	0\\
0 & 0 &	 1 &	0
\end{array} \right].
\end{equation}
The pivot 1's are in positions $(2,1), (3,2), (4,3), (1,4)$.

Consider the following random experiment on a computer.  Fill an $8
\times 8$ matrix $M$ with random real numbers between 0 and 1.  Since
the determinant is a degree $8$ polynomial in the 64 entries, $\det(M)$
will almost surely not be zero.  Therefore, $M$ is almost surely
invertible so it will correspond to a flag $F_\bullet$ in $\Fl(8)$.
Furthermore, the $(8,1)$ entry of the chosen matrix is almost surely
not 0, so $F_{1} \not \subset E_{i}$ for any $i<8$ and the column
reduced matrix will have a $1$ in position $(8,1)$ after rescaling and
0's to its right.  Similarly, the lower left $2\times 2$ submatrix of
the matrix will almost surely be invertible so the $(7,2)$ entry of
the canonical matrix of $F_{\bullet}$ will be a pivot 1 with 0's to
its right.  In fact, every lower left submatrix of $M$ is almost
surely invertible.  Thus, the canonical matrix of $F_{\bullet }$ will
almost surely be in the column-reduced form \medskip
\begin{equation}\label{eq:w0matrix}
 \left[
\begin{matrix}
* & * & * & * & * & * & * & 1\\
* & * & * & * & * & * & 1 & 0\\
* & * & * & * & * & 1 & 0 & 0\\
* & * & * & * & 1 & 0 & 0 & 0\\
* & * & * & 1 & 0 & 0 & 0 & 0\\
* & * & 1 & 0 & 0 & 0 & 0 & 0\\
* & 1 & 0 & 0 & 0 & 0 & 0 & 0\\
1 & 0 & 0 & 0 & 0 & 0 & 0 & 0\\
\end{matrix}
\right] \medskip
\end{equation}
with the stars replaced by real numbers.
If we could choose an $8 \times 8$ matrix of complex numbers
uniformly at random, we would see the same type of canonical
matrix with the stars replaced by complex numbers.  Such an
experiment, if it could be done, would produce a
\textit{generically chosen flag}. However, we cannot choose a
flag uniformly at random this way since there is no way to even
choose a single complex number uniformly at random.  See
\Cref{ex:random.flag} for one approach to randomly generating
flags.  On the other hand, a less generic flag may have canonical
matrix of the form

\medskip

\begin{equation}\label{eq:46287351}
\left[
\begin{matrix}
* & * & * & * & * & * & * & 1\\
* & * & 1 & 0 & 0 & 0 & 0 & 0\\
* & * & 0 & * & * & 1 & 0 & 0\\
1 & 0 & 0 & 0 & 0 & 0 & 0 & 0\\
0 & * & 0 & * & * & 0 & 1 & 0\\
0 & 1 & 0 & 0 & 0 & 0 & 0 & 0\\
0 & 0 & 0 & * & 1 & 0 & 0 & 0\\
0 & 0 & 0 & 1 & 0 & 0 & 0 & 0\\
\end{matrix}
 \right].
\end{equation}

Observe from the definition of the canonical matrix representing
$F_{\ci}$ that the pivot 1's in $M(F_{\ci})$ form a permutation
matrix.  Which permutation is associated to a given set of pivot 1's?
Again we have some choices about how we label the corresponding
permutation matrix using an ordered list of the numbers ${1,2,\dots ,
n}$.  Let's take a step back to review some terminology on
permutations before we make that choice.

\begin{Exercise}\label{ex:perms}
There are at least 8 natural ways to label a permutation matrix using
an ordered list on $\{1,2,\dots ,n \}$.  These labels correspond with
taking any combination of the three bijections $w \rightarrow w^{-1}$,
\ $w \rightarrow w_{0}w$, \ and $w\rightarrow w w_{0}$.  What does
each of the corresponding 8 bijections do to a permutation matrix in
terms of rotation, reflection, etc?  How would you label the
permutation matrix with 1's in the positions shown in the matrix in
\eqref{eq:46287351}?
\end{Exercise}

\begin{Exercise}\label{ex:unitary}
Prove that every complete flag in $\mathbb{C}^{n}$ can be represented
by an $n \times n$ unitary matrix.  Furthermore, there is a bijection
between $GL_{n}(\mathbb{C})/B$ and $\mathrm{U}_{n}/T$.  Here $\mathrm{U}_{n}$ is the set
of $n \times n$ unitary matrices and $T$ is the $n\times n$ 
diagonal unitary matrices.   Therefore, $\flags$ naturally forms
a \emph{compact} topological space.
\end{Exercise}

\begin{Exercise}\cite{Diaconis.2003}\label{ex:random.flag} Let $A,B$
be independent random variables with the standard normal distribution.
Let $M$ be an $n \times n$ matrix with entries given by $n^2$
independent samples drawn from $A+Bi$.  Let $F_\bullet(M)$ be the
random flag determined by $M$. What is the corresponding distribution
on $\flags$?
\end{Exercise}  %%% answer: uniform on flags

%%% What about A,B uniform on [-1,1]??? 

\subsection{Permutations}\label{sub:perms}

Permutations are fundamentally bijections on the set of numbers
$[n] \coloneqq \{1,2,\dotsc , n \}$ to itself, or shuffles of a deck of cards,
or seating assignments on a full airplane given a numbering of the
seats and passengers, etc.  One may represent a permutation
$w:[n]\to [n]$ in \textit{one-line} notation as
$w=[w(1),w(2),\dotsc , w(n)]=[w_{1},\dotsc , w_{n}]$ or just
$w_{1}w_{2}\ldots w_{n}$ in examples where $n<10$ so the commas are not needed.  An
\textit{ascent} in $w$ is a position $i$ such that $w(i)<w(i+1)$, and
a \textit{descent} in $w$ is a position $i$ such that $w(i)>w(i+1)$.
More generally, an \textit{inversion} in $w$ is a pair of positions
$(i,j)$ such that $i<j$ and $w(i)>w(j)$.  If $i<j$ and $w(i)<w(j)$,
then $(i,j)$ is a \textit{non-inversion}.  Let
\begin{itemize}
\item $\inv (w)=\#\{(i,j)\in [n]^{2} \given i<j, w(i)>w(j)\}$ be the
\emph{number of inversions} in $w$,  and
\item $\des (w)=\#\{i\in [n-1] \given w(i)>w(i+1)\}$ be the \emph{number of
descents} in $w$.
\end{itemize}
The set of all permutations on $[n]$ is the \textit{symmetric group}
$S_{n}$. By an easy counting argument, there are $n!$ permutations in
$S_{n}$.  There are two $q$-analogs of $n!$ that have inspired a wealth of
research over the past century going back to Euler and MacMahon,

\[
A_{n}(q)  \coloneqq \sum_{w \in S_{n}} q^{\des(w)} 
\]
and 
\begin{equation}\label{eq:q-analog-n!}
[n]_{q}! \coloneqq \sum_{w \in S_{n}} q^{\inv(w)} =
\prod_{k=1}^{n-1}(1+q+q^{2}+ \cdots + q^{k}).
\end{equation}
The Eulerian polynomials $A_{n}(q)$ are so rich in structure that Kyle
Petersen \cite{Petersen} wrote a 456 page book about them!  The
standard $q$-analog of $n!$ given by $[n]_{q}!$ is much more important
in the context of flags as we shall see.

Multiplication in the symmetric group is given by function
composition $ v \circ w: [n]\to[n]$, so
\begin{equation}\label{eq:perm.products}
vw =[v(w(1)), v(w(2)), \dotsc , v(w(n))]
\end{equation}
in one-line notation for permutations $v,w \in S_{n}$.  This rule is
not commutative.  For example $[1,3,2][2,1,3]=[3,1,2]$ while
$[2,1,3][1,3,2]=[2,3,1]$.  Let $t_{ij}$ be the permutation in $S_{n}$
that swaps $i$ and $j$ and fixes all other elements in $[n]$ as a
bijection.  Then, the one-line notation of $vt_{ij}$ for $i<j$ is 
\[
vt_{ij}=[v_{1},\dots, v_{i-1},v_{j},v_{i+1}, \dots , v_{j-1},v_{i},
v_{j+1},\dots , v_{n}]. 
\]
The permutations $t_{ij}$ are called \textit{transpositions} and they
play a critical role in the theory that follows.  Among the
transpositions are the \textit{simple transpositions}
$s_{i}=t_{i,i+1}$ for $1\leq i\leq n$.  The simple transpositions are
a minimal set of generators for all of $S_{n}$ with the relations

\begin{equation}\label{eq:simple.relations}
\begin{array}{lll}
s_i^2=\id, & &\\
s_is_j=s_js_i  &\text{for all} &  |i-j|\geq2,\\
s_{i}s_{i+1}s_i=s_{i+1}s_is_{i+1} & \text{for all} & 1\leq i\leq n-2.
\end{array}
\end{equation}
We encourage the reader to check these relations hold.  Proving these
relations generate all possible relations among the simple
transpositions is more difficult.  See \cite{b-b,Hum} for details.

Let $M_{w}$ be the \textit{permutation matrix} for $w \in S_{n}$, so
$M_{w}$ has a 1 in position $(w_{j},j)$ for each $j\in [n]$ and 0's
everywhere else.  For example,
\[
M_{[4,2,6,1,3,5]}= \left[
\begin{matrix}
0 &	0 &	0 &	1 &	0 &	0 	\\ %1
0 &	1 &	0 &	0 &	0 &	0 	\\ %2
0 &	0 &	0 &	0 &	1 &	0 	\\ %3
1 &	0 &	0 &	0 &	0 &	0 	\\ %4
0 &	0 &	0 &	0 &	0 &	1 	\\ %5
0 &	0 &	1 &	0 &	0 &	0 	\\ %6
\end{matrix}
\right].
\]
The canonical matrices shown in \eqref{eq:46287351} have the pivot
1's in the same positions as the permutation matrix  $M_{[4,6,2,8,7,3,5,1]}$.

The transpose of $M_{w}$ is the permutation matrix for $w^{-1}$.
Matrix multiplication agrees with permutation multiplication in the
sense that $M_v M_w = M_{vw}$.  The identity matrix corresponds with
the identity permutation mapping each value $i$ to itself.  The reader
should verify each of these claims about the group structure of
matrices and permutations.

\begin{Remark}\label{rem:columns}
We note that some authors use the transpose matrix instead of
$M_{w}$. We chose this convention so that 
 matrix multiplication agrees with permutation multiplication,
while $M_{w}^{-1}M_{v}^{-1}=M_{x}^{-1}$ flips the order if $vw=x$. 
\end{Remark}

Permutations can be described in many other ways in addition to
one-line notation and permutation matrices. Here are some of the
most important ways in the context of flag varieties.  You may want to
look ahead at \Cref{ex.2341.table} and
\Cref{tab:permutation-representations} for a preview.

The \textit{northwest rank table} $rk(w)$ is obtained from $w$ by setting
\begin{equation}\label{eq:rank.def}
\rk(w)[i,j] = \# \{h \in [j] \given w(h) \leq i \}.
\end{equation}
The value $\rk(w)[i,j]$ is the rank of the \textit{northwest
submatrix} of $M_{w}$ with upper left corner $(1,1)$ and lower right
corner $(i,j)$.  One could similarly compute the southwest, northeast,
or southeast rank table for a permutation.

A \textit{string diagram} of a permutation $w$ is a braid with the
strings proceeding from the initial order to the permuted order given
by $w=w_1w_2\ldots w_n$ in such a way that no three strings cross at
any point; see Example~\ref{ex.2341.table}.  A \textit{wiring diagram} is a more rigidly drawn string
diagram with exactly one crossing in each row and piecewise straight strings,
or you may see such diagrams transposed to save space on paper as in
\Cref{fig:wiring.diag}.  Reading the rows from top
to bottom, we can associate a \textit{reduced word} to a wiring
diagram by labeling each crossing with one more than the number of
strings to its left.  A permutation can have many reduced words.  For
$321$, there are two reduced words $(1,2,1)$ and $(2,1,2)$
corresponding to the two minimal length expressions for
$321=s_{1}s_{2}s_{1}=s_{2}s_{1}s_{2}$ using the simple transpositions.
All reduced words for $w$ have the same \textit{length}, denoted
$\ell(w)$.  Furthermore, the length of $w$ is equal to the number of
\textit{inversions} for $w$,
\begin{equation}\label{eq:inv}
\ell(w)=\inv (w)=\#\{(i,j)\in [n]^{2} \given i<j, w(i)>w(j)\}.
\end{equation}
Note, $\ell(w)$ is a number between 0 and $\binom{n}{2}$ for any
permutation in $S_{n}$.  The \textit{longest permutation} in $S_{n}$
is the unique permutation $w_{0}=[n,n-1,\dots ,2, 1]$ with
$\ell(w_{0})=\binom{n}{2}$.  The identity permutation, denoted
$\id=[1,2,\dots , n]$, is the unique permutation in $S_{n}$ with
length equal to 0. In the context of Schubert varieties, $\ell(w)$ is more
common than $\inv(w)$ so we will use that notation in this chapter
unless we want to emphasize the connection with inversions.

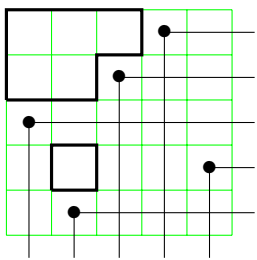
\begin{figure}[h!]
\begin{minipage}{0.45\linewidth}
\centering
\begin{tikzpicture}[scale=0.6]
\draw[step=1.0,green,thin] (0,0) grid (5,5);
\draw[very thick] (1,1)--(2,1)--(2,2)--(1,2)--(1,1);
\draw[very thick] (0,3)--(2,3)--(2,4)--(3,4)--(3,5)--(0,5)--(0,3);
\node at (0.5,2.5) {$\bullet$};
\node at (1.5,0.5) {$\bullet$};
\node at (2.5,3.5) {$\bullet$};
\node at (3.5,4.5) {$\bullet$};
\node at (4.5,1.5) {$\bullet$};
\draw(0.5,-0.5)--(0.5,2.5)--(5.5,2.5);
\draw(1.5,-0.5)--(1.5,0.5)--(5.5,0.5);
\draw(2.5,-0.5)--(2.5,3.5)--(5.5,3.5);
\draw(3.5,-0.5)--(3.5,4.5)--(5.5,4.5);
\draw(4.5,-0.5)--(4.5,1.5)--(5.5,1.5);
\node at (1.5,1.5) {};
\node at (0.5,3.5) {};
\node at (0.5,4.5) {};
\node at (1.5,3.5) {};
\node at (1.5,4.5) {};
\node at (2.5,4.5) {};
\end{tikzpicture}
\end{minipage}
\caption{The diagram of $w=43152$ is the set of outlined cells,
so \\
$D(w)=\{(1,1), (1,2), (1,3), (2,1), (2,2), (4,2) \}$. }
\label{fig:diagram.w}
\end{figure}

The \textit{diagram of a permutation} $w$ is obtained from the matrix
of $w^{-1}$ by removing all cells in an $n\times n$ array which are
weakly to the right or below a 1 in $w^{-1}$.  The remaining cells
form the diagram $D(w)$.  See \Cref{fig:diagram.w}.  Thus,
\begin{equation}\label{eq:diagram}
D(w) = \{(i,j)\in [n]^{2} \given j< w(i) \text{ and } i< w^{-1}(j) \}.
\end{equation}
Equivalently, the cells of $D(w)$ are in bijection with the inversions
of $w$, so 
\begin{equation}\label{eq:diagram.2}
D(w) = \{(i,w(j))\in [n]^{2} \given i<j \text{ and } w(i)> w(j) \}.
\end{equation}
It is unfortunate that the diagram is defined in terms of $w^{-1}$,
but that is the most common convention in the literature \cite{M2}.
Sometimes $D(w)$ is called the \textit{Rothe diagram} of $w$.

One more notation for a permutation $w \in S_n$ is its \emph{Lehmer
code}, or simply its \emph{code}, which is the $n$-tuple
\begin{equation}\label{eq:code}
c(w)=    (c(w)_1, c(w)_2, \ldots, c(w)_n)
\end{equation}
where $c(w)_i$ denotes the number of inversions $(i,j)$ with first
coordinate $i$, or equivalently, $c(w)_i$ is the number of boxes on
row $i$ of $D(w)$.  Note, $0 \leq c(w)_i \leq n-i$ for all $1 \leq i
\leq n$.  Hence, the Lehmer code gives a natural bijection between
$S_{n}$ and the product of sets $[n]\times [n-1] \times \cdots \times
[2] \times [1]$.

\vspace{0.1in}

\begin{Example}\label{ex.2341.table}  For the permutation $w=[2,3,4,1]$, we have the following
equivalent representations of $w$ shown in Table~\ref{tab:permutation-representations}.
\begin{table}[h!]
\centering
\begin{tabular}{|c|c|c|c|}
\hline
$2341$ & $\begin{bmatrix}1&2&3&4\\2&3&4&1\end{bmatrix}$ & $(1,2,3)$ & $(1,1,1,0)$  \\
one-line notation & two-line notation & reduced word & code \\\hline
$\begin{bmatrix}0&0&0&1\\1&0&0&0\\0&1&0&0\\0&0&1&0\end{bmatrix}$ & $\begin{bmatrix}0&0&0&1\\1&1&1&2\\1&2&2&3\\1&2&3&4\end{bmatrix}$ & \adjustbox{valign=m}{\begin{tikzpicture}[scale=0.5]
\draw[step=1.0,green,thin] (0,0) grid (4,4);
\node at (0.5,0.5) {$\bullet$};
\node at (1.5,3.5) {$\bullet$};
\node at (2.5,2.5) {$\bullet$};
\node at (3.5,1.5) {$\bullet$};
\draw(0.5,-0.5)--(0.5,0.5)--(4.5,0.5);
\draw(1.5,-0.5)--(1.5,3.5)--(4.5,3.5);
\draw(2.5,-0.5)--(2.5,2.5)--(4.5,2.5);
\draw(3.5,-0.5)--(3.5,1.5)--(4.5,1.5);
\draw[very thick] (0,1)--(1,1)--(1,4)--(0,4)--(0,1);
\node at (0.5,4.2) {};
\end{tikzpicture}} & \adjustbox{valign=m}{\begin{tikzpicture}[scale=0.5]
\node at (0,0) {$\bullet$};
\node[below] at (0,0) {$2$};
\node at (1,0) {$\bullet$};
\node[below] at (1,0) {$3$};
\node at (2,0) {$\bullet$};
\node[below] at (2,0) {$4$};
\node at (3,0) {$\bullet$};
\node[below] at (3,0) {$1$};
\node at (0,2) {$\bullet$};
\node[above] at (0,2) {$1$};
\node at (1,2) {$\bullet$};
\node[above] at (1,2) {$2$};
\node at (2,2) {$\bullet$};
\node[above] at (2,2) {$3$};
\node at (3,2) {$\bullet$};
\node[above] at (3,2) {$4$};
\draw(1,2)--(0,0);
\draw(2,2)--(1,0);
\draw(3,2)--(2,0);
\draw[bend right=30](0,2) to (1.5,1);
\draw[bend left=30](1.5,1) to (3,0);
\end{tikzpicture}} \\
matrix notation & rank table & Rothe diagram & string diagram \\\hline
\end{tabular}
\caption{Eight different ways to represent the same permutation $2341$}
\label{tab:permutation-representations}
\end{table}
\end{Example}

\begin{Example}\label{ex:big.8}
Try drawing the permutation diagram for $83617254$ for yourself.
Then observe that the stars in the canonical matrices given in
\eqref{eq:46287351} appear exactly in the positions corresponding with
the diagram of $D(8 3 6 1 7 2 5 4)$ and $[4, 6, 2, 8, 7, 3, 5,
1]^{-1}= [8, 3, 6, 1, 7, 2, 5, 4]$.
\end{Example}

Another crucial concept that appears everywhere in Schubert calculus
is \emph{pattern avoidance}.
See \Cref{sub:Singularloci} and several of the other chapters in this
book for more applications of pattern avoidance and variations on that
theme.

\begin{Definition}\label{def:pattern-avoidance}
A permutation $w\in S_n$ \emph{contains} the pattern $\pi\in S_k$ if
there exists indices $1\leq a_1<a_2<\cdots<a_k\leq n$ such that
$w(a_i)<w(a_j)$ if and only if $\pi(i)<\pi(j)$ for all $1\leq i<j\leq
k$. We say that $w\in S_n$ \emph{avoids} $\pi\in S_k$ if $w$ does not
contain the pattern $\pi$.
\end{Definition}

In other words, $w$ contains the pattern $\pi$ if in one-line
notation, there exists a (not necessarily consecutive) subsequence of
$w$ with the same relative ordering as $\pi$.  One can also think of a
pattern $\pi$ in $w$ as a submatrix of the matrix $M_{w}$.  

\begin{Example}
The permutation $w=6 2 5 4 3 1$ contains the subsequence $6231$ which
has the same inversion set as $4231$, so $w$ contains a $4231$
pattern. However, $w$ avoids the patterns $2143$ and $3412$.
\end{Example}

The early history of permutation pattern avoidance goes back to
MacMahon \cite{MacMahon}, Knuth \cite{MR0378456}, Tarjan
\cite{tarjan}, and Pratt \cite{Pratt:1973}.  An early result shows
that for every $\pi\in S_3$, the number of permutations in $S_n$
avoiding $\pi$ equals $\frac{1}{n+1}{2n\choose n}$, the $n^{th}$
Catalan number \cite{MacMahon,SS}.

\begin{Exercise}[\textbf{Exchange Lemma}]\label{ex:exchange}
Use string diagrams to prove that if $(a_{1},a_{2},\dots , a_{p})$ and
$(b_{1},b_{2},\dots , b_{p})$ are both reduced words for $w \in
S_{n}$, then there exists an $i$ such that $(a_{1}, b_{1},\dots,
\widehat{b_{i}},\dots , b_{p})$ is also a reduced word for $w$.
\end{Exercise}

\begin{Exercise}[\textbf{Tits/Matsumoto Theorem}]\label{ex:tits}
Let $G_{w}$ be the graph on all reduced words for $w$ with edges
connecting two words if they differ by a \emph{commutation relation}
$s_{i}s_{j}=s_{j}s_{i}$ or \emph{braid relation}
$s_{i}s_{i+1}s_i=s_{i+1}s_is_{i+1}$ \eqref{eq:simple.relations}.
Prove $G_{w}$ is connected.   
\end{Exercise}

\begin{Exercise}\cite[P0062]{dppa}\label{ex:grassmannian.def}
A permutation $w \in S_n$ is \emph{Grassmannian} if it has at most one
descent.   What set of patterns characterize the Grassmannian
permutations?
  \end{Exercise}
  
\begin{Exercise}\label{ex:bigrassmannians} \cite[P0063]{dppa}  %% noted by Josh Swanson
We say $w$ is bigrassmannian if both $w$ and $w^{-1}$ are
Grassmannian.  What set of permutation patterns characterizes the
bigrassmannian permutations?  How can you characterize them in terms
of their Rothe diagrams?  See \Cref{rem:presentaiton.cohomology} for
applications.  \end{Exercise}

\begin{Exercise}\label{ex:perm.bijections}
Describe how $w \in S_{n}$ can be recovered from its rank table,
string diagram, Rothe diagram, or its inversion set.
\end{Exercise}

\begin{Exercise}
Prove \Cref{eq:q-analog-n!} using the bijection between permutations
and their Lehmer codes, so
\[
[n]_{q}! \coloneqq \sum_{w \in S_{n}} q^{\inv(w)} = 
\prod_{k=1}^{n-1}(1+q+q^{2}+ \cdots + q^{k}).
\]
\end{Exercise}

\begin{Exercise}\label{ex:red.word}
For $w \in S_{n}$, fill the cells in $D(w)$ with positive integers as
follows.  For each $1\leq i\leq n$, starting with the leftmost cell in
row $i$, fill the cells from left to right
consecutively with values $i,i+1,i+2, \dots $.  Prove that the word
obtained by reading along the rows from top to bottom, right to
left is a reduced word for $w$.  In fact, it is the largest reduced
word for $w$ in reverse lexicographic order.
\end{Exercise}

One fun fact about Lehmer codes is that they determine the positions
of the permutations in $S_{n}$ when written out in lexicographic
order.  We leave this as an exercise for the reader.  

\begin{Exercise}\label{ex:code-to-lex-order}
Let $L_{n}=(w^{(0)}, w^{(1)},\dots , w^{(n!-1)})$ be the list of all
permutations in $S_{n}$ in lexicographic order, so $w^{(0)}$ is the
identity and $w^{(n!-1)}$ is the longest permutation $w_{0}$.  For $w
\in S_{n}$ with Lehmer code $c(w)=(c_{1},c_{2},\dots, c_{n})$, show that
$w=w^{(j)}$ for $j = \sum c_{i}(n-i)!$.
\end{Exercise}

\subsection{Schubert Cells and Schubert Varieties}\label{sub:Schubert.Varieties}

Recall that we have chosen a fixed ordered basis $(e_{1},
e_{2},\ldots, e_{n})$ for $\mathbb{C}^{n}$ to bijectively identify
$\Fl(n)$ with $GL_{n}/B$.  Let $E_{\ci}=(E_{1}, E_{2},\ldots, E_{n})$
be the flag that corresponds with this fixed ordered basis, so $E_{i}$
is spanned by $\{e_{1},e_{2},\dots , e_{i} \}$.  Each flag $F_{\ci}$
has a \textit{position} relative to this fixed flag, which is
determined by the matrix of values $\mathrm{dim} (E_{i} \cap F_{j})$.
By considering the canonical matrix representative of $F_{\ci}$, we
observe that such a table is always the rank table of some permutation
$w \in S_{n}$, so we write $\mathrm{position}(E_{\ci},F_{\ci})=w$.
The flags $E_{\ci},F_{\ci}$ are in \textit{transverse position} or
\textit{general position} if
$\mathrm{position}(E_{\ci},F_{\ci})=w_{0}$.

A table of intersection dimensions is reminiscent of the Schubert
problems we discussed in the introduction.  We could ask which other
flags intersect the subspaces $E_{1}, E_{2},\dots , E_{n}$ in the same
table of intersection dimensions?  This gives rise to the concept of a
\emph{Schubert cell} in the complete flag variety and its closure which is called a
\emph{Schubert variety}.  We spell out the details and give examples below.

\begin{Definition} \label{def:Schubert.cell}
For a permutation $w \in S_n$, the \textit{Schubert cell}
$C_{w}(E_{\ci})\subset \Fl(n)$ is the set of all flags $F_{\ci}$ with
$\mathrm{position}(E_{\ci},F_{\ci})=w$. Equivalently, we can write
$C_{w}(E_{\ci})$ using the northwest rank table of $w$ as 
\begin{align*}
C_{w}(E_{\ci}) = 
\{F_{\ci} \in \flags 
\given \mathrm{dim}(E_{i} \cap F_{j}) = \mathrm{rk}(w)[i,j] \text{ for all } 1\leq i,j \leq n \}.
\end{align*}
If the fixed flag is clear from context, we may write $C_{w}$ instead
of $C_{w}(E_{\ci}).$ 
\end{Definition}

Note, each permutation $w \in S_{n }$ gives rise to a special flag
$w_{\ci}$ given by the ordered basis $(e_{w(1)}, e_{w(2)},\ldots,
e_{w(n)})$.  The matrix corresponding to $w_{\ci}$ is the permutation
matrix $M_{w}$, which is in canonical form.  By definition, $w_{\ci}
\in C_{w}(E_{\ci})$, so each Schubert cell is nonempty.

\begin{Example}\label{example:n=4.continued.2} Identifying the flag from \Cref{example:n=4}
with its canonical matrix shown in \eqref{example:n=4.continued}, we have
\begin{align*}
F_{\ci} =    
\left[\begin{array}{cccc}
2 &	2 &	 7 &	 \enci 1 \\
\enci 1&0 &	0 &	0\\
0 &     \enci 1 &     0 &	0\\
0 &     0 &	 \enci 1 &	0
\end{array} \right]
\in 
\displaystyle C_{2341} = 
\left\{\left[ 
\begin{array}{cccc}
x &	 y &	 z &	 1 \\
1 &	0 &	 0  &	0\\
0 &    1 &     0 &	 0\\
 0 & 0 & 1 &	0
\end{array}
 \right]:  x,y,z\in \mathbb{C} \right\} \subset GL_{4}/B.
\end{align*}
The $x,y,z$ variables parameterize the Schubert cell $C_{2341}$: every
choice of complex numbers for these 3 variables gives rise to a
distinct flag in $C_{2341}$.  Here $F_{\ci}$ is the flag with ordered
basis $(2 e_{1}+ e_{2},\hspace{.1in} 2 e_{1} +e_{3},\hspace{.1in}
7e_{1}+e_{4},\hspace{.1in} e_{1}).$ We see $F_{1}$ is spanned by $2
e_{1}+ e_{2}$ so $E_{1} \cap F_{1} $ is just the origin and $F_{1}
\subset E_{2}$ so $\mathrm{dim} (E_{1} \cap F_{1} ) =0$ and
$\mathrm{dim} (E_{i} \cap F_{1} ) =1$ for all $i\geq 2$.  We can
similarly compute $\mathrm{dim}(E_{i} \cap F_{j} )$ for all $1 \leq
i,j\leq 4$ which is the rank table for $w=2341$,
\begin{equation}\label{eq:rank.table.2341}
\left[\begin{array}{cccc}
0 &	0 &	0 &	1 \\
1 &    1 &	1 &	2\\
1 &    2 &      2 &	3\\
1 &    2 &      3 &	4
\end{array} \right]. 
\end{equation}
Every other flag with the same table of intersection dimensions with
$E_{\ci}$ is also in the Schubert cell $C_{2341}(E_{\ci})$.  One could
think of $C_{w}(E_{\ci})$ as the equivalence class of flags in the
same position $w$ with respect to $E_{\ci}$.  Each equivalence class is
determined by a set of incidence relations.  We can visualize the
equivalence class $C_{2341}(E_{\ci})$ using shoebox diagrams.  Draw
the standard flag in black in a shoebox.  Then add a red flag
$R_{\bullet}$ in $C_{2341}(E_{\ci})$ to the picture.  The red
plane and the black plane meet in a line, which is the red line since
$\mathrm{dim}(E_{3}\cap R_{2})=2$.  The red point is at the point of
intersection between the red line and the black line since
$\mathrm{dim}(E_{2}\cap R_{1})=1$ and $R_{1} \subset R_{2}$.
\end{Example}

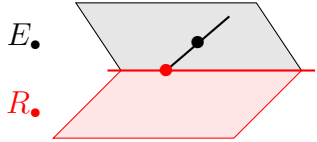
\begin{figure}[h]
\begin{tikzpicture}[scale=0.6]
\draw[red, fill=red, fill opacity=0.1](1.5,1.5)--(0,0)--(4,0)--(5.5,1.5);
\draw[fill=black, fill opacity=0.1](1.5,1.5)--(0.5,3)--(4.5,3)--(5.5,1.5);
\draw[thick](2.5,1.5)--(3.9,2.7);
\node at (3.2,2.1) {$\bullet$};
\draw[thick,red](1.2,1.5)--(5.8,1.5);
\node[red] at (2.5,1.5) {$\bullet$};
\node[left,red] at (0,0.75) {$R_{\bullet}$};
\node[left] at (0,2.25) {$E_{\bullet}$};
\end{tikzpicture}
\caption{A shoebox diagram of two flags in $\Fl(4)$.  The black flag represents
the standard flag $E_{\bullet}$.  The red flag $R_{\bullet}$ is an
element of  $C_{2341}(E_{\bullet}).$}
\end{figure}

\begin{center}

\end{center}
\bigskip

From the intersection conditions for $C_{w}(E_{\ci})$ in
\Cref{def:Schubert.cell}, one can observe that there is a close
connection between the diagram $D(w^{-1})$ and the canonical matrix
representatives for $C_{w}(E_{\ci})$.  A flag $F_{\bullet}$ is in $C_{w}(E_{\ci})$
if and only if its canonical matrix representative $M(F_{\bullet})$
has 1's in the same positions as the permutation matrix $M_{w}$.  Any
other nonzero entries in $M(F_{\bullet})$ must appear in positions
$(i,j) \in D(w^{-1})$.  Therefore, the elements in $D(w^{-1})$
determine a parametrization of the Schubert cell $C_{w}(E_{\ci})$.
Since $\inv (w)=\inv (w^{-1})$, we conclude that the complex dimension
of a Schubert cell is the number of inversions of the permutation,
\begin{equation}\label{eq:dim.Schubert.cell}
\mathrm{dim}(C_{w}(E_{\ci}))= \ell(w) = \inv (w).
\end{equation}
The dimension of $\flags$ is $\binom{n}{2}$ since
$\mathrm{dim}(GL_{n})=n^{2}$ and $\mathrm{dim}(B)=\binom{n+1}{2}$.
Thus, the codimension of $C_{w}(E_{\ci})$, denoted $\codim
(C_{w}(E_{\ci}))=\binom{n}{2}-\ell (w)$, is the number of noninversions.  We
will denote this \textit{coinversion} statistic by $\coinv(w):=\binom{n}{2}-\ell(w)$
combining codimension and noninversion.  

Recall that if a matrix $M$ represents a flag, then every matrix $Mb$
for $b\in B$ represents the same flag.  What about multiplication on
the left?  Multiplication on the left by upper triangular matrices in
$B$ can be described as adding scalar multiples of row $j$ in $M$ to
row $i$ for $i\leq j$.  Assume $M$ is a canonical matrix so it is
column-reduced.  If a column of $M$ represents a vector in $E_{i}$,
then it has zeros below row $i$.  Thus, every matrix of the form $bM$
for $b \in B$, will also have zeros below row $i$ in that column. If
the flag represented by $M$ satisfies the intersection conditions for
$C_{w}(E_{\ci})$ then so does every other flag represented by $bM$
for some $b \in B$.  The matrices in the orbit $BM$ may not be in column-reduced form.  However, since we know that every flag in
$C_{w}(E_{\ci})$ can be represented by a matrix in canonical form with
pivot 1's in positions corresponding with the permutation matrix
$M_{w}$ and any other nonzero entries only occurring in positions in $D(w)$, 
we conclude that every flag in $C_{w}(E_{\ci})$ can be represented by
a matrix in the orbit $BM_{w}$.
Hence, $C_{w}(E_{\ci})$ is exactly
the $B$-orbit of the flag $w_{\bullet}=(e_{w(1)},e_{w(2)},\ldots ,
e_{w(n)})$ under the action of left multiplication induced from the
matrix representatives.  Equating flags with their matrix
representatives, one may write $C_{w}(E_{\ci}) = BM_{w}$.

%\vspace{0.1in}

\begin{Example}\label{ex:bs}
For arbitrary $b_{i,j}$'s with $b_{1,1},b_{2,2},b_{3,3}, b_{4,4}$ nonzero, we have
\begin{small}
\[
 \left[ \begin {array}{cccc} b_{{1,1}}&b_{{1,2}}&b_{{1,3}}&b_{{1,4}}
\\ \noalign{\medskip}0&b_{{2,2}}&b_{{2,3}}&b_{{2,4}}
\\ \noalign{\medskip}0&0&b_{{3,3}}&b_{{3,4}}\\ \noalign{\medskip}0&0&0
&b_{{4,4}}\end {array} \right] 
 \left[ \begin {array}{cccc} 0&0&0&1\\ \noalign{\medskip}1&0&0&0
\\ \noalign{\medskip}0&1&0&0\\ \noalign{\medskip}0&0&1&0\end {array}
 \right]
=
\left[ \begin {array}{cccc} b_{{1,2}}&b_{{1,3}}&b_{{1,4}}&b_{{1,1}}
\\ \noalign{\medskip}b_{{2,2}}&b_{{2,3}}&b_{{2,4}}&0
\\ \noalign{\medskip}0&b_{{3,3}}&b_{{3,4}}&0\\ \noalign{\medskip}0&0&b
_{{4,4}}&0\end {array} \right] \in C_{2341}.
\]
\end{small}
\end{Example}
\bigskip

The \textit{Bruhat decomposition} of $\flags$ is the disjoint union of
the space into Schubert cells, each of which is a nonempty $B$-orbit:
\begin{equation}\label{eq:BruhatDecomposition}
\flags =  \bigcup_{w \in S_{n}} C_{w}(E_{\ci}).
\end{equation}
This decomposition was initially studied by Bruhat in the 1950s and
extended to other Lie groups $G$ and Borel subgroups $B$ \cite{Bruhat}.
In 2010, Lusztig wrote a nice overview of the contributions of Bruhat
in \cite{lusztig2010bruhat}, and says ``Bruhat decomposition is
indispensable for the understanding of both the structure and
representations of $G$.''

\begin{Definition}\label{def:Schubert.variety}
The \textit{Schubert variety} $X_{w}(E_{\ci})$ of a permutation $w$ is
defined to be the closure of $C_{w}(E_{\ci})$ under the Euclidean
topology, or equivalently the Zariski topology (see below). As in the case for Schubert cells, $X_{w}(E_{\ci})$ can be
expressed in terms of intersection dimensions as 
\begin{equation}\label{eq:def.Schubert.variety}
X_{w}(E_{\ci}) = \{F_{\ci} \in \flags
\given \mathrm{dim}(E_{i} \cap F_{j}) \geq \mathrm{rk}(w)[i,j] \text{ for all } 1\leq i,j \leq n \}.
\end{equation}
If the fixed flag is clear from context, we may write $X_{w}$ instead
of $X_{w}(E_{\ci}).$ 
\end{Definition}

\begin{Example}\label{ex:isomorphism.type.2314}
By definition $\mathrm{dim}(E_{2} \cap F_{1}) \geq 1$ and
$\mathrm{dim}(E_{3} \cap F_{3}) \geq 3$ so $F_{3}=E_{3}$ for all
$F_{\bullet}\in X_{2314}(E_{\ci})$.  These are the only binding
conditions on the flags in $X_{2314}(E_{\ci})$.  Therefore, every flag
in $X_{2314}(E_{\ci})$ is determined by a choice of the 1-dimensional
subpace $F_{1}$ in $E_{3}$.
\end{Example}

As described above, each Schubert cell is a $B$-orbit.  Similarly, one
can observe from \eqref{eq:def.Schubert.variety} that $BX_{w}=X_{w}$
since left multiplication by $B$ never decreases the northwest rank
table of a matrix. So $X_{w}$ is the union of a finite number of
$B$-orbits.

\begin{Example}\label{ex:1324.2341} By checking the intersection conditions, we see that
all flags in the Schubert cell $C_{1342}(E_{\ci})$ are in $X_{2341}(E_{\ci})$,
so
\begin{align*}
C_{1342}(E_{\ci}) = \left\{\left[\begin{array}{cccc}
 1 &	 0 &	 0 &	 0 \\
0 &	* &	 * &	 1\\
0 &     1 &     0 &	 0\\
0 &	0 &  1 &	 0 
\end{array} \right]  \right\}
\subset 
\displaystyle X_{2341}(E_{\ci}) = 
\overline{
\left\{
\left[\begin{array}{cccc}
* &	 * &	 * &	 1 \\
1 &	 0 &	 0 &	0\\
0 &      1 &     0 &	 0\\
0 &      0 &	 1 &	0
\end{array} \right]
 \right\}
},
\end{align*}
where each star represents a free parameter of the Schubert cell as
we saw with $x,y,z$ in \Cref{example:n=4.continued.2}.  We can also
see $C_{1342}$ is contained in the closure of $C_{2341}$ inherited
from the Euclidean topology as follows. For $a \neq 0$, the following
matrices all represent the same flag in $C_{2341}$:
\begin{equation*}
  \begin{bmatrix}
    a & 0 & 0 & 1 \\ 
    1 & 0 & 0 & 0\\
    0 & 1 & 0 & 0\\
    0 & 0 & 1 & 0
  \end{bmatrix} \sim 
  \begin{bmatrix}
    1 & 0 & 0 & 1 \\ 
    a^{-1} & 0 & 0 & 0\\
    0 & 1 & 0 & 0\\
    0 & 0 & 1 & 0
  \end{bmatrix} \sim \begin{bmatrix}
    1 & 0 & 0 & 0 \\ 
    a^{-1} & 0 & 0 & -a^{-1}\\
    0 & 1 & 0 & 0\\
    0 & 0 & 1 & 0
  \end{bmatrix}
\sim \begin{bmatrix}
    1 & 0 & 0 & 0 \\ 
    a^{-1} & 0 & 0 & 1\\
    0 & 1 & 0 & 0\\
    0 & 0 & 1 & 0
  \end{bmatrix}.
\end{equation*}
Letting $a \to \infty$ gives the permutation matrix $M_{1342}$, which
therefore lies in the (Euclidean) closure $\overline{{C}_{2341}}$. Since
 $X_{2341}=\overline{{C}_{2341}}$ is $B$-stable, $M_{1342} \in X_{2341}$ implies  $BM_{1342}=C_{1342} \subseteq \overline{{C}_{2341}} = X_{2341}$.
\end{Example}

The Zariski topology is used in algebraic geometry as the topology on
affine and projective varieties.  In that context, a variety
is defined as the set of solutions to a system of polynomial
equations, and such sets form the closed sets of a topology.  For instance, the statement of \Cref{ex:1324.2341} that $C_{1342} \subseteq \overline{{C}_{2341}}$ also holds when the closure is taken in the Zariski topology, where it means that every polynomial vanishing on ${C}_{2341}$ also vanishes on $C_{1342}$.

The equations defining a Schubert variety come from the southwest rank conditions and setting certain determinants of
submatrices equal to 0. However, note that a general polynomial function of the entries of a matrix $M$ is not well-defined when we view $M$ as an element of $GL_{n}/B$. We postpone further discussion of this issue and of the equations of a Schubert variety to \S\ref{sub:Grassmannians.intro}, and for now just view Schubert varieties as sets.

There was nothing special about the use of the flag $E_\bullet$ in
Definition~\ref{def:Schubert.variety}: for every pair of flags
$B_{\bullet}$ and $R_{\bullet }$ in $\flags$ we can consider the table
of intersection dimensions $\mathrm{dim}(B_{i}\cap R_{j})$.  This
table will again correspond with a permutation in $S_{n}$.  The
concept of a Schubert cell carries over as well, so
$C_{w}(B_{\bullet})$ is the set of all flags $R_{\bullet }$ in
position $w$ with respect to the base flag $B_{\bullet}$, equivalently
$\mathrm{dim}(B_{i}\cap R_{j})=\mathrm{rk}(w)[i,j]$ for all $i,j$.
Schubert varieties can be generalized similarly.  Changing the base
flag gives a convenient way to move a Schubert cell or Schubert
variety around inside $\flags$ while preserving most aspects of the
geometry.

The standard flag $E_{\bullet}=(e_{1},e_{2},\dots , e_{n})$ is most
commonly used.  The next most common flag is obtained by reversing the
standard basis, say $\oppositeE_{\bullet}=(e_{n},e_{n-1},\dots ,
e_{1})$.  Schubert varieties with respect to $\oppositeE_{\bullet}$
are called \textit{opposite Schubert varieties}.  You may see these
designated in the literature as $X^{w}$ or $\Omega_{w}$.  Every
ordered basis $(f_{1},\dots , f_{n})$ gives rise to a flag
$F_{\bullet}$ and an opposite flag with ordered basis $(f_{n},\dots ,
f_{1})$, which we could denote by $\oppositeF_{\bullet}$.  The key
observation is that $F_{\bullet}\in X_{w_{0}}(\oppositeF_{\bullet})$
and $\oppositeF_{\bullet}\in X_{w_{0}}(F_{\bullet})$  so opposite flags are in transverse position.

The position of a flag with respect to a given flag $F_{\bullet}$ can
be different from the position of the flag with respect to
$\oppositeF_{\bullet}$.  For example, the flag with ordered basis
$(e_{3}+7e_{4}, e_{2}-e_{4}, e_{1}+5e_{4},e_{4})$ is in position
$4321$ with respect to $E_{\bullet }$ and in position $2341$ with
respect to $\oppositeE_{\bullet}$.

\begin{Exercise}\label{ex:opposite.canonical}
Prove that a flag $G_{\bullet} \in C_{w}(\oppositeE_{\bullet})$ can be
represented by an ordered basis $(g_{1},g_{2},\dots , g_{n})$ such
that $g_{j}$ is in the span of $\{e_{n},e_{n-1},\dots , e_{n-w_{j}+1}
\}$ for each $j \in [n]$.  By considering the possible elementary
column operations, describe a canonical representation of a flag with
respect to $\oppositeE_{\bullet}=(e_{n},e_{n-1},\dots , e_{1})$.
\end{Exercise}

\begin{Exercise}\label{ex:transitive} 
Show $GL_{n}$ acts transitively on the complete flags in the flag variety
via left multiplication, so $\flags$ is a smooth manifold.
\end{Exercise}

\begin{Exercise}\label{ex:flags.4}
Return to the projective drawings we saw in \Cref{fig:rel.position}.
With respect to the black flag $B_{\bullet}$, which Schubert cell
$C_{w}(B_{\bullet})$ does each of the red flags lie in?  In the other
direction, try drawing out a typical flag in $C_{1432}(B_{\bullet})$
and $C_{4123}(R_{\bullet})$ using a shoebox drawing.  How many flags
are in $C_{1432}(B_{\bullet}) \cap C_{4123}(R_{\bullet})$ if
$B_{\bullet}$ and $R_{\bullet}$ are chosen generically?  If there are
any, draw one in the shoebox picture.  
\end{Exercise}

\begin{Exercise}\label{exericse:local.coords}
Prove that the reverse standard flag $\oppositeE_{\bullet}$ has an
affine neighborhood $C_{w_{0}}$ of dimension $\binom{n}{2}$ and a
local coordinate system.  Similarly, prove that a flag with canonical
matrix $g$ has an affine neighborhood $g w_{0}C_{w_{0}}$ in $\flags$.
\end{Exercise}

\subsection{Bruhat Order on Permutations}\label{sub:bruhat}

Since Schubert cells are $B$-orbits, Schubert varieties are invariant
under left multiplication by $B$ as well. So, each Schubert variety is
equal to a disjoint union of Schubert cells.  Which cells are
contained in $X_{w}$? The containment relation on Schubert varieties
$C_v \subset X_w$ defines a partial order on permutations $v\leq w$.
We write 
\begin{align}\label{prelim 100}
  X_{w} = \bigcup_{v\leq w} C_{v} 
\end{align}
similar to the Bruhat decomposition of $\flags$ in
\eqref{eq:BruhatDecomposition}.  This poset on $S_{n}$ is called
\textit{Bruhat order}.  One way to verify if $v\leq w$ is to compare
their rank tables by \eqref{eq:def.Schubert.variety}.  In particular,
$v\leq w$ in Bruhat order if and only if $\mathrm{rk}(v)[i,j]\geq
\mathrm{rk}(w)[i,j]$ for all $1\leq i,j\leq n$.  Note that before
Bruhat's work on Lie groups, this partial order on $S_{n}$ was studied
by Ehresmann around 1934 \cite{ehresmann.1934}, and generalized to
other Weyl groups by Chevalley in the 1950's \cite{chevalley.1958}.
Bruhat did not invent this partial order, but this naming convention
is prevalent in the literature and it is closely aligned with the
Bruhat decomposition.

Bruhat order has a nice description in terms of a set of generating
relations.  For a permutation $w \in S_{n}$ and integers $1\leq
i<j\leq n$, one can check by using the rank tables that $w(i)<w(j)$ if
and only if $w < w t_{ij}$.  On the other hand, if $v<w$ then there
exists a pair $c<d$ such that $v(c)<v(d)$ and
$\mathrm{rk}(vt_{cd})[i,j]\geq \mathrm{rk}(w)[i,j]$ for all $i,j$.
For example, take $c$ to be the first column where the rank tables of
$v$ and $w$ differ. Let $r$ be the first row where they differ in
column $c$.  Then $\mathrm{rk}(v)[r,c]> \mathrm{rk}(w)[r,c]$ since
$v<w$,\ $v(b)=w(b)$ for $1\leq b<c$,\ $v(c)=r$ and $w(c)>r$.  Since
$w(c)\not \in \{v(1),v(2),\dots ,v(c) \}$, we know there exists a
minimal $d>c$ such that $v(c)<v(d)$.  By considering rank tables
again, one can verify $v<vt_{cd}\leq w$.  Furthermore, $\ell(vt_{cd})
=\ell(v)+1$ since $vt_{cd}$ has exactly one more inversion than $v$.
By induction on the length difference $\ell(w)-\ell(v)$, there exists
a sequence of transpositions each moving up in Bruhat order
$v<vt_{c,d}<\dots \leq w$.  Hence, Bruhat order is equivalently
defined to be the transitive closure of the relations determined by
transpositions: $w < w t_{ij}$ for all $i<j$ and $w(i)<w(j)$.

The \textit{Hasse diagram} of a poset is a directed graph that encodes
a minimal set of generators for the poset under the transitive
closure.  The edges in the Hasse diagram are the \textit{covering
relations} for the poset.  In general, any directed acyclic graph
gives rise to a partial order by taking the transitive closure of its
edges, but the Hasse diagram is the unique minimal directed acyclic
graph generating the poset, also called its \textit{transitive
reduction}.

To draw the Hasse diagram for Bruhat order, note that not all of the
transposition relations, $w < w t_{ij}$ if $w(i)<w(j)$ and $i<j$, are
required.  For example $123<321$ in Bruhat order on $S_{3}$ and
$321=123\cdot t_{13}$, but we also have $123<213<231<321$, so the
Hasse diagram of Bruhat order for $S_{3}$ does not contain an edge
connecting $123$ and $321$.  It is a good exercise for the reader to
show that the minimal set of transposition relations also requires
$\ell(wt_{ij})=\ell (w)+1$.  Thus, Bruhat order is a ranked poset with
rank function given by the number of inversions of a permutation.  By
MacMahon's theorem, we know the rank generating function is given by
the standard $t$-analog of $n!$, namely
\[\sum_{w\in S_{n}} t^{\inv(w)}= \sum_{w\in S_{n}} t^{\ell(w)} = \prod_{k=1}^{n-1}(1+t+t^{2}+\dots +t^{k}) = [n]_{t}!.
\]
The polynomial $[n]_{t}!$ is also called the \textit{length generating
function} of $S_{n}$.

\begin{Example}
The following is the Hasse diagram of the Bruhat order on $S_{3}$.
There is an edge connecting $u$ and $w$ if $u\leq w$ and no $v$ exists
such that $u<v<w$.  We omit the arrows in this directed graph because
they are implied by the height of the elements.  A larger element is
above a smaller element in a Hasse diagram.

\[\begin{tikzpicture}[scale=0.6]
\def\a{0.7};
\def\b{0.4};
\def\h{2.4};
\newcommand\Rec[3]{
\node at (#1,#2) {#3};
\draw(#1-\a,#2-\b)--(#1-\a,#2+\b)--(#1+\a,#2+\b)--(#1+\a,#2-\b)--(#1-\a,#2-\b);
}
\Rec{0}{0}{$123$}
\Rec{-3*\a}{\h}{$132$}
\Rec{3*\a}{\h}{$213$}
\Rec{-3*\a}{2*\h}{$231$}
\Rec{3*\a}{2*\h}{$312$}
\Rec{0}{3*\h}{$321$}
\draw(0,\b)--(-3*\a,\h-\b);
\draw(0,\b)--(3*\a,\h-\b);
\draw(-3*\a,\h+\b)--(-3*\a,2*\h-\b);
\draw(-3*\a,\h+\b)--(3*\a,2*\h-\b);
\draw(3*\a,\h+\b)--(-3*\a,2*\h-\b);
\draw(3*\a,\h+\b)--(3*\a,2*\h-\b);
\draw(-3*\a,2*\h+\b)--(0,3*\h-\b);
\draw(3*\a,2*\h+\b)--(0,3*\h-\b);
\end{tikzpicture}\]
%\vspace{-0.1in}
%\setlength{\unitlength}{1cm}
%\newcommand{\perm}[1]{{\raisebox{-.1\unitlength}{\makebox(0,0)[b]{${#1}$}}}}
%\begin{center}
%\raisebox{1ex}{\begin{picture}(0,2.3)
%\put(1.5,0){\perm{132}}\put(1.5,1){\perm{231}}
%\put(.5,-1){\perm{123}}\put(.5,1.8){\perm{321}}
%\put(-.6,0){\perm{213}}\put(-.6,1){\perm{312}}
%\put(0,0){\line(0,1){1}}
%\put(1,0){\line(0,1){1}}
%\put(0,0){\line(1,1){1}}
%\put(.5,1.5){\line(1,-1){.5}}
%\put(.5,1.5){\line(-1,-1){.5}}
%\put(.5,-.5){\line(-1,1){.5}}
%\put(.5,-.5){\line(1,1){.5}}
%\put(0,1){\line(1,-1){1}}
%\end{picture}}
%\end{center}
%\vspace{.4in}
\end{Example}

\begin{Example}
The Hasse diagram of $S_{4}$ is drawn in Figure \ref{Hasse of S4}.
\begin{figure}[h]
\begin{tikzpicture}[scale=0.6]
\def\a{0.9};
\def\b{0.4};
\def\h{3.0};
\newcommand\Rec[3]{
\node at (#1,#2) {#3};
\draw(#1-\a,#2-\b)--(#1-\a,#2+\b)--(#1+\a,#2+\b)--(#1+\a,#2-\b)--(#1-\a,#2-\b);
}
\Rec{0}{0}{$1234$}
\Rec{0}{\h}{$1324$}
\Rec{-4*\a}{\h}{$1243$}
\Rec{4*\a}{\h}{$2134$}
\Rec{-8*\a}{2*\h}{$1342$}
\Rec{-4*\a}{2*\h}{$1423$}
\Rec{0*\a}{2*\h}{$2314$}
\Rec{4*\a}{2*\h}{$2143$}
\Rec{8*\a}{2*\h}{$3124$}
\Rec{-10*\a}{3*\h}{$1432$}
\Rec{-6*\a}{3*\h}{$2341$}
\Rec{-2*\a}{3*\h}{$2413$}
\Rec{2*\a}{3*\h}{$3142$}
\Rec{6*\a}{3*\h}{$3214$}
\Rec{10*\a}{3*\h}{$4123$}
\Rec{-8*\a}{4*\h}{$2431$}
\Rec{-4*\a}{4*\h}{$3412$}
\Rec{0*\a}{4*\h}{$3241$}
\Rec{4*\a}{4*\h}{$4132$}
\Rec{8*\a}{4*\h}{$4213$}
\Rec{-4*\a}{5*\h}{$3421$}
\Rec{0*\a}{5*\h}{$4231$}
\Rec{4*\a}{5*\h}{$4312$}
\Rec{0}{6*\h}{$4321$}
\draw(0,\b)--(-4*\a,\h-\b);
\draw(0,\b)--(0*\a,\h-\b);
\draw(0,\b)--(4*\a,\h-\b);
\draw(-4*\a,\h+\b)--(-8*\a,2*\h-\b);
\draw(-4*\a,\h+\b)--(-4*\a,2*\h-\b);
\draw(-4*\a,\h+\b)--(4*\a,2*\h-\b);
\draw(0*\a,\h+\b)--(-8*\a,2*\h-\b);
\draw(0*\a,\h+\b)--(-4*\a,2*\h-\b);
\draw(0*\a,\h+\b)--(0*\a,2*\h-\b);
\draw(0*\a,\h+\b)--(8*\a,2*\h-\b);
\draw(4*\a,\h+\b)--(0*\a,2*\h-\b);
\draw(4*\a,\h+\b)--(4*\a,2*\h-\b);
\draw(4*\a,\h+\b)--(8*\a,2*\h-\b);
\draw(-8*\a,2*\h+\b)--(-10*\a,3*\h-\b);
\draw(-8*\a,2*\h+\b)--(-6*\a,3*\h-\b);
\draw(-8*\a,2*\h+\b)--(2*\a,3*\h-\b);
\draw(-4*\a,2*\h+\b)--(-10*\a,3*\h-\b);
\draw(-4*\a,2*\h+\b)--(-2*\a,3*\h-\b);
\draw(-4*\a,2*\h+\b)--(10*\a,3*\h-\b);
\draw(0*\a,2*\h+\b)--(-6*\a,3*\h-\b);
\draw(0*\a,2*\h+\b)--(-2*\a,3*\h-\b);
\draw(0*\a,2*\h+\b)--(6*\a,3*\h-\b);
\draw(4*\a,2*\h+\b)--(-6*\a,3*\h-\b);
\draw(4*\a,2*\h+\b)--(-2*\a,3*\h-\b);
\draw(4*\a,2*\h+\b)--(2*\a,3*\h-\b);
\draw(4*\a,2*\h+\b)--(10*\a,3*\h-\b);
\draw(8*\a,2*\h+\b)--(2*\a,3*\h-\b);
\draw(8*\a,2*\h+\b)--(6*\a,3*\h-\b);
\draw(8*\a,2*\h+\b)--(10*\a,3*\h-\b);
\draw(-10*\a,3*\h+\b)--(-8*\a,4*\h-\b);
\draw(-10*\a,3*\h+\b)--(-4*\a,4*\h-\b);
\draw(-10*\a,3*\h+\b)--(4*\a,4*\h-\b);
\draw(-6*\a,3*\h+\b)--(-8*\a,4*\h-\b);
\draw(-6*\a,3*\h+\b)--(0*\a,4*\h-\b);
\draw(-2*\a,3*\h+\b)--(-8*\a,4*\h-\b);
\draw(-2*\a,3*\h+\b)--(-4*\a,4*\h-\b);
\draw(-2*\a,3*\h+\b)--(8*\a,4*\h-\b);
\draw(2*\a,3*\h+\b)--(-4*\a,4*\h-\b);
\draw(2*\a,3*\h+\b)--(0*\a,4*\h-\b);
\draw(2*\a,3*\h+\b)--(4*\a,4*\h-\b);
\draw(6*\a,3*\h+\b)--(-4*\a,4*\h-\b);
\draw(6*\a,3*\h+\b)--(0*\a,4*\h-\b);
\draw(6*\a,3*\h+\b)--(8*\a,4*\h-\b);
\draw(10*\a,3*\h+\b)--(4*\a,4*\h-\b);
\draw(10*\a,3*\h+\b)--(8*\a,4*\h-\b);
\draw(-8*\a,4*\h+\b)--(-4*\a,5*\h-\b);
\draw(-8*\a,4*\h+\b)--(0*\a,5*\h-\b);
\draw(-4*\a,4*\h+\b)--(-4*\a,5*\h-\b);
\draw(-4*\a,4*\h+\b)--(4*\a,5*\h-\b);
\draw(0*\a,4*\h+\b)--(-4*\a,5*\h-\b);
\draw(0*\a,4*\h+\b)--(0*\a,5*\h-\b);
\draw(4*\a,4*\h+\b)--(0*\a,5*\h-\b);
\draw(4*\a,4*\h+\b)--(4*\a,5*\h-\b);
\draw(8*\a,4*\h+\b)--(0*\a,5*\h-\b);
\draw(8*\a,4*\h+\b)--(4*\a,5*\h-\b);
\draw(-4*\a,5*\h+\b)--(0*\a,6*\h-\b);
\draw(0*\a,5*\h+\b)--(0*\a,6*\h-\b);
\draw(4*\a,5*\h+\b)--(0*\a,6*\h-\b);
\end{tikzpicture}
\caption{The Hasse diagram of $S_{4}$}\label{Hasse of S4}
\end{figure}
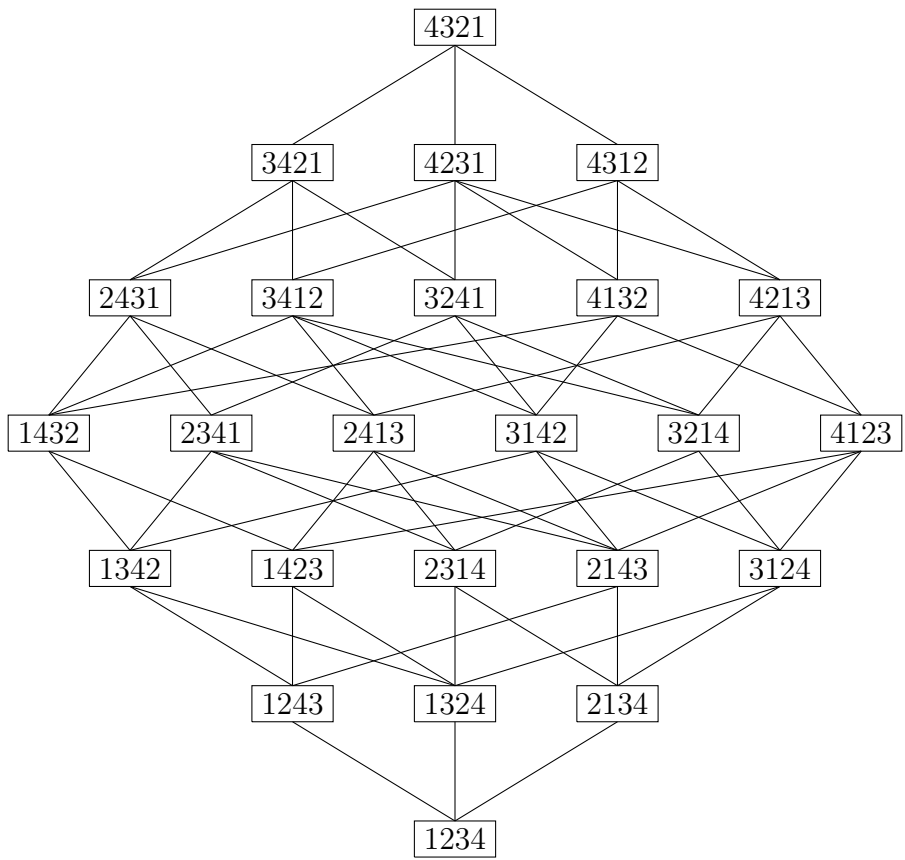
\end{Example}

\begin{Exercise}
Prove that the Hasse diagram of $S_{n}$ is self-dual, rank-symmetric and
rank-unimodal.
\end{Exercise}

Testing if $v\leq w$ in Bruhat order via the transposition relations
can be cumbersome since one would need to consider many possible paths
from $v$ to $w$.  A more efficient criteria was given by Ehresmann
\cite{ehresmann.1934} and improved upon by Bj\"orner and Brenti \cite{BB2}.  This
criteria depends on another representation of a permutation $w
=[w_{1},\dots , w_{n}]$ as a tableau (table in English) of numbers
where each row is the sorted list of initial values $\{w_{1},\dots ,
w_{i} \}$.  Since we only care about the numbers in the rows and not
any particular order, one might left justify or right justify such a
tableau as we have done below in an example.  The criteria then
compares $v$ and $w$ in Bruhat order by comparing the rows of the
tableau in what is sometimes now referred to as the Gale order \cite{Gale}
because of his studies in matroid theory in the 1960s, well after
Ehresmann.
 
\begin{Definition}\label{def:dominated.multisets.different.size}(\textbf{Gale order}) If $A=\{a_1\leq\cdots\leq a_p\}$ and $B=\{b_1 \leq \cdots \leq
b_p\}$ are subsets of $[n]$ with their elements listed in increasing
order, we say $A \trianglelefteq B$ if and only if $a_{i}\leq b_{i}$
for all $1\leq i\leq p.$
\end{Definition}

\begin{Theorem}[\textbf{The Ehresmann Tableau Criterion}]\cite{ehresmann.1934} \label{thm:Ehresmann}
For $v, w \in S_n$, we have $v \leq w$ in Bruhat order if and only if
$\{v_1,\ldots,v_j\} \trianglelefteq \{w_1,\ldots,w_j\}$ for each $j
\in [n]$.\end{Theorem}

\begin{Theorem} \cite[Cor. 5]{BB2} \label{thm:bb.improved.criterion}
For $v, w \in S_n$, we have $v \leq w$ in Bruhat order if and only if
$\{v_1,\ldots,v_j\} \trianglelefteq \{w_1,\ldots,w_j\}$ whenever $v_j
> v_{j+1}$.\end{Theorem}

Continuing \Cref{ex:1324.2341}, we can test $1324< 2341$ by checking
the tableau conditions

\[
\begin{array}{rcl}
1 & \trianglelefteq & 2\\
1 3 & \trianglelefteq & 2 3\\
1 2 3 & \trianglelefteq & 2 3 4\\
1 2 3 4 & \trianglelefteq & 1 2 3 4.
\end{array}
\]
However, $1423 \not \leq 2314$ since $14 \not \trianglelefteq 23$.

The proof of Ehresmann's tableau criterion follows from the rank
tables and intersection conditions for Schubert varieties.  We can
think of the rows of the tableau for $w$ as another way of stating the
defining rank conditions on the flags in $C_{w}(E_{\bullet})$.  By
construction, row $j$ contains $\{w_{1},\dots , w_{j} \}$ in
increasing order, say $\{w_{1},\dots , w_{j} \}=\{a_{1}<a_{2}<\dots <
a_{j} \}$. Considering the canonical form matrix representing a flag
$F_{\bullet} \in C_{w}(E_{\bullet})$, one can observe that
$\mathrm{dim}(E_{i} \cap F_{j})=0$ for $i<a_{1}$, the dimension jumps
up to $\mathrm{dim}(E_{i} \cap F_{j})=1$ for $a_{1}\leq i<a_{2}$, and
then jumps again $\mathrm{dim}(E_{i} \cap F_{j})=2$ for $a_{2}\leq
i<a_{3}$, etc.  The $a_{i}$'s determine the ``jumps'' in dimension as
we consider the list of intersections of $F_{j}$ with
$E_{1},E_{2},\dots , E_{n}$.  Some authors refer to these required
jumps as \textit{Schubert conditions}.  Ehresmann's criterion just
encodes the fact that for each $j$, the jumps for $F_{j}$ must come at
least as early in $v$ as in $w$ in order to have $v\leq w$.  The
Bj\"orner-Brenti improvement (Theorem~\ref{thm:bb.improved.criterion})
encodes the fact that the rows $j$ for which $v_{j}>v_{j+1}$ are more
binding than the others, so the Schubert conditions from other rows
follow from these.

One of the benefits of Bruhat order is a description of the
\textit{Poincar\'e polynomials} of Schubert varieties.  Since a
Schubert variety is the disjoint union of Schubert cells, there is a
way to obtain a basis for the associated cohomology ring of $X_{w}$,
denoted $H^{*}(X_{w})$, indexed by the permutations $v\leq w$.  See
\Cref{sub:Monk} for more discussion on cohomology rings.  For now, we
just note that the cohomology ring of $X_{w}$ is a graded ring and its
Hilbert series is the same as the Poincar\'e polynomial of $X_{w}$,
namely
\begin{align*}
P_{w}(t)= \sum_{v \leq w} t^{2 \ell (v)}.
\end{align*}
Because only even exponents appear in the Poincar\'e polynomials
above, we often abuse notation and define 
\begin{align*}
P_{w}(t)= \sum_{v \leq w} t^{\ell (v)}
\end{align*}
using complex dimension instead of real dimension.

\begin{Example}\label{ex:3412.Hasse.diagram}
For $w=3412$, the following permutations shown in
\Cref{fig:3412-lower-interval} are in the interval below $3412$ in
Bruhat order, denoted $[\id,3412]=\{v\in S_{n} \given v\leq 3412 \}$.
The number on the left is the number of inversions for all
permutations on that row.
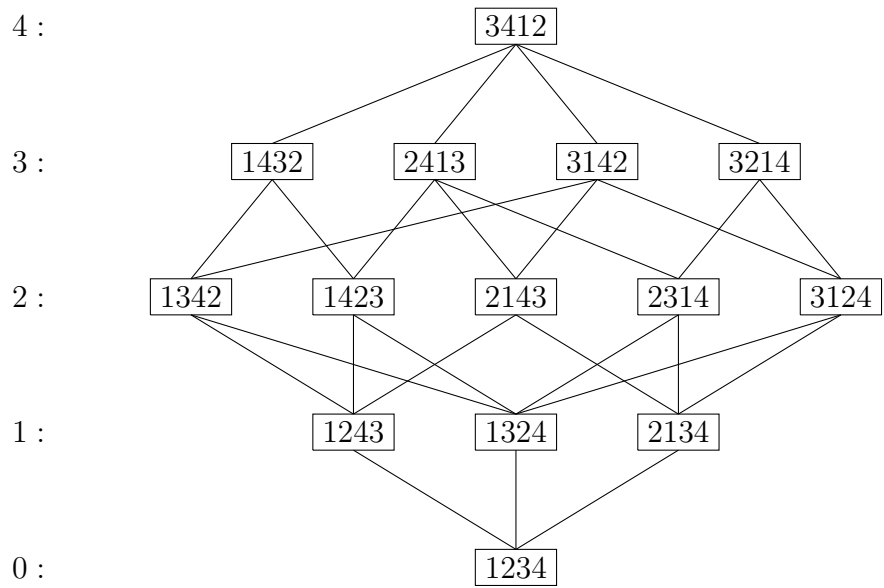
\begin{figure}[h!]
\centering
\begin{tikzpicture}[scale=0.6]
\def\a{0.9};
\def\b{0.4};
\def\h{3.0};
\newcommand\Rec[3]{
\node at (#1,#2) {#3};
\draw(#1-\a,#2-\b)--(#1-\a,#2+\b)--(#1+\a,#2+\b)--(#1+\a,#2-\b)--(#1-\a,#2-\b);
}
\Rec{0}{0}{$1234$}
\Rec{0}{\h}{$1324$}
\Rec{-4*\a}{\h}{$1243$}
\Rec{4*\a}{\h}{$2134$}
\Rec{-8*\a}{2*\h}{$1342$}
\Rec{-4*\a}{2*\h}{$1423$}
\Rec{0*\a}{2*\h}{$2143$}
\Rec{4*\a}{2*\h}{$2314$}
\Rec{8*\a}{2*\h}{$3124$}
\Rec{-6*\a}{3*\h}{$1432$}
\Rec{-2*\a}{3*\h}{$2413$}
\Rec{2*\a}{3*\h}{$3142$}
\Rec{6*\a}{3*\h}{$3214$}
\Rec{0*\a}{4*\h}{$3412$}

\draw(0,\b)--(-4*\a,\h-\b);
\draw(0,\b)--(0*\a,\h-\b);
\draw(0,\b)--(4*\a,\h-\b);
\draw(-4*\a,\h+\b)--(-8*\a,2*\h-\b);
\draw(-4*\a,\h+\b)--(-4*\a,2*\h-\b);
\draw(-4*\a,\h+\b)--(0*\a,2*\h-\b);
\draw(0*\a,\h+\b)--(-8*\a,2*\h-\b);
\draw(0*\a,\h+\b)--(-4*\a,2*\h-\b);
\draw(0*\a,\h+\b)--(4*\a,2*\h-\b);
\draw(0*\a,\h+\b)--(8*\a,2*\h-\b);
\draw(4*\a,\h+\b)--(0*\a,2*\h-\b);
\draw(4*\a,\h+\b)--(4*\a,2*\h-\b);
\draw(4*\a,\h+\b)--(8*\a,2*\h-\b);
\draw(-8*\a,2*\h+\b)--(-6*\a,3*\h-\b);
\draw(-8*\a,2*\h+\b)--(2*\a,3*\h-\b);
\draw(-4*\a,2*\h+\b)--(-6*\a,3*\h-\b);
\draw(-4*\a,2*\h+\b)--(-2*\a,3*\h-\b);
\draw(0*\a,2*\h+\b)--(-2*\a,3*\h-\b);
\draw(0*\a,2*\h+\b)--(2*\a,3*\h-\b);
\draw(4*\a,2*\h+\b)--(-2*\a,3*\h-\b);
\draw(4*\a,2*\h+\b)--(6*\a,3*\h-\b);
\draw(8*\a,2*\h+\b)--(2*\a,3*\h-\b);
\draw(8*\a,2*\h+\b)--(6*\a,3*\h-\b);
\draw(-6*\a,3*\h+\b)--(0*\a,4*\h-\b);
\draw(-2*\a,3*\h+\b)--(0*\a,4*\h-\b);
\draw(2*\a,3*\h+\b)--(0*\a,4*\h-\b);
\draw(6*\a,3*\h+\b)--(0*\a,4*\h-\b);

\node at (-12*\a,0*\h) {$0:$};
\node at (-12*\a,1*\h) {$1:$};
\node at (-12*\a,2*\h) {$2:$};
\node at (-12*\a,3*\h) {$3:$};
\node at (-12*\a,4*\h) {$4:$};
\end{tikzpicture}
\caption{The lower interval in the Bruhat order below $3412$}
\label{fig:3412-lower-interval}
\end{figure}
So $P_{3412}(t) = 1 + 3t + 5 t^{2} + 4 t^{3} +t^{4}$. One can see that
the Schubert variety $X_{3412}$ is not a smooth manifold since its
Poincar\'e polynomial is not symmetric (palindromic) which implies
that Poincar\'e duality does not hold for $H^*(X_{3412})$.
\end{Example}
\bigskip

There are many additional interesting facts about Bruhat order as a
partial order on $S_{n}$.  Here are some examples. We will encounter
some more of them in the rest of this book.  See \cite[Ch 2]{b-b} and
\cite{Hum} for more background and generalizations to all Coxeter groups.

%%%%%%
\begin{enumerate}
\item  Bruhat order is $k$-Sperner for all $k$ \cite{S7}.
\item  It has the nicest possible M\"obius function:  $\mu (v,w)
=(-1)^{\ell(w)-\ell(v)}$ \cite{Verma.71}.
\item The intervals $[x,y]$ determine the composition series for
Verma modules \cite{Ve1}.
\item  The interval $[\id,w]$ is rank symmetric if and only if $X_{w}$
is nonsingular as a variety \cite{carrell94}.  See also
\Cref{sub:Singularloci}. 
\item The order complex of $(u,v)$ is shellable \cite{b-w-82}.
\item  It is amenable to pattern avoidance \cite{Abe-Billey}.
\end{enumerate}

\bigskip

\begin{Exercise}\label{ex:order-reversing}
Let $v,w \in S_{n}$.  The following are equivalent:
\begin{enumerate}
\item $v\leq w$,
\item $v^{-1}\leq w^{-1}$,
\item $ww_{0}\leq vw_{0}$,
\item $w_{0}w \leq w_{0}v$.
\end{enumerate}
\end{Exercise}

\begin{Exercise}\label{ex:boundary} Prove that the boundary of $X_{w}$
has irreducible components given by the Schubert varieties $X_v$ such
that $v$ is covered by $w$ in Bruhat order.  Show $w$ covers $v$ if
and only if $v=wt_{ij} <w$ for some $i<j$ with $w_{i}>w_{j}$ and $\ell
(v)=\ell (w)-1$.
\end{Exercise}

\begin{Exercise}\label{ex:dense}
Prove that $C_{w}(F_{\ci})$ is a dense open set in $X_{w}(F_{\ci})$
for any $F_{\ci} \in \flags$.
\end{Exercise}

\begin{Exercise}\label{ex:flatten.bruhat} \cite[Lemma 2.1]{B3} Given
any sequence of distinct real numbers $r_1\ldots r_m$, define the
\textit{flattening function} $\fl(r_1\ldots r_m)$ to be the permutation
$v \in S_m$ such that $r_i < r_j$ if and only if $v_i<v_j$.  Recall
that a permutation is uniquely defined by its inversion set, so this
condition uniquely defines $v$.  Prove that if two permutations $v,w
\in S_{n}$ agree in position $i$, then $v \leq w$ in Bruhat order on
$S_{n}$ if and only if $\fl(v_1\ldots \widehat{v_i} \ldots v_n) \leq
\fl(w_1\ldots \widehat{w_i} \ldots w_n)$ in Bruhat order on $S_{n-1}$.
\end{Exercise}

\begin{Exercise}
What is the Poincar\'e polynomial for $H^{*}(\flags)$ for $n=3,4,5$?
\end{Exercise}

\begin{Exercise}\label{ex:richardson.varieites}
Let $\oppositeE_{\bullet}=(e_{n},\dots , e_{1})$. Which permutation
matrices represent flags in $X_{w}(\oppositeE_{\bullet})$ for $w \in S_{n}$?  
\end{Exercise}

\begin{Exercise}\label{ex:Subword.Property}    
For $v,w \in S_{n}$, show $v\leq w$  if and only if some reduced word
for $v$ is a subword of some (any) reduced word for $w$. 
\end{Exercise}

\begin{Exercise}\label{ex:bigrass.bruhat}
Recall the bigrassmannian permutations defined in
\Cref{ex:bigrassmannians}.  Prove that the bigrassmannian permutations
are the join-irreducible elements in Bruhat order, and so $v \leq w$
in Bruhat order if and only if the set of all bigrassmannian
permutations below $v$ is a subset of the bigrassmannian permutations
below $w$ \cite{LS-MacNeille}.
\end{Exercise}

\begin{Exercise}\cite[Prop. 4.2]{RWY}\label{ex:bigrass.not.below}
For $w \in S_{n}$, prove that the minimal elements in $\{y \in S_{n}
\given y \not \leq w \}$ are all bigrassmannian permutations.  Can you
explicitly determine which bigrassmannian permutations are the minimal
elements in the complement of $[\mathrm{id},w]$?
%%%% answer: the ones corresponding with essential set elements.
%%%% Namely, if (i,j) is in the essential set of w  and
%%%% $rk_{w}(s,t)=r$, then the bigrassmmannian perm
%%%% (1,2,\dots, r-1,
%%%             s+1, s+2, \dots, s+t-r+1,
%%%                      r, r+1, \dots,s,
%%%                            s+t-r+2, \dots , n)
%%% is minimal in the complement of the interval and conversely.
\end{Exercise}

\begin{Problem}\label{prob:enumerate.pairs}
If $v,w$ are chosen uniformly in $S_{n}$, what is the probability that
$v\leq w$?  See \cite{hammett2006likelihood}.
\end{Problem}

Bj\"orner showed that each interval $[u,w]=\{v \in S_{n}\given u\leq v\leq
w \}$ of rank $k$ in Bruhat order is the order poset for a CW-complex
on the sphere of dimension $k-2$ \cite{Bjorner.1984}. For $k=2$, each
interval is a diamond.  For $k=3$, each interval is a $k$-crown.
Hultman identified each possible interval type for $k=4$
\cite{hultman.4b}. For each $k$, there are only a finite number of
possible interval types according to a theorem of Dyer
\cite{Dyer.1991.Bruhat.graph}.

\begin{Problem}\label{Problem:Bruhat}
Can the set of all possible Bruhat intervals be characterized for
arbitrary $k$?  Conversely, given an abstract poset $P$ with a unique
minimal and maximal element, what is the computational complexity of
determining if $P$ is isomorphic to an interval in Bruhat order in
$S_{n}$ for some $n$?
\end{Problem}

The following conjecture was made independently by Dyer and Lusztig,
see also  Brenti-Casselli-Martinelli \cite{BCM.2006}.  The Kazhdan-Lusztig
polynomials are an important family of polynomials in $\mathbb{Z}[q]$
which are indexed by two permutations $u\leq v$ in Bruhat order.  We
will give their classical definition via an existence statement in
\Cref{thm:KL.polynomials}.  More generally, these polynomials can be
defined for any two elements in a Coxeter group related by Bruhat
order.

\begin{Problem}\label{KL-intervals}\textbf{The Combinatorial
Invariance Conjecture}. The Kazhdan-Lusztig polynomials $P_{u,v}(q) $
and $P_{w,x}(q) $ are equal if the intervals $[u,v]$ and $[w,x]$ are
isomorphic as abstract posets.
\end{Problem}

The converse of this statement does not hold. For example
$P_{123,321}=1=P_{123,123}$ but the corresponding intervals are
different.

\subsection{Grassmannians and Partial Flag Varieties}\label{sub:Grassmannians.intro}

Where did the notion of a complete flag come from?  This was in some
sense a natural extension of projective space.  The \textit{projective
space} of $\mathbb{C}^{n}$, denoted $\mathbb{P}^{n-1}$, is the
collection of complex lines in $\mathbb{C}^{n}$ through the origin.  Each
such line can be represented by a nonzero $n$-dimensional vector.
Rescaling the vector by a nonzero complex number results in another
representation of the same line. This space has dimension $n-1$, hence
the superscript on $\mathbb{P}^{n-1}$.  A point in $\mathbb{P}^{n-1}$
is represented by \emph{homogeneous coordinates} $[p_{1}:p_2: \dots :p_{n}]$
where $[p_{1}:p_2: \dots :p_{n}]=[\lambda p_{1}:\lambda p_2: \dots
:\lambda p_{n}]$ for all $\lambda \in \mathbb{C}$ such that $\lambda
\neq 0.$ We call $\mathbb{P}^{1}$ the projective line, and
$\mathbb{P}^{2}$ is the projective plane. One can think of the
projective plane as an open disk (points with $p_2 \neq 0$), plus the circle on the boundary
minus a point (points with $p_2 = 0, p_1 \neq 0$), plus that last point $[1:0:0]$.  This gives a cell decomposition
of $\mathbb{P}^{2}$, similar to the cell decomposition of the flag
variety into Schubert cells.

Projective spaces are better behaved than complex spaces in some ways.  A finite
dimensional projective space is compact.  It can be covered by open
charts, so a projective space is a manifold.  Furthermore, any two
subspaces of complementary dimension intersect,
perhaps at $\infty$.  For example, two parallel lines in a plane meet
somewhere in the closure of the plane, ``off at $\infty$''.   

Why stop at 1-dimensional subspaces? The \textit{Grassmannian
variety} $\Gr(k,n)$ is the collection of $k$-dimensional subspaces of
$\mathbb{C}^{n}$.  Such subspaces can be represented by $k$ linearly
independent vectors, or equivalently by a $n \times k$ matrix of full
rank $k$.  We can get a canonical representation of each $k$-dimensional
subspace by using the reduced column echelon form of any matrix $M$
representing it.  The determinantal minors of a full rank $n \times k$
matrix $M$ can be used to embed $\Gr(k,n)$ into a projective space as
follows.  For a size $k$ subset $I \subset [n]$, let $\Delta_{I}(M)$
represent the determinant of the submatrix of $M$ in rows $I$ and
columns $[k]$.  If we organize the size $k$ minors of $M$ into a list
\begin{equation}\label{eq:plucker.coords.gnk}
P_{k}(M)=[\Delta_{I_{1}}(M):\cdots :\Delta_{I_{\binom{n}{k}}}(M)]
\end{equation}
where the subsets of size $k$ appear in lexicographic order, we get a
homogeneous coordinate for some point in $\mathbb{P}^{\binom{n}{k}-1}$,
called the \textit{Pl\"ucker coordinates}.  Two matrices $M$ and $N$
represent the same $k$-dimensional subspace $V \in \Gr(k,n)$ if and
only if $N=MA$ for some $A \in GL_{k}$, in which case the homogeneous
coordinates $P_{k}(M)$ and $P_{k}(N)$ agree because multiplication by $A$ just
rescales the homogeneous coordinate $P_{k}(M)$ by the determinant of $A$.
Therefore, we call $P_{k}(M)$ the \textit{Pl\"ucker coordinates} for $V$.
In fact, one can recover the reduced column echelon form of the matrix
representing $V$ from $P_{k}(M)$, so the Pl\"ucker coordinates give an
embedding of the Grassmannian variety $\Gr(k,n)$ into the projective
space $\mathbb{P}^{\binom{n}{k}-1}$.  The image of the Pl\"ucker
embedding is all points in $\mathbb{P}^{\binom{n}{k}-1}$ satisfying
certain quadratic polynomials called the \textit{Pl\"ucker relations}.  The
Pl\"ucker relations play a central role in the theory of
Grassmannians, $S_{n}$ representation theory, and $GL_{n}$
representation theory.  Fulton's book \emph{Young Tableaux} does a
wonderful job of telling this story so we refer the reader to
\cite{Fulton-book} for more details.  We will return to these
relations in \Cref{sec:CAGofSchubertVarieties}.  

Recall from \Cref{def:flag} that a \emph{partial flag} with dimensions
${\bf d}$ is a sequence of subspaces $F_1 \subseteq \cdots \subseteq
F_m \subseteq \C^n$ with $\dim F_i = d_i$ for any given subset ${\bf
d} = \{d_1 < \cdots < d_m\} \subseteq [n-1]$.  The set of all such
partial flags is the \emph{partial flag variety} $\Fl(n; {\bf d})$.
For example, the Grassmannian is $\Gr(k,n) = \Fl(n; \{k\})$.  Using
minors of different sizes, one can give $\Fl(n; {\mathbf d})$ the
structure of a compact smooth manifold or a projective complex variety
in more or less the same way as $\Gr(k,n)$ or $\Fl(n)$.  Some authors
say ``Grassmannian manifold'' or ``flag manifold'' instead of
Grassmannian variety and flag variety when they wish to emphasize the
manifold structure with local coordinate charts similar to what we saw
in \Cref{exericse:local.coords}. Note also that $\GL_n$ still acts
transitively on $\Fl(n; {\mathbf d})$ by left multiplication, as in
Exercise~\ref{ex:transitive}.  Furthermore, $\Fl(n; {\mathbf d})$ is
isomorphic to a quotient $\GL_n/P$ where $P$ is the set of invertible,
block upper triangular matrices with zeros strictly southwest of the
entries $(d_{i},d_{i})$ for $1\leq i\leq m$.

There is an obvious projection map $\pi : \Fl(n) \to \Fl(n; {\mathbf
d})$ which ``forgets'' those components of a complete flag with
dimensions not in ${\mathbf d}$. It can be shown that the partial flag
varieties $\Fl(n; {\mathbf d})$ account for \emph{all} of the compact
quotients of $\GL_n$, and one reason for the ubiquity of the complete flag variety is that it is maximal among them, in the sense that it projects to all others via these projections $\pi$. In particular, Schubert varieties in the flag variety project to varieties in $\Fl(n; {\mathbf d})$. The image
of such a projection in a partial flag variety or Grassmannian is also
called a \textit{Schubert variety}.  So, take some care to note where
your Schubert varieties live.  In this chapter, we will stay focused
on Schubert varieties in the complete flag variety as they present a very rich
structure from which most results on partial flag varieties can be
deduced.

We now return to the topic of equations defining Schubert varieties,
as promised in \S\ref{sub:Schubert.Varieties}. Generalizing the
notation used for Grassmannians, let $\Delta_{I,J}(M)$ be the
determinant of the submatrix of $M$ in rows $I$ and columns $J$ for
$I,J \subset [n]$ of the same size.  All such minors $\Delta_{I,J}$
can be expressed as homogeneous polynomial functions in the entries of
a matrix using variables $\{z_{ij}\given 1\leq i,j\leq n \}$.  A
\emph{flag minor} is one of the form $\Delta_{I,[j]}$ for $j \leq
n$. Viewing $\flags$ as a subset of $\Gr(1,n) \times \Gr(2,n) \times
\cdots \times \Gr(n-1,n)$ and then applying the Pl\"ucker embedding
for Grassmannians gives an embedding $\flags \hookrightarrow \bP^{{n
\choose 1}-1} \times \bP^{{n \choose 2}-1} \times \cdots \times \bP^{{n-1
\choose 1}-1}$. The flag given by a matrix $M$ embeds as a point with
coordinates on the $j$\textsuperscript{th} factor of this product
given by the flag minors $\Delta_{I,[j]}(M)$ for all $j$-subsets $I
\subseteq [n]$. If $f$ is a polynomial in the flag minors which is homogeneous in a
suitable sense, then the equation $f = 0$ defines a subset of $\bP^{{n
\choose 1}-1} \times \bP^{{n \choose 2}-1} \times \cdots \times
\bP^{{n-1 \choose 1}-1}$, which is itself a projective variety via the
Segre embedding (see Exercise~\ref{ex:segre}).  We will return to projective varieties in
\Cref{sec:CAGofSchubertVarieties}.

\begin{Exercise}\label{ex:segre}
Show that the image of the \emph{Segre embedding} $\mathbb{P}^{m}
\times \mathbb{P}^{n} \longrightarrow \mathbb{P}^{(n+1)(m+1)-1}$
defined by $[x_{0}:\cdots:x_{m}][y_{0}:\dots :y_{n}]\to [x_{i}y_{j}:
i\in \{0,1,\dots ,m \}, j\in \{0,1,\dots ,n \}]$ is a projective variety in
$\mathbb{P}^{(n+1)(m+1)-1}$.
\end{Exercise}

A Schubert variety is the solution set to a collection of polynomial
equations in the ring generated by the flag minors. This is because the Schubert conditions $\dim(E_i \cap F_j) \geq \rk(w)[i,j]$ are equivalent to bounding the \emph{southwest} rank table of a matrix representing $F_\bullet$ by the southwest rank table of $M_w$, which in turn is equivalent to requiring certain minors to vanish.

\begin{Example}\label{ex:equations.2341}
For $n=4$, consider polynomials in the entries of a matrix using
variables $\{z_{ij}\given 1\leq i,j\leq 4 \}$.  For example,
$\Delta_{\{3,4 \},\{1,2 \}}=z_{31}z_{42}- z_{32}z_{41}$.  If $M$ is a
matrix representing a flag in $X_{2341}(E_{\bullet })$ then its
southwest rank conditions are bounded above by the table

\begin{equation}\label{eq:swranks.2341}
\left[\begin{array}{cccc}
1 &	2 &	3 &	4 \\
1 &    2 &	3 &	3\\
0 &    1 &      2 &	2\\
0 &    0 &      1 &	1
\end{array} \right]. 
\end{equation}
Therefore, the flag minors 
\begin{equation*}
\begin{array}{ll}
  \Delta_{\{3\},\{1\}} = z_{31}, & \Delta_{\{1,4\},\{1,2\}} = z_{11}z_{42}- z_{12}z_{41},\\
  \Delta_{\{4\},\{1\}} = z_{41},  & \Delta_{\{2,4\},\{1,2\}} = z_{21}z_{42}- z_{22}z_{41},\\
                                & \Delta_{\{3,4\},\{1,2\}} = z_{31}z_{42}- z_{32}z_{41},
\end{array}
\end{equation*}
evaluate to 0 on $M$. Conversely, the vanishing of these minors
implies the southwest rank conditions bounded above by
\eqref{eq:swranks.2341} all hold, so
$X_{2341}(E_\bullet)$ is the set of flags that have matrix
representations which are solutions to the five equations
$\Delta_{\{3\},\{1\}} = \Delta_{\{4\},\{1\}} =
\Delta_{\{1,4\},\{1,2\}} =  \Delta_{\{2,4\},\{1,2\}} = 
\Delta_{\{3,4\},\{1,2\}} =  0$.

It can be fruitful to consider more general minors vanishing on $M$. The set of minors vanishing on a general $M$ as above is 
\begin{equation}\label{eq:2341.minors}
\{z_{31},\ z_{41},\ z_{42},\ z_{11}z_{42}- z_{12}z_{41},\ z_{21}z_{42}- z_{22}z_{41},\ z_{31}z_{42}- z_{32}z_{41}   \}.
\end{equation}
However, $\Delta_{\{1,4\},\{1,2\}}, \ \Delta_{\{2,4\},\{1,2\}} $ and
$\Delta_{\{3,4\},\{1,2\}}$ are in the ideal generated by
$\{z_{31},\ z_{41},\ z_{42}\}$, so it suffices to check the three
entries in positions $(3,1), (4,1), (4,2)$ are 0 in any matrix $M$
representing a flag in $X_{2341}$. One should be a bit careful here:
$z_{42}$ is not a flag minor, and does not correspond to any
coordinate of the Pl\"ucker embedding of $\flags$ itself, even in a
projective sense.  However, $X_{2341}(E_\bullet)$ as a set can be also
be determined as the flags represented by the subset of $GL_{4}$ where $\{z_{31},\ z_{41},\ z_{42}\}$ vanish.

Note that if $b \in B$, then the southwest rank tables of $M$ and
$Mb$ agree everywhere, so it does not matter which representative of
the coset $MB$ we choose to represent a flag.  Conversely, if every
minor in \eqref{eq:2341.minors} vanishes on a matrix $M \in GL_{4}$,
then its southwest rank table is bounded above by
\eqref{eq:swranks.2341}, or equivalently its northwest rank table is
bounded below by \eqref{eq:rank.table.2341}.  Hence, the flag defined
by $M$ is in $X_{2341}(E_{\bullet })$. We will return to the equations
defining a Schubert variety both in projective space and in the affine
space of all $n \times n$ matrices in \S\ref{sec:CAGofSchubertVarieties}.  
\end{Example}
\bigskip

\begin{Exercise}\label{ex:Grass.cells}
For a $k$-subset $I \subset [n]$, let $C_{I}$ be the set of $k$-planes
in $\Gr (k,n)$ such that $\Delta_{I}$ is the first nonvanishing
Pl\"ucker coordinate in lex order.  Show that $C_{I}$ is a cell in the
sense that it has a parametrization as $\mathbb{C}^{d}$ for some $d$,
and these cells partition $\Gr (k,n)$.    What canonical matrices
represent the subspaces in $C_{I}$? (The sets $C_I$ are again called \emph{Schubert cells}, and their closures \emph{Schubert varieties}; in \S\ref{sub:123.StepFlags} we will explain how they relate to Schubert cells and varieties in $\flags$.)
\end{Exercise}

\begin{Exercise}\label{ex:Grass.Schubert.varieties}
Let $X_{I}=\overline{C}_{I}$ be the Schubert variety in $\Gr(k,n)$, defined
by the closure of the cell $C_{I}$, for any subset $I$ of size $k$ of
$[n]$.  Describe the containment relation on Grassmannian Schubert
varieties using a variation on the Ehresmann criterion for Bruhat
order from \Cref{thm:Ehresmann}.
\end{Exercise}

\begin{Exercise}\label{ex:find.minors}
  Identify all minors $\Delta_{I,J}$ which evaluate to be zero on every
  matrix representing a flag in $X_{4132}(E_{\bullet})$. Also, identify all \emph{flag} minors $\Delta_{I,[j]}$ which are zero on $X_{4132}(E_{\bullet})$.
  \end{Exercise}

\begin{Exercise}\label{ex:rank2.isom.types}
Prove that every one dimensional Schubert variety in $\Fl(n)$ is
isomorphic to $\mathbb{P}^{1}$ and that every two dimensional Schubert
variety in $\Fl(n)$ is isomorphic to $\mathbb{P}^{2}$ or
$\mathbb{P}^{1}\times \mathbb{P}^{1}$. 
\end{Exercise}

\begin{Exercise}\label{ex:rank2.isom.types.partials}
Does the classification of isomorphism types for Schubert varieties of
dimensions 1 and 2 from \Cref{ex:rank2.isom.types} hold in all partial
flag varieties?
\end{Exercise}

\subsection{Permutation Arrays and Modern Schubert Problems}\label{sub:PermutationArrays}

Now that we have some of the vocabulary of Schubert cells and
varieties in our repertoire, let's return to the topic of modern
Schubert problems.  Schubert was interested in enumerative geometry
problems where an intersection of geometric objects has a finite
number of possibilities in the generic case.  \textit{Modern Schubert
calculus is the study of the intersection numbers that arise for the
intersections of Schubert varieties in the generic case.}

Each Schubert variety is defined with respect to a fixed reference
flag.  When we consider intersecting Schubert varieties, we do so by
first moving them into general position.  This minimizes the overlap
as one would expect in the generic case.
Moving a Schubert variety is as simple as changing the reference flag.
Minimizing overlap of two varieties means the codimensions of the
varieties add up to the codimension of the intersection.  If the
codimensions of the subvarieties being intersected add up to the
dimension of the ambient variety containing them and the subvarieties
are in general position, then the intersection is 0-dimensional.  A
0-dimensional variety is a finite number of points.  Thus a
\textit{modern Schubert problem} on the flag variety would ask:
\medskip
\begin{quote}
How many flags are in the intersection $X_{u}(E_{\bullet}) \cap
X_{v}(F_{\bullet}) \cap \dots \cap X_{w}(G_{\bullet})$ assuming the
reference flags are chosen generically and the intersection is
0-dimensional?
\end{quote}
\medskip
Such a Schubert problem would require us to specify all of the
required dimensions of each intersection of subspaces over all the
reference flags and flags in the intersection.  This data structure is
a higher dimensional array of intersection conditions, beyond just a
matrix.

In \cite{ELcombinatorial} and \cite{ELdecomposition}, Eriksson and
Linusson developed a $d$-dimensional analog of a permutation matrix
toward characterizing all possible tables of intersection dimensions
for $d$ flags in $\flags$.  One way to generalize permutation matrices
is to consider all $d$-dimensional arrays of $0$'s and $1$'s with a
single $1$ in each hyperplane with a single fixed coordinate.  They
claim that a better way is to consider a permutation matrix to be a
two-dimensional array of 0's and 1's such that the rank of any
northwest submatrix is equal to the number of occupied rows in that
submatrix or equivalently equal to the number of occupied columns in
that submatrix.  The locations of the 1's in a permutation matrix will
be the elements in the corresponding permutation array.  We will
summarize their work here and refer the reader to their well-written
papers for further details.

Let $P$ be any collection of $d$-tuples in $[n]^{d}:=\{1,2,\dots,
n\}^{d}$.  We think of these $d$-tuples as the locations of \textit{dots} in an
$[n]^{d}$-\textit{dot array}.  Define a partial order on $[n]^{d}$ by
$$x=(x_{1},\dots, x_{d}) \preceq y=(y_{1}, \dots, y_{d}),$$ read ``$x$
is \textit{dominated} by $y$'', if $x_{i} \leq y_{i}$ for all $1\leq i
\leq d$.  This poset is a lattice with meet and join operation defined
by
\begin{eqnarray*}
x \vee y &=& z \quad \quad \text{ if  } z_{i} = \mathrm{max} (x_{i},
y_{i}) \text{ for all $i$}, 
\\
x \wedge y &=& z \quad \quad \text{ if } z_{i} = \mathrm{min} (x_{i}, y_{i})
\text{ for all $i$.}
\end{eqnarray*}
These operations extend to any set of points $R$ by taking $\bigvee
R=z$ where $z_{i}$ is the  maximum value in coordinate $i$ over the
whole set, and similarly for $\bigwedge R$.

For each $y \in [n]^{d}$, let $P[y]=\{x \in P \given x \preceq y \}$
be the \textit{principal subarray} of $P$ containing all points of $P$
which are dominated by $y$.  Define
\[
\rkj P= \# \{1\leq k\leq n \given \text{there exists } x \in P \text{ with } x_{j}=k \}.
\]
$P$ is \textit{rankable} of \textit{rank} $r$ if $\rkj P =r$ for all
$1\leq j \leq d$.  $P$ is \textit{totally rankable} if every principal
subarray of $P$ is rankable.  

\begin{Example}\label{ex:not.rankable}
The array $Q=\{(1,2),(2,2)\} \subseteq [2]^2$ is not rankable since
$\rk_1(Q)=2$ and $\rk_{2}(Q)=1$.
\end{Example}

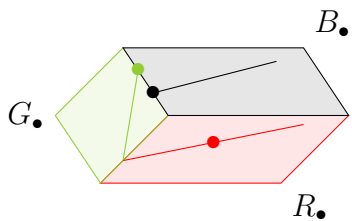
\begin{figure}[h]
\begin{tikzpicture}[scale=0.6]
\definecolor{MyGreen}{RGB}{145,205,50}
\draw[red, fill=red, fill opacity=0.1](0,0)--(-1.5,-1.5)--(2.5,-1.5)--(4,0);
\draw[MyGreen, fill=MyGreen, fill opacity=0.1](0,0)--(-1.5,-1.5)--(-2.5,0)--(-1,1.5);
\draw[black, fill=black, fill opacity=0.1](0,0)--(4,0)--(3,1.5)--(-1,1.5)--(0,0);
\draw[red] (-1,-1)--(3,-0.2);
\node[red] at (1,-0.6) {$\bullet$};
\draw[MyGreen] (-1,-1)--(-2/3,1);
\node[MyGreen] at (-2/3,1) {$\bullet$};
\draw (-1/3,0.5)--(2.4,1.2);
\node at (-1/3,0.5) {$\bullet$};
\node[left] at (-2.5,0) {$G_\bullet$};
\node[above right] at (3,1.5) {$B_\bullet$};
\node[below right] at (2.5,-1.5) {$R_\bullet$};
\end{tikzpicture}
\caption{Flags in position determined by the permutation array $P$
from \eqref{ex:permarray.P}.}
\label{fig:shoebox.perm.array}
\end{figure}

\begin{Example}\label{example.permarray.P}
For example with $n=4$, $d=3$, the following subset of $[4]^{3}$ is a totally
rankable dot array: 
\begin{equation}\label{ex:permarray.P}
P=\{(3,4,1), (4,2,2), (1,4,3), (3,3,3), (2,3,4),
(3,2,4),(4,1,4)\}.
\end{equation}
We picture this dot array as four 2-dimensional slices according to the last
coordinate, where the first one is ``slice $1$'' and the last is
``slice $4$''.  See \Cref{fig:slices.p}.  Here, $(3,4,1)$ corresponds
to the dot in the first slice from the left in \Cref{fig:slices.p}, and
$(4,2,2)$ corresponds to the dot in the second slice.  The two dots in the
third slice in correspond to $(1,4,3)$ and
$(3,3,3)$, etc.

\begin{figure}
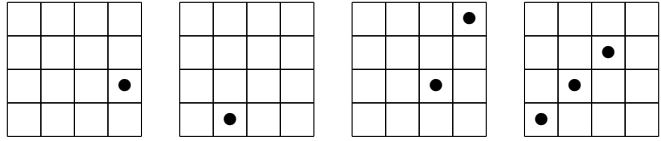

$$
\tableau{\emt & \emt & \emt & \emt \\
         \emt & \emt & \emt & \emt \\
         \emt & \emt & \emt & \ci \\
         \emt & \emt & \emt & \emt }
\hspace{.2in}
\tableau{\emt & \emt & \emt & \emt \\
         \emt & \emt & \emt & \emt \\
         \emt & \emt & \emt & \emt \\
         \emt & \ci & \emt & \emt }
\hspace{.2in}
\tableau{\emt & \emt & \emt & \ci \\
         \emt & \emt & \emt & \emt \\
         \emt & \emt & \ci & \emt \\
         \emt & \emt & \emt & \emt }
\hspace{.2in}
\tableau{\emt & \emt & \emt & \emt \\
         \emt & \emt & \ci      & \emt \\
         \emt & \ci      & \emt & \emt \\
         \ci      & \emt & \emt & \emt }
$$
\caption{Slices of a permutation array.}
\label{fig:slices.p}
\end{figure}

To verify this set $P \subset [4]^{3}$ is totally rankable, we need to
verify the rank is well defined for each principal subarray
$P[i,j,k]$ with $i,j,k\in [4]$.  For example, the points in $P$
dominated by $(3,4,4)$ are the principal subarray
\[
P[3,4,4]=\{(3,4,1), (1,4,3), (3,3,3), (2,3,4),(3,2,4)\} \subset P.
\]
These are the dots in $P$ that are visible to $(3,4,4)$ when looking
north, west, and down in \Cref{fig:slices.p} since the slices should
be considered as a $3$-dimensional transparent stack.

Note, that the first coordinate among these 3-tuples takes on values
$\{1,2,3 \}$, so ${\mathrm{rk}_{1}} P[3,4,4]=3$.  Similarly, the
second coordinate takes on values $\{2,3,4 \}$ and the third
coordinate takes on values $\{1,3,4 \}$, so ${\mathrm{rk}_{2}}
P[3,4,4]={\mathrm{rk}_{3}} P[3,4,4]=3$ also.  Since $\rkj P[3,4,4]=3$
for all $j \in [3]$, $P[3,4,4]$ is rankable and the rank of $P[3,4,4]$
equal to 3.  For the drawings below, we think of $\{1,2,3 \}$ as the
\textit{occupied rows} for the dots in $P[3,4,4]$, $\{2,3,4 \}$ are
the \textit{occupied columns}, and $\{1,3,4 \}$ are the
\textit{occupied slices} for $P[3,4,4]$. Continuing in the same way
checking that the number of occupied rows equals the number of
occupied columns equals the number of occupied slices for all
$P[i,j,k]$, we find $P$ is totally rankable. The full 3-dimensional rank table for
$P$ is represented by the list of tables in \Cref{fig:Pranktable}
showing one layer at a time.

\begin{figure}[h]
$$
\tableau{{ } & { } & { } & { } \\
         { } & { } & { } & { } \\
         { } & { } & { } & 1 \\
         { } & { } & { } & 1 }
\hspace{.2in}
\tableau{{ } & { } & { } & { } \\
         { } & { } & { } & { } \\
         { } & { } & { } & 1 \\
         { } & 1 & 1 & 2 }
\hspace{.2in}
\tableau{{ } & { } & { } & 1 \\
         { } & { } & { } & 1 \\
         { } & { } & 1 & 2 \\
         { } & 1 & 2 & 3 }
\hspace{.2in}
\tableau{{ } & { } & { } & 1 \\
         { } & { } & 1      & 2 \\
         { } & 1      & 2 & 3 \\
         1      & 2 & 3 & 4 }
$$
\caption{Rank table for $P=\{(3,4,1), (4,2,2), (1,4,3), (3,3,3),$ $(2,3,4),$
$(3,2,4),(4,1,4)\}$, where the empty boxes mean the rank is 0 for readability. }
\label{fig:Pranktable}
\end{figure}
\bigskip

The rank table in \Cref{fig:Pranktable} is realizable as the
intersection table for 3 flags in $\mathbb{C}^{4}$, call them
$B_{\bullet}, R_{\bullet}, G_{\bullet}$ for black, red, and green this
time.  Here we will think of the rows of the tables as being indexed
by $B_{1},B_{2},B_{3},B_{4}$ and the columns of each table as
$R_{1},R_{2},R_{3},R_{4}$, as before.  The $k^{th}$ slice from left to
right determines $\mathrm{dim}(B_{i}\cap R_{j}\cap G_{k})$ for $1\leq i,j\leq 4$.

A shoebox diagram of $B_{\bullet}, G_{\bullet}, R_{\bullet}$ with this
intersection data is shown in \Cref{fig:shoebox.perm.array}.  Let's
compare the intersection table with the picture.  We see
$\mathrm{dim}(B_{3}\cap R_{4}\cap G_{1})=1$ and we know $G_{1} \subset
G_{3}$ so the green point in the shoebox diagram corresponding to
$G_{1}$ is contained in the line at the intersection of the black
plane corresponding to $B_{3}$ and the green plane corresponding to
$G_{3}$.  The green dot is not coincident to either the black point or
the black line according to the data in the first two slices for
$G_{1}$ and $G_{2}$.  Since $R_{4}=\mathbb{C}^{4}$, no additional
constraints are imposed on $G_{1}$ by the red flag.  However,
$\mathrm{dim}(B_{4}\cap R_{2}\cap G_{2})=1$ so the green line
intersects the red line in the point where the red line meets the
green plane.  What is the dimension of the set of flags $B_{\bullet},
R_{\bullet}, G_{\bullet}$ that satisfy exactly the constraints imposed
by $P$?  What if $B_{\bullet}$ and $R_{\bullet}$ are fixed?
\end{Example}

\bigskip

\begin{Example}\label{example:totally.rankable}
A subset of a totally rankable array may or may not be rankable.
Consider the subset of $P$ above given by the array $A=\{(3,4,1),
(4,2,2), (2,3,4)\}$; then one can check $A$ is totally rankable.  On
the other hand, the array $B=\{(3,4,1), (4,2,2),(1,4,3) \} \subset P$,
is not totally rankable since the rank function is not well-defined at
every position in $[4]^{3}$.  For example, $B=B[4,4,4]$ is not
rankable. Indeed, the tuples in $B$ only have two distinct values
appearing in the second index $\{2,4 \}$ (occupied columns) and three
in the first and third, $\{1,3,4 \}$ and $\{1,2,3 \}$ (occupied rows
and slices) respectively. So ${\mathrm{rk}_{1}}
B[4,4,4]={\mathrm{rk}_{3}} B[4,4,4]=3,$ while ${\mathrm{rk}_{2}}
B[4,4,4]=2$.
\end{Example}

\bigskip

Many pairs $P, P'$ of totally rankable dot arrays are \emph{rank
equivalent}, i.e.\ $\rkj P[x] =\rkj P'[x]$, for all $x$ and $j$.
However, among all rank equivalent dot arrays there is a unique one
with a minimal number of dots \cite[Prop.\ 4.1]{ELcombinatorial}.  In
order to characterize the minimal totally rankable dot arrays, we give
the following two definitions.

\begin{Definition}\label{def:redundant.covered}
A position $y \in [n]^{d}$ is
\emph{redundant} in $P$ if there exists a collection of dots $R
\subset P\setminus \{y \}$ such that $y = \bigvee R$, $\# R >1$, and
every $x \in R$ has at least one $x_{i} =y_{i}$.  A position $y$ is
\emph{covered} by dots in $P$ if $y$ is redundant for some $R
\subset P \setminus \{y\}$, and for each $1 \leq j\leq d$ there exists
some $x \in R$ such that $x_{j}<y_{j}$.
\end{Definition}

For example, consider the totally rankable array $P$ in
\Cref{example.permarray.P} again.  The position $(3,4,3)$ is
redundant in $P$ since $(3,4,3)=(1,4,3) \vee (3,3,3)$ and
$(1,4,3),(3,3,3) \in P$.  However, the set $\{(1,4,3),(3,3,3) \}$ does
not cover $(3,4,3)$ since both of these dots have the same third
coordinate as $(1,4,3)$.  The set $R=\{(1,4,3), (3,3,3), (3,4,1)\}$
has join equal to $x=(3,4,3)$ and $R$ does contain an element with first
coordinate $1<x_{1}=3$, second coordinate $3<x_{2}=4$, and third coordinate $1<x_{3}=3$,
so the position $(3,4,3)$ is both redundant and covered in $P$.  On
the other hand, the position $(2,4,3)$ is not redundant in $P$, and
hence not covered.

\begin{Theorem}\label{t:EL}\cite[Theorem 3.2]{ELdecomposition} Let $P$ be a dot array.  The
following are equivalent:
\begin{enumerate}
\item $P$ is totally rankable.
\item Every two dimensional projection of every principal subarray is
totally rankable.
\item Every redundant position is covered by dots in $P$.
\item If there exist dots in $P$ in positions $y$ and $z$ and integers
$i,j$ such that $y_{i}<z_{i}$ and $y_{j}=z_{j}$, then there exists a
dot in some position  
$x \preceq (y \vee z)$ such that
$x_{i} = z_{i}$ and $x_{j} < z_{j}$.
\end{enumerate}
\end{Theorem}

Define a \textit{permutation array} in $[n]^{d}$ to be a totally
rankable dot array of rank $n$ with no redundant dots (or
equivalently, no covered dots).  The permutation arrays are the unique
representatives of each rank equivalence class of totally rankable dot
arrays with no redundant dots.  These arrays are Eriksson and
Linusson's analogs of permutation matrices.

For example, if $n=1$ and $d$ is any positive integer, there is a
unique permutation array $P=\{(1,1,..,1) \}=[1]^{d}$.  If $d=1$ and
$n$ is a positive integer, there is a unique permutation array
$P'=\{(i)\given i \in [n] \}=[n]^{1}$.  For $d=2$, the permutation
arrays are in bijection with permutations.  Just as permutations
determine all possible tables of intersection dimensions for two
flags, the definition of permutation arrays was motivated because they
include all possible relative configurations of flags such as the one
depicted in \Cref{fig:shoebox.perm.array}.

\begin{Theorem}  \label{t:EL-strata} 
\cite[Thm.\ 3.1]{ELdecomposition} Given $d$ complete flags
$E_{\ci}^{1}, E_{\ci}^{2}, \dots , E_{\ci}^{d}$ in $\mathbb{C}^{n}$,
there exists an $[n]^{d}$-permutation array $P$ with rank table equal
to  the table
of all intersection dimensions as follows.  For each $x \in [n]^{d}$,
\begin{equation}\label{e:EL-strata}
\rk(P[x]) = \mathrm{dim} \left(E_{x_{1}}^{1}\cap E_{x_{2}}^{2}\cap
\dots \cap E_{x_{d}}^{d} \right).
\end{equation}
\end{Theorem}

The feature of permutation arrays is that they are much
more manageable as data sets than the full table of intersection
dimensions.  The elements in a permutation array determine the minimal
jumps in dimension in the table of intersection dimensions of flags,
and therefore naturally correspond to critical vectors in the
geometry.

Based on many examples, Eriksson and Linusson
\cite[Conj.~3.2]{ELdecomposition} asked

\bigskip

\begin{quote}
Can every permutation array be realized by flags?
\end{quote}

\bigskip

\noindent We refer to this question as the \textit{Realizability
Conjecture}.  This question is motivated by more than curiosity.  A
fundamental question is: {\em what are the possible relative
configurations of $d$ flags?}  In other words: {\em what rank tables
(intersection dimension) are possible?}  For $d=2$, the answer leads
to the theory of Schubert cells and varieties discussed earlier.

The Realizability Conjecture is true for $d=1, 2, 3$.  For $d=1$, the
only permutation array is $[n]$.  The rank table for $[n]$ encodes the
dimensions of the subspaces in every flag in $\flags$.  For $d=2$, the
permutation arrays are in bijection with the permutation matrices,
hence they are all realizable for some pair of flags.  Realizability
of all permutation arrays in the case $d=3$ follows from \cite{SSV} (as
described in \cite[\S 3.2]{ELdecomposition}), see also \cite[\S
4.8]{Vakil-A}.  The case $n \leq 2$ is fairly clear, involving only
one-dimensional subspaces of a two-dimensional vector space (or
projectively, points on $\proj^1$), cf.\
\cite[Lemma~4.3]{ELdecomposition}.  Nonetheless, the conjecture is
false.  Counterexamples based on the Fano plane and Pappus' Theorem
were first given in \cite{Billey-Vakil}.  It is interesting that the
combinatorics of permutation arrays prevent some naive attempts at
counterexamples from working; somehow, permutation arrays see some
subtle linear algebraic information, but not all.

\begin{Problem}\label{prob:realizability}
Which permutation arrays in $[n]^{d}$ can be realized by $d$ complete
flags $E_{\ci}^{1}, E_{\ci}^{2}, \dots , E_{\ci}^{d}$ in
$\mathbb{C}^{n}$?
\end{Problem}

Eriksson and Linusson also gave an algorithm for producing all
permutation arrays in $[n]^{d}$ recursively from the permutation
arrays in $[n]^{d-1}$ \cite[Sect. 2.3]{ELdecomposition}.  We review
their algorithm, which we call the \textit{EL-algorithm} below, as
this is key to an algorithm for intersecting Schubert varieties.
Warning: this may get a bit technical!  It will not be used again in
later subsections, so it may be skipped by the reader without losing
insight in later material.

Let $A$ be any antichain of dots in $P$ under the dominance order.
Let $C(A)$ be the set of positions  covered by dots in 
$A$.  Define the \textit{downsizing} operator 
$D(A,P)$
with respect to $A$
on $P$  to be the result of the following process.
\begin{enumerate}
\item Set $Q_{1} = P \setminus A$. 
\item Set $Q_{2} = Q_{1} \cup C(A)$.
\item Set $D(A,P) = Q_{2} \setminus R(Q_{2})$ where $R(Q)$ is the set
of redundant positions of $Q$.
\end{enumerate}
The downsizing of a totally rankable array $P$ is \textit{successful}
if the resulting array is again totally rankable and has rank $\rk(P)-1$.

\begin{Theorem} \textbf{(The EL-Algorithm)}
\label{t:EL-algorithm} For positive integers $d,n>1$, the set of all
permutation arrays in $[n]^{d}$ can be obtained recursively by the
following depth-first search algorithm.

\begin{enumerate}
\item Choose a permutation array $P_{n}$ in $[n]^{d-1}$.  
\item If $P_n, \ldots , P_{i}$ have been defined and $i>2$, then
choose a nonempty antichain $A_{i}$ of dots in $P_{i}$ such that the
downsizing $D(A_{i},P_{i})$ is successful, set $P_{i-1}=D(A_{i},
P_{i})$.
\item Set $A_{1}=P_{1}$.  
\item Set $P = \bigcup_{i=1}^{n}\{(x_{1},\dots, x_{d-1},i) \given (x_{1},\dots,x_{d-1})
\in A_{i} \}$.  Add $P$ to the list of permutation arrays in $[n]^{d}$
constructed thus far.
\end{enumerate}
Furthermore, each permutation array $P$ is constructed from a unique
$P_{n}$ in $[n]^{d-1}$ and a unique sequence of nonempty antichains.
\end{Theorem}

\begin{Example} Let $d=2$ and $n=4$. 
Start with the $2$-dimensional permutation array 
\[
P_{4}=\{(1,4),(2,3),(3,1), (4,2) \}  \in [4]^{2}
\]
corresponding to the 1's in the permutation matrix for $w=3421$.  We
run through the algorithm as follows.  In the figure below, dots
correspond to elements in $P_{i}$ and circled dots correspond to
elements in $A_{i}$.
\[
\begin{array}{ll}
P_{4}=\{(1,4), (2,3),(3,1), (4,2) \} &  A_{4} = \{(1,4), (2,3) \}\\
P_{3} = \{(2,4),(3,1), (4,2) \} &       A_{3}= \{(3,1) \}\\
P_{2} = \{(2,4), (4,2) \} &     A_{2}= \{(2,4), (4,2) \}\\
P_{1} = \{(4,4)\} &     A_{1}= \{(4,4)\} 
\end{array}
\]

$$
\tableau{\emt & \emt & \emt & \emt \\
         \emt & \emt & \emt & \emt \\
         \emt & \emt & \emt & \emt \\
         \emt & \emt & \emt & \doubleci }
\hspace{.2in}
\tableau{\emt & \emt & \emt & \emt \\
         \emt & \emt & \emt & \doubleci \\
         \emt & \emt & \emt & \emt \\
         \emt & \doubleci & \emt & \emt }
\hspace{.2in}
\tableau{\emt & \emt & \emt & \emt \\
         \emt & \emt & \emt & \ci \\
         \doubleci & \emt & \emt & \emt \\
         \emt & \ci & \emt & \emt }
\hspace{.2in}
\tableau{\emt & \emt & \emt & \doubleci \\
         \emt & \emt & \doubleci & \emt \\
         \ci & \emt & \emt & \emt \\
         \emt & \ci & \emt & \emt }
$$
These choices lead to the 3-dimensional array 
\begin{equation}\label{eq:3-d.perm.array}
P= \{(4,4,1), (2,4,2), (4,2,2), (3,1,3), (1,4,4), (2,3,4) \}, 
\end{equation}
which the reader should check is a totally rankable array with no
redundant dots of rank 4.  

We prefer to display 3-dimensional dot-arrays as 2-dimensional
number-arrays as in \cite{ELdecomposition,Vakil-A} where a square
$(i,j)$ contains the number $k$ if $(i,j,k) \in P$.  The number-array 
$$
\tableau{\emt & \emt & \emt & 4 \\
         \emt & \emt & 4 & 2 \\
         3 & \emt & \emt & \emt \\
         \emt & 2 & \emt & 1 }
$$
also encodes $P$ from \eqref{eq:3-d.perm.array}, or equivalently the antichains
$A_{1},A_{2},A_{3},A_{4}$ above.  Note that there is at most one
number in any square of the number-array representing the permutation
array.  We claim this holds in general, not just this example. Why?
By Theorem~\ref{t:EL} Part~4, if two dots $y$, $z$ in a totally
rankable array $P$ existed such that $y_{1}=z_{1}, y_{2}=z_{2},
y_{3}<z_{3}$, then there exists a third dot $x \preceq (y \vee z) = z$
in $P$ with $x_{3}=z_{3}$ and $x_{1}<z_{1}$, but this
implies that $z$ is redundant for the set $R=\{x,y \}$, contradicting
the fact that $P$ is a permutation array.
\end{Example}

\begin{Corollary}\label{c:checkerboards}
Using the notation defined in Theorem~\ref{t:EL-algorithm}, each $P_{i}$ is a totally rankable
array  of rank $i$ in $[n]^{d-1}$.  Furthermore, if $P$
determines the rank table for flags $E^{1}_{\ci},\dots , E^{d}_{\ci}$,
then $P_{i}$ determines the rank table for $E^{1}_{\ci},\dots ,
E^{d-1}_{\ci}$ intersecting the subspace $E_{i}^{d}$, i.e.\
\[
\mathrm{rk}\left(P_{i}[x] \right) = \mathrm{dim}\left(E_{x_{1}}^{1}\cap
E_{x_{2}}^{2}\cap \cdots \cap E_{x_{d-1}}^{d-1} \cap E_{i}^{d} \right).
\]
\end{Corollary}

\begin{proof}
Observe $P_{i}$ is the totally rankable array obtained from
the projection
\[
\{(x_{1},\dots, x_{d}) \given (x_{1},\dots , x_{d}, x_{d+1}) \in P \text{
and } x_{d+1} \leq i \}
\]
by removing all covered elements.
\end{proof}

To represent a 4-dimensional permutation array $P$, we often draw the $n$
3-dimensional permutation arrays $P_{1}, \ldots, P_{n}$ from the
EL-algorithm.   For example, when $n=4$, the number-arrays 
$$\tableau{\emt & \emt & \emt & \emt \\
         \emt & \emt & \emt & \emt \\
         \emt & \emt & \emt & \emt \\
         \emt & 4 & \emt & \emt }
\hspace{.3in}
\tableau{\emt & \emt & \emt & \emt \\
         \emt & \emt & \emt & 4 \\
         \emt & \emt & \emt & \emt\\
         \emt & 4 & \emt & 3 }
\hspace{.3in}
\tableau{\emt & \emt & \emt & \emt \\
         \emt & \emt & \emt & 4\\
         \emt & \emt & 4 & 3 \\
%        4 & \emt & 3 & 2 } %%% Sara changed this 11.2.2023
         \emt & 4 & 3 & 2 }
\hspace{.3in}
\tableau{\emt & \emt & \emt & 1 \\
         4 & \emt & \emt & \emt \\
         \emt & \emt & 3 & \emt \\
         \emt & 2 & \emt & \emt } \vspace{.3in} 
$$
represent the 4-dimensional permutation array $P$ with entries
\[
\begin{matrix}
(4,2,4,1), (2,4,4,2), (4,4,3,2), (3,3,4,3),(3,4,3,3), (4,3,3,3), \\
  (4,4,2,3), (1,4,1,4), (2,1,4,4),(3,3,3,4), (4,2,2,4).
\end{matrix} 
\]
Here, $(4,2,4,1)$ encodes the fact that in position $(4,2)$
there is a 4 in the slice $P_{1}$ above.  Note, there are fewer
elements in $P$ than there are numbers on the four slices above.  This
is because the permutation array $P$ has no redundant dots by
definition, while the union of $P_{1},P_{2},P_{3},P_{4}$ does have
redundant dots in this case.
\bigskip

As an application, consider again a modern Schubert problem in
$\flags$.  Let $E^{1}_{\ci},\dots , E^{d}_{\ci}$ be fixed flags in
$\mathbb{C}^{n}$, and let $w^{1},w^{2},\ldots, w^{d}$ be permutations
in $S_{n}$.  The permutations are not required to be distinct.  Let

\begin{equation}\label{eq:X.SchubertProblem}
X=X_{w^{1}}(E^{1}_{\ci})\cap X_{w^{2}}(E^{2}_{\ci}) \cap \dotsb \cap
X_{w^{d}}(E^{d}_{\ci}).
\end{equation}
If $\sum_i \coinv(w^{i}) = \binom{n}{2}$, then $X$ is 0-dimensional.
The associated \textit{Schubert problem} is to find the maximum number
of flags in $X$ over all choices of the $d$ fixed flags such that $X$
is a finite set.  This maximum number is the \textit{intersection number} for
$w^{1},w^{2},\ldots, w^{d}$.  The intersection number is the answer to
the statement of a modern Schubert problem in beginning of
\Cref{sub:SchubertProblems2000}.

By \Cref{def:Schubert.variety}, a flag $F_{\bullet}$ is in
$X=X_{w^{1}}(E^{1}_{\ci})\cap X_{w^{2}}(E^{2}_{\ci}) \cap \dotsb \cap
X_{w^{d}}(E^{d}_{\ci})$ if and only if
\[
\mathrm{dim} \left(E_{i}^{a}\cap F_{j}\right) \geq
\mathrm{rk}(w^{a})[i,j]
\]
 for all $ 1\leq i,j \leq n$ and all $1\leq a \leq d$. Such a flag
also gives rise to a $(d+1)$-dimensional table of intersection dimensions
\begin{equation}\label{eq:dims}
\mathrm{dim} \left(E_{x_{1}}^{1}\cap E_{x_{2}}^{2}\cap
\cdots \cap E_{x_{d}}^{d} \cap F_{x_{d+1}}\right).
\end{equation}
By \Cref{t:EL-strata}, this table of intersection dimensions is
encoded by some permutation array $P \subset [n]^{d+1}$.  In this
case, we say $F_{\bullet}$ has \textit{relative position} $P$ with
respect to $E^{1}_{\ci},\dots , E^{d}_{\ci}$.  Since  $E_{n}^{a} =
\mathbb{C}^{n}$ for all $a\in [d]$, one can recover the $2$-dimensional tables
$\mathrm{dim} \left(E_{i}^{a}\cap F_{j}\right)$ from the rank table
for $P$, or equivalently from the table \eqref{eq:dims}, by taking all
$x_{i}=n$ except for $i=a$.  Therefore, if any other flag
$G_{\bullet}$ and $E^{1}_{\ci},\dots , E^{d}_{\ci}$ also have relative
position $P$, then $G_{\bullet} \in X$ as well.

For a permutation array $P\subset [n]^{d} $, let
$C_{P}(E^{1}_{\ci},\dots , E^{d}_{\ci}) \subset \flags$ be the set of
all flags $F_{\bullet }$ with relative position $P$ with respect to
$E^{1}_{\ci},\dots , E^{d}_{\ci}$.  Again, if $E^{1}_{\ci},\dots ,
E^{d}_{\ci}$ are fixed, we can suppress the list and just write
$C_{P}$.   The $C_{P}$'s are generalizations of
Schubert cells.  It follows that $X$ decomposes as a union of
$C_{P}$'s.  This holds even when $X$ is not a 0-dimensional
intersection or if the reference flags are not in generic position.

Billey-Vakil \cite[Sect. 5]{Billey-Vakil} showed that in the case
$E^{1}_{\ci},\dots , E^{d}_{\ci}$ are generically chosen, there is a
unique permutation array $P$ such that $X=C_{P}$ provided $X$ is
nonempty, and furthermore $P$ can be identified by an explicit
recursive algorithm.  They also showed how to use $P$ to write down
equations for $X$.  These equations can also be used to determine if
$E^{1}_{\ci},\dots , E^{d}_{\ci}$ are sufficiently general for
computing the generic number of flags in the intersection assuming the
$d$ permutations in $S_{n}$ are fixed and the reference flags are
allowed to vary.  The number of solutions will always be either
infinite or no greater than the expected number.  The expected number
is achieved on a dense open subset of $\flags^d$.  Following
\cite{Billey-Vakil}, there are two key statements needed to identify
the unique permutation array $P$ for this problem when $X$ is
nonempty.  \Cref{lem:generic} below could be a good exercise for the
reader.

\begin{Lemma}\label{lem:generic}
The permutation array corresponding to $d$ generically chosen flags
$E^{1}_{\ci},\dotsc , E^{d}_{\ci}$ is given by the {\bf transverse
permutation array}
\begin{equation}\label{e:transverse}
T_{n,d} = \left\{(x_{1},\dots , x_{d}) \in [n]^{d}\given \sum x_{i}
=(d-1)n+1 \right\}.
\end{equation}
This permutation array corresponds to the transverse rank table
$$
\rk(T_{n,d} [x]) =
\max \left( 0,  n - \sum_{i=1}^d (n-x_i) \right)
$$
for all $x=(x_{1},\dots , x_{d}) \in [n]^{d}$.  
\end{Lemma}

\begin{Remark}\label{rem:transversality}
The nomenclature comes from the geometry of linear spaces.  Two subspaces $U$ and $V$
in $\mathbb{C}^{n}$ are said to be in \textit{transverse} position
provided their intersection is just the origin unless
$\mathrm{dim}(U)+\mathrm{dim}(V)\geq n$, and if
$\mathrm{dim}(U)+\mathrm{dim}(V)\geq n$, then the dimension of $U
\cap V$ is $ \mathrm{dim}(U)+\mathrm{dim}(V)-n$, which is between 0
and $n$.  The terminology carries over to manifolds and varieties by
considering their tangent spaces at any point.  The concept
generalizes the notion of two lines in ``general position'' in
$\mathbb{R}^{3}$.  Two generically chosen flags $F_{\bullet}$ and
$G_{\bullet}$ will have subspaces $F_{i}$ and $G_{j}$ in transverse
position.  It is very difficult to specify a generic choice of flags
in practice since we are restricted by the limits of computation with
modern computers.  However, testing for transversality is easy in
comparison.  See \cite[Cor 5.2]{Billey-Vakil} for a sufficient
criterion for genericity.  
\end{Remark}

\begin{Example}\label{ex:transverse.perm.array}
In the case $n=4, d=3$, the transverse permutation array is
\[T_{4,3}=
\tableau{\emt & \emt & \emt & 4 \\
         \emt & \emt & 4 & 3 \\
	 \emt & 4 & 3 & 2 \\
	 4 & 3 & 2 & 1 }.
\]
\end{Example}
\bigskip

To describe the unique permutation array $P$ corresponding to a
Schubert problem with permutations $w^{1},\ldots , w^{d}$, note that
each intersection of the form
\[
E_{x_{1}}^{1}\cap E_{x_{2}}^{2}\cap
\cdots \cap E_{x_{d}}^{d} \cap F_{x_{d+1}}
\]
can be simplified if $x_{i}=n$ since $E_{n}^{i}=\mathbb{C}^{n}$.
Furthermore, by construction $\mathrm{dim} \left(E_{i}^{a}\cap
F_{j}\right)$ is determined by the rank table corresponding to the
permutation $w^{a}$.  Thus, to determine $\rk(P[x])$ for $x=(x_{1},\dots ,
x_{d},j)$, it suffices to consider the dimensions of intersections of
the form
$$
d_{j}(s_{1},s_{2},\dots , s_{k}) := \mathrm{dim}
\left(E_{x_{s_{1}}}^{s_{1}}\cap E_{x_{s_{2}}}^{s_{2}}\cap \dotsc \cap
E_{x_{s_{k}}}^{s_{k}} \cap F_{j}\right)=\rk(P[x])
$$
where $1\leq s_{1}<s_{2}<\dotsb <s_{k}\leq d$,\, $k\geq 2$,\ and $1
\leq x_{s_{j}} \leq n-1$ for each $1\leq j\leq k$, and $x_{i}=n$ for
all $i \not \in \{s_{1},s_{2},\ldots ,s_{k} \}$.  \Cref{t:unique} can be derived from the algorithm given in
\cite[Sect. 5]{Billey-Vakil}, but did not appear explicitly there.

\begin{Proposition}\label{t:unique}\cite[Sect. 5]{Billey-Vakil}
The unique permutation array $P$ corresponding to a Schubert problem
$X$ as in \eqref{eq:X.SchubertProblem} with permutations $w^{1},\ldots
, w^{d}$ satisfying $\sum \mathrm{codim}(w^{a}) = \binom{n }{2}$ is
determined by the table of intersection dimensions satisfying the
recurrence
$$
d_{j}(s_{1},s_{2},\dots , s_{k}) = \mathrm{max}
\begin{cases}
d_{j}(s_{1}) + d_{j}(s_{2},\dots , s_{k})-j
\\
d_{j-1}(s_{1},s_{2},\dots , s_{k})
\end{cases}
$$
where $1\leq s_{1}<s_{2}<\dotsb <s_{k}\leq d$,\, $k\geq 2$, \, $1 \leq
x_{s_{i}} \leq n-1$ for each $1\leq i\leq k$, $
\mathrm{dim}\left(F_{0}\right) :=0$, and for all $F_{\bullet}\in X$
the value $d_{j}(s_{i})=\mathrm{dim} \left(E_{x_{s_{i}}}^{s_{i}}\cap F_{j}\right)$ is
determined by the rank table corresponding to the permutation $w^{s_{i}}$.
\end{Proposition}

\begin{proof}
One can determine $\rk(P[x])$ at any fixed $x=(x_{1},\dots , x_{d},j)$ with
$1 \leq x_{s_{j}} \leq n-1$ for each $1\leq j\leq k$, and $x_{i}=n$
for all $i \not \in \{s_{1}<s_{2}<\ldots <s_{k} \}$ by induction on
$j,k\geq 0$. The base case when $j=0$ is given by
$d_{0}(s_{1},s_{2},\dots , s_{k}) =0$ for all $1\leq
s_{1}<s_{2}<\dotsb <s_{k}\leq d$ since $\mathrm{dim}(F_{0})=0$.  If
$k=0$, we have $\rk(P[x])=j$.  If $k=1$, then $\rk(P[x])$ is determined by
the hypothesis $d_{j}(s_{1})=\mathrm{dim}
\left(E_{x_{s_{1}}}^{s_{1}}\cap F_{j}\right)=\rk(P[x])$.

For $j>0$ and $k>1$, assume $d_{j}(s_{1}), d_{j}(s_{2},\dots ,
s_{k}),$ and $d_{j-1}(s_{1},s_{2},\dots , s_{k})$ are known by
induction.  Since $F_{j-1}\subset F_{j}$, we know
\[
d_{j}(s_{1},s_{2},\dots , s_{k}) \geq d_{j-1}(s_{1},s_{2},\dots , s_{k}).
\]
For any two subspaces $U,V$ in $F_{j}$, we know $\dim(U \cap V)\geq
\dim(U)+\dim(V)-j$. Therefore,
\[
d_{j}(s_{1},s_{2},\dots , s_{k}) \geq d_{j}(s_{1}) + d_{j}(s_{2},\dots , s_{k})-j.
\]
Since the base flags $E_{\bullet}^{s_{1}},E_{\bullet}^{s_{2}},\dots,
E_{\bullet}^{s_{k}} $ are generic, $\rk(P[x])=d_{j}(s_{1},s_{2},\dots , s_{k})$
is the minimal value satisfying these two inequalities.  
\end{proof}

\begin{Remark}\label{rem:twisted}
Note in \cite{Billey-Vakil}, the Schubert variety labeled by $w$ is our $X_{ww_0}$,
 but the difference is only in the label.  The
rank conditions are written the same way.  So, in their work, $\codim
(X_{w})=\ell(w)$, and a $0$-dimensional intersection has $\sum \ell(w^{(i)})=\binom{n }{2}$.  
\end{Remark}

\begin{figure}
\begin{center}
\includegraphics[scale=0.8]{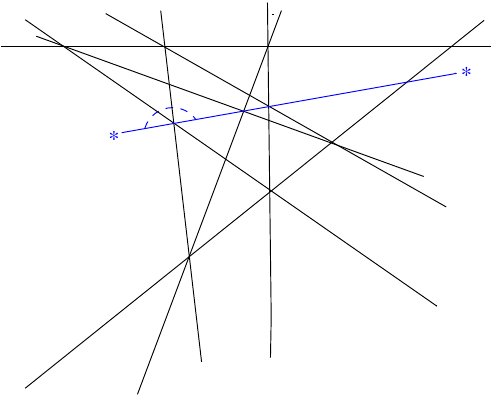}
\end{center}
\caption{The Pappus line configuration on 9 lines.  The arc represents
an impossible ``hop'' by the *blue line* over an intersection.  This
leads to an unrealizable permutation array.}  \label{fig:pappus}
\end{figure}

In \cite{Vakil-A}, Vakil gave an algorithm to determine intersection
numbers for Schubert varieties in Grassmannian varieties using a
``checkers game''.  Each step of this game corresponds geometrically to a torus
degeneration that takes the moving flag one step closer to the base
flag, in a specified way.  In each step, the corresponding
intersection either corresponds with one checker configuration or it
branches into two components.  The black checkers in the game play the
role of the permutation array determining the table of intersections
between the two complete flags.  The white checkers determine another
permutation array that determines the table of intersections of a
subspace $V$ in the intersection of two Grassmannian Schubert
varieties with the subspaces formed by intersecting the subspaces in
the two flags.  In the end, the moving flag and the base flag become
the same, and at that point, each white checker configuration
determines exactly one irreducible component of the intersection. The
intersection numbers are given by counting the different types of
terminal checkers games.

% AlgorithmsInAGVolume.pdf
% Flags.fpsac.final.pdf 
\begin{Exercise}\label{ex:rank.to.perm.array}
For a permutation array $P$ with rank table $\mathrm{rk}(P[x])$ for all
$x\in [n]^{d}$, describe an algorithm to reconstruct $P$ from the
data in the rank table.
\end{Exercise}

\begin{Exercise}\label{ex:fano.plane}
The Pappus configuration of lines is shown in \Cref{fig:pappus}.
First, find a permutation array that goes with the Pappus
configuration by choosing a special point on each line. Second, find
the permutation array that goes with the line configuration shown in
\Cref{fig:pappus} where the blue line somehow does not intersect the
leftmost crossing.  Note, this picture is drawn projectively so each
line represents a 2-dimensional plane, which all contain the origin.
\end{Exercise}

\begin{Problem}\label{Prob:degeneration}
Is there a method to compute intersection numbers for $\flags$ via
moving the dots in a permutation array in some specified way?
\end{Problem}

\subsection{Schubert Calculus and Cohomology of Flag Varieties }\label{sub:Monk}

Intersection of sets is a commutative relation, just like
multiplication in a commutative ring.  What if we could convert
Schubert problems in enumerative geometry into problems in polynomial
rings?  This question leads to the topics of intersection theory,
Schubert calculus, cohomology rings, and Chow rings.  We will survey
this approach here and refer the interested reader to the two books by
William Fulton on the topic \cite{Fulton,Fulton6}, which laid out the
rigorous foundation of the theory.  See also the book by
Eisenbud-Harris called ``3264 \& All That: Intersection Theory in
Algebraic Geometry'' \cite{eisenbud_harris_2016} and the recent book
by Anderson-Fulton ``Equivariant Cohomology''
\cite{Anderson-Fulton}.

The integral cohomology ring of a topological space $X$ is a
graded-commutative ring, denoted $H^{*}(X)= \bigoplus H^{i}(X)$.  We
will consider cohomology only over the integers, but rational numbers
and other choices of coefficient rings are important as well. If two
spaces are homotopy equivalent, then their cohomology rings are
isomorphic.  There are different ways to define cohomology such as
singular cohomology or de Rham cohomology. The Chow ring is another
variation on this theme.  It is notable that for the flag varieties $\flags$, 
these different cohomology theories agree.  In each case, one could
ask for an algebraic model of the cohomology ring of $X$, for instance
as a polynomial ring modulo an ideal.  For example, complex projective
space $\mathbb{P}^{n}$, consisting of complex lines through the origin
in $\mathbb{C}^{n+1},$ has integral cohomology ring isomorphic to
$\mathbb{Z}[x]/\langle x^{n+1} \rangle$, where $x$ is an indeterminate
of degree $2$.  Here $x$ has degree 2 because $\mathbb{C}\simeq
\mathbb{R}^{2}.$

In the 1950's, Armand Borel was a visionary in algebraic topology and
algebraic groups, not to be confused with Emile Borel of Borel measure
fame.  One of his celebrated results relates the cohomology ring of
the flag varieties to the coinvariant algebra from invariant theory.
His work generalizes to all flag varieties $G/B$ for semisimple Lie
groups $G$ and Borel subgroups $B$ \cite{Borel.1953}.  A key fact used
to prove his result is that $G/B$ can also be represented by a compact
Lie group modulo a maximal torus, $U_{n}/T$; see \Cref{ex:unitary}.
Borel was a member of the Bourbaki group, along with Bruhat, Cartan, 
Chevelley and Ehresmann \cite{borel.obit}.  The Bourbaki text
``Groupes et Alg\`ebres de Lie'' is a classic text related to
generalized flag varieties, Weyl groups, root systems, and invariant
theory \cite{bou1}.  These days much of the material in \cite{bou1} is
also covered by \cite{b-b} and \cite{Hum}, but some of the motivation
from flag varieties has been omitted.

\begin{Definition}\label{def:coinvariant.algebra}
The \emph{coinvariant algebra} $R_{n}$ is the quotient of the
polynomial ring $\mathbb{Z}[x_{1},x_{2}, \dots , x_{n}]$ modulo the
ideal $I_{n}^{+}$ generated by all homogeneous polynomials of positive
degree that are invariant under the natural $S_{n}$ action permuting
the variables.
\end{Definition}

The Fundamental Theorem of Invariant Theory says that the ring of
symmetric polynomials in $\mathbb{Z}[x_{1},x_{2}, \dots , x_{n}]$ is
exactly the freely generated polynomial ring $\mathbb{Z}[h_{1},h_{2},
\dots , h_{n}],$
where the \textit{homogeneous symmetric polynomial} of
degree $k$ is 
\begin{equation}\label{eq:hom.sym.polys}
h_{k}=h_{k}(x_{1},\dots , x_{n})= \sum x_{1}^{i_{1}}x_{2}^{i_{2}}\cdots x_{n}^{i_{n}}
\end{equation}
and the sum is over all monomials of degree $k=i_{1}+...+i_n$.
Therefore, the ideal $I_{n}^{+}$ is minimally generated by
$h_{1},h_{2},\ldots , h_{n}$.

\begin{Theorem}\cite{Borel.1953}\textbf{The Borel Presentation.}\label{thm;borel}
The cohomology ring of the flag variety is isomorphic to the
coinvariant algebra,
\[
H^{*}(\flags, \mathbb{Z}) \simeq \mathbb{Z}[x_{1},x_{2}, \dots ,
x_{n}]/I_{n}^{+} = R_{n}
\]
\end{Theorem}
\medskip

\begin{Exercise}\label{ex:artin.monomials}
Use the fact that $I_{n}^{+}=\langle h_{1},h_{2},\ldots , h_{n}
\rangle$ to prove $R_{n}$ has a linear basis of monomials of the form
$x_{1}^{c_{1}}x_{2}^{c_{2}}\cdots x_{n}^{c_{n}}$ where $0\leq
c_{i}\leq n-i$.  Hence, $\mathrm{dim}(R_{n})=n!$.
\end{Exercise}

To connect the Borel presentation with modern Schubert problems, we
give a brief introduction to the Chow ring of a variety, instead of
using singular cohomology.  Since the two rings are isomorphic for the
flag variety, one can use either variation.  Given a smooth variety
$Z$, each of its irreducible subvarieties $X \subset Z$ gives rise to
an element in the \textit{Chow ring} of $Z$, called its \textit{Chow
class} and denoted $[X]$.  We consider two irreducible subvarieties
$X,Y \subset Z$ to give rise to the same element in the Chow ring if
they are rationally equivalent, denoted $X \sim Y$.  Roughly speaking,
$X \sim Y$ means we can smoothly deform $X$ from its current position
in $Z$ to align with $Y$.
%For example,
%%any two lines in $\mathbb{R}^{2}$ are rationally equivalent, but also
%line is rationally equivalent to a parabola.
Let $A^{*}(Z) = \bigoplus A^{d}(Z)$ denote the Chow ring of $Z$ with
elements given by linear combinations of cohomology classes modulo
rational equivalence.  The graded component $ A^{d}(Z)$ is spanned by
classes $[X]$ where $\codim (X)=d$.

Every ring has a rule for addition and a rule for multiplication.
Addition in the Chow ring is just formal addition. If $X,Y$ are
subvarieties of $Z$, the sum of the corresponding classes is denoted
$[X]+[Y]$.  The additive identity is the class of the empty set, which
is a subvariety of $Z$ itself.  We write $0=[\emptyset]$.
Furthermore, if $Y \subset Z $ is the union of irreducible
subvarieties $X_{1} \cup \cdots \cup X_{k}$, then the class of $Y$
decomposes as $[Y]=[X_{1}]+\cdots+ [X_{k}]$.  Some of these components
$X_i$ may be rationally equivalent, so coefficients other
than 0 or 1 can appear in the expansion of $[Y]$.

Multiplication in the Chow ring is modeled on the intersection
problems discussed in this chapter and later in this book.  To
multiply the classes $[X]\cdot [Y]$, first we ``move'' $Y$ into
transverse position with respect to $X$ via rational equivalence to
obtain $Y'$, and then take the class of the intersection.  So by
definition $[X]\cdot [Y] =[X \cap Y']$. The class $[X \cap Y']$ is
well defined and any generic choice of $Y' \sim Y$ suffices to
identify the class of the intersection.  Furthermore, if $[X] \in
A^{d}(Z)$ and $[Y] \in A^{e}(Z)$, then $[X][Y] \in A^{d+e}$ reflects
the fact that $\codim (X)+\codim (Y)=\codim (X\cap Y)$ provided the
intersection is nonempty and $X,Y$ are in transverse position.  Thus,
the Chow ring $A^{*}(Z)$ is a graded, commutative, associative ring
with multiplicative identity given by $[Z]$.  We should note that the
process of ``moving'' $Y$ to an appropriate $Y'$ is a subtle
construction for arbitrary varieties $Z$, but in the case of flag
varieties it is straightforward, as discussed below.

Defining rational equivalence explicitly is beyond the scope of this
chapter because it is a subtle concept in general, and we can simply
state the facts we need in the context of Schubert varieties.
Instead, we offer some examples to build some intuition for the
concept. Any two hypersurfaces defined by polynomial equations $f
= 0$ and $g = 0$ in $\mathbb{R}^n$ are rationally equivalent, because
we can let $t$ vary from $0$ to $1$ in the family $tf + (1-t)g$ =
0. So every line in $\mathbb{R}^2$ is rationally equivalent to every other line, and
they are also equivalent to a parabola. Rational
equivalence is more interesting in projective space.  For instance,
suppose $X$ is the subvariety $\{[x:y] \in \mathbb{P}^1 : y = 0\} = \{[1:0]\}$
and $Y$ is the subvariety defined by $x^2 - y^2=0$ consisting of the
two points $[1:1], [1:-1]$. We cannot consider the family of
hypersurfaces defined by $t(x^2-y^2) + (1-t)y = 0$ since this equation
is not homogeneous in $x,y$.  We could homogenize the equation to get
$t(x^2-y^2) + (1-t)y^2 = 0$, but letting $t$ vary from 0 to 1 now
shows that $Y$ is rationally equivalent to the subvariety defined by $y^2 = 0$. In the algebraic
setting of the Chow ring, we must view this as defining the point
$[1:0]$ with multiplicity 2, which is not the same as $X$. With a
little more work it can be shown that $X$ and $Y$ are not rationally
equivalent because the multiplicity does matter.

As described above, the Chow ring $A^{*}(Z)$ may be infinitely
generated as a ring.  However, it has a very nice structure as an
abelian group in the following special but important case. We say $Z$
has a \emph{cellular decomposition} or \emph{affine paving} if there
is a chain of subvarieties $\emptyset = Z_0 \subseteq Z_1 \subseteq
\cdots \subseteq Z_m = Z$ such that each $Z_{i} \setminus Z_{i-1}$ is
a disjoint union  $\bigcup_j C_{ij}$ where each $C_{ij}$ is
isomorphic, as a complex variety, to some affine space $\C^d$. The
sets $C_{ij}$ are called the \emph{cells} of $Z$. Cellular
decompositions are algebraic analogs of the notion of a
CW-decomposition in topology.

If the cells of $Z$ are $\{C_{1}, \ldots, C_{k} \}$, then their
Zariski closures $\overline{C}_{i}$ are subvarieties in $Z$ by
definition. The corresponding classes $\{[\overline{C}_{i}] \given i
\in [k] \}$ form a linear basis for the Chow ring of $Z$, which is free abelian of rank $k$.  More specifically,
 the $d^{th}$ graded component $A^{d}(Z)$ is spanned by
all of the classes $[\overline{C}_{i}]$ such that $\codim (C_{i})=d$.
One can further show that $A^*(Z)$ is isomorphic to the singular
cohomology ring $H^*(Z; \bZ)$ when $Z$ is smooth and has a cellular
decomposition; see \cite[Ch. 1 and Ch. 19]{Fulton}.

\begin{Exercise}\label{ex:projective.space}
Describe the Chow ring of $\mathbb{P}^{n}$ by identifying a cellular
decomposition of the space so that the additive and multiplicative
structure is isomorphic to $\mathbb{Z}[x]/\langle x^{n+1}
\rangle.$
\end{Exercise}

\begin{Exercise}\label{ex:affine.space}
Describe the Chow ring $A^{*}(\mathbb{C}^{n})$.
\end{Exercise}

Related varieties have related Chow rings. To be precise, if $f : Y \to Z$ is a map between smooth varieties, there is a ring homomorphism $f^* : A^*(Z) \to A^*(Y)$ induced in the opposite direction \cite[Ch. 6]{Fulton}. This assignment is functorial in the sense that $(f \circ g)^* = g^* \circ f^*$ and the identity map on a variety induces the identity map on its Chow ring. In general the action of $f^*$ can be non-trivial to describe, but we use only the following special cases.
\begin{Proposition} \label{prop:pullback} \hfill 
  \begin{enumerate}[(a)]
    \item Suppose $f$ is a closed embedding. Given $[X] \in A^*(Z)$,
    choose a subvariety $X' \subseteq Z$ rationally equivalent to $X$
    and transverse to $f(Y)$, so that $[X][f(Y)] = [X' \cap
    f(Y)]$. Then $f^*([X]) = [f^{-1}(X' \cap f(Y))]$.
    \item Suppose $Y$ and $Z$ are smooth and irreducible, and $\dim
    f^{-1}(z)$ is constant for $z \in Z$. Then $f^*([X]) =
    [f^{-1}(X)]$.
  \end{enumerate}
\end{Proposition}
See \cite[Ch. 6]{Fulton} and \cite[\S 1.7, \S 8.3]{Fulton}. For example, suppose $f : \bP^n \to \bP^{n+1}$ is the embedding sending $[z_0 : \cdots : z_n]$ to $[z_0 : \cdots : z_n : 0]$. Using the isomorphism from Exercise~\ref{ex:projective.space}, $f^*$ is the quotient map $\bZ[x]/\langle x^{n+1} \rangle \to \bZ[x]/\langle x^{n} \rangle$.

The complete flag variety $\flags$ is a smooth variety and it has a
nice cell decomposition, so its Chow ring $A^{*}(\flags)$ serves as a
fundamental example for this theory.  There is a smooth path through
$GL_{n}(\mathbb{C})$ connecting every invertible matrix with any other
by left multiplication, and this action carries over to flags in $G/B
\simeq \flags$.  Therefore, $X_{w}(E_{\bullet})$ is rationally
equivalent to $X_{w}(F_{\bullet})$ for all $F_{\bullet} \in \flags$.  The choice
of base flag is irrelevant in terms of computing cohomology classes in
the Chow ring.  So, we denote the \textit{Schubert classes}
$[X_{w}(E_{\bullet})]=[X_{w}(F_{\bullet})]$ simply by $[X_{w}]$.

To identify an affine paving of $\flags$, one can use the rank
function on Bruhat order and equation~\eqref{prelim 100}.  Set
\[
Z_i = \bigcup_{\substack{w \in S_{n}:\\
 \ell(w) = i}} X_{w} = \bigcup_{\substack{v \in S_{n}:\\
 \ell(v) \leq  i}} C_{v}
\]
for $0\leq i \leq \binom{n}{2}$. It is straightforward to verify that
$\emptyset = Z_0 \subseteq Z_1 \subseteq \cdots \subseteq Z_m = Z$ and
each $Z_{i} \setminus Z_{i-1}$ is a disjoint union of cells in this
case, namely the Schubert cells $C_{w}$ for $w \in S_{n}$ with $i$
inversions.  Therefore, the Schubert classes $\{[X_{w}]\given w\in
S_{n} \}$ form a basis of the Chow ring of $\flags$, and hence, every
subvariety of the flag variety is rationally equivalent to a union of
translates of Schubert varieties $X_{w}(F_{\bullet})$ for possibly
different base flags. Since the complex codimension of $C_{w}$ is the
number of non-inversions of $w$, denoted $\coinv (w)$, we know the
$2d^{th}$ graded component of the Chow ring of $\flags$ is spanned by
all of the Schubert classes $[X_{w}]$ such that $\coinv (w)=d$.
Furthermore, $[X_{u}]=[X_{v}]$ in $A^{*}(\flags)$ if and only if $u=v$,
by linear independence.   

In summary, the Chow ring of $\flags$, denoted
\begin{equation}\label{eq:Chow}
A^{*}(\flags)= \bigoplus_{d=0}^{\binom{n}{2}} A^{2d}
\end{equation}
is a commutative graded ring such that $A^{2d}$ is spanned by
$\{[X_{w}]\given \coinv (w) =d \}$ and $A^{2d+1}$ is 0.  Therefore,
its \textit{Poincar\'e-Hilbert} series is
\begin{equation}\label{eq:poincare}
\sum_{d\geq 0} \mathrm{dim}(A^{d})q^{d} = [n]_{q^{2}}!
\coloneqq \sum_{w \in S_{n}} q^{2\coinv(w)} =
\prod_{k=1}^{n-1}(1+q^{2}+q^{4}+ \cdots + q^{2k}).
\end{equation}
Since the Chow ring and the cohomology ring of $\flags$ are isomorphic,
we know \eqref{eq:poincare} is the Poincar\'e polynomial of $\flags$
and the Euler characteristic of $\flags$ is $n!$ since 
\[
\chi (\flags)=
\sum_{d\geq 0}(-1)^{d} \mathrm{dim}(A^{d})=
\left[\prod_{k=1}^{n-1}(1+q^{2}+q^{4}+ \cdots +
q^{2k})\right]_{q=-1} = n!.  
\]

\begin{Remark}
Warning, many authors dispense with the technicality that $\mathbb{C}$
is 2-dimensional and just work with complex dimensions. In that case,
$A^{d}(\flags)$ is spanned by the set of all Schubert classes
$[X_{w}]$ for $w \in S_{n}$ such that $\coinv (w)=d$. Perhaps the only
place this technicality matters is when determining the Euler
characteristic.  The functions $[n]_{q^{2}}!$ and $[n]_{q}!$ are not
equal, and they remain not equal when $q=-1$.  In fact, $[n]_{q=-1}!=0$
for $n\geq 2$.  
\end{Remark}

Next we consider multiplication in the Chow ring of the flag variety.  
Recall $X_{w_{0}}$ is the whole flag variety itself, so
$[X_{w_{0}}][X_{w}]=[X_{w}]$ in the Chow ring.  In other words,
$[X_{w_{0}}]$ plays the role of the identity element in the ring
$A^{*}(\flags)$.

The Schubert class $[X_{\id}]$ also plays a special role in the Chow
ring.  The intersection conditions in \Cref{def:Schubert.variety}
imply $X_{\id}(E_{\bullet})$ contains exactly the one flag $E_{\bullet
}$.  We know $GL_{n}(\mathbb{C})$ acts transitively on $\flags$, so
every flag is rationally equivalent to every other flag as a point in
the flag variety.  Hence, multiplication in the Chow ring is closely
related to the 0-dimensional Schubert problems discussed in
\Cref{sub:PermutationArrays}.  The following theorem equates the
geometry of Schubert intersection problems to the algebraic problem of computing
product expansions of Schubert classes in the Schubert basis of Chow
rings. \textbf{This is the fundamental advantage of using Chow rings!}
The proof follows directly from the discussion.

\begin{Theorem}[\textbf{Intersection to Chow Ring Multiplication}]\label{thm:reduction.to.Chow.products}
The number of flags in the intersection
\[
X=X_{w^{1}}(E^{1}_{\ci})\cap X_{w^{2}}(E^{2}_{\ci}) \cap \dotsb \cap
X_{w^{d}}(E^{d}_{\ci}),
\]
assuming the reference flags are chosen generically and the
intersection is 0-dimensional, is exactly the coefficient of
$[X_{\id}]$ in the class $[X]=[X_{w^{1}}][X_{w^{2}}]\cdots [X_{w^{d}}]$
when expanded in the basis of Schubert classes.  
\end{Theorem}

\subsection{Schubert Structure Constants}\label{sub:cuvw}

The \textit{Schubert structure constants} $c_{uv}^{w}$ in the Chow ring
of the flag variety are the coefficients that arise when we expand
$[X_{u}][X_{v}]$ into the basis of Schubert classes $\{[X_{w}]\given w \in
S_{n} \}$.  Define the constants $c_{uv}^{w}$ by

\begin{equation}\label{eq:structure.constants}
[X_{u}][X_{v}] = \sum c_{uv}^{w} [X_{w}].
\end{equation}
As we will show, each structure constant $c_{uv}^{w}$ is a nonnegative
integer.  This is a special property of structure constants in the
Chow ring of the flag variety!  As we will see, the product
$[X_{u}][X_{v}]$ in the Chow ring is determined by the variety
$X_{u}(E_{\bullet}) \cap X_{v}(\oppositeE_{\bullet})$, known as the
\emph{Richardson variety} for $u,v \in S_{n}$. Despite the name, these
varieties had previously been studied by Kazhdan-Lusztig \cite{K-L2}.
Richardson varieties are irreducible subvarieties of the flag variety
\cite{richardson}.  The sum on the right side of
\eqref{eq:structure.constants} is not necessarily a single Schubert
class, so these irreducible intersections can be rationally equivalent
to a nontrivial union of Schubert varieties in general position since
the Chow ring of the flag manifold is spanned by Schubert classes.

Recall from \Cref{sub:PermutationArrays} that modern Schubert calculus
is the study of the intersection numbers that arise in the
intersections of Schubert varieties in the generic case.  Once one
knows the Schubert structure constants, one also knows how to expand
any product $[X_{w^{1}}][X_{w^{2}}]\cdots [X_{w^{d}}]$ in the basis of
Schubert classes.  Furthermore, we can use the following lemma to show
each Schubert structure constant $c_{uv}^{w}$ is determined by a
0-dimensional Schubert problem depending on intersecting three
Schubert varieties in sufficiently generic position.  

\begin{Lemma}\label{lem:duality}
Let $E_{\bullet}$ be the standard base flag, and let
$\oppositeE_{\bullet}=(e_{n},e_{n-1},\ldots , e_{1})$ be the flag
associated to the standard basis in reverse order.  Then
$X_{w}(E_{\bullet}) \cap X_{w_{0}w}(\oppositeE_{\bullet})$ contains
exactly one flag, namely $w_{\bullet}=(e_{w_{1}},e_{w_{2}},\ldots ,
e_{w_{n}})$.  More generally, $X_{w}(F_{\bullet}) \cap
X_{w_{0}w}(G_{\bullet})$ contains exactly one flag for any pair of
transverse flags $F_{\bullet}$ and $G_{\bullet}$.
\end{Lemma}

\begin{proof}
Consider the Schubert cells defined with respect to
$\oppositeE_{\bullet}=(e_{n},e_{n-1},\ldots , e_{1})$.  By
\Cref{def:Schubert.cell} and \Cref{ex:opposite.canonical}, a flag
$G_{\bullet} \in C_{w}(\oppositeE_{\bullet})$ can be represented by an
ordered basis $(g_{1},g_{2},\dots , g_{n})$ such that $g_{j}$ is in
the span of $\{e_{n},e_{n-1},\dots , e_{n-w_{j}+1} \}$ for each $j \in
[n]$.  Such a flag can be represented by a variation on the canonical
matrices from \Cref{def:canonical.mat} that are upside down from those
defined with respect to $E_{\bullet}$.  Each pivot 1 in such a
canonical representative has 0's to its right and above.  A flag is in
$C_{w}(\oppositeE_{\bullet})$ if the pivot 1's in its canonical matrix
with respect to $\oppositeE_{\bullet}$ are in positions
$(n-w_{j}+1,j)$, or equivalently in position $w$ if we count up from
the bottom instead of the usual matrix notation.  By matrix multiplication, one can check that the
permutation $[n-w_{1}+1,n-w_{2}+1, \dots,n-w_{n}+1]$ is equal to
$w_{0}w$.  Therefore, the flag $w_{\bullet }$ is the only flag in
$C_{w}(E_{\bullet}) \cap C_{w_{0}w}(\oppositeE_{\bullet})$.  See
\Cref{fig:S8.intersection} for a specific example.

\begin{figure}
$$
\left[
\begin{matrix}
* & * & * & * & * & * & * & 1\\
* & * & 1 & 0 & 0 & 0 & 0 & 0\\
* & * & 0 & * & * & 1 & 0 & 0\\
1 & 0 & 0 & 0 & 0 & 0 & 0 & 0\\
0 & * & 0 & * & * & 0 & 1 & 0\\
0 & 1 & 0 & 0 & 0 & 0 & 0 & 0\\
0 & 0 & 0 & * & 1 & 0 & 0 & 0\\
0 & 0 & 0 & 1 & 0 & 0 & 0 & 0\\
\end{matrix} 
 \right] \bigcap
\left[
\begin{matrix}
0 & 0 & 0 & 0 & 0 & 0 & 0 & 1\\
0 & 0 & 1 & 0 & 0 & 0 & 0 & 0\\
0 & 0 & * & 0 & 0 & 1 & 0 & 0\\
1 & 0 & 0 & 0 & 0 & 0 & 0 & 0\\
* & 0 & * & 0 & 0 & * & 1 & 0\\
* & 1 & 0 & 0 & 0 & 0 & 0 & 0\\
* & * & * & 0 & 1 & 0 & 0 & 0\\
* & * & * & 1 & 0 & 0 & 0 & 0\\
\end{matrix} 
 \right] =
\left[
\begin{matrix}
0 & 0 & 0 & 0 & 0 & 0 & 0 & 1\\
0 & 0 & 1 & 0 & 0 & 0 & 0 & 0\\
0 & 0 & 0 & 0 & 0 & 1 & 0 & 0\\
1 & 0 & 0 & 0 & 0 & 0 & 0 & 0\\
0 & 0 & 0 & 0 & 0 & 0 & 1 & 0\\
0 & 1 & 0 & 0 & 0 & 0 & 0 & 0\\
0 & 0 & 0 & 0 & 1 & 0 & 0 & 0\\
0 & 0 & 0 & 1 & 0 & 0 & 0 & 0\\
\end{matrix} 
 \right].
 $$
\caption{$C_{46287351}(E_{\bullet}) \cap C_{53712648}(\oppositeE_{\bullet})$
contains only the flag corresponding with the permutation matrix of
$46287351$.}  \label{fig:S8.intersection}
\end{figure}

We claim the cells of $X_{w}(E_{\bullet})$ of dimension strictly
smaller than $C_{w}$ do not intersect
$X_{w_{0}w}(\oppositeE_{\bullet})$.  To prove the claim, recall from
\Cref{ex:boundary} that the boundary of $X_{w}(E_{\bullet})$ is the
union of Schubert varieties $X_{v}(E_{\bullet})$ where $v=wt_{ij}$ for
some $i<j$ with $w_{i}>w_{j}$ and $\ell (v)=\ell (w)-1$. In
particular, any flag $G_{\bullet}\in X_{v}(E_{\bullet})$ would have
\[
\dim (E_{v_{i}} \cap G_{i}) \geq \mathrm{rk}(v)[v_{i},i] >
\mathrm{rk}(w)[v_{i},i], 
\]
so any matrix $M$ representing a flag $G_{\bullet} \in
X_{v}(E_{\bullet})$ would have rank strictly larger than
$\mathrm{rk}(w)[v_{i},i]$ in the northwest submatrix with lower right
corner at $(v_{i},i)$.  On the other hand, any flag in
$X_{w_{0}w}(\oppositeE_{\bullet})$ must be represented by a matrix
with rank at most $\mathrm{rk}(w)[v_{i},i]$ in the northwest submatrix
with lower right corner at $(v_{i},i)$ by considering the canonical
matrix representatives with 0's above and to the right of the pivot
1's.  So $X_{v}(E_{\bullet}) \cap X_{w_{0}w}(\oppositeE_{\bullet})
=\emptyset$.  This completes the proof of the first statement.

More generally, if $F_{\bullet }$ and $G_{\bullet }$ are transverse,
then the dimension of $F_{i}\cap G_{n-i+1}$ is exactly 1 for each
$1\leq i\leq n$.  Let $f_{i}$ be a nonzero vector in $F_{i}\cap
G_{n-i+1}$.  Then $F_{\bullet}=(f_{1},\ldots , f_{n})$ and $G_{\bullet
}=(f_{n},\ldots,f_1)$ so with respect to the $v$-basis, the generic
situation reduces to the case above.
\end{proof}

\begin{Exercise}\label{lem:nonintersect}
Say $w,y \in S_{n}$ and $\ell(w)=\ell (y)$, but $w \neq y$.  Prove 
$X_{w}(E_{\bullet}) \cap X_{w_{0}y}(\oppositeE_{\bullet}) = \emptyset$, where
$\oppositeE_{\bullet}=(e_{n},.., e_{1})$ is the reverse standard flag.
Conclude that $[X_{w}] [X_{w_{0}y}] =0$ in $A^{*}(\flags)$.
Hint: consider the relationship between $w$ and $y$ in Bruhat order.
\end{Exercise}

\begin{Lemma}\label{cor:perfect.pairing}
Assume $w,y \in S_{n}$ and $\ell (w)+\ell (y)=\binom{n }{2}$. Then we
have the simple product formula
\[
[X_{w}][X_{y}] = \begin{cases}
[X_{\id}]&  y=w_{0}w \\
0& y \neq w_{0}w. 
\end{cases}
\]
\end{Lemma}

\begin{proof}
If $y=w_{0}w $, then $X_{w}(F_{\bullet}) \cap X_{w_{0}w}(G_{\bullet})$
contains exactly one flag for any pair of transverse flags
$F_{\bullet}$ and $G_{\bullet}$ by \Cref{lem:duality}.  Every flag in
$\flags$ is an irreducible variety consisting of one point that is
rationally equivalent to the unique Schubert variety containing
exactly one point, namely $X_{\id}(E_{\bullet})$.  So, by the
definition of multiplication in the Chow ring
\[
[X_{w}][X_{w_{0}w}] = [X_{w}(F_{\bullet}) \cap X_{w_{0}w}(G_{\bullet})]=[X_{\id}].
\]
On the other hand, if $y \neq w_{0}w$, then $[X_{w}] [X_{y}] =0$
by \Cref{lem:nonintersect}.  
\end{proof}

\begin{Remark}\label{rem:positivity}
Recall in \Cref{thm:reduction.to.Chow.products} we observed that any
Schubert problem can be equated with a product of Schubert classes in
the Chow ring. The problem of expanding a product of many Schubert
classes reduces to the problem of expanding the product of any two
Schubert classes into the Schubert basis using the Schubert structure
constants recursively.  The following theorem interprets the Schubert
structure constants as solutions to Schubert problems themselves,
which means they are always nonnegative integers.  Thus, the
0-dimensional triple intersections/products are the key for all of
Schubert calculus related to the flag variety since these
intersections determine all of the structure constants for the Chow
ring.  In fact, these structure constants determine the structure
constants for the Chow rings of all partial flag varieties as well,
including the Grassmannian varieties, as we will explain in \Cref{sub:123.StepFlags}.  We say the Schubert basis in
the Chow ring of $\Fl(n)$ and any partial flag variety has the
\textit{positivity property} because $[X_{u}][X_{v}]$ always expands
into a positive integral sum of Schubert classes.
\end{Remark}

\begin{Theorem}[\textbf{Geometry Implies Positivity}]\label{cor:structure.constants}
For $u,v,w \in S_{n}$, the structure constant $c_{uv}^{w}=0$ unless
$\coinv (w)=\coinv (u)+\coinv (v)$.  Furthermore, if $\coinv
(u)+\coinv (v)=\coinv(w)$, then $c_{uv}^{w}$ is the number of flags in
the 0-dimensional intersection $X_{u}(E_{\bullet}) \cap
X_{v}(F_{\bullet}) \cap X_{w_{0}w}(G_{\bullet })$ whenever flags
$E_{\bullet}, F_{\bullet}, G_{\bullet } \in \flags$ are generic.
\end{Theorem}

\begin{proof}
Recall from \eqref{eq:structure.constants} that $c_{uv}^{w}$ is
defined by expanding the product $[X_{u}][X_{v}]$ into the Schubert
basis:
\begin{equation}\label{eq:structure.constants.2}
[X_{u}][X_{v}] = \sum c_{uv}^{w} [X_{w}].
\end{equation}
Since the Chow ring is a graded ring, $[X_{u}] \in A^{2\coinv (u)}$,
and $[X_{v}] \in A^{2\coinv (v)}$, we know the product $[X_{u}][X_{v}]$ is in 
$A^{2\coinv (u)+2\coinv (v)}$.  Every $[X_{w}]$ that appears with
nonzero coefficient in the expansion of $[X_{u}][X_{v}]$ must have
$2\coinv (w)=2\coinv (u)+2\coinv (v)$.  Therefore, $c_{uv}^{w}=0$ unless
$\coinv (w)=\coinv (u)+\coinv (v)$.

Assume $y \in S_{n}$ and $\coinv (y)=\coinv (u)+\coinv (v)$.  Multiply
both sides of \eqref{eq:structure.constants.2} by $[X_{w_{0}y}]$ to get
\[
[X_{u}][X_{v}][X_{w_{0}y}] = \sum c_{uv}^{w} [X_{w}][X_{w_{0}y}].
\]
By \Cref{lem:nonintersect}, $[X_{w}][X_{w_{0}y}]=0$ if $w \neq y$.  If
$w=y$, then \Cref{lem:duality} implies
$[X_{w}][X_{w_{0}w}]=[X_{\id}]$.  Therefore,
\begin{equation}\label{eq:triple.product}
[X_{u}][X_{v}][X_{w_{0}y}] = c_{uv}^{y} [X_{\id}].
\end{equation}
The triple product $[X_{u}][X_{v}][X_{w_{0}y}]$ can also be computed
as the class of the intersection $X_{u}(E_{\bullet}) \cap
X_{v}(F_{\bullet}) \cap X_{w_{0}y}(G_{\bullet})$ in the Chow ring,
provided the flags are chosen generically.

Since $\coinv (y)=\coinv (u)+\coinv (v)$ by assumption, we know
$\coinv (u)+\coinv (v)+\coinv (w_{0}y) =\binom{n}{2}$ and 
$X_{u}(E_{\bullet}) \cap X_{v}(F_{\bullet}) \cap
X_{w_{0}y}(G_{\bullet})$ is 0-dimensional.  So $X_{u}(E_{\bullet})
\cap X_{v}(F_{\bullet}) \cap X_{w_{0}y}(G_{\bullet})$ is a finite
number of points, which must equal $c_{uv}^{y}$ by
\eqref{eq:triple.product} and the fact that the Schubert classes form
a basis for the Chow ring of flag variety so the expansion is unique.
\end{proof}

The following symmetries among Schubert structure constants follow
from the fact that intersection is a commutative relation.  Other
symmetry relations will be presented in \Cref{ex:descent.cycling}.  

\begin{Corollary}\label{cor:sym.structure.constants}
For $u,v,w \in S_{n}$ such that 
$\coinv (w)=\coinv (u)+\coinv (v)$, we have
\[
c_{uv}^{w} = c_{vu}^{w} = c_{u,w_{0}w}^{w_{0}v}.
\]
\end{Corollary}

In the case $n=4$, all of the structure constants can be computed by
using the shoebox pictures described in \Cref{sub:flags} and the Hasse
diagram of $S_{4}$ shown in \Cref{Hasse of S4}.  Let's dig into these
computations with pictures and compute $c_{uv}^{w}$ for some
triples of permutations.  By \Cref{cor:structure.constants}, we know
the identity
\[
c_{uv}^{w}[X_{\id}]=[X_{u}(B_{\bullet}) \cap 
X_{v}(R_{\bullet}) \cap X_{w_{0}w}(G_{\bullet })]
\]
holds in $A^{*}(\Fl(4))$ for three transverse reference flags
$B_{\bullet}$, $R_{\bullet}$, and $G_{\bullet}$ in $\mathbb{C}^{4}$.
We can't draw $\mathbb{C}^{4}$, so we consider a shadow of this
elusive space by considering the projective picture of
$\mathbb{R}^{4}$ into $\mathbb{R}^{3}$, which we try to draw on this
2-dimensional platform.  Hopefully you can see how these pictures
suffice to make the necessary calculations. 

To compute Schubert structure constants for $n=4$, draw $B_{\bullet}$
in black, $R_{\bullet}$ in red, and $G_{\bullet}$ in green.  We saw in
\Cref{fig:one.blue.flag} how to draw two transverse flags.  Extend
this to three flags in transverse position as in
\Cref{fig:three.transverse.flags}.  Here each flag has its special
point on its special line in its special plane, and they all sit in
the same shoebox.  The planes of three flags $B_{3}, R_{3}, G_{3}$ are
not necessarily orthogonal, and they do not necessarily meet at the
origin in this projective picture.  In fact, if they were orthogonal
or met at the origin, that would not be generic, but they would be
transverse.

\begin{figure}[h]
\begin{tikzpicture}[scale=0.6]
\definecolor{MyGreen}{RGB}{145,205,50}
\draw[red, fill=red, fill opacity=0.1](0,0)--(-1.5,-1.5)--(2.5,-1.5)--(4,0);
\draw[MyGreen, fill=MyGreen, fill opacity=0.1](0,0)--(-1.5,-1.5)--(-2.5,0)--(-1,1.5);
\draw[black, fill=black, fill opacity=0.1](0,0)--(4,0)--(3,1.5)--(-1,1.5)--(0,0);
\draw[red] (-1,-1)--(3,-0.2);
\node[red] at (1,-0.6) {$\bullet$};
\draw[MyGreen] (-1.5,-1)--(-1,1);
\node[MyGreen] at (-1.3,-0.2) {$\bullet$};
\draw (2,0.2)--(0.5,1.3);
\node at (1.25,0.75) {$\bullet$};
\end{tikzpicture}
\caption{Projective representation of three transverse flags in $\mathbb{C}^{4}$ looking into a shoebox bounded by the three planes.}
\label{fig:three.transverse.flags}
\end{figure}
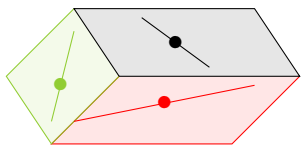

\begin{figure}[h]
\[
\begin{array}{ccccc}
\left[
\begin{array}{rrrr}
0 & 0 & 0 &	1\\
0 & 0 & 1 &	2\\
1 & 1 & 2 &	3 \\
1 & 2 & 3 &	4
\end{array}
 \right]
& \hspace{.3in} & 
\left[
\begin{array}{rrrr}
0 & 0 & 0 &	1\\
0 & 1 & 1 &	2\\
0 & 1 & 2 &	3 \\
1 & 2 & 3 &	4
\end{array}
 \right]
& \hspace{.3in} & 
\left[
\begin{array}{rrrr}
0 & 1 & 1 &	1\\
0 & 1 & 2 &	2\\
1 & 2 & 3 &	3 \\
1 & 2 & 3 &	4
\end{array}
 \right]\\
{\vspace{0in}} &	{} &	{} &	{} &	{} 
\\
\mathrm{dim}(B_{i}\cap F_{j}) &	&	\mathrm{dim}(R_{i}\cap F_{j})
&		       &      \mathrm{dim}(G_{i}\cap F_{j})
\end{array}
\]
\caption{Tables of intersection dimensions for $F_{\bullet} \in
X_{3421}(B_{\bullet})$ on the left, 
$F_{\bullet} \in X_{4231}(R_{\bullet})$ in the middle, and 
$F_{\bullet} \in X_{3124}(G_{\bullet })$ on the right from \Cref{ex:cuvw.calc}.}
\label{fig:tables.in.dim}
\end{figure}

\begin{Example}\label{ex:cuvw.calc}
Let $u=3421$, $v=4231$, and $w=2431$.  We want to prove $c_{uv}^{w}=1$
in this case. Observe that $2=\coinv (3421)+\coinv (4231)=\coinv
(2431)$, so it is possible $c_{uv}^{w}$ is not zero by
\Cref{cor:structure.constants}.  Compute $w_{0}w=3124$.  Any 
flag 
\begin{equation}\label{eq:3421.4231.3124}
F_{\bullet} \in X_{3421}(B_{\bullet}) \cap X_{4231}(R_{\bullet}) \cap
X_{3124}(G_{\bullet })
\end{equation}
must satisfy the intersection conditions given by the tables in
\Cref{fig:tables.in.dim}.

\bigskip

\noindent 
Let's translate that into conditions on the possible shoebox
drawings containing $F_{\bullet}$, along with the black, red, and
green flags.
\begin{enumerate}
\item From $\mathrm{dim}(B_{3}\cap F_{1})=1$ we see $F_{1}$ is
special, it must be a one dimensional subspace of $B_{3}$.  In the
drawing, this means the point of any $F_{\bullet}$ must be drawn in
the black plane.
\medskip
\item From $\mathrm{dim}(R_{2}\cap F_{2})=1$ we see $F_{2}$ is
special, it must intersect $R_{2}$ in a one dimensional subspace.  So,
the line for $F_{\bullet}$ will intersect the red line.
\medskip
\item From $\mathrm{dim}(G_{1}\cap F_{2})=1$, we see $F_{2}$ contains
$G_{1}$.  From $\mathrm{dim}(G_{3}\cap F_{3})=3$, we see
$F_{3}=G_{3}$.  So, the plane for $F_{\bullet}$ is exactly the green
plane, and the line for $F_{\bullet}$ must lie in the green plane and
pass through the green point.
\end{enumerate}
These conditions together uniquely determine a flag, hence
$c_{uv}^{w}=1$.  We draw this unique flag $F_{\bullet}$ in purple in
\Cref{fig:3421cap4231cap2431}.  The purple plane agrees with the green
plane.  The purple line representing $F_{2}$ must be in the green
plane, go through the green dot, and go through the point of
intersection between the red line with the green plane.  The purple
line meets the black plane in one point, which must be the purple
point representing $F_{1}$.

\begin{figure}[h] 
\begin{tikzpicture}[scale=0.6]
\definecolor{MyGreen}{RGB}{145,205,50}
\definecolor{MyPurple}{RGB}{142,69,133}
\definecolor{SomeColor}{RGB}{143,137,91}
\draw[red, fill=red, fill opacity=0.1](0,0)--(-1.5,-1.5)--(2.5,-1.5)--(4,0);
\draw[MyGreen, fill=SomeColor, fill opacity=0.1](0,0)--(-1.5,-1.5)--(-2.5,0)--(-1,1.5);
\draw[black, fill=black, fill opacity=0.1](0,0)--(4,0)--(3,1.5)--(-1,1.5)--(0,0);
\draw[red] (-1,-1)--(3,-0.2);
\node[red] at (1,-0.6) {$\bullet$};
\draw[MyPurple] (-0.1,0)--(-1.5,-1.4)--(-2.4,0)--(-1,1.4)--(-0.1,0);
\draw[MyPurple] (-1,-1)--(-2/3,1);
\node[MyPurple] at (-2/3,1) {$\bullet$};
\draw[MyGreen] (-5/6+0.5,0.2)--(-5/6-1,-0.4);
\node[MyGreen] at (-5/6,0) {$\bullet$};
\draw (2,0.2)--(0.5,1.3);
\node at (1.25,0.75) {$\bullet$};
\end{tikzpicture}
\caption{The purple flag represents the unique flag
$X_{3421}(B_{\bullet}) \cap X_{4231}(R_{\bullet}) \cap
X_{3124}(G_{\bullet })$.}
\label{fig:3421cap4231cap2431}
\end{figure}
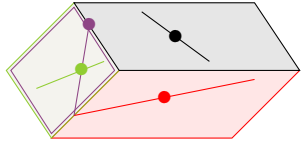
\end{Example}
\bigskip

\begin{Example}\label{ex:3421.4231.3124}
Let $u=3421$, $v=4231$, and $w=4132$.  We want to prove $c_{uv}^{w}=0$
in this case. Observe that $2=\coinv (3421)+\coinv (4231)=\coinv
(4132)$, so it is possible $c_{uv}^{w}$ is not zero by
\Cref{cor:structure.constants}.  Compute $w_{0}w=1423$.  Any 
flag $F_{\bullet}$ in 
\begin{equation}\label{eq:3421.4231.1423}
X_{3421}(B_{\bullet}) \cap X_{4231}(R_{\bullet}) \cap
X_{1423}(G_{\bullet })
\end{equation}
must satisfy the intersection conditions given by the following
tables, 

\[
\begin{array}{ccccc}
\left[
\begin{array}{rrrr}
0 & 0 & 0 &	1\\
0 & 0 & 1 &	2\\
1 & 1 & 2 &	3 \\
1 & 2 & 3 &	4
\end{array}
 \right]
& \hspace{.3in} & 
\left[
\begin{array}{rrrr}
0 & 0 & 0 &	1\\
0 & 1 & 1 &	2\\
0 & 1 & 2 &	3 \\
1 & 2 & 3 &	4
\end{array}
 \right]
& \hspace{.3in} & 
\left[
\begin{array}{rrrr}
1 & 1 & 1 &	1\\
1 & 1 & 2 &	2\\
1 & 1 & 2 &	3 \\
1 & 2 & 3 &	4
\end{array}
 \right]\\
% \left[      %%% Sara Changed this on 11.14.2023.  It w_{0}w of what
% it should be which is the table for 1423
% \begin{array}{rrrr}
% 0 & 1 & 1 &	1\\
% 0 & 1 & 2 &	2\\
% 0 & 2 & 3 &	3 \\
% 1 & 2 & 3 &	4
% \end{array}
%  \right]\\
{\vspace{0in}} &	{} &	{} &	{} &	{} 
\\
\mathrm{dim}(B_{i}\cap F_{j}) &	&	\mathrm{dim}(R_{i}\cap F_{j})
&		       &      \mathrm{dim}(G_{i}\cap F_{j}).
\end{array}
\]
\bigskip

\noindent 
Now, $\mathrm{dim}(B_{3}\cap F_{1})=1$ and $\mathrm{dim}(G_{1}\cap
F_{1})=1$ is not possible if $B_{\bullet}$ and $G_{\bullet}$ are in
transverse position.  That would imply the purple dot is both in the
black plane and equal to the green dot, which is not the case in the
drawing of three flags in transverse position.  Hence
$c_{3421,4231}^{4132}=0$.   The same type of argument also shows
$c_{3421,4231}^{4213}=0$. 
\end{Example}
\bigskip

We turn now to the computation of $c_{uv}^{w}$ in some special cases,
but for arbitrary $n$.  These are the celebrated Monk's formula and
 Pieri rule.  David Monk's paper from 1959 is a gem \cite{MONK}. His work
builds on important contributions from Borel, Ehresmann, Hirzebruch, and 
Hodge. It follows a line of reasoning counting flags in 0-dimensional
triple intersections.  Along similar lines, the Pieri formula for the
Grassmannian was generalized to the flag variety by Frank Sottile in
1996 \cite{sottile}.  This was a huge leap in Schubert calculus that
led to many further developments including the development by
Bergeron-Sottile of the $k$-Bruhat order \cite{BergSott}, a
generalization to quantum cohomology by Ciocan-Fontanine
\cite{Ciocan-Fontanine} and Postnikov \cite{Post}, and to more exotic
cohomology theories \cite{lenart-sottile,Li-et-al,Liu.22}.  We will
return to the Pieri formula in \Cref{thm:sottile} below.

Monk considered the family of Schubert varieties $X_{w}(E_{\bullet})$
in the special case $w=w_{0}s_{i}$ for some $1\leq i<n$.  Recall,
$s_{i}$ is the simple transposition that interchanges $i,i+1$.  The
varieties $X_{w_{0}s_{i}}$ are hypersurfaces in $\flags$.  In one-line
notation $w_{0}s_{i}$ is a decreasing sequence interrupted by one
ascent after position $i$,
\[
w_{0}s_{i}=[n,n-1,\dots , n-i+2, n-i,n-i+1, n-i-1,\dots, 1]. 
\]
We know $X_{w_{0}}$ is the whole flag variety, and so each
$X_{w_{0}s_{i}}$ has codimension 1. To begin, we consider the
question, 

\medskip
\begin{quote}
Question: What flags are in $X_{w_{0}s_{i}}(E_{\bullet})$?
\end{quote}

\medskip
\begin{Example}\label{ex:n=5}
Let $n=5$, \ and $w=w_{0}s_{2}=53421$.  A flag $F_{\bullet}$ in
$X_{w}(E_{\bullet})$ must have $\mathrm{dim}(E_{3}\cap F_{2})\geq 1$
by definition.  This condition is also sufficient to prove
$F_{\bullet} \in X_{w}(E_{\bullet})$, because this is the only binding
rank condition.  Thus, for every $F_{\bullet} \in X_{w}(E_{\bullet})$
there exists an ordered basis
$F_{\bullet}=(f_{1},f_{2},f_{3},f_{4},f_{5})$ such that $f_{1}$ or
$f_{2}$ in $E_{3}$.
\end{Example}
\bigskip

Monk describes the flags in $X_{w_{0}s_{i}}(E_{\bullet})$ as the set
of flags with $F_{i}$ chosen from a \textit{pencil} of $i$-planes over
$E_{n-i}$. This evokes the image of a pencil of lines through a point,
meaning all lines that go through one point in the plane.  If we
choose a different flag $G_{\bullet}$ as the reference flag, then
$X_{w_{0}s_{i}}(G_{\bullet})$ is still the set of flags such that
$F_{i}$ intersects $G_{n-i}$ in a subspace of dimension at least 1,
which would not happen if $F_{i}$ was chosen generically.  Since
$X_{w_{0}s_{i}}(G_{\bullet})$ is relatively easy to describe, Monk was
able to give a simple formula for certain products of Schubert classes
in the Chow ring of the flag variety.

\begin{Theorem}[\textbf{Monk's Formula}] \cite[Thm.3]{MONK}
\label{thm:monk}
Let $w \in S_{n}$ and $i\in [n-1]$, then 

\begin{equation}\label{eq:monk}
[X_{w_{0}s_{i}}][X_{w}] = \sum_{\substack{
h\leq i<j:
\\
\ell(w t_{hj}) = \ell(w) -1
}} [X_{wt_{hj}}].
\end{equation}
\end{Theorem}

It is quite remarkable that the coefficients on the right side of
Monk's formula in \eqref{eq:monk} are all 0 or 1.  In general, the
Schubert structure constants can be large positive integers.

\begin{Example}\label{ex:45132.2}
Let $n=5$,\ $w=45132$, and $i=2$.  Monk's formula says to consider
$w=45|132$ separated after the second position.  Find all pairs of
positions $h\leq 2<j$ such that $w_{h}>w_{j}$ and all values in
between positions $h$ and $j$ are not in the range $[w_{j},w_{h}]$.  This
ensures $\ell(wt_{h,j})=\ell(w)-1$.  For every such pair, we get a
term in the expansion so 
\[
[X_{53421}][X_{45|132}] =
[X_{15|432}]+[X_{35|142}]+[X_{41|532}]+[X_{43|152}].
\]
Note the value 2 is not swapped with $4$ or $5$ since the 3 is in
between.  
\end{Example}
\bigskip

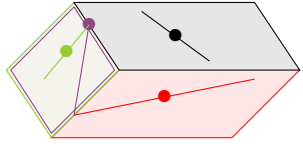
\begin{figure}
\begin{tikzpicture}[scale=0.6]
\definecolor{MyGreen}{RGB}{145,205,50}
\definecolor{MyPurple}{RGB}{142,69,133}
\definecolor{SomeColor}{RGB}{143,137,91}
\draw[red, fill=red, fill opacity=0.1](0,0)--(-1.5,-1.5)--(2.5,-1.5)--(4,0);
\draw[MyGreen, fill=SomeColor, fill opacity=0.1](0,0)--(-1.5,-1.5)--(-2.5,0)--(-1,1.5);
\draw[black, fill=black, fill opacity=0.1](0,0)--(4,0)--(3,1.5)--(-1,1.5)--(0,0);
\draw[red] (-1,-1)--(3,-0.2);
\node[red] at (1,-0.6) {$\bullet$};
\draw[MyPurple] (-0.1,0)--(-1.5,-1.4)--(-2.4,0)--(-1,1.4)--(-0.1,0);
\draw[MyPurple] (-1,-1)--(-2/3,1);
\node[MyPurple] at (-2/3,1) {$\bullet$};
\draw[MyGreen] (-2/3,1)--(-5/3,-0.2);
\node[MyGreen] at (-7/6,0.4) {$\bullet$};
\draw (2,0.2)--(0.5,1.3);
\node at (1.25,0.75) {$\bullet$};
\end{tikzpicture}
\caption{The purple flag represents the unique flag in position $3421$
with respect to the black flag, in position $4231$ with respect to the
red flag, and in position $2314$ with respect to the green flag. }
\label{fig:2314}
\end{figure}

\begin{Example}\label{ex:3421.4231}
Let $n=4$,\ $w=3421$, and $i=2$.  Monk's formula implies 
\[
[X_{4231}][X_{34|21}]= [X_{24|31}]+[X_{32|41}].
\]
Note, $3421=w_{0}s_{1}$, so one can compute the same expansion by
looking for transpositions reducing the number of inversions by 1 over
the bar in $4|231$ as well.  We can confirm this expansion
by considering shoebox diagrams.  We only need to consider terms in
the expansion of the form $[X_{v}]$ with
$\ell(v)=\ell(w)-1=5-1=4$. The permutations with 4 inversions are 
\[
4213,\ 4132,\ 3412,\ 3241,\ 2431.
\]
We know from \Cref{ex:cuvw.calc} that $c_{3421,4231}^{2431}=1$.  Since
$A^{*}(\flags)$ is commutative, $c_{4231,3421}^{2431}=1$. We also know
from \Cref{ex:3421.4231.3124} that $c_{3421,4231}^{4132}=0$ and
$c_{3421,4231}^{4213}=0$.  It remains to show $c_{3421,4231}^{3241}=1$
and $c_{3421,4231}^{3412}=0$.

To see $c_{3421,4231}^{3241}=1$, we follow the same procedure as in
\Cref{ex:cuvw.calc}.  Compute $w_{0}3241=2314$.  Any 
flag $F_{\bullet}$ in 
\begin{equation}\label{eq:3421.4231.2314}
X_{3421}(B_{\bullet}) \cap X_{4231}(R_{\bullet}) \cap
X_{2314}(G_{\bullet })
\end{equation}
must satisfy the intersection conditions given by the following
tables, 

\[
\begin{array}{ccccc}
\left[
\begin{array}{rrrr}
0 & 0 & 0 &	1\\
0 & 0 & 1 &	2\\
1 & 1 & 2 &	3 \\
1 & 2 & 3 &	4
\end{array}
 \right]
& \hspace{.3in} & 
\left[
\begin{array}{rrrr}
0 & 0 & 0 &	1\\
0 & 1 & 1 &	2\\
0 & 1 & 2 &	3 \\
1 & 2 & 3 &	4
\end{array}
 \right]
& \hspace{.3in} & 
\left[
\begin{array}{rrrr}
0 & 0 & 1 &	1\\
1 & 1 & 2 &	2\\
1 & 2 & 3 &	3 \\
1 & 2 & 3 &	4
\end{array}
 \right]\\
{\vspace{0in}} &	{} &	{} &	{} &	{} 
\\
\mathrm{dim}(B_{i}\cap F_{j}) &	&	\mathrm{dim}(R_{i}\cap F_{j})
&		       &      \mathrm{dim}(G_{i}\cap F_{j}).
\end{array}
\]
This translates into exactly one flag $F_{\bullet}$ drawn in purple in
\Cref{fig:2314}.  Here $F_{1}=B_{3}\cap G_{2}$ so the purple dot is at
the intersection of the green line and the black plane. The purple
line goes through the purple dot and the intersection of the red line
with the green plane since $F_{1} \subset F_{2}$ and
$\mathrm{dim}(R_{2}\cap F_{2})\geq 1$.  Finally, $\mathrm{dim}(G_{3}\cap
F_{3})\geq 3$, so the green plane and the purple plane agree.

The final calculation for $c_{3421,4231}^{3412}=0$ follows the same
procedure. Here $w_{0}3412=2143.$ Writing out the
intersection tables for $F_{\bullet}$ in the triple intersection we
see $F_{2}=G_{2}$ and $\mathrm{dim}(F_{2}\cap R_{2})\geq 1$.  Since
$G_{\bullet}$ and $R_{\bullet}$ are assumed to be in transverse
position, these two conditions cannot happen.  So no flags are in the
triple intersection $ X_{3421}(B_{\bullet}) \cap X_{4231}(R_{\bullet})
\cap X_{2143}(G_{\bullet })$.  
\end{Example}
\bigskip

To outline the proof of Monk's formula, we need to find the number of
flags in the appropriate triple intersections of Schubert varieties
with respect to 3 generic flags by \Cref{cor:structure.constants}.
For any 3 generic flags, we can always find a change of basis for
$\mathbb{C}^{n}$ to move two of the flags so that they become
$E_{\bullet}$ and $\oppositeE_{\bullet}$.  Thus, it suffices to
compute $c_{w,w_{0}s_{i}}^{v}$ for all $v\in S_{n}$ by finding the
number of flags in a 0-dimensional triple intersection of the form
\begin{equation}\label{eq:monk.proof}
X_{w}(E_{\bullet}) \cap X_{w_{0}v}(\oppositeE_{\bullet}) \cap
X_{w_{0}s_{i}}(G_{\bullet})
\end{equation}
assuming $G_{\bullet}$ is generic and $\ell(v)=\ell(w)-1$.  The proof
then follows from the following statements and properties of the
Bruhat order from \Cref{sub:bruhat}.  Note, \Cref{ex:Monk.nonempties}
is significantly harder than the other exercises.  Monk uses a
characterization of the irreducible hyperplane sections of $\flags$
and an interesting identity he proves on multinomial coefficients.

\begin{Exercise}\label{ex:Monk.v<=w}
For $u,w \in S_{n}$, the Richardson variety $X_{u}^{w}=
X_{w}(E_{\bullet}) \cap X_{w_{0}u}(\oppositeE_{\bullet})$ is nonempty
if and only if $u\leq w$.  In particular, $v_{\bullet} \in X_{u}^{w}$
if and only if $v \in [u,w]$.
\end{Exercise}

\begin{Exercise}\cite[\S 9, Prop 5]{MONK} \label{ex:Monk.Prop.5}
For $w \in S_{n}$ and $v=wt_{hj}$ covered by $w$ in Bruhat order,
the intersection $X_{w}(E_{\bullet}) \cap
X_{w_{0}v}(\oppositeE_{\bullet})$ consists of flags $F_{\bullet}$ such
that
\begin{enumerate}
\item the subspaces $F_{i}$ for $h\leq i<j$  must contain the span of
$\{e_{1},\dots , e_{h-1},e_{h+1},\dots , e_{i} \}$ and must be
contained in the span of $\{e_{1},\dots , e_{k},e_{j} \}$, and 
\item the initial and final subspaces in $F_{\bullet}$ are given by
$F_{i}=E_{i}$ for $1\leq i<h$ and $j\leq i\leq n$.
\end{enumerate}
\end{Exercise}

\begin{Exercise}\label{ex:Monk.empties}
For any generic flag $G_{\bullet}$ and $v=wt_{h,j}$ such that
$\ell(v)=\ell(w)-1$, the intersection in \eqref{eq:monk.proof} is
empty unless $h\leq i<j$.
\end{Exercise}

\begin{Exercise}\cite[Thm. 1]{MONK}\label{ex:Monk.degree}
Prove
\[
 \left(\sum_{k=1}^{n-1} [X_{w_{0}s_{k}}] \right)^{\binom{n}{2}}= \binom{n}{2}! \ [X_{\id}].
\]
\end{Exercise}

\begin{Exercise}\cite[\S 11]{MONK}\label{ex:Monk.nonempties}
For any generic flag $G_{\bullet}$ and $v=wt_{h,j}$ such that $h\leq
i<j$ and $\ell(v)=\ell(w)-1$, the intersection in
\eqref{eq:monk.proof} contains exactly one flag. Hint: Positivity in
this case is implied by \Cref{ex:Monk.v<=w}, so the task is to prove
there is not more than 1 flag in the intersection.  Assume there are
at least 2 flags in some such intersection so the coefficients in the
expansion of $[X_{w_{0}s_{i}}][X_{w}]$ in \eqref{eq:monk} are not
multiplicity free, then prove that then the coefficient would be
strictly larger than $\binom{n}{2}!$, contradicting
\Cref{ex:Monk.degree}.  
\end{Exercise}

\begin{Remark}
Around the same time as Monk was working on his formula, Claude
Chevalley found another proof.  His work went unpublished until 1994
and was posthumously published by Borel in \cite{chevalley}.  This
paper was highly influential even though only a few copies existed
before it was published.  The reader may wish to compare Monk's proof
to Chevalley's proof for the final details of proving
\Cref{ex:Monk.nonempties}.  
\end{Remark}

% Note the terms of the identities on multinomial coeffs and
% the determinant of $(i^{j})$. Where else does this last part of the
% proof appear in the literature? Ehrenborg's Hankel Determinant paper
% from 2000 for example.  See A000178 Superfactorials: product of first
% n factorials, but Monk's paper not mentioned there.  Relate to the
% degrees of Schubs in Postnikov-Stanley, Prop. 2.1.  }

Miraculously, Monk's rule for multiplying by the special Schubert
classes suffices to compute the structure constants for the Chow ring
of the full flag variety.  Every Schubert class can be expanded as a
sum of products of special classes.

\begin{Theorem}\label{cor:special.classes.generate.ring}
The set $\{[X_{w_{0}s_{i}}]\given 1\leq i<n\}$ generates the Chow ring
$A^{*}(\flags)$.
\end{Theorem}

\begin{proof}
Consider the subring $R$ generated by the special classes
$[X_{w_{0}s_{i}}]$.  Since $\{[X_{w}] \given w \in S_{n}\}$ is a linear
basis for $A^{*}(\flags)$, we only need to show that every $[X_{w}]$ is in $R$
to prove the statement. The empty product of generators in a ring is
defined to be the multiplicative identity.  In $A^{*}(\flags)$, we
know $1=[X_{w_{0}}]$, so $[X_{w_{0}}] \in R$ and each
$[X_{w_{0}s_{i}}] \in R$.  Assume by induction that $[X_{w}] \in R$
for all $w$ such that $\binom{n}{2}\geq \ell (w)> k$ for some $k$.   

Let $v \in S_{n}$ be a permutation with $\ell(v)=k$. Since
$k<\binom{n}{2}$, $v$ has an ascent in one-line notation.  Let $r$ be
the position of the last ascent in $v$, so
$v_{r}<v_{r+1}>v_{r+2}>\dots >v_{n}$.  Let $s$ be the largest value
such that $v_{r}<v_{s}$, so $r<s\leq n$.  The pair $(r,s)$ is the
lexicographically largest (lex) non-inversion for $v$. Let
$w=vt_{rs}$, then $\ell(w)=k+1$.  By Monk's formula,
\begin{equation}\label{eq:generate.v}
[X_{w_{0}s_{r}}][X_{w}] =[X_{v}] + \sum_{\substack{ h\leq  r<j:
\\
\ell(w t_{hj}) = k\\
(r,s) \neq (h,j) }} [X_{wt_{hj}}].
\end{equation}

If $r>1$, the product $[X_{w_{0}s_{r-1}}][X_{w}]$ has an
expansion that is very similar to \eqref{eq:generate.v} but does not
include $[X_{v}]$.  After cancellation and rearranging terms, one can
conclude
\begin{equation}\label{eq:generate.v.1}
[X_{v}] = [X_{w_{0}s_{r}}][X_{w}] - [X_{w_{0}s_{r-1}}][X_{w}]  + \sum_{\substack{ h< r:
\\
\ell(w t_{hr}) = k }} [X_{wt_{hr}}] - \sum_{\substack{ r<j\neq s:
\\
\ell(w t_{rj}) = k }} [X_{wt_{rj}}].
\end{equation}
If $s<j$, then all permutations of the form $wt_{rj}$ have
$\ell(wt_{rj})<k$ since  $w$ is
decreasing after position $r$ so $w_{r}>w_{s}=v_{r}>w_{j}=v_{j}$.  If $r<j<s$, then
$w_{r}=v_{s}<v_{j}=w_{j}$, so $\ell(wt_{rj})=\ell(w)+1=k+2$.
Therefore the negative sum is empty, and
\begin{equation}\label{eq:generate.v.2}
[X_{v}] = [X_{w_{0}s_{r}}][X_{w}] - [X_{w_{0}s_{r-1}}][X_{w}]  + \sum_{\substack{ h< r:
\\
\ell(w t_{hr}) = k }} [X_{wt_{hr}}].
\end{equation}
Furthermore, if $(r',s')$ is the
largest pair of positions corresponding with a non-inversion in
$wt_{hr}$ indexing a term in the positive sum, then either $(r'<r)$ or
($r'=r$ and $s'<s$).  Hence, by a second induction over lex order on
the largest non-inversions of permutations with $k$ inversions, all
terms on the right side of \eqref{eq:generate.v.2} are in $R$, hence
$[X_{v}] \in R$.

If $r=1$, then the sum on the right of \eqref{eq:generate.v} is empty
by an argument similar to the negative sum above.  Hence, $[X_{v}] \in R$ since
$[X_{w_{0}s_{i}}][X_{w}]$ is in $R$ by induction.
\end{proof}

Monk clearly recognized that the special Schubert classes
$[X_{w_{0}s_{i}}]$ generated the cohomology ring of the flag variety.
He also shows constructively that the ring is isomorphic to the
coinvariant algebra by relating the special classes $[X_{w_{0}s_{i}}]$
with polynomial representatives in the quotient modulo $I_{n}^{+}$.
Monk argued that if $x_{i}=[X_{w_{0}s_{i-1}}]-[X_{w_{0}s_{i}}]$ for
$1<i<n$,\ $x_{1}=-[X_{w_{0}s_{1}}]$, and $x_{n}=[X_{w_{0}s_{n-1}}]$,
then the corresponding elementary symmetric polynomials
$e_{k}(x_{1},\dots , x_{n})$ vanish, which is well known to be another
generating set for $I_{n}^{+}$ in symmetric function theory. This
reproves Borel's theorem if by a great leap of faith one assumes the
Chow ring and the singular cohomology ring of the flag variety are
isomorphic, or one studies intersection theory carefully.  In a sense,
Monk's formula determines the entire ring structure for the cohomology
ring of a flag variety!

The real power of Monk's work is that it provides a concrete path
toward identifying polynomials in the coinvariant algebra that
represent the Schubert classes.  In general, it is much easier to
multiply polynomials and expand them in the basis of polynomials than
it is to work out the details of finding the expansion of a class into
Schubert classes via rational equivalence.  Working modulo the ideal
$I_{n}^{+}$ provides some challenges, but is possible using the theory
of Gr\"obner bases.  \bigskip

How should the map from $A^{*}(\flags)$ to $R_{n}$ be constructed?  We
need to identify polynomials that represent the Schubert classes and
form an independent set in the coinvariant algebra
$R_{n}=\mathbb{Z}[x_{1},x_{2}, \dots , x_{n}]/I_{n}^{+}$.  Plus, the
map should respect the grading on both rings. It is very instructive
to try this on your own before reading the construction below.

\bigskip

\begin{figure}[h]
\centering
\includegraphics[height=15cm]{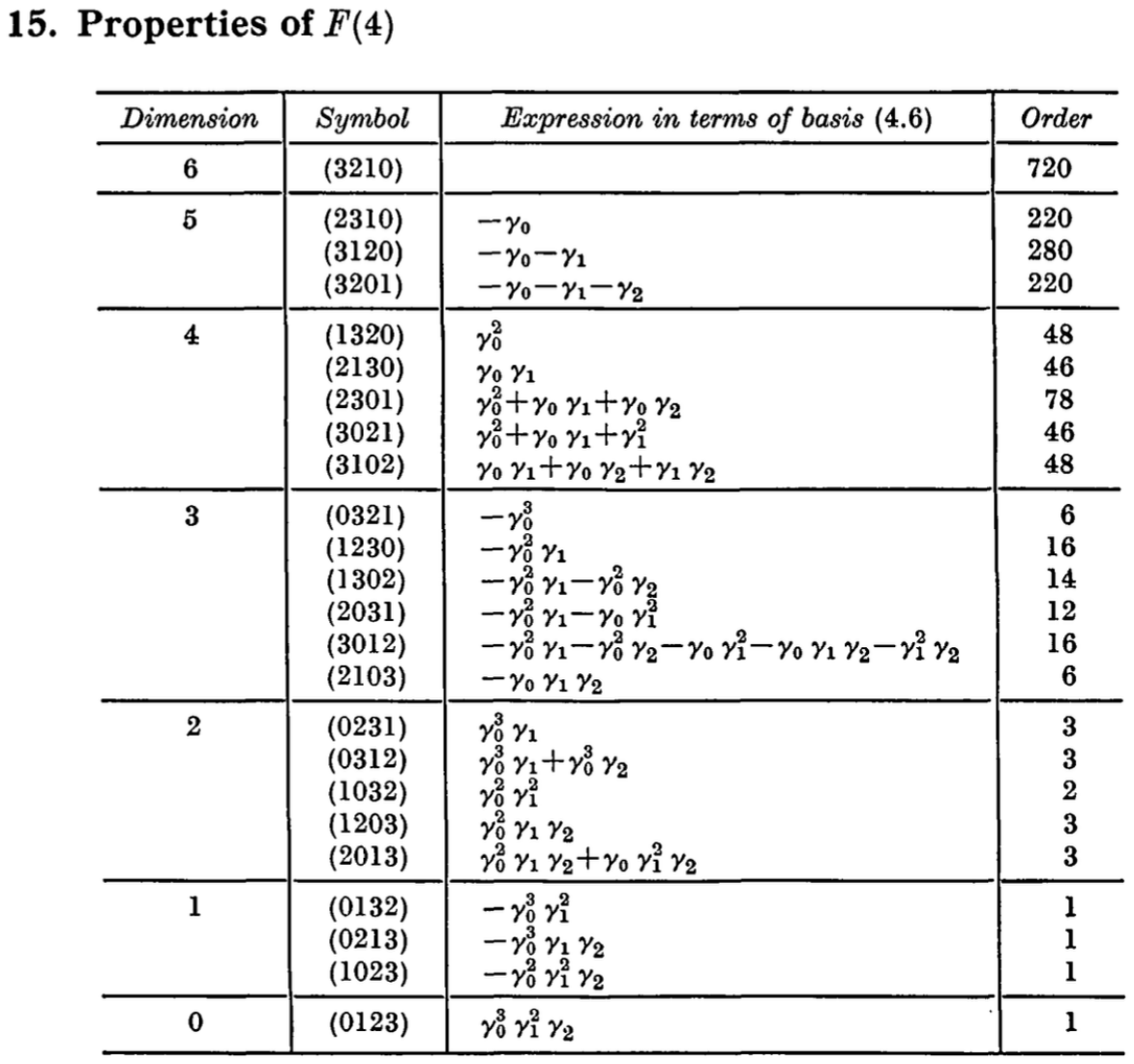}
\caption{Monk's representatives of Schubert classes in
$\mathbb{Z}[\gamma_{1},\dots , \gamma_{n}].$
He writes permutations in $S_{n}$ as bijections on the set
$\{0,1,2,\dots, n-1 \}$. The ``order'' column gives the degree of the
corresponding Schubert variety as a projective variety when one embeds
$\flags$ into a product of projective spaces via the Pl\"ucker
embedding (see \Cref{sub:Plucker}), and then into a single projective space via the Segre embedding.
Formulas for these degrees can be found in  \cite{duan-degree} and \cite{postnikov-stanley}.
}
\label{fig:Monk.Table}
\end{figure}

Given the way multiplication works in any graded ring, a Schubert
class $[X_{w}] \in A^{2d}$ should map to a polynomial of degree $d$ in
$\mathbb{Z}[x_{1},x_{2},\ldots , x_{n}]$.  We consider
$\mathrm{deg}(x_{i})=1$ for each $i$.  So, the complex codimension of
the variety should be the degree of the polynomial representing its
class.  Monk suggested mapping
\[
[X_{w_{0}s_{i}}] \mapsto -(x_{1}+x_{2}+ \cdots + x_{i}), 
\]
but he notes there is an ``ambiguity of sign'' \cite[Sect. 10]{MONK}.
His computations in the case of $n=4$ are given in
\cite[Sect. 15]{MONK}.  See \Cref{fig:Monk.Table}.  He notes that the
stated polynomial representatives alternate in sign depending on the
degree.  However, then Monk abandons this path toward realizing the
full power of the polynomial map because of the negative signs.

\bigskip

\textbf{What if Monk had made a different choice when he noticed the
``ambiguity of sign''?  Could it lead to a nicer family of
polynomials?  }

\bigskip

We now outline another natural choice for the polynomial
representatives for the special classes $[X_{w_{0}s_{i}}]$, which has
opened many doors in mathematical research.  In a ``revisionist
history'', we explore what could have been done in the 1950's.

Assume $f:A^{*}(\flags) \longrightarrow R_{n}=\mathbb{Z}[x_{1},\dots ,
x_{n}]/I_{n}^{+}$ is a ring homomorphism defined so that $[X_{w_{0}}]=1$ and for each
$i \in [n-1]$,
\[
[X_{w_{0}s_{i}}] \mapsto x_{1}+x_{2}+ \cdots + x_{i}.
\]
Then, by \eqref{eq:generate.v.2} in the proof of
\Cref{cor:special.classes.generate.ring}, the image of $[X_{v}]$ under
the map $f$ can be computed by the recurrence

\begin{equation}\label{eq:Monk.transition}
f([X_{v}]) = x_{r} f([X_{w}]) + \sum_{\substack{ h< r:
\\ \ell(v) =
\ell(w t_{hr})  }} f([X_{wt_{hr}}]), \medskip
\end{equation}
where $(r,s)$ is the lex largest non-inversion for $v$ and
$w=vt_{rs}$.  One can observe by induction that the polynomial
representative for $f([X_{v}])$ determined by the recurrence in
\eqref{eq:Monk.transition} will be a homogeneous nonnegative integer
combination of monomials in the $x_{i}$'s of degree $\coinv (v)$.
From the proof of Monk's formula \Cref{thm:monk} and the details in
his paper, we know that if $x_{i}=[X_{w_{0}s_{i-1}}]-[X_{w_{0}s_{i}}]$
for $1<i<n$,\ $x_{1}=-[X_{w_{0}s_{1}}]$, and
$x_{n}=[X_{w_{0}s_{n-1}}]$, then the corresponding elementary
symmetric polynomial $e_{k}(x_{1},\dots , x_{n})$ vanishes. This still holds
if we negate the variables.  Therefore, if the set $\{f([X_{v}])
\given v \in S_{n}\}$ also spans the coinvariant algebra $R_{n}$, then
$f:A^{*}(\flags) \longrightarrow R_{n}$ is a ring isomorphism.  We
will prove this claim and more in \Cref{cor:Schub.basis}.

The family of polynomials determined by the recurrence in
\eqref{eq:Monk.transition} and the map $f([X_{w_{0}s_{i}})]=x_{1}+\dots
+x_{i}$ have become known as the \emph{Schubert polynomials}, denoted
$\fS_{w}(x_{1},\dots , x_{n})=\fS_{w}$ for $w \in S_{n}$.  They were
not systematically studied until 20 years after Monk's work by Alain
Lascoux and Marcel-Paul Sch\"utzenberger.  They also added the twist
$f([X_{w}])=\fS_{w_{0}w}$, so $\ell (w)=\mathrm{deg}(\fS_{w})$.  We
will systematically study the Schubert polynomials below and relate
them back to solving Schubert problems in the rest of this section.

We conclude this subsection by returning to the cohomology ring of the
flag variety $H^{*}(\flags,\mathbb{Z})$.  Thanks to the work of Borel,
Fulton, Monk and others, we know $H^{*}(\flags,\mathbb{Z})$ is
isomorphic to the Chow ring $A^{*}(\flags)$ and the coinvariant
algebra $R_{n} =\mathbb{Z}[x_{1},\dots , x_{n}]/I_{n}^{+} $. These
rings have linear bases given by Schubert classes
$f([X_{w}])=\fS_{w_{0}w}$ for $w \in S_{n}$. Algebraically, they are
generated by the special classes of the form $[X_{w_{0}s_{i}}]$, or
equivalently $\fS_{s_{i}}$, via Monk's formula and the transition
equations. By \Cref{cor:structure.constants}, their structure constants are
determined by 0-dimensional triple intersections.

\begin{Remark}
From here on out, we will refer to this family of isomorphic rings as
the cohomology ring of the flag variety $H^{*}(\flags,\mathbb{Z})$ in
keeping with the literature, unless we want to emphasize some
particular aspect from the point of view in the Chow ring or the
polynomial quotient.  We will also dispense with the notation for the
map $f$, so we will just say $[X_{w}]=\fS_{w_{0}w}$, meaning that the
polynomial $\fS_{w_{0}w}$ represents the Schubert class $[X_{w}]$
modulo the ideal $I_{n}^{+}$ in $R_{n}$.  
\end{Remark}

\begin{Remark}\label{rem:presentaiton.cohomology}
The cohomology ring of a Schubert variety, or equivalently its Chow
ring, has a nice presentation as a quotient of a polynomial ring as
well.  For $w \in S_{n}$, we have an inclusion $i_w : X_{w}
\hookrightarrow \flags$.  The induced ring homomorphism $i^* :
H^*(\flags) \to H^*(X_w)$ is surjective, and $H^*(X_w)$ is isomorphic
to $R_{n}/I_{w}$ where $I_{w}$ is generated by the Schubert
polynomials $\fS_{v}$ for $v \not\leq w$ in $S_{n}$.  This theorem is
due to Aky\i ld\i z-Lascoux-Pragacz \cite{ALP} and independently to
Carrell \cite[Cor. 4.4]{carrell92}.  See also Reiner-Woo-Yong's
results on finding an efficient presentation of these rings and the
problem they pose of finding a minimal generating set for the ideal
$J_{w}$ such that $H^{*}(X_{w}) \simeq R_{n}/J_{w}$ \cite{RWY}.
Their conjecture was proved by St.\ Dizier and Yong in
\cite[Thm. 1.2]{dizier.yong2023presenting} by explicit construction of
a minimal generating set for the ideal generated by all Schubert
polynomials $\fS_{w}$ such that $v \not \leq w$.  See
\Cref{sub:isomorphisms} for more details on the cohomology of Schubert
varieties and detecting isomorphisms between Schubert varieties.  
\end{Remark}

\subsection{Transition Equations and Schubert
Polynomials}\label{sub:TransitionEqns}

The celebrated Schubert polynomials were defined by Lascoux and
Sch\"utzenberger in their 1982 paper ``Polyn\^omes de
Schubert'' \cite{LS1}.  They form a basis for all polynomials in
countably many variables $\mathbb{Z}[x_{1}, x_{2},\dots ]$.
Therefore, the product of two Schubert polynomials expands in the
basis of Schubert polynomials to produce expansion coefficients on par
with the Schubert structure constants.  Remarkably, Lascoux and
Sch\"utzenberger built on work of Bernstein-Gelfand-Gelfand \cite{BGG}
and Demazure \cite{Dem} to prove their polynomials represent Schubert classes in
the coinvariant algebra, so the Schubert structure constants are
determined by Schubert polynomials.

The following definition for Schubert polynomials is not the original
definition due to Lascoux and Sch\"utzenberger, but they did show it
is equivalent shortly after their introduction using Monk's formula
and the ring homomorphism \eqref{eq:Monk.transition} \cite{LS2}. See
also \cite[4.16]{M2}. We present this definition first because it
follows directly from Monk's work from the 1950's and it is more
efficient for computations than the original formula based on divided
differences.  Furthermore, the definition below using the transition equation
 is directly related to the geometry of Schubert problems
since it is derived from Monk's formula leading quickly to the
Inherited Positivity \Cref{cor:Schub.basis}.  This way we bypass the
work of Bernstein-Gelfand-Gelfand showing how divided differences
applied to the Vandermonde determinant modulo $I_{n}^{+}$ can produce
a family of representatives for Schubert classes in the coinvariant algebra.
While their approach has been highly influential over the past half
century, it is not as direct as Monk's approach and the transition
equation.

The definition below is stated in terms of a recurrence based on the
lexicographically (lex) largest inversion $(r,s)$ for $w$ assuming $w
\neq \id$, where as usual an inversion means $r<s$ and $w(r)>w(s)$.
Note that $r$ is the position of the largest descent in $w$, and $s$
is the largest value such that $w(r) > w(s)$.

\begin{samepage} %%def:SchubertPolynomials
\begin{Definition}[\textbf{Transition Equation for Schubert polynomials} \cite{LS2}]
 \label{t:transitionA} For $w \in S_{n}$, define the
\textit{Schubert polynomial} $\fS_{w} \in \mathbb{Z}[x_{1},x_{2},\dots
, x_{n}]$ by the recurrence relation
\begin{equation}\label{e:transA}
\fS_{w}= x_{r}\fS_{v} + \sum_{\substack{
h<r:
\\
\ell(w)=\ell(v t_{hr})
}}
\fS_{v t_{hr}}
\end{equation}
where $(r,s)$ is the lex largest inversion in
$w$, and $v=w t_{rs}$.  The base case of the
recurrence is $\fS_{\id}=1$.
\end{Definition}
\end{samepage}

\begin{Example}\label{ex:1432}
For $w=1432$, the lex largest inversion is between positions $(3,4)$, so
set $r=3, s=4$, and $v=wt_{3,4}=1423$.  By the transition equation, we
have 
$$
\fS_{1432} = x_3 \fS_{1423} + \fS_{2413}.
$$
The lex largest  inversion for $2413$ is $(2,4)$ so
$$
\fS_{2413} =  x_2 \fS_{2314} + \fS_{3214}.
$$
If we continue to use the transition equation, we find
$\fS_{2314}= x_1 x_2$  and $\fS_{3214} = x^2_1 x_2$ so 
\[
\fS_{2413} = x_2 \fS_{2314} + \fS_{3214} =
x_{1}x_{2}^{2} + x_{1}^{2}x_{2}. 
\]
Similarly, we also have $\fS_{1423}= x_2 \fS_{1324}+
\fS_{3124}$ $= x_2( x_1 + x_2) +x_1^2$.  Therefore,
\begin{equation}\label{ex:1}
\fS_{1432} =  x_1^2 x_2 + x_1^2x_3 + x_1x_2^2 + x_1x_2x_3 + x_2^2x_3.
\end{equation}
\end{Example}

Let's analyze the details in the transition equations in terms of the
diagram of the permutation for $w$.  The lex largest inversion of $w$,
called $(r,s)$ above, corresponds with the rightmost cell in $D(w)$ on
the lowest occupied row of the diagram.  It has coordinates
$(r,w_{s})$.  Removing $(r,w_{s})$ from $D(w)$ we obtain the diagram
of $v=wt_{rs}$.  So, $\ell (v)=\ell (w)-1$.

Every permutation $w'=vt_{hr}$ with $h<r$ and $\ell(w)=\ell(v t_{hr})$
appears on the right side of \eqref{e:transA}.  This set may be empty.
This happens whenever $v_{h}>v_{r}$ for all $h<r$.  To find all $h<r$
with $\ell(w)=\ell(v t_{hr})$, consider the pivot dot in position
$(r,v_{r})$ in $D(v)$.  From $(r,v_{r})$ look northwest to find all
pivot dots $(h,v_{h})$ such that the rectangle formed by $(h,v_{h})$ and
$(r,v_{r})$ don't contain any other pivot dots.  In this case,
$v_{h}<v_{r}$ and $\ell (w)=\ell(v)+1=\ell(v t_{hr})$ so $vt_{hr}$
contributes a term to the right side of the transition equation.

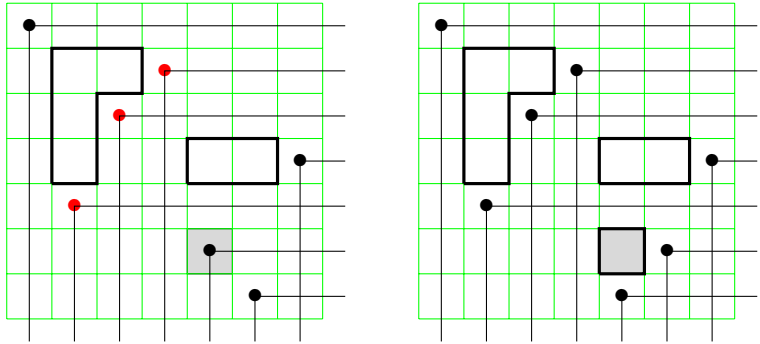
\begin{figure}[h!]
\begin{tikzpicture}[scale=0.6]
\draw[fill=gray!30] (4,1)--(5,1)--(5,2)--(4,2)--(4,1);
\draw[step=1.0,green,thin] (0,0) grid (7,7);
\draw[very thick] (1,3)--(2,3)--(2,5)--(3,5)--(3,6)--(1,6)--(1,3);
\draw[very thick] (4,3)--(6,3)--(6,4)--(4,4)--(4,3);
\node at (0.5,6.5) {$\bullet$};
\node at (1.5,2.5) {$\textcolor{red}{\bullet}$};
\node at (2.5,4.5) {$\textcolor{red}{\bullet}$};
\node at (3.5,5.5) {$\textcolor{red}{\bullet}$};
\node at (4.5,1.5) {$\bullet$};
\node at (5.5,0.5) {$\bullet$};
\node at (6.5,3.5) {$\bullet$};
\draw(0.5,-0.5)--(0.5,6.5)--(7.5,6.5);
\draw(1.5,-0.5)--(1.5,2.5)--(7.5,2.5);
\draw(2.5,-0.5)--(2.5,4.5)--(7.5,4.5);
\draw(3.5,-0.5)--(3.5,5.5)--(7.5,5.5);
\draw(4.5,-0.5)--(4.5,1.5)--(7.5,1.5);
\draw(5.5,-0.5)--(5.5,0.5)--(7.5,0.5);
\draw(6.5,-0.5)--(6.5,3.5)--(7.5,3.5);
\end{tikzpicture}
\qquad
\begin{tikzpicture}[scale=0.6]
\draw[fill=gray!30] (4,1)--(5,1)--(5,2)--(4,2)--(4,1);
\draw[step=1.0,green,thin] (0,0) grid (7,7);
\draw[very thick] (1,3)--(2,3)--(2,5)--(3,5)--(3,6)--(1,6)--(1,3);
\draw[very thick] (4,3)--(6,3)--(6,4)--(4,4)--(4,3);
\draw[very thick] (4,1)--(5,1)--(5,2)--(4,2)--(4,1);
\node at (0.5,6.5) {$\bullet$};
\node at (1.5,2.5) {$\bullet$};
\node at (2.5,4.5) {$\bullet$};
\node at (3.5,5.5) {$\bullet$};
\node at (4.5,0.5) {$\bullet$};
\node at (5.5,1.5) {$\bullet$};
\node at (6.5,3.5) {$\bullet$};
\draw(0.5,-0.5)--(0.5,6.5)--(7.5,6.5);
\draw(1.5,-0.5)--(1.5,2.5)--(7.5,2.5);
\draw(2.5,-0.5)--(2.5,4.5)--(7.5,4.5);
\draw(3.5,-0.5)--(3.5,5.5)--(7.5,5.5);
\draw(4.5,-0.5)--(4.5,0.5)--(7.5,0.5);
\draw(5.5,-0.5)--(5.5,1.5)--(7.5,1.5);
\draw(6.5,-0.5)--(6.5,3.5)--(7.5,3.5);
\end{tikzpicture}
\caption{The diagrams of $v=1437256$  and $w=1437265$ with the
difference highlighted in gray. }
\label{fig:v.w.diagrams}
\end{figure}

\begin{Example}\label{ex:4317265}
For $w=1437265 \in S_{7}$, we have $r=6, s=7$ so $v=wt_{rs}=1437256$.
See \Cref{fig:v.w.diagrams}.  The set of permutations appearing in the
sum of the transition equation correspond with the 3 dots in the
diagram of $v$ that are northwest of the dot in position
$(r,w_{s})=(6,5)$ that form a rectangle with $(6,5)$ not containing
any other dots (these are colored red in the figure).  Therefore, 
\[
\fS_{1437265}=x_{6}\fS_{1437256}+ \fS_{1537246} + \fS_{1457236} +
\fS_{1437526}.  
\]
\end{Example}
\bigskip

The permutations in $S_{n}$ naturally embed into $S_{n+1}$ by adding a
final fixed point.  The inversion set for $w \in S_{n}$, the diagram
of $w$, and the reduced words for $w$ are stable under this embedding.
Identifying two permutations if they have the same inversion set, let
\begin{equation}\label{eq:S.infinity}
S_{\infty}=\bigcup_{n\geq 1} S_{n}
\end{equation}
be the set of all permutations on
$\mathbb{Z}_{+}$ which fix all but a finite number of the positive
integers.

Schubert polynomials are also stable under the natural embedding 
$S_{n}\hookrightarrow S_{n+1}$.  This follows from the
transition equation since if $w \in S_{n}$ then all of the
permutations that appear on the right side of the transition equation
also are in $S_{n}$. So it makes sense to extend the definition of
$\fS_{w}$ to all $w \in S_{\infty}$.

\begin{Definition}\label{defn:lex.largest.inversion.order}
We define the \emph{inversion order} $\lessthan$ on permutations in
$S_{\infty}$ as follows. Given $w \in S_{\infty}$, let $\Inv(w)$ be
the ordered list of inversions in antilex order.  The antilex order is
like reading a dictionary backwards.  Therefore, $\Inv(w)$ begins with the
lex largest inversion of $u$.  Then, for $v \in S_{\infty}$, we say $v
\lessthan w $ provided $\Inv(v ) < \Inv(w) $ in lex order as lists.
For example, $\Inv(1432)=((3,4),(2,4),(2,3))$ and
$\Inv(2413)=((2,4),(2,3),(1,3))$, so $2413\prec 1432$ in inversion
order.
\end{Definition}

\begin{Definition}\label{defn:revlex}
The \textit{revlex} order on monomials is defined so that
$x_{n}<x_{n-1}<\cdots < x_{1}$ and we consider a monomial
$x_{1}^{a_{1}}x_{2}^{a_{2}}\cdots x_{n}^{a_{n}}$ to be a list of
variables with $x_{n}$ listed $a_{n}$ times, followed by $x_{n-1}$
listed $a_{n-1}$ times, etc., and put these lists in antilex order.  So, $x_{4}^{2}x_{3}=x_{4}x_{4}x_{3}$
comes before $x_{4}x_{3}^{2}=x_{4}x_{3}x_{3}$.
\end{Definition}

\begin{Lemma}\label{lem:homogeneous}  Let $w \in S_{\infty}$.
\begin{enumerate}
\item The Schubert polynomial $\fS_{w}$ is a homogeneous polynomial of
degree $\ell (w)$ with nonnegative integer coefficients.
\medskip
\item The leading monomial in revlex order in $\fS_{w}$ is
$x^{c(w)}$ where $c(w)$ is the Lehmer code
defined in \eqref{eq:code}, and all monomials that occur with nonzero
coefficient are in the set $\{x_{1}^{c_{1}}x_{2}^{c_{2}}\cdots
x_{n}^{c_{n}} \given 0\leq c_{i}\leq n-i, \ \forall 1\leq i\leq n \}$
provided $w \in S_{n}$.
\medskip
\item If $w \neq \id$, then $ \fS_{w} \in \mathbb{Z}[x_{1},x_{2},\dots
, x_{r}]$, where $r$ is the position of the last descent in $w$.  
\end{enumerate}
\end{Lemma}

\begin{proof}
The statements hold for the base case $\fS_{\id}=1$.  All of the
permutations on the right hand side of \eqref{e:transA} are strictly
smaller than $w$ in inversion order by construction.  Furthermore, the
permutations on the right hand side of \eqref{e:transA} are in $S_{n}$
provided $w \in S_{n}$. Hence there are only a finite number of terms
in the expansion of a Schubert polynomial. Apply induction over this
finite set to complete the proofs, noting how the transitions work on
the diagram of the permutation of $w$.
\end{proof}

\begin{Exercise}\label{exercise:si.transitions}
For any simple transposition $s_{r} \in S_{\infty}$, use the
Transition Equation to prove the Schubert polynomial 
$\fS_{s_{r}}$ is $x_{1}+x_{2}+\cdots + x_{r}$.
\end{Exercise}

\begin{Exercise}\label{ex:elementary.are.Schubs}
Let $w=s_{i}s_{i+1}\cdots s_{j}$ for $i\leq j$.  Show $\fS_{w}$ is the
elementary symmetric polynomial $e_{k}(x_{1},x_{2},\dots , x_{j})$
where $k=\ell (w) = j-i+1$.  Here $e_{k}(x_{1},x_{2},\dots , x_{j})$
is the sum over all square-free monomials in $x_{1},x_{2},\dots ,
x_{j}$ of degree $k$ with coefficient 1.
\end{Exercise}

\begin{Exercise}\label{ex:homogeneous.are.Schubs}
Recall, $h_{k}(x_{1},x_{2},\dots , x_{n})$ is the sum over all
monomials in $x_{1},x_{2},\dots , x_{n}$ of degree $k$ with
coefficient 1 from \eqref{eq:hom.sym.polys}.  Let
$w=s_{j}s_{j-1}\cdots s_{i}$ for $i\leq j$.  Show $\fS_{w}$ is the
homogeneous symmetric polynomial $h_{k}(x_{1},x_{2},\dots , x_{i})$
where $k=\ell (w) = j-i+1$.
\end{Exercise}

\begin{Exercise}\label{ex:dominant.1}
 Show $\fS_{w_{0}}=x^{n-1}_{1}x_{2}^{n-2}\cdots x_{n-1}^{1}x_{n}^{0}$,
where $w_{0}=[n,n-1,\dots ,2,1]$ is the longest permutation in
$S_{n}$.  More generally, show the following are equivalent for any $w
\in S_{\infty}$.  Such permutations are called \emph{dominant}.
\begin{enumerate}
\item The code of $w$ is a weakly decreasing sequence.
\item The Schubert polynomial has exactly one term, so $\fS_{w}=x^{c(w)}$.
\item The permutation $w$ avoids the pattern $132$.
\end{enumerate}

\end{Exercise}

One special feature of Schubert polynomials as representatives of the Schubert classes $[X_w]$ is that calculations which \emph{a priori} should be done in the quotient ring $R_n \simeq A^*(\flags)$ can actually be done in the polynomial ring $\bZ[x_1, x_2, \ldots]$ (Theorem~\ref{cor:Schub.basis} below). The price to be paid for this pleasant property is the following technical lemma relating the rings $A^*(\flags)$ for different values of $n$. We could also prove this for the rings $R_n$ and $R_{n+1}$ algebraically, but in keeping with our geometric theme we will give a proof by calculating explicit intersections.
\begin{Lemma} \label{lem:schubert-restrict}
    Let $g_n : \bZ[x_1, x_2, \ldots] \to A^*(\Fl(n))$ be the ring homomorphism defined by  
    \begin{equation*}
        g_n(x_i) = \begin{cases}
            [X_{w_0 s_i}] - [X_{w_0 s_{i-1}}] & \text{if $i \leq n$}\\
            0 & \text{otherwise}
        \end{cases},
    \end{equation*}
    where the expression $[X_{w_0 s_i}]$ is to be interpreted as 0 if $i < 1$ or $i > n-1$. Let $w \in S_\infty$. Then $g_n(\fS_{w}) = [X_{w_0 w}]$ if $w \in S_n$ and $g_n(\fS_w) = 0$ otherwise.
\end{Lemma}

Identifying the variable $x_{i}$ with the difference $[X_{w_0 s_i}] -
[X_{w_0 s_{i-1}}]$ is natural in the study of cohomology rings in
terms of line bundles.  The variables represent the first Chern class
of a quotient bundle determined by comparing flags in $X_{w_0 s_i}$
versus $X_{w_0 s_{i-1}}$. We will discuss these ideas somewhat informally in \Cref{sub:DegeneracyLoci.DoubleSchubs}. See also \cite[Sect. 10 and Appendix
B]{Fulton-book} for more details, in particular Lemma 9 in Section
10.3.

\begin{proof}
    The natural embedding $\C^n \hookrightarrow \C^{n+1}$ sending
$(x_{1},\dots , x_{n}) \mapsto (x_{1},\dots , x_{n},0)$ lets us view a
flag in $\C^n$ as a flag in $\C^{n+1}$, so defines an embedding $i :
\Fl(n) \hookrightarrow \Fl(n+1)$. Recall that this induces a ring homomorphism
$i^* : A^*(\Fl(n+1)) \to A^*(\Fl(n))$ satisfying $i^*([Z])
= [i^{-1}(Z \cap i(\Fl(n)))]$ whenever $Z$ and $i(\Fl(n))$ are transverse
. Write $w_0^n$ and $w_0^{n+1}$ to distinguish the
long elements in $S_n$ and $S_{n+1}$, and likewise
$\oppositeE^{\,n}_\bullet \in \Fl(n)$ and $\oppositeE^{\,n+1}_\bullet
\in \Fl(n+1)$.

    We claim that
    \begin{equation} \label{eq:level-reduction}
        X_{w_0^{n+1} w}(\oppositeE^{\,n+1}_\bullet) \cap i(\Fl(n)) = \begin{cases}
            \emptyset & \text{if $w \in S_{n+1} \setminus S_n$}\\
            i(X_{w_0^{n} w}(\oppositeE^{\,n}_\bullet)) & \text{if $w
\in S_n$} \end{cases}, \end{equation} from which we'll get a formula
for $i^*$. Assume $F_\bullet$ is in $X_{w_0^{n+1}
w}(\oppositeE_\bullet) \cap i(\Fl(n))$ from the left side of
\eqref{eq:level-reduction}. First we consider the case where $w \in
S_{n+1} \setminus S_n$. Setting $j = w^{-1}(n+1)$, we have $w_0^{n+1}
w(j) = 1$ and $j < n+1$, so the rank conditions for $X_{w_0^{n+1} w}$
say $\dim(\oppositeE_1 \cap F_j) \geq 1$. But this says $e_{n+1} \in
F_j \subseteq F_n$, while $F_\bullet \in i(\Fl(n))$ requires $F_n =
\C^n$, a contradiction. Therefore $X_{w_0^{n+1} w}(\oppositeE_\bullet)
\cap i(\Fl(n)) = \emptyset$ as claimed.

    Next, suppose $w \in S_n$. Since $F_\bullet \in i(\Fl(n))$, we
have $F_n = \langle e_1, \ldots, e_n \rangle$, so
$\oppositeE^{\,n+1}_{i+1} \cap F_j = \oppositeE^{\,n}_{i} \cap F_j$.
Since $F_\bullet \in X_{w_0^{n+1} w}(\oppositeE_\bullet)$ also, we
know 
\begin{equation}\label{eq:rev.flag.reduction}
\dim(\oppositeE^{\,n}_{i} \cap F_j) = \dim(\oppositeE^{\,n+1}_{i+1}
\cap F_j) \geq \rk(w_0^{n+1} w)[i+1,j]
\end{equation}
for all $i,j \leq n$. But $\rk(w_0^{n+1} w)[i+1,j] = \rk(w_0^{n}
w)[i,j]$ (draw the permutation matrices!).  Therefore, $F_\bullet \in
X_{w_0^{n+1} w}(\oppositeE_\bullet)$ if and only if the truncated flag
$i^{-1}(F_\bullet) = (F_1 \subseteq \cdots \subseteq F_n) \in \Fl(n)$
is in $X_{w_0^n w}(\oppositeE^{\,n}_{\bullet
})$. Thus $X_{w_0^{n+1} w}(\oppositeE^{\,n+1}_\bullet) \cap i(\Fl(n))
= i(X_{w_0^{n+1} w}(\oppositeE^{\,n}_\bullet))$, proving
\eqref{eq:level-reduction}.
    
    The intersection in \eqref{eq:level-reduction} is transverse because $i(\Fl(n)) = X_{w_0^n}(E^{n+1}_{\bullet})$, so we conclude that 
    \begin{equation*}
        i^*[X_{w_0^{n+1} w}] = \begin{cases}
            0 & \text{if $w \in S_{n+1} \setminus S_n$}\\
            [X_{w_0^{n} w}] & \text{if $w \in S_n$.}
        \end{cases}
    \end{equation*}
    In particular, $i^*[X_{w_0^{n+1} s_i}] = [X_{w_0^{n} s_i}]$ if $i
    \leq n-1$ and $i^*[X_{w_0^{n+1} s_n}] = 0$. Therefore $i^* \circ  g_{n+1} = g_n$, and so if $w \in S_{n+1} \setminus S_n$ then $g_n(\fS_w) = (i^* \circ g_{n+1})(\fS_w) = i^*[X_{w_0^{n+1} w}] = 0$. This argument can be iterated to conclude the same for any $w \in S_\infty \setminus S_n$.
 \end{proof}

We are now in a position to relate the Schubert polynomials with the
Schubert structure constants $c_{uv}^{w}$ defined in \Cref{sub:cuvw}
by
\begin{equation} \label{eq:schub-structure-constants}
[X_{u}][X_{v}] = \sum c_{uv}^{w} [X_{w}].  
\end{equation}
Recall each $c_{uv}^{w}$ is a nonnegative integer by the Geometry
Implies Positivity \Cref{cor:structure.constants}.

\begin{Theorem}[\textbf{Inherited Positivity}]  Let $n$ be a positive integer.
\label{cor:Schub.basis} \hfill
\begin{enumerate}[(a)]
\item  The set of all Schubert polynomials $\{\fS_{w} \given w \in
S_{\infty}\}$ forms a basis for $\mathbb{Z}[x_{1},x_{2},\ldots ]$.
\item The set $\{\fS_{w} \given w \in S_{n}\}$ represents a basis for
the coinvariant algebra $R_{n}=\mathbb{Z}[x_{1},\dots ,
x_{n}]/I_{n}^{+} $.
\item The linear map $f:A^{*}(\flags) \longrightarrow R_{n}$ sending
$[X_{w}]$ to $\fS_{w_{0}w}$ is a ring isomorphism.
\item \label{part:inclusion} The natural inclusion $\bZ[x_1, \ldots, x_n] \hookrightarrow \bZ[x_1,
x_2, \ldots]$ induces a ring isomorphism
\[
R_n \simeq \bZ[x_1, x_2, \ldots] / (\fS_w : w \in S_\infty \setminus S_n).
\]

\item For $u,v,w \in S_{n}$, the Schubert structure constant
$c_{uv}^{w}$ is the coefficient of $\fS_{w_{0}w}$ in the product
$\fS_{w_{0}u}\fS_{w_{0}v}$ where $w_{0}$ is the longest permutation in
$S_{n}$. This is true regardless of whether the product is computed in $\bZ[x_1, x_2, \ldots]$ or in the quotient ring $R_n$.
\item Products of Schubert polynomials expand in the basis of Schubert
polynomials with positive integer coefficients.
\end{enumerate}
\end{Theorem}

\begin{proof} \hfill
    \begin{enumerate}[(a)]
        \item

        Let $\mathcal{P}_{\infty} =\mathbb{Z}[x_{1},x_{2}\dots ]$, and
denote the standard monomial basis by $\{x^{\alpha } \}$, where
$x^{\alpha }=x_{1}^{\alpha_{1}}x_{2}^{\alpha_{2}}\cdots$.  Recall the
(Lehmer) code gives a natural bijection between $S_{n}$ and the
product of sets $[n]\times [n-1] \times \cdots \times [2] \times [1]$.
The code of a permutation is invariant under the embedding $S_{n}
\hookrightarrow S_{n+1}$, so the code also gives a bijection between
permutations in $S_{\infty}$ and monomials in
$\mathcal{P}_{\infty}$. By \Cref{lem:homogeneous}, the leading term of
$\fS_{w}$ is $x^{c(w)}$ appearing with coefficient 1, so each Schubert
polynomial is determined by its leading term.  Hence, Schubert
polynomials form a basis of $\mathcal{P}_{\infty}$ with an invertible
triangular change of basis matrix to the standard monomial basis
$\{x^{\alpha } \}$.

        \item By \Cref{ex:artin.monomials}, we know the monomials
        \begin{equation}
        \mathcal{B}=\{x_{1}^{c_{1}}x_{2}^{c_{2}}\cdots x_{n}^{c_{n}} \given 0\leq
        c_{i}\leq n-i \ \forall 1\leq i \leq n \}
        \end{equation}
        form a basis for the coinvariant algebra $R_{n}$. These are the leading terms of the Schubert polynomials $\fS_w$ for $w \in S_n$. Furthermore, all
        monomials appearing with nonzero coefficient in $\fS_{w}$ are in $\mathcal{B}$. Therefore the same argument as in (a) works.

        \item Let $J_n$ be the ideal $(\fS_w : w \in S_\infty
\setminus S_n) \subseteq \mathcal{P}_{\infty}$. By
Lemma~\ref{lem:schubert-restrict}, there is a ring homomorphism $g :
\mathcal{P}_{\infty} / J_n \to A^*(\flags)$ with $g(\fS_w) = [X_{w_0
w}]$ for $w \in S_n$. We know the Schubert classes $[X_{w_0 w}]$ for
$w \in S_{n}$ form a basis for $A^*(\Fl(n))$ and the Schubert
polynomials form a basis for $\mathcal{P}_{\infty}$, so the Schubert
polynomials $\fS_w$ must represent a basis for the quotient
$\mathcal{P}_{\infty} / J_n$ and $g$ is a ring isomorphism.  
        
        We claim that the inclusion $\bZ[x_1, \ldots, x_n] \to
	\mathcal{P}_{\infty}$ descends to an isomorphism $R_n \simeq
	\mathcal{P}_{\infty} / J_n$. Indeed,
	Exercise~\ref{ex:elementary.are.Schubs} says $e_k(x_1, \ldots,
	x_n) = \fS_{s_{n-k+1}s_{n-k+2} \cdots s_n} \in J_n$, so $I_n^+
	\subseteq J_n$. Therefore $R_n \to \mathcal{P}_{\infty} / J_n$  is a
	well-defined ring homomorphism, and sends the basis of Schubert polynomials to the basis of Schubert polynomials. The desired isomorphism $f$ is the inverse of the composition $R_n \to \mathcal{P}_{\infty} / J_n \xrightarrow{g} A^*(\flags)$.

        \item Proved in the proof of (c).

        \item Since Schubert polynomials form a $\bZ$-basis of $\mathcal{P}_{\infty}$, we can write $\fS_{w_0 u} \fS_{w_0 v} = \sum_{z \in S_\infty} a_{uv}^{z} \fS_z$. with $a_{uv}^{z} \in \bZ$. By (d), this expression is equal to $\sum_{z \in S_n} a_{uv}^{z} \fS_z$ in the quotient ring $R_n$. On the other hand, applying the map $f$ from (c) to \eqref{eq:schub-structure-constants} shows $\sum_{z \in S_n} a_{uv}^{z} \fS_z = \sum_{w \in S_n} c_{uv}^{w} \fS_{w_0 w}$, so $a_{uv}^{w_0 w} = c_{uv}^w$ since $\{\fS_w : w \in S_n\}$ is a basis of $R_n$ by (b).
        
        \item This follows from (e) since each $c^w_{uv}$ is a nonnegative
	integer as stated above.
    \end{enumerate}
\end{proof}

Monk's formula (\Cref{thm:monk}) gives rise to a product
formula for Schubert polynomials in the special case of multiplication
by a degree 1 Schubert polynomial.  This is a fundamental result in
the literature on Schubert polynomials.  The expansion depends only on
the covering relations in Bruhat order.  A permutation $w'$
\textit{covers} $w$ in Bruhat order if there is no $v$ with $w' > v >
w$, or equivalently, $w' > w$ and $\ell(w') = \ell(w)+1$. The covering
relations correspond with the edges in the Hasse diagram of Bruhat
order as shown in  in \Cref{Hasse of S4}.

\begin{Theorem}[\textbf{Monk's Formula for Schubert Polynomials}]\label{thm:Schubert.Monk}
For $w \in S_{\infty}$ and a positive integer $r$, 
\begin{equation}\label{eq:Monk-v1}
\fS_{s_r}\fS_{w}=(x_1+\cdots+x_{r})\fS_{w}=\sum_{\substack{k\leq
r<l\\ w t_{k,l}\gtrdot w}}\fS_{w t_{k,l}}.
\end{equation}
Here, $w' \gtrdot w$ means that $w'$ covers $w$ in Bruhat order.  
\end{Theorem}

\medskip

\begin{proof}
Let $n$ be large enough so that all permutations appearing in
\eqref{eq:Monk-v1} have only fixed points after $n$.  By Inherited
Positivity \Cref{cor:Schub.basis}, it suffices to prove
\eqref{eq:Monk-v1} holds in $R_{n}$.  
The polynomials $\fS_{w}$ for $w \in S_{n}$ are defined by the
transition equation in \Cref{t:transitionA}.  By
\Cref{cor:Schub.basis}, $\{\fS_{w} \given w \in S_{n}\}$ forms a basis
for $R_{n}$.  By \Cref{exercise:si.transitions}, we know that
$\fS_{s_{r}}=x_{1}+x_{2}+\cdots + x_{r}$ for all $1\leq r<n$.  Recall
the ring isomorphism $f:A^{*}(\flags) \longrightarrow
R_{n}=\mathbb{Z}[x_{1},\dots , x_{n}]/I_{n}^{+}$ defined by
\eqref{eq:Monk.transition}.  So,
\[
\fS_{s_{r}}=x_{1}+x_{2}+\cdots + x_{r}=f([X_{w_{0}s_{r}}])
\]
and for $\ell(w)>1$, we have $\fS_{w}=f([X_{w_{0}w}])$ by comparing
\eqref{eq:Monk.transition} and \Cref{t:transitionA}.  Therefore, the
Schubert polynomials multiply as in Monk's formula \Cref{thm:monk}.
Therefore, \eqref{eq:Monk-v1} follows by applying the map $w\to
w_{0}w$.   
\end{proof}

\bigskip
\begin{figure}
\begin{center}
\begin{tabular}{|c|l|}
\hline
degree &	Schubert polynomial\\
\hline
0& $\fS_{1 2 3 4} = 1$\\
\hline
1& $\fS_{2 1 3 4} = x_1 $\\
& $\fS_{1 3 2 4} = x_2  + x_1 $\\
& $\fS_{1 2 4 3} = x_3  + x_2  + x_1 $\\
\hline
2& $\fS_{3 1 2 4} = x_1^2 $\\
& $\fS_{2 3 1 4} = x_1 x_2 $\\
& $\fS_{2 1 4 3} = x_1 x_3  + x_1 x_2  + x_1^2 $\\
& $\fS_{1 4 2 3} = x_2^2  + x_1 x_2  + x_1^2 $\\
& $\fS_{1 3 4 2} = x_2 x_3  + x_1 x_3  + x_1 x_2 $\\
\hline
3& $\fS_{4 1 2 3} = x_1^3 $\\
& $\fS_{3 2 1 4} = x_1^2 x_2 $\\
& $\fS_{2 4 1 3} = x_1 x_2^2  + x_1^2 x_2 $\\
& $\fS_{3 1 4 2} = x_1^2 x_3  + x_1^2 x_2 $\\
& $\fS_{1 4 3 2} = x_2^2 x_3  + x_1 x_2 x_3  + x_1^2 x_3  + x_1 x_2^2  + x_1^2 x_2 $\\
& $\fS_{2 3 4 1} = x_1 x_2 x_3 $\\
\hline
4& $\fS_{4 2 1 3} = x_1^3 x_2 $\\
& $\fS_{4 1 3 2} = x_1^3 x_3  + x_1^3 x_2 $\\
& $\fS_{3 4 1 2} = x_1^2 x_2^2 $\\
& $\fS_{3 2 4 1} = x_1^2 x_2 x_3 $\\
& $\fS_{2 4 3 1} = x_1 x_2^2 x_3  + x_1^2 x_2 x_3 $\\
\hline
5& $\fS_{4 3 1 2} = x_1^3 x_2^2 $\\
& $\fS_{4 2 3 1} = x_1^3 x_2 x_3 $\\
& $\fS_{3 4 2 1} = x_1^2 x_2^2 x_3 $\\
\hline
6& $\fS_{4 3 2 1} = x_1^3 x_2^2 x_3 $\\
\hline 
\end{tabular}
\end{center}
\caption{Schubert polynomials for permutations in $S_{4}$.}
\label{fig:Schubs.n=4}
\end{figure}

\begin{Example}\label{example:monk}
If $n=4$, we have $\fS_{1324}\fS_{2134} = \fS_{3124}+\fS_{2314}$ by
Monk's formula.  We also can verify this example using
\Cref{fig:Schubs.n=4} and/or from \Cref{ex:3421.4231}.  Note, unlike
Monk's original formula in \Cref{thm:monk}, we need to consider
permutations to always be in $S_{\infty}$ for products of Schubert
polynomials.  For example, $\fS_{s_{2}} \fS_{41|32} =
\fS_{51|324}+\fS_{42|31}+\fS_{43|12}$.  
\end{Example}

The Schubert polynomials $\{\fS_w\:|\: w\in S_n\}$ can also be
recursively computed by the recurrence
\begin{equation}\label{eq:divided.difference} 
\fS_w =\begin{cases}x_1^{n-1}x_2^{n-2}\cdots x_{n-1} &\text{if }w=[n,
n-1, \ldots, 1] \\ \partial_i\fS_{ws_i}(x)&\text{if }w(i)<w(i+1), \end{cases}
\end{equation}
where the \emph{divided difference operator} $\partial_i$ acts on polynomials by
\[\partial_i f=\frac{f-s_if}{x_i-x_{i+1}}.\]
This is the approach Lascoux and Sch\"utzenberger originally used to
define Schubert polynomials building on the divided difference
operators used by Bernstein-Gelfand-Gelfand \cite{BGG} and Demazure
\cite{Dem}.  See Lemma 9 and the proof of Proposition 4 in
\cite[Section 10.3]{Fulton-book} for the cohomological interpretation
of the divided difference operators in terms of a composition of a
pullback and a pushforward of natural projection operations on flags.
We will prove this recurrence using the nil-Coxeter algebra in
\Cref{sub:Fomin-Stanley} as part of our revisionist history of
Schubert calculus.  See also \Cref{ex:div.diffs.to.transition} for an
outline of the classical proof.  

\begin{Example}\label{ex:div.diffs}
To compute $\fS_{4132}$ via divided differences, we observe $4132=w_{0}s_{3}s_{2}$.  Hence,

\[
\fS_{4132}= \partial_{2}\partial_3 x_{1}^{3}x_{2}^{2}x_{3}.
\]
In steps, first apply $\partial_{3}$, 
\[
\partial_3 x_{1}^{3}x_{2}^{2}x_{3} =
\frac{x_{1}^{3}x_{2}^{2}x_{3}-x_{1}^{3}x_{2}^{2}x_{4}}{x_{3}-x_{4}} = x_{1}^{3}x_{2}^{2},
\]
and then apply  $\partial_{2}$ to the result to get 
\[
\partial_2 x_{1}^{3}x_{2}^{2} =
\frac{x_{1}^{3}x_{2}^{2}-x_{1}^{3}x_{3}^{2}}{x_{2}-x_{3}} = x_{1}^{3}(x_{2}+x_{3}).
\]
Therefore, $\fS_{4132}=x_{1}^{3}x_{2}+x_{1}^{3}x_{3},$   which agrees
with the data in \Cref{fig:Schubs.n=4}.
\end{Example}
\bigskip

\begin{Exercise}
Prove the divided difference operators act linearly on the polynomial
ring of countably many variables
$P_{\infty}=\mathbb{Z}[x_{1},x_{2},\dots ]$.
\end{Exercise}

\begin{Exercise}\label{ex:divided-difference}
Prove that for each positive integer $i$, the kernel and the image of
$\partial_{i}$ is the subalgebra of polynomials symmetric in $x_{i}$
and $x_{i+1}$.  In particular, $\fS_{w}$ is symmetric in $x_{i}$ and
$x_{i+1}$ if and only if $w_{i}<w_{i+1}$.  Furthermore, the following
relations on the divided difference operators hold:
\begin{itemize}
\item $\partial_i^2=0$,
\item $\partial_i\partial_j=\partial_j\partial_i$ if $|i-j|\geq2$,
\item
$\partial_i\partial_{i+1}\partial_i=\partial_{i+1}\partial_i\partial_{i+1}$
for $1\leq i\leq n-2$.
\item Leibniz rule:  For polynomials $f,g$, we have
$\partial_i(fg)=(\partial_i f)g + (s_{i}f) (\partial_i g)$.
\end{itemize}
\end{Exercise}

\begin{Exercise} \label{item:partial.w} Given a reduced word
$(r_{1},r_{2},\ldots , r_{p})$ for $w \in S_{\infty}$, define
$\partial_{w}:P_{\infty} \longrightarrow P_{\infty}$ by the
composition of divided difference operators
$\partial_{w}=\partial_{r_{1}}\partial_{r_{2}}\cdots
\partial_{r_{p}}$.  Prove $\partial_{w}$ is well-defined as an
operator on $P_{\infty}$, and hence $\fS_{w}$ is well-defined by
\eqref{eq:divided.difference}.
\end{Exercise}

A more challenging exercise is the following.
\begin{Exercise}\cite[Prop. 3.1]{Knutson.2012}\label{ex:descent.cycling}
Say $u,v,w \in S_{n}$ such that $\coinv (w)=\coinv (u)+\coinv (v)$.
If there exists an $i$ such that $us_{i}>u$ and $ws_{i}<w$, then use
properties of divided differences to prove the following
\textit{descent cycling symmetries} for Schubert structure constants.
\begin{enumerate}
\item If $vs_{i}>v$, then $c_{uv}^{w}=0$.
\item If $vs_{i}<v$, then $c_{uv}^{w}=c_{u,vs_{i}}^{ws_{i}}$.
\end{enumerate}
\end{Exercise}

Some properties of Schubert polynomials are more readily apparent via
the transition equation (such as monomial positivity) and some via the
divided differences (such as symmetry properties and descent cycling).
We will do a deep dive into the rich combinatorics of Schubert
polynomials in \Cref{sec:CombinatoricsofSchubertPolys}.  In
particular, we explore some symmetry properties of Schubert
polynomials in \Cref{sub:Grassmannians} and identify the set of Schur
polynomials from symmetric function theory and the representation
theory of $S_{n}$ and $GL_{n}$ as particular Schubert polynomials.
For example, we saw in \Cref{ex:elementary.are.Schubs} and
\Cref{ex:homogeneous.are.Schubs} that the elementary symmetric
polynomials and the homogeneous symmetric polynomials in any finite
alphabet are Schubert polynomials.  These are special cases of Schur
polynomials from $S_{n}$ representation theory \cite{Fulton-book,Sag}. It is well known that the Schur polynomials indexed by
partitions represent the Schubert classes in the cohomology ring of
the Grassmannian manifold $\Gr(k,n)$.  See the Math Monthly article by
Kleiman-Laksov from 1972 for a history up to that point \cite{KL}.
The structure constants for Schur polynomials give rise to the
Littlewood-Richardson coefficients. Therefore, an important philosophy
guiding the research in Schubert polynomials over the past 40 years is
as follows.

\bigskip
\begin{Remark}
\textbf{Every tool used to study Schur polynomials has a potential
analog for Schubert polynomials.}
\end{Remark}

We conclude this subsection with one beautiful example of this
philosophy: a generalization of the Pieri formula to Schubert
polynomials due to Frank Sottile \cite{sottile}.  His proof has a
similar feel to the proof of Monk's formula outlined above.  See
\cite[Sect. 5]{sottile} for the ``Geometry of Intersections'' in this
case.  His paper also outlines the connections to prior work on Pieri
formulas.

Let $b,d$ be positive integers.  Let
\begin{equation}\label{eq:one.row.perm}
r[b,d]=[1,2,\ldots,b-1,b+d,b,b+1,\ldots]=s_{b+d-1}\cdots s_{b+1} s_{b}
\end{equation}
where the indices on the simple transpositions are consecutive and
decreasing.  Note that the permutation diagram $D(r[b,d])$ is entirely
contained in row $b$, so we call $r[b,d]$ a \textit{one-row
permutation}.  \Cref{ex:homogeneous.are.Schubs} implies the Schubert
polynomial ${\fS}_{r[b,d]}=h_d(x_1, x_2, \ldots, x_b)$, the
homogeneous symmetric polynomial of degree $d$.  The permutation
$r[b,d]$ is of the special type called Grassmannian from
\Cref{ex:grassmannian.def} since it has has one descent. These will
be important in \Cref{sub:StanleySymmetrics}.

Similarly, define the \textit{one-column permutation}
\begin{equation}\label{eq:one.col}
c[b,d]=[1,2,\ldots ,b-d,b-d+2,\ldots ,b+1,b-d+1,b+2,b+3,\ldots
]=s_{b-d+1}\cdots s_{b-1}s_{b}
\end{equation}
for any positive integers $b,d$ such that $b\geq d$.  Here the indices on the simple
transposition are consecutive and increasing.  The diagram of the
permutation $D(c[b,d])$ is contained in column $b-d+1$.  The permutation $c[b,d]$
is also a Grassmannian permutation, and the Schubert polynomial
${\fS}_{c[b,d]}=e_d(x_1, x_2, \ldots, x_b)$, the elementary symmetric
polynomial of degree $d$ as seen in \Cref{ex:elementary.are.Schubs}.

\begin{Example}\label{ex:one.row.col}
From the diagram in \Cref{ex.2341.table}, one can observe that $2341$
is a one-column permutation, namely $c[3,3]$.  The permutation
$13452=s_{2}s_{3}s_{4}=c[4,3]$, and
$(13452)^{-1}=15234=s_{4}s_{3}s_{2}=r[2,3]$.
\end{Example}

\begin{Theorem}[\textbf{Pieri Formula
for Schubert Polynomials}]\label{thm:sottile}\cite{sottile}
Given $v\in S_{\infty}$ and any one-column permutation $c[b,d]$,

\begin{equation}\label{eq:sottile.elems}
\fS_{c[b,d]}\fS_{v}=e_{d}(x_{1},\dots , x_{b})\ \fS_{v}=\sum \fS_{w}
\end{equation}
\medskip where the sum runs over all distinct $w=vt_{k_1l_1}t_{k_2l_2}
\cdots t_{k_dl_d}$ with $k_{1},k_{2},\dots , k_{d}$ distinct, each
$k_i\leq b<l_i$, 
and $\ell(vt_{k_1l_1} \cdots t_{k_il_i})=\ell (v)+i$ for $1\leq i\leq
d$.  \bigskip

\noindent 
Similarly, given any one-row permutation $r[b,d]$,

\begin{equation}\label{eq:sottile.hom}
\fS_{r[b,d]}\fS_{v}= h_{d}(x_{1},\dots , x_{b}) \ \fS_{v} = \sum
{\fS}_{w} 
\end{equation}
\medskip where the sum indexing is the same as above except that
$l_{1},l_{2},\dots , l_{d}$ are distinct this time.
\end{Theorem}

\begin{Example}\label{ex:pieri}
We saw $\fS_{s_{2}} \fS_{41|32} =
\fS_{51|324}+\fS_{42|31}+\fS_{43|12}$ in \Cref{example:monk}.  The
permutation $s_{2}=132=r[2,1]=c[2,1]$.  Next consider the one-row
permutation $r[2,2]=s_{3}s_{2}=1423$.  By the Pieri formula for
Schubert polynomials, we obtain the expansion of $\fS_{r[2,2]}
\fS_{41|32}$ by starting with $\fS_{s_{2}} \fS_{41|32} =
\fS_{51|324}+\fS_{42|31}+\fS_{43|12}$ and applying all transpositions
$t_{kl}$ such that $k\leq 2<l$ to the right side while ensuring the
position $l$  was not moved in the first application of
Monk's formula.  By Monk's formula,
\[
\fS_{s_{2}} \fS_{51|324} = \fS_{61|3245}+\fS_{53|124}+\fS_{52|314}, 
\]
\[
\fS_{s_{2}} \fS_{42|31} = \fS_{52|314}+ \fS_{43|21},
\]
and 
\[
\fS_{s_2} \fS_{43|12} = \fS_{53|124} + \fS_{45|123}.
\]
Each term on the right side in the three equations above will
contribute to the expansion of $\fS_{r[2,2]} \fS_{41|32}$ since each
transposed position to the right of the line is distinct. Note,
$\fS_{52|314}$ appears in two of the products above, but we only need
to know the set of terms appearing in the expansion above.  Therefore,
% \begin{equation}\label{eq:pieri.example}
% \fS_{c[2,2]}\ \fS_{41|32} = \fS_{53|124}+  \fS_{52|314}.
% \end{equation}
% %(<1.PERM:(5 3 1 2 4)> <1.PERM:(5 2 3 1 4)>)
\begin{equation}\label{eq:pieri.example}
\fS_{r[2,2]}\ \fS_{41|32} = \fS_{61|3245}+\fS_{53|124}+\fS_{52|314}
+ \fS_{45|123} +\fS_{43|21}.
\end{equation}
\end{Example}

\begin{Exercise}\label{ex:div.diffs.to.transition}
Recall, Lascoux and Sch\"utzenberger defined Schubert polynomials for
$w \in S_{n}$ by the formula
\[
\fS_{w} = \begin{cases}
x_{1}^{n-1}x_{2}^{n-2}\cdots x_{n}^{0}& w=w_{0}  \\
\partial_{i} \fS_{ws_{i}} &  w_{i}< w_{i+1}.
\end{cases}
\]
The following steps outline a proof that shows these polynomials
satisfy the Transition Equation \eqref{e:transA}.  Prove each step
using the following notation. Let
$P_{\infty}=\mathbb{Z}[x_{1},x_{2},\dots ]$ and
$P_{n}=\mathbb{Z}[x_{1},x_{2},\dots,x_{n} ]$.  Let $A_{n}$ be the
$\mathbb{Z}$-span of the \emph{Artin monomials} $\{x^{\alpha} \given
0\leq \alpha_{i}\leq n-i \}$ as a finite dimensional subspace of
$P_{n}\subset P_{\infty}$.  These steps follow the outline in
\cite{M2}.
\begin{enumerate}

 \item \label{item:partial.vw}
 Prove
 \[
 \partial_{v}\partial_{w} =
 \begin{cases}
 \partial_{vw}&  \text{if }\ell(v)+\ell(w)=\ell(vw) \\
 0 &   \text{otherwise}.
 \end{cases}
 \]
 Hence, $\partial_{v}\fS_{w} = \fS_{wv^{-1}}$ if $\ell(wv^{-1}) =
 \ell(w)-\ell(v)$ or 0 otherwise.

\item \label{item:expansion.coef.constantterm.partialw}
If $f=\sum c_{w} \fS_{w} \in P_{\infty}$, prove $c_{w}$ is the
constant term of $\partial_{w}(f)$.

\item \label{item:An.preserving} For $1\leq i<n$, show $\partial_{i}$
acts on $A_{n}$ and $P_{n}$.

\item \label{item:dominant.case}
If $w \in S_{n}$ is dominant of shape $\lambda$, then
$\fS_{w}=x^{\lambda}$ including the cases $\fS_{\id}=1$ and
$\fS_{(n-1,n-2,\ldots,1,n)}=x_{1}^{n-2}x_{2}^{n-3}\cdots x_{n-1}^{0}$.
Hence, $\fS_{w}$ is stable as $S_{n} \hookrightarrow
S_{n+1}\hookrightarrow \cdots.$

\item \label{item:Schub.basis}
Prove $\{\fS_{w} \given w \in S_{n} \}$ is a basis for $A_{n}$,\
$\{\fS_{w} \given w \in S_{\infty} \}$ is a basis for $P_{\infty}$,\
and $\{\fS_{w} \given w_{n+1}<w_{n+2}< \cdots \}$ is a basis for
$P_{n}$.  

\item \label{item:linear.product}
Let $f=\sum c_{i} x_{i} \in P_{\infty}$ be a linear polynomial, and let $w \in S_{\infty}$. Then
the product of $f$ and $\fS_{w}$ expands in the Schubert basis as 
\[
f \fS_{w} = \sum_{\ell(wt_{ij}) = \ell(w)+1} (c_{i} - c_{j}) \fS_{wt_{ij}}.
\]
\item Prove
\[
x_{r} \fS_{w} = \sum_{r<j: \ell(wt_{rj}) = \ell(w)+1}  \fS_{wt_{rj}} - \sum_{i<r: \ell(wt_{ir}) = \ell(w)+1}  \fS_{wt_{ir}}.
\]

\item  For $\id \neq w\in S_{\infty}$, let $(r,s)$ be its lex largest
inversion.  Let $v=wt_{rs}$.  Prove
\[
x_{r}\fS_{v} = \fS_{w} - \sum_{i<r:\ell(vt_{ir})=\ell(w)}
\fS_{vt_{ir}}, 
\]
and use this identity to complete the proof that the Schubert
polynomials as defined by divided difference operators satisfy
$\fS_{id}=1$ and the
Transition Equation for $w \neq \id$. 
\end{enumerate}
\end{Exercise}

\bigskip

\bigskip
\subsection{Solving Schubert Problems in 2000 (Reprise)}\label{sub:SchubertProblems2000.reprise}

Going back to the beginning of the story, Hilbert asked mathematicians
to put Schubert's calculus on a rigorous foundation.  In the context
of flag varieties, partial flag varieties and Grassmannians, this task
is complete and rigorous.

Recall from \Cref{sub:SchubertProblems1900} that Hermann C\"asar
Hannibal Schubert (1848--1911) was interested in counting
configurations of subspaces in given arrangements.  The classic
``Schubert Problem'' is:

\medskip
\begin{quote}
How many lines in $\mathbb{R}^{3}$ intersect 4 given lines?
\end{quote}
\medskip

\noindent The possible answers are 0, 1, 2, or infinity.  In the
generic situation, that answer is always 2.  We have phrased
Schubert's question over $\mathbb{R}$ for readers to visualize this
problem. In the following discussions, we will, however, stick with
$\mathbb{C}$, as before in our discussion of the Chow ring and
cohomology. It turns out that in this particular question, the generic
answer is the same in both $\mathbb{R}$ and $\mathbb{C}$ \cite[Thm
C]{sottile96}. Real Schubert calculus is another massive field of its
own and we will not discuss it in this chapter except briefly in \Cref{sub:fields}.  See
\cite{Sottile.book.2011} for details.  

Let's calculate the answer 2 to the classic Schubert problem with the
modern tools of Schubert calculus.  We have been drawing pictures of
flags projectively.  A flag in $\Fl(4)$ is
represented projectively as a point, on a line, in a plane, in a
shoebox.  The shoebox represents an affine hyperplane $V$ in
$\mathbb{C}^{4}$. If one specifies 4 lines $L_{1},L_{2},L_{3},L_{4}$
in $V$, let $\tilde{L}_{1},\tilde{L}_{2},\tilde{L}_{3},\tilde{L}_{4}$
represent the corresponding planes in $\mathbb{C}^{4}$.  Consider the
set of all flags $F_{\bullet}$ such that $\mathrm{dim} (F_{2}\cap
\tilde{L}_{i})\geq 1$ for all $i$.  Such flags would be drawn
projectively with $F_{2}$ represented as a line $L \subset  V$ that intersects
all 4 given lines in the shoebox if such an $F_{\bullet}$ exists.  How
can we describe these flags as a set or as a subvariety of $\Fl(4)$?

For a fixed affine line $L$ in $V$, choose a flag $H_{\bullet} \in
\Fl(4)$ with $H_{2}=\tilde{L}$, the span of the vectors intersecting
$L$.  The set of all flags $F_{\bullet}$ such that $F_{2}$ ``meets''
the given plane $\tilde{L}=H_{2}$, is determined by the single binding
intersection condition $\mathrm{dim}(H_{2} \cap F_{2})\geq 1$, along
with the necessary conditions $\mathrm{dim}(H_{i} \cap F_{j})\geq
i+j-n$ for any $i,j \in [4]$.  The required rank conditions are
bounded below by
\[
\rk(4231)=
\tableau{0 & 0 & 0 & 1 \\
         0 & 1 & 1 & 2 \\
         0 & 1 & 2 & 3 \\
         1 & 2 & 3 & 4 }
\]
and these bounds are tight, so the set of all flags $F_{\bullet}$ such
that $F_{2}$ ``meets'' the given plane $\tilde{L}=H_{2}$ is exactly the
Schubert variety
\[
X_{4231}(H_{\ci}) = \{F_{\ci} \in \Fl(4)
\given \mathrm{dim}(H_{i} \cap F_{j}) \geq \mathrm{rk}(4231)[i,j] \text{ for all } 1\leq i,j \leq 4 \}.
\]

Similarly, the set of all flags $F_{\bullet}$ such that $F_{2}$ equals
the fixed plane $\tilde{L}=H_{2}$, is determined by the single binding
intersection condition $\mathrm{dim}(H_{2} \cap F_{2})\geq 2$, along
with $\mathrm{dim}(H_{i} \cap F_{j})\geq i+j-n$ for all $i,j \in [4]$.
Observe, $\mathrm{dim}(H_{2} \cap F_{2})\geq 2$ forces
$\mathrm{dim}(H_{1} \cap F_{2})\geq 1$ and $\mathrm{dim}(H_{2} \cap
F_{1})\geq 1$ by the nesting conditions on the subspaces in flags.
These rank conditions are tightly bounded below by
\[
\rk(2143)=\tableau{0 & 1 & 1 & 1 \\
         1 & 2 & 2 & 2 \\
         1 & 2 & 2 & 3 \\
         1 & 2 & 3 & 4 }
\]
so the set of all flags $F_{\bullet}$ such that $F_{2}=\tilde{L}$ is precisely
the Schubert variety $X_{2143}(H_{\bullet})$.

Putting these computations together with the classic Schubert problem,
let $G_{\bullet}^{(1)}, G_{\bullet}^{(2)},$
$G_{\bullet}^{(3)}$, $G_{\bullet}^{(4)}$ be 4 flags chosen so that
$G_{2}^{(i)}=L_{i}$ for each $i \in [4]$.  Then

\[
F_{\bullet}\in  X_{4231}(G_{\bullet}^{(1)})\cap
X_{4231}(G_{\bullet}^{(2)}) \cap X_{4231}(G_{\bullet}^{(3)}) \cap
X_{4231}(G_{\bullet}^{(4)})
\]

\bigskip 

\noindent if and only if $\mathrm{dim}(F_{2}\cap L_{i})\geq 1$ for
each $i$, that is the line for $F_{2}$ will be drawn intersecting
$L_{i}$ projectively.  We don't want to overcount flags since there
are an infinite number of flags with $F_{2}$ as their second
component.  Collecting them up according to their second component is
the same as determining how many irreducible components of
$X_{4231}(G_{\bullet}^{(1)})\cap X_{4231}(G_{\bullet}^{(2)}) \cap
X_{4231}(G_{\bullet}^{(3)}) \cap X_{4231}(G_{\bullet}^{(4)})$ are
rationally equivalent to $X_{2143}(H_{\bullet})$.  This is easy to do
via Monk's formula for Schubert classes in the Chow ring, and the
answer does not depend on any choice of flags made above.  Since the
lines $L_{1},L_{2},L_{3}, L_{4}$ are assumed to be generically chosen,
we can assume $G_{\bullet}^{(1)},G_{\bullet}^{(2)},
G_{\bullet}^{(3)},G_{\bullet}^{(4)}$ are 4 generic flags.  Then,

\[
[X_{4231}(G_{\bullet}^{(1)})\cap X_{4231}(G_{\bullet}^{(2)}) \cap
X_{4231}(G_{\bullet}^{(3)}) \cap X_{4231}(G_{\bullet}^{(4)})] =
[X_{4231}]^{4}.
\]

\bigskip

\noindent Applying Monk's formula for multiplication by the special
class $[X_{4231}]$ three times, we have
\begin{align*}
    [X_{4231}]^{4} &= [X_{4231}]^{2} ([X_{3241}] + [X_{4132}])\\
    &= [X_{4231}] ([X_{3142}] + [X_{3142}])\\
    &= ([X_{2143}] + [X_{2143}])\\
    &= 2[X_{2143}].
\end{align*}
Therefore, 4 generic 2-dimensional subspaces $\tilde{L}_{1},
\tilde{L}_{2},\tilde{L}_{3},\tilde{L}_{4}$ in $\mathbb{C}^{4}$ will
all intersect nontrivially with exactly 2 additional 2-dimensional
subspaces.  Hence, we have used the power of the cohomology ring of
the flag variety to completely rigorously compute the generic answer
``2'' to Schubert's classic problem of ``How many lines meet 4 given
lines in 3-dimensional space?''.

To solve the same classic Schubert problem using Schubert polynomials,
we just need to compute
$[X_{4231}]^{4}=\fS_{1324}^{4}=(x_{1}+x_{2})^{4}$ and expand this
polynomial in the basis of Schubert polynomials by
Inherited Positivity \Cref{cor:Schub.basis}.  The coefficient of $[X_{2143}]=\fS_{3412}$
should be 2 in the expansion using the calculations above.  Indeed,
\begin{equation}\label{ex:1324.4}
(\fS_{1324})^{4}=(x_{1}+x_{2})^{4} =
2\fS_{3412}+3\fS_{25134}+ \fS_{162345}.
\end{equation}

\begin{Exercise}\label{ex:3components}
What is the geometric interpretation of the coefficient 3 in
\eqref{ex:1324.4}?
\end{Exercise}

More generally, Schubert was interested in testing when a family of
linear spaces in certain relative positions was intersected by only a
finite number of other linear spaces and determining the number of
generic solutions. These types of relative positions give rise to
problems about intersecting Schubert varieties with respect to certain
flags.  Each of these intersection problems can in principle be done by
expanding the corresponding product of Schubert polynomials and
expanding in the Schubert basis to find the specific
coefficients in the product that determine the generic multiplicity
for the given Schubert problem.  Thus, the generic solutions can be
found by polynomial arithmetic and some linear algebra.  The theory is
completely rigorous.  However, the problem then becomes one of
computational complexity.

Finding the most efficient possible algorithm to compute Schubert
structure constants remains an open problem. Both time and memory
become an issue for large $n$, meaning $n\geq 10$ on a typical
home/office computer of today.  This is closely related to
combinatorial rules for computing Littlewood-Richardson coefficients,
which are the structure constants for Schur functions.  The best-known
combinatorial interpretations of the Littlewood-Richardson
coefficients use the skew semistandard tableaux with reverse lattice
reading word which correspond with paths in Young's lattice, jeu de
taquin, Mondrian tableaux \cite{coskun:LR}, the Remmel-Whitney rule
\cite{RW}, and Vakil's checkers game \cite{Vakil-A}.  The first two
can be found in standard textbooks such as \cite{Fulton-book,Sag}.  We
will describe two more Littlewood-Richardson rules in
\Cref{sub:123.StepFlags} and \Cref{sub:StanleySymmetrics} using
puzzles and leaves in the transition tree.

\begin{Problem}\label{prob:structure.constants}
Find an interpretation for the Schubert structure constants in terms
of counting some sort of combinatorial objects such as paths in Bruhat
order, Mondrian tableaux, labeled diagrams, permutation arrays, or
$n$-dimensional chess games.
\end{Problem}

Bergeron and Sottile refined this problem to looking at paths in what
they call $k$-Bruhat order \cite{BergSott}.  This is a subposet of
Bruhat order on $S_{n}$ that can be naturally obtained from the
conditions in Monk's formula: $w <_{k} wt_{ij}$ if $i\leq k<j$.  See
\cite[Equation (1)]{BergSott}.  The Pieri formula is a perfect example
of the sort of rule one would like for all products of Schubert
classes.

\begin{Remark}\label{rem:count.points}
Note, the Schubert structure constants already count the number of
points in a certain type of generic 0-dimensional intersection by
Geometry Implies Positivity \Cref{cor:structure.constants}.  Perhaps
one could call this a combinatorial interpretation, since they do
count something!  However, it is very difficult to test if flags are
truly in generic position, even though presumably almost anything you
could choose would suffice. Furthermore, solving the equations for a
0-dimensional intersection of varieties in anything but small
dimensions such as $n\leq 6$ is prohibitive.  Therefore, developing
further combinatorial tools to avoid both the genericity problem and
the solving of equations would be useful.  This is why the problem
above is still considered a major open problem in this field.
\end{Remark}

%% From Collen Robichaux: 
\begin{Remark}\label{rem:LRwarning}
A warning about this open problem is in order.  The Schubert structure
constants include the Littlewood-Richardson coefficients for Schur
polynomials as a proper subset \cite{M2}.  In 2006, Narayanan showed
that computing the Littlewood-Richardson coefficients is
$\#P$-complete \cite{Narayanan} for binary inputs.  However, for
arbitrary Schubert problems, one might prefer unary inputs, which is
an alternative notation for the same objects using the unary numeral
system \cite{wiki:unary}. Conjecturally, computing Schubert structure
constants remains $\#P$-hard for unary inputs
\cite{pak2022combinatorialinterpretation}.  Recent work of Pak and
Robichaux has proven that deciding the vanishing of Schubert structure
constants lies in the second level of the polynomial hierarchy,
assuming the Generalized Riemann Hypothesis, that all nontrivial zeros
of $L$-functions have real part $1/2$
\cite{pak.robichaux.2024vanishingschubertcoefficients}.  They provide
an overview of many aspects of the theory of computation related to
the Schubert problem.
\end{Remark}

As mentioned above, there are some ideas for improving on the approach
to computing Schubert structure constants by using Schubert
polynomials and linear algebra, but they all might necessarily include
special cases which are known to be as hard or harder than counting
the number of Hamiltonian paths in a graph, assuming unary inputs.  Is
computing Schubert structure constants on par with traveling salesman
problems \cite{TSP-book} or computing irreducible characters for
$S_{n}$ \cite{Ikenmeyer-Pak-Panova} or finding a perfect strategy for
an $n\times n$ board chess game \cite{EXP-chess}?  Maybe you will help
answer this question.

%% file: section4.tex
\section{Combinatorics of Schubert Polynomials}\label{sec:CombinatoricsofSchubertPolys}

Schubert polynomials are special multivariate polynomials with
nonnegative integer coefficients invented by Lascoux and
Sch\"utzenberger in the early 1980s \cite{LS1}.  They generalize the
Schur polynomials, which play a critical role in $S_{n}$ and $GL_{n}(\C)$
representation theory as well as the cohomology of Grassmannian
varieties.  Schubert polynomials have been widely used and studied
over the past 40 years.  An excellent summary of the early work on
these polynomials appears in Macdonald's notes \cite{M2}; see
Manivel's book \cite{manivel-book} for a more recent treatment.  In
this section, we take an alternative approach to defining Schubert
polynomials using pipe dreams, and then use the nil-Coxeter algebra to
prove they satisfy the original recurrence due to Lascoux and
Sch\"utzenberger using divided difference operators.  The algorithms
described here are reminiscent of games, including chutes/ladders,
Little bumps, and mitosis.  

This section can be read independently from the previous sections.
Just know that the main motivation comes from Schubert calculus.
Specifically, we will show that the generating functions for reduced
pipe dreams equal the Schubert polynomials by bijectively proving that
they satisfy the transition equation of \Cref{t:transitionA}.  The map
itself is quite simple, but there are many details to check.  Once
this first combinatorial interpretation for the monomials in a Schubert
polynomial is established, several others follow easily.

\subsection{Games: Pipe Dreams and Little Bumps}\label{sub:Games}

We will review some basic notation and definitions relating to
permutations, in addition to the material covered in \Cref{sub:perms}.
Recall $S_n$ is the symmetric group of all permutations on
$[n]=\{1,\ldots,n\}$. Multiplication of permutations is defined via
composition, $vw(i):=v(w(i))$.  We write $t_{ij}$ for the
transposition $(i\,j)$ which swaps $i$ and $j$, and we write $s_i = t_{i, i+1}$
$(1 \leq i \leq n-1)$.  The $s_i$ are called \emph{simple
transpositions}; they generate $S_n$ as a Coxeter group.  An
\emph{inversion} of $w \in S_n$ is an ordered pair $(i,j)$ such that
$ i < j $ and $w(i) > w(j)$. The \emph{length} $\ell(w)$ is the number
of inversions of $w$.

Let $w \in S_n$ be a permutation.  A \emph{word for $w$} is a list of
positive integers 
$\mathbf{a}= (a_1, \ldots, a_{k})$ such that 
\[
    s_{a_1} s_{a_2} \ldots s_{a_{k}} = w.
\]
If $k = \ell(w)$, then we say that $\mathbf{a}$ is a \emph{reduced
word} for $w$.  The reduced words are precisely the minimum-length
ways of representing $w$ as a product of simple transpositions.  For
instance, the permutation $[3,2,1] \in S_3$ has two reduced words:
$(1,2,1)$ and $(2,1,2)$.  The empty word $()$ is the unique reduced
word for the identity permutation $[1,2,\ldots, n] \in S_{n}$.  As
with permutations, the \emph{ascent set} of a word $\mathbf{a}=(a_1,
\ldots, a_{k})$ is $\{i\given a_{i}<a_{i+1} \} \subseteq \{1,\ldots,
k-1 \}$.  The \emph{descent set} of $\mathbf{a}$ is the complement in
$[k-1]$.

Write $R(w)$ for the set of all reduced words of the permutation $w$.
The set $R(w)$ has been extensively studied, in part due to interest
in Bott-Samelson varieties and Schubert calculus.  Its size has an
interpretation in terms of counting standard tableaux and the Stanley
symmetric functions~\cite{LS1,little2003combinatorial,Sta84}, which will
be presented in \Cref{sub:StanleySymmetrics}.

\begin{figure}
%\vspace{.2in}
\centering
\raisebox{2mm}{\includegraphics[width=1.7in]{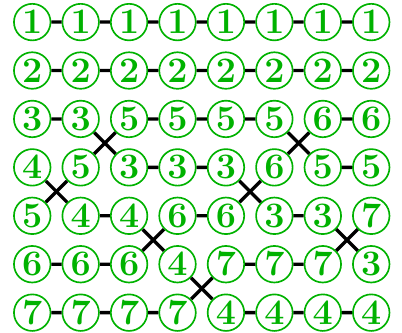}}
\hspace{.2in}
\includegraphics[width=1.7in]{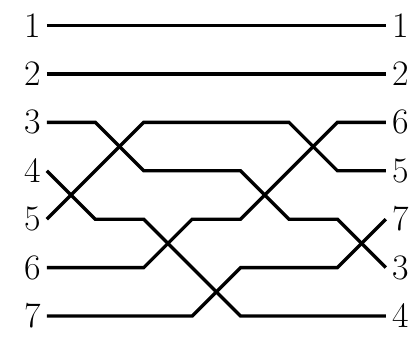}
\hspace{.2in}
\includegraphics[width=1.7in]{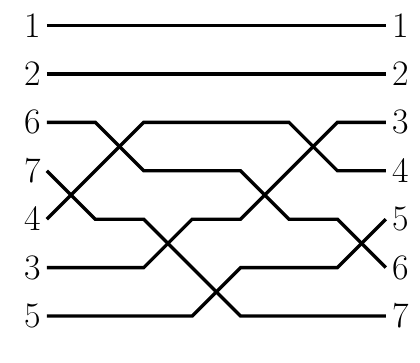}

\caption{The wiring diagram for the reduced word $(4, 3,5, 6, 4, 3,5)
  \in R(1265734)$ notated in three different ways: with the
  intermediate permutations $w^{(t)}$ shown on the left,
   the left-labelling in the middle, and the right-labelling on
   the right. The crossings in columns 2 and 6 are both at row 3.}
\label{fig:wiring.diag}
\end{figure}

We visualize a word as a dynamical system using its wiring
diagram. Wiring diagrams are more rigidly drawn versions of the
string diagrams introduced in \Cref{sub:perms}.  See
Figure~\ref{fig:wiring.diag} for an example.

\begin{Definition}
Let $\mathbf{a}=(a_1,\ldots,a_k)$ be a word for $w$.   For each $0\leq t\leq k$, define the permutation $w^{(t)} \in S_n$ \emph{at time} $t$
by
\[
    w^{(t)} =  s_{a_1}s_{a_2}\cdots s_{a_t},
\]
with the convention that $w^{(0)}$ is the identity, while
$w^{(k)}=w$. The $i$-\emph{wire} (or the wire with \emph{label} $i$)
of $\mathbf{a}$ is defined to be the piecewise linear path joining the
points $(w^{(t)}_{i},t)$ for $0 \leq t \leq k$. We use ``matrix
coordinates", that is, $(w^{(t)}_{i},t)$ refers to row $w^{(t)}(i)$ (numbered from the
top of the diagram) and column $t$ (numbered from the left). The
\emph{wiring diagram} is the union of these $n$ wires.  The diagram is
\textit{left-labeled} because the wires are labeled in order
$1,2,\dots ,n$ down the left side and they retain the label as they
proceed to the right where they appear in order $w_{1},w_{2},\dots ,
w_{n}$ reading down the right side.  The \textit{right-labeled} wiring
diagram for $\mathbf{a}$ is the same union of piecewise linear paths,
but it has wires labeled $1,2,\dots , n$ down the right side of the
wiring diagram, so the labels appear in the order
$w^{-1}(1),w^{-1}(2),\dots , w^{-1}(n)$ reading down the left side.
\end{Definition}

For all $t \geq 1$, observe that between columns $t-1$ and
$t$ in the wiring diagram for $\bfa$, precisely two wires
$i$ and $j$ intersect. This intersection is called a
\emph{crossing}.  One can identify a crossing by its
\emph{column} $t$.  We call $a_t$ the \emph{row} of the
crossing at column $t$.  When the word $\mathbf{a}$ is
reduced, the minimality of the length of $\mathbf{a}$
ensures that any two wires cross at most once.  In this
case, we can also identify a crossing by the unordered pair
of wire labels that are involved, i.e.\ the pair
$\{w^{(t)}(a_t),w^{(t)}(a_{t+1})\}$.

Note that the terms row and column have slightly different
meaning when we refer to a crossing versus a wire.  The
upper left corner of a wiring diagram is at $(1,0)$.  When
we say a crossing is in row $i$ column $j$ it means the
intersection of the crossing is at
$(i+\tfrac12,j-\tfrac12)$.  When we say wire $r$ is in row
$i$ at column $j$, we mean that $w^{(j)}(i)=r$, so that the
$r$-wire passes through the point $(i,j)$.

Observe that for $i<j$, wires $w(i)$ and $w(j)$ cross in the wiring
diagram for $\bfa \in R(w)$ if and only if $w(i)>w(j)$.  This occurs
if and only if $(i,j)$ is an inversion of $w$, which in turn is
equivalent to the wire labels $(w(j),w(i))$ being an inversion of
$w^{-1}$.  Let 
\begin{equation}\label{eq:inversion.def}
\Inv(w)=\{(i,j)\given i<j,w_{i}>w_{j} \} 
\end{equation}
denote the \emph{inversion set} of $w$ indexed by positions.  Note
$\inv (w)=\ell(w)$ is the number of inversions, while $\Inv(w)$ is
the set of pairs indexing the inversions.

Sorry, dear readers! Defining inversion sets as $\Inv(w)=\{(i,j)\given
i<j,w_{i}>w_{j} \}$ is another choice that we just need to live with.
Sometimes one may prefer $\Inv(w^{-1})=\{(w_{j},w_{i})\given
i<j,w_{i}>w_{j} \}$, which one may also call the inversion set of $w$
indexed ``by values''.  The advantage of $\Inv(w^{-1})$ is that the
$i$-wire and the $j$-wires cross in the left-labeled wiring diagram of
$\mathbf{a}$ if and only if $(i,j) \in \Inv(w^{-1})$.  Check that in
\Cref{fig:wiring.diag}!  However, reversing any word for $w$ gives a
word for $w^{-1}$. Thus, if we label the wires $1,2,3,\ldots$ in
increasing order down the right side of a wiring diagram instead of
the left, then the corresponding wires travel right to left, and
appear in order $w^{-1}$ down the left side. Also, the $i$-wire
and the $j$-wire cross in the right-labeled wiring diagram for
$\bfa \in R(w)$ if and only if $(i,j) \in \Inv (w)$ using positions.
Hence, it can be helpful to have both the right-labeled and
left-labeled wiring diagram for $\mathbf{a}$ in mind.  Either way, the
wiring diagram of $\mathbf{a}$ is the same union of unlabeled
piecewise linear paths.

We advise the reader to choose a favorite permutation which is not its
own inverse to keep in mind throughout this chapter and use it to test
your understanding of the notation.  A simple example is $2413 \in
S_{4}$.  It has two reduced words: $(1,3,2)$ and $(3,1,2)$.

\subsubsection{Pipe Dreams}\label{sub:pipes}
Let's consider a variation on the idea of wiring diagrams.  
\begin{Definition}\label{def:pipedream}
A \emph{pipe dream} $D$ is a finite subset of $\mathbb{Z}_+\times
\mathbb{Z}_+$. We usually draw a pipe dream by placing a \+-tile,
called a ``cross", at every $(i,j)\in D$, and a \elbow-tile, called a
``bump" or an ``elbow", at every $(i,j)\in
\mathbb{Z}_+\times\mathbb{Z}_+\setminus D$, creating pipes (or wires)
connecting the left boundary to the top boundary (in matrix
coordinates).

If the wires are numbered $1, 2, 3, \ldots$ across the top, then the
corresponding wires read along the left side of the diagram from
top to bottom form a permutation $w$ of the positive integers that
fixes all but finitely many values.  We call $w$ the \emph{permutation
of $D$} following the literature. Following the terminology for reduced words, we say that $D$ is
\emph{reduced} if $w$ is the permutation of $D$ and $\ell(w) = |D|$.
We write $\rp(w)$ for the set of all reduced pipe dreams for $w$. 
\end{Definition}

\begin{figure}
\begin{center}
\begin{tikzpicture}[scale=0.3]
\begin{scope}[xshift=40em, scale=0.9, thick]
\begin{scope}[transparency group, opacity=0.75]
\draw[-stealth,line width=5pt, orange!50!yellow] (11, 10)-- (-3, 10);
\draw[-stealth,line width=5pt, orange!50!yellow]  (9, 8) -- (-3, 8);
\draw[-stealth,line width=5pt, orange!50!yellow]  (7, 6) -- (-3, 6);
\draw[-stealth,line width=5pt, orange!50!yellow]  (5, 4) -- (-3, 4);
\draw[-stealth,line width=5pt, orange!50!yellow]  (3, 2) -- (-3, 2);
\draw[-stealth,line width=5pt, orange!50!yellow]  (1, 0) -- (-3, 0);
\end{scope}\draw (-0.5,8) arc (270:360:0.5)
(0,8.5) -- (0,11.5);
\draw (-0.5,0) arc (270:360:0.5)
(0,0.5) -- (0,3.5)
(0,3.5) arc (180:90:0.5)
(0.5,4) -- (1.5,4)
(1.5,4) arc (270:360:0.5)
(2,4.5) -- (2,11.5);
\draw (-0.5,10) -- (3.5,10)
(3.5, 10) arc (270:360:0.5)
(4,10.5) -- (4,11.5);
\draw (-0.5,6) arc (270:360:0.5)
(0,6.5) -- (0,7.5)
(0,7.5)  arc (180:90:0.5)
(0.5,8) -- (3.5,8)
(3.5,8) arc (270:360:0.5)
(4,8.5) -- (4,9.5)
(4,9.5) arc (180:90:0.5)
(4.5,10) -- (5.5,10)
(5.5,10) arc (270:360:0.5)
(6,10.5) -- (6,11.5);
\draw (-0.5,2) -- (1.5,2)
(1.5,2) arc (270:360:0.5)
(2,2.5) -- (2,3.5)
(2,3.5) arc (180:90:0.5)
(2.5,4) -- (3.5,4)
(3.5,4) arc (270:360:0.5)
(4,4.5) -- (4,5.5)
(4,5.5) arc (180:90:0.5)
(4.5,6) -- (5.5,6)
(5.5,6) arc (270:360:0.5)
(6,6.5) -- (6,7.5)
(6,7.5) arc (180:90:0.5)
(6.5,8) -- (7.5,8)
(7.5,8) arc (270:360:0.5)
(8,8.5) -- (8,11.5);
\draw (-0.5,4) arc (270:360:0.5)
(0,4.5) -- (0,5.5)
(0,5.5) arc (180:90:0.5)
(0.5,6) -- (3.5,6)
(3.5,6) arc (270:360:0.5)
(4,6.5) -- (4,7.5)
(4,7.5) arc (180:90:0.5)
(4.5,8) -- (5.5,8)
(5.5,8) arc (270:360:0.5)
(6,8.5) -- (6,9.5)
(6,9.5) arc (180:90:0.5)
(6.5,10) -- (9.5,10)
(9.5,10) arc (270:360:0.5)
(10,10.5) -- (10,11.5);
\draw (8.6,10.9) node {\textbf{1}};
\draw (2.6,10.9) node {\textbf{2}};
\draw (0.6,10.9) node {\textbf{3}};
\draw (2.6, 8.9) node {\textbf{4}};
\draw (2.6, 6.9) node {\textbf{5}};
\draw (0.6, 2.9) node {\textbf{6}};
\end{scope}

\begin{scope}[scale=0.9,thick]
\draw (-0.5,8) arc (270:360:0.5)
(0,8.5) -- (0,11.5);
\draw (-0.5,0) arc (270:360:0.5)
(0,0.5) -- (0,3.5)
(0,3.5) arc (180:90:0.5)
(0.5,4) -- (1.5,4)
(1.5,4) arc (270:360:0.5)
(2,4.5) -- (2,11.5);
\draw (-0.5,10) -- (3.5,10)
(3.5, 10) arc (270:360:0.5)
(4,10.5) -- (4,11.5);
\draw (-0.5,6) arc (270:360:0.5)
(0,6.5) -- (0,7.5)
(0,7.5)  arc (180:90:0.5)
(0.5,8) -- (3.5,8)
(3.5,8) arc (270:360:0.5)
(4,8.5) -- (4,9.5)
(4,9.5) arc (180:90:0.5)
(4.5,10) -- (5.5,10)
(5.5,10) arc (270:360:0.5)
(6,10.5) -- (6,11.5);
\draw (-0.5,2) -- (1.5,2)
(1.5,2) arc (270:360:0.5)
(2,2.5) -- (2,3.5)
(2,3.5) arc (180:90:0.5)
(2.5,4) -- (3.5,4)
(3.5,4) arc (270:360:0.5)
(4,4.5) -- (4,5.5)
(4,5.5) arc (180:90:0.5)
(4.5,6) -- (5.5,6)
(5.5,6) arc (270:360:0.5)
(6,6.5) -- (6,7.5)
(6,7.5) arc (180:90:0.5)
(6.5,8) -- (7.5,8)
(7.5,8) arc (270:360:0.5)
(8,8.5) -- (8,11.5);
\draw (-0.5,4) arc (270:360:0.5)
(0,4.5) -- (0,5.5)
(0,5.5) arc (180:90:0.5)
(0.5,6) -- (3.5,6)
(3.5,6) arc (270:360:0.5)
(4,6.5) -- (4,7.5)
(4,7.5) arc (180:90:0.5)
(4.5,8) -- (5.5,8)
(5.5,8) arc (270:360:0.5)
(6,8.5) -- (6,9.5)
(6,9.5) arc (180:90:0.5)
(6.5,10) -- (9.5,10)
(9.5,10) arc (270:360:0.5)
(10,10.5) -- (10,11.5);

\draw (-1,0) node {\small{2}};
\draw (-1,2) node {\small{5}};
\draw (-1,4) node {\small{6}};
\draw (-1,6) node {\small{4}};
\draw (-1,8) node {\small{1}};
\draw (-1,10) node {\small{3}};
\draw (0,12.3) node {\small{1}};
\draw (2,12.3) node {\small{2}};
\draw (4,12.3) node {\small{3}};
\draw (6,12.3) node {\small{4}};
\draw (8,12.3) node {\small{5}};
\draw (10,12.3) node {\small{6}};
% \draw (-1,0) node {\small{6}};
% \draw (-1,2) node {\small{5}};
% \draw (-1,4) node {\small{4}};
% \draw (-1,6) node {\small{3}};
% \draw (-1,8) node {\small{2}};
% \draw (-1,10) node {\small{1}};
% \draw (0,12.3) node {\small{2}};
% \draw (2,12.3) node {\small{6}};
% \draw (4,12.3) node {\small{1}};
% \draw (6,12.3) node {\small{3}};
% \draw (8,12.3) node {\small{5}};
% \draw (10,12.3) node {\small{4}};
\end{scope}

\begin{scope}[xshift=80em, yshift=2ex, scale=1.75, thick]

\node at (-0.3,5) {$1$}; \node at (-0.3,4) {$2$}; \node at
(-0.3,3) {$3$}; \node at (-0.3,2) {$4$}; \node at (-0.3,1)
{$5$}; \node at (-0.3,0) {$6$};

\node at (6.3,5) {$3$}; \node at (6.3,4) {$1$}; \node at
(6.3,3) {$4$}; \node at (6.3,2) {$6$}; \node at (6.3,1)
{$5$}; \node at (6.3,0) {$2$};

\coordinate (0/6) at (0, 0); \coordinate (0/5) at (0, 1);
\coordinate (0/4) at (0, 2); \coordinate (0/3) at (0, 3);
\coordinate (0/2) at (0, 4); \coordinate (0/1) at (0, 5);

\coordinate (1/5) at (1, 0); \coordinate (1/6) at (1, 1);
\coordinate (1/4) at (1, 2); \coordinate (1/3) at (1, 3);
\coordinate (1/2) at (1, 4); \coordinate (1/1) at (1, 5);

\coordinate (2/5) at (2, 0); \coordinate (2/6) at (2, 1);
\coordinate (2/4) at (2, 2); \coordinate (2/2) at (2, 3);
\coordinate (2/3) at (2, 4); \coordinate (2/1) at (2, 5);

\coordinate (3/5) at (3, 0); \coordinate (3/6) at (3, 1);
\coordinate (3/4) at (3, 2); \coordinate (3/2) at (3, 3);
\coordinate (3/1) at (3, 4); \coordinate (3/3) at (3, 5);

\coordinate (4/5) at (4, 0); \coordinate (4/6) at (4, 1);
\coordinate (4/2) at (4, 2); \coordinate (4/4) at (4, 3);
\coordinate (4/1) at (4, 4); \coordinate (4/3) at (4, 5);

\coordinate (5/5) at (5, 0); \coordinate (5/2) at (5, 1);
\coordinate (5/6) at (5, 2); \coordinate (5/4) at (5, 3);
\coordinate (5/1) at (5, 4); \coordinate (5/3) at (5, 5);

\coordinate (6/2) at (6, 0); \coordinate (6/5) at (6, 1);
\coordinate (6/6) at (6, 2); \coordinate (6/4) at (6, 3);
\coordinate (6/1) at (6, 4); \coordinate (6/3) at (6, 5);

\draw (0/1) -- (1/1) -- (2/1) -- (3/1) -- (4/1) -- (5/1) --
(6/1); \draw (0/2) -- (1/2) -- (2/2) -- (3/2) -- (4/2) --
(5/2) -- (6/2); \draw (0/3) -- (1/3) -- (2/3) -- (3/3) --
(4/3) -- (5/3) -- (6/3); \draw (0/4) -- (1/4) -- (2/4) --
(3/4) -- (4/4) -- (5/4) -- (6/4); \draw (0/5) -- (1/5) --
(2/5) -- (3/5) -- (4/5) -- (5/5) -- (6/5); \draw (0/6) --
(1/6) -- (2/6) -- (3/6) -- (4/6) -- (5/6) -- (6/6);

\end{scope}
\end{tikzpicture}
 \caption{Left: a reduced pipe dream $D$ for
   $w=[3,1,4,6,5,2]$.
Middle: the reading
   order for the crossings, with numbers indicating
   position in the order.  The associated reduced word is
   $\mathbf{r}_D =(5,2,1,3,4,5) \in R(w)$.  Reading the sequence of row
   numbers and column numbers in the same order we obtain $\mathbf{i}_D=(1,1,1,2,3,5)$ and
   $\mathbf{j}_D=(5,2,1,2,2,1)$ respectively.  Right: the left-labeled
   wiring diagram of the
   associated reduced word
   $\mathbf{r}_D =(5,2,1,3,4,5) \in R(w)$.
   }
 \label{fig:rcgraphs}
\end{center}
\end{figure}
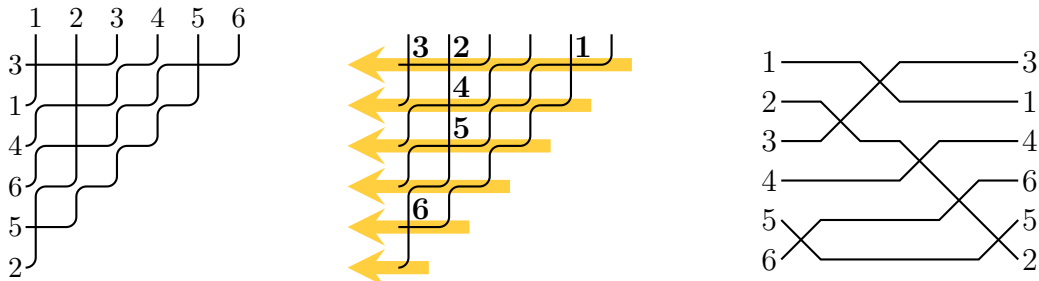

We only need to draw a finite
number of wires in a triangular array to represent a pipe dream since
it necessarily contains only a finite number of crossings.  See
Figure~\ref{fig:rcgraphs} for an example of a pipe dream for
$w=314652$.

We call the elements of a pipe dream $D \subset \mathbb{Z}_+ \times
\mathbb{Z}_+ $ its \emph{crossings} or \emph{occupied positions}, and
the elements $(i,j)$ of $ \mathbb{Z}_+ \times \mathbb{Z}_+ \setminus
D$ the \emph{unoccupied positions}.  Each crossing involves two wires,
which are said to \emph{enter the crossing horizontally} and
\emph{vertically}.

Two
wires labeled $i<j$ cross somewhere in a pipe dream $D \in \rp(w)$ if and only if
$(i,j) \in \Inv(w^{-1})$.  Observe that the larger labeled
wire necessarily enters the crossing horizontally in a reduced pipe
dream as it proceeds down and left.

One can easily read off a reduced word $\mathbf{r}_D$ for $w$ from a
reduced pipe dream $D \in \rp(w)$.  Order the crossings in $D$ in the order given
by reading rows from top to bottom, and from right to left within each
row, as in the middle diagram of Figure~\ref{fig:rcgraphs}.  We call this the \emph{reading order} on $D$.  We construct
three words from the ordered list of crossings: the row numbers of the
crossings $\mathbf{i}_{D} = (i_{1}, i_{2},\ldots, i_{p})$, the column
numbers $\mathbf{j}_{D} = (j_{1}, j_{2},\ldots, j_{p})$ and the
diagonal numbers $\mathbf{r}_{D} =
(i_{1}+j_{1}-1,i_{2}+j_{2}-1,\ldots, i_{p}+j_{p}-1) =
\mathbf{j}_D+\mathbf{i}_D -\mathbf{1}$, which is the corresponding
reduced word of $D$.  Any two of
$\mathbf{i}_D,\mathbf{j}_D,\mathbf{r}_D$ suffice to determine $D$.

For example, Figure~\ref{fig:rcgraphs} shows the pipe dream $D$ with
the reduced word $\mathbf{r}_{D}=(5,2,1,3,4,5)$, row numbers
$\mathbf{i}_{D}=(1,1,1,2,3,5)$, and column numbers
$\mathbf{j}_{D}=(5,2,1,2,2,1)$. Observe that the associated wiring diagram looks
significantly different than the pipe dream.  Can you see how to
reflect, rotate, and stretch the pipe dream to get the crossings to
 appear in the same order?

From the reading order of a reduced pipe dream $D$, one can observe
that $(\mathbf{r}_{D},\mathbf{i}_{D}) = ((r_1, \ldots, r_p), (i_1, \ldots, i_p))$ satisfies three
\textit{compatibility conditions}
\begin{equation}\label{eq:compatible}
\begin{array}{lll}
\bullet & i_{1}\leq  i_{2}\leq \cdots \leq  i_{p}, &\\
\bullet & i_{k}\leq  r_{k}, & \forall \ 1\leq k \leq p,\\
\bullet & \text{if $r_{k}<r_{k+1}$, then } i_{k}<i_{k+1}& \forall \ 1\leq k \leq p.
\end{array}
\end{equation}
For any reduced word $\mathbf{r}$, say the sequence $\mathbf{i}$ is \textit{compatible} with
$\mathbf{r}$ if it satisfies the conditions above.   A reduced word can have many compatible sequences
or it may have none.  For example, the reduced word $(2,3,1) \in R(3142)$ has
none, while $(2,1,3)$ has compatible sequences $(1,1,2)$ and
$(1,1,3)$.

\begin{Definition}\cite{b-j-s}\label{def:compatible-sequence}
For a permutation $w$, a pair of positive integer sequences of the
same length
$(\mathbf{a},\mathbf{i})=\big((a_1,\ldots,a_{\ell}),(i_1,\ldots,i_{\ell})\big)$
is called a \emph{compatible pair} for $w$ if $\mathbf{a}\in R(w)$ is
a reduced word for $w$ and $(\mathbf{a},\mathbf{i})$ satisfies the
compatibility conditions in \eqref{eq:compatible}.  
\end{Definition}

\begin{Exercise}\label{ex:compatible.pair}
Prove that the reduced pipe dreams for $w$ are in bijection with the
compatible pairs for $w$. The forward direction was done in the
observation above.  For the converse statement, one must show that any
biword $((a_{1}, a_{2},\ldots, a_{p}), (i_{1}, i_{2},\ldots, i_{p}))$
such that $\mathbf{a}=(a_{1}, a_{2},\ldots, a_{p})$ is a reduced word
and $(i_{1}, i_{2},\ldots, i_{p})$ is compatible with
$\mathbf{a}$ gives rise to a reduced pipe dream for the permutation
$s_{a_{1}}s_{a_{2}}\cdots s_{a_{p}}$.
\end{Exercise}

The \emph{monomial weight} of a pipe dream $D$ is given by the product
over row numbers of the crossings
\[
 x^D := \prod_{(i,j) \in D} x_i = x_{i_{1}}x_{i_{2}}\cdots x_{i_{p}}
\]
where $x_1, x_2, \ldots$ are formal commuting variables.  Adding a second set of
commuting variables, $y_{1},y_{2},\ldots$ one can also record the column
numbers.  By convention, we associate to $D$ the polynomial
\[
(x-y)^{D}=\prod_{(i,j) \in D} (x_i-y_{j}).  
\]
For example, the pipe dream of \Cref{fig:rcgraphs} has weight
$x^D=x_1^3x_2x_3x_5$
and \[(x-y)^D=(x_1-y_1)(x_1-y_2)(x_1-y_5)(x_2-y_2)(x_3-y_2)(x_{5}-y_{1}).\]

The following important theorem in this theory shows that each
Schubert polynomial is the generating function for reduced pipe dreams
weighted by their monomials as follows.  This formula could serve as
the definition of Schubert polynomials in analogy with the definition
of Schur functions as the sum over semistandard tableaux in symmetric
function theory.  

\begin{Theorem}\label{thm:pipedreams}
The \emph{Schubert polynomial} of $w \in S_n$ is 
\[
    \fS_{w}=\fS_{w}(x_1,x_2,\ldots,x_n) := \sum_{D \in \rp(w)} x^D.
\]
\end{Theorem}

\begin{Corollary}[\textbf{The BJS Formula}]\label{cor:bjs}
The \emph{Schubert polynomial} of $w \in S_n$ is 
\[
    \fS_{w}=\fS_{w}(x_1,x_2,\ldots,x_n) := \sum
    x_{i_{1}}x_{i_{2}}\cdots x_{i_{p}}
\]
summed over all compatible pairs $(\mathbf{a},\mathbf{i}) = ((a_1,\ldots,a_p),(i_1,\ldots,i_p))$ for $w$.  
\end{Corollary}

\begin{Remark}\label{rem:history.Schubs}
This theorem and corollary build on work of many people, so we pause
for a historical note.  Recall, Lascoux-Sch\"utzenberger \cite{LS1}
invented the Schubert polynomials as a solution to the divided
difference equations using the specific choice of representative for
the cohomology class of a point given by
$[X_{\id}]=\fS_{w_{0}}=x_{1}^{n-1}x_{2}^{n-2}\cdots x_{n-1}$ as we saw
in \Cref{sub:TransitionEqns} in the early 1980s.  In the early 1990s,
Stanley conjectured the monomial expansion for the Schubert polynomial
$\fS_{w}$ in terms of biwords consisting of reduced word and
compatible sequence pairs for $w$, \Cref{cor:bjs}.  He had used
compatible pairs to study the number of reduced expressions for
$w_{0}$ in his celebrated paper \cite{Sta84}.  His conjecture was
first proved in 1992 by Billey-Jockusch-Stanley \cite{b-j-s} using a
different recurrence than the transition equation.  Shortly after,
Fomin-Stanley \cite{FS} gave another proof related to the nil-Coxeter
algebra, which we will return to in \Cref{sub:Fomin-Stanley}.
Fomin-Kirillov \cite{FK} gave an interpretation of the compatible
pairs as \emph{pseudo-line} arrangements, which are more free-form
than the wiring diagrams drawn so far.  Billey-Bergeron were the first
to use pipe dreams as we have drawn them above \cite{billey-bergeron},
though they called them \emph{RC-graphs} since they are a
visualization of the reduced word and compatible sequence pairs
$(\mathbf{r}_{D},\mathbf{i}_{D})$.  They used them to get a bijective
proof of Monk's formula using an insertion algorithm, and developed
the chute and ladder moves connecting all reduced pipe dreams for $w$.
Knutson-Miller coined the name ``pipe dreams'', which seems to have
stuck, in their paper \cite{knutson-miller-2005} where they give other
geometric and algebraic interpretations of individual pipe dreams. See
also Kogan's work related to toric varieties and pipe dreams
\cite{kogan.phd}. Recently, Nadeau-Spink-Tewari have surveyed many of
the combinatorial objects related to Schubert polynomials and given an
all positive monomial expansion using new creation operators starting
at $\fS_{\id}=1$ instead of divided difference operators starting at
$\fS_{w_{0}}$ \cite{nadeau.spink.tewari.2024}.  
\end{Remark}

There are many proofs now of \Cref{thm:pipedreams} or its equivalent
formulation in \Cref{cor:bjs}.  Our goal is to present a proof that
could have come shortly after Monk's work in the 1950's.  This
approach is a bijective proof based on the transition equation in
\Cref{t:transitionA}.  The base case says $\fS_{\id}=1$, which
corresponds with the unique pipe dream with no crossings.  To prove
the recurrence holds for the weight generating function over pipe
dreams of $w$, we will use a variation on Little's bumping
algorithm~\cite{little2003combinatorial}, also known as a ``Little
bump''.  It is a map on the set of all possible reduced words.  It was
introduced to study the decomposition of Stanley symmetric functions
into Schur functions in a bijective way.  Later, the Little algorithm
was found to be related to the Robinson-Schensted-Knuth
map~\cite{little2} and the Edelman-Greene map~\cite{hamaker-young}; it
has been extended to signed
permutations~\cite{billey-hamaker-roberts-young}, affine permutations
\cite{lam-shimozono}, and the subset of involutions in $S_{n}$
\cite{Hamaker-Marberg-Pawlowski.2016}.

Before we move on to Little bumps, you might be asking: how does one
find all reduced pipe dreams for $w \in S_{n}$?  Drawing one at random
is a challenge.  Luckily, there are a couple that are straightforward
to find, and these can be used to find all of the others.  Recall the
code of a permutation is $c(w)= (c(w)_1, c(w)_2, \ldots, c(w)_n)$ where $c(w)_i =
|\{j \in [n] \given (i,j) \in \Inv(w) \}|$.  We construct the
\textit{bottom pipe dream} for $w$, denoted $\Dbot(w)$, by placing $c_{i}(w)$ left-justified
cross tiles on row $i$ for each $1\leq i\leq n$.  Since $0\leq
c_{i}(w)\leq n-i$, the pipe dream $\Dbot(w)$ can be represented by a
pipe dream with $n$ wires. The $c_{1}(w)$ \+-tiles on row 1 will lead
the wire labeled $w_{1}=c_{1}(w)+1$ along at the top to turn left and exit on the
first row.  By induction, wires $w_{2},\dots , w_{n}$ will exit on the
left in rows $2,3,\dots , n$, proving $\Dbot(w)\in \rp(w)$.

We get a second pipe dream for $w$ by transposing $\Dbot(w^{-1})$.
This is called the \textit{top pipe dream}, denoted $\Dtop(w)$.  It will have its \+-tiles all top-justified in the columns.

For example, if $w=25143$, then $c(w)=(1,3,0,1,0)$ and
$c(w^{-1})=(2,0,2,1,0)$, so 
$$
\Dbot(w)=\begin{matrix}
+&.&.&.\\
+&+&+\\
.&.\\
+
\end{matrix}
\hspace{.5in}
\Dtop(w)=\begin{matrix}
+&.&+&+\\
+&.&+\\
.&.\\
.
\end{matrix}
$$
by definition.  Note that we may omit any \elbow-tiles, since all the
information of a pipe dream is contained in the \+-tiles.  The reader
should compare these pipe dreams to the diagram $D(w)$ from
\eqref{eq:diagram} since the code of $w$ also counts the number of
elements in $D(w)$ on each row.

\begin{figure}[h!]
\centering
\begin{tikzpicture}[scale=0.4]
\node at (0,0) {$\cdot$};
\node at (1,0) {$\cdot$};
\node at (0,-1) {$+$};
\node at (1,-1) {$+$};
\node at (0,-2) {$+$};
\node at (1,-2) {$+$};
\node at (0,-3) {$+$};
\node at (1,-3) {$\cdot$};
\def\a{5};
\node at (\a+0,0) {$\cdot$};
\node at (\a+1,0) {$+$};
\node at (\a+0,-1) {$+$};
\node at (\a+1,-1) {$+$};
\node at (\a+0,-2) {$+$};
\node at (\a+1,-2) {$+$};
\node at (\a+0,-3) {$\cdot$};
\node at (\a+1,-3) {$\cdot$};
\draw[very thick,->] (2,-1.5)--(4,-1.5);
\end{tikzpicture}
\qquad\qquad
\begin{tikzpicture}[scale=0.2]
\def\a{0};
\def\b{0};
\draw[step=2, color=green](\a+2,-10+\b)grid(\a+6,-2+\b);
\draw[very thick](\a+3,-2+\b)--(\a+3,-2.5+\b);
\draw[very thick](\a+3,-3.5+\b)--(\a+3,-4+\b);
\draw[very thick](\a+2,-3+\b)--(\a+2.5,-3+\b);
\draw[very thick](\a+3.5,-3+\b)--(\a+4,-3+\b);
\draw[very thick](\a+2.5,-3+\b)arc(270:360:0.5);
\draw[very thick](\a+3.5,-3+\b)arc(90:180:0.5);
\draw[very thick](\a+5,-2+\b)--(\a+5,-2.5+\b);
\draw[very thick](\a+5,-3.5+\b)--(\a+5,-4+\b);
\draw[very thick](\a+4,-3+\b)--(\a+4.5,-3+\b);
\draw[very thick](\a+5.5,-3+\b)--(\a+6,-3+\b);
\draw[very thick](\a+4.5,-3+\b)arc(270:360:0.5);
\draw[very thick](\a+5.5,-3+\b)arc(90:180:0.5);
\draw[very thick](\a+2,-5+\b)--(\a+4,-5+\b);
\draw[very thick](\a+3,-4+\b)--(\a+3,-6+\b);
\draw[very thick](\a+4,-5+\b)--(\a+6,-5+\b);
\draw[very thick](\a+5,-4+\b)--(\a+5,-6+\b);
\draw[very thick](\a+2,-7+\b)--(\a+4,-7+\b);
\draw[very thick](\a+3,-6+\b)--(\a+3,-8+\b);
\draw[very thick](\a+4,-7+\b)--(\a+6,-7+\b);
\draw[very thick](\a+5,-6+\b)--(\a+5,-8+\b);
\draw[very thick](\a+2,-9+\b)--(\a+4,-9+\b);
\draw[very thick](\a+3,-8+\b)--(\a+3,-10+\b);
\draw[very thick](\a+5,-8+\b)--(\a+5,-8.5+\b);
\draw[very thick](\a+5,-9.5+\b)--(\a+5,-10+\b);
\draw[very thick](\a+4,-9+\b)--(\a+4.5,-9+\b);
\draw[very thick](\a+5.5,-9+\b)--(\a+6,-9+\b);
\draw[very thick](\a+4.5,-9+\b)arc(270:360:0.5);
\draw[very thick](\a+5.5,-9+\b)arc(90:180:0.5);

\def\a{10};
\def\b{0};
\draw[step=2, color=green](\a+2,-10+\b)grid(\a+6,-2+\b);
\draw[very thick](\a+3,-2+\b)--(\a+3,-2.5+\b);
\draw[very thick](\a+3,-3.5+\b)--(\a+3,-4+\b);
\draw[very thick](\a+2,-3+\b)--(\a+2.5,-3+\b);
\draw[very thick](\a+3.5,-3+\b)--(\a+4,-3+\b);
\draw[very thick](\a+2.5,-3+\b)arc(270:360:0.5);
\draw[very thick](\a+3.5,-3+\b)arc(90:180:0.5);
\draw[very thick](\a+4,-3+\b)--(\a+6,-3+\b);
\draw[very thick](\a+5,-2+\b)--(\a+5,-4+\b);
\draw[very thick](\a+2,-5+\b)--(\a+4,-5+\b);
\draw[very thick](\a+3,-4+\b)--(\a+3,-6+\b);
\draw[very thick](\a+4,-5+\b)--(\a+6,-5+\b);
\draw[very thick](\a+5,-4+\b)--(\a+5,-6+\b);
\draw[very thick](\a+2,-7+\b)--(\a+4,-7+\b);
\draw[very thick](\a+3,-6+\b)--(\a+3,-8+\b);
\draw[very thick](\a+4,-7+\b)--(\a+6,-7+\b);
\draw[very thick](\a+5,-6+\b)--(\a+5,-8+\b);
\draw[very thick](\a+3,-8+\b)--(\a+3,-8.5+\b);
\draw[very thick](\a+3,-9.5+\b)--(\a+3,-10+\b);
\draw[very thick](\a+2,-9+\b)--(\a+2.5,-9+\b);
\draw[very thick](\a+3.5,-9+\b)--(\a+4,-9+\b);
\draw[very thick](\a+2.5,-9+\b)arc(270:360:0.5);
\draw[very thick](\a+3.5,-9+\b)arc(90:180:0.5);
\draw[very thick](\a+5,-8+\b)--(\a+5,-8.5+\b);
\draw[very thick](\a+5,-9.5+\b)--(\a+5,-10+\b);
\draw[very thick](\a+4,-9+\b)--(\a+4.5,-9+\b);
\draw[very thick](\a+5.5,-9+\b)--(\a+6,-9+\b);
\draw[very thick](\a+4.5,-9+\b)arc(270:360:0.5);
\draw[very thick](\a+5.5,-9+\b)arc(90:180:0.5);

\draw[very thick,->] (7,-6)--(11,-6);
\end{tikzpicture}
\caption{Ladder moves on pipe dreams move one cross to the right and
preserve the permutation.}
\label{fig:ladder-move-PD}
\end{figure}
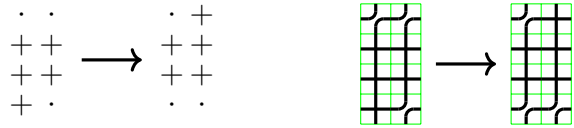

There are two types of moves on pipe dreams that preserve the
permutation: chute moves (\Cref{fig:chute-move-PD}) and ladder moves
(\Cref{fig:ladder-move-PD}).  A chute move can be thought of as
swapping a \+-tile with the first \elbow-tile on the row below and to
the left in such a way as to preserve the permutation.  Formally,
$(i,j+k+1)\in D$ moves to $(i+1,j)\notin D$ if $(i+a,j+b)\in D$ for
all $a\in\{0,1\}$, $b\in\{1,2,\ldots,k\}$ and $(i,j),(i+1,j+k+1)\notin
D$, for some $k\geq0$. Similarly, a ladder move swaps a \+-tile with
the first \elbow-tile on the column to the right and above in such a
way as to preserve the permutation.  Formally, $(i+k+1,j)\in D$ moves
to $(i,j+1)\notin D$ if $(i+a,j+b)\in D$ for all
$a\in\{1,2,\ldots,k\}$, $b\in\{0,1\}$ and $(i,j),(i+k+1,j+1)\notin D$
for some $k\geq0$. It is easy to see from the wires in each case that
the chute and ladder moves on pipe dreams preserve the underlying
permutation.
\begin{figure}[h!]
\centering
\begin{tikzpicture}[scale=0.4]
\node at (0,0) {$\cdot$};
\node at (1,0) {$+$};
\node at (2,0) {$+$};
\node at (3,0) {$+$};
\node at (4,0) {$\cdot$};
\node at (0,-1) {$+$};
\node at (1,-1) {$+$};
\node at (2,-1) {$+$};
\node at (3,-1) {$+$};
\node at (4,-1) {$\cdot$};
\def\a{-8};
\node at (\a+0,0) {$\cdot$};
\node at (\a+1,0) {$+$};
\node at (\a+2,0) {$+$};
\node at (\a+3,0) {$+$};
\node at (\a+4,0) {$+$};
\node at (\a+0,-1) {$\cdot$};
\node at (\a+1,-1) {$+$};
\node at (\a+2,-1) {$+$};
\node at (\a+3,-1) {$+$};
\node at (\a+4,-1) {$\cdot$};
\draw[very thick,->] (-3,-0.5)--(-1,-0.5);
\end{tikzpicture}
\qquad
\begin{tikzpicture}[scale=0.2]
\def\a{0};
\def\b{0};
\draw[step=2, color=green](\a+2,-6+\b)grid(\a+12,-2+\b);
\draw[very thick](\a+3,-2+\b)--(\a+3,-2.5+\b);
\draw[very thick](\a+3,-3.5+\b)--(\a+3,-4+\b);
\draw[very thick](\a+2,-3+\b)--(\a+2.5,-3+\b);
\draw[very thick](\a+3.5,-3+\b)--(\a+4,-3+\b);
\draw[very thick](\a+2.5,-3+\b)arc(270:360:0.5);
\draw[very thick](\a+3.5,-3+\b)arc(90:180:0.5);
\draw[very thick](\a+4,-3+\b)--(\a+6,-3+\b);
\draw[very thick](\a+5,-2+\b)--(\a+5,-4+\b);
\draw[very thick](\a+6,-3+\b)--(\a+8,-3+\b);
\draw[very thick](\a+7,-2+\b)--(\a+7,-4+\b);
\draw[very thick](\a+8,-3+\b)--(\a+10,-3+\b);
\draw[very thick](\a+9,-2+\b)--(\a+9,-4+\b);
\draw[very thick](\a+11,-2+\b)--(\a+11,-2.5+\b);
\draw[very thick](\a+11,-3.5+\b)--(\a+11,-4+\b);
\draw[very thick](\a+10,-3+\b)--(\a+10.5,-3+\b);
\draw[very thick](\a+11.5,-3+\b)--(\a+12,-3+\b);
\draw[very thick](\a+10.5,-3+\b)arc(270:360:0.5);
\draw[very thick](\a+11.5,-3+\b)arc(90:180:0.5);
\draw[very thick](\a+2,-5+\b)--(\a+4,-5+\b);
\draw[very thick](\a+3,-4+\b)--(\a+3,-6+\b);
\draw[very thick](\a+4,-5+\b)--(\a+6,-5+\b);
\draw[very thick](\a+5,-4+\b)--(\a+5,-6+\b);
\draw[very thick](\a+6,-5+\b)--(\a+8,-5+\b);
\draw[very thick](\a+7,-4+\b)--(\a+7,-6+\b);
\draw[very thick](\a+8,-5+\b)--(\a+10,-5+\b);
\draw[very thick](\a+9,-4+\b)--(\a+9,-6+\b);
\draw[very thick](\a+11,-4+\b)--(\a+11,-4.5+\b);
\draw[very thick](\a+11,-5.5+\b)--(\a+11,-6+\b);
\draw[very thick](\a+10,-5+\b)--(\a+10.5,-5+\b);
\draw[very thick](\a+11.5,-5+\b)--(\a+12,-5+\b);
\draw[very thick](\a+10.5,-5+\b)arc(270:360:0.5);
\draw[very thick](\a+11.5,-5+\b)arc(90:180:0.5);

\def\a{-16};
\def\b{0};
\draw[step=2, color=green](\a+2,-6+\b)grid(\a+12,-2+\b);
\draw[very thick](\a+3,-2+\b)--(\a+3,-2.5+\b);
\draw[very thick](\a+3,-3.5+\b)--(\a+3,-4+\b);
\draw[very thick](\a+2,-3+\b)--(\a+2.5,-3+\b);
\draw[very thick](\a+3.5,-3+\b)--(\a+4,-3+\b);
\draw[very thick](\a+2.5,-3+\b)arc(270:360:0.5);
\draw[very thick](\a+3.5,-3+\b)arc(90:180:0.5);
\draw[very thick](\a+4,-3+\b)--(\a+6,-3+\b);
\draw[very thick](\a+5,-2+\b)--(\a+5,-4+\b);
\draw[very thick](\a+6,-3+\b)--(\a+8,-3+\b);
\draw[very thick](\a+7,-2+\b)--(\a+7,-4+\b);
\draw[very thick](\a+8,-3+\b)--(\a+10,-3+\b);
\draw[very thick](\a+9,-2+\b)--(\a+9,-4+\b);
\draw[very thick](\a+10,-3+\b)--(\a+12,-3+\b);
\draw[very thick](\a+11,-2+\b)--(\a+11,-4+\b);
\draw[very thick](\a+3,-4+\b)--(\a+3,-4.5+\b);
\draw[very thick](\a+3,-5.5+\b)--(\a+3,-6+\b);
\draw[very thick](\a+2,-5+\b)--(\a+2.5,-5+\b);
\draw[very thick](\a+3.5,-5+\b)--(\a+4,-5+\b);
\draw[very thick](\a+2.5,-5+\b)arc(270:360:0.5);
\draw[very thick](\a+3.5,-5+\b)arc(90:180:0.5);
\draw[very thick](\a+4,-5+\b)--(\a+6,-5+\b);
\draw[very thick](\a+5,-4+\b)--(\a+5,-6+\b);
\draw[very thick](\a+6,-5+\b)--(\a+8,-5+\b);
\draw[very thick](\a+7,-4+\b)--(\a+7,-6+\b);
\draw[very thick](\a+8,-5+\b)--(\a+10,-5+\b);
\draw[very thick](\a+9,-4+\b)--(\a+9,-6+\b);
\draw[very thick](\a+11,-4+\b)--(\a+11,-4.5+\b);
\draw[very thick](\a+11,-5.5+\b)--(\a+11,-6+\b);
\draw[very thick](\a+10,-5+\b)--(\a+10.5,-5+\b);
\draw[very thick](\a+11.5,-5+\b)--(\a+12,-5+\b);
\draw[very thick](\a+10.5,-5+\b)arc(270:360:0.5);
\draw[very thick](\a+11.5,-5+\b)arc(90:180:0.5);
\draw[very thick,->] (-3,-4)--(1,-4);
\end{tikzpicture}
\caption{Chute moves on pipe dreams move one cross to the left and
preserve the permutation.}
\label{fig:chute-move-PD}
\end{figure}
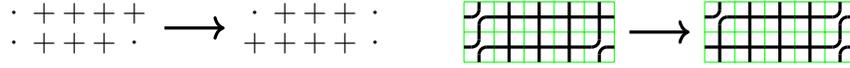

\begin{Theorem}\cite[Thm 3.7]{billey-bergeron}\label{thm:chutes.and.ladders}
For any $w \in S_{n}$, every reduced pipe dream for $w$ can be
obtained from a sequence of ladder moves on $\Dbot(w)$, and every reduced pipe
dream for $w$ can be obtained from a sequence of chute moves on $\Dtop(w)$.
\end{Theorem}

\begin{Example}\label{ex:1432.1}
The reduced pipe dreams for $w=1432$ are shown in
\Cref{fig:1432.ladders}. Note that the code is $c(1432) = (0,2,1,0)$, so $\Dbot(1432)$ is the leftmost pipe dream.  Therefore by Theorem~\ref{thm:chutes.and.ladders}, as computed in \eqref{ex:1} of \Cref{ex:1432}, 
\begin{equation}
\fS_{1432} =  x_1^2 x_2 + x_1^2x_3 + x_1x_2^2 + x_1x_2x_3 + x_2^2x_3.
\end{equation}

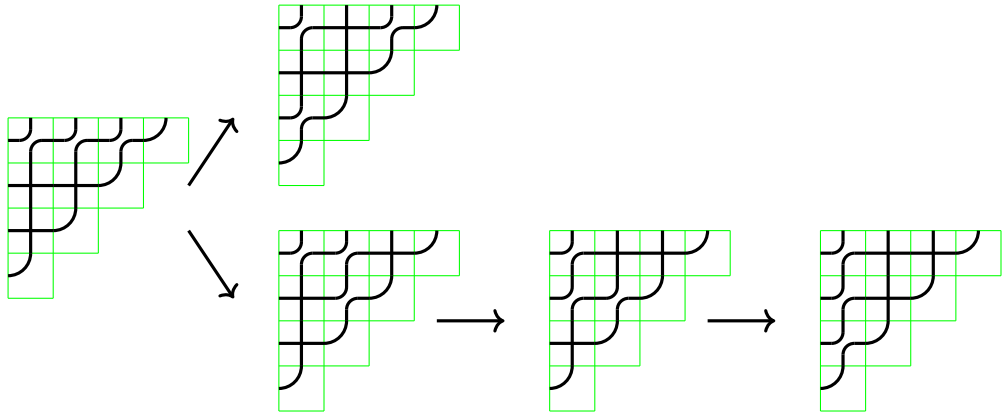
\begin{figure}[h]
\centering
\begin{tikzpicture}[scale=0.300000000000000]
\def\a{0};
\def\b{0};
\draw[color=green](\a+2,-2+\b)--(\a+10,-2+\b);
\draw[color=green](\a+2,-2+\b)--(\a+2,-10+\b);
\draw[color=green](\a+2,-4+\b)--(\a+10,-4+\b);
\draw[color=green](\a+4,-2+\b)--(\a+4,-10+\b);
\draw[color=green](\a+2,-6+\b)--(\a+8,-6+\b);
\draw[color=green](\a+6,-2+\b)--(\a+6,-8+\b);
\draw[color=green](\a+2,-8+\b)--(\a+6,-8+\b);
\draw[color=green](\a+8,-2+\b)--(\a+8,-6+\b);
\draw[color=green](\a+2,-10+\b)--(\a+4,-10+\b);
\draw[color=green](\a+10,-2+\b)--(\a+10,-4+\b);
\draw[very thick](\a+3,-2+\b)--(\a+3,-2.50000000000000+\b);
\draw[very thick](\a+3,-3.50000000000000+\b)--(\a+3,-4+\b);
\draw[very thick](\a+2,-3+\b)--(\a+2.50000000000000,-3+\b);
\draw[very thick](\a+3.50000000000000,-3+\b)--(\a+4,-3+\b);
\draw[very thick](\a+2.50000000000000,-3+\b)arc(270:360:0.500000000000000);
\draw[very thick](\a+3.50000000000000,-3+\b)arc(90:180:0.500000000000000);
\draw[very thick](\a+5,-2+\b)--(\a+5,-2.50000000000000+\b);
\draw[very thick](\a+5,-3.50000000000000+\b)--(\a+5,-4+\b);
\draw[very thick](\a+4,-3+\b)--(\a+4.50000000000000,-3+\b);
\draw[very thick](\a+5.50000000000000,-3+\b)--(\a+6,-3+\b);
\draw[very thick](\a+4.50000000000000,-3+\b)arc(270:360:0.500000000000000);
\draw[very thick](\a+5.50000000000000,-3+\b)arc(90:180:0.500000000000000);
\draw[very thick](\a+7,-2+\b)--(\a+7,-2.50000000000000+\b);
\draw[very thick](\a+7,-3.50000000000000+\b)--(\a+7,-4+\b);
\draw[very thick](\a+6,-3+\b)--(\a+6.50000000000000,-3+\b);
\draw[very thick](\a+7.50000000000000,-3+\b)--(\a+8,-3+\b);
\draw[very thick](\a+6.50000000000000,-3+\b)arc(270:360:0.500000000000000);
\draw[very thick](\a+7.50000000000000,-3+\b)arc(90:180:0.500000000000000);
\draw[very thick](\a+2,-5+\b)--(\a+4,-5+\b);
\draw[very thick](\a+3,-4+\b)--(\a+3,-6+\b);
\draw[very thick](\a+4,-5+\b)--(\a+6,-5+\b);
\draw[very thick](\a+5,-4+\b)--(\a+5,-6+\b);
\draw[very thick](\a+2,-7+\b)--(\a+4,-7+\b);
\draw[very thick](\a+3,-6+\b)--(\a+3,-8+\b);
\draw[very thick](\a+8,-3+\b)arc(270:360:1);
\draw[very thick](\a+6,-5+\b)arc(270:360:1);
\draw[very thick](\a+4,-7+\b)arc(270:360:1);
\draw[very thick](\a+2,-9+\b)arc(270:360:1);

\def\a{12};
\def\b{5};
\draw[color=green](\a+2,-2+\b)--(\a+10,-2+\b);
\draw[color=green](\a+2,-2+\b)--(\a+2,-10+\b);
\draw[color=green](\a+2,-4+\b)--(\a+10,-4+\b);
\draw[color=green](\a+4,-2+\b)--(\a+4,-10+\b);
\draw[color=green](\a+2,-6+\b)--(\a+8,-6+\b);
\draw[color=green](\a+6,-2+\b)--(\a+6,-8+\b);
\draw[color=green](\a+2,-8+\b)--(\a+6,-8+\b);
\draw[color=green](\a+8,-2+\b)--(\a+8,-6+\b);
\draw[color=green](\a+2,-10+\b)--(\a+4,-10+\b);
\draw[color=green](\a+10,-2+\b)--(\a+10,-4+\b);
\draw[very thick](\a+3,-2+\b)--(\a+3,-2.50000000000000+\b);
\draw[very thick](\a+3,-3.50000000000000+\b)--(\a+3,-4+\b);
\draw[very thick](\a+2,-3+\b)--(\a+2.50000000000000,-3+\b);
\draw[very thick](\a+3.50000000000000,-3+\b)--(\a+4,-3+\b);
\draw[very thick](\a+2.50000000000000,-3+\b)arc(270:360:0.500000000000000);
\draw[very thick](\a+3.50000000000000,-3+\b)arc(90:180:0.500000000000000);
\draw[very thick](\a+4,-3+\b)--(\a+6,-3+\b);
\draw[very thick](\a+5,-2+\b)--(\a+5,-4+\b);
\draw[very thick](\a+7,-2+\b)--(\a+7,-2.50000000000000+\b);
\draw[very thick](\a+7,-3.50000000000000+\b)--(\a+7,-4+\b);
\draw[very thick](\a+6,-3+\b)--(\a+6.50000000000000,-3+\b);
\draw[very thick](\a+7.50000000000000,-3+\b)--(\a+8,-3+\b);
\draw[very thick](\a+6.50000000000000,-3+\b)arc(270:360:0.500000000000000);
\draw[very thick](\a+7.50000000000000,-3+\b)arc(90:180:0.500000000000000);
\draw[very thick](\a+2,-5+\b)--(\a+4,-5+\b);
\draw[very thick](\a+3,-4+\b)--(\a+3,-6+\b);
\draw[very thick](\a+4,-5+\b)--(\a+6,-5+\b);
\draw[very thick](\a+5,-4+\b)--(\a+5,-6+\b);
\draw[very thick](\a+3,-6+\b)--(\a+3,-6.50000000000000+\b);
\draw[very thick](\a+3,-7.50000000000000+\b)--(\a+3,-8+\b);
\draw[very thick](\a+2,-7+\b)--(\a+2.50000000000000,-7+\b);
\draw[very thick](\a+3.50000000000000,-7+\b)--(\a+4,-7+\b);
\draw[very thick](\a+2.50000000000000,-7+\b)arc(270:360:0.500000000000000);
\draw[very thick](\a+3.50000000000000,-7+\b)arc(90:180:0.500000000000000);
\draw[very thick](\a+8,-3+\b)arc(270:360:1);
\draw[very thick](\a+6,-5+\b)arc(270:360:1);
\draw[very thick](\a+4,-7+\b)arc(270:360:1);
\draw[very thick](\a+2,-9+\b)arc(270:360:1);

\def\a{12};
\def\b{-5};
\draw[color=green](\a+2,-2+\b)--(\a+10,-2+\b);
\draw[color=green](\a+2,-2+\b)--(\a+2,-10+\b);
\draw[color=green](\a+2,-4+\b)--(\a+10,-4+\b);
\draw[color=green](\a+4,-2+\b)--(\a+4,-10+\b);
\draw[color=green](\a+2,-6+\b)--(\a+8,-6+\b);
\draw[color=green](\a+6,-2+\b)--(\a+6,-8+\b);
\draw[color=green](\a+2,-8+\b)--(\a+6,-8+\b);
\draw[color=green](\a+8,-2+\b)--(\a+8,-6+\b);
\draw[color=green](\a+2,-10+\b)--(\a+4,-10+\b);
\draw[color=green](\a+10,-2+\b)--(\a+10,-4+\b);
\draw[very thick](\a+3,-2+\b)--(\a+3,-2.50000000000000+\b);
\draw[very thick](\a+3,-3.50000000000000+\b)--(\a+3,-4+\b);
\draw[very thick](\a+2,-3+\b)--(\a+2.50000000000000,-3+\b);
\draw[very thick](\a+3.50000000000000,-3+\b)--(\a+4,-3+\b);
\draw[very thick](\a+2.50000000000000,-3+\b)arc(270:360:0.500000000000000);
\draw[very thick](\a+3.50000000000000,-3+\b)arc(90:180:0.500000000000000);
\draw[very thick](\a+5,-2+\b)--(\a+5,-2.50000000000000+\b);
\draw[very thick](\a+5,-3.50000000000000+\b)--(\a+5,-4+\b);
\draw[very thick](\a+4,-3+\b)--(\a+4.50000000000000,-3+\b);
\draw[very thick](\a+5.50000000000000,-3+\b)--(\a+6,-3+\b);
\draw[very thick](\a+4.50000000000000,-3+\b)arc(270:360:0.500000000000000);
\draw[very thick](\a+5.50000000000000,-3+\b)arc(90:180:0.500000000000000);
\draw[very thick](\a+6,-3+\b)--(\a+8,-3+\b);
\draw[very thick](\a+7,-2+\b)--(\a+7,-4+\b);
\draw[very thick](\a+2,-5+\b)--(\a+4,-5+\b);
\draw[very thick](\a+3,-4+\b)--(\a+3,-6+\b);
\draw[very thick](\a+5,-4+\b)--(\a+5,-4.50000000000000+\b);
\draw[very thick](\a+5,-5.50000000000000+\b)--(\a+5,-6+\b);
\draw[very thick](\a+4,-5+\b)--(\a+4.50000000000000,-5+\b);
\draw[very thick](\a+5.50000000000000,-5+\b)--(\a+6,-5+\b);
\draw[very thick](\a+4.50000000000000,-5+\b)arc(270:360:0.500000000000000);
\draw[very thick](\a+5.50000000000000,-5+\b)arc(90:180:0.500000000000000);
\draw[very thick](\a+2,-7+\b)--(\a+4,-7+\b);
\draw[very thick](\a+3,-6+\b)--(\a+3,-8+\b);
\draw[very thick](\a+8,-3+\b)arc(270:360:1);
\draw[very thick](\a+6,-5+\b)arc(270:360:1);
\draw[very thick](\a+4,-7+\b)arc(270:360:1);
\draw[very thick](\a+2,-9+\b)arc(270:360:1);

\def\a{24};
\def\b{-5};
\draw[color=green](\a+2,-2+\b)--(\a+10,-2+\b);
\draw[color=green](\a+2,-2+\b)--(\a+2,-10+\b);
\draw[color=green](\a+2,-4+\b)--(\a+10,-4+\b);
\draw[color=green](\a+4,-2+\b)--(\a+4,-10+\b);
\draw[color=green](\a+2,-6+\b)--(\a+8,-6+\b);
\draw[color=green](\a+6,-2+\b)--(\a+6,-8+\b);
\draw[color=green](\a+2,-8+\b)--(\a+6,-8+\b);
\draw[color=green](\a+8,-2+\b)--(\a+8,-6+\b);
\draw[color=green](\a+2,-10+\b)--(\a+4,-10+\b);
\draw[color=green](\a+10,-2+\b)--(\a+10,-4+\b);
\draw[very thick](\a+3,-2+\b)--(\a+3,-2.50000000000000+\b);
\draw[very thick](\a+3,-3.50000000000000+\b)--(\a+3,-4+\b);
\draw[very thick](\a+2,-3+\b)--(\a+2.50000000000000,-3+\b);
\draw[very thick](\a+3.50000000000000,-3+\b)--(\a+4,-3+\b);
\draw[very thick](\a+2.50000000000000,-3+\b)arc(270:360:0.500000000000000);
\draw[very thick](\a+3.50000000000000,-3+\b)arc(90:180:0.500000000000000);
\draw[very thick](\a+4,-3+\b)--(\a+6,-3+\b);
\draw[very thick](\a+5,-2+\b)--(\a+5,-4+\b);
\draw[very thick](\a+6,-3+\b)--(\a+8,-3+\b);
\draw[very thick](\a+7,-2+\b)--(\a+7,-4+\b);
\draw[very thick](\a+3,-4+\b)--(\a+3,-4.50000000000000+\b);
\draw[very thick](\a+3,-5.50000000000000+\b)--(\a+3,-6+\b);
\draw[very thick](\a+2,-5+\b)--(\a+2.50000000000000,-5+\b);
\draw[very thick](\a+3.50000000000000,-5+\b)--(\a+4,-5+\b);
\draw[very thick](\a+2.50000000000000,-5+\b)arc(270:360:0.500000000000000);
\draw[very thick](\a+3.50000000000000,-5+\b)arc(90:180:0.500000000000000);
\draw[very thick](\a+5,-4+\b)--(\a+5,-4.50000000000000+\b);
\draw[very thick](\a+5,-5.50000000000000+\b)--(\a+5,-6+\b);
\draw[very thick](\a+4,-5+\b)--(\a+4.50000000000000,-5+\b);
\draw[very thick](\a+5.50000000000000,-5+\b)--(\a+6,-5+\b);
\draw[very thick](\a+4.50000000000000,-5+\b)arc(270:360:0.500000000000000);
\draw[very thick](\a+5.50000000000000,-5+\b)arc(90:180:0.500000000000000);
\draw[very thick](\a+2,-7+\b)--(\a+4,-7+\b);
\draw[very thick](\a+3,-6+\b)--(\a+3,-8+\b);
\draw[very thick](\a+8,-3+\b)arc(270:360:1);
\draw[very thick](\a+6,-5+\b)arc(270:360:1);
\draw[very thick](\a+4,-7+\b)arc(270:360:1);
\draw[very thick](\a+2,-9+\b)arc(270:360:1);

\def\a{36};
\def\b{-5};
\draw[color=green](\a+2,-2+\b)--(\a+10,-2+\b);
\draw[color=green](\a+2,-2+\b)--(\a+2,-10+\b);
\draw[color=green](\a+2,-4+\b)--(\a+10,-4+\b);
\draw[color=green](\a+4,-2+\b)--(\a+4,-10+\b);
\draw[color=green](\a+2,-6+\b)--(\a+8,-6+\b);
\draw[color=green](\a+6,-2+\b)--(\a+6,-8+\b);
\draw[color=green](\a+2,-8+\b)--(\a+6,-8+\b);
\draw[color=green](\a+8,-2+\b)--(\a+8,-6+\b);
\draw[color=green](\a+2,-10+\b)--(\a+4,-10+\b);
\draw[color=green](\a+10,-2+\b)--(\a+10,-4+\b);
\draw[very thick](\a+3,-2+\b)--(\a+3,-2.50000000000000+\b);
\draw[very thick](\a+3,-3.50000000000000+\b)--(\a+3,-4+\b);
\draw[very thick](\a+2,-3+\b)--(\a+2.50000000000000,-3+\b);
\draw[very thick](\a+3.50000000000000,-3+\b)--(\a+4,-3+\b);
\draw[very thick](\a+2.50000000000000,-3+\b)arc(270:360:0.500000000000000);
\draw[very thick](\a+3.50000000000000,-3+\b)arc(90:180:0.500000000000000);
\draw[very thick](\a+4,-3+\b)--(\a+6,-3+\b);
\draw[very thick](\a+5,-2+\b)--(\a+5,-4+\b);
\draw[very thick](\a+6,-3+\b)--(\a+8,-3+\b);
\draw[very thick](\a+7,-2+\b)--(\a+7,-4+\b);
\draw[very thick](\a+3,-4+\b)--(\a+3,-4.50000000000000+\b);
\draw[very thick](\a+3,-5.50000000000000+\b)--(\a+3,-6+\b);
\draw[very thick](\a+2,-5+\b)--(\a+2.50000000000000,-5+\b);
\draw[very thick](\a+3.50000000000000,-5+\b)--(\a+4,-5+\b);
\draw[very thick](\a+2.50000000000000,-5+\b)arc(270:360:0.500000000000000);
\draw[very thick](\a+3.50000000000000,-5+\b)arc(90:180:0.500000000000000);
\draw[very thick](\a+4,-5+\b)--(\a+6,-5+\b);
\draw[very thick](\a+5,-4+\b)--(\a+5,-6+\b);
\draw[very thick](\a+3,-6+\b)--(\a+3,-6.50000000000000+\b);
\draw[very thick](\a+3,-7.50000000000000+\b)--(\a+3,-8+\b);
\draw[very thick](\a+2,-7+\b)--(\a+2.50000000000000,-7+\b);
\draw[very thick](\a+3.50000000000000,-7+\b)--(\a+4,-7+\b);
\draw[very thick](\a+2.50000000000000,-7+\b)arc(270:360:0.500000000000000);
\draw[very thick](\a+3.50000000000000,-7+\b)arc(90:180:0.500000000000000);
\draw[very thick](\a+8,-3+\b)arc(270:360:1);
\draw[very thick](\a+6,-5+\b)arc(270:360:1);
\draw[very thick](\a+4,-7+\b)arc(270:360:1);
\draw[very thick](\a+2,-9+\b)arc(270:360:1);

\draw[very thick, ->](10,-5)--(12,-2);
\draw[very thick, ->](10,-7)--(12,-10);
\draw[very thick, ->](21,-11)--(24,-11);
\draw[very thick, ->](33,-11)--(36,-11);
\end{tikzpicture}
\caption{Construction of $\rp(w)$ by ladder moves on $\Dbot(1432)$.}
\label{fig:1432.ladders}
\end{figure}
\end{Example}

Lascoux-Sch\"utzenberger \cite{LS1} also defined the \textit{double
Schubert polynomials} using divided difference operators, similar
to \eqref{eq:divided.difference}, by the formula
\begin{equation}\label{eq:divided.difference.for.doubleSchubs}
\fS_{w}(X;Y)= \fS_{w}(x_{1},\dots ,x_{n};y_{1},\dots , y_{n})=
\begin{cases}\prod_{i+j\leq n}(x_{i}-y_{j}) &\text{if }w=[n,
n-1, \ldots, 1] \\ \partial_i\fS_{ws_i}(X;Y)&\text{if }w(i)<w(i+1). \end{cases}
\end{equation}
Here the divided difference operators act on the $x$-variables, and
treat polynomials in the $y$-variables alone as scalars.
The double Schubert polynomials can be written in terms of the
single Schubert polynomials as 
\[
\fS_{w}(X;Y)= \sum (-1)^{\ell(v)}\fS_{u}(X)\fS_{v}(Y)
\]
where the sum is taken over all factorizations $v^{-1}u=w$ such that
$\ell(u)+\ell(v)=\ell(w)$ \cite[(6.3)]{M2}.  Fomin-Kirillov \cite{FK} also gave the
monomial expansion of the double Schubert polynomials,
\begin{equation}\label{eq:double.schubs}
 \fS_{w}(x_{1},\dots ,x_{n};y_{1},\dots ,
y_{n})= \sum_{D \in \rp(w)} (x-y)^D.
\end{equation}
See also Knutson's co-transition formula in
\cite{knutson.2022} and a variation for the
Grothendieck polynomials in $K$-theory.  Bergeron-Billey also gave a
double Schubert polynomial formula using ladder moves applied to
$\Dbot(w)$ allowing the crossings to move above the first row
\cite[Thm. 4.1]{billey-bergeron}.  The double Schubert polynomials play
a key role in the equivariant Schubert calculus discussed in successive
chapters of this book.  We will return to the double Schubert
polynomials in \Cref{sub:MatrixSchubs} when we discuss matrix Schubert
varieties and the ``naturality'' of pipe dreams.

\begin{Exercise}\label{ex:1432.2}  Use pipe dreams to show 
$\fS_{2413} = x_1 x_2^{2} + x_1^2x_{2}$.
\end{Exercise}

\begin{Exercise}\label{ex:w_{0}}
Prove there is exactly one reduced pipe dream for $w_{0}=[n,n-1,\dots
, 1]$.  
\end{Exercise}

\begin{Exercise}
Prove that the reduced word $\mathrm{r}_{\Dbot(w)}$ is largest in
reverse lexicographic order among all reduced words for $w$.
Similarly, $\mathrm{r}_{\Dtop(w)}$ is the smallest reduced word for
$w$ in lexicographic order.   
\end{Exercise}

\begin{Exercise}
Recall the dominant permutations from \Cref{ex:dominant.1} index the
Schubert polynomials with exactly one monomial in their expansion,
hence they have exactly one reduced pipe dream.  Why are there no
valid ladder moves on $\Dbot$ in this case?  What can you say about
Schubert polynomials with at most $k$ reduced pipe dreams?
\end{Exercise}

\begin{Exercise}
    Use pipe dreams to compute the double Schubert polynomial $\fS_{1432}(X;Y)$.
\end{Exercise}

\begin{Exercise} 
Given a reduced pipe dream $D$ for a permutation in $S_{n}$ and $k \in
[n]$, color the pipes red that exit on rows $1,2,\ldots , k$ and color
the remaining pipes blue.  Remove all blue pipes leaving a tiling
using red elbows, red crosses, and red horizontal line segments.
Contract all of the red horizontal line segments.  Prove the
resulting tiling is a  reduced pipe dream for a permutation in
$S_{k}$.  This algorithm is used by Bergeron-Ceballos-Pilaud to
construct the coproduct that they use to construct a Hopf algebra on pipe
dreams \cite{Bergeron-Ceballos-Pilaud}.
\end{Exercise}

\subsubsection{Little Bumps}\label{sub:Little.Bumps}

Recall that Schubert polynomials were defined in this chapter using
the transition equation of \Cref{t:transitionA}.  Therefore, in order
to prove \Cref{thm:pipedreams} we need to show that the generating
functions for reduced pipe dreams satisfy the same recurrence via a
bijection. We will give a bijection that proves
this recurrence using the Little bump algorithm and a variation on
that theme.  The map itself is quite simple, but there are many
details to check to complete the proof.

\begin{figure}
\centering
\includegraphics[width=1.65in]{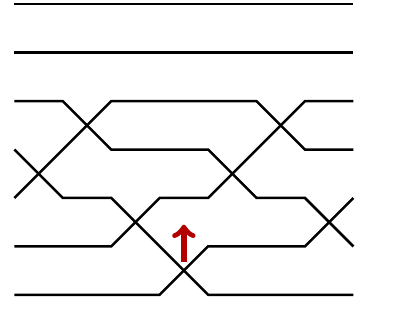}
\includegraphics[width=1.65in]{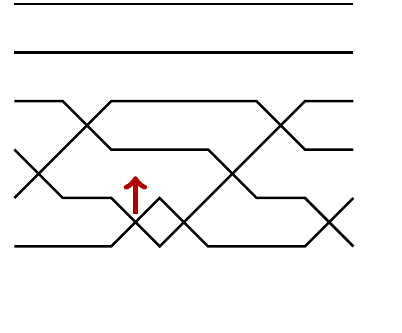}
\includegraphics[width=1.65in]{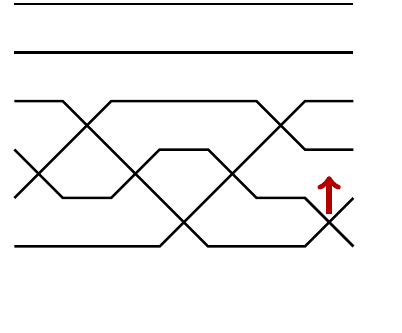} \\
\vspace{.2in}
\includegraphics[width=1.65in]{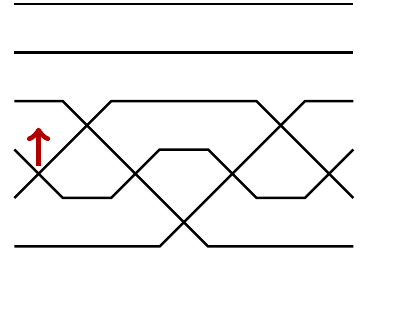}
\includegraphics[width=1.65in]{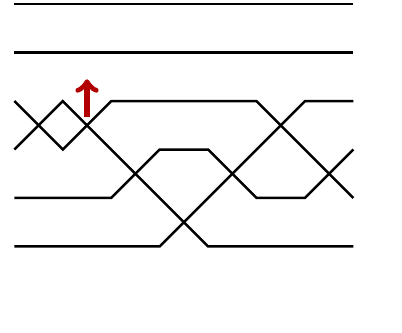}
\includegraphics[width=1.65in]{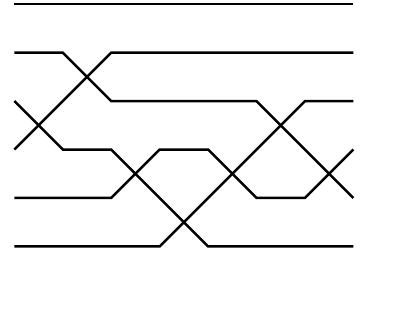}
\caption{An example of the sequence of wiring diagrams for the words
  $\bfa'$ which appear when running the bounded bump algorithm on
  input
  $\bfa =(4, 3, 5, 6, 4, 3, 5), \ \bfb = (2,2,2,2,2,2,2),\ t_0=4,
  \text{ and } \d=-.$ The arrows indicate which crossing will move in
  the next step.  After the first step, row $7$ contains a wire with
  no swaps, which is therefore not shown.}
\label{fig:little.bump}
\end{figure}

The idea of a Little bump is to push a crossing in a reduced wiring
diagram up or down, depending on the chosen direction of travel, and
then to iteratively try to correct the resulting wiring diagram if it
is not reduced by pushing up or down again at another specified
crossing.  \Cref{fig:little.bump} shows an example.  This will
translate to pushing a crossing in a pipe dream left or right on each
step so the monomial weight of the pipe dream $x^{D}$ will be
preserved throughout the process.  We will now discuss the specific
details of the algorithms.

\begin{Definition}
    Let $\mathbf{a} = (a_1, \ldots, a_k)$ be a word.  Define the
\emph{decrement-push}, \emph{increment-push},
\emph{deletion} and \emph{insertion} of $\mathbf{a}$ at
column $t$, respectively, to be \begin{align*}
        \Push^-_t \bfa &= (a_1, \ldots, a_{t-1}, a_t-1, a_{t+1}, \ldots, a_k); \\
        \Push^+_t \bfa &= (a_1, \ldots, a_{t-1}, a_t+1, a_{t+1}, \ldots, a_k); \\
        \Delete_t \bfa &= (a_1, \ldots, a_{t-1},  a_{t+1}, \ldots,
        a_k); \\
        \Insert_t^{x} \bfa &= (a_1, \ldots, a_{t-1},x,  a_{t}, \ldots, a_k).
    \end{align*}
\end{Definition}
	
\begin{Definition} \label{defn:nearly-reduced}
    Let $\mathbf{a}=(a_1, \ldots, a_k)$ be a word, and assume  $1\leq t\leq k$.  If $\Delete_t \bfa$ is reduced, then we say that $\mathbf{a}$ is \emph{nearly reduced at $t$}.
\end{Definition}

The term ``nearly reduced'' is inspired by Lam~\cite[Chapter
  3]{LLMSSZ}, who uses the terminology ``$t$-marked nearly reduced''.  Words that are
nearly reduced at $t$ may or may not also be reduced, and words that
are reduced may not be nearly reduced at $t$. 
However, every reduced word $\bfa$ is nearly reduced at
some index $t$.  For instance, a reduced word $\bfa$ of
length $k$ is nearly reduced at $1$ and at $k$.

\begin{Exercise}\cite{little2003combinatorial,lam-shimozono}
\label{lem:nearly_reduced} Assume $\bfa$ is not reduced, but is nearly
reduced at $t$. Prove that in the wiring diagram of $\bfa$, the two wires
crossing in column $t$ cross in exactly one other column $t'$.
Furthermore,  $\Delete_{t'} \bfa $ is 
reduced.
\end{Exercise}

\begin{Definition}
In the situation of \Cref{lem:nearly_reduced}, we say
that $t'$ \emph{forms a defect with} $t$ in $\bfa$, and
write $\defect_t(\bfa) = t'$.
\end{Definition}

A crucial point is that the definitions of ``reduced'', ``nearly
reduced'', and the ``$\defect$ map'' make sense even if we are given
only the word $\bfa$, but not the corresponding permutation $w\in
S_n$, nor even its size $n$. Indeed, we can take $n$ to be any integer
greater than the largest element of $\bfa$.  We typically draw the
wiring diagram of $\bfa$ using the minimal number of wires.

\begin{Definition}
A word $\bfb$ is a \emph{bounded word} for another word
$\bfa$ if the words have the same length and $1\leq b_i\leq
a_i$ for all $i$. A \emph{bounded pair} (for a permutation
$w$) is an ordered pair $(\bfa,\bfb)$ such that $\bfa$ is
a reduced word (for $w$) and $\bfb$ is a bounded word for
$\bfa$.  Let $\BoundedPairs(w)$ be the set of all bounded
pairs for $w$.
\end{Definition}

For example, for the simple transposition $s_k$,
the set is
$$\BoundedPairs(s_k) = \bigl\{\bigl((k),(i)\bigr): 1\leq i \leq
k\bigr\}.$$ For $w=321$, the set $\BoundedPairs(321)$ has 6 elements, 
\[
((1,2,1),(1,1,1)), ((1,2,1),(1,2,1)), ((2,1,2),(1,1,1)),
\]
\[
((2,1,2),(1,1,2)), ((2,1,2),(2,1,1)), ((2,1,2),(2,1,2)).
\]

\begin{algorithm}[\textbf{Bounded Bump Algorithm},\cite{Billey-Holroyd-Young}]
  \label{algorithm:little bump}\

\noindent \textbf{Input}: $(\bfa, \bfb, t_0, \d)$, where
$\bfa$ is a word that is nearly reduced at $t_0$,
$\bfb$ is a bounded word for $\bfa$, and $\d \in \{-,+\}=\{-1,+1 \}$ is a direction.

\noindent \textbf{Output}: $\Bump^{\d}_{t_0}(\bfa, \bfb) = (\bfa',
\bfb', i, j, \outcome)$, where $\bfa'$ is a reduced word, $\bfb'$ is a
bounded word for $\bfa'$, $i$ is the row and $j$ is the column of the
last crossing pushed in the algorithm, and $\outcome$ is a binary
indicator explained below.

\begin{enumerate}
    \item Initialize $\bfa'\leftarrow \bfa,\, \bfb'
        \leftarrow \bfb, \, t \leftarrow t_0$.
\item Push in direction
$\d$ at column $t$, i.e.\ set $\bfa' \leftarrow
    \Push^{\d}_{t}\bfa'$ and $\bfb' \leftarrow
    \Push^{\d}_{t}\bfb'$.
\item If $b'_t = 0$, return $(\Delete_t \bfa', \Delete_t
    \mathbf{b}', \bfa'_{t}, t, \deleted)$ and
    \textbf{stop}.
\item If $\bfa'$ is reduced, return $(\bfa', \bfb',
    \bfa'_{t}, t, \bumped)$ and \textbf{stop}.
\item Otherwise, set $t \leftarrow \defect_t(\bfa')$ and
    \textbf{return to step 2.}
\end{enumerate}
\end{algorithm}

\begin{Example}\label{ex:little.bump}
Consider the bounded bump algorithm on 
$$\bfa =(4, 3, 5, 6, 4, 3, 5), \ \bfb = (2,2,2,2,2,2,2),\ t_0=4, \text{ and } \d=-.$$
The result is
$$\Bump^{-}_{4}(\bfa, \bfb) = \bigl((3, 2, 4, 5, 4, 3, 4),
(1,1,1,1,2,2,1),2,2, \bumped\bigl).$$ The sequence of pushes used to
obtain the output reduced word $(3, 2, 4, 5, 4, 3, 4)$ is shown in
Figure~\ref{fig:little.bump}.  The corresponding bounded word
$(1,1,1,1,2,2,1)$ is obtained from $\bfb = (2,2,2,2,2,2,2)$ by
decrementing each position corresponding to a column that was pushed
by the bump algorithm.  On the other hand, with input
$\widetilde{\mathbf{b}} = (2,2,2,2,2,2,1)$ and the same $\bfa$ the
bounded bump algorithm stops after the third push in the sequence
because $\widetilde{b}_7=1$, so
$$\Bump^{-}_{4}(\bfa, \widetilde{\mathbf{b}}) =
\bigl((4, 3, 4, 5, 4, 3), (2,2,1,1,2,2), 4, 7, \deleted\bigr).$$
\end{Example}

The original Little bump algorithm, defined by David Little in \cite{little2003combinatorial}, inspired
the bounded bump algorithm.  In his algorithm, the bounded word is essentially fixed to be the same as the input
reduced word.  The key difference is in Step 3.  Instead of deleting a
letter when $b'_t = 0$, add 1 to each letter in $\bfa'$ and return the
resulting reduced word.  Let $\bfa'+1$ denote the word obtained from
$\bfa'$ by adding 1 to each entry.

\begin{algorithm}[\textbf{Little Bump Algorithm},\cite{little2003combinatorial}]
  \label{algorithm:little bump original}\

\noindent \textbf{Input}: $(\bfa,t_{0},\d)$, where $\bfa$ is a reduced
word that is nearly reduced at $t_0$, and $\d \in \{-,+\}$.

\noindent \textbf{Output}: $\LBump^{\d}_{t_0}(\bfa) = \bfa'$, where $\bfa'$ is a reduced word.

\begin{enumerate}
    \item Initialize $\bfa'\leftarrow \bfa,\, t \leftarrow t_0$.
\item Set $\bfa' \leftarrow \Push^{\d}_{t}\bfa'$.
\item If $a'_t = 0$, return $\bfa'+1$ and \textbf{stop}.
\item If $\bfa'$ is reduced, return $\bfa'$ and \textbf{stop}.
\item Otherwise, set $t \leftarrow \defect_t(\bfa')$ and
    \textbf{return to step 2.}
\end{enumerate}
\end{algorithm}

We now make some remarks about the Little bump and bounded bump
algorithms.  Since we will use the bounded bump algorithm below, we
will focus on that variation.  The initial input word $\bfa$ may or
may not be reduced, but, if we reach step $5$ then $\bfa'$ is always
not reduced but nearly reduced at $t$, so the $\defect$ map is defined
on $\bfa'$.

Suppose that the input word $\bfa$ is a word for a permutation $w\in
S_n$. Pushes may in general result in words with elements outside the
interval $[1,n-1]$. Specifically, in the case $\d=+$, step 2 may
result in a word $\bfa'$ with an element $a_t'=n$. As mentioned above,
this can be interpreted as a word for a permutation in $S_{n+1}$.  In
fact, in this case the algorithm will immediately stop at step 4,
since this new word is necessarily reduced.  On the other hand, in the
case $\d=-$, if step 2 ever results in a word with $a_t'=0$, we must
have $b_t'=0$ as well, so the algorithm will immediately stop at step
3, and the $0$ will be deleted.  Note that it is also possible for a
non-zero element of $\bfa$ to be deleted at step 3, since
$b_{i}'<a_{i}'$ is possible.  Thus, the bounded bump algorithm clearly
terminates in a finite number of steps.  In fact, each position gets
pushed at most once during either the bounded bump algorithm or the
Little bump algorithm.

The exercises below collect several technical facts about the bounded
bump algorithm that are analogous to facts proved by Little about
his algorithm \cite{little2003combinatorial}. These statements may be
checked by essentially the same arguments as in
\cite{little2003combinatorial} -- the inclusion of $\bfb$ has scant
effect here.  We encourage readers to find their own proof and peek at
the article as needed.  In each exercise, we assume $\bfa$ is a word
that is nearly reduced at $t$, and $\bfb$ is a bounded word for
$\bfa$, and let $\d \in \{ +, - \}$.  Assume $\Bump^{\d}_{t}(\bfa,
\bfb) = (\bfa',\bfb',i,j,\outcome)$.

\begin{Exercise}\label{item:1}
Suppose $\bfa$ is reduced, not just nearly reduced. Then,
Algorithm~\ref{algorithm:little bump} is reversible in the sense that
we can recover the inputs by negating the direction $\d$ and using the
output parameters $i,j,\outcome$ to reverse process. More
specifically, say $\Bump^{\d}_{t_0}(\bfa, \bfb) = (\bfa', \bfb', i, j,
\outcome)$, and consider the two possible cases for $\outcome$.  If
$\outcome=\deleted$, then it must be that $\d=-1$, so reverse the bump
by inserting an $i$ into position $j$ and initiating the bumping
algorithm there starting with a positive push to recover $(\bfa,\bfb)$
and $t$, denoted
\[
\Bump^{-\d}_{j}(\Insert_{j}^{i} \bfa', \Insert_{j}^{0} \bfb') =
(\bfa, \bfb, \bfa_{t}, t, \bumped).
\]
If $\outcome=\bumped$, then initiate the bumping algorithm starting at
$j$ in the reverse direction
\[
\Bump^{-\d}_j(\bfa', \bfb') = (\bfa, \bfb, \bfa_{t}, t, \bumped).
\]
In particular, both result in a final push which is a $\bumped$ and not
a $\deleted$ outcome in the reverse direction.  
\end{Exercise}

\begin{Exercise}\label{item:2}
If $\bfa \in R(w)$, then $\Delete_t \bfa \in R(w t_{k,l})$, where
$(k<l)$ is the inversion of $w$ whose wires cross in column $t$ of the
right-labeled wiring diagram for $\bfa$.  If $\outcome = \bumped$,
then $\bfa' \in R(w t_{k,l} t_{x,y})$ where $\{x<y\}$ is the crossing
in column $j$ of the word $\bfa'$ for $w t_{k,l}
t_{x,y}$. Furthermore, if $\d=+$, then $l=x$. If $\d=-$, then $k=y$.
\end{Exercise}

\begin{Exercise} \label{item:3} 
Suppose $\Delete_t \mathbf{a} \in R(v)$.  After every iteration of
step 2 in the bounded bump algorithm computing
$\Bump^{\d}_{t}(\bfa, \bfb)$, the pair $(\Delete_t \bfa', \Delete_t
\mathbf{b}')$ is also a bounded pair for $v$.  In particular, if
$\outcome = \deleted$, then $\bfa' \in R(v)$.
\end{Exercise}

\begin{Exercise}
If $\outcome=\bumped$, then the input
and output words $\bfa$ and $\bfa'$ have the same ascent set.  If
$\outcome=\deleted$, then the ascent set of $\Insert_{j}^{i}(\bfa')$
is the same as the ascent set of $\bfa$.
\end{Exercise}

A pipe dream for a permutation $w$ may be interpreted as a bounded
pair of a special type for the same $w$.  To make this precise, we use
the biword $(\mathbf{r}_{D},\mathbf{j}_{D})$ to encode a pipe dream
$D$. Recall, $\mathbf{r}_{D}$ is the corresponding word for $w$, and
$\mathbf{j}_{D}$ is the sequence of column numbers read along rows top
to bottom, right to left as shown in \Cref{fig:rcgraphs}.  Since
$r_{k}=i_{k}+j_{k}-1$ for each $k$, it is clear that $r_{k}\geq j_{k}$
so $\mathbf{j}_{D}$ is a bounded word for $\mathbf{r}_{D}$.  Note,
$\mathbf{j}_{D}$ is a bounded word for $\mathbf{r}_{D}$ even if $D$ is not a reduced pipe dream.   

\begin{Exercise}\label{ex:bounded.pair}
Not every bounded pair corresponds to a pipe dream.  Consider for
example $(121,121)$.  Prove that a bounded pair $(\mathbf{a},
\mathbf{b})=((a_{1},\ldots , a_{p}), (b_{1},\ldots, b_{p}))$
corresponds to a pipe dream if and only if the list
$[(i_1,b_1),\ldots,(i_p,b_p)]$ has $p$ distinct elements which are
increasing in the reading order shown in \Cref{fig:rcgraphs}.
Here $i_k=a_k-b_k+1$ for $1\leq k\leq p$ is the sequence of row
numbers, and $(b_{1},\ldots, b_{p})$ is the sequence of column
numbers.  Equivalently, $(\mathbf{a}, \mathbf{b})$ corresponds to a
pipe dream if and only if the pairs $(i_1,-b_1),\ldots,(i_p,-b_p)$ are
in strictly increasing lex order.
\end{Exercise}

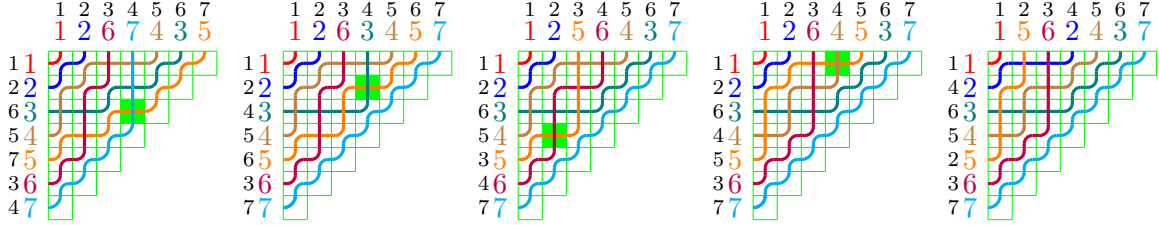
\begin{figure}
\begin{center}
\begin{tikzpicture}[scale=0.16000000000000]
\filldraw[color=green](8,-6)--(10,-6)--(10,-8)--(8,-8);
\draw[color=green](2,-2)--(16,-2);
\draw[color=green](2,-2)--(2,-16);
\draw[color=green](2,-4)--(16,-4);
\draw[color=green](4,-2)--(4,-16);
\draw[color=green](2,-6)--(14,-6);
\draw[color=green](6,-2)--(6,-14);
\draw[color=green](2,-8)--(12,-8);
\draw[color=green](8,-2)--(8,-12);
\draw[color=green](2,-10)--(10,-10);
\draw[color=green](10,-2)--(10,-10);
\draw[color=green](2,-12)--(8,-12);
\draw[color=green](12,-2)--(12,-8);
\draw[color=green](2,-14)--(6,-14);
\draw[color=green](14,-2)--(14,-6);
\draw[color=green](2,-16)--(4,-16);
\draw[color=green](16,-2)--(16,-4);
\draw[very thick, color=orange](14,-3)arc(270:360:1);
\draw[very thick, color=teal](13,-2)--(13,-2.30000000000000);
\draw[very thick, color=orange](13,-3.70000000000000)--(13,-4);
\draw[very thick, color=teal](12,-3)--(12.3000000000000,-3);
\draw[very thick, color=orange](13.7000000000000,-3)--(14,-3);
\draw[very thick, color=teal](12.3000000000000,-3)arc(270:360:0.700000000000000);
\draw[very thick, color=orange](13.7000000000000,-3)arc(90:180:0.700000000000000);
\draw[very thick, color=brown](11,-2)--(11,-2.30000000000000);
\draw[very thick, color=teal](11,-3.70000000000000)--(11,-4);
\draw[very thick, color=brown](10,-3)--(10.3000000000000,-3);
\draw[very thick, color=teal](11.7000000000000,-3)--(12,-3);
\draw[very thick, color=brown](10.3000000000000,-3)arc(270:360:0.700000000000000);
\draw[very thick, color=teal](11.7000000000000,-3)arc(90:180:0.700000000000000);
\draw[very thick, color=brown](8,-3)--(10,-3);
\draw[very thick, color=cyan](9,-2)--(9,-4);
\draw[very thick, color=brown](6,-3)--(8,-3);
\draw[very thick, color=purple](7,-2)--(7,-4);
\draw[very thick, color=blue](5,-2)--(5,-2.30000000000000);
\draw[very thick, color=brown](5,-3.70000000000000)--(5,-4);
\draw[very thick, color=blue](4,-3)--(4.30000000000000,-3);
\draw[very thick, color=brown](5.70000000000000,-3)--(6,-3);
\draw[very thick, color=blue](4.30000000000000,-3)arc(270:360:0.700000000000000);
\draw[very thick, color=brown](5.70000000000000,-3)arc(90:180:0.700000000000000);
\draw[very thick, color=red](3,-2)--(3,-2.30000000000000);
\draw[very thick, color=blue](3,-3.70000000000000)--(3,-4);
\draw[very thick, color=red](2,-3)--(2.30000000000000,-3);
\draw[very thick, color=blue](3.70000000000000,-3)--(4,-3);
\draw[very thick, color=red](2.30000000000000,-3)arc(270:360:0.700000000000000);
\draw[very thick, color=blue](3.70000000000000,-3)arc(90:180:0.700000000000000);
\draw[very thick, color=orange](12,-5)arc(270:360:1);
\draw[very thick, color=teal](11,-4)--(11,-4.30000000000000);
\draw[very thick, color=orange](11,-5.70000000000000)--(11,-6);
\draw[very thick, color=teal](10,-5)--(10.3000000000000,-5);
\draw[very thick, color=orange](11.7000000000000,-5)--(12,-5);
\draw[very thick, color=teal](10.3000000000000,-5)arc(270:360:0.700000000000000);
\draw[very thick, color=orange](11.7000000000000,-5)arc(90:180:0.700000000000000);
\draw[very thick, color=teal](8,-5)--(10,-5);
\draw[very thick, color=cyan](9,-4)--(9,-6);
\draw[very thick, color=purple](7,-4)--(7,-4.30000000000000);
\draw[very thick, color=teal](7,-5.70000000000000)--(7,-6);
\draw[very thick, color=purple](6,-5)--(6.30000000000000,-5);
\draw[very thick, color=teal](7.70000000000000,-5)--(8,-5);
\draw[very thick, color=purple](6.30000000000000,-5)arc(270:360:0.700000000000000);
\draw[very thick, color=teal](7.70000000000000,-5)arc(90:180:0.700000000000000);
\draw[very thick, color=brown](5,-4)--(5,-4.30000000000000);
\draw[very thick, color=purple](5,-5.70000000000000)--(5,-6);
\draw[very thick, color=brown](4,-5)--(4.30000000000000,-5);
\draw[very thick, color=purple](5.70000000000000,-5)--(6,-5);
\draw[very thick, color=brown](4.30000000000000,-5)arc(270:360:0.700000000000000);
\draw[very thick, color=purple](5.70000000000000,-5)arc(90:180:0.700000000000000);
\draw[very thick, color=blue](3,-4)--(3,-4.30000000000000);
\draw[very thick, color=brown](3,-5.70000000000000)--(3,-6);
\draw[very thick, color=blue](2,-5)--(2.30000000000000,-5);
\draw[very thick, color=brown](3.70000000000000,-5)--(4,-5);
\draw[very thick, color=blue](2.30000000000000,-5)arc(270:360:0.700000000000000);
\draw[very thick, color=brown](3.70000000000000,-5)arc(90:180:0.700000000000000);
\draw[very thick, color=orange](10,-7)arc(270:360:1);
\draw[very thick, color=orange](8,-7)--(10,-7);
\draw[very thick, color=cyan](9,-6)--(9,-8);
\draw[very thick, color=teal](7,-6)--(7,-6.30000000000000);
\draw[very thick, color=orange](7,-7.70000000000000)--(7,-8);
\draw[very thick, color=teal](6,-7)--(6.30000000000000,-7);
\draw[very thick, color=orange](7.70000000000000,-7)--(8,-7);
\draw[very thick, color=teal](6.30000000000000,-7)arc(270:360:0.700000000000000);
\draw[very thick, color=orange](7.70000000000000,-7)arc(90:180:0.700000000000000);
\draw[very thick, color=teal](4,-7)--(6,-7);
\draw[very thick, color=purple](5,-6)--(5,-8);
\draw[very thick, color=teal](2,-7)--(4,-7);
\draw[very thick, color=brown](3,-6)--(3,-8);
\draw[very thick, color=cyan](8,-9)arc(270:360:1);
\draw[very thick, color=orange](7,-8)--(7,-8.30000000000000);
\draw[very thick, color=cyan](7,-9.70000000000000)--(7,-10);
\draw[very thick, color=orange](6,-9)--(6.30000000000000,-9);
\draw[very thick, color=cyan](7.70000000000000,-9)--(8,-9);
\draw[very thick, color=orange](6.30000000000000,-9)arc(270:360:0.700000000000000);
\draw[very thick, color=cyan](7.70000000000000,-9)arc(90:180:0.700000000000000);
\draw[very thick, color=orange](4,-9)--(6,-9);
\draw[very thick, color=purple](5,-8)--(5,-10);
\draw[very thick, color=brown](3,-8)--(3,-8.30000000000000);
\draw[very thick, color=orange](3,-9.70000000000000)--(3,-10);
\draw[very thick, color=brown](2,-9)--(2.30000000000000,-9);
\draw[very thick, color=orange](3.70000000000000,-9)--(4,-9);
\draw[very thick, color=brown](2.30000000000000,-9)arc(270:360:0.700000000000000);
\draw[very thick, color=orange](3.70000000000000,-9)arc(90:180:0.700000000000000);
\draw[very thick, color=cyan](6,-11)arc(270:360:1);
\draw[very thick, color=purple](5,-10)--(5,-10.3000000000000);
\draw[very thick, color=cyan](5,-11.7000000000000)--(5,-12);
\draw[very thick, color=purple](4,-11)--(4.30000000000000,-11);
\draw[very thick, color=cyan](5.70000000000000,-11)--(6,-11);
\draw[very thick, color=purple](4.30000000000000,-11)arc(270:360:0.700000000000000);
\draw[very thick, color=cyan](5.70000000000000,-11)arc(90:180:0.700000000000000);
\draw[very thick, color=orange](3,-10)--(3,-10.3000000000000);
\draw[very thick, color=purple](3,-11.7000000000000)--(3,-12);
\draw[very thick, color=orange](2,-11)--(2.30000000000000,-11);
\draw[very thick, color=purple](3.70000000000000,-11)--(4,-11);
\draw[very thick, color=orange](2.30000000000000,-11)arc(270:360:0.700000000000000);
\draw[very thick, color=purple](3.70000000000000,-11)arc(90:180:0.700000000000000);
\draw[very thick, color=cyan](4,-13)arc(270:360:1);
\draw[very thick, color=purple](3,-12)--(3,-12.3000000000000);
\draw[very thick, color=cyan](3,-13.7000000000000)--(3,-14);
\draw[very thick, color=purple](2,-13)--(2.30000000000000,-13);
\draw[very thick, color=cyan](3.70000000000000,-13)--(4,-13);
\draw[very thick, color=purple](2.30000000000000,-13)arc(270:360:0.700000000000000);
\draw[very thick, color=cyan](3.70000000000000,-13)arc(90:180:0.700000000000000);
\draw[very thick, color=cyan](2,-15)arc(270:360:1);
\def\a{1.5};
\def\b{2};
\node[left, color=red] at (2,-3) {\small $1$};
\node[left, color=blue] at (2,-5) {\small $2$};
\node[left, color=teal] at (2,-7) {\small $3$};
\node[left, color=brown] at (2,-9) {\small $4$};
\node[left, color=orange] at (2,-11) {\small $5$};
\node[left, color=purple] at (2,-13) {\small $6$};
\node[left, color=cyan] at (2,-15) {\small $7$};
\node[above, color=red] at (3,-2) {\small $1$};
\node[above, color=blue] at (5,-2) {\small $2$};
\node[above, color=purple] at (7,-2) {\small $6$};
\node[above, color=cyan] at (9,-2) {\small $7$};
\node[above, color=brown] at (11,-2) {\small $4$};
\node[above, color=teal] at (13,-2) {\small $3$};
\node[above, color=orange] at (15,-2) {\small $5$};
\node[left, color=black] at (2-\a,-3) {\tiny $1$};
\node[left, color=black] at (2-\a,-5) {\tiny $2$};
\node[left, color=black] at (2-\a,-7) {\tiny $6$};
\node[left, color=black] at (2-\a,-9) {\tiny $5$};
\node[left, color=black] at (2-\a,-11) {\tiny $7$};
\node[left, color=black] at (2-\a,-13) {\tiny $3$};
\node[left, color=black] at (2-\a,-15) {\tiny $4$};
\node[above, color=black] at (3,-2+\b) {\tiny $1$};
\node[above, color=black] at (5,-2+\b) {\tiny $2$};
\node[above, color=black] at (7,-2+\b) {\tiny $3$};
\node[above, color=black] at (9,-2+\b) {\tiny $4$};
\node[above, color=black] at (11,-2+\b) {\tiny $5$};
\node[above, color=black] at (13,-2+\b) {\tiny $6$};
\node[above, color=black] at (15,-2+\b) {\tiny $7$};
\end{tikzpicture}
\begin{tikzpicture}[scale=0.160000000000000]
\filldraw[color=green](8,-4)--(10,-4)--(10,-6)--(8,-6);
\draw[color=green](2,-2)--(16,-2);
\draw[color=green](2,-2)--(2,-16);
\draw[color=green](2,-4)--(16,-4);
\draw[color=green](4,-2)--(4,-16);
\draw[color=green](2,-6)--(14,-6);
\draw[color=green](6,-2)--(6,-14);
\draw[color=green](2,-8)--(12,-8);
\draw[color=green](8,-2)--(8,-12);
\draw[color=green](2,-10)--(10,-10);
\draw[color=green](10,-2)--(10,-10);
\draw[color=green](2,-12)--(8,-12);
\draw[color=green](12,-2)--(12,-8);
\draw[color=green](2,-14)--(6,-14);
\draw[color=green](14,-2)--(14,-6);
\draw[color=green](2,-16)--(4,-16);
\draw[color=green](16,-2)--(16,-4);
\draw[very thick, color=cyan](14,-3)arc(270:360:1);
\draw[very thick, color=orange](13,-2)--(13,-2.30000000000000);
\draw[very thick, color=cyan](13,-3.70000000000000)--(13,-4);
\draw[very thick, color=orange](12,-3)--(12.3000000000000,-3);
\draw[very thick, color=cyan](13.7000000000000,-3)--(14,-3);
\draw[very thick, color=orange](12.3000000000000,-3)arc(270:360:0.700000000000000);
\draw[very thick, color=cyan](13.7000000000000,-3)arc(90:180:0.700000000000000);
\draw[very thick, color=brown](11,-2)--(11,-2.30000000000000);
\draw[very thick, color=orange](11,-3.70000000000000)--(11,-4);
\draw[very thick, color=brown](10,-3)--(10.3000000000000,-3);
\draw[very thick, color=orange](11.7000000000000,-3)--(12,-3);
\draw[very thick, color=brown](10.3000000000000,-3)arc(270:360:0.700000000000000);
\draw[very thick, color=orange](11.7000000000000,-3)arc(90:180:0.700000000000000);
\draw[very thick, color=brown](8,-3)--(10,-3);
\draw[very thick, color=teal](9,-2)--(9,-4);
\draw[very thick, color=brown](6,-3)--(8,-3);
\draw[very thick, color=purple](7,-2)--(7,-4);
\draw[very thick, color=blue](5,-2)--(5,-2.30000000000000);
\draw[very thick, color=brown](5,-3.70000000000000)--(5,-4);
\draw[very thick, color=blue](4,-3)--(4.30000000000000,-3);
\draw[very thick, color=brown](5.70000000000000,-3)--(6,-3);
\draw[very thick, color=blue](4.30000000000000,-3)arc(270:360:0.700000000000000);
\draw[very thick, color=brown](5.70000000000000,-3)arc(90:180:0.700000000000000);
\draw[very thick, color=red](3,-2)--(3,-2.30000000000000);
\draw[very thick, color=blue](3,-3.70000000000000)--(3,-4);
\draw[very thick, color=red](2,-3)--(2.30000000000000,-3);
\draw[very thick, color=blue](3.70000000000000,-3)--(4,-3);
\draw[very thick, color=red](2.30000000000000,-3)arc(270:360:0.700000000000000);
\draw[very thick, color=blue](3.70000000000000,-3)arc(90:180:0.700000000000000);
\draw[very thick, color=cyan](12,-5)arc(270:360:1);
\draw[very thick, color=orange](11,-4)--(11,-4.30000000000000);
\draw[very thick, color=cyan](11,-5.70000000000000)--(11,-6);
\draw[very thick, color=orange](10,-5)--(10.3000000000000,-5);
\draw[very thick, color=cyan](11.7000000000000,-5)--(12,-5);
\draw[very thick, color=orange](10.3000000000000,-5)arc(270:360:0.700000000000000);
\draw[very thick, color=cyan](11.7000000000000,-5)arc(90:180:0.700000000000000);
\draw[very thick, color=orange](8,-5)--(10,-5);
\draw[very thick, color=teal](9,-4)--(9,-6);
\draw[very thick, color=purple](7,-4)--(7,-4.30000000000000);
\draw[very thick, color=orange](7,-5.70000000000000)--(7,-6);
\draw[very thick, color=purple](6,-5)--(6.30000000000000,-5);
\draw[very thick, color=orange](7.70000000000000,-5)--(8,-5);
\draw[very thick, color=purple](6.30000000000000,-5)arc(270:360:0.700000000000000);
\draw[very thick, color=orange](7.70000000000000,-5)arc(90:180:0.700000000000000);
\draw[very thick, color=brown](5,-4)--(5,-4.30000000000000);
\draw[very thick, color=purple](5,-5.70000000000000)--(5,-6);
\draw[very thick, color=brown](4,-5)--(4.30000000000000,-5);
\draw[very thick, color=purple](5.70000000000000,-5)--(6,-5);
\draw[very thick, color=brown](4.30000000000000,-5)arc(270:360:0.700000000000000);
\draw[very thick, color=purple](5.70000000000000,-5)arc(90:180:0.700000000000000);
\draw[very thick, color=blue](3,-4)--(3,-4.30000000000000);
\draw[very thick, color=brown](3,-5.70000000000000)--(3,-6);
\draw[very thick, color=blue](2,-5)--(2.30000000000000,-5);
\draw[very thick, color=brown](3.70000000000000,-5)--(4,-5);
\draw[very thick, color=blue](2.30000000000000,-5)arc(270:360:0.700000000000000);
\draw[very thick, color=brown](3.70000000000000,-5)arc(90:180:0.700000000000000);
\draw[very thick, color=cyan](10,-7)arc(270:360:1);
\draw[very thick, color=teal](9,-6)--(9,-6.30000000000000);
\draw[very thick, color=cyan](9,-7.70000000000000)--(9,-8);
\draw[very thick, color=teal](8,-7)--(8.30000000000000,-7);
\draw[very thick, color=cyan](9.70000000000000,-7)--(10,-7);
\draw[very thick, color=teal](8.30000000000000,-7)arc(270:360:0.700000000000000);
\draw[very thick, color=cyan](9.70000000000000,-7)arc(90:180:0.700000000000000);
\draw[very thick, color=teal](6,-7)--(8,-7);
\draw[very thick, color=orange](7,-6)--(7,-8);
\draw[very thick, color=teal](4,-7)--(6,-7);
\draw[very thick, color=purple](5,-6)--(5,-8);
\draw[very thick, color=teal](2,-7)--(4,-7);
\draw[very thick, color=brown](3,-6)--(3,-8);
\draw[very thick, color=cyan](8,-9)arc(270:360:1);
\draw[very thick, color=orange](7,-8)--(7,-8.30000000000000);
\draw[very thick, color=cyan](7,-9.70000000000000)--(7,-10);
\draw[very thick, color=orange](6,-9)--(6.30000000000000,-9);
\draw[very thick, color=cyan](7.70000000000000,-9)--(8,-9);
\draw[very thick, color=orange](6.30000000000000,-9)arc(270:360:0.700000000000000);
\draw[very thick, color=cyan](7.70000000000000,-9)arc(90:180:0.700000000000000);
\draw[very thick, color=orange](4,-9)--(6,-9);
\draw[very thick, color=purple](5,-8)--(5,-10);
\draw[very thick, color=brown](3,-8)--(3,-8.30000000000000);
\draw[very thick, color=orange](3,-9.70000000000000)--(3,-10);
\draw[very thick, color=brown](2,-9)--(2.30000000000000,-9);
\draw[very thick, color=orange](3.70000000000000,-9)--(4,-9);
\draw[very thick, color=brown](2.30000000000000,-9)arc(270:360:0.700000000000000);
\draw[very thick, color=orange](3.70000000000000,-9)arc(90:180:0.700000000000000);
\draw[very thick, color=cyan](6,-11)arc(270:360:1);
\draw[very thick, color=purple](5,-10)--(5,-10.3000000000000);
\draw[very thick, color=cyan](5,-11.7000000000000)--(5,-12);
\draw[very thick, color=purple](4,-11)--(4.30000000000000,-11);
\draw[very thick, color=cyan](5.70000000000000,-11)--(6,-11);
\draw[very thick, color=purple](4.30000000000000,-11)arc(270:360:0.700000000000000);
\draw[very thick, color=cyan](5.70000000000000,-11)arc(90:180:0.700000000000000);
\draw[very thick, color=orange](3,-10)--(3,-10.3000000000000);
\draw[very thick, color=purple](3,-11.7000000000000)--(3,-12);
\draw[very thick, color=orange](2,-11)--(2.30000000000000,-11);
\draw[very thick, color=purple](3.70000000000000,-11)--(4,-11);
\draw[very thick, color=orange](2.30000000000000,-11)arc(270:360:0.700000000000000);
\draw[very thick, color=purple](3.70000000000000,-11)arc(90:180:0.700000000000000);
\draw[very thick, color=cyan](4,-13)arc(270:360:1);
\draw[very thick, color=purple](3,-12)--(3,-12.3000000000000);
\draw[very thick, color=cyan](3,-13.7000000000000)--(3,-14);
\draw[very thick, color=purple](2,-13)--(2.30000000000000,-13);
\draw[very thick, color=cyan](3.70000000000000,-13)--(4,-13);
\draw[very thick, color=purple](2.30000000000000,-13)arc(270:360:0.700000000000000);
\draw[very thick, color=cyan](3.70000000000000,-13)arc(90:180:0.700000000000000);
\draw[very thick, color=cyan](2,-15)arc(270:360:1);
\draw[very thick, color=cyan](2,-15)arc(270:360:1);
\node[left, color=red] at (2,-3) {\small $1$};
\node[left, color=blue] at (2,-5) {\small $2$};
\node[left, color=teal] at (2,-7) {\small $3$};
\node[left, color=brown] at (2,-9) {\small $4$};
\node[left, color=orange] at (2,-11) {\small $5$};
\node[left, color=purple] at (2,-13) {\small $6$};
\node[left, color=cyan] at (2,-15) {\small $7$};
\node[above, color=red] at (3,-2) {\small $1$};
\node[above, color=blue] at (5,-2) {\small $2$};
\node[above, color=purple] at (7,-2) {\small $6$};
\node[above, color=teal] at (9,-2) {\small $3$};
\node[above, color=brown] at (11,-2) {\small $4$};
\node[above, color=orange] at (13,-2) {\small $5$};
\node[above, color=cyan] at (15,-2) {\small $7$};
\def\a{1.5};
\def\b{2};
\node[left, color=black] at (2-\a,-3) {\tiny $1$};
\node[left, color=black] at (2-\a,-5) {\tiny $2$};
\node[left, color=black] at (2-\a,-7) {\tiny $4$};
\node[left, color=black] at (2-\a,-9) {\tiny $5$};
\node[left, color=black] at (2-\a,-11) {\tiny $6$};
\node[left, color=black] at (2-\a,-13) {\tiny $3$};
\node[left, color=black] at (2-\a,-15) {\tiny $7$};
\node[above, color=black] at (3,-2+\b) {\tiny $1$};
\node[above, color=black] at (5,-2+\b) {\tiny $2$};
\node[above, color=black] at (7,-2+\b) {\tiny $3$};
\node[above, color=black] at (9,-2+\b) {\tiny $4$};
\node[above, color=black] at (11,-2+\b) {\tiny $5$};
\node[above, color=black] at (13,-2+\b) {\tiny $6$};
\node[above, color=black] at (15,-2+\b) {\tiny $7$};
\end{tikzpicture}
\begin{tikzpicture}[scale=0.160000000000000]
\filldraw[color=green](4,-8)--(6,-8)--(6,-10)--(4,-10);
\draw[color=green](2,-2)--(16,-2);
\draw[color=green](2,-2)--(2,-16);
\draw[color=green](2,-4)--(16,-4);
\draw[color=green](4,-2)--(4,-16);
\draw[color=green](2,-6)--(14,-6);
\draw[color=green](6,-2)--(6,-14);
\draw[color=green](2,-8)--(12,-8);
\draw[color=green](8,-2)--(8,-12);
\draw[color=green](2,-10)--(10,-10);
\draw[color=green](10,-2)--(10,-10);
\draw[color=green](2,-12)--(8,-12);
\draw[color=green](12,-2)--(12,-8);
\draw[color=green](2,-14)--(6,-14);
\draw[color=green](14,-2)--(14,-6);
\draw[color=green](2,-16)--(4,-16);
\draw[color=green](16,-2)--(16,-4);
\draw[very thick, color=cyan](14,-3)arc(270:360:1);
\draw[very thick, color=teal](13,-2)--(13,-2.30000000000000);
\draw[very thick, color=cyan](13,-3.70000000000000)--(13,-4);
\draw[very thick, color=teal](12,-3)--(12.3000000000000,-3);
\draw[very thick, color=cyan](13.7000000000000,-3)--(14,-3);
\draw[very thick, color=teal](12.3000000000000,-3)arc(270:360:0.700000000000000);
\draw[very thick, color=cyan](13.7000000000000,-3)arc(90:180:0.700000000000000);
\draw[very thick, color=brown](11,-2)--(11,-2.30000000000000);
\draw[very thick, color=teal](11,-3.70000000000000)--(11,-4);
\draw[very thick, color=brown](10,-3)--(10.3000000000000,-3);
\draw[very thick, color=teal](11.7000000000000,-3)--(12,-3);
\draw[very thick, color=brown](10.3000000000000,-3)arc(270:360:0.700000000000000);
\draw[very thick, color=teal](11.7000000000000,-3)arc(90:180:0.700000000000000);
\draw[very thick, color=brown](8,-3)--(10,-3);
\draw[very thick, color=purple](9,-2)--(9,-4);
\draw[very thick, color=brown](6,-3)--(8,-3);
\draw[very thick, color=orange](7,-2)--(7,-4);
\draw[very thick, color=blue](5,-2)--(5,-2.30000000000000);
\draw[very thick, color=brown](5,-3.70000000000000)--(5,-4);
\draw[very thick, color=blue](4,-3)--(4.30000000000000,-3);
\draw[very thick, color=brown](5.70000000000000,-3)--(6,-3);
\draw[very thick, color=blue](4.30000000000000,-3)arc(270:360:0.700000000000000);
\draw[very thick, color=brown](5.70000000000000,-3)arc(90:180:0.700000000000000);
\draw[very thick, color=red](3,-2)--(3,-2.30000000000000);
\draw[very thick, color=blue](3,-3.70000000000000)--(3,-4);
\draw[very thick, color=red](2,-3)--(2.30000000000000,-3);
\draw[very thick, color=blue](3.70000000000000,-3)--(4,-3);
\draw[very thick, color=red](2.30000000000000,-3)arc(270:360:0.700000000000000);
\draw[very thick, color=blue](3.70000000000000,-3)arc(90:180:0.700000000000000);
\draw[very thick, color=cyan](12,-5)arc(270:360:1);
\draw[very thick, color=teal](11,-4)--(11,-4.30000000000000);
\draw[very thick, color=cyan](11,-5.70000000000000)--(11,-6);
\draw[very thick, color=teal](10,-5)--(10.3000000000000,-5);
\draw[very thick, color=cyan](11.7000000000000,-5)--(12,-5);
\draw[very thick, color=teal](10.3000000000000,-5)arc(270:360:0.700000000000000);
\draw[very thick, color=cyan](11.7000000000000,-5)arc(90:180:0.700000000000000);
\draw[very thick, color=purple](9,-4)--(9,-4.30000000000000);
\draw[very thick, color=teal](9,-5.70000000000000)--(9,-6);
\draw[very thick, color=purple](8,-5)--(8.30000000000000,-5);
\draw[very thick, color=teal](9.70000000000000,-5)--(10,-5);
\draw[very thick, color=purple](8.30000000000000,-5)arc(270:360:0.700000000000000);
\draw[very thick, color=teal](9.70000000000000,-5)arc(90:180:0.700000000000000);
\draw[very thick, color=purple](6,-5)--(8,-5);
\draw[very thick, color=orange](7,-4)--(7,-6);
\draw[very thick, color=brown](5,-4)--(5,-4.30000000000000);
\draw[very thick, color=purple](5,-5.70000000000000)--(5,-6);
\draw[very thick, color=brown](4,-5)--(4.30000000000000,-5);
\draw[very thick, color=purple](5.70000000000000,-5)--(6,-5);
\draw[very thick, color=brown](4.30000000000000,-5)arc(270:360:0.700000000000000);
\draw[very thick, color=purple](5.70000000000000,-5)arc(90:180:0.700000000000000);
\draw[very thick, color=blue](3,-4)--(3,-4.30000000000000);
\draw[very thick, color=brown](3,-5.70000000000000)--(3,-6);
\draw[very thick, color=blue](2,-5)--(2.30000000000000,-5);
\draw[very thick, color=brown](3.70000000000000,-5)--(4,-5);
\draw[very thick, color=blue](2.30000000000000,-5)arc(270:360:0.700000000000000);
\draw[very thick, color=brown](3.70000000000000,-5)arc(90:180:0.700000000000000);
\draw[very thick, color=cyan](10,-7)arc(270:360:1);
\draw[very thick, color=teal](9,-6)--(9,-6.30000000000000);
\draw[very thick, color=cyan](9,-7.70000000000000)--(9,-8);
\draw[very thick, color=teal](8,-7)--(8.30000000000000,-7);
\draw[very thick, color=cyan](9.70000000000000,-7)--(10,-7);
\draw[very thick, color=teal](8.30000000000000,-7)arc(270:360:0.700000000000000);
\draw[very thick, color=cyan](9.70000000000000,-7)arc(90:180:0.700000000000000);
\draw[very thick, color=teal](6,-7)--(8,-7);
\draw[very thick, color=orange](7,-6)--(7,-8);
\draw[very thick, color=teal](4,-7)--(6,-7);
\draw[very thick, color=purple](5,-6)--(5,-8);
\draw[very thick, color=teal](2,-7)--(4,-7);
\draw[very thick, color=brown](3,-6)--(3,-8);
\draw[very thick, color=cyan](8,-9)arc(270:360:1);
\draw[very thick, color=orange](7,-8)--(7,-8.30000000000000);
\draw[very thick, color=cyan](7,-9.70000000000000)--(7,-10);
\draw[very thick, color=orange](6,-9)--(6.30000000000000,-9);
\draw[very thick, color=cyan](7.70000000000000,-9)--(8,-9);
\draw[very thick, color=orange](6.30000000000000,-9)arc(270:360:0.700000000000000);
\draw[very thick, color=cyan](7.70000000000000,-9)arc(90:180:0.700000000000000);
\draw[very thick, color=orange](4,-9)--(6,-9);
\draw[very thick, color=purple](5,-8)--(5,-10);
\draw[very thick, color=brown](3,-8)--(3,-8.30000000000000);
\draw[very thick, color=orange](3,-9.70000000000000)--(3,-10);
\draw[very thick, color=brown](2,-9)--(2.30000000000000,-9);
\draw[very thick, color=orange](3.70000000000000,-9)--(4,-9);
\draw[very thick, color=brown](2.30000000000000,-9)arc(270:360:0.700000000000000);
\draw[very thick, color=orange](3.70000000000000,-9)arc(90:180:0.700000000000000);
\draw[very thick, color=cyan](6,-11)arc(270:360:1);
\draw[very thick, color=purple](5,-10)--(5,-10.3000000000000);
\draw[very thick, color=cyan](5,-11.7000000000000)--(5,-12);
\draw[very thick, color=purple](4,-11)--(4.30000000000000,-11);
\draw[very thick, color=cyan](5.70000000000000,-11)--(6,-11);
\draw[very thick, color=purple](4.30000000000000,-11)arc(270:360:0.700000000000000);
\draw[very thick, color=cyan](5.70000000000000,-11)arc(90:180:0.700000000000000);
\draw[very thick, color=orange](3,-10)--(3,-10.3000000000000);
\draw[very thick, color=purple](3,-11.7000000000000)--(3,-12);
\draw[very thick, color=orange](2,-11)--(2.30000000000000,-11);
\draw[very thick, color=purple](3.70000000000000,-11)--(4,-11);
\draw[very thick, color=orange](2.30000000000000,-11)arc(270:360:0.700000000000000);
\draw[very thick, color=purple](3.70000000000000,-11)arc(90:180:0.700000000000000);
\draw[very thick, color=cyan](4,-13)arc(270:360:1);
\draw[very thick, color=purple](3,-12)--(3,-12.3000000000000);
\draw[very thick, color=cyan](3,-13.7000000000000)--(3,-14);
\draw[very thick, color=purple](2,-13)--(2.30000000000000,-13);
\draw[very thick, color=cyan](3.70000000000000,-13)--(4,-13);
\draw[very thick, color=purple](2.30000000000000,-13)arc(270:360:0.700000000000000);
\draw[very thick, color=cyan](3.70000000000000,-13)arc(90:180:0.700000000000000);
\draw[very thick, color=cyan](2,-15)arc(270:360:1);
\node[left, color=red] at (2,-3) {\small $1$};
\node[left, color=blue] at (2,-5) {\small $2$};
\node[left, color=teal] at (2,-7) {\small $3$};
\node[left, color=brown] at (2,-9) {\small $4$};
\node[left, color=orange] at (2,-11) {\small $5$};
\node[left, color=purple] at (2,-13) {\small $6$};
\node[left, color=cyan] at (2,-15) {\small $7$};
\node[above, color=red] at (3,-2) {\small $1$};
\node[above, color=blue] at (5,-2) {\small $2$};
\node[above, color=orange] at (7,-2) {\small $5$};
\node[above, color=purple] at (9,-2) {\small $6$};
\node[above, color=brown] at (11,-2) {\small $4$};
\node[above, color=teal] at (13,-2) {\small $3$};
\node[above, color=cyan] at (15,-2) {\small $7$};
\def\a{1.5};
\def\b{2};
\node[left, color=black] at (2-\a,-3) {\tiny $1$};
\node[left, color=black] at (2-\a,-5) {\tiny $2$};
\node[left, color=black] at (2-\a,-7) {\tiny $6$};
\node[left, color=black] at (2-\a,-9) {\tiny $5$};
\node[left, color=black] at (2-\a,-11) {\tiny $3$};
\node[left, color=black] at (2-\a,-13) {\tiny $4$};
\node[left, color=black] at (2-\a,-15) {\tiny $7$};
\node[above, color=black] at (3,-2+\b) {\tiny $1$};
\node[above, color=black] at (5,-2+\b) {\tiny $2$};
\node[above, color=black] at (7,-2+\b) {\tiny $3$};
\node[above, color=black] at (9,-2+\b) {\tiny $4$};
\node[above, color=black] at (11,-2+\b) {\tiny $5$};
\node[above, color=black] at (13,-2+\b) {\tiny $6$};
\node[above, color=black] at (15,-2+\b) {\tiny $7$};
\end{tikzpicture}
\begin{tikzpicture}[scale=0.160000000000000]
\filldraw[color=green](8,-2)--(10,-2)--(10,-4)--(8,-4);
\draw[color=green](2,-2)--(16,-2);
\draw[color=green](2,-2)--(2,-16);
\draw[color=green](2,-4)--(16,-4);
\draw[color=green](4,-2)--(4,-16);
\draw[color=green](2,-6)--(14,-6);
\draw[color=green](6,-2)--(6,-14);
\draw[color=green](2,-8)--(12,-8);
\draw[color=green](8,-2)--(8,-12);
\draw[color=green](2,-10)--(10,-10);
\draw[color=green](10,-2)--(10,-10);
\draw[color=green](2,-12)--(8,-12);
\draw[color=green](12,-2)--(12,-8);
\draw[color=green](2,-14)--(6,-14);
\draw[color=green](14,-2)--(14,-6);
\draw[color=green](2,-16)--(4,-16);
\draw[color=green](16,-2)--(16,-4);
\draw[very thick, color=cyan](14,-3)arc(270:360:1);
\draw[very thick, color=teal](13,-2)--(13,-2.30000000000000);
\draw[very thick, color=cyan](13,-3.70000000000000)--(13,-4);
\draw[very thick, color=teal](12,-3)--(12.3000000000000,-3);
\draw[very thick, color=cyan](13.7000000000000,-3)--(14,-3);
\draw[very thick, color=teal](12.3000000000000,-3)arc(270:360:0.700000000000000);
\draw[very thick, color=cyan](13.7000000000000,-3)arc(90:180:0.700000000000000);
\draw[very thick, color=orange](11,-2)--(11,-2.30000000000000);
\draw[very thick, color=teal](11,-3.70000000000000)--(11,-4);
\draw[very thick, color=orange](10,-3)--(10.3000000000000,-3);
\draw[very thick, color=teal](11.7000000000000,-3)--(12,-3);
\draw[very thick, color=orange](10.3000000000000,-3)arc(270:360:0.700000000000000);
\draw[very thick, color=teal](11.7000000000000,-3)arc(90:180:0.700000000000000);
\draw[very thick, color=orange](8,-3)--(10,-3);
\draw[very thick, color=brown](9,-2)--(9,-4);
\draw[very thick, color=orange](6,-3)--(8,-3);
\draw[very thick, color=purple](7,-2)--(7,-4);
\draw[very thick, color=blue](5,-2)--(5,-2.30000000000000);
\draw[very thick, color=orange](5,-3.70000000000000)--(5,-4);
\draw[very thick, color=blue](4,-3)--(4.30000000000000,-3);
\draw[very thick, color=orange](5.70000000000000,-3)--(6,-3);
\draw[very thick, color=blue](4.30000000000000,-3)arc(270:360:0.700000000000000);
\draw[very thick, color=orange](5.70000000000000,-3)arc(90:180:0.700000000000000);
\draw[very thick, color=red](3,-2)--(3,-2.30000000000000);
\draw[very thick, color=blue](3,-3.70000000000000)--(3,-4);
\draw[very thick, color=red](2,-3)--(2.30000000000000,-3);
\draw[very thick, color=blue](3.70000000000000,-3)--(4,-3);
\draw[very thick, color=red](2.30000000000000,-3)arc(270:360:0.700000000000000);
\draw[very thick, color=blue](3.70000000000000,-3)arc(90:180:0.700000000000000);
\draw[very thick, color=cyan](12,-5)arc(270:360:1);
\draw[very thick, color=teal](11,-4)--(11,-4.30000000000000);
\draw[very thick, color=cyan](11,-5.70000000000000)--(11,-6);
\draw[very thick, color=teal](10,-5)--(10.3000000000000,-5);
\draw[very thick, color=cyan](11.7000000000000,-5)--(12,-5);
\draw[very thick, color=teal](10.3000000000000,-5)arc(270:360:0.700000000000000);
\draw[very thick, color=cyan](11.7000000000000,-5)arc(90:180:0.700000000000000);
\draw[very thick, color=brown](9,-4)--(9,-4.30000000000000);
\draw[very thick, color=teal](9,-5.70000000000000)--(9,-6);
\draw[very thick, color=brown](8,-5)--(8.30000000000000,-5);
\draw[very thick, color=teal](9.70000000000000,-5)--(10,-5);
\draw[very thick, color=brown](8.30000000000000,-5)arc(270:360:0.700000000000000);
\draw[very thick, color=teal](9.70000000000000,-5)arc(90:180:0.700000000000000);
\draw[very thick, color=brown](6,-5)--(8,-5);
\draw[very thick, color=purple](7,-4)--(7,-6);
\draw[very thick, color=orange](5,-4)--(5,-4.30000000000000);
\draw[very thick, color=brown](5,-5.70000000000000)--(5,-6);
\draw[very thick, color=orange](4,-5)--(4.30000000000000,-5);
\draw[very thick, color=brown](5.70000000000000,-5)--(6,-5);
\draw[very thick, color=orange](4.30000000000000,-5)arc(270:360:0.700000000000000);
\draw[very thick, color=brown](5.70000000000000,-5)arc(90:180:0.700000000000000);
\draw[very thick, color=blue](3,-4)--(3,-4.30000000000000);
\draw[very thick, color=orange](3,-5.70000000000000)--(3,-6);
\draw[very thick, color=blue](2,-5)--(2.30000000000000,-5);
\draw[very thick, color=orange](3.70000000000000,-5)--(4,-5);
\draw[very thick, color=blue](2.30000000000000,-5)arc(270:360:0.700000000000000);
\draw[very thick, color=orange](3.70000000000000,-5)arc(90:180:0.700000000000000);
\draw[very thick, color=cyan](10,-7)arc(270:360:1);
\draw[very thick, color=teal](9,-6)--(9,-6.30000000000000);
\draw[very thick, color=cyan](9,-7.70000000000000)--(9,-8);
\draw[very thick, color=teal](8,-7)--(8.30000000000000,-7);
\draw[very thick, color=cyan](9.70000000000000,-7)--(10,-7);
\draw[very thick, color=teal](8.30000000000000,-7)arc(270:360:0.700000000000000);
\draw[very thick, color=cyan](9.70000000000000,-7)arc(90:180:0.700000000000000);
\draw[very thick, color=teal](6,-7)--(8,-7);
\draw[very thick, color=purple](7,-6)--(7,-8);
\draw[very thick, color=teal](4,-7)--(6,-7);
\draw[very thick, color=brown](5,-6)--(5,-8);
\draw[very thick, color=teal](2,-7)--(4,-7);
\draw[very thick, color=orange](3,-6)--(3,-8);
\draw[very thick, color=cyan](8,-9)arc(270:360:1);
\draw[very thick, color=purple](7,-8)--(7,-8.30000000000000);
\draw[very thick, color=cyan](7,-9.70000000000000)--(7,-10);
\draw[very thick, color=purple](6,-9)--(6.30000000000000,-9);
\draw[very thick, color=cyan](7.70000000000000,-9)--(8,-9);
\draw[very thick, color=purple](6.30000000000000,-9)arc(270:360:0.700000000000000);
\draw[very thick, color=cyan](7.70000000000000,-9)arc(90:180:0.700000000000000);
\draw[very thick, color=brown](5,-8)--(5,-8.30000000000000);
\draw[very thick, color=purple](5,-9.70000000000000)--(5,-10);
\draw[very thick, color=brown](4,-9)--(4.30000000000000,-9);
\draw[very thick, color=purple](5.70000000000000,-9)--(6,-9);
\draw[very thick, color=brown](4.30000000000000,-9)arc(270:360:0.700000000000000);
\draw[very thick, color=purple](5.70000000000000,-9)arc(90:180:0.700000000000000);
\draw[very thick, color=brown](2,-9)--(4,-9);
\draw[very thick, color=orange](3,-8)--(3,-10);
\draw[very thick, color=cyan](6,-11)arc(270:360:1);
\draw[very thick, color=purple](5,-10)--(5,-10.3000000000000);
\draw[very thick, color=cyan](5,-11.7000000000000)--(5,-12);
\draw[very thick, color=purple](4,-11)--(4.30000000000000,-11);
\draw[very thick, color=cyan](5.70000000000000,-11)--(6,-11);
\draw[very thick, color=purple](4.30000000000000,-11)arc(270:360:0.700000000000000);
\draw[very thick, color=cyan](5.70000000000000,-11)arc(90:180:0.700000000000000);
\draw[very thick, color=orange](3,-10)--(3,-10.3000000000000);
\draw[very thick, color=purple](3,-11.7000000000000)--(3,-12);
\draw[very thick, color=orange](2,-11)--(2.30000000000000,-11);
\draw[very thick, color=purple](3.70000000000000,-11)--(4,-11);
\draw[very thick, color=orange](2.30000000000000,-11)arc(270:360:0.700000000000000);
\draw[very thick, color=purple](3.70000000000000,-11)arc(90:180:0.700000000000000);
\draw[very thick, color=cyan](4,-13)arc(270:360:1);
\draw[very thick, color=purple](3,-12)--(3,-12.3000000000000);
\draw[very thick, color=cyan](3,-13.7000000000000)--(3,-14);
\draw[very thick, color=purple](2,-13)--(2.30000000000000,-13);
\draw[very thick, color=cyan](3.70000000000000,-13)--(4,-13);
\draw[very thick, color=purple](2.30000000000000,-13)arc(270:360:0.700000000000000);
\draw[very thick, color=cyan](3.70000000000000,-13)arc(90:180:0.700000000000000);
\draw[very thick, color=cyan](2,-15)arc(270:360:1);
\draw[very thick, color=cyan](2,-15)arc(270:360:1);
\node[left, color=red] at (2,-3) {\small $1$};
\node[left, color=blue] at (2,-5) {\small $2$};
\node[left, color=teal] at (2,-7) {\small $3$};
\node[left, color=brown] at (2,-9) {\small $4$};
\node[left, color=orange] at (2,-11) {\small $5$};
\node[left, color=purple] at (2,-13) {\small $6$};
\node[left, color=cyan] at (2,-15) {\small $7$};
\node[above, color=red] at (3,-2) {\small $1$};
\node[above, color=blue] at (5,-2) {\small $2$};
\node[above, color=purple] at (7,-2) {\small $6$};
\node[above, color=brown] at (9,-2) {\small $4$};
\node[above, color=orange] at (11,-2) {\small $5$};
\node[above, color=teal] at (13,-2) {\small $3$};
\node[above, color=cyan] at (15,-2) {\small $7$};
\def\a{1.5};
\def\b{2};
\node[left, color=black] at (2-\a,-3) {\tiny $1$};
\node[left, color=black] at (2-\a,-5) {\tiny $2$};
\node[left, color=black] at (2-\a,-7) {\tiny $6$};
\node[left, color=black] at (2-\a,-9) {\tiny $4$};
\node[left, color=black] at (2-\a,-11) {\tiny $5$};
\node[left, color=black] at (2-\a,-13) {\tiny $3$};
\node[left, color=black] at (2-\a,-15) {\tiny $7$};
\node[above, color=black] at (3,-2+\b) {\tiny $1$};
\node[above, color=black] at (5,-2+\b) {\tiny $2$};
\node[above, color=black] at (7,-2+\b) {\tiny $3$};
\node[above, color=black] at (9,-2+\b) {\tiny $4$};
\node[above, color=black] at (11,-2+\b) {\tiny $5$};
\node[above, color=black] at (13,-2+\b) {\tiny $6$};
\node[above, color=black] at (15,-2+\b) {\tiny $7$};
\end{tikzpicture}
\begin{tikzpicture}[scale=0.160000000000000]
\draw[color=green](2,-2)--(16,-2);
\draw[color=green](2,-2)--(2,-16);
\draw[color=green](2,-4)--(16,-4);
\draw[color=green](4,-2)--(4,-16);
\draw[color=green](2,-6)--(14,-6);
\draw[color=green](6,-2)--(6,-14);
\draw[color=green](2,-8)--(12,-8);
\draw[color=green](8,-2)--(8,-12);
\draw[color=green](2,-10)--(10,-10);
\draw[color=green](10,-2)--(10,-10);
\draw[color=green](2,-12)--(8,-12);
\draw[color=green](12,-2)--(12,-8);
\draw[color=green](2,-14)--(6,-14);
\draw[color=green](14,-2)--(14,-6);
\draw[color=green](2,-16)--(4,-16);
\draw[color=green](16,-2)--(16,-4);
\draw[very thick, color=cyan](14,-3)arc(270:360:1);
\draw[very thick, color=teal](13,-2)--(13,-2.30000000000000);
\draw[very thick, color=cyan](13,-3.70000000000000)--(13,-4);
\draw[very thick, color=teal](12,-3)--(12.3000000000000,-3);
\draw[very thick, color=cyan](13.7000000000000,-3)--(14,-3);
\draw[very thick, color=teal](12.3000000000000,-3)arc(270:360:0.700000000000000);
\draw[very thick, color=cyan](13.7000000000000,-3)arc(90:180:0.700000000000000);
\draw[very thick, color=brown](11,-2)--(11,-2.30000000000000);
\draw[very thick, color=teal](11,-3.70000000000000)--(11,-4);
\draw[very thick, color=brown](10,-3)--(10.3000000000000,-3);
\draw[very thick, color=teal](11.7000000000000,-3)--(12,-3);
\draw[very thick, color=brown](10.3000000000000,-3)arc(270:360:0.700000000000000);
\draw[very thick, color=teal](11.7000000000000,-3)arc(90:180:0.700000000000000);
\draw[very thick, color=blue](9,-2)--(9,-2.30000000000000);
\draw[very thick, color=brown](9,-3.70000000000000)--(9,-4);
\draw[very thick, color=blue](8,-3)--(8.30000000000000,-3);
\draw[very thick, color=brown](9.70000000000000,-3)--(10,-3);
\draw[very thick, color=blue](8.30000000000000,-3)arc(270:360:0.700000000000000);
\draw[very thick, color=brown](9.70000000000000,-3)arc(90:180:0.700000000000000);
\draw[very thick, color=blue](6,-3)--(8,-3);
\draw[very thick, color=purple](7,-2)--(7,-4);
\draw[very thick, color=blue](4,-3)--(6,-3);
\draw[very thick, color=orange](5,-2)--(5,-4);
\draw[very thick, color=red](3,-2)--(3,-2.30000000000000);
\draw[very thick, color=blue](3,-3.70000000000000)--(3,-4);
\draw[very thick, color=red](2,-3)--(2.30000000000000,-3);
\draw[very thick, color=blue](3.70000000000000,-3)--(4,-3);
\draw[very thick, color=red](2.30000000000000,-3)arc(270:360:0.700000000000000);
\draw[very thick, color=blue](3.70000000000000,-3)arc(90:180:0.700000000000000);
\draw[very thick, color=cyan](12,-5)arc(270:360:1);
\draw[very thick, color=teal](11,-4)--(11,-4.30000000000000);
\draw[very thick, color=cyan](11,-5.70000000000000)--(11,-6);
\draw[very thick, color=teal](10,-5)--(10.3000000000000,-5);
\draw[very thick, color=cyan](11.7000000000000,-5)--(12,-5);
\draw[very thick, color=teal](10.3000000000000,-5)arc(270:360:0.700000000000000);
\draw[very thick, color=cyan](11.7000000000000,-5)arc(90:180:0.700000000000000);
\draw[very thick, color=brown](9,-4)--(9,-4.30000000000000);
\draw[very thick, color=teal](9,-5.70000000000000)--(9,-6);
\draw[very thick, color=brown](8,-5)--(8.30000000000000,-5);
\draw[very thick, color=teal](9.70000000000000,-5)--(10,-5);
\draw[very thick, color=brown](8.30000000000000,-5)arc(270:360:0.700000000000000);
\draw[very thick, color=teal](9.70000000000000,-5)arc(90:180:0.700000000000000);
\draw[very thick, color=brown](6,-5)--(8,-5);
\draw[very thick, color=purple](7,-4)--(7,-6);
\draw[very thick, color=orange](5,-4)--(5,-4.30000000000000);
\draw[very thick, color=brown](5,-5.70000000000000)--(5,-6);
\draw[very thick, color=orange](4,-5)--(4.30000000000000,-5);
\draw[very thick, color=brown](5.70000000000000,-5)--(6,-5);
\draw[very thick, color=orange](4.30000000000000,-5)arc(270:360:0.700000000000000);
\draw[very thick, color=brown](5.70000000000000,-5)arc(90:180:0.700000000000000);
\draw[very thick, color=blue](3,-4)--(3,-4.30000000000000);
\draw[very thick, color=orange](3,-5.70000000000000)--(3,-6);
\draw[very thick, color=blue](2,-5)--(2.30000000000000,-5);
\draw[very thick, color=orange](3.70000000000000,-5)--(4,-5);
\draw[very thick, color=blue](2.30000000000000,-5)arc(270:360:0.700000000000000);
\draw[very thick, color=orange](3.70000000000000,-5)arc(90:180:0.700000000000000);
\draw[very thick, color=cyan](10,-7)arc(270:360:1);
\draw[very thick, color=teal](9,-6)--(9,-6.30000000000000);
\draw[very thick, color=cyan](9,-7.70000000000000)--(9,-8);
\draw[very thick, color=teal](8,-7)--(8.30000000000000,-7);
\draw[very thick, color=cyan](9.70000000000000,-7)--(10,-7);
\draw[very thick, color=teal](8.30000000000000,-7)arc(270:360:0.700000000000000);
\draw[very thick, color=cyan](9.70000000000000,-7)arc(90:180:0.700000000000000);
\draw[very thick, color=teal](6,-7)--(8,-7);
\draw[very thick, color=purple](7,-6)--(7,-8);
\draw[very thick, color=teal](4,-7)--(6,-7);
\draw[very thick, color=brown](5,-6)--(5,-8);
\draw[very thick, color=teal](2,-7)--(4,-7);
\draw[very thick, color=orange](3,-6)--(3,-8);
\draw[very thick, color=cyan](8,-9)arc(270:360:1);
\draw[very thick, color=purple](7,-8)--(7,-8.30000000000000);
\draw[very thick, color=cyan](7,-9.70000000000000)--(7,-10);
\draw[very thick, color=purple](6,-9)--(6.30000000000000,-9);
\draw[very thick, color=cyan](7.70000000000000,-9)--(8,-9);
\draw[very thick, color=purple](6.30000000000000,-9)arc(270:360:0.700000000000000);
\draw[very thick, color=cyan](7.70000000000000,-9)arc(90:180:0.700000000000000);
\draw[very thick, color=brown](5,-8)--(5,-8.30000000000000);
\draw[very thick, color=purple](5,-9.70000000000000)--(5,-10);
\draw[very thick, color=brown](4,-9)--(4.30000000000000,-9);
\draw[very thick, color=purple](5.70000000000000,-9)--(6,-9);
\draw[very thick, color=brown](4.30000000000000,-9)arc(270:360:0.700000000000000);
\draw[very thick, color=purple](5.70000000000000,-9)arc(90:180:0.700000000000000);
\draw[very thick, color=brown](2,-9)--(4,-9);
\draw[very thick, color=orange](3,-8)--(3,-10);
\draw[very thick, color=cyan](6,-11)arc(270:360:1);
\draw[very thick, color=purple](5,-10)--(5,-10.3000000000000);
\draw[very thick, color=cyan](5,-11.7000000000000)--(5,-12);
\draw[very thick, color=purple](4,-11)--(4.30000000000000,-11);
\draw[very thick, color=cyan](5.70000000000000,-11)--(6,-11);
\draw[very thick, color=purple](4.30000000000000,-11)arc(270:360:0.700000000000000);
\draw[very thick, color=cyan](5.70000000000000,-11)arc(90:180:0.700000000000000);
\draw[very thick, color=orange](3,-10)--(3,-10.3000000000000);
\draw[very thick, color=purple](3,-11.7000000000000)--(3,-12);
\draw[very thick, color=orange](2,-11)--(2.30000000000000,-11);
\draw[very thick, color=purple](3.70000000000000,-11)--(4,-11);
\draw[very thick, color=orange](2.30000000000000,-11)arc(270:360:0.700000000000000);
\draw[very thick, color=purple](3.70000000000000,-11)arc(90:180:0.700000000000000);
\draw[very thick, color=cyan](4,-13)arc(270:360:1);
\draw[very thick, color=purple](3,-12)--(3,-12.3000000000000);
\draw[very thick, color=cyan](3,-13.7000000000000)--(3,-14);
\draw[very thick, color=purple](2,-13)--(2.30000000000000,-13);
\draw[very thick, color=cyan](3.70000000000000,-13)--(4,-13);
\draw[very thick, color=purple](2.30000000000000,-13)arc(270:360:0.700000000000000);
\draw[very thick, color=cyan](3.70000000000000,-13)arc(90:180:0.700000000000000);
\draw[very thick, color=cyan](2,-15)arc(270:360:1);
\node[left, color=red] at (2,-3) {\small $1$};
\node[left, color=blue] at (2,-5) {\small $2$};
\node[left, color=teal] at (2,-7) {\small $3$};
\node[left, color=brown] at (2,-9) {\small $4$};
\node[left, color=orange] at (2,-11) {\small $5$};
\node[left, color=purple] at (2,-13) {\small $6$};
\node[left, color=cyan] at (2,-15) {\small $7$};
\node[above, color=red] at (3,-2) {\small $1$};
\node[above, color=orange] at (5,-2) {\small $5$};
\node[above, color=purple] at (7,-2) {\small $6$};
\node[above, color=blue] at (9,-2) {\small $2$};
\node[above, color=brown] at (11,-2) {\small $4$};
\node[above, color=teal] at (13,-2) {\small $3$};
\node[above, color=cyan] at (15,-2) {\small $7$};
\def\a{1.5};
\def\b{2};
\node[left, color=black] at (2-\a,-3) {\tiny $1$};
\node[left, color=black] at (2-\a,-5) {\tiny $4$};
\node[left, color=black] at (2-\a,-7) {\tiny $6$};
\node[left, color=black] at (2-\a,-9) {\tiny $5$};
\node[left, color=black] at (2-\a,-11) {\tiny $2$};
\node[left, color=black] at (2-\a,-13) {\tiny $3$};
\node[left, color=black] at (2-\a,-15) {\tiny $7$};
\node[above, color=black] at (3,-2+\b) {\tiny $1$};
\node[above, color=black] at (5,-2+\b) {\tiny $2$};
\node[above, color=black] at (7,-2+\b) {\tiny $3$};
\node[above, color=black] at (9,-2+\b) {\tiny $4$};
\node[above, color=black] at (11,-2+\b) {\tiny $5$};
\node[above, color=black] at (13,-2+\b) {\tiny $6$};
\node[above, color=black] at (15,-2+\b) {\tiny $7$};
\end{tikzpicture}
\caption{If $D$ is
the pipe dream on the left with reduced word ${\bf r}_D = (4,3,5,6,4,3,5)$, then $\TransitionMap(D)$ is the
pipe dream on the right. In between we show the stack pushes in the
bounded bump algorithm.  The crossing initiating a stack push is
highlighted for each step. Here, $w=[1265734]$, hence $w^{-1} =
[1267435]$, $r=5, s=7$, $v = [1265437]$.  In this case,
$\TransitionMap(D)$ is a pipe dream for $v t_{25}=[1465237]$ so $q=2$.
\label{fig:pipe dream bump}}
\end{center}
\end{figure}

Given any pipe dream $D$, we can apply the bounded bump algorithm to
the corresponding bounded pair $(\mathbf{r}_D,\mathbf{j}_D)$ at any
position $t_{0}$ such that $\mathbf{r}_D$ is nearly reduced at $t_{0}$
and in either direction $\epsilon \in \{+,- \}$. One may observe that
the bounded pairs encountered during the steps of the bounded bump
algorithm do \emph{not} all encode pipe dreams, but it will turn out
that the departures from ``pipe dream encoding status'' are temporary,
and have a straightforward structure that will be analyzed in the
proof of \Cref{l:stack.push} below where we develop the notion of a
\textit{stack push}. It may be helpful to look ahead at \Cref{fig:pipe dream bump} to understand the proof of \Cref{l:stack.push}.

\begin{Lemma}\label{l:stack.push}
Let $D$ be a reduced pipe dream, and suppose that $\mathbf{r}_{D}$ is
nearly reduced at position $t$.  Let $\d\in \{+,- \}$ and suppose 
$$\Bump^{\d}_{t}(\mathbf{r}_{D}, \mathbf{j}_{D})=
(\bfa',\bfb', i, j, \outcome).$$ Then, the bounded pair
$(\bfa',\bfb')$ also encodes a reduced pipe dream.
\end{Lemma}

\begin{figure}
\begin{center}
\begin{tikzpicture}[scale=0.300000000000000] \def\a{24};
\filldraw[color=yellow](\a+6.5,-5)--(\a+8,-4)--(\a+8,-4.5)--(\a+13,-4.5)--(\a+13,-5.5)--(\a+8,-5.5)--(\a+8,-6);
\draw[color=green](2,-2)--(20,-2); \draw[color=green](2,-2)--(2,-10);
\draw[color=green](2,-4)--(20,-4); \draw[color=green](4,-2)--(4,-10);
\draw[color=green](2,-6)--(18,-6); \draw[color=green](6,-2)--(6,-10);
\draw[color=green](2,-8)--(16,-8); \draw[color=green](8,-2)--(8,-10);
\draw[color=green](2,-10)--(14,-10);
\draw[color=green](10,-2)--(10,-10);
\draw[color=green](2,-10)--(12,-10);
\draw[color=green](12,-2)--(12,-10);
\draw[color=green](2,-10)--(10,-10);
\draw[color=green](14,-2)--(14,-10);
\draw[color=green](2,-10)--(8,-10);
\draw[color=green](16,-2)--(16,-8);
\draw[color=green](2,-10)--(6,-10);
\draw[color=green](18,-2)--(18,-6);
\draw[color=green](2,-10)--(4,-10);
\draw[color=green](20,-2)--(20,-4); \draw[very thick](2,-3)--(4,-3);
\draw[very thick](3,-2)--(3,-4); \draw[very
thick](5,-2)--(5,-2.50000000000000); \draw[very
thick](5,-3.50000000000000)--(5,-4); \draw[very
thick](4,-3)--(4.50000000000000,-3); \draw[very
thick](5.50000000000000,-3)--(6,-3); \draw[very
thick](4.50000000000000,-3)arc(270:360:0.500000000000000); \draw[very
thick](5.50000000000000,-3)arc(90:180:0.500000000000000); \draw[very
thick](6,-3)--(8,-3); \draw[very thick](7,-2)--(7,-4); \draw[very
thick](8,-3)--(10,-3); \draw[very thick](9,-2)--(9,-4); \draw[very
thick](10,-3)--(12,-3); \draw[very thick](11,-2)--(11,-4); \draw[very
thick](12,-3)--(14,-3); \draw[very thick](13,-2)--(13,-4); \draw[very
thick](14,-3)--(16,-3); \draw[very thick](15,-2)--(15,-4); \draw[very
thick](17,-2)--(17,-2.50000000000000); \draw[very
thick](17,-3.50000000000000)--(17,-4); \draw[very
thick](16,-3)--(16.5000000000000,-3); \draw[very
thick](17.5000000000000,-3)--(18,-3); \draw[very
thick](16.5000000000000,-3)arc(270:360:0.500000000000000); \draw[very
thick](17.5000000000000,-3)arc(90:180:0.500000000000000); \draw[very
thick](3,-4)--(3,-4.50000000000000); \draw[very
thick](3,-5.50000000000000)--(3,-6); \draw[very
thick](2,-5)--(2.50000000000000,-5); \draw[very
thick](3.50000000000000,-5)--(4,-5); \draw[very
thick](2.50000000000000,-5)arc(270:360:0.500000000000000); \draw[very
thick](3.50000000000000,-5)arc(90:180:0.500000000000000); \draw[very
thick](5,-4)--(5,-4.50000000000000); \draw[very
thick](5,-5.50000000000000)--(5,-6); \draw[very
thick](4,-5)--(4.50000000000000,-5); \draw[very
thick](5.50000000000000,-5)--(6,-5); \draw[very
thick](4.50000000000000,-5)arc(270:360:0.500000000000000); \draw[very
thick](5.50000000000000,-5)arc(90:180:0.500000000000000); \draw[very
thick](7,-4)--(7,-4.50000000000000); \draw[very
thick](7,-5.50000000000000)--(7,-6); \draw[very
thick](6,-5)--(6.50000000000000,-5); \draw[very
thick](7.50000000000000,-5)--(8,-5); \draw[very
thick](6.50000000000000,-5)arc(270:360:0.500000000000000); \draw[very
thick](7.50000000000000,-5)arc(90:180:0.500000000000000); \draw[very
thick](9,-4)--(9,-4.50000000000000); \draw[very
thick](9,-5.50000000000000)--(9,-6); \draw[very
thick](8,-5)--(8.50000000000000,-5); \draw[very
thick](9.50000000000000,-5)--(10,-5); \draw[very
thick](8.50000000000000,-5)arc(270:360:0.500000000000000); \draw[very
thick](9.50000000000000,-5)arc(90:180:0.500000000000000); \draw[very
thick](10,-5)--(12,-5); \draw[very thick](11,-4)--(11,-6); \draw[very
thick](12,-5)--(14,-5); \draw[very thick](13,-4)--(13,-6); \draw[very
thick](14,-5)--(16,-5); \draw[very thick](15,-4)--(15,-6); \draw[very
thick](3,-6)--(3,-6.50000000000000); \draw[very
thick](3,-7.50000000000000)--(3,-8); \draw[very
thick](2,-7)--(2.50000000000000,-7); \draw[very
thick](3.50000000000000,-7)--(4,-7); \draw[very
thick](2.50000000000000,-7)arc(270:360:0.500000000000000); \draw[very
thick](3.50000000000000,-7)arc(90:180:0.500000000000000); \draw[very
thick](4,-7)--(6,-7); \draw[very thick](5,-6)--(5,-8); \draw[very
thick](6,-7)--(8,-7); \draw[very thick](7,-6)--(7,-8); \draw[very
thick](9,-6)--(9,-6.50000000000000); \draw[very
thick](9,-7.50000000000000)--(9,-8); \draw[very
thick](8,-7)--(8.50000000000000,-7); \draw[very
thick](9.50000000000000,-7)--(10,-7); \draw[very
thick](8.50000000000000,-7)arc(270:360:0.500000000000000); \draw[very
thick](9.50000000000000,-7)arc(90:180:0.500000000000000); \draw[very
thick](11,-6)--(11,-6.50000000000000); \draw[very
thick](11,-7.50000000000000)--(11,-8); \draw[very
thick](10,-7)--(10.5000000000000,-7); \draw[very
thick](11.5000000000000,-7)--(12,-7); \draw[very
thick](10.5000000000000,-7)arc(270:360:0.500000000000000); \draw[very
thick](11.5000000000000,-7)arc(90:180:0.500000000000000); \draw[very
thick](13,-6)--(13,-6.50000000000000); \draw[very
thick](13,-7.50000000000000)--(13,-8); \draw[very
thick](12,-7)--(12.5000000000000,-7); \draw[very
thick](13.5000000000000,-7)--(14,-7); \draw[very
thick](12.5000000000000,-7)arc(270:360:0.500000000000000); \draw[very
thick](13.5000000000000,-7)arc(90:180:0.500000000000000); \draw[very
thick](3,-8)--(3,-8.50000000000000); \draw[very
thick](3,-9.50000000000000)--(3,-10); \draw[very
thick](2,-9)--(2.50000000000000,-9); \draw[very
thick](3.50000000000000,-9)--(4,-9); \draw[very
thick](2.50000000000000,-9)arc(270:360:0.500000000000000); \draw[very
thick](3.50000000000000,-9)arc(90:180:0.500000000000000); \draw[very
thick](5,-8)--(5,-8.50000000000000); \draw[very
thick](5,-9.50000000000000)--(5,-10); \draw[very
thick](4,-9)--(4.50000000000000,-9); \draw[very
thick](5.50000000000000,-9)--(6,-9); \draw[very
thick](4.50000000000000,-9)arc(270:360:0.500000000000000); \draw[very
thick](5.50000000000000,-9)arc(90:180:0.500000000000000); \draw[very
thick](7,-8)--(7,-8.50000000000000); \draw[very
thick](7,-9.50000000000000)--(7,-10); \draw[very
thick](6,-9)--(6.50000000000000,-9); \draw[very
thick](7.50000000000000,-9)--(8,-9); \draw[very
thick](6.50000000000000,-9)arc(270:360:0.500000000000000); \draw[very
thick](7.50000000000000,-9)arc(90:180:0.500000000000000); \draw[very
thick](9,-8)--(9,-8.50000000000000); \draw[very
thick](9,-9.50000000000000)--(9,-10); \draw[very
thick](8,-9)--(8.50000000000000,-9); \draw[very
thick](9.50000000000000,-9)--(10,-9); \draw[very
thick](8.50000000000000,-9)arc(270:360:0.500000000000000); \draw[very
thick](9.50000000000000,-9)arc(90:180:0.500000000000000); \draw[very
thick](11,-8)--(11,-8.50000000000000); \draw[very
thick](11,-9.50000000000000)--(11,-10); \draw[very
thick](10,-9)--(10.5000000000000,-9); \draw[very
thick](11.5000000000000,-9)--(12,-9); \draw[very
thick](10.5000000000000,-9)arc(270:360:0.500000000000000); \draw[very
thick](11.5000000000000,-9)arc(90:180:0.500000000000000); \draw[very
thick](18,-3)arc(270:360:1); \draw[very thick](16,-5)arc(270:360:1);
\draw[very thick](14,-7)arc(270:360:1); \draw[very
thick](12,-9)arc(270:360:1); \node at (0,-6) {$D$}; \draw[very thick,
color=red](15,-4.5)arc(90:450:0.5);

\node at (\a+0,-3) {$+$};
\node at (\a+2,-3) {$\cdot$};
\node at (\a+4,-3) {$+$};
\node at (\a+6,-3) {$+$};
\node at (\a+8,-3) {$+$};
\node at (\a+10,-3) {$+$};
\node at (\a+12,-3) {$+$};
\node at (\a+0,-5) {$\cdot$};
\node at (\a+2,-5) {$\cdot$};
\node at (\a+4,-5) {$\cdot$};
\node at (\a+6,-5) {$\cdot$};
\node at (\a+8,-5) {$+$};
\node at (\a+10,-5) {$+$};
\node at (\a+12,-5) {$+$};
\node at (\a+0,-7) {$\cdot$};
\node at (\a+2,-7) {$+$};
\node at (\a+4,-7) {$+$};
\node at (\a+6,-7) {$\cdot$};
\node at (\a+8,-7) {$\cdot$};
\node at (\a+10,-7) {$\cdot$};
\node at (\a+12,-7) {$\cdot$};
\end{tikzpicture}

\

\begin{tikzpicture}[scale=0.300000000000000]
\def\a{24};
\filldraw[color=yellow](\a+0.5,-7)--(\a+2,-6)--(\a+2,-6.5)--(\a+5,-6.5)--(\a+5,-7.5)--(\a+2,-7.5)--(\a+2,-8);

\draw[color=green](2,-2)--(20,-2);
\draw[color=green](2,-2)--(2,-10);
\draw[color=green](2,-4)--(20,-4);
\draw[color=green](4,-2)--(4,-10);
\draw[color=green](2,-6)--(18,-6);
\draw[color=green](6,-2)--(6,-10);
\draw[color=green](2,-8)--(16,-8);
\draw[color=green](8,-2)--(8,-10);
\draw[color=green](2,-10)--(14,-10);
\draw[color=green](10,-2)--(10,-10);
\draw[color=green](2,-10)--(12,-10);
\draw[color=green](12,-2)--(12,-10);
\draw[color=green](2,-10)--(10,-10);
\draw[color=green](14,-2)--(14,-10);
\draw[color=green](2,-10)--(8,-10);
\draw[color=green](16,-2)--(16,-8);
\draw[color=green](2,-10)--(6,-10);
\draw[color=green](18,-2)--(18,-6);
\draw[color=green](2,-10)--(4,-10);
\draw[color=green](20,-2)--(20,-4);
\draw[very thick](2,-3)--(4,-3);
\draw[very thick](3,-2)--(3,-4);
\draw[very thick](5,-2)--(5,-2.50000000000000);
\draw[very thick](5,-3.50000000000000)--(5,-4);
\draw[very thick](4,-3)--(4.50000000000000,-3);
\draw[very thick](5.50000000000000,-3)--(6,-3);
\draw[very thick](4.50000000000000,-3)arc(270:360:0.500000000000000);
\draw[very thick](5.50000000000000,-3)arc(90:180:0.500000000000000);
\draw[very thick](6,-3)--(8,-3);
\draw[very thick](7,-2)--(7,-4);
\draw[very thick](8,-3)--(10,-3);
\draw[very thick](9,-2)--(9,-4);
\draw[very thick](10,-3)--(12,-3);
\draw[very thick](11,-2)--(11,-4);
\draw[very thick](12,-3)--(14,-3);
\draw[very thick](13,-2)--(13,-4);
\draw[very thick](14,-3)--(16,-3);
\draw[very thick](15,-2)--(15,-4);
\draw[very thick](17,-2)--(17,-2.50000000000000);
\draw[very thick](17,-3.50000000000000)--(17,-4);
\draw[very thick](16,-3)--(16.5000000000000,-3);
\draw[very thick](17.5000000000000,-3)--(18,-3);
\draw[very thick](16.5000000000000,-3)arc(270:360:0.500000000000000);
\draw[very thick](17.5000000000000,-3)arc(90:180:0.500000000000000);
\draw[very thick](3,-4)--(3,-4.50000000000000);
\draw[very thick](3,-5.50000000000000)--(3,-6);
\draw[very thick](2,-5)--(2.50000000000000,-5);
\draw[very thick](3.50000000000000,-5)--(4,-5);
\draw[very thick](2.50000000000000,-5)arc(270:360:0.500000000000000);
\draw[very thick](3.50000000000000,-5)arc(90:180:0.500000000000000);
\draw[very thick](5,-4)--(5,-4.50000000000000);
\draw[very thick](5,-5.50000000000000)--(5,-6);
\draw[very thick](4,-5)--(4.50000000000000,-5);
\draw[very thick](5.50000000000000,-5)--(6,-5);
\draw[very thick](4.50000000000000,-5)arc(270:360:0.500000000000000);
\draw[very thick](5.50000000000000,-5)arc(90:180:0.500000000000000);
\draw[very thick](7,-4)--(7,-4.50000000000000);
\draw[very thick](7,-5.50000000000000)--(7,-6);
\draw[very thick](6,-5)--(6.50000000000000,-5);
\draw[very thick](7.50000000000000,-5)--(8,-5);
\draw[very thick](6.50000000000000,-5)arc(270:360:0.500000000000000);
\draw[very thick](7.50000000000000,-5)arc(90:180:0.500000000000000);
\draw[very thick](8,-5)--(10,-5);
\draw[very thick](9,-4)--(9,-6);
\draw[very thick](10,-5)--(12,-5);
\draw[very thick](11,-4)--(11,-6);
\draw[very thick](12,-5)--(14,-5);
\draw[very thick](13,-4)--(13,-6);
\draw[very thick](15,-4)--(15,-4.50000000000000);
\draw[very thick](15,-5.50000000000000)--(15,-6);
\draw[very thick](14,-5)--(14.5000000000000,-5);
\draw[very thick](15.5000000000000,-5)--(16,-5);
\draw[very thick](14.5000000000000,-5)arc(270:360:0.500000000000000);
\draw[very thick](15.5000000000000,-5)arc(90:180:0.500000000000000);
\draw[very thick](3,-6)--(3,-6.50000000000000);
\draw[very thick](3,-7.50000000000000)--(3,-8);
\draw[very thick](2,-7)--(2.50000000000000,-7);
\draw[very thick](3.50000000000000,-7)--(4,-7);
\draw[very thick](2.50000000000000,-7)arc(270:360:0.500000000000000);
\draw[very thick](3.50000000000000,-7)arc(90:180:0.500000000000000);
\draw[very thick](4,-7)--(6,-7);
\draw[very thick](5,-6)--(5,-8);
\draw[very thick](6,-7)--(8,-7);
\draw[very thick](7,-6)--(7,-8);
\draw[very thick](9,-6)--(9,-6.50000000000000);
\draw[very thick](9,-7.50000000000000)--(9,-8);
\draw[very thick](8,-7)--(8.50000000000000,-7);
\draw[very thick](9.50000000000000,-7)--(10,-7);
\draw[very thick](8.50000000000000,-7)arc(270:360:0.500000000000000);
\draw[very thick](9.50000000000000,-7)arc(90:180:0.500000000000000);
\draw[very thick](11,-6)--(11,-6.50000000000000);
\draw[very thick](11,-7.50000000000000)--(11,-8);
\draw[very thick](10,-7)--(10.5000000000000,-7);
\draw[very thick](11.5000000000000,-7)--(12,-7);
\draw[very thick](10.5000000000000,-7)arc(270:360:0.500000000000000);
\draw[very thick](11.5000000000000,-7)arc(90:180:0.500000000000000);
\draw[very thick](13,-6)--(13,-6.50000000000000);
\draw[very thick](13,-7.50000000000000)--(13,-8);
\draw[very thick](12,-7)--(12.5000000000000,-7);
\draw[very thick](13.5000000000000,-7)--(14,-7);
\draw[very thick](12.5000000000000,-7)arc(270:360:0.500000000000000);
\draw[very thick](13.5000000000000,-7)arc(90:180:0.500000000000000);
\draw[very thick](3,-8)--(3,-8.50000000000000);
\draw[very thick](3,-9.50000000000000)--(3,-10);
\draw[very thick](2,-9)--(2.50000000000000,-9);
\draw[very thick](3.50000000000000,-9)--(4,-9);
\draw[very thick](2.50000000000000,-9)arc(270:360:0.500000000000000);
\draw[very thick](3.50000000000000,-9)arc(90:180:0.500000000000000);
\draw[very thick](5,-8)--(5,-8.50000000000000);
\draw[very thick](5,-9.50000000000000)--(5,-10);
\draw[very thick](4,-9)--(4.50000000000000,-9);
\draw[very thick](5.50000000000000,-9)--(6,-9);
\draw[very thick](4.50000000000000,-9)arc(270:360:0.500000000000000);
\draw[very thick](5.50000000000000,-9)arc(90:180:0.500000000000000);
\draw[very thick](7,-8)--(7,-8.50000000000000);
\draw[very thick](7,-9.50000000000000)--(7,-10);
\draw[very thick](6,-9)--(6.50000000000000,-9);
\draw[very thick](7.50000000000000,-9)--(8,-9);
\draw[very thick](6.50000000000000,-9)arc(270:360:0.500000000000000);
\draw[very thick](7.50000000000000,-9)arc(90:180:0.500000000000000);
\draw[very thick](9,-8)--(9,-8.50000000000000);
\draw[very thick](9,-9.50000000000000)--(9,-10);
\draw[very thick](8,-9)--(8.50000000000000,-9);
\draw[very thick](9.50000000000000,-9)--(10,-9);
\draw[very thick](8.50000000000000,-9)arc(270:360:0.500000000000000);
\draw[very thick](9.50000000000000,-9)arc(90:180:0.500000000000000);
\draw[very thick](11,-8)--(11,-8.50000000000000);
\draw[very thick](11,-9.50000000000000)--(11,-10);
\draw[very thick](10,-9)--(10.5000000000000,-9);
\draw[very thick](11.5000000000000,-9)--(12,-9);
\draw[very thick](10.5000000000000,-9)arc(270:360:0.500000000000000);
\draw[very thick](11.5000000000000,-9)arc(90:180:0.500000000000000);
\draw[very thick](18,-3)arc(270:360:1);
\draw[very thick](16,-5)arc(270:360:1);
\draw[very thick](14,-7)arc(270:360:1);
\draw[very thick](12,-9)arc(270:360:1);
\node at (0,-6) {$D'$};
\draw[very thick, color=red](9,-4.5)arc(90:450:0.5);
\draw[very thick, color=red](7,-6.5)arc(90:450:0.5);

\node at (\a+0,-3) {$+$};
\node at (\a+2,-3) {$\cdot$};
\node at (\a+4,-3) {$+$};
\node at (\a+6,-3) {$+$};
\node at (\a+8,-3) {$+$};
\node at (\a+10,-3) {$+$};
\node at (\a+12,-3) {$+$};
\node at (\a+0,-5) {$\cdot$};
\node at (\a+2,-5) {$\cdot$};
\node at (\a+4,-5) {$\cdot$};
\node at (\a+6,-5) {$+$};
\node at (\a+8,-5) {$+$};
\node at (\a+10,-5) {$+$};
\node at (\a+12,-5) {$\cdot$};
\node at (\a+0,-7) {$\cdot$};
\node at (\a+2,-7) {$+$};
\node at (\a+4,-7) {$+$};
\node at (\a+6,-7) {$\cdot$};
\node at (\a+8,-7) {$\cdot$};
\node at (\a+10,-7) {$\cdot$};
\node at (\a+12,-7) {$\cdot$};
\end{tikzpicture}

\

\begin{tikzpicture}[scale=0.300000000000000]
\def\a{24};
\filldraw[color=yellow](\a+2.5,-3)--(\a+4,-2)--(\a+4,-2.5)--(\a+13,-2.5)--(\a+13,-3.5)--(\a+4,-3.5)--(\a+4,-4);

\draw[color=green](2,-2)--(20,-2);
\draw[color=green](2,-2)--(2,-10);
\draw[color=green](2,-4)--(20,-4);
\draw[color=green](4,-2)--(4,-10);
\draw[color=green](2,-6)--(18,-6);
\draw[color=green](6,-2)--(6,-10);
\draw[color=green](2,-8)--(16,-8);
\draw[color=green](8,-2)--(8,-10);
\draw[color=green](2,-10)--(14,-10);
\draw[color=green](10,-2)--(10,-10);
\draw[color=green](2,-10)--(12,-10);
\draw[color=green](12,-2)--(12,-10);
\draw[color=green](2,-10)--(10,-10);
\draw[color=green](14,-2)--(14,-10);
\draw[color=green](2,-10)--(8,-10);
\draw[color=green](16,-2)--(16,-8);
\draw[color=green](2,-10)--(6,-10);
\draw[color=green](18,-2)--(18,-6);
\draw[color=green](2,-10)--(4,-10);
\draw[color=green](20,-2)--(20,-4);
\draw[very thick](2,-3)--(4,-3);
\draw[very thick](3,-2)--(3,-4);
\draw[very thick](5,-2)--(5,-2.50000000000000);
\draw[very thick](5,-3.50000000000000)--(5,-4);
\draw[very thick](4,-3)--(4.50000000000000,-3);
\draw[very thick](5.50000000000000,-3)--(6,-3);
\draw[very thick](4.50000000000000,-3)arc(270:360:0.500000000000000);
\draw[very thick](5.50000000000000,-3)arc(90:180:0.500000000000000);
\draw[very thick](6,-3)--(8,-3);
\draw[very thick](7,-2)--(7,-4);
\draw[very thick](8,-3)--(10,-3);
\draw[very thick](9,-2)--(9,-4);
\draw[very thick](10,-3)--(12,-3);
\draw[very thick](11,-2)--(11,-4);
\draw[very thick](12,-3)--(14,-3);
\draw[very thick](13,-2)--(13,-4);
\draw[very thick](14,-3)--(16,-3);
\draw[very thick](15,-2)--(15,-4);
\draw[very thick](17,-2)--(17,-2.50000000000000);
\draw[very thick](17,-3.50000000000000)--(17,-4);
\draw[very thick](16,-3)--(16.5000000000000,-3);
\draw[very thick](17.5000000000000,-3)--(18,-3);
\draw[very thick](16.5000000000000,-3)arc(270:360:0.500000000000000);
\draw[very thick](17.5000000000000,-3)arc(90:180:0.500000000000000);
\draw[very thick](3,-4)--(3,-4.50000000000000);
\draw[very thick](3,-5.50000000000000)--(3,-6);
\draw[very thick](2,-5)--(2.50000000000000,-5);
\draw[very thick](3.50000000000000,-5)--(4,-5);
\draw[very thick](2.50000000000000,-5)arc(270:360:0.500000000000000);
\draw[very thick](3.50000000000000,-5)arc(90:180:0.500000000000000);
\draw[very thick](5,-4)--(5,-4.50000000000000);
\draw[very thick](5,-5.50000000000000)--(5,-6);
\draw[very thick](4,-5)--(4.50000000000000,-5);
\draw[very thick](5.50000000000000,-5)--(6,-5);
\draw[very thick](4.50000000000000,-5)arc(270:360:0.500000000000000);
\draw[very thick](5.50000000000000,-5)arc(90:180:0.500000000000000);
\draw[very thick](7,-4)--(7,-4.50000000000000);
\draw[very thick](7,-5.50000000000000)--(7,-6);
\draw[very thick](6,-5)--(6.50000000000000,-5);
\draw[very thick](7.50000000000000,-5)--(8,-5);
\draw[very thick](6.50000000000000,-5)arc(270:360:0.500000000000000);
\draw[very thick](7.50000000000000,-5)arc(90:180:0.500000000000000);
\draw[very thick](8,-5)--(10,-5);
\draw[very thick](9,-4)--(9,-6);
\draw[very thick](10,-5)--(12,-5);
\draw[very thick](11,-4)--(11,-6);
\draw[very thick](12,-5)--(14,-5);
\draw[very thick](13,-4)--(13,-6);
\draw[very thick](15,-4)--(15,-4.50000000000000);
\draw[very thick](15,-5.50000000000000)--(15,-6);
\draw[very thick](14,-5)--(14.5000000000000,-5);
\draw[very thick](15.5000000000000,-5)--(16,-5);
\draw[very thick](14.5000000000000,-5)arc(270:360:0.500000000000000);
\draw[very thick](15.5000000000000,-5)arc(90:180:0.500000000000000);
\draw[very thick](2,-7)--(4,-7);
\draw[very thick](3,-6)--(3,-8);
\draw[very thick](4,-7)--(6,-7);
\draw[very thick](5,-6)--(5,-8);
\draw[very thick](7,-6)--(7,-6.50000000000000);
\draw[very thick](7,-7.50000000000000)--(7,-8);
\draw[very thick](6,-7)--(6.50000000000000,-7);
\draw[very thick](7.50000000000000,-7)--(8,-7);
\draw[very thick](6.50000000000000,-7)arc(270:360:0.500000000000000);
\draw[very thick](7.50000000000000,-7)arc(90:180:0.500000000000000);
\draw[very thick](9,-6)--(9,-6.50000000000000);
\draw[very thick](9,-7.50000000000000)--(9,-8);
\draw[very thick](8,-7)--(8.50000000000000,-7);
\draw[very thick](9.50000000000000,-7)--(10,-7);
\draw[very thick](8.50000000000000,-7)arc(270:360:0.500000000000000);
\draw[very thick](9.50000000000000,-7)arc(90:180:0.500000000000000);
\draw[very thick](11,-6)--(11,-6.50000000000000);
\draw[very thick](11,-7.50000000000000)--(11,-8);
\draw[very thick](10,-7)--(10.5000000000000,-7);
\draw[very thick](11.5000000000000,-7)--(12,-7);
\draw[very thick](10.5000000000000,-7)arc(270:360:0.500000000000000);
\draw[very thick](11.5000000000000,-7)arc(90:180:0.500000000000000);
\draw[very thick](13,-6)--(13,-6.50000000000000);
\draw[very thick](13,-7.50000000000000)--(13,-8);
\draw[very thick](12,-7)--(12.5000000000000,-7);
\draw[very thick](13.5000000000000,-7)--(14,-7);
\draw[very thick](12.5000000000000,-7)arc(270:360:0.500000000000000);
\draw[very thick](13.5000000000000,-7)arc(90:180:0.500000000000000);
\draw[very thick](3,-8)--(3,-8.50000000000000);
\draw[very thick](3,-9.50000000000000)--(3,-10);
\draw[very thick](2,-9)--(2.50000000000000,-9);
\draw[very thick](3.50000000000000,-9)--(4,-9);
\draw[very thick](2.50000000000000,-9)arc(270:360:0.500000000000000);
\draw[very thick](3.50000000000000,-9)arc(90:180:0.500000000000000);
\draw[very thick](5,-8)--(5,-8.50000000000000);
\draw[very thick](5,-9.50000000000000)--(5,-10);
\draw[very thick](4,-9)--(4.50000000000000,-9);
\draw[very thick](5.50000000000000,-9)--(6,-9);
\draw[very thick](4.50000000000000,-9)arc(270:360:0.500000000000000);
\draw[very thick](5.50000000000000,-9)arc(90:180:0.500000000000000);
\draw[very thick](7,-8)--(7,-8.50000000000000);
\draw[very thick](7,-9.50000000000000)--(7,-10);
\draw[very thick](6,-9)--(6.50000000000000,-9);
\draw[very thick](7.50000000000000,-9)--(8,-9);
\draw[very thick](6.50000000000000,-9)arc(270:360:0.500000000000000);
\draw[very thick](7.50000000000000,-9)arc(90:180:0.500000000000000);
\draw[very thick](9,-8)--(9,-8.50000000000000);
\draw[very thick](9,-9.50000000000000)--(9,-10);
\draw[very thick](8,-9)--(8.50000000000000,-9);
\draw[very thick](9.50000000000000,-9)--(10,-9);
\draw[very thick](8.50000000000000,-9)arc(270:360:0.500000000000000);
\draw[very thick](9.50000000000000,-9)arc(90:180:0.500000000000000);
\draw[very thick](11,-8)--(11,-8.50000000000000);
\draw[very thick](11,-9.50000000000000)--(11,-10);
\draw[very thick](10,-9)--(10.5000000000000,-9);
\draw[very thick](11.5000000000000,-9)--(12,-9);
\draw[very thick](10.5000000000000,-9)arc(270:360:0.500000000000000);
\draw[very thick](11.5000000000000,-9)arc(90:180:0.500000000000000);
\draw[very thick](18,-3)arc(270:360:1);
\draw[very thick](16,-5)arc(270:360:1);
\draw[very thick](14,-7)arc(270:360:1);
\draw[very thick](12,-9)arc(270:360:1);
\node at (0,-6) {$D''$};
\draw[very thick, color=red](15,-2.5)arc(90:450:0.5);
\draw[very thick, color=red](3,-6.5)arc(90:450:0.5);

\node at (\a+0,-3) {$+$};
\node at (\a+2,-3) {$\cdot$};
\node at (\a+4,-3) {$+$};
\node at (\a+6,-3) {$+$};
\node at (\a+8,-3) {$+$};
\node at (\a+10,-3) {$+$};
\node at (\a+12,-3) {$+$};
\node at (\a+0,-5) {$\cdot$};
\node at (\a+2,-5) {$\cdot$};
\node at (\a+4,-5) {$\cdot$};
\node at (\a+6,-5) {$+$};
\node at (\a+8,-5) {$+$};
\node at (\a+10,-5) {$+$};
\node at (\a+12,-5) {$\cdot$};
\node at (\a+0,-7) {$+$};
\node at (\a+2,-7) {$+$};
\node at (\a+4,-7) {$\cdot$};
\node at (\a+6,-7) {$\cdot$};
\node at (\a+8,-7) {$\cdot$};
\node at (\a+10,-7) {$\cdot$};
\node at (\a+12,-7) {$\cdot$};
\end{tikzpicture}

\

\begin{tikzpicture}[scale=0.300000000000000]
\def\a{24};

\draw[color=green](2,-2)--(20,-2);
\draw[color=green](2,-2)--(2,-10);
\draw[color=green](2,-4)--(20,-4);
\draw[color=green](4,-2)--(4,-10);
\draw[color=green](2,-6)--(18,-6);
\draw[color=green](6,-2)--(6,-10);
\draw[color=green](2,-8)--(16,-8);
\draw[color=green](8,-2)--(8,-10);
\draw[color=green](2,-10)--(14,-10);
\draw[color=green](10,-2)--(10,-10);
\draw[color=green](2,-10)--(12,-10);
\draw[color=green](12,-2)--(12,-10);
\draw[color=green](2,-10)--(10,-10);
\draw[color=green](14,-2)--(14,-10);
\draw[color=green](2,-10)--(8,-10);
\draw[color=green](16,-2)--(16,-8);
\draw[color=green](2,-10)--(6,-10);
\draw[color=green](18,-2)--(18,-6);
\draw[color=green](2,-10)--(4,-10);
\draw[color=green](20,-2)--(20,-4);
\draw[very thick](2,-3)--(4,-3);
\draw[very thick](3,-2)--(3,-4);
\draw[very thick](4,-3)--(6,-3);
\draw[very thick](5,-2)--(5,-4);
\draw[very thick](6,-3)--(8,-3);
\draw[very thick](7,-2)--(7,-4);
\draw[very thick](8,-3)--(10,-3);
\draw[very thick](9,-2)--(9,-4);
\draw[very thick](10,-3)--(12,-3);
\draw[very thick](11,-2)--(11,-4);
\draw[very thick](12,-3)--(14,-3);
\draw[very thick](13,-2)--(13,-4);
\draw[very thick](15,-2)--(15,-2.50000000000000);
\draw[very thick](15,-3.50000000000000)--(15,-4);
\draw[very thick](14,-3)--(14.5000000000000,-3);
\draw[very thick](15.5000000000000,-3)--(16,-3);
\draw[very thick](14.5000000000000,-3)arc(270:360:0.500000000000000);
\draw[very thick](15.5000000000000,-3)arc(90:180:0.500000000000000);
\draw[very thick](17,-2)--(17,-2.50000000000000);
\draw[very thick](17,-3.50000000000000)--(17,-4);
\draw[very thick](16,-3)--(16.5000000000000,-3);
\draw[very thick](17.5000000000000,-3)--(18,-3);
\draw[very thick](16.5000000000000,-3)arc(270:360:0.500000000000000);
\draw[very thick](17.5000000000000,-3)arc(90:180:0.500000000000000);
\draw[very thick](3,-4)--(3,-4.50000000000000);
\draw[very thick](3,-5.50000000000000)--(3,-6);
\draw[very thick](2,-5)--(2.50000000000000,-5);
\draw[very thick](3.50000000000000,-5)--(4,-5);
\draw[very thick](2.50000000000000,-5)arc(270:360:0.500000000000000);
\draw[very thick](3.50000000000000,-5)arc(90:180:0.500000000000000);
\draw[very thick](5,-4)--(5,-4.50000000000000);
\draw[very thick](5,-5.50000000000000)--(5,-6);
\draw[very thick](4,-5)--(4.50000000000000,-5);
\draw[very thick](5.50000000000000,-5)--(6,-5);
\draw[very thick](4.50000000000000,-5)arc(270:360:0.500000000000000);
\draw[very thick](5.50000000000000,-5)arc(90:180:0.500000000000000);
\draw[very thick](7,-4)--(7,-4.50000000000000);
\draw[very thick](7,-5.50000000000000)--(7,-6);
\draw[very thick](6,-5)--(6.50000000000000,-5);
\draw[very thick](7.50000000000000,-5)--(8,-5);
\draw[very thick](6.50000000000000,-5)arc(270:360:0.500000000000000);
\draw[very thick](7.50000000000000,-5)arc(90:180:0.500000000000000);
\draw[very thick](8,-5)--(10,-5);
\draw[very thick](9,-4)--(9,-6);
\draw[very thick](10,-5)--(12,-5);
\draw[very thick](11,-4)--(11,-6);
\draw[very thick](12,-5)--(14,-5);
\draw[very thick](13,-4)--(13,-6);
\draw[very thick](15,-4)--(15,-4.50000000000000);
\draw[very thick](15,-5.50000000000000)--(15,-6);
\draw[very thick](14,-5)--(14.5000000000000,-5);
\draw[very thick](15.5000000000000,-5)--(16,-5);
\draw[very thick](14.5000000000000,-5)arc(270:360:0.500000000000000);
\draw[very thick](15.5000000000000,-5)arc(90:180:0.500000000000000);
\draw[very thick](2,-7)--(4,-7);
\draw[very thick](3,-6)--(3,-8);
\draw[very thick](4,-7)--(6,-7);
\draw[very thick](5,-6)--(5,-8);
\draw[very thick](7,-6)--(7,-6.50000000000000);
\draw[very thick](7,-7.50000000000000)--(7,-8);
\draw[very thick](6,-7)--(6.50000000000000,-7);
\draw[very thick](7.50000000000000,-7)--(8,-7);
\draw[very thick](6.50000000000000,-7)arc(270:360:0.500000000000000);
\draw[very thick](7.50000000000000,-7)arc(90:180:0.500000000000000);
\draw[very thick](9,-6)--(9,-6.50000000000000);
\draw[very thick](9,-7.50000000000000)--(9,-8);
\draw[very thick](8,-7)--(8.50000000000000,-7);
\draw[very thick](9.50000000000000,-7)--(10,-7);
\draw[very thick](8.50000000000000,-7)arc(270:360:0.500000000000000);
\draw[very thick](9.50000000000000,-7)arc(90:180:0.500000000000000);
\draw[very thick](11,-6)--(11,-6.50000000000000);
\draw[very thick](11,-7.50000000000000)--(11,-8);
\draw[very thick](10,-7)--(10.5000000000000,-7);
\draw[very thick](11.5000000000000,-7)--(12,-7);
\draw[very thick](10.5000000000000,-7)arc(270:360:0.500000000000000);
\draw[very thick](11.5000000000000,-7)arc(90:180:0.500000000000000);
\draw[very thick](13,-6)--(13,-6.50000000000000);
\draw[very thick](13,-7.50000000000000)--(13,-8);
\draw[very thick](12,-7)--(12.5000000000000,-7);
\draw[very thick](13.5000000000000,-7)--(14,-7);
\draw[very thick](12.5000000000000,-7)arc(270:360:0.500000000000000);
\draw[very thick](13.5000000000000,-7)arc(90:180:0.500000000000000);
\draw[very thick](3,-8)--(3,-8.50000000000000);
\draw[very thick](3,-9.50000000000000)--(3,-10);
\draw[very thick](2,-9)--(2.50000000000000,-9);
\draw[very thick](3.50000000000000,-9)--(4,-9);
\draw[very thick](2.50000000000000,-9)arc(270:360:0.500000000000000);
\draw[very thick](3.50000000000000,-9)arc(90:180:0.500000000000000);
\draw[very thick](5,-8)--(5,-8.50000000000000);
\draw[very thick](5,-9.50000000000000)--(5,-10);
\draw[very thick](4,-9)--(4.50000000000000,-9);
\draw[very thick](5.50000000000000,-9)--(6,-9);
\draw[very thick](4.50000000000000,-9)arc(270:360:0.500000000000000);
\draw[very thick](5.50000000000000,-9)arc(90:180:0.500000000000000);
\draw[very thick](7,-8)--(7,-8.50000000000000);
\draw[very thick](7,-9.50000000000000)--(7,-10);
\draw[very thick](6,-9)--(6.50000000000000,-9);
\draw[very thick](7.50000000000000,-9)--(8,-9);
\draw[very thick](6.50000000000000,-9)arc(270:360:0.500000000000000);
\draw[very thick](7.50000000000000,-9)arc(90:180:0.500000000000000);
\draw[very thick](9,-8)--(9,-8.50000000000000);
\draw[very thick](9,-9.50000000000000)--(9,-10);
\draw[very thick](8,-9)--(8.50000000000000,-9);
\draw[very thick](9.50000000000000,-9)--(10,-9);
\draw[very thick](8.50000000000000,-9)arc(270:360:0.500000000000000);
\draw[very thick](9.50000000000000,-9)arc(90:180:0.500000000000000);
\draw[very thick](11,-8)--(11,-8.50000000000000);
\draw[very thick](11,-9.50000000000000)--(11,-10);
\draw[very thick](10,-9)--(10.5000000000000,-9);
\draw[very thick](11.5000000000000,-9)--(12,-9);
\draw[very thick](10.5000000000000,-9)arc(270:360:0.500000000000000);
\draw[very thick](11.5000000000000,-9)arc(90:180:0.500000000000000);
\draw[very thick](18,-3)arc(270:360:1);
\draw[very thick](16,-5)arc(270:360:1);
\draw[very thick](14,-7)arc(270:360:1);
\draw[very thick](12,-9)arc(270:360:1);
\node at (0,-6) {$D'''$};
\node at (\a+0,-3) {$+$};
\node at (\a+2,-3) {$+$};
\node at (\a+4,-3) {$+$};
\node at (\a+6,-3) {$+$};
\node at (\a+8,-3) {$+$};
\node at (\a+10,-3) {$+$};
\node at (\a+12,-3) {$\cdot$};
\node at (\a+0,-5) {$\cdot$};
\node at (\a+2,-5) {$\cdot$};
\node at (\a+4,-5) {$\cdot$};
\node at (\a+6,-5) {$+$};
\node at (\a+8,-5) {$+$};
\node at (\a+10,-5) {$+$};
\node at (\a+12,-5) {$\cdot$};
\node at (\a+0,-7) {$+$};
\node at (\a+2,-7) {$+$};
\node at (\a+4,-7) {$\cdot$};
\node at (\a+6,-7) {$\cdot$};
\node at (\a+8,-7) {$\cdot$};
\node at (\a+10,-7) {$\cdot$};
\node at (\a+12,-7) {$\cdot$};
\end{tikzpicture}
\end{center}
\caption{Computing $\Bump^{-}_{6}(\mathbf{r}_{D},j_{D})$ for $D$ starting in position $(2,7)$ with
$r_{D}=(7,6,5,4,3,1,8,7,6,5,4)$ and $i_{D}=(1,1,1,1,1,1,2,2,2,3,3)$
via 3 stack pushes. } \label{fig:stack.push}
\end{figure}
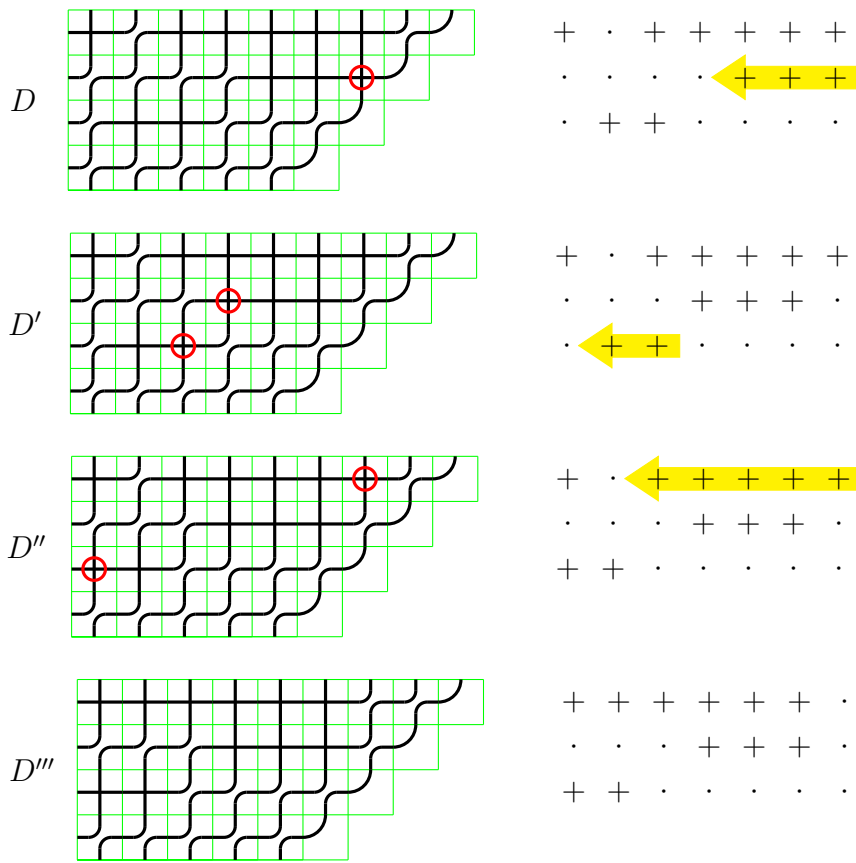

\begin{proof}
Consider the effect of the bounded bumping algorithm in terms of pipe
dreams.  To be concrete, assume $\d=-$; the case $\d=+$ being similar.
Recalling Definition~\ref{defn:nearly-reduced}, $\mathbf{r}_D$ is nearly reduced at position $t$ if and only
if removing the $t$\textsuperscript{th} crossing of $D$ in reading order results in a reduced pipe
dream.  Observe that when we initially decrement-push
$(\mathbf{r}_{D},\mathbf{j}_{D}) $ in column $t$, it has the effect of
moving the $t^{th}$ crossing in the reading order on $D$, say in
position $(i,j)\in D$, one column to the left to position
$(i,j-1)$. If this location is already occupied, that is if $(i,j-1) \in D$, then
$\Push_{t}^{-} \mathbf{r}_{D}$ returns a nearly reduced word with letter $i+j-2$ in both positions $t$ and $t+1$.  The resulting bounded pair does
not encode a pipe dream.  Then, the next step of the bounded bumping
algorithm will decrement-push at $t+1$. If $(i,j-2) \in D$ also, then
$\bfa'= \Push_{t+1}^{-}\Push_{t}^{-} \mathbf{r}_{D}$ will again have duplicate
copies of the letter $i+j-3$ in positions $t+1$ and $t+2$ so the next
decrement-push will be in position $t+2$, and so on.  Note that since
the algorithm decrement-pushes both of the words in the bounded pair
in the same position at each iteration, the entrywise differences
$\bfa'-\bfb' = \mathbf{r}_{D} - \mathbf{j}_{D}$ agree, so the original
row numbers $\mathbf{i}_{D}$ are maintained throughout the bounded
bumping algorithm unless a deletion occurs.

We can group the push steps along one row so a decrement-push in
position $(i,j)$ pushes all of the adjacent $+$'s to its left over by
one as a stack. See \Cref{fig:stack.push} for example.   Thus, the
effect of the bounded bumping algorithm on the pipe dream amounts to a
sequence of such ``stack pushes''.  If at the end of a stack push, a
$+$ in column 1 of the pipe dream is decrement-pushed, the bounded
bump algorithm terminates by deleting that position because there will
be a 0 in the bounded word.  Otherwise, a stack push ends with a
bounded pair that corresponds to a pipe dream, which may or may not be
reduced. If it is reduced, the algorithm stops and returns
$\outcome=\bumped$.  Otherwise, we find the defect and continue with
another stack push in a different row.  In either case, the final
bounded pair $(\bfa',\bfb')$ encodes a reduced pipe dream.
\end{proof}

\begin{Definition}\label{def:stack.push}
A \textit{stack push} on a (not necessarily reduced) pipe dream $D$
starting at $(i,j)$ is the result of applying the consecutive sequence
of pushes on adjacent letters in $\mathbf{r}_{D}$ that correspond with
crossings on row $i$ when applying the bounded bump algorithm until
the defect is resolved or moves to another row.  Such stack pushes
necessarily correspond with a substring of $\mathbf{r}_{D}$ of
consecutive increasing or decreasing integers.
\end{Definition}

\begin{Exercise}\label{ex:stack.push}
For the pipe dream $D$ with $r_{D}=(8,7,6)$, $i_{D}=(2,2,2)$, and
$j_{D}=(7,6,5)$, verify $\Bump^{-}_{1}(\mathbf{r}_{D},
\mathbf{j}_{D})=((7,6,5), (6,5,4))$ can be done with one stack push on
row 2.  In particular, show the bounded pair $((7,6,5), (6,5,4))$ is
the reduced word and column word for a reduced pipe dream $D'$ and
$i_{D'}=(2,2,2)$, so all crossings are on the same row as $D$.
\end{Exercise}

Next, we give the promised bijective proof of \Cref{thm:pipedreams}
that proves the generating function for reduced pipe dreams for a
permutation $w$ satisfies the recurrence given by the transition equation (\Cref{t:transitionA}) for Schubert polynomials.  The idea is
to take any $D\in\rp(w)$, find the lex
largest inversion $(r,s) \in \Inv(w)$, and initiate the decrementing bounded bump
algorithm (so $\d = -1$) in the position of $D$ where wires $w_r$ and $w_s$ cross.  The first push initiates a stack push to
the left so a \+-tile either gets pushed out of the diagram into
column 0 or it moves to the first empty position, denoted by a bump
\elbow-tile, to its left if one exists.  In the first case, this only
happens if the initial crossing is on row $r$ and the remaining
diagram is a reduced pipe dream for $v=wt_{rs}$.  In the second case,
the resulting pipe dream may not be reduced, but we know it will be
nearly reduced so it has a well defined defect crossing, and we try to
rectify the situation by stack pushing the defect crossing one step to
its left and continue similarly until the defect is resolved.  The key point is that the $\d=-1$ bounded bump
applied to the crossing of the $w_{r}$-wire and the $w_{s}$-wire in
$D$ always terminates in a reduced pipe dream either for $v$ or for
one of the permutations $vt_{hr}$ appearing in the transition equation
from \Cref{t:transitionA}, and the process is reversible provided
$h,r$ are known.  See \Cref{fig:pipe dream bump} for an
example. Pinning down the details of this bijection does take some
work.  The proof given here follows the exposition in
\cite{Billey-Holroyd-Young}, where interested readers can find further
technical details.  Essentially the same construction has been given
by Buch according to \cite[p.11]{Knutson.2012}.  We begin by
describing how the bounded bump algorithm acts on reduced pipe dreams
to formally define the transition map and its inverse.

Warning! For this algorithm and proofs below, it is useful to relabel
the wires in a pipe dream by listing the positive integers in
increasing order from top to bottom on the left side and in order of
$w^{-1}$ on the top.  Then $(i,j) \in \Inv(w)$ if and only if
the $\{i,j \}$-wires cross.  This relabelling does not change the
monomial weight.  Similarly, we will use right-labeled wiring diagrams
of reduced words for the same reason.

\begin{proof}[Proof of \Cref{thm:pipedreams}] In the case $w =
\id$, we know $\fS_{w}=1$ and the only reduced pipe dream for the
identity is the empty pipe dream with no crossings, so the theorem
holds trivially.  So, assume $w \neq \id$, $(r,s)$ is the lex largest
inversion for $w$, and  $v = w t_{r,s}$.  Let
\begin{equation}\label{eq:Tset}
  \Tset(w):= \rp(v) \cup \bigcup_{\substack{ h<r:
      \\
      \ell(w)=\ell(v t_{hr}) }}
  \rp(v t_{hr}).
\end{equation}
We think of $v = v t_{r,r}$ so each pipe dream in $\Tset(w)$ is for a
permutation of the form $v t_{h,r}$ with $1\leq h\leq r$, though not
all such $v t_{h,r}$ necessarily occur.  We will give a bijection
$\TransitionMap: \rp(w) \longrightarrow \Tset(w)$ that preserves
weight, except in the cases $\TransitionMap(D)=E \in \rp(v)$, where
the weight will change by $x_r$, so $x^D = x_r x^E$.  Therefore, this
bijection will show that $\sum_{D\in\rp(w)}x^D$ satisfies the
transition equation for all permutations $w \in S_{\infty}$ with
$\ell(w)>0$. Since the transition equation along with the initial
condition $\fS_{\mathrm{id}}=1$ was used initially to define the
Schubert polynomial $\fS_w$ in \Cref{t:transitionA}, we will obtain
the desired equality.

\begin{algorithm}[\textbf{Transition Map}] \label{algorithm:TransitionMap}\
Suppose $w \neq \id$ is given, and let $(r,s)$ and $v$
be defined as above.

 \noindent \textbf{Input}: $D$, a non-empty reduced pipe dream for
 $w$ encoded as the biword $(\mathbf{r}_{D},\mathbf{j}_{D})$.

\noindent \textbf{Output}: $\TransitionMap(D) = E \in \Tset(w)$.

\begin{enumerate}

\item[(1)]  Let $t_0$ be the unique column containing the $\{r,s\}$-wiring
   crossing in the right-labeled wiring diagram for $\mathbf{r}_D$.

\item[(2)] Compute $\Bump^{-}_{t_0}(\mathbf{r}_{D}, \mathbf{j}_{D}) =
  (\bfa',\bfb', i, j, \outcome).$

\item[(3a)] If $\outcome=\deleted$, then $i=r-1$, $j=\ell(w)$, and
  $(\mathbf{a}', \bfb')$ encodes a pipe dream $E \in\rp(v) \subset
  \Tset(w)$.  Return $E$ and \textbf{stop}.

\item[(3b)] If $\outcome=\bumped$, then $(\mathbf{a}', \bfb')$ encodes a
    pipe dream $E \in\rp(v t_{hr})$ for some $h<r$ with
    $\ell(w)=\ell(v t_{hr})$.  Thus, $E \in \Tset(w)$.  Return
    $E$ and \textbf{stop}.
\end{enumerate}
 \end{algorithm}

\bigskip

 We claim that $\TransitionMap$ is a weight preserving bijection.
This claim follows from \Cref{item:1} since the bounded bump
algorithm is reversible.  We leave the details to the reader as an
exercise. 
 \end{proof}

\begin{Example}
  In Figure~\ref{fig:pipe dream bump}, we give an example computing
$\TransitionMap(D)$ for $w=1265734$.  Note that a defect can occur
either above or below the pushed crossing.  Going from the fourth to
the fifth pipe dream, two consecutive pushes on the same row are
combined into one step.  This is an example of a nontrivial stack
push (cf. Definition~\ref{def:stack.push}).  Here the curves in the elbow tiles have been straightened out,
which can be helpful visually.  Coloring the wires can also be
helpful.
\end{Example}

\bigskip

We note that the transition map on pipe dreams used in the proof of
\Cref{thm:pipedreams} is very close to, but different than, an
algorithm used earlier by Billey-Bergeron in \cite{billey-bergeron} to
give a bijective proof of Monk's formula
\[
\fS_{s_r}\fS_{w}=(x_1+\cdots+x_{r})\fS_{w}=\sum_{\substack{k\leq r<l\\
w t_{k,l}\gtrdot w}}\fS_{w t_{k,l}}.
\]
The key difference is the need to map pipe dreams for $w$ to those for
$wt_{kl}$ where $k\leq r <l$ with $\ell(wt_{kl})=\ell(w)+1$ when
inserting $x_{i}$ for $1\leq i\leq r$. This has to be built into the
algorithm.  We briefly review the BB-insertion algorithm for
comparison.

\medskip

\begin{algbreak}[\textbf{BB-insertion}] \label{algorithm:BB-insertion} 
\noindent \textbf{Input}:  An integer $1 \leq i
   \leq r$ and $D \in \rp(w)$ with wires labeled in increasing order
   down the left side.

 \noindent \textbf{Output}: $I(D,i,r)=E$, a reduced pipe dream for $wt_{kl}$
for some $k\leq r <l$ such that $\ell(wt_{k,l})=\ell(w)+1$.

\begin{enumerate}
\item Identify the smallest column $j$ on row $i$ such that there is a
bump tile in position $(i,j)$ and  the wire $k$ from the left and the
wire $l$
from below satisfy the \textit{insertion condition}
\begin{equation}\label{eq:cross.condition}
k\leq r<l.  
\end{equation}
Change the tile at $(i,j)$ to a cross to obtain the pipe dream $E$.  Note, such a pair of wires
always exists assuming one can use the invisible wire labeled $n+1$ if
necessary since the $i$ wire enters row $i$ on the left edge and
some wire on row $i$ (possibly invisible) has larger label than $r$.

\item If $E$ is reduced, stop and return $E \in \rp(wt_{kl})$.
\item If $E$ is not reduced, then wires $(k,l)$ cross again say in
entry $(i',j')$ such that $i'<i$.  Remove the crossing at $(i',j')$
from $E$ and add another crossing on row $i'$ in the largest column
$j''<j'$ such that two wires $k'<l'$ come together but don't cross
satisfying the insertion condition $k'\leq r < l'$.  Again such an
entry must exist.  Update $E$ to be the resulting pipe dream,
$i \leftarrow i'$, $j\leftarrow j''$, and repeat step 2. 
\end{enumerate}
\end{algbreak}

For example, consider the reduced pipe dream for $21534$ 
\[
D=\begin{matrix}
+&.&.&+\\
.&+&.\\
.&.\\
.
\end{matrix}
\]
The reader may want to draw the wires into the pipe dreams above and below.
To compute $I(D,2,3)$, we insert a cross at position $(2,3)$ where the
$3$-wire and $5$-wire bump, 
\[
\begin{matrix}
+&.&.&+\\
.&+&+\\
.&.\\
.
\end{matrix}
\]
but this is not reduced.  Remove the $+$ at $(1,4)$ and reinsert it further left on
row 1, as follows.  Position $(1,3)$ is a bump tile with wires $4$ and $5$, so we skip
that and consider position $(1,2)$ where wires $1$ and $4$ bump.
Since $1\leq r=3<4$, we place a plus tile at $(1,2)$ to get 
\[
E=\begin{matrix}
+&+&.&.\\
.&+&+\\
.&.\\
.
\end{matrix}
\]
Now $E$ is reduced for $v = 31524 = wt_{14}$, so $I(D,2,3)=E$.  Note, the BB-insertion algorithm
doesn't just stack push the crossing on row 1, but instead removes it
and reinserts in a way that maintains the insertion condition.  It
also always finds the defect in a smaller row as it proceeds through
the steps.

\begin{Exercise}\label{ex:BB-inverse}
Write out the inverse map for BB-insertion.
\end{Exercise}

\subsection{The nil-Coxeter Algebra and Divided
Differences}\label{sub:Fomin-Stanley}

Fomin and Stanley \cite{FS} utilized the nil-Coxeter algebra to study
Schubert polynomials.  In their setup, pipe dreams are ``natural''
combinatorial objects to study.  In our revisionist history, we will
use their proof in reverse to prove that Schubert polynomials satisfy
the divided difference recurrence given in Equation \eqref{eq:divided.difference}.
See the later chapter \cite{richmond2023nilheckeringsschubertcalculus} in this handbook for more on Nil-Hecke rings.

\begin{Definition}
The \emph{nil-Coxeter algebra} $\mathcal{N}_n$ has generators $u_1,u_2,\ldots,u_{n-1}$ with coefficients in $\mathbb{Z}[x_1,x_2,\ldots,x_n]$ and the following relations:
\begin{itemize}
\item $u_i^2=0$,
\item $u_iu_j=u_ju_i$ if $|i-j|\geq2$,
\item $u_{i}u_{i+1}u_i=u_{i+1}u_iu_{i+1}$ for $1\leq i\leq n-2$.
\end{itemize}
\end{Definition}
Note that these relations are reminiscent of those for the symmetric
group generated by the simple transpositions $s_i$'s in \eqref{eq:simple.relations}. Consequently, we
know that $\mathcal{N}_n$ has a linear basis given by $\{u_w\}_{w\in
S_n}$ where \[u_w=u_{i_1}u_{i_2}\cdots u_{i_{\ell(w)}}\] and
$s_{i_1}s_{i_2}\cdots s_{i_{\ell(w)}}$ is a reduced word for
$w$. Since $u_{\id}$ is the multiplicative identity in $\mathcal{N}_{n}$, we
sometimes represent it by the identity $1$ in
$\mathbb{Z}[x_1,x_2,\ldots,x_n]$. For example, $(x_{1}+x_{2})u_{231} +
x_{3}^{2}u_{123}=(x_{1}+x_{2})u_{231} + x_{3}^{2}\ \in \mathcal{N}_{3}$.

With $n$ fixed, define the following expressions: 
\begin{itemize}
\item $h_i(x)=1+xu_i$ for $i=1,\ldots,n-1$,
\item $A_i(x)=h_{n-1}(x)h_{n-2}(x)\cdots h_i(x)$ for $i=1,\ldots,n-1$, and
\item $\fS(\bfx) = \fS(x_{1},\ldots, x_{n}) =A_1(x_1)A_2(x_2)\cdots A_{n-1}(x_{n-1})$.
\end{itemize}
In other words, $\fS(\bfx)$ is a product of the $\binom{n}{2}$ binomials
\begin{equation}\label{eq:schub}
\begin{split}
\fS(\bfx)=(1+x_1u_{n-1})(1+x_1u_{n-2})\cdots (1+x_1u_2)(1+x_1u_1)&\\
(1+x_2u_{n-1})\cdots (1+x_2u_3)(1+x_2u_2)&\\
\ddots\ \ \ \ \ \ \ \ \ \ \ \ \ \ \ \ \vdots&\\
(1+x_{n-1}u_{n-1}).&
\end{split}
\end{equation}

Pipe dreams are hidden in plain sight in \eqref{eq:schub}. Each factor in \Cref{eq:schub} consists of two terms: $1$ which means a ``bump''-tile and $x_iu_j$ which means a ``cross''-tile on row $i$ that contributes $s_j$ to the resulting
reduced expression.  

\begin{Example}
Let $n=4$. Consider the following underlined term in $\fS(w)$ with its corresponding pipe dream on the side in \Cref{fig:FS-pipe-dream-example}. This pipe dream has weight $x_1x_2^2$ that contributes to the permutation $w=s_1s_3s_2=2413$.
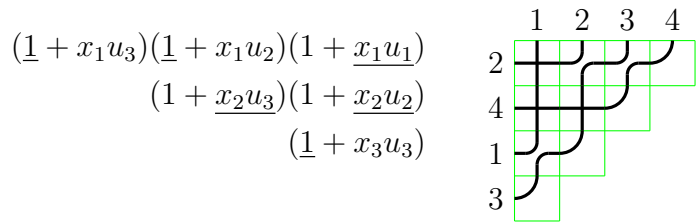
\begin{figure}[h!]
\centering
\begin{tikzpicture}[scale=0.5]
\node at (0,0) {$\begin{aligned}
(\underline{1}+x_1u_3)(\underline{1}+x_1u_2)(1+\underline{x_1u_1})&\\
(1+\underline{x_2u_3})(1+\underline{x_2u_2})&\\
(\underline{1}+x_3u_3)&
\end{aligned}$};
\node at (0,-3) {};
\end{tikzpicture}
\quad
\begin{tikzpicture}[scale=0.300000000000000]
\draw[color=green](2,-2)--(10,-2);
\draw[color=green](2,-2)--(2,-10);
\draw[color=green](2,-4)--(10,-4);
\draw[color=green](4,-2)--(4,-10);
\draw[color=green](2,-6)--(8,-6);
\draw[color=green](6,-2)--(6,-8);
\draw[color=green](2,-8)--(6,-8);
\draw[color=green](8,-2)--(8,-6);
\draw[color=green](2,-10)--(4,-10);
\draw[color=green](10,-2)--(10,-4);
\draw[very thick](2,-3)--(4,-3);
\draw[very thick](3,-2)--(3,-4);
\draw[very thick](5,-2)--(5,-2.50000000000000);
\draw[very thick](5,-3.50000000000000)--(5,-4);
\draw[very thick](4,-3)--(4.50000000000000,-3);
\draw[very thick](5.50000000000000,-3)--(6,-3);
\draw[very thick](4.50000000000000,-3)arc(270:360:0.500000000000000);
\draw[very thick](5.50000000000000,-3)arc(90:180:0.500000000000000);
\draw[very thick](7,-2)--(7,-2.50000000000000);
\draw[very thick](7,-3.50000000000000)--(7,-4);
\draw[very thick](6,-3)--(6.50000000000000,-3);
\draw[very thick](7.50000000000000,-3)--(8,-3);
\draw[very thick](6.50000000000000,-3)arc(270:360:0.500000000000000);
\draw[very thick](7.50000000000000,-3)arc(90:180:0.500000000000000);
\draw[very thick](2,-5)--(4,-5);
\draw[very thick](3,-4)--(3,-6);
\draw[very thick](4,-5)--(6,-5);
\draw[very thick](5,-4)--(5,-6);
\draw[very thick](3,-6)--(3,-6.50000000000000);
\draw[very thick](3,-7.50000000000000)--(3,-8);
\draw[very thick](2,-7)--(2.50000000000000,-7);
\draw[very thick](3.50000000000000,-7)--(4,-7);
\draw[very thick](2.50000000000000,-7)arc(270:360:0.500000000000000);
\draw[very thick](3.50000000000000,-7)arc(90:180:0.500000000000000);
\draw[very thick](8,-3)arc(270:360:1);
\draw[very thick](6,-5)arc(270:360:1);
\draw[very thick](4,-7)arc(270:360:1);
\draw[very thick](2,-9)arc(270:360:1);
\node[left] at (2,-3) {$2$};
\node[left] at (2,-5) {$4$};
\node[left] at (2,-7) {$1$};
\node[left] at (2,-9) {$3$};
\node[above] at (3,-2) {$1$};
\node[above] at (5,-2) {$2$};
\node[above] at (7,-2) {$3$};
\node[above] at (9,-2) {$4$};
\end{tikzpicture}
\caption{A term in $\fS(w)$ and its corresponding pipe dream}
\label{fig:FS-pipe-dream-example}
\end{figure}
There is a horizontal reflection between the two pictures due to the order that we read the crossings in a pipe dream as shown in \Cref{fig:rcgraphs}.
\end{Example}

As a result of the observation that each factor in \Cref{eq:schub}
consists of two terms corresponding with the bump and cross tiles, we
deduced the following lemma.

\begin{Lemma} \cite{FS}
The coefficients of $\fS(\bfx)$ expanded in the $\{u_{w} \}$ basis are
the Schubert polynomials, 
\[\fS(\bfx)=\sum_{w\in S_n}\fS_w(x_{1},\dots , x_{n})u_w.\]
\end{Lemma}

Recall, a central player in Schubert calculus is the \emph{divided difference operator} $\partial_i$, for $i=1,2,\ldots,n-1$ that acts on polynomials by
\[\partial_i f=\frac{f-s_if}{x_i-x_{i+1}},\]
where $s_i$ acts by swapping $x_i$ and $x_{i+1}$.
We extend the action of $\partial_i$ to the nil-Coxeter algebra where
$\partial_i$ acts on the coefficients in
$\mathbb{Z}[x_1,\ldots,x_n]$.  Recall from
\Cref{ex:divided-difference} that these operators satisfy the
following relations: 
\begin{itemize}
\item $\partial_i^2=0$,
\item $\partial_i\partial_j=\partial_j\partial_i$ if $|i-j|\geq2$,
\item $\partial_i\partial_{i+1}\partial_i=\partial_{i+1}\partial_i\partial_{i+1}$ for $1\leq i\leq n-2$.
\end{itemize}
The key algebraic computation in \cite{FS} is the following lemma.

\begin{Lemma} \cite[Lemma 3.5]{FS}\label{lem:FS-divided-difference-right-multiplication} For every $i=1,\ldots,n-1$,
\begin{equation}\label{eq:schub-calculation}
\partial_i\fS(\bfx)=\fS(\bfx)u_i.
\end{equation}
\end{Lemma}
Establishing \Cref{lem:FS-divided-difference-right-multiplication}
boils down to proving the identity 
\[\partial_i\big(A_i(x_i)A_{i+1}(x_{i+1})\big)=A_i(x_i)A_{i+1}(x_{i+1})u_i,\] which is a straightforward algebraic manipulation. We omit the details here. Interested readers can check out \cite{FS} and \cite{manivel-book}.

By comparing the coefficients of $u_w$ on both sides of
\Cref{eq:schub-calculation}, we obtain the following recurrence for the
Schubert polynomials. This formula is more widely known as the
definition of Schubert polynomials by Lascoux and Sch\"utzenberger
\cite{l-s}. However, in our approach the following theorem is a
consequence of the proof of \cite[Thm 2.2]{FS}.

\begin{Theorem} \cite{FS,l-s} \label{cor:schubert-divided-difference}
For $w\in S_n$,
\[\fS_w(x_1,\ldots,x_n)=\begin{cases}x_1^{n-1}x_2^{n-2}\cdots x_{n-1}
&\text{if }w=w_0=[n,\ldots,1],\\
\partial_i\fS_{ws_i}(x_{1},\ldots,x_{n} )&\text{if }w(i)<w(i+1).
\end{cases}\]
\end{Theorem}

The framework of the nil-Coxeter algebra has more applications.  Let
$\fS(x,\ldots,x)$ be the result of substituting a single indeterminate $x$
in for all $x_1,x_2, \ldots,x_{n}$ in $\fS(\bfx)$.

\begin{Lemma} \cite[Lemma 5.1]{FS}\label{lem:FS-exponential-product}
For indeterminates $x,y$, we have 
$$\fS(x,\ldots,x)\fS(y,\ldots,y)=\fS(x+y,\ldots,x+y).$$

\end{Lemma}
The proof of \Cref{lem:FS-exponential-product} can be done via
explicit calculation.  We encourage the reader to try the cases
$n=3,4$.

\begin{Lemma}\cite{FS}\label{lem:FS-exponential-formula}
We have
$$\fS(x,\ldots,x)=\exp(x(u_1+2u_2+\cdots+(n-1)u_{n-1})).$$
\end{Lemma}
\begin{proof}
By \Cref{lem:FS-exponential-product}, one may assume there exists
$f\in\mathcal{N}_n$ such that $\fS(x,\ldots,x)=\exp(xf)$.
To solve for $f$, observe that we must have 
\[f=\frac{d}{dx}\exp(xf)|_{x=0}=\frac{d}{dx}\fS(x,\ldots,x)|_{x=0}=u_1+2u_2+\cdots+(n-1)u_{n-1}\]
since $u_i$ appears in $i$ linear factors in \Cref{eq:schub}.
\end{proof}

\begin{Theorem}[\textbf{Macdonald's identity}]\cite[(6.11)]{M2} 
\label{cor:Macdonald-identity-q-1}
For $w\in S_n$, the number of reduced pipe dreams of $w$
equals \[|\rp(w)|=\fS_{w}(1,\ldots,1)= \frac{1}{\ell!}\sum_{(a_1,\ldots,a_{\ell})\in R(w)}a_1 a_2\cdots a_{\ell}.\]
\end{Theorem}

\begin{proof}
The first equality follows from \Cref{thm:pipedreams}. The second
equality follows from \Cref{lem:FS-exponential-formula} by setting
$x=1$ comparing the coefficient of $u_w$ on both sides of the
resulting identity.  
\end{proof}

We will outline an additional proof of
\Cref{cor:Macdonald-identity-q-1} that has been inspirational in the
literature. With \Cref{lem:FS-exponential-formula}, Fomin and Stanley
in the 90's hinted at taking derivatives of Schubert structures from a
certain perspective. This concept was further developed by Hamaker, Pechenik,
Speyer and Weigandt \cite{hamaker2020derivatives} in 2020, who were
motivated by Stanley's operator as an attempt to establish the Sperner
property of the weak Bruhat order. Curiously, Stanley presented his
operator during FominFest in 2018, a celebration for Fomin's 60
birthday. The Sperner property of the weak Bruhat order on the
symmetric group is first proved by Gaetz and Gao
\cite{gaetzgao-sperner-2,gaetzgao-sperner-1}, while Hamaker et
al. \cite{hamaker2020derivatives} followed up to prove the full
determinantal conjecture by Stanley \cite{stanley2017shenanigans}.
\begin{Theorem}\cite{hamaker2020derivatives}\label{thm:derivative-schubert}
For any fixed positive integer $n$, 
\[\nabla=\frac{\partial}{\partial x_1}+\frac{\partial}{\partial x_2}+\cdots+\frac{\partial}{\partial x_{n}}.\]
Then applying $\nabla$ to the Schubert polynomial for $w \in S_{n}$, we have
\[\nabla\fS_w(x_{1},\dots ,
x_{n})=\sum_{k:s_kw<w}k\fS_{s_kw}(x_{1},\dots , x_{n}).\]
\end{Theorem}
\begin{Exercise}
Use \Cref{thm:derivative-schubert} to prove Macdonald's identity \Cref{cor:Macdonald-identity-q-1}.
\end{Exercise}

We end this section by remarking that the setup here is much more powerful than what's presented. In particular, Fomin and Kirillov \cite{FK} are able to derive the pipe dream formula for double Schubert polynomials from the divided difference operators. Furthermore, a $q$-analog of Macdonald's identity can be proved.
\begin{Theorem}[Theorem~2.4 of \cite{FS}]\label{thm:FS}
For a permutation $w \in S_{n}$, \[\fS_w(1,q,\ldots,q^{n-1})=\frac{1}{[\ell]_q!}\sum_{(a_1,\ldots,a_{\ell})\in R(w)}[a_1]_q\cdots [a_\ell]_qq^{\sum_{a_i<a_{i+1}}i}\]
where $[a]_q=1+q+\cdots+q^{a-1}$, and $[m]_q!=[1]_q\cdots[m]_q$.
\end{Theorem}

A bijective proof of \Cref{thm:FS} was given in
\cite{Billey-Holroyd-Young} using the bounded bump algorithm
\Cref{algorithm:little bump} and a generalization of the Transition
Equation. Holroyd proposed studying the distribution on reduced
expressions for the longest permutation $w_{0}=[n,\ldots,1] \in S_{n}$ where $(a_{1},\dots , a_{\ell})$ is
chosen with probability proportional to the product $a_{1}\cdots
a_{\ell}$.  Does this distribution have nice properties as with the
uniform distribution studied in \cite{angel2007random,Dauvergne}?  The bounded
bump algorithm can be used to select random reduced words according to
this distribution.

Recently, Nadeau and Tewari recognized an interesting similarity
between Macdonald's identity and an identity due to Klyachko in his
study of the closure of the orbit of a generic flag in $\flags$ under
the action of the invertible diagonal matrices
\cite{nadeau.tewari.2021,Klyachko.1985}. 
In addition, they describe a $q$-deformation of what they call the
Klyachko algebra on noncommuting indeterminates $u_{1}',\dots , u_{n}'$
in a style similar to the nil-Coxeter algebra in order to give a
$q$-Klyachko-Macdonald identity
\[
\fS_{w}(u_{1}', u_{2}'-u_{1}',\dots ) =
\frac{1}{[\ell]_q!}\sum_{(a_1,\ldots,a_{\ell})\in R(w)}u_{a_{1}}'\cdots
u_{a_\ell}' q^{\sum_{a_i<a_{i+1}}i}
\]
They connect this identity to the geometry of certain Deligne-Lusztig
varieties connected to the representation of finite group of Lie type
using results of Kim \cite{Kim.2020}.

\subsection{More Games: Mitosis}\label{sub:Mitosis}

If Schubert calculus had started with Monk's formula, the transition
equation, and pipe dreams, would we have discovered the divided
difference recurrence for Schubert polynomials?  Good question, maybe
not.  But, in this revisionist history of Schubert polynomials, one
could ask if there is a combinatorial game style procedure that
relates to the divided difference recurrence
\eqref{eq:divided.difference}?  Yes!  This was given by Knutson-Miller
in their ``mitosis'' algorithm \cite{knutson-miller-2005}.  We follow
\cite{Miller.2003} in this exposition.  See also \cite{B,lenart.2004} for other variations on this theme.

Given a pipe dream $D$, let $\mathrm{start}_i(D)$ be the column index
of the leftmost \elbow-tile in row $i$, and let
$C_i(D)=\{q<\mathrm{start}_i(D)\:|\: (i+1,q)\text{ is an \elbow-tile
in }D\}$. For each $q\in C_i(D)$, the \emph{offspring} $D_{i,q}$ is
obtained from $D$ by deleting the cross tile at $(i,q)$ and then
moving all crosses $(i,p)$ in row $i$ one step down to $(i+1,p)$ for
all $p<q$ and $p\in C_i(D)$. See Figure~\ref{fig:mitosis-example} for example.
\begin{Definition}
The $i$th \emph{mitosis operator} acts on each pipe dream $D$ by producing
its set of $i$-offspring 
\[\mathrm{mitosis}_i(D)=\{D_{i,q}\:|\: q\in C_i(D)\}.\]
We also write $\mathrm{mitosis}_i(\mathcal{D})=\bigcup_{D\in\mathcal{D}}\mathrm{mitosis}_i(D)$ for a set of pipe dreams $\mathcal{D}$.
\end{Definition}
\begin{Exercise}
Prove the following statements. \begin{enumerate}
\item If $i$ is an ascent of $w$, then $\mathrm{mitosis}_i(\rp(w))=\emptyset$.
\item If $i$ is a descent of $w$, then every
$D'\in\mathrm{mitosis}_i(\rp(w))$ is a pipe dream of $ws_i$.
\end{enumerate}

\end{Exercise}

The main theorem of this section is that all pipe dreams can be
generated via mitosis from the unique pipe dream of the longest
permutation $w_0\in S_n$. Readers are referred to \cite{Miller.2003}
for a detailed proof.

\begin{Theorem}\cite{knutson-miller-2005,Miller.2003}\label{thm:mitosis}
For any $w\in S_n$ with $\ell(ws_i)<\ell(w)$, there is a partition of
$\rp(ws_i)$ into a disjoint union of offspring given by  \[\rp(ws_i)=\bigsqcup_{D\in\rp(w)}\mathrm{mitosis}_i(D).\]
As a result, let $D_0$ be the unique pipe dream for the longest permutation $w_0\in S_n$, then \[\rp(w)=\mathrm{mitosis}_{i_{\ell}}\cdots\mathrm{mitosis}_{i_1}(D_0)\]
where $s_{i_1}\cdots s_{i_\ell}$ is any reduced word for $w_0w$. 
\end{Theorem}

\begin{Example}
Consider the following pipe dream $D\in\rp(261453)$ in Figure~\ref{fig:mitosis-example}. We see that
$\mathrm{start}_2(D)=5$ and $C_2(D)=\{1,3,4\}$ with the three
offspring in 
$\mathrm{mitosis}_2(D)$ shown on the right in Figure~\ref{fig:mitosis-example}, all
of which are pipe dreams of $216453=261453s_2$.
\begin{figure}[h!]
\centering
\begin{tikzpicture}[scale=0.2]
\def\a{-30};
\def\b{0};
\draw[color=green](\a+2,-2+\b)--(\a+14,-2+\b);
\draw[color=green](\a+2,-2+\b)--(\a+2,-14+\b);
\draw[color=green](\a+2,-4+\b)--(\a+14,-4+\b);
\draw[color=green](\a+4,-2+\b)--(\a+4,-14+\b);
\draw[color=green](\a+2,-6+\b)--(\a+12,-6+\b);
\draw[color=green](\a+6,-2+\b)--(\a+6,-12+\b);
\draw[color=green](\a+2,-8+\b)--(\a+10,-8+\b);
\draw[color=green](\a+8,-2+\b)--(\a+8,-10+\b);
\draw[color=green](\a+2,-10+\b)--(\a+8,-10+\b);
\draw[color=green](\a+10,-2+\b)--(\a+10,-8+\b);
\draw[color=green](\a+2,-12+\b)--(\a+6,-12+\b);
\draw[color=green](\a+12,-2+\b)--(\a+12,-6+\b);
\draw[color=green](\a+2,-14+\b)--(\a+4,-14+\b);
\draw[color=green](\a+14,-2+\b)--(\a+14,-4+\b);
\draw[very thick](\a+2,-3+\b)--(\a+4,-3+\b);
\draw[very thick](\a+3,-2+\b)--(\a+3,-4+\b);
\draw[very thick](\a+5,-2+\b)--(\a+5,-2.5+\b);
\draw[very thick](\a+5,-3.5+\b)--(\a+5,-4+\b);
\draw[very thick](\a+4,-3+\b)--(\a+4.5,-3+\b);
\draw[very thick](\a+5.5,-3+\b)--(\a+6,-3+\b);
\draw[very thick](\a+4.5,-3+\b)arc(270:360:0.5);
\draw[very thick](\a+5.5,-3+\b)arc(90:180:0.5);
\draw[very thick](\a+7,-2+\b)--(\a+7,-2.5+\b);
\draw[very thick](\a+7,-3.5+\b)--(\a+7,-4+\b);
\draw[very thick](\a+6,-3+\b)--(\a+6.5,-3+\b);
\draw[very thick](\a+7.5,-3+\b)--(\a+8,-3+\b);
\draw[very thick](\a+6.5,-3+\b)arc(270:360:0.5);
\draw[very thick](\a+7.5,-3+\b)arc(90:180:0.5);
\draw[very thick](\a+9,-2+\b)--(\a+9,-2.5+\b);
\draw[very thick](\a+9,-3.5+\b)--(\a+9,-4+\b);
\draw[very thick](\a+8,-3+\b)--(\a+8.5,-3+\b);
\draw[very thick](\a+9.5,-3+\b)--(\a+10,-3+\b);
\draw[very thick](\a+8.5,-3+\b)arc(270:360:0.5);
\draw[very thick](\a+9.5,-3+\b)arc(90:180:0.5);
\draw[very thick](\a+11,-2+\b)--(\a+11,-2.5+\b);
\draw[very thick](\a+11,-3.5+\b)--(\a+11,-4+\b);
\draw[very thick](\a+10,-3+\b)--(\a+10.5,-3+\b);
\draw[very thick](\a+11.5,-3+\b)--(\a+12,-3+\b);
\draw[very thick](\a+10.5,-3+\b)arc(270:360:0.5);
\draw[very thick](\a+11.5,-3+\b)arc(90:180:0.5);
\draw[very thick](\a+2,-5+\b)--(\a+4,-5+\b);
\draw[very thick](\a+3,-4+\b)--(\a+3,-6+\b);
\draw[very thick](\a+4,-5+\b)--(\a+6,-5+\b);
\draw[very thick](\a+5,-4+\b)--(\a+5,-6+\b);
\draw[very thick](\a+6,-5+\b)--(\a+8,-5+\b);
\draw[very thick](\a+7,-4+\b)--(\a+7,-6+\b);
\draw[very thick](\a+8,-5+\b)--(\a+10,-5+\b);
\draw[very thick](\a+9,-4+\b)--(\a+9,-6+\b);
\draw[very thick](\a+3,-6+\b)--(\a+3,-6.5+\b);
\draw[very thick](\a+3,-7.5+\b)--(\a+3,-8+\b);
\draw[very thick](\a+2,-7+\b)--(\a+2.5,-7+\b);
\draw[very thick](\a+3.5,-7+\b)--(\a+4,-7+\b);
\draw[very thick](\a+2.5,-7+\b)arc(270:360:0.5);
\draw[very thick](\a+3.5,-7+\b)arc(90:180:0.5);
\draw[very thick](\a+4,-7+\b)--(\a+6,-7+\b);
\draw[very thick](\a+5,-6+\b)--(\a+5,-8+\b);
\draw[very thick](\a+7,-6+\b)--(\a+7,-6.5+\b);
\draw[very thick](\a+7,-7.5+\b)--(\a+7,-8+\b);
\draw[very thick](\a+6,-7+\b)--(\a+6.5,-7+\b);
\draw[very thick](\a+7.5,-7+\b)--(\a+8,-7+\b);
\draw[very thick](\a+6.5,-7+\b)arc(270:360:0.5);
\draw[very thick](\a+7.5,-7+\b)arc(90:180:0.5);
\draw[very thick](\a+3,-8+\b)--(\a+3,-8.5+\b);
\draw[very thick](\a+3,-9.5+\b)--(\a+3,-10+\b);
\draw[very thick](\a+2,-9+\b)--(\a+2.5,-9+\b);
\draw[very thick](\a+3.5,-9+\b)--(\a+4,-9+\b);
\draw[very thick](\a+2.5,-9+\b)arc(270:360:0.5);
\draw[very thick](\a+3.5,-9+\b)arc(90:180:0.5);
\draw[very thick](\a+4,-9+\b)--(\a+6,-9+\b);
\draw[very thick](\a+5,-8+\b)--(\a+5,-10+\b);
\draw[very thick](\a+3,-10+\b)--(\a+3,-10.5+\b);
\draw[very thick](\a+3,-11.5+\b)--(\a+3,-12+\b);
\draw[very thick](\a+2,-11+\b)--(\a+2.5,-11+\b);
\draw[very thick](\a+3.5,-11+\b)--(\a+4,-11+\b);
\draw[very thick](\a+2.5,-11+\b)arc(270:360:0.5);
\draw[very thick](\a+3.5,-11+\b)arc(90:180:0.5);
\draw[very thick](\a+12,-3+\b)arc(270:360:1);
\draw[very thick](\a+10,-5+\b)arc(270:360:1);
\draw[very thick](\a+8,-7+\b)arc(270:360:1);
\draw[very thick](\a+6,-9+\b)arc(270:360:1);
\draw[very thick](\a+4,-11+\b)arc(270:360:1);
\draw[very thick](\a+2,-13+\b)arc(270:360:1);
\node[left] at (\a+2,-3+\b) {$2$};
\node[left] at (\a+2,-5+\b) {$6$};
\node[left] at (\a+2,-7+\b) {$1$};
\node[left] at (\a+2,-9+\b) {$4$};
\node[left] at (\a+2,-11+\b) {$5$};
\node[left] at (\a+2,-13+\b) {$3$};
\node[above] at (\a+3,-2+\b) {$1$};
\node[above] at (\a+5,-2+\b) {$2$};
\node[above] at (\a+7,-2+\b) {$3$};
\node[above] at (\a+9,-2+\b) {$4$};
\node[above] at (\a+11,-2+\b) {$5$};
\node[above] at (\a+13,-2+\b) {$6$};

\def\a{0};
\def\b{0};
\draw[color=green](\a+2,-2+\b)--(\a+14,-2+\b);
\draw[color=green](\a+2,-2+\b)--(\a+2,-14+\b);
\draw[color=green](\a+2,-4+\b)--(\a+14,-4+\b);
\draw[color=green](\a+4,-2+\b)--(\a+4,-14+\b);
\draw[color=green](\a+2,-6+\b)--(\a+12,-6+\b);
\draw[color=green](\a+6,-2+\b)--(\a+6,-12+\b);
\draw[color=green](\a+2,-8+\b)--(\a+10,-8+\b);
\draw[color=green](\a+8,-2+\b)--(\a+8,-10+\b);
\draw[color=green](\a+2,-10+\b)--(\a+8,-10+\b);
\draw[color=green](\a+10,-2+\b)--(\a+10,-8+\b);
\draw[color=green](\a+2,-12+\b)--(\a+6,-12+\b);
\draw[color=green](\a+12,-2+\b)--(\a+12,-6+\b);
\draw[color=green](\a+2,-14+\b)--(\a+4,-14+\b);
\draw[color=green](\a+14,-2+\b)--(\a+14,-4+\b);
\draw[very thick](\a+2,-3+\b)--(\a+4,-3+\b);
\draw[very thick](\a+3,-2+\b)--(\a+3,-4+\b);
\draw[very thick](\a+5,-2+\b)--(\a+5,-2.5+\b);
\draw[very thick](\a+5,-3.5+\b)--(\a+5,-4+\b);
\draw[very thick](\a+4,-3+\b)--(\a+4.5,-3+\b);
\draw[very thick](\a+5.5,-3+\b)--(\a+6,-3+\b);
\draw[very thick](\a+4.5,-3+\b)arc(270:360:0.5);
\draw[very thick](\a+5.5,-3+\b)arc(90:180:0.5);
\draw[very thick](\a+7,-2+\b)--(\a+7,-2.5+\b);
\draw[very thick](\a+7,-3.5+\b)--(\a+7,-4+\b);
\draw[very thick](\a+6,-3+\b)--(\a+6.5,-3+\b);
\draw[very thick](\a+7.5,-3+\b)--(\a+8,-3+\b);
\draw[very thick](\a+6.5,-3+\b)arc(270:360:0.5);
\draw[very thick](\a+7.5,-3+\b)arc(90:180:0.5);
\draw[very thick](\a+9,-2+\b)--(\a+9,-2.5+\b);
\draw[very thick](\a+9,-3.5+\b)--(\a+9,-4+\b);
\draw[very thick](\a+8,-3+\b)--(\a+8.5,-3+\b);
\draw[very thick](\a+9.5,-3+\b)--(\a+10,-3+\b);
\draw[very thick](\a+8.5,-3+\b)arc(270:360:0.5);
\draw[very thick](\a+9.5,-3+\b)arc(90:180:0.5);
\draw[very thick](\a+11,-2+\b)--(\a+11,-2.5+\b);
\draw[very thick](\a+11,-3.5+\b)--(\a+11,-4+\b);
\draw[very thick](\a+10,-3+\b)--(\a+10.5,-3+\b);
\draw[very thick](\a+11.5,-3+\b)--(\a+12,-3+\b);
\draw[very thick](\a+10.5,-3+\b)arc(270:360:0.5);
\draw[very thick](\a+11.5,-3+\b)arc(90:180:0.5);
\draw[very thick](\a+3,-4+\b)--(\a+3,-4.5+\b);
\draw[very thick](\a+3,-5.5+\b)--(\a+3,-6+\b);
\draw[very thick](\a+2,-5+\b)--(\a+2.5,-5+\b);
\draw[very thick](\a+3.5,-5+\b)--(\a+4,-5+\b);
\draw[very thick](\a+2.5,-5+\b)arc(270:360:0.5);
\draw[very thick](\a+3.5,-5+\b)arc(90:180:0.5);
\draw[very thick](\a+4,-5+\b)--(\a+6,-5+\b);
\draw[very thick](\a+5,-4+\b)--(\a+5,-6+\b);
\draw[very thick](\a+6,-5+\b)--(\a+8,-5+\b);
\draw[very thick](\a+7,-4+\b)--(\a+7,-6+\b);
\draw[very thick](\a+8,-5+\b)--(\a+10,-5+\b);
\draw[very thick](\a+9,-4+\b)--(\a+9,-6+\b);
\draw[very thick](\a+3,-6+\b)--(\a+3,-6.5+\b);
\draw[very thick](\a+3,-7.5+\b)--(\a+3,-8+\b);
\draw[very thick](\a+2,-7+\b)--(\a+2.5,-7+\b);
\draw[very thick](\a+3.5,-7+\b)--(\a+4,-7+\b);
\draw[very thick](\a+2.5,-7+\b)arc(270:360:0.5);
\draw[very thick](\a+3.5,-7+\b)arc(90:180:0.5);
\draw[very thick](\a+4,-7+\b)--(\a+6,-7+\b);
\draw[very thick](\a+5,-6+\b)--(\a+5,-8+\b);
\draw[very thick](\a+7,-6+\b)--(\a+7,-6.5+\b);
\draw[very thick](\a+7,-7.5+\b)--(\a+7,-8+\b);
\draw[very thick](\a+6,-7+\b)--(\a+6.5,-7+\b);
\draw[very thick](\a+7.5,-7+\b)--(\a+8,-7+\b);
\draw[very thick](\a+6.5,-7+\b)arc(270:360:0.5);
\draw[very thick](\a+7.5,-7+\b)arc(90:180:0.5);
\draw[very thick](\a+3,-8+\b)--(\a+3,-8.5+\b);
\draw[very thick](\a+3,-9.5+\b)--(\a+3,-10+\b);
\draw[very thick](\a+2,-9+\b)--(\a+2.5,-9+\b);
\draw[very thick](\a+3.5,-9+\b)--(\a+4,-9+\b);
\draw[very thick](\a+2.5,-9+\b)arc(270:360:0.5);
\draw[very thick](\a+3.5,-9+\b)arc(90:180:0.5);
\draw[very thick](\a+4,-9+\b)--(\a+6,-9+\b);
\draw[very thick](\a+5,-8+\b)--(\a+5,-10+\b);
\draw[very thick](\a+3,-10+\b)--(\a+3,-10.5+\b);
\draw[very thick](\a+3,-11.5+\b)--(\a+3,-12+\b);
\draw[very thick](\a+2,-11+\b)--(\a+2.5,-11+\b);
\draw[very thick](\a+3.5,-11+\b)--(\a+4,-11+\b);
\draw[very thick](\a+2.5,-11+\b)arc(270:360:0.5);
\draw[very thick](\a+3.5,-11+\b)arc(90:180:0.5);
\draw[very thick](\a+12,-3+\b)arc(270:360:1);
\draw[very thick](\a+10,-5+\b)arc(270:360:1);
\draw[very thick](\a+8,-7+\b)arc(270:360:1);
\draw[very thick](\a+6,-9+\b)arc(270:360:1);
\draw[very thick](\a+4,-11+\b)arc(270:360:1);
\draw[very thick](\a+2,-13+\b)arc(270:360:1);

\def\a{14};
\def\b{0};
\draw[color=green](\a+2,-2+\b)--(\a+14,-2+\b);
\draw[color=green](\a+2,-2+\b)--(\a+2,-14+\b);
\draw[color=green](\a+2,-4+\b)--(\a+14,-4+\b);
\draw[color=green](\a+4,-2+\b)--(\a+4,-14+\b);
\draw[color=green](\a+2,-6+\b)--(\a+12,-6+\b);
\draw[color=green](\a+6,-2+\b)--(\a+6,-12+\b);
\draw[color=green](\a+2,-8+\b)--(\a+10,-8+\b);
\draw[color=green](\a+8,-2+\b)--(\a+8,-10+\b);
\draw[color=green](\a+2,-10+\b)--(\a+8,-10+\b);
\draw[color=green](\a+10,-2+\b)--(\a+10,-8+\b);
\draw[color=green](\a+2,-12+\b)--(\a+6,-12+\b);
\draw[color=green](\a+12,-2+\b)--(\a+12,-6+\b);
\draw[color=green](\a+2,-14+\b)--(\a+4,-14+\b);
\draw[color=green](\a+14,-2+\b)--(\a+14,-4+\b);
\draw[very thick](\a+2,-3+\b)--(\a+4,-3+\b);
\draw[very thick](\a+3,-2+\b)--(\a+3,-4+\b);
\draw[very thick](\a+5,-2+\b)--(\a+5,-2.5+\b);
\draw[very thick](\a+5,-3.5+\b)--(\a+5,-4+\b);
\draw[very thick](\a+4,-3+\b)--(\a+4.5,-3+\b);
\draw[very thick](\a+5.5,-3+\b)--(\a+6,-3+\b);
\draw[very thick](\a+4.5,-3+\b)arc(270:360:0.5);
\draw[very thick](\a+5.5,-3+\b)arc(90:180:0.5);
\draw[very thick](\a+7,-2+\b)--(\a+7,-2.5+\b);
\draw[very thick](\a+7,-3.5+\b)--(\a+7,-4+\b);
\draw[very thick](\a+6,-3+\b)--(\a+6.5,-3+\b);
\draw[very thick](\a+7.5,-3+\b)--(\a+8,-3+\b);
\draw[very thick](\a+6.5,-3+\b)arc(270:360:0.5);
\draw[very thick](\a+7.5,-3+\b)arc(90:180:0.5);
\draw[very thick](\a+9,-2+\b)--(\a+9,-2.5+\b);
\draw[very thick](\a+9,-3.5+\b)--(\a+9,-4+\b);
\draw[very thick](\a+8,-3+\b)--(\a+8.5,-3+\b);
\draw[very thick](\a+9.5,-3+\b)--(\a+10,-3+\b);
\draw[very thick](\a+8.5,-3+\b)arc(270:360:0.5);
\draw[very thick](\a+9.5,-3+\b)arc(90:180:0.5);
\draw[very thick](\a+11,-2+\b)--(\a+11,-2.5+\b);
\draw[very thick](\a+11,-3.5+\b)--(\a+11,-4+\b);
\draw[very thick](\a+10,-3+\b)--(\a+10.5,-3+\b);
\draw[very thick](\a+11.5,-3+\b)--(\a+12,-3+\b);
\draw[very thick](\a+10.5,-3+\b)arc(270:360:0.5);
\draw[very thick](\a+11.5,-3+\b)arc(90:180:0.5);
\draw[very thick](\a+3,-4+\b)--(\a+3,-4.5+\b);
\draw[very thick](\a+3,-5.5+\b)--(\a+3,-6+\b);
\draw[very thick](\a+2,-5+\b)--(\a+2.5,-5+\b);
\draw[very thick](\a+3.5,-5+\b)--(\a+4,-5+\b);
\draw[very thick](\a+2.5,-5+\b)arc(270:360:0.5);
\draw[very thick](\a+3.5,-5+\b)arc(90:180:0.5);
\draw[very thick](\a+4,-5+\b)--(\a+6,-5+\b);
\draw[very thick](\a+5,-4+\b)--(\a+5,-6+\b);
\draw[very thick](\a+7,-4+\b)--(\a+7,-4.5+\b);
\draw[very thick](\a+7,-5.5+\b)--(\a+7,-6+\b);
\draw[very thick](\a+6,-5+\b)--(\a+6.5,-5+\b);
\draw[very thick](\a+7.5,-5+\b)--(\a+8,-5+\b);
\draw[very thick](\a+6.5,-5+\b)arc(270:360:0.5);
\draw[very thick](\a+7.5,-5+\b)arc(90:180:0.5);
\draw[very thick](\a+8,-5+\b)--(\a+10,-5+\b);
\draw[very thick](\a+9,-4+\b)--(\a+9,-6+\b);
\draw[very thick](\a+2,-7+\b)--(\a+4,-7+\b);
\draw[very thick](\a+3,-6+\b)--(\a+3,-8+\b);
\draw[very thick](\a+4,-7+\b)--(\a+6,-7+\b);
\draw[very thick](\a+5,-6+\b)--(\a+5,-8+\b);
\draw[very thick](\a+7,-6+\b)--(\a+7,-6.5+\b);
\draw[very thick](\a+7,-7.5+\b)--(\a+7,-8+\b);
\draw[very thick](\a+6,-7+\b)--(\a+6.5,-7+\b);
\draw[very thick](\a+7.5,-7+\b)--(\a+8,-7+\b);
\draw[very thick](\a+6.5,-7+\b)arc(270:360:0.5);
\draw[very thick](\a+7.5,-7+\b)arc(90:180:0.5);
\draw[very thick](\a+3,-8+\b)--(\a+3,-8.5+\b);
\draw[very thick](\a+3,-9.5+\b)--(\a+3,-10+\b);
\draw[very thick](\a+2,-9+\b)--(\a+2.5,-9+\b);
\draw[very thick](\a+3.5,-9+\b)--(\a+4,-9+\b);
\draw[very thick](\a+2.5,-9+\b)arc(270:360:0.5);
\draw[very thick](\a+3.5,-9+\b)arc(90:180:0.5);
\draw[very thick](\a+4,-9+\b)--(\a+6,-9+\b);
\draw[very thick](\a+5,-8+\b)--(\a+5,-10+\b);
\draw[very thick](\a+3,-10+\b)--(\a+3,-10.5+\b);
\draw[very thick](\a+3,-11.5+\b)--(\a+3,-12+\b);
\draw[very thick](\a+2,-11+\b)--(\a+2.5,-11+\b);
\draw[very thick](\a+3.5,-11+\b)--(\a+4,-11+\b);
\draw[very thick](\a+2.5,-11+\b)arc(270:360:0.5);
\draw[very thick](\a+3.5,-11+\b)arc(90:180:0.5);
\draw[very thick](\a+12,-3+\b)arc(270:360:1);
\draw[very thick](\a+10,-5+\b)arc(270:360:1);
\draw[very thick](\a+8,-7+\b)arc(270:360:1);
\draw[very thick](\a+6,-9+\b)arc(270:360:1);
\draw[very thick](\a+4,-11+\b)arc(270:360:1);
\draw[very thick](\a+2,-13+\b)arc(270:360:1);

\def\a{28};
\def\b{0};
\draw[color=green](\a+2,-2+\b)--(\a+14,-2+\b);
\draw[color=green](\a+2,-2+\b)--(\a+2,-14+\b);
\draw[color=green](\a+2,-4+\b)--(\a+14,-4+\b);
\draw[color=green](\a+4,-2+\b)--(\a+4,-14+\b);
\draw[color=green](\a+2,-6+\b)--(\a+12,-6+\b);
\draw[color=green](\a+6,-2+\b)--(\a+6,-12+\b);
\draw[color=green](\a+2,-8+\b)--(\a+10,-8+\b);
\draw[color=green](\a+8,-2+\b)--(\a+8,-10+\b);
\draw[color=green](\a+2,-10+\b)--(\a+8,-10+\b);
\draw[color=green](\a+10,-2+\b)--(\a+10,-8+\b);
\draw[color=green](\a+2,-12+\b)--(\a+6,-12+\b);
\draw[color=green](\a+12,-2+\b)--(\a+12,-6+\b);
\draw[color=green](\a+2,-14+\b)--(\a+4,-14+\b);
\draw[color=green](\a+14,-2+\b)--(\a+14,-4+\b);
\draw[very thick](\a+2,-3+\b)--(\a+4,-3+\b);
\draw[very thick](\a+3,-2+\b)--(\a+3,-4+\b);
\draw[very thick](\a+5,-2+\b)--(\a+5,-2.5+\b);
\draw[very thick](\a+5,-3.5+\b)--(\a+5,-4+\b);
\draw[very thick](\a+4,-3+\b)--(\a+4.5,-3+\b);
\draw[very thick](\a+5.5,-3+\b)--(\a+6,-3+\b);
\draw[very thick](\a+4.5,-3+\b)arc(270:360:0.5);
\draw[very thick](\a+5.5,-3+\b)arc(90:180:0.5);
\draw[very thick](\a+7,-2+\b)--(\a+7,-2.5+\b);
\draw[very thick](\a+7,-3.5+\b)--(\a+7,-4+\b);
\draw[very thick](\a+6,-3+\b)--(\a+6.5,-3+\b);
\draw[very thick](\a+7.5,-3+\b)--(\a+8,-3+\b);
\draw[very thick](\a+6.5,-3+\b)arc(270:360:0.5);
\draw[very thick](\a+7.5,-3+\b)arc(90:180:0.5);
\draw[very thick](\a+9,-2+\b)--(\a+9,-2.5+\b);
\draw[very thick](\a+9,-3.5+\b)--(\a+9,-4+\b);
\draw[very thick](\a+8,-3+\b)--(\a+8.5,-3+\b);
\draw[very thick](\a+9.5,-3+\b)--(\a+10,-3+\b);
\draw[very thick](\a+8.5,-3+\b)arc(270:360:0.5);
\draw[very thick](\a+9.5,-3+\b)arc(90:180:0.5);
\draw[very thick](\a+11,-2+\b)--(\a+11,-2.5+\b);
\draw[very thick](\a+11,-3.5+\b)--(\a+11,-4+\b);
\draw[very thick](\a+10,-3+\b)--(\a+10.5,-3+\b);
\draw[very thick](\a+11.5,-3+\b)--(\a+12,-3+\b);
\draw[very thick](\a+10.5,-3+\b)arc(270:360:0.5);
\draw[very thick](\a+11.5,-3+\b)arc(90:180:0.5);
\draw[very thick](\a+3,-4+\b)--(\a+3,-4.5+\b);
\draw[very thick](\a+3,-5.5+\b)--(\a+3,-6+\b);
\draw[very thick](\a+2,-5+\b)--(\a+2.5,-5+\b);
\draw[very thick](\a+3.5,-5+\b)--(\a+4,-5+\b);
\draw[very thick](\a+2.5,-5+\b)arc(270:360:0.5);
\draw[very thick](\a+3.5,-5+\b)arc(90:180:0.5);
\draw[very thick](\a+4,-5+\b)--(\a+6,-5+\b);
\draw[very thick](\a+5,-4+\b)--(\a+5,-6+\b);
\draw[very thick](\a+7,-4+\b)--(\a+7,-4.5+\b);
\draw[very thick](\a+7,-5.5+\b)--(\a+7,-6+\b);
\draw[very thick](\a+6,-5+\b)--(\a+6.5,-5+\b);
\draw[very thick](\a+7.5,-5+\b)--(\a+8,-5+\b);
\draw[very thick](\a+6.5,-5+\b)arc(270:360:0.5);
\draw[very thick](\a+7.5,-5+\b)arc(90:180:0.5);
\draw[very thick](\a+9,-4+\b)--(\a+9,-4.5+\b);
\draw[very thick](\a+9,-5.5+\b)--(\a+9,-6+\b);
\draw[very thick](\a+8,-5+\b)--(\a+8.5,-5+\b);
\draw[very thick](\a+9.5,-5+\b)--(\a+10,-5+\b);
\draw[very thick](\a+8.5,-5+\b)arc(270:360:0.5);
\draw[very thick](\a+9.5,-5+\b)arc(90:180:0.5);
\draw[very thick](\a+2,-7+\b)--(\a+4,-7+\b);
\draw[very thick](\a+3,-6+\b)--(\a+3,-8+\b);
\draw[very thick](\a+4,-7+\b)--(\a+6,-7+\b);
\draw[very thick](\a+5,-6+\b)--(\a+5,-8+\b);
\draw[very thick](\a+6,-7+\b)--(\a+8,-7+\b);
\draw[very thick](\a+7,-6+\b)--(\a+7,-8+\b);
\draw[very thick](\a+3,-8+\b)--(\a+3,-8.5+\b);
\draw[very thick](\a+3,-9.5+\b)--(\a+3,-10+\b);
\draw[very thick](\a+2,-9+\b)--(\a+2.5,-9+\b);
\draw[very thick](\a+3.5,-9+\b)--(\a+4,-9+\b);
\draw[very thick](\a+2.5,-9+\b)arc(270:360:0.5);
\draw[very thick](\a+3.5,-9+\b)arc(90:180:0.5);
\draw[very thick](\a+4,-9+\b)--(\a+6,-9+\b);
\draw[very thick](\a+5,-8+\b)--(\a+5,-10+\b);
\draw[very thick](\a+3,-10+\b)--(\a+3,-10.5+\b);
\draw[very thick](\a+3,-11.5+\b)--(\a+3,-12+\b);
\draw[very thick](\a+2,-11+\b)--(\a+2.5,-11+\b);
\draw[very thick](\a+3.5,-11+\b)--(\a+4,-11+\b);
\draw[very thick](\a+2.5,-11+\b)arc(270:360:0.5);
\draw[very thick](\a+3.5,-11+\b)arc(90:180:0.5);
\draw[very thick](\a+12,-3+\b)arc(270:360:1);
\draw[very thick](\a+10,-5+\b)arc(270:360:1);
\draw[very thick](\a+8,-7+\b)arc(270:360:1);
\draw[very thick](\a+6,-9+\b)arc(270:360:1);
\draw[very thick](\a+4,-11+\b)arc(270:360:1);
\draw[very thick](\a+2,-13+\b)arc(270:360:1);

\def\a{-30};
\def\b{-14};
\draw[step=2, color=green](\a+2,-14+\b)grid(\a+14,-2+\b);
\node at (\a+3,-3+\b) {$+$};
\node at (\a+3,-5+\b) {$+$};
\node at (\a+5,-5+\b) {$+$};
\node at (\a+7,-5+\b) {$+$};
\node at (\a+9,-5+\b) {$+$};
\node at (\a+5,-7+\b) {$+$};
\node at (\a+5,-9+\b) {$+$};
\node[above] at (\a+11,-1+\b) {$\mathrm{start}_2$};
\draw[->](\a+11,-0.5+\b)--(\a+11,-2+\b);

\def\a{0};
\def\b{-14};
\draw[step=2, color=green](\a+2,-14+\b)grid(\a+14,-2+\b);
\node at (\a+3,-3+\b) {$+$};
\node at (\a+5,-5+\b) {$+$};
\node at (\a+7,-5+\b) {$+$};
\node at (\a+9,-5+\b) {$+$};
\node at (\a+5,-7+\b) {$+$};
\node at (\a+5,-9+\b) {$+$};

\def\a{14};
\def\b{-14};
\draw[step=2, color=green](\a+2,-14+\b)grid(\a+14,-2+\b);
\node at (\a+3,-3+\b) {$+$};
\node at (\a+5,-5+\b) {$+$};
\node at (\a+3,-7+\b) {$+$};
\node at (\a+9,-5+\b) {$+$};
\node at (\a+5,-7+\b) {$+$};
\node at (\a+5,-9+\b) {$+$};

\def\a{28};
\def\b{-14};
\draw[step=2, color=green](\a+2,-14+\b)grid(\a+14,-2+\b);
\node at (\a+3,-3+\b) {$+$};
\node at (\a+5,-5+\b) {$+$};
\node at (\a+3,-7+\b) {$+$};
\node at (\a+7,-7+\b) {$+$};
\node at (\a+5,-7+\b) {$+$};
\node at (\a+5,-9+\b) {$+$};

\draw[very thick,->](-14,-15)--(-1,-15);
\node[above] at (-7.5,-15) {$\mathrm{mitosis}_2(D)$};
\node at (-18,-10) {$D$};
\node at (10,-14) {$D_{2,1}$};
\node at (24,-14) {$D_{2,3}$};
\node at (38,-14) {$D_{2,4}$};
\end{tikzpicture}
\caption{Example of the mitosis operator}
\label{fig:mitosis-example}
\end{figure}
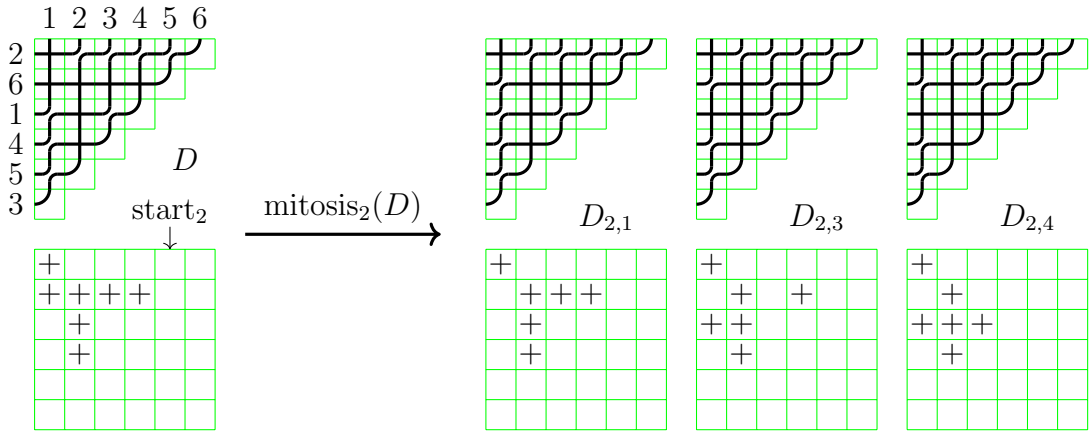

\noindent 
The monomial weight of the pipe dream $D$ is $x^D=x_1x_2^4x_3x_4$, and
we can calculate
that \[\partial_2x^D=x_1x_2^3x_3x_4+x_1x_2^2x_3^2x_4+x_1x_2x_3^3x_4=x^{D_{2,1}}+x^{D_{2,3}}+x^{D_{2,4}}.\]
This calculation may be a little bit misleading in general, but it is
the main idea of how mitosis is related to the divided difference
operators. 
\end{Example}

\begin{Exercise}
Let $D'$ be the pipe dream obtained from $D$ in
\Cref{fig:mitosis-example} by 3 chute moves from row 2 to row 1, so 
\[
D'=\{(1,1),(1,3),(1,4),(1,5),(2,1),(3,2),(4,2) \}.
\]
Compute $\partial_2x^{D'}$, $C_2(D')$, and $\mathrm{mitosis}_2(D')$.
How do $\partial_2x^D$ and $\mathrm{mitosis}_2(D')$ contribute to the
computation of $\fS_{261453s_{2}}$ in this case?
\end{Exercise}

\subsection{Balanced Tableaux and Balanced Labellings}\label{sub:Balanced}

In this section, we introduce the notion of \emph{balanced labellings}
of permutation diagrams, which also provides a combinatorial
interpretation for the Schubert polynomials, due to Fomin, Greene,
Reiner, Shimozono \cite{fomin1997balanced} . Balanced labellings can
be viewed as generalizations of \emph{balanced tableaux} used by
Edelman and Greene \cite{EG}. This celebrated paper \cite{EG},
originally aimed at solving enumeration problems related to
\emph{Stanley symmetric functions} \cite{Sta84}, introduced the
\emph{Edelman-Greene insertion algorithm} with great importance
towards many aspects of algebraic combinatorics. This section will not
focus on the insertion algorithm, but rather the more general setting
of balanced labellings.

\begin{Definition}\label{def:diagram}
A \emph{diagram} $D$ is a finite subset of boxes of the
$\mathbb{Z}_+\times\mathbb{Z}_+$ grid, drawn using matrix
coordinates.  
A \emph{tableau} (or a \emph{labelling}) of \emph{shape} $D$ is a
filling of the cells of $D$ by positive integers, or equivalently, a
map $D\rightarrow\mathbb{Z}_{>0}$. A tableau $T$ has \emph{content}
$\alpha=(\alpha_1,\alpha_2,\ldots)$ if there are $\alpha_i$ copies of
$i$, and we write $x^T=x_1^{\alpha_1}x_2^{\alpha_{2}}\cdots$.
\end{Definition}

The concept of a tableau of partition shape is prevalent in the
representation theory of the symmetric group \cite{Fulton-book,Sag}.
Permutation diagrams defined in \eqref{eq:diagram} fit this definition
of a diagram, so it is natural to consider fillings in that context as
well.  Using the crossings in a wiring diagram or pipe dream, they too
can be interpreted as diagrams.

The key definition from \cite{EG} is that of a \emph{balanced hook}.
``Balanced'' in this context refers to a stability property.

\begin{Definition}\label{def:hook}
To each cell $(i,j)$ of a diagram $D$, we associate the \emph{hook}
$H_{i,j}=H_{i,j}(D)$ consisting of cells $(i',j')$ of $D$ such that
either $i'=i$ and $j'\geq j$ or $i'\geq i$ and $j'=j$.  A labelling of the hook $H_{i,j}$ is \emph{balanced} if when one rearranges the labels within the hook so that the values increase weakly from right to left and from top to bottom, the corner entry stays the same.
\end{Definition}

\begin{Example}
On the left is a balanced hook, with the rearrangement of its values
into increasing order moving left and down shown on the right:
\[\ytableausetup{boxsize=1.3em}
\begin{ytableau}
3 & \none & 6 & 1 & 2 & \none & 2\\
1 \\
\none \\ 
5
\end{ytableau}\qquad\qquad
\begin{ytableau}
3 & \none & 2 & 2 & 1 & \none & 1\\
5 \\
\none \\
6
\end{ytableau}.\]
\end{Example}
\begin{Definition}\cite{EG,fomin1997balanced}\label{def:balanced-flagged-column-injective}
Let $D$ be a diagram with $\ell$ cells. A labelling of $D$ is
\begin{itemize}
\item \emph{balanced} if each hook $H_{i,j}(D)$ is balanced;
\item \emph{bijective} if each of the labels $1,2,\ldots,\ell$ appears exactly once;
\item \emph{column-injective} if there is at most one copy of $i$ in each column, for all $i$;
\item \emph{flagged} if all numbers in row $i$ do not exceed $i$, for all $i$.
\end{itemize}
\end{Definition}

\begin{figure}[h!]
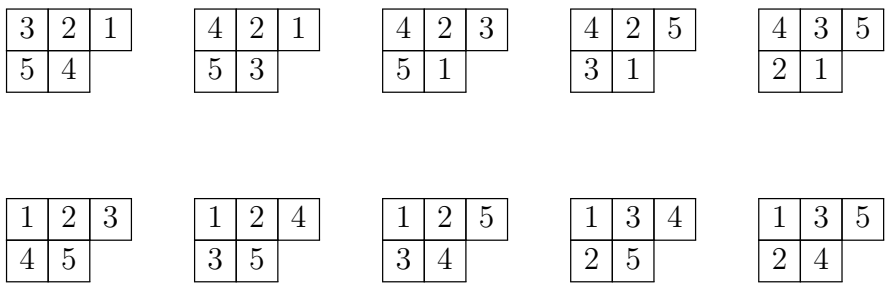

\centering
\[\ytableausetup{boxsize=1.3em}
\begin{ytableau}
3&2&1 \\ 5&4
\end{ytableau}\qquad
\begin{ytableau}
4&2&1\\5&3
\end{ytableau}\qquad
\begin{ytableau}
4&2&3\\5&1
\end{ytableau}\qquad
\begin{ytableau}
4&2&5\\3&1
\end{ytableau}\qquad
\begin{ytableau}
4&3&5\\2&1
\end{ytableau}\]

\

\

\[\ytableausetup{boxsize=1.3em}
\begin{ytableau}
1&2&3 \\ 4&5
\end{ytableau}\qquad
\begin{ytableau}
1&2&4\\3&5
\end{ytableau}\qquad
\begin{ytableau}
1&2&5\\3&4
\end{ytableau}\qquad
\begin{ytableau}
1&3&4\\2&5
\end{ytableau}\qquad
\begin{ytableau}
1&3&5\\2&4
\end{ytableau}\]
\caption{Balanced labellings (top) and standard Young tableaux (bottom) for $\lambda=(3,2)$}
\label{fig:BL-biject-SYT}
\end{figure}

Denote the set of bijective balanced labellings for a diagram $D$ as
$\mathrm{BBL}(D)$. An example of a bijective balanced labelling of a
diagram is shown in \Cref{fig:balanced-labellings}.  If $D=D(w)$ is
the Rothe diagram of a permutation $w$, we use the shorthand
$\mathrm{BBL}(w)$, and if $D$ is a partition shape $\lambda$, we use
$\mathrm{BBL}(\lambda)$. Similarly, denote the set of column-injective
flagged balanced labellings for a diagram $D$ as $\mathrm{CFBL}(D)$,
and write $\mathrm{CFBL}(w)$ for $\mathrm{CFBL}(D(w))$.  If $D$ is the
Young diagram of a partition, then the bijective balanced labellings
of $D$ have a nice bijection with standard tableaux.  See
Figure~\ref{fig:BL-biject-SYT} for an example of the bijection.

\begin{Theorem}\cite{EG}\label{thm:EG-balanced=standard}
Let $\lambda$ be a partition. Then $|\mathrm{BBL}(\lambda)|=|\mathrm{SYT}(\lambda)|$
where $\mathrm{SYT}$ is the set of standard Young tableaux of shape $\lambda$.
\end{Theorem} 

One intuition towards the balanced condition on permutation diagrams
is in connection with \emph{reflection orders}.  Recall that $R(w)$
denotes the set of reduced words for a permutation $w\in S_n$. For
$\bfa=(a_1,\ldots,a_{\ell})\in R(w)$, define an ordered pair
$\gamma_i=\big(s_{a_1}\cdots s_{a_{i-1}}(a_i),s_{a_1}\cdots
s_{a_{i-1}}(a_i+1)\big)$ which represents two values that are adjacent
and necessarily in increasing order in the intermediate permutation
$w^{(i-1)}=s_{a_1}\cdots s_{a_{i-1}}$.  In $w^{(i)}=s_{a_1}\cdots
s_{a_{i-1}}s_{a_{i}}$ these two values are switched.  Since $\bfa$ is
reduced, these two values will not switch again in the successive
intermediate permutations $w^{(i+1)}, \ldots, w^{(\ell)}=w$ so they
correspond with an inversion in $w$, but written in terms of values
instead of positions.  The set $\{(w_{j},w_{i}) \given i<j,
w_{i}>w_{j}\}$ is the inversion set of $w^{-1}$ as defined in
\eqref{eq:inversion.def}, and there is an obvious bijection with
$\Inv(w)$.

This sequence $\gamma_1,\ldots,\gamma_{\ell}$ is a total order of the
pairs $\{(w_{j},w_{i}) \given i<j, w_{i}>w_{j}\}$.  By considering the
wiring diagram of $\bfa$, one can observe that if there exists a
\textit{321-pattern} $i<j<k$ such that $w_{i}>w_{j}>w_{k}$ then
$(w_{k},w_{i})$ must appear between $(w_{j},w_{i})$ and
$(w_{k},w_{j})$ in the total order $\gamma_1,\ldots,\gamma_{\ell}$.
Conversely, one can show that any total order on $\{(w_{j},w_{i})
\given i<j, w_{i}>w_{j}\}$ that can be extended to a total order on
all pairs $(i,j)$ for $1 \leq i<j\leq n$ satisfying the rule for all
$321$-patterns corresponds with some reduced word for $w$. Such an
order is called a \emph{reflection order}.  See \cite{b-b} for more
background on reflection orders.  The point here is that balanced
labellings are another way of recording reduced words, reflection
orders, and the condition on $321$-patterns.

\begin{Definition}\cite{fomin1997balanced}\label{def:canonical-labelling}
For a reduced word $\bfa=(a_1,\ldots,a_{\ell})\in R(w)$, define a tableau
$T_{\bfa}$ of shape $D(w)$ such that $T_{\bfa}(i,j)=k$ if the simple
transposition $s_{a_k}$ switches the values $j<w_i$ in the
intermediate permutation $w^{(k-1)}=s_{a_1}\cdots s_{a_{k-1}}$.
\end{Definition}

\begin{Exercise}\label{ex:diagram.balanced}
Given the definition of $D(w)$ in \eqref{eq:diagram}, show each $(i,j)
\in D(w)$ is assigned a value by the map $T_{\bfa}$.
\end{Exercise}

\begin{Theorem}\cite{fomin1997balanced}\label{thm:caonical-labelling}
The map $a\mapsto T_{\bfa}$ is a bijection between $R(w)$ and $\mathrm{BBL}(w)$. 
\end{Theorem} 

Edelman and Greene used \Cref{def:canonical-labelling} and
\Cref{thm:caonical-labelling} in the case where $w_0=[n,\ldots,1]$ is
the longest permutation in $S_n$, whose Rothe diagram is the staircase
partition shape $(n{-}1,n{-}2,\ldots,1)$. Together with the celebrated
Edelman-Greene insertion algorithm which provides a bijection between
$R(w_0)$ and the standard Young tableaux of the staircase shape, they
established \Cref{thm:EG-balanced=standard} for the staircase shape,
which serves as the base case for any partition shape $\lambda$. The
current presentation of \Cref{def:canonical-labelling} and
\Cref{thm:caonical-labelling} is due to Fomin-Greene-Reiner-Shimozono
\cite{fomin1997balanced}.

\begin{Exercise}\label{ex:canonical-labelling}
Prove \Cref{thm:caonical-labelling} by constructing the inverse map.
\end{Exercise}

\begin{Example}

Consider $w=43152$ and its reduced word $\bfa=(2,3,2,1,4,2)\in R(w)$. The procedure of
generating $T_{\bfa}$ and the final result is shown in \Cref{fig:balanced-labellings}. One can indeed check that $T_{\bfa}$ is balanced.

\begin{figure}[h!]
\begin{minipage}{0.45\linewidth}
\centering
\begin{tabular}{c|c|c|c} 
$k$ & $a_k$ & $w^{(k)}$ & \text{cell} \\\hline
$1$ & $2$ & $1\underline{\textbf{32}}45$ & $(2,2)$ \\
$2$ & $3$ & $13\underline{\textbf{42}}5$ & $(1,2)$ \\
$3$ & $2$ & $1\underline{\textbf{43}}25$ & $(1,3)$ \\
$4$ & $1$ & $\underline{\textbf{41}}325$ & $(1,1)$ \\
$5$ & $4$ & $413\underline{\textbf{52}}$ & $(4,2)$ \\
$6$ & $2$ & $4\underline{\textbf{31}}52$ & $(2,1)$
\end{tabular}
\end{minipage}
\hspace{0.25cm}
\begin{minipage}{0.45\linewidth}
\centering
\begin{tikzpicture}[scale=0.6]
\draw[step=1.0,green,thin] (0,0) grid (5,5);
\draw[very thick] (1,1)--(2,1)--(2,2)--(1,2)--(1,1);
\draw[very thick] (0,3)--(2,3)--(2,4)--(3,4)--(3,5)--(0,5)--(0,3);
\node at (0.5,2.5) {$\bullet$};
\node at (1.5,0.5) {$\bullet$};
\node at (2.5,3.5) {$\bullet$};
\node at (3.5,4.5) {$\bullet$};
\node at (4.5,1.5) {$\bullet$};
\draw(0.5,-0.5)--(0.5,2.5)--(5.5,2.5);
\draw(1.5,-0.5)--(1.5,0.5)--(5.5,0.5);
\draw(2.5,-0.5)--(2.5,3.5)--(5.5,3.5);
\draw(3.5,-0.5)--(3.5,4.5)--(5.5,4.5);
\draw(4.5,-0.5)--(4.5,1.5)--(5.5,1.5);
\node at (1.5,1.5) {$5$};
\node at (0.5,3.5) {$6$};
\node at (0.5,4.5) {$4$};
\node at (1.5,3.5) {$1$};
\node at (1.5,4.5) {$2$};
\node at (2.5,4.5) {$3$};
\end{tikzpicture}
\end{minipage}
\caption{Construction of the tableau $T_{\bfa}$ for  $\bfa=(2,3,2,1,4,2)\in R(43152)$}
\label{fig:balanced-labellings}
\end{figure}
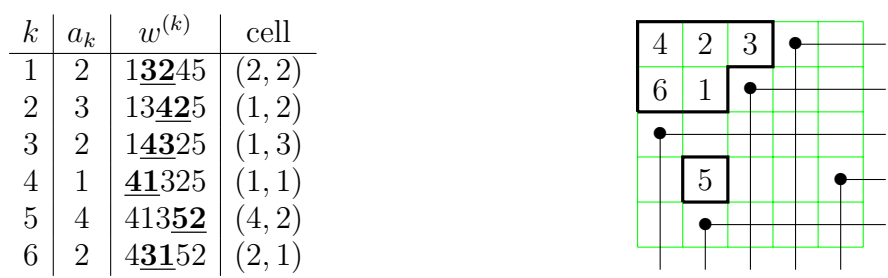

\end{Example}
Schubert polynomials can also be expanded using the idea of balanced
labellings. To be precise, Schubert polynomials are generating
functions for column-injective flagged balanced labellings of
permutation diagrams.
\begin{Theorem}\cite{fomin1997balanced}\label{thm:balanced-schubert}
For a permutation $w$, $\fS_w=\sum_{T\in \mathrm{CFBL}(w)}x^T$.
\end{Theorem}
\begin{Example}
Consider $w=1432$ with $\fS_w=x_1^2x_2+x_1^2x_3+x_1x_2^2+x_1x_2x_3+x_2^2x_3$. The column-injective flagged balanced labellings of shape $D(w)$ are shown in \Cref{fig:column-injective-balanced-labellings}.
\begin{figure}[h!]
\centering
\begin{tikzpicture}[scale=0.5]
\draw[step=1.0,green,thin] (0,0) grid (3,3);
\draw[very thick] (1,2)--(3,2)--(3,1)--(1,1);
\draw[very thick] (1,2)--(1,0)--(2,0)--(2,2);
\node at (1.5,1.5) {$1$};
\node at (2.5,1.5) {$1$};
\node at (1.5,0.5) {$2$};
\end{tikzpicture}
\quad
\begin{tikzpicture}[scale=0.5]
\draw[step=1.0,green,thin] (0,0) grid (3,3);
\draw[very thick] (1,2)--(3,2)--(3,1)--(1,1);
\draw[very thick] (1,2)--(1,0)--(2,0)--(2,2);
\node at (1.5,1.5) {$1$};
\node at (2.5,1.5) {$1$};
\node at (1.5,0.5) {$3$};
\end{tikzpicture}
\quad
\begin{tikzpicture}[scale=0.5]
\draw[step=1.0,green,thin] (0,0) grid (3,3);
\draw[very thick] (1,2)--(3,2)--(3,1)--(1,1);
\draw[very thick] (1,2)--(1,0)--(2,0)--(2,2);
\node at (1.5,1.5) {$2$};
\node at (2.5,1.5) {$1$};
\node at (1.5,0.5) {$3$};
\end{tikzpicture}
\quad
\begin{tikzpicture}[scale=0.5]
\draw[step=1.0,green,thin] (0,0) grid (3,3);
\draw[very thick] (1,2)--(3,2)--(3,1)--(1,1);
\draw[very thick] (1,2)--(1,0)--(2,0)--(2,2);
\node at (1.5,1.5) {$2$};
\node at (2.5,1.5) {$2$};
\node at (1.5,0.5) {$1$};
\end{tikzpicture}
\quad
\begin{tikzpicture}[scale=0.5]
\draw[step=1.0,green,thin] (0,0) grid (3,3);
\draw[very thick] (1,2)--(3,2)--(3,1)--(1,1);
\draw[very thick] (1,2)--(1,0)--(2,0)--(2,2);
\node at (1.5,1.5) {$2$};
\node at (2.5,1.5) {$2$};
\node at (1.5,0.5) {$3$};
\end{tikzpicture}
\caption{Column-injective flagged balanced labellings of shape $D(1432)$}
\label{fig:column-injective-balanced-labellings}
\end{figure}
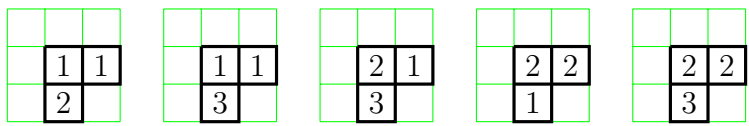
\end{Example}

\subsection{Bumpless Pipe Dreams}\label{sub:Bumpless}
Lam, Lee, and Shimozono \cite{LLS} introduced \emph{bumpless pipe dreams} (BPDs) in their work on the infinite flag variety and back-stable Schubert calculus and used them to give a formula for (double) Schubert polynomials. Since then, the rich combinatorial, algebraic and geometric properties of bumpless pipe dreams have been extensively explored by the community, and it appears that there are still many exciting developments on the horizon.

\begin{Definition}
A \emph{bumpless pipe dream} (abbreviated as $\BPD$) $D$ is a tiling
of the $\mathbb{Z}_+\times\mathbb{Z}_+$ grid with matrix coordinates
using the following six types of tiles:
\begin{center}
\includegraphics[scale=0.5]{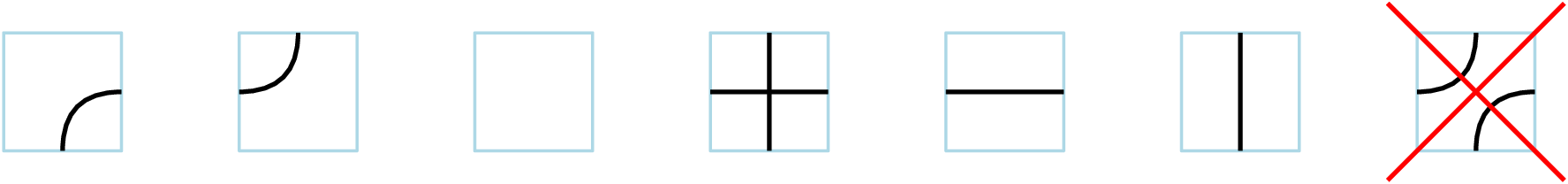}
\end{center}
forming pipes that travel in the northeast direction. A bumpless pipe
dream is \emph{reduced} if no two pipes cross twice. A reduced
bumpless pipe dream corresponds to a permutation $w\in S_{\infty}$ if
pipe $i$ goes from $(\infty,i)$ to $(w(i),\infty)$ for all $i$. Let
$\mathrm{blank}(D)$ and $\mathrm{cross}(D)$ be the coordinates of
\emph{blank tiles} \bl\ and \emph{cross tiles} \+\ in a bumpless pipe
dream $D$, respectively. The weight of a bumpless pipe dream $D$ is
given by the product of variables indexed by the row numbers of its blank tiles \[x^D:=\prod_{(i,j)\in\mathrm{blank}(D)}x_i.\] For a permutation $w$, let $\BPD(w)$ be the set of reduced bumpless pipe dreams of $w$.
\end{Definition}
The name ``bumpless" comes from the fact that the \emph{bump tile}
$\elbow$ that was used to construct the original pipe dreams is not
allowed in such a tiling. Note that unlike the original pipe dreams, a
bumpless pipe dream is not determined by either $\mathrm{blank}(D)$ or
$\mathrm{cross}(D)$. For a permutation $w\in S_n$, we typically draw a
bumpless pipe dream of $w$ in a $n\times n$ grid.   See \Cref{fig:BPDs-2143} for
example. 

\begin{Exercise}
Show that in a bumpless pipe dream $D$, we have $|\mathrm{blank}(D)|=|\mathrm{cross}(D)|$.
\end{Exercise}

\begin{Exercise}\label{exer:BPD-blank-cross}
Show that a bumpless pipe dream $D$ is determined by $\mathrm{blank}(D)$ and $\mathrm{cross}(D)$ together.
\end{Exercise}

Remarkably, Lam-Lee-Shimozono proved that the
reduced bumpless pipe dreams for $w$ enumerate monomials in the
corresponding Schubert polynomial, just like the classic pipe dreams.
Their proof was via considerations on the cohomology of the flag
variety and special positroid varieties known as \emph{graph
Schubert varieties}, not via an explicit bijection.  See
\Cref{exer:graph-schubert} below for more on these varieties.

\begin{Theorem}\cite{LLS}\label{thm:bumpless-schubert} For all $w
\in S_{n}$, the Schubert polynomial satisfies
\[\fS_{w}(x_1,x_2,\ldots,x_{n})=\sum_{D\in\BPD(w)}x^D.\]
\end{Theorem}
\begin{Example}\label{ex:PD-BPD-2143}
We have $\fS_{2143}=x_1^2+x_1x_2+x_1x_3$ with the following reduced
bumpless pipe dreams shown in \Cref{fig:BPDs-2143}. We also provide
the pipe dreams in \Cref{fig:PDs-2143}.  Try to identify a natural
bijection between these two sets.
\begin{figure}[h!]
\centering
\begin{tikzpicture}[scale=0.300000000000000]
\draw[step=2, color=green](2,-10)grid(10,-2);

\draw[very thick](5,-3.50000000000000)--(5,-4);
\draw[very thick](5.50000000000000,-3)--(6,-3);
\draw[very thick](5.50000000000000,-3)arc(90:180:0.500000000000000);

\draw[very thick](6,-3)--(8,-3);

\draw[very thick](8,-3)--(10,-3);

\draw[very thick](3,-5.50000000000000)--(3,-6);
\draw[very thick](3.50000000000000,-5)--(4,-5);
\draw[very thick](3.50000000000000,-5)arc(90:180:0.500000000000000);

\draw[very thick](4,-5)--(6,-5);
\draw[very thick](5,-4)--(5,-6);

\draw[very thick](6,-5)--(8,-5);

\draw[very thick](8,-5)--(10,-5);

\draw[very thick](3,-6)--(3,-8);

\draw[very thick](5,-6)--(5,-8);

\draw[very thick](9,-7.50000000000000)--(9,-8);
\draw[very thick](9.50000000000000,-7)--(10,-7);
\draw[very thick](9.50000000000000,-7)arc(90:180:0.500000000000000);

\draw[very thick](3,-8)--(3,-10);

\draw[very thick](5,-8)--(5,-10);

\draw[very thick](7,-9.50000000000000)--(7,-10);
\draw[very thick](7.50000000000000,-9)--(8,-9);
\draw[very thick](7.50000000000000,-9)arc(90:180:0.500000000000000);

\draw[very thick](8,-9)--(10,-9);
\draw[very thick](9,-8)--(9,-10);

\node[right] at (10,-3) {$2$};
\node[right] at (10,-5) {$1$};
\node[right] at (10,-7) {$4$};
\node[right] at (10,-9) {$3$};
\node[below] at (3,-10) {$1$};
\node[below] at (5,-10) {$2$};
\node[below] at (7,-10) {$3$};
\node[below] at (9,-10) {$4$};
\end{tikzpicture}
\quad
\begin{tikzpicture}[scale=0.300000000000000]
\draw[step=2, color=green](2,-10)grid(10,-2);

\draw[very thick](5,-3.50000000000000)--(5,-4);
\draw[very thick](5.50000000000000,-3)--(6,-3);
\draw[very thick](5.50000000000000,-3)arc(90:180:0.500000000000000);

\draw[very thick](6,-3)--(8,-3);

\draw[very thick](8,-3)--(10,-3);

\draw[very thick](5,-4)--(5,-6);

\draw[very thick](7,-5.50000000000000)--(7,-6);
\draw[very thick](7.50000000000000,-5)--(8,-5);
\draw[very thick](7.50000000000000,-5)arc(90:180:0.500000000000000);

\draw[very thick](8,-5)--(10,-5);

\draw[very thick](3,-7.50000000000000)--(3,-8);
\draw[very thick](3.50000000000000,-7)--(4,-7);
\draw[very thick](3.50000000000000,-7)arc(90:180:0.500000000000000);

\draw[very thick](4,-7)--(6,-7);
\draw[very thick](5,-6)--(5,-8);

\draw[very thick](7,-6)--(7,-6.50000000000000);
\draw[very thick](6,-7)--(6.50000000000000,-7);
\draw[very thick](6.50000000000000,-7)arc(270:360:0.500000000000000);

\draw[very thick](9,-7.50000000000000)--(9,-8);
\draw[very thick](9.50000000000000,-7)--(10,-7);
\draw[very thick](9.50000000000000,-7)arc(90:180:0.500000000000000);

\draw[very thick](3,-8)--(3,-10);

\draw[very thick](5,-8)--(5,-10);

\draw[very thick](7,-9.50000000000000)--(7,-10);
\draw[very thick](7.50000000000000,-9)--(8,-9);
\draw[very thick](7.50000000000000,-9)arc(90:180:0.500000000000000);

\draw[very thick](8,-9)--(10,-9);
\draw[very thick](9,-8)--(9,-10);

\node[right] at (10,-3) {$2$};
\node[right] at (10,-5) {$1$};
\node[right] at (10,-7) {$4$};
\node[right] at (10,-9) {$3$};
\node[below] at (3,-10) {$1$};
\node[below] at (5,-10) {$2$};
\node[below] at (7,-10) {$3$};
\node[below] at (9,-10) {$4$};
\end{tikzpicture}
\quad
\begin{tikzpicture}[scale=0.300000000000000]
\draw[step=2, color=green](2,-10)grid(10,-2);

\draw[very thick](7,-3.50000000000000)--(7,-4);
\draw[very thick](7.50000000000000,-3)--(8,-3);
\draw[very thick](7.50000000000000,-3)arc(90:180:0.500000000000000);

\draw[very thick](8,-3)--(10,-3);

\draw[very thick](3,-5.50000000000000)--(3,-6);
\draw[very thick](3.50000000000000,-5)--(4,-5);
\draw[very thick](3.50000000000000,-5)arc(90:180:0.500000000000000);

\draw[very thick](4,-5)--(6,-5);

\draw[very thick](6,-5)--(8,-5);
\draw[very thick](7,-4)--(7,-6);

\draw[very thick](8,-5)--(10,-5);

\draw[very thick](3,-6)--(3,-8);

\draw[very thick](5,-7.50000000000000)--(5,-8);
\draw[very thick](5.50000000000000,-7)--(6,-7);
\draw[very thick](5.50000000000000,-7)arc(90:180:0.500000000000000);

\draw[very thick](7,-6)--(7,-6.50000000000000);
\draw[very thick](6,-7)--(6.50000000000000,-7);
\draw[very thick](6.50000000000000,-7)arc(270:360:0.500000000000000);

\draw[very thick](9,-7.50000000000000)--(9,-8);
\draw[very thick](9.50000000000000,-7)--(10,-7);
\draw[very thick](9.50000000000000,-7)arc(90:180:0.500000000000000);

\draw[very thick](3,-8)--(3,-10);

\draw[very thick](5,-8)--(5,-10);

\draw[very thick](7,-9.50000000000000)--(7,-10);
\draw[very thick](7.50000000000000,-9)--(8,-9);
\draw[very thick](7.50000000000000,-9)arc(90:180:0.500000000000000);

\draw[very thick](8,-9)--(10,-9);
\draw[very thick](9,-8)--(9,-10);

\node[right] at (10,-3) {$2$};
\node[right] at (10,-5) {$1$};
\node[right] at (10,-7) {$4$};
\node[right] at (10,-9) {$3$};
\node[below] at (3,-10) {$1$};
\node[below] at (5,-10) {$2$};
\node[below] at (7,-10) {$3$};
\node[below] at (9,-10) {$4$};
\end{tikzpicture}
\caption{Reduced bumpless pipe dreams for $2143$.}
\label{fig:BPDs-2143}
\end{figure}

\begin{figure}[h!]
\centering
\begin{tikzpicture}[scale=0.300000000000000]
\draw[color=green](2,-2)--(10,-2);
\draw[color=green](2,-2)--(2,-10);
\draw[color=green](2,-4)--(10,-4);
\draw[color=green](4,-2)--(4,-10);
\draw[color=green](2,-6)--(8,-6);
\draw[color=green](6,-2)--(6,-8);
\draw[color=green](2,-8)--(6,-8);
\draw[color=green](8,-2)--(8,-6);
\draw[color=green](2,-10)--(4,-10);
\draw[color=green](10,-2)--(10,-4);
\draw[very thick](2,-3)--(4,-3);
\draw[very thick](3,-2)--(3,-4);
\draw[very thick](5,-2)--(5,-2.50000000000000);
\draw[very thick](5,-3.50000000000000)--(5,-4);
\draw[very thick](4,-3)--(4.50000000000000,-3);
\draw[very thick](5.50000000000000,-3)--(6,-3);
\draw[very thick](4.50000000000000,-3)arc(270:360:0.500000000000000);
\draw[very thick](5.50000000000000,-3)arc(90:180:0.500000000000000);
\draw[very thick](7,-2)--(7,-2.50000000000000);
\draw[very thick](7,-3.50000000000000)--(7,-4);
\draw[very thick](6,-3)--(6.50000000000000,-3);
\draw[very thick](7.50000000000000,-3)--(8,-3);
\draw[very thick](6.50000000000000,-3)arc(270:360:0.500000000000000);
\draw[very thick](7.50000000000000,-3)arc(90:180:0.500000000000000);
\draw[very thick](3,-4)--(3,-4.50000000000000);
\draw[very thick](3,-5.50000000000000)--(3,-6);
\draw[very thick](2,-5)--(2.50000000000000,-5);
\draw[very thick](3.50000000000000,-5)--(4,-5);
\draw[very thick](2.50000000000000,-5)arc(270:360:0.500000000000000);
\draw[very thick](3.50000000000000,-5)arc(90:180:0.500000000000000);
\draw[very thick](5,-4)--(5,-4.50000000000000);
\draw[very thick](5,-5.50000000000000)--(5,-6);
\draw[very thick](4,-5)--(4.50000000000000,-5);
\draw[very thick](5.50000000000000,-5)--(6,-5);
\draw[very thick](4.50000000000000,-5)arc(270:360:0.500000000000000);
\draw[very thick](5.50000000000000,-5)arc(90:180:0.500000000000000);
\draw[very thick](2,-7)--(4,-7);
\draw[very thick](3,-6)--(3,-8);
\draw[very thick](8,-3)arc(270:360:1);
\draw[very thick](6,-5)arc(270:360:1);
\draw[very thick](4,-7)arc(270:360:1);
\draw[very thick](2,-9)arc(270:360:1);
\node[left] at (2,-3) {$2$};
\node[left] at (2,-5) {$1$};
\node[left] at (2,-7) {$4$};
\node[left] at (2,-9) {$3$};
\node[above] at (3,-2) {$1$};
\node[above] at (5,-2) {$2$};
\node[above] at (7,-2) {$3$};
\node[above] at (9,-2) {$4$};
\end{tikzpicture}
\quad
\begin{tikzpicture}[scale=0.300000000000000]
\draw[color=green](2,-2)--(10,-2);
\draw[color=green](2,-2)--(2,-10);
\draw[color=green](2,-4)--(10,-4);
\draw[color=green](4,-2)--(4,-10);
\draw[color=green](2,-6)--(8,-6);
\draw[color=green](6,-2)--(6,-8);
\draw[color=green](2,-8)--(6,-8);
\draw[color=green](8,-2)--(8,-6);
\draw[color=green](2,-10)--(4,-10);
\draw[color=green](10,-2)--(10,-4);
\draw[very thick](2,-3)--(4,-3);
\draw[very thick](3,-2)--(3,-4);
\draw[very thick](5,-2)--(5,-2.50000000000000);
\draw[very thick](5,-3.50000000000000)--(5,-4);
\draw[very thick](4,-3)--(4.50000000000000,-3);
\draw[very thick](5.50000000000000,-3)--(6,-3);
\draw[very thick](4.50000000000000,-3)arc(270:360:0.500000000000000);
\draw[very thick](5.50000000000000,-3)arc(90:180:0.500000000000000);
\draw[very thick](7,-2)--(7,-2.50000000000000);
\draw[very thick](7,-3.50000000000000)--(7,-4);
\draw[very thick](6,-3)--(6.50000000000000,-3);
\draw[very thick](7.50000000000000,-3)--(8,-3);
\draw[very thick](6.50000000000000,-3)arc(270:360:0.500000000000000);
\draw[very thick](7.50000000000000,-3)arc(90:180:0.500000000000000);
\draw[very thick](3,-4)--(3,-4.50000000000000);
\draw[very thick](3,-5.50000000000000)--(3,-6);
\draw[very thick](2,-5)--(2.50000000000000,-5);
\draw[very thick](3.50000000000000,-5)--(4,-5);
\draw[very thick](2.50000000000000,-5)arc(270:360:0.500000000000000);
\draw[very thick](3.50000000000000,-5)arc(90:180:0.500000000000000);
\draw[very thick](4,-5)--(6,-5);
\draw[very thick](5,-4)--(5,-6);
\draw[very thick](3,-6)--(3,-6.50000000000000);
\draw[very thick](3,-7.50000000000000)--(3,-8);
\draw[very thick](2,-7)--(2.50000000000000,-7);
\draw[very thick](3.50000000000000,-7)--(4,-7);
\draw[very thick](2.50000000000000,-7)arc(270:360:0.500000000000000);
\draw[very thick](3.50000000000000,-7)arc(90:180:0.500000000000000);
\draw[very thick](8,-3)arc(270:360:1);
\draw[very thick](6,-5)arc(270:360:1);
\draw[very thick](4,-7)arc(270:360:1);
\draw[very thick](2,-9)arc(270:360:1);
\node[left] at (2,-3) {$2$};
\node[left] at (2,-5) {$1$};
\node[left] at (2,-7) {$4$};
\node[left] at (2,-9) {$3$};
\node[above] at (3,-2) {$1$};
\node[above] at (5,-2) {$2$};
\node[above] at (7,-2) {$3$};
\node[above] at (9,-2) {$4$};
\end{tikzpicture}
\quad
\begin{tikzpicture}[scale=0.300000000000000]
\draw[color=green](2,-2)--(10,-2);
\draw[color=green](2,-2)--(2,-10);
\draw[color=green](2,-4)--(10,-4);
\draw[color=green](4,-2)--(4,-10);
\draw[color=green](2,-6)--(8,-6);
\draw[color=green](6,-2)--(6,-8);
\draw[color=green](2,-8)--(6,-8);
\draw[color=green](8,-2)--(8,-6);
\draw[color=green](2,-10)--(4,-10);
\draw[color=green](10,-2)--(10,-4);
\draw[very thick](2,-3)--(4,-3);
\draw[very thick](3,-2)--(3,-4);
\draw[very thick](5,-2)--(5,-2.50000000000000);
\draw[very thick](5,-3.50000000000000)--(5,-4);
\draw[very thick](4,-3)--(4.50000000000000,-3);
\draw[very thick](5.50000000000000,-3)--(6,-3);
\draw[very thick](4.50000000000000,-3)arc(270:360:0.500000000000000);
\draw[very thick](5.50000000000000,-3)arc(90:180:0.500000000000000);
\draw[very thick](6,-3)--(8,-3);
\draw[very thick](7,-2)--(7,-4);
\draw[very thick](3,-4)--(3,-4.50000000000000);
\draw[very thick](3,-5.50000000000000)--(3,-6);
\draw[very thick](2,-5)--(2.50000000000000,-5);
\draw[very thick](3.50000000000000,-5)--(4,-5);
\draw[very thick](2.50000000000000,-5)arc(270:360:0.500000000000000);
\draw[very thick](3.50000000000000,-5)arc(90:180:0.500000000000000);
\draw[very thick](5,-4)--(5,-4.50000000000000);
\draw[very thick](5,-5.50000000000000)--(5,-6);
\draw[very thick](4,-5)--(4.50000000000000,-5);
\draw[very thick](5.50000000000000,-5)--(6,-5);
\draw[very thick](4.50000000000000,-5)arc(270:360:0.500000000000000);
\draw[very thick](5.50000000000000,-5)arc(90:180:0.500000000000000);
\draw[very thick](3,-6)--(3,-6.50000000000000);
\draw[very thick](3,-7.50000000000000)--(3,-8);
\draw[very thick](2,-7)--(2.50000000000000,-7);
\draw[very thick](3.50000000000000,-7)--(4,-7);
\draw[very thick](2.50000000000000,-7)arc(270:360:0.500000000000000);
\draw[very thick](3.50000000000000,-7)arc(90:180:0.500000000000000);
\draw[very thick](8,-3)arc(270:360:1);
\draw[very thick](6,-5)arc(270:360:1);
\draw[very thick](4,-7)arc(270:360:1);
\draw[very thick](2,-9)arc(270:360:1);
\node[left] at (2,-3) {$2$};
\node[left] at (2,-5) {$1$};
\node[left] at (2,-7) {$4$};
\node[left] at (2,-9) {$3$};
\node[above] at (3,-2) {$1$};
\node[above] at (5,-2) {$2$};
\node[above] at (7,-2) {$3$};
\node[above] at (9,-2) {$4$};
\end{tikzpicture}
\caption{Reduced pipe dreams for $2143$.}
\label{fig:PDs-2143}
\end{figure}
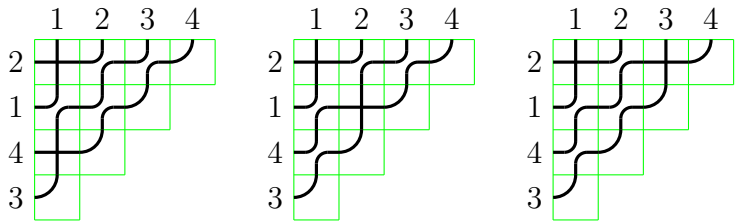
\end{Example}

\begin{Remark}
It is shown in \cite{LLS} that bumpless pipe dreams can also compute double Schubert polynomials. After adding variables $y_1,y_2,\ldots$, the weight of a bumpless pipe dream $D$ becomes \[(x-y)^D:=\prod_{(i,j)\in\mathrm{blank}(D)}(x_i-y_j).\] Lam-Lee-Shimozono showed that
\begin{equation}
 \fS_{w}(x_{1},\dots ,x_{n};y_{1},\dots ,
y_{n})= \sum_{D \in \BPD(w)} (x-y)^D.
\end{equation}
Although both pipe dream models provide
combinatorial formulae for double Schubert polynomials, the
factored sums are very different. Continuing \Cref{ex:PD-BPD-2143},
we observe from the  bumpless pipe dreams in \Cref{fig:BPDs-2143} that 
\[\fS_{2143}(X;Y)=(x_1-y_1)(x_3-y_3)+(x_1-y_1)(x_2-y_1)+(x_1-y_1)(x_1-y_2),\]
while from the pipe dreams in \Cref{fig:PDs-2143} and
\eqref{eq:double.schubs}, we have 
\[\fS_{2143}(X;Y)=(x_1-y_1)(x_1-y_3)+(x_1-y_1)(x_2-y_2)+(x_1-y_1)(x_3-y_1).\]
Hence, $\fS_{2143}(X;Y)$ can be expressed as a sum of products of
differences in different ways.
\end{Remark}

Weigandt \cite{weigandt-BPD-ASM} observed that the (not necessarily
reduced) bumpless pipe dreams on the $n \times n$ grid are in natural bijection
with \emph{alternating sign matrices}  of size $n$.  Let's
recall the definition of an ASM.

\begin{Definition}
An \emph{alternating sign matrix} (ASM) of size $n$ is an $n\times n$ matrix with entries in $\{-1,0,1\}$ such that within each row and each column, the nonzero entries sum up to $1$ and alternate between $1$ and $-1$. 
\end{Definition}
There is an incredibly rich literature and history around ASMs in
enumerative and algebraic combinatorics. Famously, the number of ASMs
of size $n$ has a beautiful product formula,
\[
\prod_{k=0}^{n-1} \frac{(3k+1)!}{(n+k)!}.
\]
The formula was conjectured by Mills-Robbins-Rumsey in 1983
\cite{MRR83} and proved in 1992 by Zeilberger \cite{zeilberger-ASM},
connected to the six vertex state model by Kuperberg
\cite{kuperberg-ASM}, and an operator formula on monotone triangles by
Fischer \cite{Fischer.2007}. On the other hand, the number of ASMs
of size $n$ equals the number of \emph{totally symmetric
self-complementary plane partitions} (TSSCPP) of size $n$
\cite{andrews-TSSCPP}. It is still open to find an explicit bijection
between ASMs and TSSCPPs, while some progress has been made with the
perspective of bumpless pipe dreams \cite{huang-TSSCPP}. The ASMs also
arise as the \emph{Dedekind–MacNeille completion} of the Bruhat order
to a lattice \cite{LS-MacNeille}.

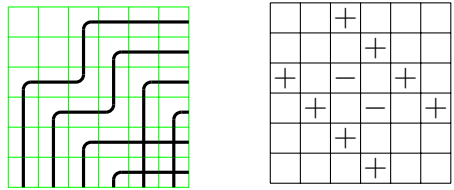
\begin{figure}[h!]
\centering
\begin{tikzpicture}[scale=0.200000000000000]
\def\a{0};
\def\b{0};
\draw[step=2, color=green](\a+2,-14+\b)grid(\a+14,-2+\b);
\draw[very thick](\a+7,-3.50000000000000+\b)--(\a+7,-4+\b);
\draw[very thick](\a+7.50000000000000,-3+\b)--(\a+8,-3+\b);
\draw[very thick](\a+7.50000000000000,-3+\b)arc(90:180:0.500000000000000);
\draw[very thick](\a+8,-3+\b)--(\a+10,-3+\b);
\draw[very thick](\a+10,-3+\b)--(\a+12,-3+\b);
\draw[very thick](\a+12,-3+\b)--(\a+14,-3+\b);
\draw[very thick](\a+7,-4+\b)--(\a+7,-6+\b);
\draw[very thick](\a+9,-5.50000000000000+\b)--(\a+9,-6+\b);
\draw[very thick](\a+9.50000000000000,-5+\b)--(\a+10,-5+\b);
\draw[very thick](\a+9.50000000000000,-5+\b)arc(90:180:0.500000000000000);
\draw[very thick](\a+10,-5+\b)--(\a+12,-5+\b);
\draw[very thick](\a+12,-5+\b)--(\a+14,-5+\b);
\draw[very thick](\a+3,-7.50000000000000+\b)--(\a+3,-8+\b);
\draw[very thick](\a+3.50000000000000,-7+\b)--(\a+4,-7+\b);
\draw[very thick](\a+3.50000000000000,-7+\b)arc(90:180:0.500000000000000);
\draw[very thick](\a+4,-7+\b)--(\a+6,-7+\b);
\draw[very thick](\a+7,-6+\b)--(\a+7,-6.50000000000000+\b);
\draw[very thick](\a+6,-7+\b)--(\a+6.50000000000000,-7+\b);
\draw[very thick](\a+6.50000000000000,-7+\b)arc(270:360:0.500000000000000);
\draw[very thick](\a+9,-6+\b)--(\a+9,-8+\b);
\draw[very thick](\a+11,-7.50000000000000+\b)--(\a+11,-8+\b);
\draw[very thick](\a+11.5000000000000,-7+\b)--(\a+12,-7+\b);
\draw[very thick](\a+11.5000000000000,-7+\b)arc(90:180:0.500000000000000);
\draw[very thick](\a+12,-7+\b)--(\a+14,-7+\b);
\draw[very thick](\a+3,-8+\b)--(\a+3,-10+\b);
\draw[very thick](\a+5,-9.50000000000000+\b)--(\a+5,-10+\b);
\draw[very thick](\a+5.50000000000000,-9+\b)--(\a+6,-9+\b);
\draw[very thick](\a+5.50000000000000,-9+\b)arc(90:180:0.500000000000000);
\draw[very thick](\a+6,-9+\b)--(\a+8,-9+\b);
\draw[very thick](\a+9,-8+\b)--(\a+9,-8.50000000000000+\b);
\draw[very thick](\a+8,-9+\b)--(\a+8.50000000000000,-9+\b);
\draw[very thick](\a+8.50000000000000,-9+\b)arc(270:360:0.500000000000000);
\draw[very thick](\a+11,-8+\b)--(\a+11,-10+\b);
\draw[very thick](\a+13,-9.50000000000000+\b)--(\a+13,-10+\b);
\draw[very thick](\a+13.5000000000000,-9+\b)--(\a+14,-9+\b);
\draw[very thick](\a+13.5000000000000,-9+\b)arc(90:180:0.500000000000000);
\draw[very thick](\a+3,-10+\b)--(\a+3,-12+\b);
\draw[very thick](\a+5,-10+\b)--(\a+5,-12+\b);
\draw[very thick](\a+7,-11.5000000000000+\b)--(\a+7,-12+\b);
\draw[very thick](\a+7.50000000000000,-11+\b)--(\a+8,-11+\b);
\draw[very thick](\a+7.50000000000000,-11+\b)arc(90:180:0.500000000000000);
\draw[very thick](\a+8,-11+\b)--(\a+10,-11+\b);
\draw[very thick](\a+10,-11+\b)--(\a+12,-11+\b);
\draw[very thick](\a+11,-10+\b)--(\a+11,-12+\b);
\draw[very thick](\a+12,-11+\b)--(\a+14,-11+\b);
\draw[very thick](\a+13,-10+\b)--(\a+13,-12+\b);
\draw[very thick](\a+3,-12+\b)--(\a+3,-14+\b);
\draw[very thick](\a+5,-12+\b)--(\a+5,-14+\b);
\draw[very thick](\a+7,-12+\b)--(\a+7,-14+\b);
\draw[very thick](\a+9,-13.5000000000000+\b)--(\a+9,-14+\b);
\draw[very thick](\a+9.50000000000000,-13+\b)--(\a+10,-13+\b);
\draw[very thick](\a+9.50000000000000,-13+\b)arc(90:180:0.500000000000000);
\draw[very thick](\a+10,-13+\b)--(\a+12,-13+\b);
\draw[very thick](\a+11,-12+\b)--(\a+11,-14+\b);
\draw[very thick](\a+12,-13+\b)--(\a+14,-13+\b);
\draw[very thick](\a+13,-12+\b)--(\a+13,-14+\b);
\end{tikzpicture}
\qquad
\begin{tikzpicture}[scale=0.2]
\draw[step=2](0,0)grid(12,12);
\node at (7,1) {$+$};
\node at (5,3) {$+$};
\node at (3,5) {$+$};
\node at (7,5) {$-$};
\node at (11,5) {$+$};
\node at (1,7) {$+$};
\node at (5,7) {$-$};
\node at (9,7) {$+$};
\node at (7,9) {$+$};
\node at (5,11) {$+$};
\end{tikzpicture}
\caption{The bijection between BPDs and ASMs}
\label{fig:BPD-ASM}
\end{figure}

Given a (not necessarily reduced) bumpless pipe dream, we can obtain
an ASM by translating a \rt-tile to $1$, a \jt-tile to $-1$ and all
other tiles to $0$. Knowing the set of up and down elbow tiles in a
bumpless pipe dream determines the entire set of tiles since the pipes
progress up and to the right. Therefore, it
is not hard to see that this map establishes a bijection. An example
is seen in \Cref{fig:BPD-ASM}.
It is not very easy to tell from the data of an ASM whether its corresponding BPD is reduced. Not necessarily reduced BPDs (and also not necessarily reduced PDs) are utilized to compute \emph{Grothendieck polynomials} \cite{weigandt-BPD-ASM}, the $K$-theoretic analogs of Schubert polynomials.  It is worth noting that Lascoux \cite{lascoux2008chern} was the first to write down Grothendieck polynomials as weighted sums over ASMs. 

Pipe dreams and bumpless pipe dreams have a lot of similarities
pictorially, as both can be viewed as certain versions of wiring
diagrams for a fixed permutation.  We highlight some notable
comparisons between these two combinatorial objects.  For the rest of
this section, a (bumpless) pipe dream is assumed to be reduced, but
some of the material will carry over to the non-reduced case as well.

\begin{itemize}
\item Both PDs and BPDs can be used to immediately show the stability of Schubert polynomials by extending the diagrams in the SE direction. 
\item One can extend a BPD in the NW direction and
consider permutations in $S_{\mathbb{Z}}\supset S_{\infty}$ to obtain
\emph{back stable Schubert polynomials} \cite{LLS}.  
 A similar backward extension for pipe dreams gives another formula
 for double Schubert polynomials \cite[Thm 4.1]{billey-bergeron}.
\item A PD is determined by the positions of its \+-tiles. However, a BPD cannot be determined by the data of its \+-tiles or \bl-tiles alone, but can be determined by both data together (\Cref{exer:BPD-blank-cross}).
\item Both PDs and BPDs of a fixed permutation are connected by local moves, called ``chute" / ``ladder" moves (discussed in \Cref{sub:pipes}) and ``droop" moves (\Cref{def:droop}) respectively. Moreover, one can start with the \emph{bottom} (bumpless) pipe dream and apply these local moves to obtain all other (bumpless) pipe dreams (\Cref{thm:chutes.and.ladders}, \Cref{lem:BPD-obtained-from-Rothe}). 
\end{itemize}

Let us discuss these local moves in more details. The following definition generalizes both chute moves and ladder moves that we have seen in \Cref{sub:pipes}.
\begin{Definition}\cite{billey-bergeron}\label{def:generalized-chute}
For a pipe dream $D$, a \emph{generalized chute move} from $(i,j)$ to $(i-a,j+b)$ with $a,b\in\mathbb{Z}_{+}$ is a local change that removes $(i,j)$ from $D$ and adds $(i-a,j+b)$ to $D$ provided that $(i',j')\in D$ for all $i-a\leq i'\leq i$ and $j\leq j'\leq j+b$ except $(i-a,j),(i,j+b),(i-a,j+b)\notin D$. 
\end{Definition}
See \Cref{fig:generalized-chute} for a much clearer explanation of
\Cref{def:generalized-chute} with diagrams. One can easily check that
a generalized chute move preserves the corresponding permutation by
inspection after drawing in the corresponding pipes.
\begin{figure}[h!]
\centering
\begin{tikzpicture}[scale=0.4]
\draw[fill=gray!30] (-.5,-2.5)--(.5,-2.5)--(.5,-1.5)--(-.5,-1.5)--(-.5,-2.5);
\node at (0,0) {$\cdot$};
\node at (1,0) {$+$};
\node at (2,0) {$+$};
\node at (3,0) {$\cdot$};
\node at (0,-1) {$+$};
\node at (1,-1) {$+$};
\node at (2,-1) {$+$};
\node at (3,-1) {$+$};
\node at (0,-2) {$+$};
\node at (1,-2) {$+$};
\node at (2,-2) {$+$};
\node at (3,-2) {$\cdot$};
\node at (-1.5,0) {$i{-}a$};
\node at (-1,-2) {$i$};
\node at (0,1) {$j$};
\node at (3,1) {$j{+}b$};

\def\a{10};
\draw[fill=gray!30] (\a+3.5,-.5)--(\a + 2.5,-.5)--(\a + 2.5,.5)--(\a+3.5,.5)--(\a+3.5,-.5);
\node at (\a,0) {$\cdot$};
\node at (\a+1,0) {$+$};
\node at (\a+2,0) {$+$};
\node at (\a+3,0) {$+$};
\node at (\a+0,-1) {$+$};
\node at (\a+1,-1) {$+$};
\node at (\a+2,-1) {$+$};
\node at (\a+3,-1) {$+$};
\node at (\a+0,-2) {$\cdot$};
\node at (\a+1,-2) {$+$};
\node at (\a+2,-2) {$+$};
\node at (\a+3,-2) {$\cdot$};
\node at (\a-1.5,0) {$i{-}a$};
\node at (\a-1,-2) {$i$};
\node at (\a+0,1) {$j$};
\node at (\a+3,1) {$j{+}b$};

\draw[very thick,->] (5,-1)--(6.5,-1);
\end{tikzpicture}
\caption{Generalized chute moves}
\label{fig:generalized-chute}
\end{figure}

Recall that the \emph{bottom pipe dream} of a permutation $w$ is defined to be the pipe dream that contains $c(w)_i:=|\{j>i\:|\:w(j)<w(i)\}|$ left-adjusted \+-tiles in row $i$, for each $i$, and that any pipe dream of $w$ can be obtained from the bottom pipe dream of $w$ using (generalized) chute moves (\Cref{thm:chutes.and.ladders}). A parallel story exists for bumpless pipe dreams.

\begin{Definition}\cite{LLS}\label{def:droop}
For a (reduced) bumpless pipe dream $D$, a \emph{droop move} from a
\rt-tile at $(i,j)$ to a \bl-tile at $(i+a,j+b)$ with
$a,b\in\mathbb{Z}_+$ is a local modification of a pipe $p$ that travels
$(i+a,j)\rightarrow (i,j)\rightarrow (i,j+b)$ to $(i+a,j)\rightarrow
(i+a,j+b)\rightarrow (i,j+b)$ such that the end result is still a
(reduced) bumpless pipe dream.  See \Cref{fig:droop-example} for
visualization.  A droop move is called a \emph{min-droop} if all tiles
in the rectangle $[i,i+a]\times [j,j+b]$ are \+-tiles except the four
corners.
\end{Definition}

\begin{figure}[h!]
\centering
\begin{tikzpicture}[scale=0.200000000000000]
\def\a{0};
\def\b{0};
\def\a{0};
\def\b{0};
\draw[step=2, color=green](\a+2,-8+\b)grid(\a+10,-2+\b);
\draw[very thick](\a+3,-3.50000000000000+\b)--(\a+3,-4+\b);
\draw[very thick](\a+3.50000000000000,-3+\b)--(\a+4,-3+\b);
\draw[very thick](\a+3.50000000000000,-3+\b)arc(90:180:0.500000000000000);
\draw[very thick](\a+4,-3+\b)--(\a+6,-3+\b);
\draw[very thick](\a+5,-2+\b)--(\a+5,-4+\b);
\draw[very thick](\a+6,-3+\b)--(\a+8,-3+\b);
\draw[very thick](\a+9,-2+\b)--(\a+9,-2.50000000000000+\b);
\draw[very thick](\a+8,-3+\b)--(\a+8.50000000000000,-3+\b);
\draw[very thick](\a+8.50000000000000,-3+\b)arc(270:360:0.500000000000000);
\draw[very thick](\a+2,-5+\b)--(\a+4,-5+\b);
\draw[very thick](\a+3,-4+\b)--(\a+3,-6+\b);
\draw[very thick](\a+4,-5+\b)--(\a+6,-5+\b);
\draw[very thick](\a+5,-4+\b)--(\a+5,-6+\b);
\draw[very thick](\a+6,-5+\b)--(\a+8,-5+\b);
\draw[very thick](\a+8,-5+\b)--(\a+10,-5+\b);
\draw[very thick](\a+3,-6+\b)--(\a+3,-8+\b);
\draw[very thick](\a+5,-6+\b)--(\a+5,-8+\b);

\def\a{16};
\def\b{0};
\draw[step=2, color=green](\a+2,-8+\b)grid(\a+10,-2+\b);
\draw[very thick](\a+5,-2+\b)--(\a+5,-4+\b);
\draw[very thick](\a+9,-2+\b)--(\a+9,-4+\b);
\draw[very thick](\a+2,-5+\b)--(\a+4,-5+\b);
\draw[very thick](\a+4,-5+\b)--(\a+6,-5+\b);
\draw[very thick](\a+5,-4+\b)--(\a+5,-6+\b);
\draw[very thick](\a+6,-5+\b)--(\a+8,-5+\b);
\draw[very thick](\a+8,-5+\b)--(\a+10,-5+\b);
\draw[very thick](\a+9,-4+\b)--(\a+9,-6+\b);
\draw[very thick](\a+3,-7.50000000000000+\b)--(\a+3,-8+\b);
\draw[very thick](\a+3.50000000000000,-7+\b)--(\a+4,-7+\b);
\draw[very thick](\a+3.50000000000000,-7+\b)arc(90:180:0.500000000000000);
\draw[very thick](\a+4,-7+\b)--(\a+6,-7+\b);
\draw[very thick](\a+5,-6+\b)--(\a+5,-8+\b);
\draw[very thick](\a+6,-7+\b)--(\a+8,-7+\b);
\draw[very thick](\a+9,-6+\b)--(\a+9,-6.50000000000000+\b);
\draw[very thick](\a+8,-7+\b)--(\a+8.50000000000000,-7+\b);
\draw[very thick](\a+8.50000000000000,-7+\b)arc(270:360:0.500000000000000);

\draw[very thick,->] (12,-5)--(16,-5); 
\end{tikzpicture}
\caption{A droop move (but not a min-droop) on bumpless pipe dreams}
\label{fig:droop-example}
\end{figure}

\begin{figure}[h!]
\centering
\begin{tikzpicture}[scale=0.1500000000000000]
\def\a{0};
\def\b{0};
\draw[color=green](\a+2,-2+\b)--(\a+10,-2+\b+\b);
\draw[color=green](\a+2,-2+\b)--(\a+2,-10+\b);
\draw[color=green](\a+2,-4+\b)--(\a+10,-4+\b);
\draw[color=green](\a+4,-2+\b)--(\a+4,-10+\b);
\draw[color=green](\a+2,-6+\b)--(\a+8,-6+\b);
\draw[color=green](\a+6,-2+\b)--(\a+6,-8+\b);
\draw[color=green](\a+2,-8+\b)--(\a+6,-8+\b);
\draw[color=green](\a+8,-2+\b)--(\a+8,-6+\b);
\draw[color=green](\a+2,-10+\b)--(\a+4,-10+\b);
\draw[color=green](\a+10,-2+\b)--(\a+10,-4+\b);
\draw[thick](\a+3,-2+\b)--(\a+3,-2.50000000000000+\b);
\draw[thick](\a+3,-3.50000000000000+\b)--(\a+3,-4+\b);
\draw[thick](\a+2,-3+\b)--(\a+2.50000000000000,-3+\b);
\draw[thick](\a+3.50000000000000,-3+\b)--(\a+4,-3+\b);
\draw[thick](\a+2.50000000000000,-3+\b)arc(270:360:0.500000000000000);
\draw[thick](\a+3.50000000000000,-3+\b)arc(90:180:0.500000000000000);
\draw[thick](\a+5,-2+\b)--(\a+5,-2.50000000000000+\b);
\draw[thick](\a+5,-3.50000000000000+\b)--(\a+5,-4+\b);
\draw[thick](\a+4,-3+\b)--(\a+4.50000000000000,-3+\b);
\draw[thick](\a+5.50000000000000,-3+\b)--(\a+6,-3+\b);
\draw[thick](\a+4.50000000000000,-3+\b)arc(270:360:0.500000000000000);
\draw[thick](\a+5.50000000000000,-3+\b)arc(90:180:0.500000000000000);
\draw[thick](\a+7,-2+\b)--(\a+7,-2.50000000000000+\b);
\draw[thick](\a+7,-3.50000000000000+\b)--(\a+7,-4+\b);
\draw[thick](\a+6,-3+\b)--(\a+6.50000000000000,-3+\b);
\draw[thick](\a+7.50000000000000,-3+\b)--(\a+8,-3+\b);
\draw[thick](\a+6.50000000000000,-3+\b)arc(270:360:0.500000000000000);
\draw[thick](\a+7.50000000000000,-3+\b)arc(90:180:0.500000000000000);
\draw[thick](\a+2,-5+\b)--(\a+4,-5+\b);
\draw[thick](\a+3,-4+\b)--(\a+3,-6+\b);
\draw[thick](\a+4,-5+\b)--(\a+6,-5+\b);
\draw[thick](\a+5,-4+\b)--(\a+5,-6+\b);
\draw[thick](\a+2,-7+\b)--(\a+4,-7+\b);
\draw[thick](\a+3,-6+\b)--(\a+3,-8+\b);
\draw[thick](\a+8,-3+\b)arc(270:360:1);
\draw[thick](\a+6,-5+\b)arc(270:360:1);
\draw[thick](\a+4,-7+\b)arc(270:360:1);
\draw[thick](\a+2,-9+\b)arc(270:360:1);

\def\a{-10};
\def\b{10};
\draw[color=green](\a+2,-2+\b)--(\a+10,-2+\b);
\draw[color=green](\a+2,-2+\b)--(\a+2,-10+\b);
\draw[color=green](\a+2,-4+\b)--(\a+10,-4+\b);
\draw[color=green](\a+4,-2+\b)--(\a+4,-10+\b);
\draw[color=green](\a+2,-6+\b)--(\a+8,-6+\b);
\draw[color=green](\a+6,-2+\b)--(\a+6,-8+\b);
\draw[color=green](\a+2,-8+\b)--(\a+6,-8+\b);
\draw[color=green](\a+8,-2+\b)--(\a+8,-6+\b);
\draw[color=green](\a+2,-10+\b)--(\a+4,-10+\b);
\draw[color=green](\a+10,-2+\b)--(\a+10,-4+\b);
\draw[thick](\a+3,-2+\b)--(\a+3,-2.50000000000000+\b);
\draw[thick](\a+3,-3.50000000000000+\b)--(\a+3,-4+\b);
\draw[thick](\a+2,-3+\b)--(\a+2.50000000000000,-3+\b);
\draw[thick](\a+3.50000000000000,-3+\b)--(\a+4,-3+\b);
\draw[thick](\a+2.50000000000000,-3+\b)arc(270:360:0.500000000000000);
\draw[thick](\a+3.50000000000000,-3+\b)arc(90:180:0.500000000000000);
\draw[thick](\a+5,-2+\b)--(\a+5,-2.50000000000000+\b);
\draw[thick](\a+5,-3.50000000000000+\b)--(\a+5,-4+\b);
\draw[thick](\a+4,-3+\b)--(\a+4.50000000000000,-3+\b);
\draw[thick](\a+5.50000000000000,-3+\b)--(\a+6,-3+\b);
\draw[thick](\a+4.50000000000000,-3+\b)arc(270:360:0.500000000000000);
\draw[thick](\a+5.50000000000000,-3+\b)arc(90:180:0.500000000000000);
\draw[thick](\a+6,-3+\b)--(\a+8,-3+\b);
\draw[thick](\a+7,-2+\b)--(\a+7,-4+\b);
\draw[thick](\a+2,-5+\b)--(\a+4,-5+\b);
\draw[thick](\a+3,-4+\b)--(\a+3,-6+\b);
\draw[thick](\a+5,-4+\b)--(\a+5,-4.50000000000000+\b);
\draw[thick](\a+5,-5.50000000000000+\b)--(\a+5,-6+\b);
\draw[thick](\a+4,-5+\b)--(\a+4.50000000000000,-5+\b);
\draw[thick](\a+5.50000000000000,-5+\b)--(\a+6,-5+\b);
\draw[thick](\a+4.50000000000000,-5+\b)arc(270:360:0.500000000000000);
\draw[thick](\a+5.50000000000000,-5+\b)arc(90:180:0.500000000000000);
\draw[thick](\a+2,-7+\b)--(\a+4,-7+\b);
\draw[thick](\a+3,-6+\b)--(\a+3,-8+\b);
\draw[thick](\a+8,-3+\b)arc(270:360:1);
\draw[thick](\a+6,-5+\b)arc(270:360:1);
\draw[thick](\a+4,-7+\b)arc(270:360:1);
\draw[thick](\a+2,-9+\b)arc(270:360:1);

\def\a{-10};
\def\b{22};
\draw[color=green](\a+2,-2+\b)--(\a+10,-2+\b);
\draw[color=green](\a+2,-2+\b)--(\a+2,-10+\b);
\draw[color=green](\a+2,-4+\b)--(\a+10,-4+\b);
\draw[color=green](\a+4,-2+\b)--(\a+4,-10+\b);
\draw[color=green](\a+2,-6+\b)--(\a+8,-6+\b);
\draw[color=green](\a+6,-2+\b)--(\a+6,-8+\b);
\draw[color=green](\a+2,-8+\b)--(\a+6,-8+\b);
\draw[color=green](\a+8,-2+\b)--(\a+8,-6+\b);
\draw[color=green](\a+2,-10+\b)--(\a+4,-10+\b);
\draw[color=green](\a+10,-2+\b)--(\a+10,-4+\b);
\draw[thick](\a+3,-2+\b)--(\a+3,-2.50000000000000+\b);
\draw[thick](\a+3,-3.50000000000000+\b)--(\a+3,-4+\b);
\draw[thick](\a+2,-3+\b)--(\a+2.50000000000000,-3+\b);
\draw[thick](\a+3.50000000000000,-3+\b)--(\a+4,-3+\b);
\draw[thick](\a+2.50000000000000,-3+\b)arc(270:360:0.500000000000000);
\draw[thick](\a+3.50000000000000,-3+\b)arc(90:180:0.500000000000000);
\draw[thick](\a+4,-3+\b)--(\a+6,-3+\b);
\draw[thick](\a+5,-2+\b)--(\a+5,-4+\b);
\draw[thick](\a+6,-3+\b)--(\a+8,-3+\b);
\draw[thick](\a+7,-2+\b)--(\a+7,-4+\b);
\draw[thick](\a+3,-4+\b)--(\a+3,-4.50000000000000+\b);
\draw[thick](\a+3,-5.50000000000000+\b)--(\a+3,-6+\b);
\draw[thick](\a+2,-5+\b)--(\a+2.50000000000000,-5+\b);
\draw[thick](\a+3.50000000000000,-5+\b)--(\a+4,-5+\b);
\draw[thick](\a+2.50000000000000,-5+\b)arc(270:360:0.500000000000000);
\draw[thick](\a+3.50000000000000,-5+\b)arc(90:180:0.500000000000000);
\draw[thick](\a+5,-4+\b)--(\a+5,-4.50000000000000+\b);
\draw[thick](\a+5,-5.50000000000000+\b)--(\a+5,-6+\b);
\draw[thick](\a+4,-5+\b)--(\a+4.50000000000000,-5+\b);
\draw[thick](\a+5.50000000000000,-5+\b)--(\a+6,-5+\b);
\draw[thick](\a+4.50000000000000,-5+\b)arc(270:360:0.500000000000000);
\draw[thick](\a+5.50000000000000,-5+\b)arc(90:180:0.500000000000000);
\draw[thick](\a+2,-7+\b)--(\a+4,-7+\b);
\draw[thick](\a+3,-6+\b)--(\a+3,-8+\b);
\draw[thick](\a+8,-3+\b)arc(270:360:1);
\draw[thick](\a+6,-5+\b)arc(270:360:1);
\draw[thick](\a+4,-7+\b)arc(270:360:1);
\draw[thick](\a+2,-9+\b)arc(270:360:1);

\def\a{0};
\def\b{32};
\draw[color=green](\a+2,-2+\b)--(\a+10,-2+\b);
\draw[color=green](\a+2,-2+\b)--(\a+2,-10+\b);
\draw[color=green](\a+2,-4+\b)--(\a+10,-4+\b);
\draw[color=green](\a+4,-2+\b)--(\a+4,-10+\b);
\draw[color=green](\a+2,-6+\b)--(\a+8,-6+\b);
\draw[color=green](\a+6,-2+\b)--(\a+6,-8+\b);
\draw[color=green](\a+2,-8+\b)--(\a+6,-8+\b);
\draw[color=green](\a+8,-2+\b)--(\a+8,-6+\b);
\draw[color=green](\a+2,-10+\b)--(\a+4,-10+\b);
\draw[color=green](\a+10,-2+\b)--(\a+10,-4+\b);
\draw[thick](\a+3,-2+\b)--(\a+3,-2.50000000000000+\b);
\draw[thick](\a+3,-3.50000000000000+\b)--(\a+3,-4+\b);
\draw[thick](\a+2,-3+\b)--(\a+2.50000000000000,-3+\b);
\draw[thick](\a+3.50000000000000,-3+\b)--(\a+4,-3+\b);
\draw[thick](\a+2.50000000000000,-3+\b)arc(270:360:0.500000000000000);
\draw[thick](\a+3.50000000000000,-3+\b)arc(90:180:0.500000000000000);
\draw[thick](\a+4,-3+\b)--(\a+6,-3+\b);
\draw[thick](\a+5,-2+\b)--(\a+5,-4+\b);
\draw[thick](\a+6,-3+\b)--(\a+8,-3+\b);
\draw[thick](\a+7,-2+\b)--(\a+7,-4+\b);
\draw[thick](\a+3,-4+\b)--(\a+3,-4.50000000000000+\b);
\draw[thick](\a+3,-5.50000000000000+\b)--(\a+3,-6+\b);
\draw[thick](\a+2,-5+\b)--(\a+2.50000000000000,-5+\b);
\draw[thick](\a+3.50000000000000,-5+\b)--(\a+4,-5+\b);
\draw[thick](\a+2.50000000000000,-5+\b)arc(270:360:0.500000000000000);
\draw[thick](\a+3.50000000000000,-5+\b)arc(90:180:0.500000000000000);
\draw[thick](\a+4,-5+\b)--(\a+6,-5+\b);
\draw[thick](\a+5,-4+\b)--(\a+5,-6+\b);
\draw[thick](\a+3,-6+\b)--(\a+3,-6.50000000000000+\b);
\draw[thick](\a+3,-7.50000000000000+\b)--(\a+3,-8+\b);
\draw[thick](\a+2,-7+\b)--(\a+2.50000000000000,-7+\b);
\draw[thick](\a+3.50000000000000,-7+\b)--(\a+4,-7+\b);
\draw[thick](\a+2.50000000000000,-7+\b)arc(270:360:0.500000000000000);
\draw[thick](\a+3.50000000000000,-7+\b)arc(90:180:0.500000000000000);
\draw[thick](\a+8,-3+\b)arc(270:360:1);
\draw[thick](\a+6,-5+\b)arc(270:360:1);
\draw[thick](\a+4,-7+\b)arc(270:360:1);
\draw[thick](\a+2,-9+\b)arc(270:360:1);

\def\a{10};
\def\b{16};
\draw[color=green](\a+2,-2+\b)--(\a+10,-2+\b);
\draw[color=green](\a+2,-2+\b)--(\a+2,-10+\b);
\draw[color=green](\a+2,-4+\b)--(\a+10,-4+\b);
\draw[color=green](\a+4,-2+\b)--(\a+4,-10+\b);
\draw[color=green](\a+2,-6+\b)--(\a+8,-6+\b);
\draw[color=green](\a+6,-2+\b)--(\a+6,-8+\b);
\draw[color=green](\a+2,-8+\b)--(\a+6,-8+\b);
\draw[color=green](\a+8,-2+\b)--(\a+8,-6+\b);
\draw[color=green](\a+2,-10+\b)--(\a+4,-10+\b);
\draw[color=green](\a+10,-2+\b)--(\a+10,-4+\b);
\draw[thick](\a+3,-2+\b)--(\a+3,-2.50000000000000+\b);
\draw[thick](\a+3,-3.50000000000000+\b)--(\a+3,-4+\b);
\draw[thick](\a+2,-3+\b)--(\a+2.50000000000000,-3+\b);
\draw[thick](\a+3.50000000000000,-3+\b)--(\a+4,-3+\b);
\draw[thick](\a+2.50000000000000,-3+\b)arc(270:360:0.500000000000000);
\draw[thick](\a+3.50000000000000,-3+\b)arc(90:180:0.500000000000000);
\draw[thick](\a+4,-3+\b)--(\a+6,-3+\b);
\draw[thick](\a+5,-2+\b)--(\a+5,-4+\b);
\draw[thick](\a+7,-2+\b)--(\a+7,-2.50000000000000+\b);
\draw[thick](\a+7,-3.50000000000000+\b)--(\a+7,-4+\b);
\draw[thick](\a+6,-3+\b)--(\a+6.50000000000000,-3+\b);
\draw[thick](\a+7.50000000000000,-3+\b)--(\a+8,-3+\b);
\draw[thick](\a+6.50000000000000,-3+\b)arc(270:360:0.500000000000000);
\draw[thick](\a+7.50000000000000,-3+\b)arc(90:180:0.500000000000000);
\draw[thick](\a+2,-5+\b)--(\a+4,-5+\b);
\draw[thick](\a+3,-4+\b)--(\a+3,-6+\b);
\draw[thick](\a+4,-5+\b)--(\a+6,-5+\b);
\draw[thick](\a+5,-4+\b)--(\a+5,-6+\b);
\draw[thick](\a+3,-6+\b)--(\a+3,-6.50000000000000+\b);
\draw[thick](\a+3,-7.50000000000000+\b)--(\a+3,-8+\b);
\draw[thick](\a+2,-7+\b)--(\a+2.50000000000000,-7+\b);
\draw[thick](\a+3.50000000000000,-7+\b)--(\a+4,-7+\b);
\draw[thick](\a+2.50000000000000,-7+\b)arc(270:360:0.500000000000000);
\draw[thick](\a+3.50000000000000,-7+\b)arc(90:180:0.500000000000000);
\draw[thick](\a+8,-3+\b)arc(270:360:1);
\draw[thick](\a+6,-5+\b)arc(270:360:1);
\draw[thick](\a+4,-7+\b)arc(270:360:1);
\draw[thick](\a+2,-9+\b)arc(270:360:1);

\draw[very thick] (1,-1)--(-3,3);
\draw[very thick] (-3,9)--(-3,15);
\draw[very thick] (-3,21)--(1,25);
\draw[very thick] (11,-1)--(17,8);
\draw[very thick] (17,15)--(7,25);
\end{tikzpicture}
\qquad\qquad
\begin{tikzpicture}[scale=0.1500000000000000]
\def\a{0};
\def\b{0};
\draw[step=2, color=green](\a+2,-10+\b)grid(\a+10,-2+\b);
\draw[thick](\a+3,-3.50000000000000+\b)--(\a+3,-4+\b);
\draw[thick](\a+3.50000000000000,-3+\b)--(\a+4,-3+\b);
\draw[thick](\a+3.50000000000000,-3+\b)arc(90:180:0.500000000000000);
\draw[thick](\a+4,-3+\b)--(\a+6,-3+\b);
\draw[thick](\a+6,-3+\b)--(\a+8,-3+\b);
\draw[thick](\a+8,-3+\b)--(\a+10,-3+\b);
\draw[thick](\a+3,-4+\b)--(\a+3,-6+\b);
\draw[thick](\a+9,-5.50000000000000+\b)--(\a+9,-6+\b);
\draw[thick](\a+9.50000000000000,-5+\b)--(\a+10,-5+\b);
\draw[thick](\a+9.50000000000000,-5+\b)arc(90:180:0.500000000000000);
\draw[thick](\a+3,-6+\b)--(\a+3,-8+\b);
\draw[thick](\a+7,-7.50000000000000+\b)--(\a+7,-8+\b);
\draw[thick](\a+7.50000000000000,-7+\b)--(\a+8,-7+\b);
\draw[thick](\a+7.50000000000000,-7+\b)arc(90:180:0.500000000000000);
\draw[thick](\a+8,-7+\b)--(\a+10,-7+\b);
\draw[thick](\a+9,-6+\b)--(\a+9,-8+\b);
\draw[thick](\a+3,-8+\b)--(\a+3,-10+\b);
\draw[thick](\a+5,-9.50000000000000+\b)--(\a+5,-10+\b);
\draw[thick](\a+5.50000000000000,-9+\b)--(\a+6,-9+\b);
\draw[thick](\a+5.50000000000000,-9+\b)arc(90:180:0.500000000000000);
\draw[thick](\a+6,-9+\b)--(\a+8,-9+\b);
\draw[thick](\a+7,-8+\b)--(\a+7,-10+\b);
\draw[thick](\a+8,-9+\b)--(\a+10,-9+\b);
\draw[thick](\a+9,-8+\b)--(\a+9,-10+\b);

\def\a{0};
\def\b{12};
\draw[step=2, color=green](\a+2,-10+\b)grid(\a+10,-2+\b);
\draw[thick](\a+5,-3.50000000000000+\b)--(\a+5,-4+\b);
\draw[thick](\a+5.50000000000000,-3+\b)--(\a+6,-3+\b);
\draw[thick](\a+5.50000000000000,-3+\b)arc(90:180:0.500000000000000);
\draw[thick](\a+6,-3+\b)--(\a+8,-3+\b);
\draw[thick](\a+8,-3+\b)--(\a+10,-3+\b);
\draw[thick](\a+3,-5.50000000000000+\b)--(\a+3,-6+\b);
\draw[thick](\a+3.50000000000000,-5+\b)--(\a+4,-5+\b);
\draw[thick](\a+3.50000000000000,-5+\b)arc(90:180:0.500000000000000);
\draw[thick](\a+5,-4+\b)--(\a+5,-4.50000000000000+\b);
\draw[thick](\a+4,-5+\b)--(\a+4.50000000000000,-5+\b);
\draw[thick](\a+4.50000000000000,-5+\b)arc(270:360:0.500000000000000);
\draw[thick](\a+9,-5.50000000000000+\b)--(\a+9,-6+\b);
\draw[thick](\a+9.50000000000000,-5+\b)--(\a+10,-5+\b);
\draw[thick](\a+9.50000000000000,-5+\b)arc(90:180:0.500000000000000);
\draw[thick](\a+3,-6+\b)--(\a+3,-8+\b);
\draw[thick](\a+7,-7.50000000000000+\b)--(\a+7,-8+\b);
\draw[thick](\a+7.50000000000000,-7+\b)--(\a+8,-7+\b);
\draw[thick](\a+7.50000000000000,-7+\b)arc(90:180:0.500000000000000);
\draw[thick](\a+8,-7+\b)--(\a+10,-7+\b);
\draw[thick](\a+9,-6+\b)--(\a+9,-8+\b);
\draw[thick](\a+3,-8+\b)--(\a+3,-10+\b);
\draw[thick](\a+5,-9.50000000000000+\b)--(\a+5,-10+\b);
\draw[thick](\a+5.50000000000000,-9+\b)--(\a+6,-9+\b);
\draw[thick](\a+5.50000000000000,-9+\b)arc(90:180:0.500000000000000);
\draw[thick](\a+6,-9+\b)--(\a+8,-9+\b);
\draw[thick](\a+7,-8+\b)--(\a+7,-10+\b);
\draw[thick](\a+8,-9+\b)--(\a+10,-9+\b);
\draw[thick](\a+9,-8+\b)--(\a+9,-10+\b);

\def\a{-10};
\def\b{22};
\draw[step=2, color=green](\a+2,-10+\b)grid(\a+10,-2+\b);
\draw[thick](\a+5,-3.50000000000000+\b)--(\a+5,-4+\b);
\draw[thick](\a+5.50000000000000,-3+\b)--(\a+6,-3+\b);
\draw[thick](\a+5.50000000000000,-3+\b)arc(90:180:0.500000000000000);
\draw[thick](\a+6,-3+\b)--(\a+8,-3+\b);
\draw[thick](\a+8,-3+\b)--(\a+10,-3+\b);
\draw[thick](\a+5,-4+\b)--(\a+5,-6+\b);
\draw[thick](\a+9,-5.50000000000000+\b)--(\a+9,-6+\b);
\draw[thick](\a+9.50000000000000,-5+\b)--(\a+10,-5+\b);
\draw[thick](\a+9.50000000000000,-5+\b)arc(90:180:0.500000000000000);
\draw[thick](\a+3,-7.50000000000000+\b)--(\a+3,-8+\b);
\draw[thick](\a+3.50000000000000,-7+\b)--(\a+4,-7+\b);
\draw[thick](\a+3.50000000000000,-7+\b)arc(90:180:0.500000000000000);
\draw[thick](\a+5,-6+\b)--(\a+5,-6.50000000000000+\b);
\draw[thick](\a+4,-7+\b)--(\a+4.50000000000000,-7+\b);
\draw[thick](\a+4.50000000000000,-7+\b)arc(270:360:0.500000000000000);
\draw[thick](\a+7,-7.50000000000000+\b)--(\a+7,-8+\b);
\draw[thick](\a+7.50000000000000,-7+\b)--(\a+8,-7+\b);
\draw[thick](\a+7.50000000000000,-7+\b)arc(90:180:0.500000000000000);
\draw[thick](\a+8,-7+\b)--(\a+10,-7+\b);
\draw[thick](\a+9,-6+\b)--(\a+9,-8+\b);
\draw[thick](\a+3,-8+\b)--(\a+3,-10+\b);
\draw[thick](\a+5,-9.50000000000000+\b)--(\a+5,-10+\b);
\draw[thick](\a+5.50000000000000,-9+\b)--(\a+6,-9+\b);
\draw[thick](\a+5.50000000000000,-9+\b)arc(90:180:0.500000000000000);
\draw[thick](\a+6,-9+\b)--(\a+8,-9+\b);
\draw[thick](\a+7,-8+\b)--(\a+7,-10+\b);
\draw[thick](\a+8,-9+\b)--(\a+10,-9+\b);
\draw[thick](\a+9,-8+\b)--(\a+9,-10+\b);

\def\a{10};
\def\b{22};
\draw[step=2, color=green](\a+2,-10+\b)grid(\a+10,-2+\b);
\draw[thick](\a+7,-3.50000000000000+\b)--(\a+7,-4+\b);
\draw[thick](\a+7.50000000000000,-3+\b)--(\a+8,-3+\b);
\draw[thick](\a+7.50000000000000,-3+\b)arc(90:180:0.500000000000000);
\draw[thick](\a+8,-3+\b)--(\a+10,-3+\b);
\draw[thick](\a+3,-5.50000000000000+\b)--(\a+3,-6+\b);
\draw[thick](\a+3.50000000000000,-5+\b)--(\a+4,-5+\b);
\draw[thick](\a+3.50000000000000,-5+\b)arc(90:180:0.500000000000000);
\draw[thick](\a+4,-5+\b)--(\a+6,-5+\b);
\draw[thick](\a+7,-4+\b)--(\a+7,-4.50000000000000+\b);
\draw[thick](\a+6,-5+\b)--(\a+6.50000000000000,-5+\b);
\draw[thick](\a+6.50000000000000,-5+\b)arc(270:360:0.500000000000000);
\draw[thick](\a+9,-5.50000000000000+\b)--(\a+9,-6+\b);
\draw[thick](\a+9.50000000000000,-5+\b)--(\a+10,-5+\b);
\draw[thick](\a+9.50000000000000,-5+\b)arc(90:180:0.500000000000000);
\draw[thick](\a+3,-6+\b)--(\a+3,-8+\b);
\draw[thick](\a+7,-7.50000000000000+\b)--(\a+7,-8+\b);
\draw[thick](\a+7.50000000000000,-7+\b)--(\a+8,-7+\b);
\draw[thick](\a+7.50000000000000,-7+\b)arc(90:180:0.500000000000000);
\draw[thick](\a+8,-7+\b)--(\a+10,-7+\b);
\draw[thick](\a+9,-6+\b)--(\a+9,-8+\b);
\draw[thick](\a+3,-8+\b)--(\a+3,-10+\b);
\draw[thick](\a+5,-9.50000000000000+\b)--(\a+5,-10+\b);
\draw[thick](\a+5.50000000000000,-9+\b)--(\a+6,-9+\b);
\draw[thick](\a+5.50000000000000,-9+\b)arc(90:180:0.500000000000000);
\draw[thick](\a+6,-9+\b)--(\a+8,-9+\b);
\draw[thick](\a+7,-8+\b)--(\a+7,-10+\b);
\draw[thick](\a+8,-9+\b)--(\a+10,-9+\b);
\draw[thick](\a+9,-8+\b)--(\a+9,-10+\b);

\def\a{0};
\def\b{32};
\draw[step=2, color=green](\a+2,-10+\b)grid(\a+10,-2+\b);
\draw[thick](\a+7,-3.50000000000000+\b)--(\a+7,-4+\b);
\draw[thick](\a+7.50000000000000,-3+\b)--(\a+8,-3+\b);
\draw[thick](\a+7.50000000000000,-3+\b)arc(90:180:0.500000000000000);
\draw[thick](\a+8,-3+\b)--(\a+10,-3+\b);
\draw[thick](\a+5,-5.50000000000000+\b)--(\a+5,-6+\b);
\draw[thick](\a+5.50000000000000,-5+\b)--(\a+6,-5+\b);
\draw[thick](\a+5.50000000000000,-5+\b)arc(90:180:0.500000000000000);
\draw[thick](\a+7,-4+\b)--(\a+7,-4.50000000000000+\b);
\draw[thick](\a+6,-5+\b)--(\a+6.50000000000000,-5+\b);
\draw[thick](\a+6.50000000000000,-5+\b)arc(270:360:0.500000000000000);
\draw[thick](\a+9,-5.50000000000000+\b)--(\a+9,-6+\b);
\draw[thick](\a+9.50000000000000,-5+\b)--(\a+10,-5+\b);
\draw[thick](\a+9.50000000000000,-5+\b)arc(90:180:0.500000000000000);
\draw[thick](\a+3,-7.50000000000000+\b)--(\a+3,-8+\b);
\draw[thick](\a+3.50000000000000,-7+\b)--(\a+4,-7+\b);
\draw[thick](\a+3.50000000000000,-7+\b)arc(90:180:0.500000000000000);
\draw[thick](\a+5,-6+\b)--(\a+5,-6.50000000000000+\b);
\draw[thick](\a+4,-7+\b)--(\a+4.50000000000000,-7+\b);
\draw[thick](\a+4.50000000000000,-7+\b)arc(270:360:0.500000000000000);
\draw[thick](\a+7,-7.50000000000000+\b)--(\a+7,-8+\b);
\draw[thick](\a+7.50000000000000,-7+\b)--(\a+8,-7+\b);
\draw[thick](\a+7.50000000000000,-7+\b)arc(90:180:0.500000000000000);
\draw[thick](\a+8,-7+\b)--(\a+10,-7+\b);
\draw[thick](\a+9,-6+\b)--(\a+9,-8+\b);
\draw[thick](\a+3,-8+\b)--(\a+3,-10+\b);
\draw[thick](\a+5,-9.50000000000000+\b)--(\a+5,-10+\b);
\draw[thick](\a+5.50000000000000,-9+\b)--(\a+6,-9+\b);
\draw[thick](\a+5.50000000000000,-9+\b)arc(90:180:0.500000000000000);
\draw[thick](\a+6,-9+\b)--(\a+8,-9+\b);
\draw[thick](\a+7,-8+\b)--(\a+7,-10+\b);
\draw[thick](\a+8,-9+\b)--(\a+10,-9+\b);
\draw[thick](\a+9,-8+\b)--(\a+9,-10+\b);

\draw[very thick] (6,-1.5)--(6,1.5);
\draw[very thick] (1,7)--(-3,11);
\draw[very thick] (11,7)--(15,11);
\draw[very thick] (-3,21)--(1,25);
\draw[very thick] (11,25)--(15,21);
\end{tikzpicture}
\caption{The posets of pipe dreams and of bumpless pipe dreams for $1432$}
\label{fig:chute-versus-droop}
\end{figure}
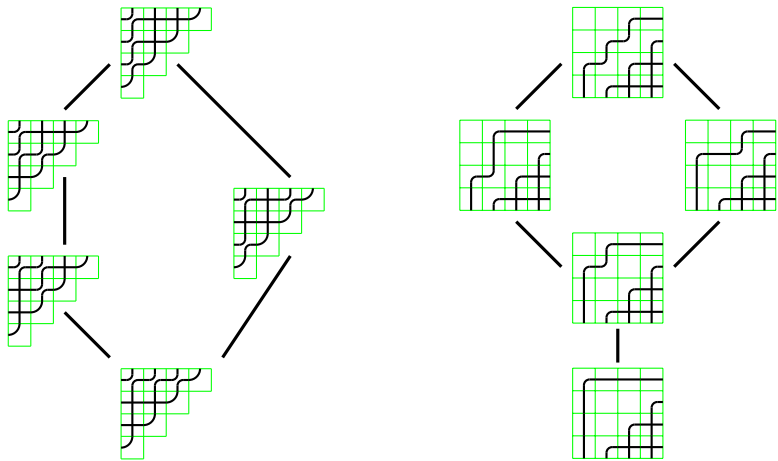

\begin{Definition}\cite{LLS}
The \emph{Rothe bumpless pipe dream} (or the \emph{bottom bumpless pipe dream}) for a permutation $w$ is the bumpless pipe dream where pipe $i$ only turns once at the \rt-tile at $(i,w(i))$, for all $i$.
\end{Definition}
\begin{Lemma}\cite{LLS}\label{lem:BPD-obtained-from-Rothe}
Any bumpless pipe dream of $w$ can be obtained from the Rothe bumpless pipe dream of $w$ using droop moves.
\end{Lemma}
Although generalized chute moves on pipe dreams and droop moves on
bumpless pipe dreams have a similar flavor, their respective poset
structures are different. To be precise, fix a permutation $w$. Then,
one can construct a poset on $\rp(w)$ under generalized chute moves
with the bottom pipe dream as the unique minimal element $\hat0$, and
one can construct a poset on $\BPD(w)$ under droop moves with the
Rothe bumpless pipe dream as $\hat0$. These two posets have the same
number of elements, but they hardly have anything else in common. Recall, there is a \emph{top pipe dream} $\Dtop(w)$ for all permutations $w$ (\Cref{sub:Games}), but there is only a top bumpless pipe dream in some cases (\Cref{ex:bpd-poset-vexillary}). See \Cref{fig:chute-versus-droop} for further comparison.

There is a fun conjecture by Rubey \cite{rubey2012maximal} on the
structure of generalized chute moves, stated below.  Certain special
cases of \Cref{conj:rubey-lattice} are known, including dominant
permutations, vexillary permutations, and $1432$-avoiding
permutations, but the conjecture is still open in general.  Recall
from \Cref{ex:dominant.1} that the
\emph{dominant permutations} are those which have exactly one reduced
pipe dream, or equivalently, the $132$-avoiding permutations. \emph{Vexillary permutations} are those avoiding $2143$ and will be discussed more extensively in \S\ref{sub:StanleySymmetrics}.

\begin{Conjecture}\cite{rubey2012maximal}\label{conj:rubey-lattice}
For  any permutation $w$, the poset on $\rp(w)$ under generalized
chute moves is a lattice.
\end{Conjecture}

\begin{Exercise}\label{ex:bpd-poset-vexillary}
Fix a permutation $w \in S_{n}$. Show that the poset on $\BPD(w)$
under droop moves has a unique maximal element if and only if $w$ is
vexillary. Furthermore, show that when $w$ is
vexillary, this poset is a lattice.
\end{Exercise}

Of course, one immediate question is to describe an explicit and
``natural" weight-preserving bijection between reduced pipe dreams and
reduced bumpless pipe dreams. Since both objects compute Schubert
polynomials with relevant theorems established independently, a
weight-preserving bijection is guaranteed to exist. It is the exact
description that is important. This question is answered by Gao and
Huang \cite{GaoHuang-bijection}, and we will mainly follow their
exposition here.  Recall we are still assuming all (bumpless and
classic) pipe dreams are reduced in what remains of this subsection.

Recall from \Cref{sub:Games} that pipe dreams $D$ are cryptomorphic to
their corresponding biwords $(\mathbf{r}_{D},\mathbf{j}_{D})$ and
$(\mathbf{r}_{D},\mathbf{i}_{D})$.  We will use the
$(\mathbf{r}_{D},\mathbf{i}_{D})$ encoding as compatible pairs with
the conditions described in \Cref{ex:bounded.pair}
to give a bijection between bumpless pipe dreams and pipe dreams.
Recall that we index each pipe of a bumpless pipe dream by the column that it starts off from the south boundary. 

\begin{Definition}\cite{GaoHuang-bijection}\label{def:pop}
Given a permutation $w \neq \mathrm{id}$ and  $D\in\BPD(w)$, the following procedure produces $\nabla
D\in\BPD(w')$ for some $w'=s_aw$ such that
$\ell(w')=\ell(w)-1$ and outputs a pair of positive integers
$\mathrm{pop}(D)=(r,i)$. 

\begin{enumerate}
\item Let $i$ be the smallest row index that contains a \bl-tile. Mark
the rightmost \bl-tile in row $i$ with a label $\times$, say in
position $(x,y)$. Here $x=i$ initially. Let $p$ be the unique pipe passing through
$(x,y+1)$, which must exist by the choice of blank tile.

\item If $p\leq y$, there is a \jt-tile at some coordinate $(x',y+1)$ with
$x'>x$ belonging to pipe $p$. Undroop this pipe from $(x',y+1)$ to
$(x,y)$, and place the label $\times$ at $(x',y+1)$. See
\Cref{fig:bpd-pop-situation1}.  

\item Next, move $\times$ to the rightmost
\bl-tile among its horizontally contiguous block of \bl-tiles, and
update  $(x,y)$ to be its coordinates. There is a unique pipe $p$ that passes through $(x,y+1)$ again. Repeat step (2),(3) with the updated $x,y,p$.

\item If $p=y+1$, pipe $y$ must intersect pipe
$p=y+1$ at coordinate $(x',y+1)$ with $x'>x$. Remove this crossing and
remove the label $\times$ by replacing the \+-tile at $(x',y+1)$ with
a \elbow-tile temporarily and undrooping its \jt-part to $(x,y)$. See
\Cref{fig:bpd-pop-situation2}. The whole process ends here. Finally, let $r=y$,
and return $(r,i)$.
\end{enumerate}
\end{Definition}

Note, this definition of $\nabla$ is different than the one defined in
\Cref{thm:derivative-schubert}.  Both uses of this notation indicate
the concept of going down in Bruhat order and decreasing the degree by
1.  At this point, the reader may guess what we are doing next to get the
desired bijection.

\begin{figure}[h!]
\centering
\begin{tikzpicture}[scale=0.200000000000000]
\def\a{0};
\def\b{0};
\draw[step=2, color=green](\a+2,-14+\b)grid(\a+6,-2+\b);
\draw[very thick](\a+5,-2+\b)--(\a+5,-4+\b);
\draw[very thick](\a+2,-5+\b)--(\a+4,-5+\b);
\draw[very thick](\a+4,-5+\b)--(\a+6,-5+\b);
\draw[very thick](\a+5,-4+\b)--(\a+5,-6+\b);
\draw[very thick](\a+5,-6+\b)--(\a+5,-8+\b);
\draw[very thick](\a+3,-9.50000000000000+\b)--(\a+3,-10+\b);
\draw[very thick](\a+3.50000000000000,-9+\b)--(\a+4,-9+\b);
\draw[very thick](\a+3.50000000000000,-9+\b)arc(90:180:0.500000000000000);
\draw[very thick](\a+4,-9+\b)--(\a+6,-9+\b);
\draw[very thick](\a+5,-8+\b)--(\a+5,-10+\b);
\draw[very thick](\a+3,-10+\b)--(\a+3,-10.5000000000000+\b);
\draw[very thick](\a+2,-11+\b)--(\a+2.50000000000000,-11+\b);
\draw[very thick](\a+2.50000000000000,-11+\b)arc(270:360:0.500000000000000);
\draw[very thick](\a+5,-10+\b)--(\a+5,-12+\b);
\draw[very thick](\a+2,-13+\b)--(\a+4,-13+\b);
\draw[very thick](\a+5,-12+\b)--(\a+5,-12.5000000000000+\b);
\draw[very thick](\a+4,-13+\b)--(\a+4.50000000000000,-13+\b);
\draw[very thick](\a+4.50000000000000,-13+\b)arc(270:360:0.500000000000000);

\def\a{8};
\def\b{0};
\draw[step=2, color=green](\a+2,-14+\b)grid(\a+6,-2+\b);
\draw[very thick](\a+3,-3.50000000000000+\b)--(\a+3,-4+\b);
\draw[very thick](\a+3.50000000000000,-3+\b)--(\a+4,-3+\b);
\draw[very thick](\a+3.50000000000000,-3+\b)arc(90:180:0.500000000000000);
\draw[very thick](\a+5,-2+\b)--(\a+5,-2.50000000000000+\b);
\draw[very thick](\a+4,-3+\b)--(\a+4.50000000000000,-3+\b);
\draw[very thick](\a+4.50000000000000,-3+\b)arc(270:360:0.500000000000000);
\draw[very thick](\a+2,-5+\b)--(\a+4,-5+\b);
\draw[very thick](\a+3,-4+\b)--(\a+3,-6+\b);
\draw[very thick](\a+4,-5+\b)--(\a+6,-5+\b);
\draw[very thick](\a+3,-6+\b)--(\a+3,-8+\b);
\draw[very thick](\a+3,-8+\b)--(\a+3,-10+\b);
\draw[very thick](\a+5,-9.50000000000000+\b)--(\a+5,-10+\b);
\draw[very thick](\a+5.50000000000000,-9+\b)--(\a+6,-9+\b);
\draw[very thick](\a+5.50000000000000,-9+\b)arc(90:180:0.500000000000000);
\draw[very thick](\a+2,-11+\b)--(\a+4,-11+\b);
\draw[very thick](\a+3,-10+\b)--(\a+3,-12+\b);
\draw[very thick](\a+5,-10+\b)--(\a+5,-10.5000000000000+\b);
\draw[very thick](\a+4,-11+\b)--(\a+4.50000000000000,-11+\b);
\draw[very thick](\a+4.50000000000000,-11+\b)arc(270:360:0.500000000000000);
\draw[very thick](\a+3,-12+\b)--(\a+3,-12.5000000000000+\b);
\draw[very thick](\a+2,-13+\b)--(\a+2.50000000000000,-13+\b);
\draw[very thick](\a+2.50000000000000,-13+\b)arc(270:360:0.500000000000000);

\node at (3,-3) {$\times$};
\node at (3,-7) {$\vdots$};
\node at (13,-13) {$\times$};
\node at (13,-7) {$\vdots$};
\draw[very thick,->](7,-8)--(9,-8);
\end{tikzpicture}
\qquad
\begin{tikzpicture}[scale=0.200000000000000]
\def\a{0};
\def\b{0};
\draw[step=2, color=green](\a+2,-14+\b)grid(\a+6,-2+\b);
\draw[very thick](\a+5,-2+\b)--(\a+5,-4+\b);
\draw[very thick](\a+2,-5+\b)--(\a+4,-5+\b);
\draw[very thick](\a+4,-5+\b)--(\a+6,-5+\b);
\draw[very thick](\a+5,-4+\b)--(\a+5,-6+\b);
\draw[very thick](\a+5,-6+\b)--(\a+5,-8+\b);
\draw[very thick](\a+3,-9.50000000000000+\b)--(\a+3,-10+\b);
\draw[very thick](\a+3.50000000000000,-9+\b)--(\a+4,-9+\b);
\draw[very thick](\a+3.50000000000000,-9+\b)arc(90:180:0.500000000000000);
\draw[very thick](\a+4,-9+\b)--(\a+6,-9+\b);
\draw[very thick](\a+5,-8+\b)--(\a+5,-10+\b);
\draw[very thick](\a+3,-10+\b)--(\a+3,-10.5000000000000+\b);
\draw[very thick](\a+2,-11+\b)--(\a+2.50000000000000,-11+\b);
\draw[very thick](\a+2.50000000000000,-11+\b)arc(270:360:0.500000000000000);
\draw[very thick](\a+5,-10+\b)--(\a+5,-12+\b);
\draw[very thick](\a+3,-13.5000000000000+\b)--(\a+3,-14+\b);
\draw[very thick](\a+3.50000000000000,-13+\b)--(\a+4,-13+\b);
\draw[very thick](\a+3.50000000000000,-13+\b)arc(90:180:0.500000000000000);
\draw[very thick](\a+5,-12+\b)--(\a+5,-12.5000000000000+\b);
\draw[very thick](\a+4,-13+\b)--(\a+4.50000000000000,-13+\b);
\draw[very thick](\a+4.50000000000000,-13+\b)arc(270:360:0.500000000000000);

\def\a{8};
\def\b{0};
\draw[step=2, color=green](\a+2,-14+\b)grid(\a+6,-2+\b);
\draw[very thick](\a+3,-3.50000000000000+\b)--(\a+3,-4+\b);
\draw[very thick](\a+3.50000000000000,-3+\b)--(\a+4,-3+\b);
\draw[very thick](\a+3.50000000000000,-3+\b)arc(90:180:0.500000000000000);
\draw[very thick](\a+5,-2+\b)--(\a+5,-2.50000000000000+\b);
\draw[very thick](\a+4,-3+\b)--(\a+4.50000000000000,-3+\b);
\draw[very thick](\a+4.50000000000000,-3+\b)arc(270:360:0.500000000000000);
\draw[very thick](\a+2,-5+\b)--(\a+4,-5+\b);
\draw[very thick](\a+3,-4+\b)--(\a+3,-6+\b);
\draw[very thick](\a+4,-5+\b)--(\a+6,-5+\b);
\draw[very thick](\a+3,-6+\b)--(\a+3,-8+\b);
\draw[very thick](\a+3,-8+\b)--(\a+3,-10+\b);
\draw[very thick](\a+5,-9.50000000000000+\b)--(\a+5,-10+\b);
\draw[very thick](\a+5.50000000000000,-9+\b)--(\a+6,-9+\b);
\draw[very thick](\a+5.50000000000000,-9+\b)arc(90:180:0.500000000000000);
\draw[very thick](\a+2,-11+\b)--(\a+4,-11+\b);
\draw[very thick](\a+3,-10+\b)--(\a+3,-12+\b);
\draw[very thick](\a+5,-10+\b)--(\a+5,-10.5000000000000+\b);
\draw[very thick](\a+4,-11+\b)--(\a+4.50000000000000,-11+\b);
\draw[very thick](\a+4.50000000000000,-11+\b)arc(270:360:0.500000000000000);
\draw[very thick](\a+3,-12+\b)--(\a+3,-14+\b);

\node at (3,-3) {$\times$};
\node at (3,-7) {$\vdots$};
\node at (13,-13) {$\times$};
\node at (13,-7) {$\vdots$};
\draw[very thick,->](7,-8)--(9,-8);
\end{tikzpicture}
\qquad
\begin{tikzpicture}[scale=0.200000000000000]
\def\a{0};
\def\b{0};
\draw[step=2, color=green](\a+2,-14+\b)grid(\a+6,-2+\b);
\draw[very thick](\a+5,-3.50000000000000+\b)--(\a+5,-4+\b);
\draw[very thick](\a+5.50000000000000,-3+\b)--(\a+6,-3+\b);
\draw[very thick](\a+5.50000000000000,-3+\b)arc(90:180:0.500000000000000);
\draw[very thick](\a+2,-5+\b)--(\a+4,-5+\b);
\draw[very thick](\a+4,-5+\b)--(\a+6,-5+\b);
\draw[very thick](\a+5,-4+\b)--(\a+5,-6+\b);
\draw[very thick](\a+5,-6+\b)--(\a+5,-8+\b);
\draw[very thick](\a+3,-9.50000000000000+\b)--(\a+3,-10+\b);
\draw[very thick](\a+3.50000000000000,-9+\b)--(\a+4,-9+\b);
\draw[very thick](\a+3.50000000000000,-9+\b)arc(90:180:0.500000000000000);
\draw[very thick](\a+4,-9+\b)--(\a+6,-9+\b);
\draw[very thick](\a+5,-8+\b)--(\a+5,-10+\b);
\draw[very thick](\a+3,-10+\b)--(\a+3,-10.5000000000000+\b);
\draw[very thick](\a+2,-11+\b)--(\a+2.50000000000000,-11+\b);
\draw[very thick](\a+2.50000000000000,-11+\b)arc(270:360:0.500000000000000);
\draw[very thick](\a+5,-10+\b)--(\a+5,-12+\b);
\draw[very thick](\a+3,-13.5000000000000+\b)--(\a+3,-14+\b);
\draw[very thick](\a+3.50000000000000,-13+\b)--(\a+4,-13+\b);
\draw[very thick](\a+3.50000000000000,-13+\b)arc(90:180:0.500000000000000);
\draw[very thick](\a+5,-12+\b)--(\a+5,-12.5000000000000+\b);
\draw[very thick](\a+4,-13+\b)--(\a+4.50000000000000,-13+\b);
\draw[very thick](\a+4.50000000000000,-13+\b)arc(270:360:0.500000000000000);

\def\a{8};
\def\b{0};
\draw[step=2, color=green](\a+2,-14+\b)grid(\a+6,-2+\b);
\draw[very thick](\a+3,-3.50000000000000+\b)--(\a+3,-4+\b);
\draw[very thick](\a+3.50000000000000,-3+\b)--(\a+4,-3+\b);
\draw[very thick](\a+3.50000000000000,-3+\b)arc(90:180:0.500000000000000);
\draw[very thick](\a+4,-3+\b)--(\a+6,-3+\b);
\draw[very thick](\a+2,-5+\b)--(\a+4,-5+\b);
\draw[very thick](\a+3,-4+\b)--(\a+3,-6+\b);
\draw[very thick](\a+4,-5+\b)--(\a+6,-5+\b);
\draw[very thick](\a+3,-6+\b)--(\a+3,-8+\b);
\draw[very thick](\a+3,-8+\b)--(\a+3,-10+\b);
\draw[very thick](\a+5,-9.50000000000000+\b)--(\a+5,-10+\b);
\draw[very thick](\a+5.50000000000000,-9+\b)--(\a+6,-9+\b);
\draw[very thick](\a+5.50000000000000,-9+\b)arc(90:180:0.500000000000000);
\draw[very thick](\a+2,-11+\b)--(\a+4,-11+\b);
\draw[very thick](\a+3,-10+\b)--(\a+3,-12+\b);
\draw[very thick](\a+5,-10+\b)--(\a+5,-10.5000000000000+\b);
\draw[very thick](\a+4,-11+\b)--(\a+4.50000000000000,-11+\b);
\draw[very thick](\a+4.50000000000000,-11+\b)arc(270:360:0.500000000000000);
\draw[very thick](\a+3,-12+\b)--(\a+3,-14+\b);

\node at (3,-3) {$\times$};
\node at (3,-7) {$\vdots$};
\node at (13,-13) {$\times$};
\node at (13,-7) {$\vdots$};
\draw[very thick,->](7,-8)--(9,-8);
\end{tikzpicture}
\qquad
\begin{tikzpicture}[scale=0.200000000000000]
\def\a{0};
\def\b{0};
\draw[step=2, color=green](\a+2,-14+\b)grid(\a+6,-2+\b);
\draw[very thick](\a+5,-3.50000000000000+\b)--(\a+5,-4+\b);
\draw[very thick](\a+5.50000000000000,-3+\b)--(\a+6,-3+\b);
\draw[very thick](\a+5.50000000000000,-3+\b)arc(90:180:0.500000000000000);
\draw[very thick](\a+2,-5+\b)--(\a+4,-5+\b);
\draw[very thick](\a+4,-5+\b)--(\a+6,-5+\b);
\draw[very thick](\a+5,-4+\b)--(\a+5,-6+\b);
\draw[very thick](\a+5,-6+\b)--(\a+5,-8+\b);
\draw[very thick](\a+3,-9.50000000000000+\b)--(\a+3,-10+\b);
\draw[very thick](\a+3.50000000000000,-9+\b)--(\a+4,-9+\b);
\draw[very thick](\a+3.50000000000000,-9+\b)arc(90:180:0.500000000000000);
\draw[very thick](\a+4,-9+\b)--(\a+6,-9+\b);
\draw[very thick](\a+5,-8+\b)--(\a+5,-10+\b);
\draw[very thick](\a+3,-10+\b)--(\a+3,-10.5000000000000+\b);
\draw[very thick](\a+2,-11+\b)--(\a+2.50000000000000,-11+\b);
\draw[very thick](\a+2.50000000000000,-11+\b)arc(270:360:0.500000000000000);
\draw[very thick](\a+5,-10+\b)--(\a+5,-12+\b);
\draw[very thick](\a+2,-13+\b)--(\a+4,-13+\b);
\draw[very thick](\a+5,-12+\b)--(\a+5,-12.5000000000000+\b);
\draw[very thick](\a+4,-13+\b)--(\a+4.50000000000000,-13+\b);
\draw[very thick](\a+4.50000000000000,-13+\b)arc(270:360:0.500000000000000);

\def\a{8};
\def\b{0};
\draw[step=2, color=green](\a+2,-14+\b)grid(\a+6,-2+\b);
\draw[very thick](\a+3,-3.50000000000000+\b)--(\a+3,-4+\b);
\draw[very thick](\a+3.50000000000000,-3+\b)--(\a+4,-3+\b);
\draw[very thick](\a+3.50000000000000,-3+\b)arc(90:180:0.500000000000000);
\draw[very thick](\a+4,-3+\b)--(\a+6,-3+\b);
\draw[very thick](\a+2,-5+\b)--(\a+4,-5+\b);
\draw[very thick](\a+3,-4+\b)--(\a+3,-6+\b);
\draw[very thick](\a+4,-5+\b)--(\a+6,-5+\b);
\draw[very thick](\a+3,-6+\b)--(\a+3,-8+\b);
\draw[very thick](\a+3,-8+\b)--(\a+3,-10+\b);
\draw[very thick](\a+5,-9.50000000000000+\b)--(\a+5,-10+\b);
\draw[very thick](\a+5.50000000000000,-9+\b)--(\a+6,-9+\b);
\draw[very thick](\a+5.50000000000000,-9+\b)arc(90:180:0.500000000000000);
\draw[very thick](\a+2,-11+\b)--(\a+4,-11+\b);
\draw[very thick](\a+3,-10+\b)--(\a+3,-12+\b);
\draw[very thick](\a+5,-10+\b)--(\a+5,-10.5000000000000+\b);
\draw[very thick](\a+4,-11+\b)--(\a+4.50000000000000,-11+\b);
\draw[very thick](\a+4.50000000000000,-11+\b)arc(270:360:0.500000000000000);
\draw[very thick](\a+3,-12+\b)--(\a+3,-12.5000000000000+\b);
\draw[very thick](\a+2,-13+\b)--(\a+2.50000000000000,-13+\b);
\draw[very thick](\a+2.50000000000000,-13+\b)arc(270:360:0.500000000000000);

\node at (3,-3) {$\times$};
\node at (3,-7) {$\vdots$};
\node at (13,-13) {$\times$};
\node at (13,-7) {$\vdots$};
\draw[very thick,->](7,-8)--(9,-8);
\end{tikzpicture}
\caption{Step (2) of \Cref{def:pop}}
\label{fig:bpd-pop-situation1}
\end{figure}

\begin{figure}[h!]
\centering
\begin{tikzpicture}[scale=0.200000000000000]
\def\a{0};
\def\b{0};
\draw[step=2, color=green](\a+2,-14+\b)grid(\a+6,-2+\b);
\draw[very thick](\a+5,-2+\b)--(\a+5,-4+\b);
\draw[very thick](\a+2,-5+\b)--(\a+4,-5+\b);
\draw[very thick](\a+4,-5+\b)--(\a+6,-5+\b);
\draw[very thick](\a+5,-4+\b)--(\a+5,-6+\b);
\draw[very thick](\a+5,-6+\b)--(\a+5,-8+\b);
\draw[very thick](\a+3,-9.50000000000000+\b)--(\a+3,-10+\b);
\draw[very thick](\a+3.50000000000000,-9+\b)--(\a+4,-9+\b);
\draw[very thick](\a+3.50000000000000,-9+\b)arc(90:180:0.500000000000000);
\draw[very thick](\a+4,-9+\b)--(\a+6,-9+\b);
\draw[very thick](\a+5,-8+\b)--(\a+5,-10+\b);
\draw[very thick](\a+3,-10+\b)--(\a+3,-10.5000000000000+\b);
\draw[very thick](\a+2,-11+\b)--(\a+2.50000000000000,-11+\b);
\draw[very thick](\a+2.50000000000000,-11+\b)arc(270:360:0.500000000000000);
\draw[very thick](\a+5,-10+\b)--(\a+5,-12+\b);
\draw[very thick](\a+3,-13.5000000000000+\b)--(\a+3,-14+\b);
\draw[very thick](\a+3.50000000000000,-13+\b)--(\a+4,-13+\b);
\draw[very thick](\a+3.50000000000000,-13+\b)arc(90:180:0.500000000000000);
\draw[very thick](\a+4,-13+\b)--(\a+6,-13+\b);
\draw[very thick](\a+5,-12+\b)--(\a+5,-14+\b);

\def\a{8};
\def\b{0};
\draw[step=2, color=green](\a+2,-14+\b)grid(\a+6,-2+\b);
\draw[very thick](\a+5,-2+\b)--(\a+5,-4+\b);
\draw[very thick](\a+2,-5+\b)--(\a+4,-5+\b);
\draw[very thick](\a+4,-5+\b)--(\a+6,-5+\b);
\draw[very thick](\a+5,-4+\b)--(\a+5,-6+\b);
\draw[very thick](\a+5,-6+\b)--(\a+5,-8+\b);
\draw[very thick](\a+3,-9.5+\b)--(\a+3,-10+\b);
\draw[very thick](\a+3.5,-9+\b)--(\a+4,-9+\b);
\draw[very thick](\a+3.5,-9+\b)arc(90:180:0.5);
\draw[very thick](\a+4,-9+\b)--(\a+6,-9+\b);
\draw[very thick](\a+5,-8+\b)--(\a+5,-10+\b);
\draw[very thick](\a+3,-10+\b)--(\a+3,-10.5+\b);
\draw[very thick](\a+2,-11+\b)--(\a+2.5,-11+\b);
\draw[very thick](\a+2.5,-11+\b)arc(270:360:0.5);
\draw[very thick](\a+5,-10+\b)--(\a+5,-12+\b);
\draw[very thick](\a+3,-13.5+\b)--(\a+3,-14+\b);
\draw[very thick](\a+3.5,-13+\b)--(\a+4,-13+\b);
\draw[very thick](\a+3.5,-13+\b)arc(90:180:0.5);
\draw[very thick](\a+5,-12+\b)--(\a+5,-12.5+\b);
\draw[very thick](\a+4,-13+\b)--(\a+4.5,-13+\b);
\draw[very thick](\a+4.5,-13+\b)arc(270:360:0.5);
\draw[very thick](\a+5,-13.5+\b)--(\a+5,-14+\b);
\draw[very thick](\a+5.5,-13+\b)--(\a+6,-13+\b);
\draw[very thick](\a+5.5,-13+\b)arc(90:180:0.5);

\def\a{16};
\def\b{0};
\draw[step=2, color=green](\a+2,-14+\b)grid(\a+6,-2+\b);
\draw[very thick](\a+3,-3.50000000000000+\b)--(\a+3,-4+\b);
\draw[very thick](\a+3.50000000000000,-3+\b)--(\a+4,-3+\b);
\draw[very thick](\a+3.50000000000000,-3+\b)arc(90:180:0.500000000000000);
\draw[very thick](\a+5,-2+\b)--(\a+5,-2.50000000000000+\b);
\draw[very thick](\a+4,-3+\b)--(\a+4.50000000000000,-3+\b);
\draw[very thick](\a+4.50000000000000,-3+\b)arc(270:360:0.500000000000000);
\draw[very thick](\a+2,-5+\b)--(\a+4,-5+\b);
\draw[very thick](\a+3,-4+\b)--(\a+3,-6+\b);
\draw[very thick](\a+4,-5+\b)--(\a+6,-5+\b);
\draw[very thick](\a+3,-6+\b)--(\a+3,-8+\b);
\draw[very thick](\a+3,-8+\b)--(\a+3,-10+\b);
\draw[very thick](\a+5,-9.50000000000000+\b)--(\a+5,-10+\b);
\draw[very thick](\a+5.50000000000000,-9+\b)--(\a+6,-9+\b);
\draw[very thick](\a+5.50000000000000,-9+\b)arc(90:180:0.500000000000000);
\draw[very thick](\a+2,-11+\b)--(\a+4,-11+\b);
\draw[very thick](\a+3,-10+\b)--(\a+3,-12+\b);
\draw[very thick](\a+5,-10+\b)--(\a+5,-10.5000000000000+\b);
\draw[very thick](\a+4,-11+\b)--(\a+4.50000000000000,-11+\b);
\draw[very thick](\a+4.50000000000000,-11+\b)arc(270:360:0.500000000000000);
\draw[very thick](\a+3,-12+\b)--(\a+3,-14+\b);
\draw[very thick](\a+5,-13.5000000000000+\b)--(\a+5,-14+\b);
\draw[very thick](\a+5.50000000000000,-13+\b)--(\a+6,-13+\b);
\draw[very thick](\a+5.50000000000000,-13+\b)arc(90:180:0.500000000000000);

\node at (3,-3) {$\times$};
\node at (11,-3) {$\times$};
\node at (3,-7) {$\vdots$};
\node at (11,-7) {$\vdots$};
\node at (21,-7) {$\vdots$};
\draw[very thick,->](7,-8)--(9,-8);
\draw[very thick,->](15,-8)--(17,-8);
\end{tikzpicture}
\qquad
\begin{tikzpicture}[scale=0.200000000000000]
\def\a{0};
\def\b{0};
\draw[step=2, color=green](\a+2,-14+\b)grid(\a+6,-2+\b);
\draw[very thick](\a+5,-3.50000000000000+\b)--(\a+5,-4+\b);
\draw[very thick](\a+5.50000000000000,-3+\b)--(\a+6,-3+\b);
\draw[very thick](\a+5.50000000000000,-3+\b)arc(90:180:0.500000000000000);
\draw[very thick](\a+2,-5+\b)--(\a+4,-5+\b);
\draw[very thick](\a+4,-5+\b)--(\a+6,-5+\b);
\draw[very thick](\a+5,-4+\b)--(\a+5,-6+\b);
\draw[very thick](\a+5,-6+\b)--(\a+5,-8+\b);
\draw[very thick](\a+3,-9.50000000000000+\b)--(\a+3,-10+\b);
\draw[very thick](\a+3.50000000000000,-9+\b)--(\a+4,-9+\b);
\draw[very thick](\a+3.50000000000000,-9+\b)arc(90:180:0.500000000000000);
\draw[very thick](\a+4,-9+\b)--(\a+6,-9+\b);
\draw[very thick](\a+5,-8+\b)--(\a+5,-10+\b);
\draw[very thick](\a+3,-10+\b)--(\a+3,-10.5000000000000+\b);
\draw[very thick](\a+2,-11+\b)--(\a+2.50000000000000,-11+\b);
\draw[very thick](\a+2.50000000000000,-11+\b)arc(270:360:0.500000000000000);
\draw[very thick](\a+5,-10+\b)--(\a+5,-12+\b);
\draw[very thick](\a+3,-13.5000000000000+\b)--(\a+3,-14+\b);
\draw[very thick](\a+3.50000000000000,-13+\b)--(\a+4,-13+\b);
\draw[very thick](\a+3.50000000000000,-13+\b)arc(90:180:0.500000000000000);
\draw[very thick](\a+4,-13+\b)--(\a+6,-13+\b);
\draw[very thick](\a+5,-12+\b)--(\a+5,-14+\b);

\def\a{8};
\def\b{0};
\draw[step=2, color=green](\a+2,-14+\b)grid(\a+6,-2+\b);
\draw[very thick](\a+5,-3.5+\b)--(\a+5,-4+\b);
\draw[very thick](\a+5.5,-3+\b)--(\a+6,-3+\b);
\draw[very thick](\a+5.5,-3+\b)arc(90:180:0.5);
\draw[very thick](\a+2,-5+\b)--(\a+4,-5+\b);
\draw[very thick](\a+4,-5+\b)--(\a+6,-5+\b);
\draw[very thick](\a+5,-4+\b)--(\a+5,-6+\b);
\draw[very thick](\a+5,-6+\b)--(\a+5,-8+\b);
\draw[very thick](\a+3,-9.5+\b)--(\a+3,-10+\b);
\draw[very thick](\a+3.5,-9+\b)--(\a+4,-9+\b);
\draw[very thick](\a+3.5,-9+\b)arc(90:180:0.5);
\draw[very thick](\a+4,-9+\b)--(\a+6,-9+\b);
\draw[very thick](\a+5,-8+\b)--(\a+5,-10+\b);
\draw[very thick](\a+3,-10+\b)--(\a+3,-10.5+\b);
\draw[very thick](\a+2,-11+\b)--(\a+2.5,-11+\b);
\draw[very thick](\a+2.5,-11+\b)arc(270:360:0.5);
\draw[very thick](\a+5,-10+\b)--(\a+5,-12+\b);
\draw[very thick](\a+3,-13.5+\b)--(\a+3,-14+\b);
\draw[very thick](\a+3.5,-13+\b)--(\a+4,-13+\b);
\draw[very thick](\a+3.5,-13+\b)arc(90:180:0.5);
\draw[very thick](\a+5,-12+\b)--(\a+5,-12.5+\b);
\draw[very thick](\a+4,-13+\b)--(\a+4.5,-13+\b);
\draw[very thick](\a+4.5,-13+\b)arc(270:360:0.5);
\draw[very thick](\a+5,-13.5+\b)--(\a+5,-14+\b);
\draw[very thick](\a+5.5,-13+\b)--(\a+6,-13+\b);
\draw[very thick](\a+5.5,-13+\b)arc(90:180:0.5);

\def\a{16};
\def\b{0};
\draw[step=2, color=green](\a+2,-14+\b)grid(\a+6,-2+\b);
\draw[very thick](\a+3,-3.50000000000000+\b)--(\a+3,-4+\b);
\draw[very thick](\a+3.50000000000000,-3+\b)--(\a+4,-3+\b);
\draw[very thick](\a+3.50000000000000,-3+\b)arc(90:180:0.500000000000000);
\draw[very thick](\a+4,-3+\b)--(\a+6,-3+\b);
\draw[very thick](\a+2,-5+\b)--(\a+4,-5+\b);
\draw[very thick](\a+3,-4+\b)--(\a+3,-6+\b);
\draw[very thick](\a+4,-5+\b)--(\a+6,-5+\b);
\draw[very thick](\a+3,-6+\b)--(\a+3,-8+\b);
\draw[very thick](\a+3,-8+\b)--(\a+3,-10+\b);
\draw[very thick](\a+5,-9.50000000000000+\b)--(\a+5,-10+\b);
\draw[very thick](\a+5.50000000000000,-9+\b)--(\a+6,-9+\b);
\draw[very thick](\a+5.50000000000000,-9+\b)arc(90:180:0.500000000000000);
\draw[very thick](\a+2,-11+\b)--(\a+4,-11+\b);
\draw[very thick](\a+3,-10+\b)--(\a+3,-12+\b);
\draw[very thick](\a+5,-10+\b)--(\a+5,-10.5000000000000+\b);
\draw[very thick](\a+4,-11+\b)--(\a+4.50000000000000,-11+\b);
\draw[very thick](\a+4.50000000000000,-11+\b)arc(270:360:0.500000000000000);
\draw[very thick](\a+3,-12+\b)--(\a+3,-14+\b);
\draw[very thick](\a+5,-13.5000000000000+\b)--(\a+5,-14+\b);
\draw[very thick](\a+5.50000000000000,-13+\b)--(\a+6,-13+\b);
\draw[very thick](\a+5.50000000000000,-13+\b)arc(90:180:0.500000000000000);

\node at (3,-3) {$\times$};
\node at (11,-3) {$\times$};
\node at (3,-7) {$\vdots$};
\node at (11,-7) {$\vdots$};
\node at (21,-7) {$\vdots$};
\draw[very thick,->](7,-8)--(9,-8);
\draw[very thick,->](15,-8)--(17,-8);
\end{tikzpicture}
\caption{Step (4) of \Cref{def:pop}}
\label{fig:bpd-pop-situation2}
\end{figure}

\begin{Definition}\cite{GaoHuang-bijection}\label{def:bijection-PD-BPD}
For $D\in\BPD(w)$ with $\ell(w)=\ell$, define \[\varphi(D)=\big(\mathbf{r}=(r_1,\ldots,r_{\ell}),\mathbf{i}=(i_1,\ldots,i_{\ell})\big)\]
where $\mathrm{pop}(\nabla^{k-1}D)=(r_k,i_k)$ for $k=1,\ldots,\ell$.
\end{Definition}
\begin{Example}
Consider $D\in\BPD(w)$ with $w=2157346$ in
\Cref{fig:pop-example}, which shows the steps of \Cref{def:pop} that
result in $\nabla D\in\BPD(2147356)$ and
$\mathrm{pop}(D)=(4,1)$ where $r=4$, $i=1$.
The edge labels in the figure are the different values of $y$ in the algorithm.
If we continue the process as in \Cref{def:bijection-PD-BPD}, we eventually end up with the compatible sequence
\[\big(\mathbf{r}=(4,1,6,5,3,4),\mathbf{i}=(1,1,2,3,3,4)\big)\]
that corresponds to the pipe dream as in \Cref{fig:bijection-example}.

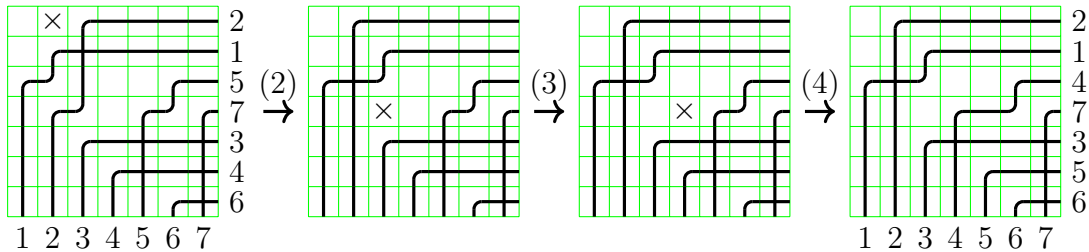
\begin{figure}[h!]
\centering
\begin{tikzpicture}[scale=0.200000000000000]
\def\a{0};
\def\b{0};
\draw[step=2, color=green](\a+2,-16+\b)grid(\a+16,-2+\b);
\draw[very thick](\a+7,-3.50000000000000+\b)--(\a+7,-4+\b);
\draw[very thick](\a+7.50000000000000,-3+\b)--(\a+8,-3+\b);
\draw[very thick](\a+7.50000000000000,-3+\b)arc(90:180:0.500000000000000);
\draw[very thick](\a+8,-3+\b)--(\a+10,-3+\b);
\draw[very thick](\a+10,-3+\b)--(\a+12,-3+\b);
\draw[very thick](\a+12,-3+\b)--(\a+14,-3+\b);
\draw[very thick](\a+14,-3+\b)--(\a+16,-3+\b);
\draw[very thick](\a+5,-5.50000000000000+\b)--(\a+5,-6+\b);
\draw[very thick](\a+5.50000000000000,-5+\b)--(\a+6,-5+\b);
\draw[very thick](\a+5.50000000000000,-5+\b)arc(90:180:0.500000000000000);
\draw[very thick](\a+6,-5+\b)--(\a+8,-5+\b);
\draw[very thick](\a+7,-4+\b)--(\a+7,-6+\b);
\draw[very thick](\a+8,-5+\b)--(\a+10,-5+\b);
\draw[very thick](\a+10,-5+\b)--(\a+12,-5+\b);
\draw[very thick](\a+12,-5+\b)--(\a+14,-5+\b);
\draw[very thick](\a+14,-5+\b)--(\a+16,-5+\b);
\draw[very thick](\a+3,-7.50000000000000+\b)--(\a+3,-8+\b);
\draw[very thick](\a+3.50000000000000,-7+\b)--(\a+4,-7+\b);
\draw[very thick](\a+3.50000000000000,-7+\b)arc(90:180:0.500000000000000);
\draw[very thick](\a+5,-6+\b)--(\a+5,-6.50000000000000+\b);
\draw[very thick](\a+4,-7+\b)--(\a+4.50000000000000,-7+\b);
\draw[very thick](\a+4.50000000000000,-7+\b)arc(270:360:0.500000000000000);
\draw[very thick](\a+7,-6+\b)--(\a+7,-8+\b);
\draw[very thick](\a+13,-7.50000000000000+\b)--(\a+13,-8+\b);
\draw[very thick](\a+13.5000000000000,-7+\b)--(\a+14,-7+\b);
\draw[very thick](\a+13.5000000000000,-7+\b)arc(90:180:0.500000000000000);
\draw[very thick](\a+14,-7+\b)--(\a+16,-7+\b);
\draw[very thick](\a+3,-8+\b)--(\a+3,-10+\b);
\draw[very thick](\a+5,-9.50000000000000+\b)--(\a+5,-10+\b);
\draw[very thick](\a+5.50000000000000,-9+\b)--(\a+6,-9+\b);
\draw[very thick](\a+5.50000000000000,-9+\b)arc(90:180:0.500000000000000);
\draw[very thick](\a+7,-8+\b)--(\a+7,-8.50000000000000+\b);
\draw[very thick](\a+6,-9+\b)--(\a+6.50000000000000,-9+\b);
\draw[very thick](\a+6.50000000000000,-9+\b)arc(270:360:0.500000000000000);
\draw[very thick](\a+11,-9.50000000000000+\b)--(\a+11,-10+\b);
\draw[very thick](\a+11.5000000000000,-9+\b)--(\a+12,-9+\b);
\draw[very thick](\a+11.5000000000000,-9+\b)arc(90:180:0.500000000000000);
\draw[very thick](\a+13,-8+\b)--(\a+13,-8.50000000000000+\b);
\draw[very thick](\a+12,-9+\b)--(\a+12.5000000000000,-9+\b);
\draw[very thick](\a+12.5000000000000,-9+\b)arc(270:360:0.500000000000000);
\draw[very thick](\a+15,-9.50000000000000+\b)--(\a+15,-10+\b);
\draw[very thick](\a+15.5000000000000,-9+\b)--(\a+16,-9+\b);
\draw[very thick](\a+15.5000000000000,-9+\b)arc(90:180:0.500000000000000);
\draw[very thick](\a+3,-10+\b)--(\a+3,-12+\b);
\draw[very thick](\a+5,-10+\b)--(\a+5,-12+\b);
\draw[very thick](\a+7,-11.5000000000000+\b)--(\a+7,-12+\b);
\draw[very thick](\a+7.50000000000000,-11+\b)--(\a+8,-11+\b);
\draw[very thick](\a+7.50000000000000,-11+\b)arc(90:180:0.500000000000000);
\draw[very thick](\a+8,-11+\b)--(\a+10,-11+\b);
\draw[very thick](\a+10,-11+\b)--(\a+12,-11+\b);
\draw[very thick](\a+11,-10+\b)--(\a+11,-12+\b);
\draw[very thick](\a+12,-11+\b)--(\a+14,-11+\b);
\draw[very thick](\a+14,-11+\b)--(\a+16,-11+\b);
\draw[very thick](\a+15,-10+\b)--(\a+15,-12+\b);
\draw[very thick](\a+3,-12+\b)--(\a+3,-14+\b);
\draw[very thick](\a+5,-12+\b)--(\a+5,-14+\b);
\draw[very thick](\a+7,-12+\b)--(\a+7,-14+\b);
\draw[very thick](\a+9,-13.5000000000000+\b)--(\a+9,-14+\b);
\draw[very thick](\a+9.50000000000000,-13+\b)--(\a+10,-13+\b);
\draw[very thick](\a+9.50000000000000,-13+\b)arc(90:180:0.500000000000000);
\draw[very thick](\a+10,-13+\b)--(\a+12,-13+\b);
\draw[very thick](\a+11,-12+\b)--(\a+11,-14+\b);
\draw[very thick](\a+12,-13+\b)--(\a+14,-13+\b);
\draw[very thick](\a+14,-13+\b)--(\a+16,-13+\b);
\draw[very thick](\a+15,-12+\b)--(\a+15,-14+\b);
\draw[very thick](\a+3,-14+\b)--(\a+3,-16+\b);
\draw[very thick](\a+5,-14+\b)--(\a+5,-16+\b);
\draw[very thick](\a+7,-14+\b)--(\a+7,-16+\b);
\draw[very thick](\a+9,-14+\b)--(\a+9,-16+\b);
\draw[very thick](\a+11,-14+\b)--(\a+11,-16+\b);
\draw[very thick](\a+13,-15.5000000000000+\b)--(\a+13,-16+\b);
\draw[very thick](\a+13.5000000000000,-15+\b)--(\a+14,-15+\b);
\draw[very thick](\a+13.5000000000000,-15+\b)arc(90:180:0.500000000000000);
\draw[very thick](\a+14,-15+\b)--(\a+16,-15+\b);
\draw[very thick](\a+15,-14+\b)--(\a+15,-16+\b);
\node[right] at (\a+16,-3+\b) {$2$};
\node[right] at (\a+16,-5+\b) {$1$};
\node[right] at (\a+16,-7+\b) {$5$};
\node[right] at (\a+16,-9+\b) {$7$};
\node[right] at (\a+16,-11+\b) {$3$};
\node[right] at (\a+16,-13+\b) {$4$};
\node[right] at (\a+16,-15+\b) {$6$};
\node[below] at (\a+3,-16+\b) {$1$};
\node[below] at (\a+5,-16+\b) {$2$};
\node[below] at (\a+7,-16+\b) {$3$};
\node[below] at (\a+9,-16+\b) {$4$};
\node[below] at (\a+11,-16+\b) {$5$};
\node[below] at (\a+13,-16+\b) {$6$};
\node[below] at (\a+15,-16+\b) {$7$};

\def\a{20};
\def\b{0};
\draw[step=2, color=green](\a+2,-16+\b)grid(\a+16,-2+\b);
\draw[very thick](\a+5,-3.50000000000000+\b)--(\a+5,-4+\b);
\draw[very thick](\a+5.50000000000000,-3+\b)--(\a+6,-3+\b);
\draw[very thick](\a+5.50000000000000,-3+\b)arc(90:180:0.500000000000000);
\draw[very thick](\a+6,-3+\b)--(\a+8,-3+\b);
\draw[very thick](\a+8,-3+\b)--(\a+10,-3+\b);
\draw[very thick](\a+10,-3+\b)--(\a+12,-3+\b);
\draw[very thick](\a+12,-3+\b)--(\a+14,-3+\b);
\draw[very thick](\a+14,-3+\b)--(\a+16,-3+\b);
\draw[very thick](\a+5,-4+\b)--(\a+5,-6+\b);
\draw[very thick](\a+7,-5.50000000000000+\b)--(\a+7,-6+\b);
\draw[very thick](\a+7.50000000000000,-5+\b)--(\a+8,-5+\b);
\draw[very thick](\a+7.50000000000000,-5+\b)arc(90:180:0.500000000000000);
\draw[very thick](\a+8,-5+\b)--(\a+10,-5+\b);
\draw[very thick](\a+10,-5+\b)--(\a+12,-5+\b);
\draw[very thick](\a+12,-5+\b)--(\a+14,-5+\b);
\draw[very thick](\a+14,-5+\b)--(\a+16,-5+\b);
\draw[very thick](\a+3,-7.50000000000000+\b)--(\a+3,-8+\b);
\draw[very thick](\a+3.50000000000000,-7+\b)--(\a+4,-7+\b);
\draw[very thick](\a+3.50000000000000,-7+\b)arc(90:180:0.500000000000000);
\draw[very thick](\a+4,-7+\b)--(\a+6,-7+\b);
\draw[very thick](\a+5,-6+\b)--(\a+5,-8+\b);
\draw[very thick](\a+7,-6+\b)--(\a+7,-6.50000000000000+\b);
\draw[very thick](\a+6,-7+\b)--(\a+6.50000000000000,-7+\b);
\draw[very thick](\a+6.50000000000000,-7+\b)arc(270:360:0.500000000000000);
\draw[very thick](\a+13,-7.50000000000000+\b)--(\a+13,-8+\b);
\draw[very thick](\a+13.5000000000000,-7+\b)--(\a+14,-7+\b);
\draw[very thick](\a+13.5000000000000,-7+\b)arc(90:180:0.500000000000000);
\draw[very thick](\a+14,-7+\b)--(\a+16,-7+\b);
\draw[very thick](\a+3,-8+\b)--(\a+3,-10+\b);
\draw[very thick](\a+5,-8+\b)--(\a+5,-10+\b);
\draw[very thick](\a+11,-9.50000000000000+\b)--(\a+11,-10+\b);
\draw[very thick](\a+11.5000000000000,-9+\b)--(\a+12,-9+\b);
\draw[very thick](\a+11.5000000000000,-9+\b)arc(90:180:0.500000000000000);
\draw[very thick](\a+13,-8+\b)--(\a+13,-8.50000000000000+\b);
\draw[very thick](\a+12,-9+\b)--(\a+12.5000000000000,-9+\b);
\draw[very thick](\a+12.5000000000000,-9+\b)arc(270:360:0.500000000000000);
\draw[very thick](\a+15,-9.50000000000000+\b)--(\a+15,-10+\b);
\draw[very thick](\a+15.5000000000000,-9+\b)--(\a+16,-9+\b);
\draw[very thick](\a+15.5000000000000,-9+\b)arc(90:180:0.500000000000000);
\draw[very thick](\a+3,-10+\b)--(\a+3,-12+\b);
\draw[very thick](\a+5,-10+\b)--(\a+5,-12+\b);
\draw[very thick](\a+7,-11.5000000000000+\b)--(\a+7,-12+\b);
\draw[very thick](\a+7.50000000000000,-11+\b)--(\a+8,-11+\b);
\draw[very thick](\a+7.50000000000000,-11+\b)arc(90:180:0.500000000000000);
\draw[very thick](\a+8,-11+\b)--(\a+10,-11+\b);
\draw[very thick](\a+10,-11+\b)--(\a+12,-11+\b);
\draw[very thick](\a+11,-10+\b)--(\a+11,-12+\b);
\draw[very thick](\a+12,-11+\b)--(\a+14,-11+\b);
\draw[very thick](\a+14,-11+\b)--(\a+16,-11+\b);
\draw[very thick](\a+15,-10+\b)--(\a+15,-12+\b);
\draw[very thick](\a+3,-12+\b)--(\a+3,-14+\b);
\draw[very thick](\a+5,-12+\b)--(\a+5,-14+\b);
\draw[very thick](\a+7,-12+\b)--(\a+7,-14+\b);
\draw[very thick](\a+9,-13.5000000000000+\b)--(\a+9,-14+\b);
\draw[very thick](\a+9.50000000000000,-13+\b)--(\a+10,-13+\b);
\draw[very thick](\a+9.50000000000000,-13+\b)arc(90:180:0.500000000000000);
\draw[very thick](\a+10,-13+\b)--(\a+12,-13+\b);
\draw[very thick](\a+11,-12+\b)--(\a+11,-14+\b);
\draw[very thick](\a+12,-13+\b)--(\a+14,-13+\b);
\draw[very thick](\a+14,-13+\b)--(\a+16,-13+\b);
\draw[very thick](\a+15,-12+\b)--(\a+15,-14+\b);
\draw[very thick](\a+3,-14+\b)--(\a+3,-16+\b);
\draw[very thick](\a+5,-14+\b)--(\a+5,-16+\b);
\draw[very thick](\a+7,-14+\b)--(\a+7,-16+\b);
\draw[very thick](\a+9,-14+\b)--(\a+9,-16+\b);
\draw[very thick](\a+11,-14+\b)--(\a+11,-16+\b);
\draw[very thick](\a+13,-15.5000000000000+\b)--(\a+13,-16+\b);
\draw[very thick](\a+13.5000000000000,-15+\b)--(\a+14,-15+\b);
\draw[very thick](\a+13.5000000000000,-15+\b)arc(90:180:0.500000000000000);
\draw[very thick](\a+14,-15+\b)--(\a+16,-15+\b);
\draw[very thick](\a+15,-14+\b)--(\a+15,-16+\b);

\def\a{38};
\def\b{0};
\draw[step=2, color=green](\a+2,-16+\b)grid(\a+16,-2+\b);
\draw[very thick](\a+5,-3.50000000000000+\b)--(\a+5,-4+\b);
\draw[very thick](\a+5.50000000000000,-3+\b)--(\a+6,-3+\b);
\draw[very thick](\a+5.50000000000000,-3+\b)arc(90:180:0.500000000000000);
\draw[very thick](\a+6,-3+\b)--(\a+8,-3+\b);
\draw[very thick](\a+8,-3+\b)--(\a+10,-3+\b);
\draw[very thick](\a+10,-3+\b)--(\a+12,-3+\b);
\draw[very thick](\a+12,-3+\b)--(\a+14,-3+\b);
\draw[very thick](\a+14,-3+\b)--(\a+16,-3+\b);
\draw[very thick](\a+5,-4+\b)--(\a+5,-6+\b);
\draw[very thick](\a+7,-5.50000000000000+\b)--(\a+7,-6+\b);
\draw[very thick](\a+7.50000000000000,-5+\b)--(\a+8,-5+\b);
\draw[very thick](\a+7.50000000000000,-5+\b)arc(90:180:0.500000000000000);
\draw[very thick](\a+8,-5+\b)--(\a+10,-5+\b);
\draw[very thick](\a+10,-5+\b)--(\a+12,-5+\b);
\draw[very thick](\a+12,-5+\b)--(\a+14,-5+\b);
\draw[very thick](\a+14,-5+\b)--(\a+16,-5+\b);
\draw[very thick](\a+3,-7.50000000000000+\b)--(\a+3,-8+\b);
\draw[very thick](\a+3.50000000000000,-7+\b)--(\a+4,-7+\b);
\draw[very thick](\a+3.50000000000000,-7+\b)arc(90:180:0.500000000000000);
\draw[very thick](\a+4,-7+\b)--(\a+6,-7+\b);
\draw[very thick](\a+5,-6+\b)--(\a+5,-8+\b);
\draw[very thick](\a+7,-6+\b)--(\a+7,-6.50000000000000+\b);
\draw[very thick](\a+6,-7+\b)--(\a+6.50000000000000,-7+\b);
\draw[very thick](\a+6.50000000000000,-7+\b)arc(270:360:0.500000000000000);
\draw[very thick](\a+13,-7.50000000000000+\b)--(\a+13,-8+\b);
\draw[very thick](\a+13.5000000000000,-7+\b)--(\a+14,-7+\b);
\draw[very thick](\a+13.5000000000000,-7+\b)arc(90:180:0.500000000000000);
\draw[very thick](\a+14,-7+\b)--(\a+16,-7+\b);
\draw[very thick](\a+3,-8+\b)--(\a+3,-10+\b);
\draw[very thick](\a+5,-8+\b)--(\a+5,-10+\b);
\draw[very thick](\a+11,-9.50000000000000+\b)--(\a+11,-10+\b);
\draw[very thick](\a+11.5000000000000,-9+\b)--(\a+12,-9+\b);
\draw[very thick](\a+11.5000000000000,-9+\b)arc(90:180:0.500000000000000);
\draw[very thick](\a+13,-8+\b)--(\a+13,-8.50000000000000+\b);
\draw[very thick](\a+12,-9+\b)--(\a+12.5000000000000,-9+\b);
\draw[very thick](\a+12.5000000000000,-9+\b)arc(270:360:0.500000000000000);
\draw[very thick](\a+15,-9.50000000000000+\b)--(\a+15,-10+\b);
\draw[very thick](\a+15.5000000000000,-9+\b)--(\a+16,-9+\b);
\draw[very thick](\a+15.5000000000000,-9+\b)arc(90:180:0.500000000000000);
\draw[very thick](\a+3,-10+\b)--(\a+3,-12+\b);
\draw[very thick](\a+5,-10+\b)--(\a+5,-12+\b);
\draw[very thick](\a+7,-11.5000000000000+\b)--(\a+7,-12+\b);
\draw[very thick](\a+7.50000000000000,-11+\b)--(\a+8,-11+\b);
\draw[very thick](\a+7.50000000000000,-11+\b)arc(90:180:0.500000000000000);
\draw[very thick](\a+8,-11+\b)--(\a+10,-11+\b);
\draw[very thick](\a+10,-11+\b)--(\a+12,-11+\b);
\draw[very thick](\a+11,-10+\b)--(\a+11,-12+\b);
\draw[very thick](\a+12,-11+\b)--(\a+14,-11+\b);
\draw[very thick](\a+14,-11+\b)--(\a+16,-11+\b);
\draw[very thick](\a+15,-10+\b)--(\a+15,-12+\b);
\draw[very thick](\a+3,-12+\b)--(\a+3,-14+\b);
\draw[very thick](\a+5,-12+\b)--(\a+5,-14+\b);
\draw[very thick](\a+7,-12+\b)--(\a+7,-14+\b);
\draw[very thick](\a+9,-13.5000000000000+\b)--(\a+9,-14+\b);
\draw[very thick](\a+9.50000000000000,-13+\b)--(\a+10,-13+\b);
\draw[very thick](\a+9.50000000000000,-13+\b)arc(90:180:0.500000000000000);
\draw[very thick](\a+10,-13+\b)--(\a+12,-13+\b);
\draw[very thick](\a+11,-12+\b)--(\a+11,-14+\b);
\draw[very thick](\a+12,-13+\b)--(\a+14,-13+\b);
\draw[very thick](\a+14,-13+\b)--(\a+16,-13+\b);
\draw[very thick](\a+15,-12+\b)--(\a+15,-14+\b);
\draw[very thick](\a+3,-14+\b)--(\a+3,-16+\b);
\draw[very thick](\a+5,-14+\b)--(\a+5,-16+\b);
\draw[very thick](\a+7,-14+\b)--(\a+7,-16+\b);
\draw[very thick](\a+9,-14+\b)--(\a+9,-16+\b);
\draw[very thick](\a+11,-14+\b)--(\a+11,-16+\b);
\draw[very thick](\a+13,-15.5000000000000+\b)--(\a+13,-16+\b);
\draw[very thick](\a+13.5000000000000,-15+\b)--(\a+14,-15+\b);
\draw[very thick](\a+13.5000000000000,-15+\b)arc(90:180:0.500000000000000);
\draw[very thick](\a+14,-15+\b)--(\a+16,-15+\b);
\draw[very thick](\a+15,-14+\b)--(\a+15,-16+\b);

\def\a{56};
\def\b{0};
\draw[step=2, color=green](\a+2,-16+\b)grid(\a+16,-2+\b);
\draw[very thick](\a+5,-3.50000000000000+\b)--(\a+5,-4+\b);
\draw[very thick](\a+5.50000000000000,-3+\b)--(\a+6,-3+\b);
\draw[very thick](\a+5.50000000000000,-3+\b)arc(90:180:0.500000000000000);
\draw[very thick](\a+6,-3+\b)--(\a+8,-3+\b);
\draw[very thick](\a+8,-3+\b)--(\a+10,-3+\b);
\draw[very thick](\a+10,-3+\b)--(\a+12,-3+\b);
\draw[very thick](\a+12,-3+\b)--(\a+14,-3+\b);
\draw[very thick](\a+14,-3+\b)--(\a+16,-3+\b);
\draw[very thick](\a+5,-4+\b)--(\a+5,-6+\b);
\draw[very thick](\a+7,-5.50000000000000+\b)--(\a+7,-6+\b);
\draw[very thick](\a+7.50000000000000,-5+\b)--(\a+8,-5+\b);
\draw[very thick](\a+7.50000000000000,-5+\b)arc(90:180:0.500000000000000);
\draw[very thick](\a+8,-5+\b)--(\a+10,-5+\b);
\draw[very thick](\a+10,-5+\b)--(\a+12,-5+\b);
\draw[very thick](\a+12,-5+\b)--(\a+14,-5+\b);
\draw[very thick](\a+14,-5+\b)--(\a+16,-5+\b);
\draw[very thick](\a+3,-7.50000000000000+\b)--(\a+3,-8+\b);
\draw[very thick](\a+3.50000000000000,-7+\b)--(\a+4,-7+\b);
\draw[very thick](\a+3.50000000000000,-7+\b)arc(90:180:0.500000000000000);
\draw[very thick](\a+4,-7+\b)--(\a+6,-7+\b);
\draw[very thick](\a+5,-6+\b)--(\a+5,-8+\b);
\draw[very thick](\a+7,-6+\b)--(\a+7,-6.50000000000000+\b);
\draw[very thick](\a+6,-7+\b)--(\a+6.50000000000000,-7+\b);
\draw[very thick](\a+6.50000000000000,-7+\b)arc(270:360:0.500000000000000);
\draw[very thick](\a+13,-7.50000000000000+\b)--(\a+13,-8+\b);
\draw[very thick](\a+13.5000000000000,-7+\b)--(\a+14,-7+\b);
\draw[very thick](\a+13.5000000000000,-7+\b)arc(90:180:0.500000000000000);
\draw[very thick](\a+14,-7+\b)--(\a+16,-7+\b);
\draw[very thick](\a+3,-8+\b)--(\a+3,-10+\b);
\draw[very thick](\a+5,-8+\b)--(\a+5,-10+\b);
\draw[very thick](\a+9,-9.50000000000000+\b)--(\a+9,-10+\b);
\draw[very thick](\a+9.50000000000000,-9+\b)--(\a+10,-9+\b);
\draw[very thick](\a+9.50000000000000,-9+\b)arc(90:180:0.500000000000000);
\draw[very thick](\a+10,-9+\b)--(\a+12,-9+\b);
\draw[very thick](\a+13,-8+\b)--(\a+13,-8.50000000000000+\b);
\draw[very thick](\a+12,-9+\b)--(\a+12.5000000000000,-9+\b);
\draw[very thick](\a+12.5000000000000,-9+\b)arc(270:360:0.500000000000000);
\draw[very thick](\a+15,-9.50000000000000+\b)--(\a+15,-10+\b);
\draw[very thick](\a+15.5000000000000,-9+\b)--(\a+16,-9+\b);
\draw[very thick](\a+15.5000000000000,-9+\b)arc(90:180:0.500000000000000);
\draw[very thick](\a+3,-10+\b)--(\a+3,-12+\b);
\draw[very thick](\a+5,-10+\b)--(\a+5,-12+\b);
\draw[very thick](\a+7,-11.5000000000000+\b)--(\a+7,-12+\b);
\draw[very thick](\a+7.50000000000000,-11+\b)--(\a+8,-11+\b);
\draw[very thick](\a+7.50000000000000,-11+\b)arc(90:180:0.500000000000000);
\draw[very thick](\a+8,-11+\b)--(\a+10,-11+\b);
\draw[very thick](\a+9,-10+\b)--(\a+9,-12+\b);
\draw[very thick](\a+10,-11+\b)--(\a+12,-11+\b);
\draw[very thick](\a+12,-11+\b)--(\a+14,-11+\b);
\draw[very thick](\a+14,-11+\b)--(\a+16,-11+\b);
\draw[very thick](\a+15,-10+\b)--(\a+15,-12+\b);
\draw[very thick](\a+3,-12+\b)--(\a+3,-14+\b);
\draw[very thick](\a+5,-12+\b)--(\a+5,-14+\b);
\draw[very thick](\a+7,-12+\b)--(\a+7,-14+\b);
\draw[very thick](\a+9,-12+\b)--(\a+9,-14+\b);
\draw[very thick](\a+11,-13.5000000000000+\b)--(\a+11,-14+\b);
\draw[very thick](\a+11.5000000000000,-13+\b)--(\a+12,-13+\b);
\draw[very thick](\a+11.5000000000000,-13+\b)arc(90:180:0.500000000000000);
\draw[very thick](\a+12,-13+\b)--(\a+14,-13+\b);
\draw[very thick](\a+14,-13+\b)--(\a+16,-13+\b);
\draw[very thick](\a+15,-12+\b)--(\a+15,-14+\b);
\draw[very thick](\a+3,-14+\b)--(\a+3,-16+\b);
\draw[very thick](\a+5,-14+\b)--(\a+5,-16+\b);
\draw[very thick](\a+7,-14+\b)--(\a+7,-16+\b);
\draw[very thick](\a+9,-14+\b)--(\a+9,-16+\b);
\draw[very thick](\a+11,-14+\b)--(\a+11,-16+\b);
\draw[very thick](\a+13,-15.5000000000000+\b)--(\a+13,-16+\b);
\draw[very thick](\a+13.5000000000000,-15+\b)--(\a+14,-15+\b);
\draw[very thick](\a+13.5000000000000,-15+\b)arc(90:180:0.500000000000000);
\draw[very thick](\a+14,-15+\b)--(\a+16,-15+\b);
\draw[very thick](\a+15,-14+\b)--(\a+15,-16+\b);
\node[right] at (\a+16,-3+\b) {$2$};
\node[right] at (\a+16,-5+\b) {$1$};
\node[right] at (\a+16,-7+\b) {$4$};
\node[right] at (\a+16,-9+\b) {$7$};
\node[right] at (\a+16,-11+\b) {$3$};
\node[right] at (\a+16,-13+\b) {$5$};
\node[right] at (\a+16,-15+\b) {$6$};
\node[below] at (\a+3,-16+\b) {$1$};
\node[below] at (\a+5,-16+\b) {$2$};
\node[below] at (\a+7,-16+\b) {$3$};
\node[below] at (\a+9,-16+\b) {$4$};
\node[below] at (\a+11,-16+\b) {$5$};
\node[below] at (\a+13,-16+\b) {$6$};
\node[below] at (\a+15,-16+\b) {$7$};

\draw[very thick,->] (19,-9)--(21,-9);
\node[above] at (20,-9) {(2)};
\draw[very thick,->] (37,-9)--(39,-9);
\node[above] at (38,-9) {(3)};
\draw[very thick,->] (55,-9)--(57,-9);
\node[above] at (56,-9) {(4)};
\node at (5,-3) {$\times$};
\node at (27,-9) {$\times$};
\node at (47,-9) {$\times$};
\end{tikzpicture}
\caption{Steps for obtaining $\nabla D$ from $D$ by the algorithm in \Cref{def:pop}.}
\label{fig:pop-example}
\end{figure}

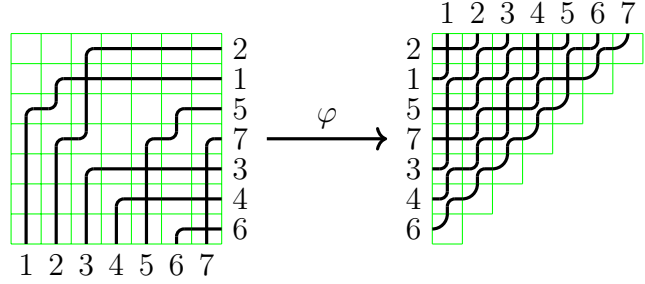
\begin{figure}[h!]
\centering
\begin{tikzpicture}[scale=0.200000000000000]
\draw[color=green](2,-2)--(16,-2);
\draw[color=green](2,-2)--(2,-16);
\draw[color=green](2,-4)--(16,-4);
\draw[color=green](4,-2)--(4,-16);
\draw[color=green](2,-6)--(14,-6);
\draw[color=green](6,-2)--(6,-14);
\draw[color=green](2,-8)--(12,-8);
\draw[color=green](8,-2)--(8,-12);
\draw[color=green](2,-10)--(10,-10);
\draw[color=green](10,-2)--(10,-10);
\draw[color=green](2,-12)--(8,-12);
\draw[color=green](12,-2)--(12,-8);
\draw[color=green](2,-14)--(6,-14);
\draw[color=green](14,-2)--(14,-6);
\draw[color=green](2,-16)--(4,-16);
\draw[color=green](16,-2)--(16,-4);
\draw[very thick](2,-3)--(4,-3);
\draw[very thick](3,-2)--(3,-4);
\draw[very thick](5,-2)--(5,-2.50000000000000);
\draw[very thick](5,-3.50000000000000)--(5,-4);
\draw[very thick](4,-3)--(4.50000000000000,-3);
\draw[very thick](5.50000000000000,-3)--(6,-3);
\draw[very thick](4.50000000000000,-3)arc(270:360:0.500000000000000);
\draw[very thick](5.50000000000000,-3)arc(90:180:0.500000000000000);
\draw[very thick](7,-2)--(7,-2.50000000000000);
\draw[very thick](7,-3.50000000000000)--(7,-4);
\draw[very thick](6,-3)--(6.50000000000000,-3);
\draw[very thick](7.50000000000000,-3)--(8,-3);
\draw[very thick](6.50000000000000,-3)arc(270:360:0.500000000000000);
\draw[very thick](7.50000000000000,-3)arc(90:180:0.500000000000000);
\draw[very thick](8,-3)--(10,-3);
\draw[very thick](9,-2)--(9,-4);
\draw[very thick](11,-2)--(11,-2.50000000000000);
\draw[very thick](11,-3.50000000000000)--(11,-4);
\draw[very thick](10,-3)--(10.5000000000000,-3);
\draw[very thick](11.5000000000000,-3)--(12,-3);
\draw[very thick](10.5000000000000,-3)arc(270:360:0.500000000000000);
\draw[very thick](11.5000000000000,-3)arc(90:180:0.500000000000000);
\draw[very thick](13,-2)--(13,-2.50000000000000);
\draw[very thick](13,-3.50000000000000)--(13,-4);
\draw[very thick](12,-3)--(12.5000000000000,-3);
\draw[very thick](13.5000000000000,-3)--(14,-3);
\draw[very thick](12.5000000000000,-3)arc(270:360:0.500000000000000);
\draw[very thick](13.5000000000000,-3)arc(90:180:0.500000000000000);
\draw[very thick](3,-4)--(3,-4.50000000000000);
\draw[very thick](3,-5.50000000000000)--(3,-6);
\draw[very thick](2,-5)--(2.50000000000000,-5);
\draw[very thick](3.50000000000000,-5)--(4,-5);
\draw[very thick](2.50000000000000,-5)arc(270:360:0.500000000000000);
\draw[very thick](3.50000000000000,-5)arc(90:180:0.500000000000000);
\draw[very thick](5,-4)--(5,-4.50000000000000);
\draw[very thick](5,-5.50000000000000)--(5,-6);
\draw[very thick](4,-5)--(4.50000000000000,-5);
\draw[very thick](5.50000000000000,-5)--(6,-5);
\draw[very thick](4.50000000000000,-5)arc(270:360:0.500000000000000);
\draw[very thick](5.50000000000000,-5)arc(90:180:0.500000000000000);
\draw[very thick](7,-4)--(7,-4.50000000000000);
\draw[very thick](7,-5.50000000000000)--(7,-6);
\draw[very thick](6,-5)--(6.50000000000000,-5);
\draw[very thick](7.50000000000000,-5)--(8,-5);
\draw[very thick](6.50000000000000,-5)arc(270:360:0.500000000000000);
\draw[very thick](7.50000000000000,-5)arc(90:180:0.500000000000000);
\draw[very thick](9,-4)--(9,-4.50000000000000);
\draw[very thick](9,-5.50000000000000)--(9,-6);
\draw[very thick](8,-5)--(8.50000000000000,-5);
\draw[very thick](9.50000000000000,-5)--(10,-5);
\draw[very thick](8.50000000000000,-5)arc(270:360:0.500000000000000);
\draw[very thick](9.50000000000000,-5)arc(90:180:0.500000000000000);
\draw[very thick](10,-5)--(12,-5);
\draw[very thick](11,-4)--(11,-6);
\draw[very thick](2,-7)--(4,-7);
\draw[very thick](3,-6)--(3,-8);
\draw[very thick](5,-6)--(5,-6.50000000000000);
\draw[very thick](5,-7.50000000000000)--(5,-8);
\draw[very thick](4,-7)--(4.50000000000000,-7);
\draw[very thick](5.50000000000000,-7)--(6,-7);
\draw[very thick](4.50000000000000,-7)arc(270:360:0.500000000000000);
\draw[very thick](5.50000000000000,-7)arc(90:180:0.500000000000000);
\draw[very thick](6,-7)--(8,-7);
\draw[very thick](7,-6)--(7,-8);
\draw[very thick](9,-6)--(9,-6.50000000000000);
\draw[very thick](9,-7.50000000000000)--(9,-8);
\draw[very thick](8,-7)--(8.50000000000000,-7);
\draw[very thick](9.50000000000000,-7)--(10,-7);
\draw[very thick](8.50000000000000,-7)arc(270:360:0.500000000000000);
\draw[very thick](9.50000000000000,-7)arc(90:180:0.500000000000000);
\draw[very thick](2,-9)--(4,-9);
\draw[very thick](3,-8)--(3,-10);
\draw[very thick](5,-8)--(5,-8.50000000000000);
\draw[very thick](5,-9.50000000000000)--(5,-10);
\draw[very thick](4,-9)--(4.50000000000000,-9);
\draw[very thick](5.50000000000000,-9)--(6,-9);
\draw[very thick](4.50000000000000,-9)arc(270:360:0.500000000000000);
\draw[very thick](5.50000000000000,-9)arc(90:180:0.500000000000000);
\draw[very thick](7,-8)--(7,-8.50000000000000);
\draw[very thick](7,-9.50000000000000)--(7,-10);
\draw[very thick](6,-9)--(6.50000000000000,-9);
\draw[very thick](7.50000000000000,-9)--(8,-9);
\draw[very thick](6.50000000000000,-9)arc(270:360:0.500000000000000);
\draw[very thick](7.50000000000000,-9)arc(90:180:0.500000000000000);
\draw[very thick](3,-10)--(3,-10.5000000000000);
\draw[very thick](3,-11.5000000000000)--(3,-12);
\draw[very thick](2,-11)--(2.50000000000000,-11);
\draw[very thick](3.50000000000000,-11)--(4,-11);
\draw[very thick](2.50000000000000,-11)arc(270:360:0.500000000000000);
\draw[very thick](3.50000000000000,-11)arc(90:180:0.500000000000000);
\draw[very thick](5,-10)--(5,-10.5000000000000);
\draw[very thick](5,-11.5000000000000)--(5,-12);
\draw[very thick](4,-11)--(4.50000000000000,-11);
\draw[very thick](5.50000000000000,-11)--(6,-11);
\draw[very thick](4.50000000000000,-11)arc(270:360:0.500000000000000);
\draw[very thick](5.50000000000000,-11)arc(90:180:0.500000000000000);
\draw[very thick](3,-12)--(3,-12.5000000000000);
\draw[very thick](3,-13.5000000000000)--(3,-14);
\draw[very thick](2,-13)--(2.50000000000000,-13);
\draw[very thick](3.50000000000000,-13)--(4,-13);
\draw[very thick](2.50000000000000,-13)arc(270:360:0.500000000000000);
\draw[very thick](3.50000000000000,-13)arc(90:180:0.500000000000000);
\draw[very thick](14,-3)arc(270:360:1);
\draw[very thick](12,-5)arc(270:360:1);
\draw[very thick](10,-7)arc(270:360:1);
\draw[very thick](8,-9)arc(270:360:1);
\draw[very thick](6,-11)arc(270:360:1);
\draw[very thick](4,-13)arc(270:360:1);
\draw[very thick](2,-15)arc(270:360:1);
\node[left] at (2,-3) {$2$};
\node[left] at (2,-5) {$1$};
\node[left] at (2,-7) {$5$};
\node[left] at (2,-9) {$7$};
\node[left] at (2,-11) {$3$};
\node[left] at (2,-13) {$4$};
\node[left] at (2,-15) {$6$};
\node[above] at (3,-2) {$1$};
\node[above] at (5,-2) {$2$};
\node[above] at (7,-2) {$3$};
\node[above] at (9,-2) {$4$};
\node[above] at (11,-2) {$5$};
\node[above] at (13,-2) {$6$};
\node[above] at (15,-2) {$7$};

\def\a{-28};
\def\b{0};
\draw[step=2, color=green](\a+2,-16+\b)grid(\a+16,-2+\b);
\draw[very thick](\a+7,-3.50000000000000+\b)--(\a+7,-4+\b);
\draw[very thick](\a+7.50000000000000,-3+\b)--(\a+8,-3+\b);
\draw[very thick](\a+7.50000000000000,-3+\b)arc(90:180:0.500000000000000);
\draw[very thick](\a+8,-3+\b)--(\a+10,-3+\b);
\draw[very thick](\a+10,-3+\b)--(\a+12,-3+\b);
\draw[very thick](\a+12,-3+\b)--(\a+14,-3+\b);
\draw[very thick](\a+14,-3+\b)--(\a+16,-3+\b);
\draw[very thick](\a+5,-5.50000000000000+\b)--(\a+5,-6+\b);
\draw[very thick](\a+5.50000000000000,-5+\b)--(\a+6,-5+\b);
\draw[very thick](\a+5.50000000000000,-5+\b)arc(90:180:0.500000000000000);
\draw[very thick](\a+6,-5+\b)--(\a+8,-5+\b);
\draw[very thick](\a+7,-4+\b)--(\a+7,-6+\b);
\draw[very thick](\a+8,-5+\b)--(\a+10,-5+\b);
\draw[very thick](\a+10,-5+\b)--(\a+12,-5+\b);
\draw[very thick](\a+12,-5+\b)--(\a+14,-5+\b);
\draw[very thick](\a+14,-5+\b)--(\a+16,-5+\b);
\draw[very thick](\a+3,-7.50000000000000+\b)--(\a+3,-8+\b);
\draw[very thick](\a+3.50000000000000,-7+\b)--(\a+4,-7+\b);
\draw[very thick](\a+3.50000000000000,-7+\b)arc(90:180:0.500000000000000);
\draw[very thick](\a+5,-6+\b)--(\a+5,-6.50000000000000+\b);
\draw[very thick](\a+4,-7+\b)--(\a+4.50000000000000,-7+\b);
\draw[very thick](\a+4.50000000000000,-7+\b)arc(270:360:0.500000000000000);
\draw[very thick](\a+7,-6+\b)--(\a+7,-8+\b);
\draw[very thick](\a+13,-7.50000000000000+\b)--(\a+13,-8+\b);
\draw[very thick](\a+13.5000000000000,-7+\b)--(\a+14,-7+\b);
\draw[very thick](\a+13.5000000000000,-7+\b)arc(90:180:0.500000000000000);
\draw[very thick](\a+14,-7+\b)--(\a+16,-7+\b);
\draw[very thick](\a+3,-8+\b)--(\a+3,-10+\b);
\draw[very thick](\a+5,-9.50000000000000+\b)--(\a+5,-10+\b);
\draw[very thick](\a+5.50000000000000,-9+\b)--(\a+6,-9+\b);
\draw[very thick](\a+5.50000000000000,-9+\b)arc(90:180:0.500000000000000);
\draw[very thick](\a+7,-8+\b)--(\a+7,-8.50000000000000+\b);
\draw[very thick](\a+6,-9+\b)--(\a+6.50000000000000,-9+\b);
\draw[very thick](\a+6.50000000000000,-9+\b)arc(270:360:0.500000000000000);
\draw[very thick](\a+11,-9.50000000000000+\b)--(\a+11,-10+\b);
\draw[very thick](\a+11.5000000000000,-9+\b)--(\a+12,-9+\b);
\draw[very thick](\a+11.5000000000000,-9+\b)arc(90:180:0.500000000000000);
\draw[very thick](\a+13,-8+\b)--(\a+13,-8.50000000000000+\b);
\draw[very thick](\a+12,-9+\b)--(\a+12.5000000000000,-9+\b);
\draw[very thick](\a+12.5000000000000,-9+\b)arc(270:360:0.500000000000000);
\draw[very thick](\a+15,-9.50000000000000+\b)--(\a+15,-10+\b);
\draw[very thick](\a+15.5000000000000,-9+\b)--(\a+16,-9+\b);
\draw[very thick](\a+15.5000000000000,-9+\b)arc(90:180:0.500000000000000);
\draw[very thick](\a+3,-10+\b)--(\a+3,-12+\b);
\draw[very thick](\a+5,-10+\b)--(\a+5,-12+\b);
\draw[very thick](\a+7,-11.5000000000000+\b)--(\a+7,-12+\b);
\draw[very thick](\a+7.50000000000000,-11+\b)--(\a+8,-11+\b);
\draw[very thick](\a+7.50000000000000,-11+\b)arc(90:180:0.500000000000000);
\draw[very thick](\a+8,-11+\b)--(\a+10,-11+\b);
\draw[very thick](\a+10,-11+\b)--(\a+12,-11+\b);
\draw[very thick](\a+11,-10+\b)--(\a+11,-12+\b);
\draw[very thick](\a+12,-11+\b)--(\a+14,-11+\b);
\draw[very thick](\a+14,-11+\b)--(\a+16,-11+\b);
\draw[very thick](\a+15,-10+\b)--(\a+15,-12+\b);
\draw[very thick](\a+3,-12+\b)--(\a+3,-14+\b);
\draw[very thick](\a+5,-12+\b)--(\a+5,-14+\b);
\draw[very thick](\a+7,-12+\b)--(\a+7,-14+\b);
\draw[very thick](\a+9,-13.5000000000000+\b)--(\a+9,-14+\b);
\draw[very thick](\a+9.50000000000000,-13+\b)--(\a+10,-13+\b);
\draw[very thick](\a+9.50000000000000,-13+\b)arc(90:180:0.500000000000000);
\draw[very thick](\a+10,-13+\b)--(\a+12,-13+\b);
\draw[very thick](\a+11,-12+\b)--(\a+11,-14+\b);
\draw[very thick](\a+12,-13+\b)--(\a+14,-13+\b);
\draw[very thick](\a+14,-13+\b)--(\a+16,-13+\b);
\draw[very thick](\a+15,-12+\b)--(\a+15,-14+\b);
\draw[very thick](\a+3,-14+\b)--(\a+3,-16+\b);
\draw[very thick](\a+5,-14+\b)--(\a+5,-16+\b);
\draw[very thick](\a+7,-14+\b)--(\a+7,-16+\b);
\draw[very thick](\a+9,-14+\b)--(\a+9,-16+\b);
\draw[very thick](\a+11,-14+\b)--(\a+11,-16+\b);
\draw[very thick](\a+13,-15.5000000000000+\b)--(\a+13,-16+\b);
\draw[very thick](\a+13.5000000000000,-15+\b)--(\a+14,-15+\b);
\draw[very thick](\a+13.5000000000000,-15+\b)arc(90:180:0.500000000000000);
\draw[very thick](\a+14,-15+\b)--(\a+16,-15+\b);
\draw[very thick](\a+15,-14+\b)--(\a+15,-16+\b);
\node[right] at (\a+16,-3+\b) {$2$};
\node[right] at (\a+16,-5+\b) {$1$};
\node[right] at (\a+16,-7+\b) {$5$};
\node[right] at (\a+16,-9+\b) {$7$};
\node[right] at (\a+16,-11+\b) {$3$};
\node[right] at (\a+16,-13+\b) {$4$};
\node[right] at (\a+16,-15+\b) {$6$};
\node[below] at (\a+3,-16+\b) {$1$};
\node[below] at (\a+5,-16+\b) {$2$};
\node[below] at (\a+7,-16+\b) {$3$};
\node[below] at (\a+9,-16+\b) {$4$};
\node[below] at (\a+11,-16+\b) {$5$};
\node[below] at (\a+13,-16+\b) {$6$};
\node[below] at (\a+15,-16+\b) {$7$};

\draw[very thick,->] (-9,-9)--(-1,-9);
\node[above] at (-5,-9) {$\varphi$};
\end{tikzpicture}
\caption{The bijection between BPDs and PDs}
\label{fig:bijection-example}
\end{figure}
\end{Example}
\begin{Theorem}\cite{GaoHuang-bijection}\label{thm:bijection-PD-BPD}
The map $\varphi$ in \Cref{def:bijection-PD-BPD} is a weight-preserving bijection from reduced bumpless pipe dreams to reduced compatible pairs of a fixed permutation.
\end{Theorem}

To prove the theorem given $D\in\BPD(w)$, it is
straightforward to check $\varphi(D)$ has the same monomial weight.
To prove $\varphi(D)$ is a reduced compatible sequence for the same
permutation is not so hard.  Given that \Cref{thm:bumpless-schubert}
has already been established, it then
suffices to show $\varphi$ is an injection. However, the ``canonical
nature'' of $\varphi$ as the right bijection is significantly harder
and richer to prove.

As shown in \Cref{sub:Games}, insertion rules on pipe dreams can be
used to bijectively prove the Transition Formula \eqref{t:transitionA}
and Monk's rule \eqref{eq:Monk-v1}.  Consider the following related
product formula for one variable times a Schubert polynomial, where
the terms shown are all positive and have coefficient 1.

\begin{Theorem}[\textbf{Variation on Monk's rule}]\label{thm:Monk}
For $w\in S_{\infty}$ and $\alpha\in\mathbb{Z}_{>0}$, 
\begin{equation}\label{eq:Monk-v2}
x_{\alpha}\fS_{w}+\sum_{\substack{k<\alpha\\w t_{k,\alpha}\gtrdot w}}\fS_{w t_{k,\alpha}}=\sum_{\substack{l>\alpha\\ w t_{\alpha,l}\gtrdot w}}\fS_{w t_{\alpha,l}}.
\end{equation}
\end{Theorem}
It is not hard to see that \Cref{eq:Monk-v2} is equivalent to Monk's rule (\Cref{eq:Monk-v1}) after a step of subtraction. 
We see that \Cref{eq:Monk-v2} is the transition equation (\Cref{eq:Monk.transition}) if
the right hand side has a single term. \Cref{eq:Monk-v2} is called a
\emph{cotransition equation} whenever the left hand side has a single term
$x_{\alpha}\fS_{w}$.

 \Cref{thm:Monk} can also be proved bijectively using bumpless pipe
dreams.  To construct the weight-preserving bijection from the objects
enumerated by the left hand side to that of the right hand side in
\eqref{eq:Monk-v2}, we will define two maps
$x_\alpha{\rightsquigarrow}$ and $m_{k,\alpha}$,
\[\begin{cases}
m_{k,\alpha}:&\BPD(w t_{k,\alpha})\rightarrow
\bigcup_{\substack{l>\alpha\\w t_{\alpha,l}\gtrdot w}}\BPD(w
t_{\alpha,l}),\ k<\alpha,w t_{k,\alpha}\gtrdot w\\
x_\alpha{\rightsquigarrow}:&\BPD(w)\rightarrow \bigcup_{\substack{l>\alpha\\w t_{\alpha,l}\gtrdot w}}\BPD(w t_{\alpha,l})
\end{cases}\] that were first described by Huang
\cite{huang2020bijective}.  The $m_{k,\alpha}$ map is a way of
``moving'' a crossing. The $x_{\alpha}\rightsquigarrow D$ map is defined via an
insertion algorithm.

\begin{Definition}\cite{huang2020bijective}\label{def:map}
Under the assumption that $k<\alpha$ and $wt_{k,\alpha}\gtrdot w$, for $D\in\BPD(w t_{k,\alpha})$, the bumpless pipe dream $m_{k,\alpha}(D)$ is constructed as follows:
\begin{enumerate}
\item Locate the \+-tile at coordinate $(x,y)$ between pipe $w^{-1}(k)$ and $w^{-1}(\alpha)$ and turn it temporarily into \elbow.
\item Do a min-droop at $(x,y)$ to $(x',y')$, allowing a temporary
\elbow\ if $(x',y')$ contained a \rt. 
\begin{enumerate}
\item[(a)] If $(x',y')$ was \bl, then it is now \jt, update $(x,y)$ to be
the coordinate of the unique \rt-tile in row $x'$ from the same pipe passing through $(x',y')$
and repeat step (2).

\item[(b)] If $(x',y')$ was \rt, then it is now \elbow, and we turn it into \+. If
these two pipes cross somewhere else, the position must be unique, say at $(x'',y'')$. Replace that tile with
\elbow, update $(x,y)=(x'',y'')$, and repeat step (2).  Otherwise
stop, and return the updated tiling.
\end{enumerate}
\end{enumerate}
\end{Definition}

\begin{Definition}
For $D\in\BPD(w)$, $x_\alpha{\rightsquigarrow}D$ is
constructed by first locating the rightmost \rt-tile in row
$\alpha$, say in position $(x,y)$. Do a min-droop at $(x,y)$ as in step (2) of \Cref{def:map} and continue as in this step until the algorithm outputs
a tiling.
\end{Definition}

\begin{Theorem}\cite{huang2020bijective}\label{thm:huang-bump}
Both $x_\alpha{\rightsquigarrow}$ and $m_{k,\alpha}$ map bumpless pipe
dreams to bumpless pipe dreams. Furthermore, together they define a
weight-preserving bijection between the objects on the left and right
side of \eqref{eq:Monk-v2}.
\end{Theorem}

\begin{figure}[h!]
\centering
\begin{tikzpicture}[scale=0.200000000000000]
\draw[fill=gray!30] (4,-2)--(6,-2)--(6,-4)--(4,-4)--(4,-2);
\draw[fill=gray!30] (24,-4)--(26,-4)--(26,-6)--(24,-6)--(24,-4);
\draw[fill=gray!30] (38,-6)--(40,-6)--(40,-8)--(38,-8)--(38,-6);
\draw[fill=gray!30] (56,-8)--(58,-8)--(58,-10)--(56,-10)--(56,-8);
\def\a{0};
\def\b{0};
\draw[step=2, color=green](\a+2,-12+\b)grid(\a+12,-2+\b);
\draw[very thick](\a+5,-3.50000000000000+\b)--(\a+5,-4+\b);
\draw[very thick](\a+5.50000000000000,-3+\b)--(\a+6,-3+\b);
\draw[very thick](\a+5.50000000000000,-3+\b)arc(90:180:0.500000000000000);
\draw[very thick](\a+6,-3+\b)--(\a+8,-3+\b);
\draw[very thick](\a+8,-3+\b)--(\a+10,-3+\b);
\draw[very thick](\a+10,-3+\b)--(\a+12,-3+\b);
\draw[very thick](\a+5,-4+\b)--(\a+5,-6+\b);
\draw[very thick](\a+7,-5.50000000000000+\b)--(\a+7,-6+\b);
\draw[very thick](\a+7.50000000000000,-5+\b)--(\a+8,-5+\b);
\draw[very thick](\a+7.50000000000000,-5+\b)arc(90:180:0.500000000000000);
\draw[very thick](\a+8,-5+\b)--(\a+10,-5+\b);
\draw[very thick](\a+10,-5+\b)--(\a+12,-5+\b);
\draw[very thick](\a+3,-7.50000000000000+\b)--(\a+3,-8+\b);
\draw[very thick](\a+3.50000000000000,-7+\b)--(\a+4,-7+\b);
\draw[very thick](\a+3.50000000000000,-7+\b)arc(90:180:0.500000000000000);
\draw[very thick](\a+4,-7+\b)--(\a+6,-7+\b);
\draw[very thick](\a+5,-6+\b)--(\a+5,-8+\b);
\draw[very thick](\a+7,-6+\b)--(\a+7,-6.50000000000000+\b);
\draw[very thick](\a+6,-7+\b)--(\a+6.50000000000000,-7+\b);
\draw[very thick](\a+6.50000000000000,-7+\b)arc(270:360:0.500000000000000);
\draw[very thick](\a+11,-7.50000000000000+\b)--(\a+11,-8+\b);
\draw[very thick](\a+11.5000000000000,-7+\b)--(\a+12,-7+\b);
\draw[very thick](\a+11.5000000000000,-7+\b)arc(90:180:0.500000000000000);
\draw[very thick](\a+3,-8+\b)--(\a+3,-10+\b);
\draw[very thick](\a+5,-8+\b)--(\a+5,-10+\b);
\draw[very thick](\a+7,-9.50000000000000+\b)--(\a+7,-10+\b);
\draw[very thick](\a+7.50000000000000,-9+\b)--(\a+8,-9+\b);
\draw[very thick](\a+7.50000000000000,-9+\b)arc(90:180:0.500000000000000);
\draw[very thick](\a+8,-9+\b)--(\a+10,-9+\b);
\draw[very thick](\a+10,-9+\b)--(\a+12,-9+\b);
\draw[very thick](\a+11,-8+\b)--(\a+11,-10+\b);
\draw[very thick](\a+3,-10+\b)--(\a+3,-12+\b);
\draw[very thick](\a+5,-10+\b)--(\a+5,-12+\b);
\draw[very thick](\a+7,-10+\b)--(\a+7,-12+\b);
\draw[very thick](\a+9,-11.5000000000000+\b)--(\a+9,-12+\b);
\draw[very thick](\a+9.50000000000000,-11+\b)--(\a+10,-11+\b);
\draw[very thick](\a+9.50000000000000,-11+\b)arc(90:180:0.500000000000000);
\draw[very thick](\a+10,-11+\b)--(\a+12,-11+\b);
\draw[very thick](\a+11,-10+\b)--(\a+11,-12+\b);
\node[right] at (\a+12,-3+\b) {$2$};
\node[right] at (\a+12,-5+\b) {$1$};
\node[right] at (\a+12,-7+\b) {$5$};
\node[right] at (\a+12,-9+\b) {$3$};
\node[right] at (\a+12,-11+\b) {$4$};
\node[below] at (\a+3,-12+\b) {$1$};
\node[below] at (\a+5,-12+\b) {$2$};
\node[below] at (\a+7,-12+\b) {$3$};
\node[below] at (\a+9,-12+\b) {$4$};
\node[below] at (\a+11,-12+\b) {$5$};

\def\a{18};
\def\b{0};
\draw[step=2, color=green](\a+2,-12+\b)grid(\a+12,-2+\b);
\draw[very thick](\a+7,-3.50000000000000+\b)--(\a+7,-4+\b);
\draw[very thick](\a+7.50000000000000,-3+\b)--(\a+8,-3+\b);
\draw[very thick](\a+7.50000000000000,-3+\b)arc(90:180:0.500000000000000);
\draw[very thick](\a+8,-3+\b)--(\a+10,-3+\b);
\draw[very thick](\a+10,-3+\b)--(\a+12,-3+\b);
\draw[very thick](\a+5,-5.50000000000000+\b)--(\a+5,-6+\b);
\draw[very thick](\a+5.50000000000000,-5+\b)--(\a+6,-5+\b);
\draw[very thick](\a+5.50000000000000,-5+\b)arc(90:180:0.500000000000000);
\draw[very thick](\a+7,-5.50000000000000+\b)--(\a+7,-6+\b);
\draw[very thick](\a+7.50000000000000,-5+\b)--(\a+8,-5+\b);
\draw[very thick](\a+7.50000000000000,-5+\b)arc(90:180:0.500000000000000);
\draw[very thick](\a+7,-4+\b)--(\a+7,-4.50000000000000+\b);
\draw[very thick](\a+6,-5+\b)--(\a+6.50000000000000,-5+\b);
\draw[very thick](\a+6.50000000000000,-5+\b)arc(270:360:0.500000000000000);
\draw[very thick](\a+8,-5+\b)--(\a+10,-5+\b);
\draw[very thick](\a+10,-5+\b)--(\a+12,-5+\b);
\draw[very thick](\a+3,-7.50000000000000+\b)--(\a+3,-8+\b);
\draw[very thick](\a+3.50000000000000,-7+\b)--(\a+4,-7+\b);
\draw[very thick](\a+3.50000000000000,-7+\b)arc(90:180:0.500000000000000);
\draw[very thick](\a+4,-7+\b)--(\a+6,-7+\b);
\draw[very thick](\a+5,-6+\b)--(\a+5,-8+\b);
\draw[very thick](\a+7,-6+\b)--(\a+7,-6.50000000000000+\b);
\draw[very thick](\a+6,-7+\b)--(\a+6.50000000000000,-7+\b);
\draw[very thick](\a+6.50000000000000,-7+\b)arc(270:360:0.500000000000000);
\draw[very thick](\a+11,-7.50000000000000+\b)--(\a+11,-8+\b);
\draw[very thick](\a+11.5000000000000,-7+\b)--(\a+12,-7+\b);
\draw[very thick](\a+11.5000000000000,-7+\b)arc(90:180:0.500000000000000);
\draw[very thick](\a+3,-8+\b)--(\a+3,-10+\b);
\draw[very thick](\a+5,-8+\b)--(\a+5,-10+\b);
\draw[very thick](\a+7,-9.50000000000000+\b)--(\a+7,-10+\b);
\draw[very thick](\a+7.50000000000000,-9+\b)--(\a+8,-9+\b);
\draw[very thick](\a+7.50000000000000,-9+\b)arc(90:180:0.500000000000000);
\draw[very thick](\a+8,-9+\b)--(\a+10,-9+\b);
\draw[very thick](\a+10,-9+\b)--(\a+12,-9+\b);
\draw[very thick](\a+11,-8+\b)--(\a+11,-10+\b);
\draw[very thick](\a+3,-10+\b)--(\a+3,-12+\b);
\draw[very thick](\a+5,-10+\b)--(\a+5,-12+\b);
\draw[very thick](\a+7,-10+\b)--(\a+7,-12+\b);
\draw[very thick](\a+9,-11.5000000000000+\b)--(\a+9,-12+\b);
\draw[very thick](\a+9.50000000000000,-11+\b)--(\a+10,-11+\b);
\draw[very thick](\a+9.50000000000000,-11+\b)arc(90:180:0.500000000000000);
\draw[very thick](\a+10,-11+\b)--(\a+12,-11+\b);
\draw[very thick](\a+11,-10+\b)--(\a+11,-12+\b);

\def\a{34};
\def\b{0};
\draw[step=2, color=green](\a+2,-12+\b)grid(\a+12,-2+\b);
\draw[very thick](\a+7,-3.50000000000000+\b)--(\a+7,-4+\b);
\draw[very thick](\a+7.50000000000000,-3+\b)--(\a+8,-3+\b);
\draw[very thick](\a+7.50000000000000,-3+\b)arc(90:180:0.500000000000000);
\draw[very thick](\a+8,-3+\b)--(\a+10,-3+\b);
\draw[very thick](\a+10,-3+\b)--(\a+12,-3+\b);
\draw[very thick](\a+5,-5.50000000000000+\b)--(\a+5,-6+\b);
\draw[very thick](\a+5.50000000000000,-5+\b)--(\a+6,-5+\b);
\draw[very thick](\a+5.50000000000000,-5+\b)arc(90:180:0.500000000000000);
\draw[very thick](\a+6,-5+\b)--(\a+8,-5+\b);
\draw[very thick](\a+7,-4+\b)--(\a+7,-6+\b);
\draw[very thick](\a+8,-5+\b)--(\a+10,-5+\b);
\draw[very thick](\a+10,-5+\b)--(\a+12,-5+\b);
\draw[very thick](\a+3,-7.50000000000000+\b)--(\a+3,-8+\b);
\draw[very thick](\a+3.50000000000000,-7+\b)--(\a+4,-7+\b);
\draw[very thick](\a+3.50000000000000,-7+\b)arc(90:180:0.500000000000000);
\draw[very thick](\a+5,-7.50000000000000+\b)--(\a+5,-8+\b);
\draw[very thick](\a+5.50000000000000,-7+\b)--(\a+6,-7+\b);
\draw[very thick](\a+5.50000000000000,-7+\b)arc(90:180:0.500000000000000);
\draw[very thick](\a+5,-6+\b)--(\a+5,-6.50000000000000+\b);
\draw[very thick](\a+4,-7+\b)--(\a+4.50000000000000,-7+\b);
\draw[very thick](\a+4.50000000000000,-7+\b)arc(270:360:0.500000000000000);
\draw[very thick](\a+7,-6+\b)--(\a+7,-6.50000000000000+\b);
\draw[very thick](\a+6,-7+\b)--(\a+6.50000000000000,-7+\b);
\draw[very thick](\a+6.50000000000000,-7+\b)arc(270:360:0.500000000000000);
\draw[very thick](\a+11,-7.50000000000000+\b)--(\a+11,-8+\b);
\draw[very thick](\a+11.5000000000000,-7+\b)--(\a+12,-7+\b);
\draw[very thick](\a+11.5000000000000,-7+\b)arc(90:180:0.500000000000000);
\draw[very thick](\a+3,-8+\b)--(\a+3,-10+\b);
\draw[very thick](\a+5,-8+\b)--(\a+5,-10+\b);
\draw[very thick](\a+7,-9.50000000000000+\b)--(\a+7,-10+\b);
\draw[very thick](\a+7.50000000000000,-9+\b)--(\a+8,-9+\b);
\draw[very thick](\a+7.50000000000000,-9+\b)arc(90:180:0.500000000000000);
\draw[very thick](\a+8,-9+\b)--(\a+10,-9+\b);
\draw[very thick](\a+10,-9+\b)--(\a+12,-9+\b);
\draw[very thick](\a+11,-8+\b)--(\a+11,-10+\b);
\draw[very thick](\a+3,-10+\b)--(\a+3,-12+\b);
\draw[very thick](\a+5,-10+\b)--(\a+5,-12+\b);
\draw[very thick](\a+7,-10+\b)--(\a+7,-12+\b);
\draw[very thick](\a+9,-11.5000000000000+\b)--(\a+9,-12+\b);
\draw[very thick](\a+9.50000000000000,-11+\b)--(\a+10,-11+\b);
\draw[very thick](\a+9.50000000000000,-11+\b)arc(90:180:0.500000000000000);
\draw[very thick](\a+10,-11+\b)--(\a+12,-11+\b);
\draw[very thick](\a+11,-10+\b)--(\a+11,-12+\b);

\def\a{50};
\def\b{0};
\draw[step=2, color=green](\a+2,-12+\b)grid(\a+12,-2+\b);
\draw[very thick](\a+7,-3.50000000000000+\b)--(\a+7,-4+\b);
\draw[very thick](\a+7.50000000000000,-3+\b)--(\a+8,-3+\b);
\draw[very thick](\a+7.50000000000000,-3+\b)arc(90:180:0.500000000000000);
\draw[very thick](\a+8,-3+\b)--(\a+10,-3+\b);
\draw[very thick](\a+10,-3+\b)--(\a+12,-3+\b);
\draw[very thick](\a+5,-5.50000000000000+\b)--(\a+5,-6+\b);
\draw[very thick](\a+5.50000000000000,-5+\b)--(\a+6,-5+\b);
\draw[very thick](\a+5.50000000000000,-5+\b)arc(90:180:0.500000000000000);
\draw[very thick](\a+6,-5+\b)--(\a+8,-5+\b);
\draw[very thick](\a+7,-4+\b)--(\a+7,-6+\b);
\draw[very thick](\a+8,-5+\b)--(\a+10,-5+\b);
\draw[very thick](\a+10,-5+\b)--(\a+12,-5+\b);
\draw[very thick](\a+3,-7.50000000000000+\b)--(\a+3,-8+\b);
\draw[very thick](\a+3.50000000000000,-7+\b)--(\a+4,-7+\b);
\draw[very thick](\a+3.50000000000000,-7+\b)arc(90:180:0.500000000000000);
\draw[very thick](\a+5,-6+\b)--(\a+5,-6.50000000000000+\b);
\draw[very thick](\a+4,-7+\b)--(\a+4.50000000000000,-7+\b);
\draw[very thick](\a+4.50000000000000,-7+\b)arc(270:360:0.500000000000000);
\draw[very thick](\a+7,-6+\b)--(\a+7,-8+\b);
\draw[very thick](\a+11,-7.50000000000000+\b)--(\a+11,-8+\b);
\draw[very thick](\a+11.5000000000000,-7+\b)--(\a+12,-7+\b);
\draw[very thick](\a+11.5000000000000,-7+\b)arc(90:180:0.500000000000000);
\draw[very thick](\a+3,-8+\b)--(\a+3,-10+\b);
\draw[very thick](\a+5,-9.50000000000000+\b)--(\a+5,-10+\b);
\draw[very thick](\a+5.50000000000000,-9+\b)--(\a+6,-9+\b);
\draw[very thick](\a+5.50000000000000,-9+\b)arc(90:180:0.500000000000000);
\draw[very thick](\a+6,-9+\b)--(\a+8,-9+\b);
\draw[very thick](\a+7,-8+\b)--(\a+7,-10+\b);
\draw[very thick](\a+8,-9+\b)--(\a+10,-9+\b);
\draw[very thick](\a+10,-9+\b)--(\a+12,-9+\b);
\draw[very thick](\a+11,-8+\b)--(\a+11,-10+\b);
\draw[very thick](\a+3,-10+\b)--(\a+3,-12+\b);
\draw[very thick](\a+5,-10+\b)--(\a+5,-12+\b);
\draw[very thick](\a+7,-10+\b)--(\a+7,-12+\b);
\draw[very thick](\a+9,-11.5000000000000+\b)--(\a+9,-12+\b);
\draw[very thick](\a+9.50000000000000,-11+\b)--(\a+10,-11+\b);
\draw[very thick](\a+9.50000000000000,-11+\b)arc(90:180:0.500000000000000);
\draw[very thick](\a+10,-11+\b)--(\a+12,-11+\b);
\draw[very thick](\a+11,-10+\b)--(\a+11,-12+\b);
\node[right] at (\a+12,-3+\b) {$3$};
\node[right] at (\a+12,-5+\b) {$1$};
\node[right] at (\a+12,-7+\b) {$5$};
\node[right] at (\a+12,-9+\b) {$2$};
\node[right] at (\a+12,-11+\b) {$4$};
\node[below] at (\a+3,-12+\b) {$1$};
\node[below] at (\a+5,-12+\b) {$2$};
\node[below] at (\a+7,-12+\b) {$3$};
\node[below] at (\a+9,-12+\b) {$4$};
\node[below] at (\a+11,-12+\b) {$5$};

\draw[very thick,->] (15,-7)--(19,-7);
\draw[very thick,->] (31,-7)--(35,-7);
\draw[very thick,->] (47,-7)--(51,-7);

\end{tikzpicture}
\caption{An example of the map $x_\alpha{\rightsquigarrow}$ on BPD with $\alpha=1$ and $w=21534$}
\label{fig:BPD-Monk}
\end{figure}
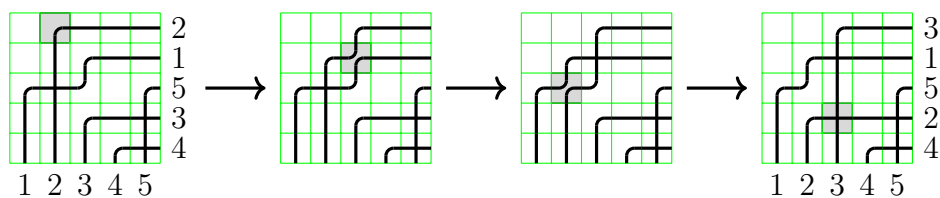

The maps described in \Cref{thm:huang-bump} are in the
spirit of adding a crossing and adjusting the crossing iteratively if
we encounter a double crossing, similar to \Cref{sub:pipes}. See \Cref{fig:BPD-Monk} for an
example.

We can define analogous maps $x_\alpha{\rightsquigarrow}$ and $m_{k,\alpha}$ on pipe dreams to prove \Cref{thm:Monk} bijectively. We will not spell out the details on these maps to avoid confusion with earlier material in \Cref{sub:Games} on little bumps. The bijection $\varphi$ in \Cref{def:bijection-PD-BPD} preserves the combinatorial proof of the variation of Monk's rule (\Cref{thm:Monk}) using both PDs and BPDs, which establishes its canonical nature. The proof is technical with many cases to analyze, so we will not include it here. Much of the difficulty comes from the fact that generalized chute moves on pipe dreams and droop moves on bumpless pipe dreams have few structural similarities. To be precise, we have the following theorem.
\begin{Theorem}\cite{GaoHuang-bijection}\label{thm:bijection-canonical}
The following diagrams commute for $w\in S_{\infty}$:
\[\begin{tikzcd}
\BPD(w)\arrow[r,"x_{\alpha}{\rightsquigarrow}"]\arrow[d,"\varphi"] & \bigcup_{\substack{l>\alpha\\w t_{\alpha,l}\gtrdot w}}\BPD(w t_{\alpha,l})\arrow[d,"\varphi"] \\
\rp(w)\arrow[r,"x_{\alpha}{\rightsquigarrow}"] & \bigcup_{\substack{l>\alpha\\w t_{\alpha,l}\gtrdot w}}\rp(w t_{\alpha,l})
\end{tikzcd},
\begin{tikzcd}
\BPD(w t_{k,\alpha})\arrow[r,"m_{k,\alpha}"]\arrow[d,"\varphi"] & \bigcup_{\substack{l>\alpha\\w t_{\alpha,l}\gtrdot w}}\BPD(w t_{\alpha,l})\arrow[d,"\varphi"] \\
\rp(w t_{k,\alpha})\arrow[r,"m_{k,\alpha}"] & \bigcup_{\substack{l>\alpha\\w t_{\alpha,l}\gtrdot w}}\rp(w t_{\alpha,l})
\end{tikzcd}\]
\end{Theorem}

There has been a lot of exciting development around bumpless pipe dreams since they were introduced. The above bijection can be generalized to \emph{marked} bumpless pipe dreams and not-necessarily-reduced bumpless pipe dreams \cite{huang2024markedbumplesspipedreamscompatible}, both of which can be used to compute the $\beta$-Grothendieck polynomials. 
Also, as an attempt to further understand the relations between pipe dreams and bumpless pipe dreams, Knutson and Udell \cite{knutson-udell-hybrid-FPSAC} defined and studied \emph{hybrid pipe dreams}, where each row of the tiling can be either \emph{ordinary} or \emph{bumpless}.
\begin{Definition}
A \emph{hybrid pipe dream} of type $\tau=\tau_1\cdots\tau_n\in\{O,B\}^n$ is a tiling of the $n\times n$ square grid such that 
\begin{itemize}
\item if $\tau_i=O$, called \emph{ordinary}, then its $i$th row is filled with tiles \+,\,\htile,\,\jt,\,\rt,\,\elbow,\,\bl,
\item if $\tau_i=B$, called \emph{bumpless}, then its $i$th row is filled with tiles \+,\,\htile,\,\rtrotatecounter,\,\rtrotateclock,\,\vtile,\,\bl,
\end{itemize}
forming pipes that go from the north boundary to the west and east boundaries, where no two pipes land on the same row and no two pipes cross twice. We label the pipes via the column numbers that they start with as before. Write the pipe labels on the endpoints on the west and east boundaries. We then obtain the corresponding permutation $w$ by reading the pipe labels in the counterclockwise direction from the west boundary to the south boundary to the east boundary. The monomial weight of a hybrid pipe dream of type $\tau$ is given by \[x^D:=\left(\prod_{\tau_i=O}\prod_{D(i,j)=\+\text{ or }\htile}x_i\right)\left(\prod_{\tau_i'=B}\prod_{D(i',j)=\bl}x_i'\right).\]
\end{Definition}
Write $\mathrm{HPD}(w,\tau)$ for all hybrid pipe dreams of $w$ of type $\tau$. 
\begin{figure}[h!]
\centering
\begin{tikzpicture}[scale=0.3]
\def\a{0};
\def\b{0};
\draw[step=2, color=green](\a+2,-12+\b)grid(\a+12,-2+\b);
\draw[very thick](\a+3,-3.5+\b)--(\a+3,-4+\b);
\draw[very thick](\a+3.5,-3+\b)--(\a+4,-3+\b);
\draw[very thick](\a+3.5,-3+\b)arc(90:180:0.5);
\draw[very thick](\a+3,-2+\b)--(\a+3,-2.5+\b);
\draw[very thick](\a+2,-3+\b)--(\a+2.5,-3+\b);
\draw[very thick](\a+2.5,-3+\b)arc(270:360:0.5);
\draw[very thick](\a+5,-3.5+\b)--(\a+5,-4+\b);
\draw[very thick](\a+5.5,-3+\b)--(\a+6,-3+\b);
\draw[very thick](\a+5.5,-3+\b)arc(90:180:0.5);
\draw[very thick](\a+5,-2+\b)--(\a+5,-2.5+\b);
\draw[very thick](\a+4,-3+\b)--(\a+4.5,-3+\b);
\draw[very thick](\a+4.5,-3+\b)arc(270:360:0.5);
\draw[very thick](\a+7,-3.5+\b)--(\a+7,-4+\b);
\draw[very thick](\a+7.5,-3+\b)--(\a+8,-3+\b);
\draw[very thick](\a+7.5,-3+\b)arc(90:180:0.5);
\draw[very thick](\a+7,-2+\b)--(\a+7,-2.5+\b);
\draw[very thick](\a+6,-3+\b)--(\a+6.5,-3+\b);
\draw[very thick](\a+6.5,-3+\b)arc(270:360:0.5);
\draw[very thick](\a+8,-3+\b)--(\a+10,-3+\b);
\draw[very thick](\a+9,-2+\b)--(\a+9,-4+\b);
\draw[very thick](\a+11,-2+\b)--(\a+11,-2.5+\b);
\draw[very thick](\a+10,-3+\b)--(\a+10.5,-3+\b);
\draw[very thick](\a+10.5,-3+\b)arc(270:360:0.5);
\draw[very thick](\a+3.5,-5+\b)--(\a+4,-5+\b);
\draw[very thick](\a+3,-4+\b)--(\a+3,-4.5+\b);
\draw[very thick](\a+3,-4.5+\b)arc(180:270:0.5);
\draw[very thick](\a+4,-5+\b)--(\a+6,-5+\b);
\draw[very thick](\a+5,-4+\b)--(\a+5,-6+\b);
\draw[very thick](\a+6,-5+\b)--(\a+8,-5+\b);
\draw[very thick](\a+7,-4+\b)--(\a+7,-6+\b);
\draw[very thick](\a+8,-5+\b)--(\a+10,-5+\b);
\draw[very thick](\a+9,-4+\b)--(\a+9,-6+\b);
\draw[very thick](\a+10,-5+\b)--(\a+12,-5+\b);
\draw[very thick](\a+2,-7+\b)--(\a+4,-7+\b);
\draw[very thick](\a+5,-7.5+\b)--(\a+5,-8+\b);
\draw[very thick](\a+5.5,-7+\b)--(\a+6,-7+\b);
\draw[very thick](\a+5.5,-7+\b)arc(90:180:0.5);
\draw[very thick](\a+5,-6+\b)--(\a+5,-6.5+\b);
\draw[very thick](\a+4,-7+\b)--(\a+4.5,-7+\b);
\draw[very thick](\a+4.5,-7+\b)arc(270:360:0.5);
\draw[very thick](\a+7,-7.5+\b)--(\a+7,-8+\b);
\draw[very thick](\a+7.5,-7+\b)--(\a+8,-7+\b);
\draw[very thick](\a+7.5,-7+\b)arc(90:180:0.5);
\draw[very thick](\a+7,-6+\b)--(\a+7,-6.5+\b);
\draw[very thick](\a+6,-7+\b)--(\a+6.5,-7+\b);
\draw[very thick](\a+6.5,-7+\b)arc(270:360:0.5);
\draw[very thick](\a+9,-6+\b)--(\a+9,-6.5+\b);
\draw[very thick](\a+8,-7+\b)--(\a+8.5,-7+\b);
\draw[very thick](\a+8.5,-7+\b)arc(270:360:0.5);
\draw[very thick](\a+5,-8+\b)--(\a+5,-10+\b);
\draw[very thick](\a+7.5,-9+\b)--(\a+8,-9+\b);
\draw[very thick](\a+7,-8+\b)--(\a+7,-8.5+\b);
\draw[very thick](\a+7,-8.5+\b)arc(180:270:0.5);
\draw[very thick](\a+8,-9+\b)--(\a+10,-9+\b);
\draw[very thick](\a+10,-9+\b)--(\a+12,-9+\b);
\draw[very thick](\a+2,-11+\b)--(\a+4,-11+\b);
\draw[very thick](\a+5,-10+\b)--(\a+5,-10.5+\b);
\draw[very thick](\a+4,-11+\b)--(\a+4.5,-11+\b);
\draw[very thick](\a+4.5,-11+\b)arc(270:360:0.5);
\node[above] at (\a+3,-2+\b) {$1$};
\node[above] at (\a+5,-2+\b) {$2$};
\node[above] at (\a+7,-2+\b) {$3$};
\node[above] at (\a+9,-2+\b) {$4$};
\node[above] at (\a+11,-2+\b) {$5$};
\node[right] at (\a+12,-3+\b) {$ $};
\node[right] at (\a+12,-5+\b) {$2$};
\node[right] at (\a+12,-7+\b) {$ $};
\node[right] at (\a+12,-9+\b) {$4$};
\node[right] at (\a+12,-11+\b) {$ $};
\node[left] at (\a+2,-3+\b) {$1$};
\node[left] at (\a+2,-5+\b) {$ $};
\node[left] at (\a+2,-7+\b) {$3$};
\node[left] at (\a+2,-9+\b) {$ $};
\node[left] at (\a+2,-11+\b) {$5$};
\end{tikzpicture}
\caption{A hybrid pipe dream of $w=13542$ and type $\tau=OBOBO$ with weight $x_1x_3x_4x_5$.} 
\end{figure}

\begin{Theorem}\cite{knutson-udell-hybrid-FPSAC} \label{thm:Knutson.Udell}
For all $w\in S_{n}$ and all $\tau\in\{O,B\}^n$, the Schubert polynomial satisfies \[\fS_{w}(x_1,x_2,\ldots,x_{n})=\sum_{D\in\mathrm{HPD}(w,\tau)}x^D.\]
\end{Theorem}

In fact, the notion of hybrid pipe
dreams gives us $2^n$ different models for Schubert polynomials!  We
have just scratched the surface here since much of the theory
around hybrid pipe dreams, including their relations with the above
defined bijection $\varphi$, is still in development. It is also worth
noting that Yu \cite{yu2024embeddingbumplesspipedreamsbruhat} has also
provided $(n-1)!$ different models for Schubert polynomials using
Bruhat chains.  Furthermore, Knutson and Zinn-Justin investigated
generating functions over \textit{generic pipe dreams} with a mixture
of the bumpless and bump tiles in their work on lower-upper varieties
and ways to interpolate between the classical and bumpless pipe dreams.

Let $D$ be a non-reduced classical or bumpless pipe dream in the
$n\times n$ grid.  Recall, we read the ``word'' of a classical pipe
dream in the reading order in \Cref{fig:rcgraphs}.  For bumpless
pipe dreams the reading order for the crossings goes up the columns
from left to right.  We can associate the \textit{natural} permutation
to $D$ in the usual way by observing the order of the exiting pipes or
we can associate the \textit{Demazure permutation} to $D$ by taking
the Demazure product of the reading word.  The \textit{Demazure
product} is defined recursively by saying $w\circ s_{i}$ is $ws_{i}$ if
$w<ws_{i}$ and $w$ if $w>ws_{i}$.  One can visualize the Demazure
permutation of $D$ by considering the crossings in reading order one
at a time, and if the pipes crossing at $(i,j)$ have previously
crossed, replace the crossing at $(i,j)$ with a bump tile.  The
\textit{Demazure permutation} for $D$ is the permutation naturally
associated with the exiting pipes in $D$.  This is the product used
for computing Grothendieck polynomials.

\begin{Exercise}
Prove that the natural permutation and the Demazure permutation for a
pipe dream $D$ agree if and only if $D$ is reduced.  
\end{Exercise}

\begin{Exercise}\cite[Lem. 7.2]{weigandt-BPD-ASM}
Prove that $w \in S_{n}$ is vexillary if and only if every
(not-necessarily-reduced) bumpless pipe dream with Demazure
permutation $w$ is reduced.  Does the same property hold for classical
pipe dreams?   %%%  Answer: no, w=1243 has non-reduced pipe dreams.
\end{Exercise}

\subsection{Schubitopes, Saturated Newton Polytopes and Vanishing Coefficients}\label{sub:VanishingCoefs}
In this section, we continue our study on the rich combinatorial structures of pipe dreams and related objects, with connections to polytopes, matroids and computational complexity. 
\begin{Definition}  
For any $\alpha=(\alpha_{1},\dots
,\alpha_{n})\in\mathbb{Z}_{\geq0}^n$, let $x^{\alpha}$ denote the
monomial $x^{\alpha}=x_{1}^{\alpha_{1}}\cdots x_{1}^{\alpha_{n}}$. The
\emph{Newton polytope} of a polynomial
$f=\sum_{\alpha\in\mathbb{Z}_{\geq0}^n} c_{\alpha}x^{\alpha }$ is the
convex hull of its exponent vectors
\[N(f):=\mathrm{conv}(\{\alpha\:|\: c_{\alpha}\neq0\})\in\mathbb{R}^n.\]
We say that a polynomial $f$ has \emph{saturated Newton polytope} (SNP) if $c_{\alpha}\neq0$ whenever $\alpha\in N(f)\cap \mathbb{Z}_{\geq0}^n.$
\end{Definition}

Many polynomials that arise naturally in algebra and combinatorics
have SNP. Unsurprisingly, we will be focusing on the Newton polytope
of the Schubert polynomials, which are called \emph{Schubitopes}
defined by Monical, Tokcan, and Yong
\cite{monical-tokcan-yong-newton-polytope}.

Recall from \Cref{def:diagram} that a \emph{diagram} is a subset of
boxes in the grid $[n]^2$. For a diagram $D$, let
\begin{equation}\label{eq:Diagram.Dj}
D_j=\{i\given (i,j)\in D\}
\end{equation}
denote the rows of the boxes in the  $j^{th}$ column of $D$, for $j=1,\ldots,n$.
\begin{Definition}\cite{monical-tokcan-yong-newton-polytope}\label{def:schubitope}
Let $D\subset[n]^2$ be a diagram. For $I\subset[n]$ and $j\in[n]$, construct $\mathrm{word}_{j,I}(D)$ by reading through column $j$ of $D$ from top to bottom and recording:
\begin{itemize}
\item $($ if $(i,j)\notin D$ and $i\in I$;
\item $)$ if $(i,j)\in D$ and $i\notin I$;
\item $*$ if $(i,j)\in D$ and $i\in I$.
\end{itemize}
Then define $\theta_D^j(I)=\#\text{ paired }()'s\text{ in
}\mathrm{word}_{j,I}(D)\ +\#*'s\text{ in }\mathrm{word}_{j,I}(D)$, and
define 
\[\theta_D(I)=\sum_{j=1}^n\theta_D^j(I).\]
The \emph{Schubitope} $\mathcal{S}_D$ associated to $D$ is the polytope
\[\mathcal{S}_D:=\left\{(\alpha_1,\ldots,\alpha_n)\in\mathbb{R}^n_{\geq0}\::\:\sum\alpha_i=\#D\text{ and }\sum_{i\in I}\alpha_i\leq \theta_D(I)\text{ for all }I\subset[n]\right\}.\]
For a permutation $w\in S_n$, its \emph{Schubitope} $\mathcal{S}_w$ is defined as $\mathcal{S}_{D(w)}$ where $D(w)$ is the Rothe diagram of $w$.
\end{Definition}

\begin{figure}[h!]
\centering
\begin{tikzpicture}[scale=0.5]
\draw[step=1.0,green,thin] (0,0) grid (4,4);
\draw[very thick] (3,0)--(4,0)--(4,1)--(3,1)--(3,0);
\draw[very thick] (0,4)--(2,4)--(2,3)--(0,3);
\draw[very thick] (0,4)--(0,2)--(1,2)--(1,4);
\end{tikzpicture}
\qquad
\begin{tikzpicture}[scale=0.5]
\draw[step=1.0,green,thin] (0,0) grid (4,4);
\draw[very thick] (3,0)--(4,0)--(4,1)--(3,1)--(3,0);
\draw[very thick] (0,4)--(2,4)--(2,3)--(0,3);
\draw[very thick] (0,4)--(0,2)--(1,2)--(1,4);
\node at (0.5,3.5) {$*$};
\node at (0.5,2.5) {$*$};
\node at (1.5,3.5) {$*$};
\node at (2.5,3.5) {$($};
\node at (3.5,3.5) {$($};
\node at (1.5,2.5) {$($};
\node at (2.5,2.5) {$($};
\node at (3.5,2.5) {$($};
\node at (3.5,0.5) {$)$};
\end{tikzpicture}
\qquad
\begin{tikzpicture}[scale=0.5]
\draw[step=1.0,green,thin] (0,0) grid (4,4);
\draw[very thick] (3,0)--(4,0)--(4,1)--(3,1)--(3,0);
\draw[very thick] (0,4)--(2,4)--(2,3)--(0,3);
\draw[very thick] (0,4)--(0,2)--(1,2)--(1,4);
\node at (0.5,1.5) {$($};
\node at (1.5,1.5) {$($};
\node at (2.5,1.5) {$($};
\node at (3.5,1.5) {$($};
\node at (3.5,0.5) {$)$};
\node at (0.5,2.5) {$)$};
\node at (0.5,3.5) {$)$};
\node at (1.5,3.5) {$)$};
\end{tikzpicture}
\caption{Left: a diagram $D$; middle: construction for $\mathrm{word}_{j,I}$ for $I=\{1,2\}$; right: construction for $\mathrm{word}_{j,I}$ for $I=\{3\}$.}
\label{fig:schubitope-example}
\end{figure}
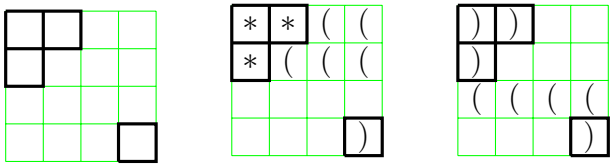

\begin{Example}
Consider the diagram $D\in[4]^2$ in \Cref{fig:schubitope-example} with
$n=4$ and $D_1=\{1,2\},D_2=\{1\},D_3=\emptyset,D_4=\{4\}$. Take
$I=\{1,2\}$ and fill in the squares according to
\Cref{def:schubitope}. We see that $\theta_D^1(\{1,2\})=2$,
$\theta_D^2(\{1,2\})=1$, $\theta_D^3(\{1,2\})=0$ and
$\theta_D^4(\{1,2\})=1$ where one pair of $()$ is formed and contributes to $\theta_D^4(\{1,2\})$. Thus,
$\theta_D(\{1,2\})=4$. For $I=\{3\}$, we see that there are no $*$'s and only one pair of
$()$ can be formed, which is in column $4$. Therefore,
$\theta_D(\{3\})=1$.  So each $ (\alpha_{1},\alpha_{2},\alpha_{3}, \alpha_{4}) \in
\mathcal{S}_D$ satisfies the equations $\alpha_{1}+\alpha_{2}\leq 4$
and $\alpha_{3}\leq 1$.  In fact, one can check that the following
selections of $I$ are sufficient to define the
Schubitope
\begin{equation}\label{eq:Schubitope.D}
\mathcal{S}_D=\left\{(\alpha_1,\alpha_2,\alpha_3,\alpha_4)\in\mathbb{R}_{\geq0}^4\:\Biggm|\: \begin{split}\alpha_1+\alpha_3+\alpha_4\leq&\,\,3\\\alpha_2+\alpha_3+\alpha_4\leq&\,\,2\\\alpha_3+\alpha_4\leq&\,\,1\\\alpha_1+\alpha_2+\alpha_3+\alpha_4=&\,\,4\end{split}\right\}.
\end{equation}
\end{Example}

\begin{Theorem}[Conjectured by Monical-Tokcan-Yong
\cite{monical-tokcan-yong-newton-polytope}, proved by
Fink-M\'{e}sz\'{a}ros-St\@.\,Dizier \cite{fink-meszaros-dizier-Schubert}]\label{thm:schubert-SNP}
The Schubert polynomial $\fS_w$ has saturated Newton polytope given by $\mathcal{S}_{w}=\mathcal{S}_{D(w)}$. 
\end{Theorem}

In fact, Fink, M\'{e}sz\'{a}ros and St\@.\,Dizier
\cite{fink-meszaros-dizier-Schubert} showed that the Newton polytopes
of Schubert polynomials, or more generally, the Newton polytopes of
\emph{dual characters of flagged Weyl modules}, can be written as a
Minkowski sum of matroid polytopes. Kra\'{s}kiewicz and Pragacz
\cite{kraskiewicz-pragacz-foncteurs} showed that Schubert polynomials
are dual characters of flagged Weyl modules, establishing the
connection. We explain here how this decomposition works, see
\Cref{thm:newton-schubert-minkowski}. See also \cite{anderson2024} for
an alternative approach.

Before stating the theorem, we need some notation and vocabulary from
the theory of matroids. See \cite{Arila.2023,Oxley} for more background on
matroids if needed.  Let $M$ be a matroid of rank $k$ on ground set
$[n]$ represented by its collection of bases. The bases of $M$ are
size $k$ subsets of $[n]$. Recall from
\Cref{def:dominated.multisets.different.size} that the Gale partial
order on subsets of $[n]$ of size $k$ is defined by $A \trianglelefteq
B$ if and only if after sorting $A=\{a_1<\cdots<a_k\}$,
$B=\{b_1<\cdots<b_k\}$, we have $a_i\leq b_i$ for all $i=1,\ldots,k$.
Hence the Gale order induces a partial order on the bases of $M$.  Let
$e_A=\sum_{i\in A}e_i\in\mathbb{R}^n$ where $e_{i}$ is the $i^{th}$
standard basis vector with a 1 in position 1 and 0's elsewhere.  Let
$P(M)$ denote the \emph{matroid polytope} of $M$, which is the convex
hull of the vectors $e_B$ for each base $B$ of $M$.  Putting these
concepts together gives rise to the Schubert matroid and polytope.

%% Derkson-Fink say Crapo originally studied Schubert matroids in his
%% 1964 paper, but Fink-Meszaros-St.Dizer cite Section 2.4 of Oriented
%% Matroids book. 
\begin{Definition}\label{def:schubert-matroid-polytope}
Let $\mathrm{SM}_n(B)$ denote the \emph{Schubert matroid} whose bases
are precisely all the sets $A$ that are less than or equal to $B$ in
the Gale order.    For $B\subset[n]$ of size $k$, the \emph{Schubert matroid
polytope} is \[P(\mathrm{SM}_n(B)):=\mathrm{conv}(\{e_A\given A
\trianglelefteq B\text{ in Gale order}\}).\]
\end{Definition}

\begin{Theorem}\cite{fink-meszaros-dizier-Schubert}\label{thm:newton-schubert-minkowski}
For $w\in S_n$, the Newton polytope of the Schubert polynomial is a
Minkowski sum of Schubert matroid polytopes. That
is, \[N(\fS_w)=\mathcal{S}_w=\sum_{j=1}^n P(\mathrm{SM}_n(D_j)),\]
where $D_{j}=\{i\given (i,j)\in D(w)\}$.
\end{Theorem}
In particular, \Cref{thm:newton-schubert-minkowski} is saying that the
Newton polytope of a Schubert polynomial is a \emph{generalized
permutahedron}, introduced by Postnikov \cite{postnikov-permutahedra}
and now widely used in combinatorics.

\begin{Example}
Consider $w=32154$ with the Schubert
polynomial
\begin{equation}\label{eq:32154}
\fS_{32154}=x_1^3x_2+x_1^2x_2^2+x_1^2x_2x_3+x_1^2x_2x_4=x_1^2x_2(x_1+x_2+x_3+x_4)
\end{equation}
with permutation diagram given in \Cref{fig:schubitope-example}.  One
can verify from the monomial expansion that its Schubitope is
$\mathcal{S}_{D(w)}$ from \Cref{fig:schubitope-example}.  Let's list the
nontrivial Minkowski summands of $\mathcal{S}_{32154}$ as in
\Cref{thm:newton-schubert-minkowski}:
\begin{align*}
D_1=\{1,2\},&\quad \mathrm{SM}_n(D_1)=\mathrm{conv}(\{e_1+e_2\}),\\
D_2=\{1\},&\quad \mathrm{SM}_n(D_2)=\mathrm{conv}(\{e_1\}),\\
D_4=\{4\},&\quad \mathrm{SM}_n(D_4)=\mathrm{conv}(\{e_1,e_2,e_3,e_4\}).
\end{align*}
Comparing \eqref{eq:32154}, one observes that
$\mathcal{S}_{32154}=P(\mathrm{SM}_n(D_1))+P(\mathrm{SM}_n(D_2))+P(\mathrm{SM}_n(D_4))$
in this case as expected from the theorem. 
\end{Example}

As an application of the framework of \Cref{thm:schubert-SNP}, Adve-Robichaux-Yong \cite{adve-robichaux-yong-efficient-vanishing-schubert} gave a polynomial time algorithm deciding whether a certain coefficient $x^{\alpha}$ of a Schubert polynomial $\fS_w$ vanishes or not. Note that there are $O(2^n)$ inequalities in the description of the Schubitope (\Cref{def:schubitope}). Many of these inequalities are redundant. The first step is to find a ``short witness" for all these inequalities.
\begin{Definition}\cite{adve-robichaux-yong-efficient-vanishing-schubert}\label{def:perfect-tableaux}
A tableau of diagram shape $D$ is called \emph{perfect} if the following two conditions are satisfied:
\begin{itemize}
\item column-injectivity: there is at most one copy of $i$ in each column, for all $i$;
\item flaggedness: all numbers in row $i$ do not exceed $i$, for all $i$.
\end{itemize}
\end{Definition}
\begin{Theorem}\cite{adve-robichaux-yong-efficient-vanishing-schubert}\label{thm:perfect-tableaux}
There exists a perfect tableau of shape $D(w)$ and content $\alpha$ if and only if the coefficient of $x^{\alpha}$ in the Schubert polynomial $\fS_w$ is positive. Moreover, in such cases, there exists a perfect tableau which increases along columns.
\end{Theorem}
\Cref{thm:perfect-tableaux} should be contrasted to
\Cref{thm:balanced-schubert}. Essentially, the balanced condition from
\Cref{def:balanced-flagged-column-injective} provides the correct coefficient, whereas column-injectivity and flagged can already indicate whether the coefficient of interest is zero or not.

The proof of \Cref{thm:perfect-tableaux} depends on the
characterizations of the Schubitope in \Cref{def:schubitope} and
\Cref{thm:schubert-SNP}.
\begin{Example}
Consider $w=1432$ with $\fS_w=x_1^2x_2+x_1^2x_3+x_1x_2^2+x_1x_2x_3+x_2^2x_3$. The perfect tableaux of shape $D(w)$ are shown in \Cref{fig:perfect-tableaux}. See \Cref{fig:column-injective-balanced-labellings} for a comparison with column-injective balanced labellings of the same shape.
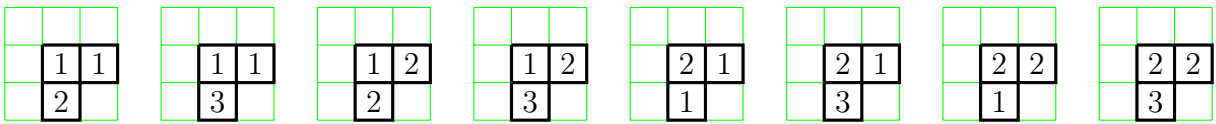
\begin{figure}[h!]
\centering
\begin{tikzpicture}[scale=0.5]
\draw[step=1.0,green,thin] (0,0) grid (3,3);
\draw[very thick] (1,2)--(3,2)--(3,1)--(1,1);
\draw[very thick] (1,2)--(1,0)--(2,0)--(2,2);
\node at (1.5,1.5) {$1$};
\node at (2.5,1.5) {$1$};
\node at (1.5,0.5) {$2$};
\end{tikzpicture}
\quad
\begin{tikzpicture}[scale=0.5]
\draw[step=1.0,green,thin] (0,0) grid (3,3);
\draw[very thick] (1,2)--(3,2)--(3,1)--(1,1);
\draw[very thick] (1,2)--(1,0)--(2,0)--(2,2);
\node at (1.5,1.5) {$1$};
\node at (2.5,1.5) {$1$};
\node at (1.5,0.5) {$3$};
\end{tikzpicture}
\quad
\begin{tikzpicture}[scale=0.5]
\draw[step=1.0,green,thin] (0,0) grid (3,3);
\draw[very thick] (1,2)--(3,2)--(3,1)--(1,1);
\draw[very thick] (1,2)--(1,0)--(2,0)--(2,2);
\node at (1.5,1.5) {$1$};
\node at (2.5,1.5) {$2$};
\node at (1.5,0.5) {$2$};
\end{tikzpicture}
\quad
\begin{tikzpicture}[scale=0.5]
\draw[step=1.0,green,thin] (0,0) grid (3,3);
\draw[very thick] (1,2)--(3,2)--(3,1)--(1,1);
\draw[very thick] (1,2)--(1,0)--(2,0)--(2,2);
\node at (1.5,1.5) {$1$};
\node at (2.5,1.5) {$2$};
\node at (1.5,0.5) {$3$};
\end{tikzpicture}
\quad
\begin{tikzpicture}[scale=0.5]
\draw[step=1.0,green,thin] (0,0) grid (3,3);
\draw[very thick] (1,2)--(3,2)--(3,1)--(1,1);
\draw[very thick] (1,2)--(1,0)--(2,0)--(2,2);
\node at (1.5,1.5) {$2$};
\node at (2.5,1.5) {$1$};
\node at (1.5,0.5) {$1$};
\end{tikzpicture}
\quad
\begin{tikzpicture}[scale=0.5]
\draw[step=1.0,green,thin] (0,0) grid (3,3);
\draw[very thick] (1,2)--(3,2)--(3,1)--(1,1);
\draw[very thick] (1,2)--(1,0)--(2,0)--(2,2);
\node at (1.5,1.5) {$2$};
\node at (2.5,1.5) {$1$};
\node at (1.5,0.5) {$3$};
\end{tikzpicture}
\quad
\begin{tikzpicture}[scale=0.5]
\draw[step=1.0,green,thin] (0,0) grid (3,3);
\draw[very thick] (1,2)--(3,2)--(3,1)--(1,1);
\draw[very thick] (1,2)--(1,0)--(2,0)--(2,2);
\node at (1.5,1.5) {$2$};
\node at (2.5,1.5) {$2$};
\node at (1.5,0.5) {$1$};
\end{tikzpicture}
\quad
\begin{tikzpicture}[scale=0.5]
\draw[step=1.0,green,thin] (0,0) grid (3,3);
\draw[very thick] (1,2)--(3,2)--(3,1)--(1,1);
\draw[very thick] (1,2)--(1,0)--(2,0)--(2,2);
\node at (1.5,1.5) {$2$};
\node at (2.5,1.5) {$2$};
\node at (1.5,0.5) {$3$};
\end{tikzpicture}
\caption{Perfect tableaux of shape $D(1432)$}
\label{fig:perfect-tableaux}
\end{figure}
\end{Example}

We can represent the data of a perfect tableau using a matrix $\{a_{ij}\}_{i,j=1}^n$ where $a_{ij}$ denotes the number of $i$'s in column $j$. Rewriting the conditions from \Cref{def:perfect-tableaux}, we arrive at the following system of linear inequalities:
\begin{equation}\label{eq:perfect-tableaux}\begin{cases}
\text{column-injectivity: }&0\leq a_{ij}\leq 1,\text{ for }1\leq i,j\leq n,\\
\text{content: }&\sum_{j=1}^n a_{ij}=\alpha_i,\text{ for }1\leq j\leq n,\\
\text{flagged: }&\sum_{i=k+1}^n a_{ij}\leq \#\{s>k\given s\in D_j\},\text{ for }1\leq j,k\leq n.
\end{cases}\end{equation}
A perfect tableau of shape $D$ and content $\alpha$ exists if and only
if the system of linear inequalities in \eqref{eq:perfect-tableaux}
has an \textbf{integral} solution. Integer linear programming is
unfortunately still NP-hard, but the good news in our case is that all
the constraints from \eqref{eq:perfect-tableaux} are \emph{totally
unimodular}. The phrase \emph{totally unimodular} means that if we
write down all the linear inequalities using $M\vec{a}\leq\vec{b}$
where $\vec{a}$ is the vector consisting of all $a_{ij}'s$, then all
the maximal minors of $M$ have determinants belonging to
$\{0,\pm1\}$. By Cramer's rule, the polytope defined by the
inequalities in \eqref{eq:perfect-tableaux} has integral
vertices. Therefore, this polytope contains a lattice point if and
only if it is nonempty.  Determining if a polytope is non-empty is a
much easier problem in general.  In fact, there are polynomial-time
algorithms to solve linear programming problems without the integer
constraint. For example, the ellipsoid algorithm \cite{khazhiyan} from optimization
works in polynomial time.  In summary, determining the nonzero
coefficients of a Schubert polynomial is relatively easy compared to
determining its coefficients explicitly.

\begin{Theorem}\cite{adve-robichaux-yong-efficient-vanishing-schubert}
There exists a perfect tableau of shape $D$ and content $\alpha$ if and only if there exists a solution for \eqref{eq:perfect-tableaux} in $\mathbb{R}^{n^2}$. Consequently, deciding whether the coefficient of $x^{\alpha}$ in $\fS_w$ equals $0$ can be done in polynomial time.
\end{Theorem}

Given that the Newton polytope of the Schubert polynomial $\fS_w$ is
saturated (\Cref{thm:schubert-SNP}), it is natural to ask when $\fS_w$
precisely equals the sum of monomials in $\mathcal{S}_w$?  In other words,
when does the monomial expansion of $\fS_w$ have coefficients in
$\{0,1\}$? Such Schubert polynomials are called \emph{zero-one
polynomials} or \emph{multiplicity-free}.  These were completely
characterized by Fink-M\'{e}sz\'{a}ros-St\@.\,Dizier
\cite{zero-one-schubert} using pattern avoidance (\Cref{def:pattern-avoidance}).

\begin{Theorem}\cite{zero-one-schubert}\label{thm:zero-one-schubert}
The Schubert polynomial $\fS_w$ is zero-one if and only if $w$ avoids
the 12 patterns
\[
12543, 13254, 13524, 13542, 21543, 125364,
\]
\[
125634,215364, 215634, 315264, 315624,315642.
\]
\end{Theorem}

A generalization of \Cref{thm:zero-one-schubert} classifying zero-one dual characters of flagged Weyl modules is conjectured by \cite{meszaros-stdizier-tanjaya} and established by \cite{guo2024zeroonedualcharactersflagged}.

\begin{Exercise}
What is the coefficient of  $x_1^5 x_2^6 x_3^5 x_4^3 x_5^2 x_6 x_7$ in
$\fS_{2 8 9 7 4 5 3 1 6}$?
\end{Exercise}

% %(sort (unique-elems  (loop for term in  (schub '(3 1 2 8 5 4 7 6) 'a) collect (first-n 2 (term-mon term)))) 'lex-order)
% ((2 0) (2 1) (2 2) (2 3) (2 4) (3 0) (3 1) (3 2) (3 3) (3 4) (4 0) (4 1) (4 2)
%  (4 3) (4 4) (5 0) (5 1) (5 2) (5 3) (6 0) (6 1) (6 2))
% CL-USER> 
\begin{Exercise}
Let $w=3 1 2 8 5 4 7 6\in S_{8}$.  Let $P_{w}(2)$ be the set of
vectors $(\alpha_1,\alpha_2)$ such that there exists a vector
$\alpha=(\alpha_1,\alpha_2,\dots , \alpha_8)$ such that $x^{\alpha}$
appears in $\fS_{w}$ with nonzero coefficient.
\begin{enumerate}
\item Sketch $P_{w}(2)$.  
\item Is $P_{w}(2)$ a Minkowski sum of Schubitopes?  If so, what are
the Minkowski summands.  If not, explain why it cannot be so.
\item For each $(i,j) \in P_{w}(2)$, find the set $R_{w}(i,j)$ of all
realizable coefficients $c_{\alpha}$ such that $\alpha_{1}=i$ and
$\alpha_{2}=j$ in the expansion $\fS_{w}= \sum c_{\alpha}x^{\alpha}$.
Add these sets to your sketch of $P_{w}(2)$.  How do the sets
$R_{w}(i,j)$ relate to each other?  Formulate a conjecture about some
phenomena you may observe.  Then, try to prove your conjecture holds
for all $w \in S_{n}$ or find a counterexample to disprove it.
\end{enumerate}
\end{Exercise}

\begin{Exercise}\cite{Castillo_Cid-Ruiz_Mohammadi_Montano_2023,monical-tokcan-yong-newton-polytope}
\label{ex:double.schub.snp}
Do the double Schubert polynomials $\fS_{w}(X,Y)$ have the saturated
Newton polytope property?
\end{Exercise}

\subsection{1-2-3 Step Flags and Puzzles}\label{sub:123.StepFlags}

In our hurry to mine the rich structure of the complete flag variety,
we have skipped over the Schubert calculus of some simpler related
spaces like projective space, Grassmannians, and partial flag
varieties.  Let's start by recalling some notation from
\Cref{sub:Grassmannians.intro} to connect the Schubert varieties in
partial flag varieties with the corresponding special Schubert
polynomials.  The corresponding Schubert structure constants can be seen
through the lens of ``puzzles,'' first introduced by Knutson and 
built upon by many authors including
\cite{Buch-Kresch-Purbhoo.2-step-rule,Buch-Kresch-Tamvakis,KnutsonTao-puzzles,
Knutson-Tao-Woodward, knutsonZin-Justin2020schubert}.  We start by
reviewing some vocabulary from \Cref{sub:Grassmannians.intro}.  

\begin{Definition}
Given a subset ${\bf d} = \{d_1 < \cdots < d_m\} \subseteq [n-1]$, a
\emph{partial flag} with dimensions ${\bf d}$ is a sequence of $m$
nested subspaces $F_1 \subseteq \cdots \subseteq F_m \subseteq \C^n$
with $\dim F_i = d_i$. The set of all such partial flags is the
\emph{partial flag variety} $\Fl(n; {\bf d})$.  We say $f \in \Fl(n;
{\bf d})$ is an $m$-\emph{step flag} if $|{\bf d}|=m$.
\end{Definition}

One can give $\Fl(n; {\mathbf d})$ the structure of a compact smooth manifold or a projective complex variety in more or less the same way as $\Fl(n)$. There is an obvious projection map $\pi : \Fl(n) \to \Fl(n; {\mathbf d})$ which ``forgets'' those components of a complete flag with dimensions not in ${\mathbf d}$. 

\begin{Example}
The \emph{Grassmannian} of $k$-planes in $\C^n$ is 
\begin{equation*}
    \Gr(k,n) = \Fl(n; \{k\}) = \{V \subseteq \C^n\text{ a linear subspace} \given \dim V = k\}.
\end{equation*}
In particular, $\Gr(1,n)$ is complex projective space $\C\bP^{n-1}$,
the space of (complex) lines through the origin in $\C^n$.  The
partial flag variety $\Fl(n; \{1,n-1\})$ is the space of pairs $(L,H)$
where $L$ is a line contained in a hyperplane $H$. Fixing an inner
product on $\C^n$ and sending $(L,H) \mapsto (L,H^\perp)$ gives an
isomorphism to the space of pairs of orthogonal lines in $\C^n$.
\end{Example}

Given $E_{\bullet} \in \Fl(n)$ and a permutation $w \in S_n$, we
define an associated Schubert variety in $\Fl(n; {\mathbf d})$ almost
exactly as in the complete flag variety, namely 
\begin{equation*}
    X_w(E_\bullet; {\mathbf d}) = \{F_\bullet \in \Fl(n; {\mathbf d}) \given \dim(E_i \cap F_{j}) \geq \rk(w)[i,j] \text{ for all $(i,j) \in [n] \times {\mathbf d}$}\}.
\end{equation*}
The Schubert cell $C_w(E_\bullet; {\mathbf d})$ is defined the same
way but replacing $\geq$ with $=$. Evidently the projection
$\pi(C_w(E_\bullet)) = C_w(E_\bullet; {\mathbf d})$. As before, we
drop the $E_\bullet$ or the ${\mathbf d}$ from the notation when clear
from context. Just as in the complete flag variety, we have:
\begin{itemize}
    \item Each $C_w$ is isomorphic to an affine space.
    \item $X_w$ is the Zariski (or Euclidean) closure of $\overline{C}_w$.
    \item The left $B$-orbits on $\Fl(n; {\mathbf d})$ are the distinct Schubert cells.
\end{itemize}

\begin{Example} \label{ex:4-lines-Gr}
Recall from \Cref{sub:SchubertProblems2000.reprise} that lines in
$\bR\bP^3$ naturally correspond to $2$-planes through the origin in
$\bR^4$. Given a line $L$ in $\bR\bP^3$, let $E_\bullet^L$ be any flag
in $\Fl(4)$ with ordered basis over $\bR$ such that $E_2^L$ is the
$2$-plane corresponding to $L$. Then
\begin{equation*}
X_{4231}(E_{\bullet}^L; \{2\}) = \{V \in \Gr_{\bR}(2,4) \given \dim(V \cap E_2^L) \geq 1\} = \{\ell \in \bR\bP^3 \given \ell \cap L \neq \emptyset\}.
\end{equation*}
Recall, the classic Schubert problem states that given any four
generic lines $L_1, \ldots, L_4$ in 3-space, there are exactly 2 lines
intersecting all of them from \Cref{sub:SchubertProblems1900} and the
reprise in \Cref{sub:SchubertProblems2000.reprise}.  We now see how to
phrase this in terms of Schubert varieties in a partial flag variety:
the intersection $\bigcap_{j=1}^4 X_{4231}(E_\bullet^{L_j}; \{2\})$
consists of two points.
\end{Example}

Our indexing is somewhat problematic as $X_w(E_\bullet; {\mathbf d})$ does not uniquely determine $w$ when ${\mathbf d} \neq [n-1]$. For instance, we could have used $w = 4213$ instead of $4231$ in Example~\ref{ex:4-lines-Gr}, given that $X_w(E_{\bullet}^L; \{2\})$ only depends on column 2 of the rank table $\rk[w]$.

\begin{Proposition} \label{prop:partial-flag} Fix a partial flag
variety $\Fl(n;\mathbf{d})$ and a flag $E_\bullet \in \Fl(n)$.  Let $v,w
\in S_{n}$.
\begin{enumerate}[(a)]
\item The (partial flag) Schubert varieties $X_v(E_\bullet; {\mathbf
d})$ and $X_w(E_\bullet; {\mathbf d})$ are equal if and only if
$v([d_j+1,d_{j+1}]) = w([d_j+1,d_{j+1}])$ for $j = 0, \ldots, m$,
where we set $d_0 = 0$ and $d_{m+1} = n$.
\item If $w$ is decreasing on each interval $[d_j+1, d_{j+1}]$, then
$\pi^{-1}(X_w(E_\bullet, {\mathbf d})) = X_w(E_\bullet) \subset  \Fl(n) $.
% \item If $w$ is as in (b), then $\dim X_w(\mathbf{d}) = \ell(w) - \sum_{j=1}^{m+1} {d_{j+1}-d_j \choose 2}$.
\end{enumerate}
\end{Proposition}

\begin{proof} \hfill
\begin{enumerate}[(a)]
\item Modifying $v$ by permuting the entries $v(d_{j}+1), v(d_j+2),
\ldots, v(d_{j+1})$ among themselves does not change the rank table
$\rk[v]$ in columns $d_1, \ldots, d_m$, which is all the data needed
to determine $X_v(E_\bullet; {\mathbf d})$. We leave it as an exercise
to check the converse: that if one of these columns is changed, then
so is $X_v(E_\bullet; {\mathbf d})$.

\item Let $X_{w}$ denote the Schubert variety in $\Fl(n)$ indexed by
$w$ with respect to $E_{\bullet}$, and let $X_w({\mathbf d})$ denote
the projection in $\Fl(n;\mathbf{d})$.   Since $\pi(X_w) = X_w({\mathbf
d})$, we have $\pi^{-1}(X_w(\mathbf{d})) = \bigcup_{v \in S} X_v$
where $S = \{v \in S_n \given X_v({\mathbf d}) = X_w({\mathbf d})\}$
by Part (a).
If $w$ is decreasing on each interval $[d_j+1, d_{j+1}]$ as in (b),
then $w \geq v$ in Bruhat order for all $v \in S$ by the Ehresmann
Tableau Criterion \Cref{thm:Ehresmann}.  Therefore, by definition of the
Bruhat order
\[
\pi^{-1}(X_w(\mathbf{d}))=\bigcup_{v \in S} X_v = X_w.
\]

\end{enumerate}
\end{proof}

Let $S_{\mathbf d}$ be the \emph{parabolic subgroup} (or \emph{Young
subgroup}) of $S_n$ whose elements map each interval $[d_j+1, d_{j+1}]$
to itself. Then Proposition~\ref{prop:partial-flag}(a) says that
$X_v({\mathbf d}) = X_w({\mathbf d})$ if and only if $vS_{\mathbf d} =
wS_{\mathbf d}$. Each right coset of $S_{\mathbf d}$ contains a unique
Bruhat-maximal element and a unique Bruhat-minimal element.  Since
\Cref{prop:partial-flag}(b) uses the Bruhat-maximal elements, we
identify the cosets in $S_n/S_{\mathbf d}$ with the set of such
maximal permutations and use them to index the distinct Schubert cells
and varieties in $\Fl(n,{\mathbf d})=\bigcup_{v \in S_n/S_{\mathbf d}}
C_v({\mathbf{d}})$.  Bruhat order on $S_{n}$ induces the partial
order given by containment of Schubert varieties for the partial flag
variety also, so
\begin{equation}\label{eq:partial.bruhat}
X_w({\mathbf d}) = \bigcup_{\substack{v \in S_n/S_{\mathbf d} \\
v\leq w}} C_v({\mathbf{d}}).
\end{equation}

Recall that the cohomology class $[X_{w}]$ can be represented in the
coinvariant algebra $R_{n}$ by $\fS_{w_{0}w}$ by the Inherited
Positivity \Cref{cor:Schub.basis}, so if $w$ is Bruhat-maximal in its
coset, then $w_{0}w$ is Bruhat-minimal in its coset.  As in the full
flag variety case, Schubert polynomials form a basis for $H^*(\Fl(n;
{\mathbf d}))$ as shown in the following theorem.

\begin{Theorem} \label{thm:partial-flag-cohom}
    The span of $\{\fS_{w_{0}w} : w \in S_n/S_{\mathbf d}\}$ in the coinvariant algebra $R_n \simeq A^*(\flags)$ is a subring isomorphic to $A^*(\Fl(n; {\mathbf d})) \simeq H^*(\Fl(n; {\mathbf d}))$ via the map $\fS_{w_{0}w} \mapsto [X_{w}(\mathbf{d})]$.
\end{Theorem}
\begin{proof}[Proof sketch] The projection map $\pi : \Fl(n) \to
\Fl(n; {\mathbf d})$ induces a ring homomorphism $\pi^* : A^*(\Fl(n;
{\mathbf d})) \to A^*(\Fl(n))$ sending the class of a variety $[Z] \in
A^*(\Fl(n; {\mathbf d}))$ to $[\pi^{-1}(Z)] \in A^*(\Fl(n))$. Just as in $\flags$, one can use Bruhat
order restricted to $S_n/S_{\mathbf d}$ to show that the Schubert
varieties $X_w({\mathbf d})$ for $w \in S_n/S_{\mathbf d}$ form an
affine paving of $\Fl(n; \mathbf d)$. Therefore $A^*(\Fl(n; {\mathbf
d}))$ is a free abelian group on $\{[X_w({\mathbf d})] \given w \in
S_n/S_{\mathbf d}\}$. This implies that $\pi^*$ is injective because
$\pi^*[X_w({\mathbf d})] = [\pi^{-1} X_w({\mathbf d})] = [X_w]$ for $w
\in S_n/S_{\mathbf d}$ by \Cref{prop:partial-flag}, and we already
know that all Schubert classes in $A^*(\flags)$ are linearly
independent from the introduction to the Chow ring in \Cref{sub:Monk}.
Therefore, $A^*(\Fl(n; {\mathbf d}))$ is isomorphic to the subring of
$A^*(\Fl(n))$ spanned by $\{[X_w({\mathbf d})] \given w \in S_n/S_{\mathbf
d}\}$.  Similarly, the affine paving also implies $A^*(\Fl(n; {\mathbf
d})) \simeq H^*(\Fl(n; {\mathbf d}))$.  Therefore, since the
cohomology class $[X_{w}]$ is represented by $\fS_{w_{0}w}$ in
$R_{n}$, the set $\{\fS_{w_{0}w} \given w \in S_n/S_{\mathbf d}\}$ in the
coinvariant algebra spans the subring in $R_n$ isomorphic to
$A^*(\Fl(n; {\mathbf d}))$.
\end{proof}

Along with Inherited Positivity \Cref{cor:Schub.basis},
\Cref{thm:partial-flag-cohom} says that Schubert polynomials for
$S_{n}$ encode all of Schubert calculus on every partial flag
variety, which is why we have focused on that more general case---but
this also makes Schubert problems for $\flags$ maximally difficult.
We will analyze the special cases of $1,2,3$-step flag varieties
where the Schubert structure constants can be interpreted as counting
certain combinatorial objects.

\subsubsection{Grassmannians}\label{sub:Grassmannians}

We now consider the case of Grassmannians more carefully.  In this
case, $\Gr (k,n)$ is the partial flag variety $\Fl(n;\mathbf{d})$ for
$\mathbf{d}=\{k \}$.  So Grassmannians are the \emph{1-step flag
varieties}.  The parabolic subgroup $S_{\{k\}}$ is generated by all
simple reflections except for $s_{k}$.  The permutations indexing
Schubert varieties and Schubert cells in $\Gr (k,n)$ are indexed by
the Bruhat-maximal elements in the cosets of $S_n / S_{\{k\}}$.  By
\Cref{thm:partial-flag-cohom}, the corresponding Schubert polynomials
are indexed by the Bruhat-minimal elements in their right $S_{\{k\}}$
cosets.  A permutation $w \in S_n$ is \emph{$k$-Grassmannian} if $w_i
> w_{i+1}$ implies $i = k$; these are the Bruhat-minimal
representatives in their cosets. Similarly, $w$ is
\emph{$k$-anti-Grassmannian} if $w_i < w_{i+1}$ implies $i = k$; these
are exactly the Bruhat-maximal elements in their cosets.  By the
bijection $w \mapsto w_{0}w$ we can toggle between the two sets.
Furthermore, sending $w \mapsto \{w_1, \ldots, w_k\}$ gives a
well-defined bijection from the cosets $S_n / S_{\{k\}}$ to the
collection ${[n] \choose k}$ of $k$-subsets of $[n]$, and from now on
we will usually index Schubert varieties in $\Gr(k,n)$ by $k$-subsets.
Therefore, if $w \in S_n $ and $I = \{w_{1},\dots , w_{k} \}$, then
the Schubert cell $C_{w}(E_{\bullet};\{k \})$ will be denoted by
$C_{I}(E_{\bullet})$ or just $C_{I}$ if the fixed flag is known from
context.  Schubert varieties $X_{I}$ will be defined similarly.

\begin{Example} \label{ex:gr-schubert}
Let $E_{\bullet}$ be the standard flag in $\C^5$. Then
$C_{\{2,3,5\}}(E_\bullet) \subseteq \Gr(3,5)$ is the set of $3$-planes
of the form
\begin{equation*}
    \operatorname{colspan} \begin{bmatrix}
    \ast & \ast & \ast \\ 
    \ast & \ast & 1\\
    \ast & 1 & 0\\
    \ast & 0 & 0 \\
    1 & 0 & 0
    \end{bmatrix}\,  = \, \operatorname{colspan} \begin{bmatrix}
    \ast & \ast & \ast \\ 
    0 & 0 & 1\\
    0 & 1 & 0\\
    \ast & 0 & 0 \\
    1 & 0 & 0
    \end{bmatrix}
 \, = \,
\operatorname{colspan} \begin{bmatrix}
    \ast & \ast & \ast \\ 
    1 & 0 & 0\\
    0 & 1 & 0\\
    0 & 0 & \ast \\
    0 & 0 & 1
    \end{bmatrix}
\end{equation*}
by Gauss elimination.  Recall from \Cref{prop:partial-flag} that
$C_{\{2,3,5 \}}$ is the image of every $C_{w}$ for $w\in S_{5}$
projected to $\Gr (3,5)$ such that $\{2,3,5 \}=\{w_{1},w_{2},w_{3}
\}$.  For example, the first matrix above looks like the first three
columns of the canonical form for matrices representing flags in the
Schubert cell $C_{53241} \subset \Fl_{5}$ and the last one looks like
the Schubert cell $C_{23514} \subset \Fl_{5}$.  Here $53241$ is a
$3$-anti-Grassmannian permutation, and $23514$ is the unique
$3$-Grassmannian permutation in the same $S_{\{3\}} \simeq S_{3}\times S_{2}$ coset.
By \Cref{thm:partial-flag-cohom}, the class $[X_{\{2,3,5 \}}]$ in
$H^*(\Gr (3,5))$ is represented by the Schubert polynomial indexed by
$w_{0}\, 53241=13425$, which is another $3$-Grassmannian
permutation. Observe that $\mathrm{dim}(C_{\{2,3,5 \}})=\ell(23514)=4$ and
$\ell(13425)=\codim(C_{\{2,3,5 \}})=2$.  
\end{Example}

The $k$-subsets indexing Grassmannian Schubert cells record the ``jump
set'' of the dimensions $\dim(V \cap E_i)$: the set $I = \{i \given
\dim(V \cap E_i) > \dim(V \cap E_{i-1})\}$ is the unique $I$ with $V
\in C_I(E_\bullet)$.  In \Cref{ex:gr-schubert}, the indexing set
$\{2,3,5 \}$ equals the row numbers with the pivots in the canonical
matrices, so we see the ``jumps'' must occur in these positions. Alternatively, we can describe Grassmannian Schubert cells and varieties as follows:
\begin{align} \label{eq:grassmannian-rank-conditions}
C_I(E_\bullet) = \{V \in \Gr(k,n) \given \dim(V \cap E_i) = |I \cap [i]| \text{ for all $i \in [n]$}\}\\
X_I(E_\bullet) = \{V \in \Gr(k,n) \given \dim(V \cap E_i) \geq |I \cap [i]| \text{ for all $i \in [n]$}\}. \nonumber 
\end{align}
Each Grassmannian Schubert variety decomposes into Grassmannian
Schubert cells as in the complete flag manifold.  From
\Cref{eq:grassmannian-rank-conditions}, we see
\begin{equation}
X_I(E_\bullet)  = \bigcup C_{H}(E_\bullet) 
\end{equation}
over all $H\trianglelefteq I$ in Gale order.  This is reminiscent of
the construction of the Schubert matroid in
\Cref{def:schubert-matroid-polytope}.   

There is also another common indexing of Grassmannian Schubert
varieties used in the literature. Consider the region in the plane
$[0,n-k] \times [0,k]$ as a $k \times (n-k)$ grid of unit
squares. Given a $k$-subset $I \subseteq [n]$, let $p$ be the lattice
path from $(0,0)$ in the lower left corner to $(n-k,k)$ in the upper
right corner which moves upward one unit on step $i$ if $i \in I$, and
rightward otherwise. For instance, here are the paths corresponding to the subsets $\{2,3,5\}$ and $\{1,3,5\}$ of $[5]$:
\begin{center}
\begin{tikzpicture}[scale=0.5]
\draw (0,0) -- (2,0);
\draw (0,1) -- (2,1);
\draw (0,2) -- (2,2);
\draw (0,3) -- (2,3);
\draw (0,0) -- (0,3);
\draw (1,0) -- (1,3);
\draw (2,0) -- (2,3);
\draw[very thick, red] (0,0) -- (1,0) -- (1,1) -- (1,2) -- (2,2) -- (2,3);
\end{tikzpicture} \qquad \raisebox{0.75cm}{\text{and}} \qquad 
\begin{tikzpicture}[scale=0.5]
    \draw (0,0) -- (2,0);
    \draw (0,1) -- (2,1);
    \draw (0,2) -- (2,2);
    \draw (0,3) -- (2,3);
    \draw (0,0) -- (0,3);
    \draw (1,0) -- (1,3);
    \draw (2,0) -- (2,3);
    \draw[very thick, red] (0,0) -- (0,1) -- (1,1) -- (1,2) -- (2,2) -- (2,3);
    \end{tikzpicture}
     \raisebox{0.75cm}{\text{.}} 
\end{center}
The set of squares weakly northwest of $p$ is the \emph{Young diagram}
of an integer partition $\lambda(I)=(\lambda_1(I) \geq \cdots \geq
\lambda_k(I) \geq 0)$, where $\lambda_i(I)$ is the number of squares
in the $i$\textsuperscript{th} row from the top. In the example above,
$\lambda(\{2,3,5 \}) = (2,1,1)$ and $\lambda(\{1,3,5 \}=(2,1,0)$. This sets up a bijection between ${[n] \choose
k}$ and the integer partitions whose Young diagram fits in a $k \times
(n-k)$ rectangle, and accordingly we will write $X_I$ or
$X_{\lambda(I)}$ as convenient.  

Write $\lambda=(\lambda_{1},\dots ,\lambda_{\ell}) \subseteq k \times
(n-k)$ to mean that (the Young diagram) of $\lambda$ fits in a $k
\times (n-k)$ rectangle, i.e. $\lambda_1 \leq n-k$ and $\lambda$ has
at most $k$ nonzero parts.  The \emph{size} of $\lambda$,
$|\lambda|=\lambda_{1}+\dots +\lambda_{\ell}$.  If $p$ is the path
defining $\lambda$, then observe that southeast of $p$ in the $k
\times (n-k)$ rectangle there is another partition shape rotated by
180\textdegree, which we will call the \textit{box complement}
$\lambda^{\vee}$.  Observe that $\lambda^{\vee}_i =
n-k-\lambda_{k-i+1}$.

\begin{Exercise}
Let $C_{I}$ be a Schubert variety in $\Gr (k,n)$.  Show that if a
general element of $C_I$ is written as the column span of a matrix $M$
in reduced column echelon form, like the third matrix in
Example~\ref{ex:gr-schubert}, then the pattern of $\ast$ entries in
$M^T$ is the Young diagram of $\lambda(I)$.  Hence, the dimension of
$X_{I}=X_{\lambda (I)}$ is  $|\lambda (I)|$.  Furthermore, the
codimension of $C_{I}$ is the size of the box complement
$\lambda^{\vee}(I)$. 
\end{Exercise}

The dimension and codimension of $C_{I}$ can also be found in terms of
permutations associated to $I$.  Furthermore, there is a close
connection between these permutations and the box complement
operation.  

\begin{Exercise}\label{exercise:dims.grass}
For a $k$-subset $I\subset [n]$, let $w$ be the $k$-Grassmannian
permutation such that $I=\{w_{1},\dots , w_{k} \}$.
\begin{enumerate}
\item Show $\mathrm{dim}(C_{I})=\ell(w)$.
\item If $I'=\{n-w_{1}+1,\dots , n-w_{k}+1 \}$ and
$w'$ is the  $k$-Grassmannian permutation $w'$ such that
$I'=\{w'_{1},\dots , w'_{k} \}$, then show $\codim(C_{I})=\ell(w')$.  
\item Prove the box complement of $\lambda (I)$ is $\lambda (I')$.
\end{enumerate}
\end{Exercise}

The Schubert polynomials representing Grassmannian Schubert classes
are particularly nice from an algebraic point of view.  These are the
cohomology classes for the Schubert varieties in Grassmannian
varieties by \Cref{thm:partial-flag-cohom}.

\begin{Proposition} \label{prop:symmetry}
If $w \in S_n$ is $k$-Grassmannian, then $\fS_w$ is a symmetric polynomial in $\bZ[x_1, \ldots, x_k]$, i.e. is invariant under all permutations of its variables $x_1, \ldots, x_k$.
\end{Proposition} 
\begin{proof}
By \Cref{lem:homogeneous} and \Cref{ex:divided-difference}, we know the Schubert polynomial $\fS_{w}$ is
symmetric in all of the variables $x_{1},x_{2},\dots , x_{k}$ and does
not depend on any other variables by \Cref{lem:homogeneous}.
\end{proof}

\begin{Definition}\label{def:shape.of.perm}
For any permutation $w \in S_{\infty}$, define an \textit{associated
partition shape} $\lambda (w)$ by sorting the code $c(w)$ into weakly
decreasing order.
\end{Definition}

\begin{Exercise}\label{ex:code.shape}
Show that if $w$ is $k$-Grassmannian, then $\lambda (w)$ equals
$\lambda(\{w_{1},\dots , w_{k} \})$.  Note, this is not true for
permutations which are not $k$-Grassmannian.
\end{Exercise}

\begin{Definition}\label{def:schur.poly}
Given an integer partition $\lambda=(\lambda_{1}\geq \dots \geq
\lambda_{k}\geq 0)$ with at most $k$ parts, the \emph{Schur
polynomial} $s_{\lambda}(x_1, \ldots, x_k)$ is $\fS_w$ where $w$ is
the unique $k$-Grassmannian permutation in $S_{\infty}$ with $\lambda
= \lambda(w)$.  
\end{Definition}

\begin{Example}\label{ex:322}
The permutation $w=346125$ is $3$-Grassmannian.  Its code is
$c(w)=(2,2,3,0,0,0)$ so its shape is $\lambda (w)=(3,2,2)$. Via the
transition equation or ladder moves on pipe dreams (\Cref{thm:chutes.and.ladders}),
one can verify that
\[
\fS_{w}=x_{1}^{2}x_{2}^{2}x_{3}^{2}(x_{1}+x_{2}+x_{3})=s_{(3,2,2)}(x_{1},x_{2},x_{3}).
\]
Similarly, $1\times w=1457236$ is $4$-Grassmannian and $\fS_{1\times
w}=s_{(3,2,2)}(x_{1},x_{2},x_{3},x_{4})$. Try drawing out the diagram
of these permutations and the corresponding bottom pipe dreams to
compare with the partition shape for $(3,2,2)$.  
\end{Example}

For readers who are familiar with the Schur polynomials as generating
functions for semistandard Young tableaux, you may be asking is there
a bijection from Grassmannian reduced pipe dreams to these tableaux?
Yes!  See \Cref{exercise:Grassmannian.Schur} below.

\begin{Exercise} \label{exer:schur-stability}
  Verify the stability property of Schur polynomials that
$s_{\lambda}(x_1, \ldots, x_k, 0) = s_{\lambda}(x_1, \ldots, x_k)$
using the transition equation for Schubert polynomials
\Cref{t:transitionA}.
\end{Exercise}

 \begin{Exercise} \label{exer:schur-variables} Show that if
$\lambda=(\lambda_{1}\geq \dots \geq \lambda_{\ell}>0)$ has exactly
$\ell$ nonzero parts, then $s_{\lambda}(x_1, \ldots, x_\ell)$ is
divisible by $x_1 x_2 \cdots x_{\ell}$. (Hint: use pipe dreams.)
\end{Exercise}

\begin{Example}
\ytableausetup{boxsize=4pt}
The nonidentity 1- and 2-Grassmannian permutations in $S_4$ and corresponding Schur polynomials are 
\begin{equation*}
\begin{array}{ll}
w = 2134 & s_{\ydiagram{1}}(x_1) = x_1\\
w = 3124 & s_{\ydiagram{2}}(x_1) = x_1^2 \\
w = 4123 & s_{\ydiagram{3}}(x_1) = x_1^3
\end{array} \qquad \qquad  \begin{array}{ll}
w = 1324 & s_{\ydiagram{1}}(x_1,x_2) = x_1 + x_2\\
w = 1423 & s_{\ydiagram{2}}(x_1,x_2) = x_1^2 + x_1 x_2 + x_2^2 \\
w = 2314 & s_{\ydiagram{1,1}}(x_1,x_2) = x_1 x_2\\ 
w = 2413 & s_{\ydiagram{2,1}}(x_1,x_2) = x_1^2 x_2 + x_1 x_2^2\\
w = 3412 & s_{\ydiagram{2,2}}(x_1,x_2) = x_1^2 x_2^2
\end{array}
\end{equation*}
\end{Example}

\begin{Exercise} \label{exer:schur-basis}
Using the fact that the leading term of $\fS_w$ in lex order is
$x_1^{c(w)_1} \cdots x_n^{c(w)_n}$ (by \Cref{lem:homogeneous} or
\Cref{thm:chutes.and.ladders}), show that the Schur polynomials
$s_{\lambda}(x_1, \ldots, x_k)$ are a $\bZ$-basis for the ring of
symmetric polynomials in $\bZ[x_1, \ldots, x_k]$ as $\lambda$ runs
over integer partitions of length $\leq k$.
\end{Exercise}

We can now describe the cohomology ring of $\Gr(k,n)$ in terms of
Schur polynomials. Let $\bZ[x_1, \ldots, x_k]^{S_k}$ be the subring of
\textit{symmetric polynomials} in $\bZ[x_1, \ldots, x_k]$.  A basis
for this ring is given by Schur polynomials by indexed by partitions
$\lambda \subseteq k \times (n-k)$ by  Exercise~\ref{exer:schur-basis}.

\begin{Theorem} \label{thm:grassmannian-cohom} For positive integers
$k\leq n$,   the linear map
\begin{equation*} H^*(\Gr(k,n)) \to \bZ[x_1, \ldots, x_k]^{S_k} /
(s_{\lambda}(x_1, \ldots, x_k) : \lambda \not \subseteq k \times
(n-k))
\end{equation*}
sending the Schubert class $[X_\lambda] \in H^*(\Gr(k,n))$ to the
Schur polynomial $s_{\lambda^{\vee}}(x_1, \ldots, x_k)$ is a ring
isomorphism.
\end{Theorem}

\begin{proof}
Let $\mathcal{S}_{k,n}=\bZ[x_1, \ldots, x_k]^{S_k} / (s_{\lambda}(x_1,
\ldots, x_k) : \lambda \not \subseteq k \times (n-k))$. By
\Cref{cor:Schub.basis}, there is a ring homomorphism $f : \bZ[x_1,
x_2, \ldots] \to H^*(\flags)$ with $f(\fS_w) = [X_{w_0 w}]$ if $w \in
S_n$ and $f(\fS_w) = 0$ if $w \in S_\infty \setminus S_n$. If $\lambda_1 > n-k$ and $\lambda$ has length
$\leq k$, then $s_{\lambda}(x_1, \ldots x_k) = \fS_w$ for some
$k$-Grassmannian permutation $w \in S_\infty \setminus S_n$, so
$f(s_{\lambda}(x_1, \ldots x_k)) = 0$.  If
$\lambda$ has more than $k$ nonzero parts, then $s_{\lambda}(x_1,
\ldots, x_k) = 0$ by Exercises~\ref{exer:schur-stability} and~\ref{exer:schur-variables}. The restriction of $f$ to $\bZ[x_1,
\ldots, x_k]^{S_k}$ therefore descends to a ring homomorphism
\[
j : \mathcal{S}_{k,n} \to H^*(\flags).
\]

For each $I \in {[n] \choose k}$, we have $j(\fS_w) = [X_{w_0 w}]$ where $w$ is the $k$-Grassmannian
permutation with $I = \{w_1, \ldots, w_k\}$. By
Theorem~\ref{thm:partial-flag-cohom}, the image of $j$ can be
identified with $H^*(\Gr(k,n))$, in which case our convention for
Grassmannian Schubert classes then identifies $[X_{w_0 w}]$ with the
Grassmannian Schubert class $[X_{\lambda(I')}]$ where $I' = \{n-w_1+1,
\ldots, n-w_k+1\}$. By definition, $\fS_w = s_{\lambda(I)}(x_1, \ldots, x_k)$,
 while $\lambda(I') = \lambda(I)^\vee$ by Exercise~\ref{exercise:dims.grass}.
 Since $I \mapsto \lambda(I)$ is a bijection between $k$-subsets and partitions
  $\lambda \subseteq k \times (n-k)$, this gives a  ring homomorphism
  $j' :\mathcal{S}_{k,n} \to H^*(\Gr(k,n))$ sending $s_{\lambda}(x_1,
\ldots, x_k)$ to $[X_{\lambda^\vee}]$. The map $j'$ sends the
spanning set $\{s_{\lambda}(x_1,\ldots,k) : \lambda \subseteq k \times
(n-k)\}$ of $\mathcal{S}_{k,n}$ to the Schubert basis of
$H^*(\Gr(k,n))$, so that spanning set is actually a basis and $j'$ is an
isomorphism.
\end{proof}

\begin{Exercise}\label{ex:grassmannian.positivity}
Review the statements in the Inherited Positivity \Cref{cor:Schub.basis}.   How do each of
these statements descend to the Grassmannian varieties?  
\end{Exercise}

We end this subsection by discussing Schubert structure constants for
the Grassmannian varieties. Recall from
\Cref{rem:positivity} that all Schubert problems can be solved by
computing the Schubert structure constants.  By the Inherited
Positivity \Cref{cor:Schub.basis}, these constants can be computed by
multiplying Schubert polynomials and expanding the product into the
basis of Schubert polynomials as a positive integral sum.  In special
cases, these constants are known to also count certain combinatorial
objects, leading to a deeper understanding of the related processes.
One well-known special case worth reviewing is the expansion
coefficients for the product of Schubert polynomials indexed by
$k$-Grassmannian permutations, which correspond to the product of two
Schur polynomials.

\begin{Definition}\label{def:LR.coeffs}
    The \emph{Littlewood-Richardson} coefficients $c_{\lambda \mu}^{\nu}$ are defined by the expression 
    \begin{equation*}
    s_{\lambda}(X_k) s_{\mu}(X_k)  = \sum_{\nu} c_{\lambda \mu}^{\nu} s_{\nu}(X_k)
    \end{equation*}
    where $\nu$ ranges over partitions with length at most $k$, and we
    have abbreviated the variable list $x_1, \ldots, x_k$ as $X_k$.
    Via the bijection from $k$-subsets to integer partitions that fit
    in the $k \times (n-k)$ rectangle, the Littlewood-Richardson
    coefficients may also be written as $c_{I,J}^{K}$.  
\end{Definition}

For example, $s_{\ydiagram{1}}(X_2)^2 = s_{\ydiagram{2}}(X_2) +
s_{\ydiagram{1,1}}(X_2)$, which shows that
$c_{\ydiagram{1},\ydiagram{1}}^{\ydiagram{2}} =
c_{\ydiagram{1},\ydiagram{1}}^{\ydiagram{1,1}} = 1$. We could also get
$c_{\ydiagram{1},\ydiagram{1}}^{\ydiagram{2}} = 1$ from the simpler
expression $s_{\ydiagram{1}}(x_1)^2 = x_1^2 = s_{\ydiagram{2}}(x_1)$,
but one should not make the mistake of concluding that
$c_{\ydiagram{1},\ydiagram{1}}^{\ydiagram{1,1}}$ is zero from this: to
compute $c_{\lambda \mu}^{\nu}$ correctly we need to work in at least
as many variables as the lengths of $\lambda$, $\mu$, and $\nu$. A
common fix for this annoyance is to define the \emph{Schur function}
\begin{equation}\label{eq:Schur.function.def}
s_{\lambda} = \lim_{k \to \infty} s_{\lambda}(x_1, \ldots, x_k),
\end{equation}
a symmetric formal power series in $x_1, x_2, \ldots$, in which case
we have $s_\lambda s_\mu = \sum_\nu c_{\lambda \mu}^{\nu} s_{\nu}$
with the sum ranging over \emph{all} partitions $\nu$ of size $|\lambda|+|\mu|$.

Schur functions and Littlewood-Richardson coefficients appear in a
remarkable variety of different contexts. We will name a few of them.

\begin{itemize}
    \item It follows from the discussion of structure constants for
Schubert varieties, Geometry Implies Positivity
\Cref{cor:structure.constants}, and the definition of Schur functions
given above that the coefficient $c_{\lambda
\mu}^{\nu}$ is the number of points in the intersection $X_{\lambda^{\vee}}(E_\bullet) \cap X_{\mu^{\vee}}(F_\bullet) \cap X_{\nu}(G_\bullet)$ with respect to generic flags $E_\bullet, F_\bullet, G_\bullet$.
\item The Schur polynomials $s_{\lambda}(X_k)$ are the characters of
irreducible polynomial representations of $\GL_k(\C)$. Upon
decomposing the tensor product of two such representations into
irreducibles, the multiplicities are the Littlewood-Richardson
coefficients \cite[Part II]{Fulton-book}.
\item As described below in
\S\ref{subsec:eigenvalues}, the set of possible spectra of triples of
Hermitian matrices $(A,B,A+B)$ is a convex polytope whose facets
correspond to those $c_{\lambda \mu}^{\nu}$ equal to $1$.
\item A
nilpotent linear operator $T : V \to V$ on a finite-dimensional vector space
over an algebraically closed field is characterized up to conjugacy by
its Jordan type $\lambda(T) = (\lambda_1 \geq \cdots \geq \lambda_k >
0)$, the list of block sizes in the Jordan normal form of $T$. A short
exact sequence of such operators $T_1, T_2, T_3$ is a commuting
diagram \begin{equation*} \begin{tikzcd}
        0 \arrow[r] & V_1 \arrow[r, "f"] \arrow[d, "T_1"] & V_2 \arrow[r, "g"] \arrow[d, "T_2"] & V_3 \arrow[r] \arrow[d, "T_3"] & 0\\
        0 \arrow[r] & V_1 \arrow[r, "f"] & V_2 \arrow[r, "g"] & V_3 \arrow[r]  & 0
         \end{tikzcd}
    \end{equation*}
     where $f$ and $g$ are linear maps and the rows are exact. Fixing
     $T_1, T_2, T_3$, there exists such a sequence if and only if
     $c_{\lambda(T_1), \lambda(T_3)}^{\lambda(T_2)} > 0$, see \cite{klein-nilpotent}.
\end{itemize}

There are many different \emph{Littlewood-Richardson rules} known,
i.e. families of combinatorial objects which are enumerated by
$c_{\lambda \mu}^{\nu}$. We describe the \emph{puzzle rule} of Knutson
and Tao \cite{KnutsonTao-puzzles}.  Here the natural encoding for the
Littlewood-Richardson coefficients coerresponding with the structure
constants in $\Gr(k,n)$ is in terms of $k$-subsets of $[n]$ encoded as
binary sequences instead of partitions.  Given a $k$-subset $I
\subseteq [n]$, let $b(I)$ be the binary word of length $n$ with
$b(I)_i = 0$ if and only if $i \in I$.

\begin{Definition}\label{def:puzzle}
A \emph{puzzle} is any assembly of the tiles
\begin{center}
    \begin{tikzpicture}
    \draw (-0.5,0) -- node {$1$} (0.5,0) -- node {$1$} (0,0.866) -- node {$1$} (-0.5,0);
    \end{tikzpicture}     \qquad 
    \begin{tikzpicture}
    \draw (-0.5,0) -- node {$0$} (0.5,0) -- node {$0$} (0,0.866) -- node {$0$} (-0.5,0);
    \end{tikzpicture}   \qquad 
    \begin{tikzpicture}
    \draw (-0.5,0) -- node {$0$} (0.5,0) -- node {$1$} (0,0.866) -- node {$10$} (-0.5,0);
    \end{tikzpicture} \raisebox{5.5mm}{and their rotations}
\end{center}
into an equilateral triangle (oriented the same way), where the labels along any edge shared by two tiles must match.
\end{Definition}

\begin{Theorem}\cite{KnutsonTao-puzzles}\label{thm:KnutsonTao-puzzles}
Let $I, J, K \subseteq [n]$ be $k$-subsets. Then $c_{I,J}^{K}=c_{\lambda(I), \lambda(J)}^{\lambda(K)}$ is the number of puzzles of the form 
\begin{center}
     \begin{tikzpicture}
    \draw (-0.5,0) -- node[below] {$b(K)$} (0.5,0) -- node[right] {$b(J)$} (0,0.866) -- node[left] {$b(I)$} (-0.5,0);
    \end{tikzpicture}
\end{center}
where the word formed by the labels along of the three edges is read left to right.
\end{Theorem}

\begin{Example} We compute the Schubert class expansion of $[X_{\ydiagram{2,1}}]^2$ in $H^*(\Gr(2,4))$ using the puzzle rule. Since $[X_{\ydiagram{2,1}}]$ corresponds to the Schur function $s_{\ydiagram{1}}$ by Theorem~\ref{thm:partial-flag-cohom}, we must compute all coefficients $c_{\ydiagram{1}, \ydiagram{1}}^{\nu}$ where the Young diagram of $\nu$ fits in a $2 \times 2$ box. The partition $(1)$ in the $2 \times 2$ box corresponds to the subset $\{1,3\} \subseteq [4]$, so the edge labels of the two non-horizontal sides are $b(\{1,3\}) = 0101$. This already forces part of a puzzle:
\begin{center}
\begin{tikzpicture} 
\draw (-2, 0) -- node {$0$} (-1.5, 0.866) -- node {$1$} (-1, 0.866*2) -- node {$0$} (-0.5, 0.866*3) -- node {$1$} (0, 0.866*4);
\draw (0,0.866*4) -- node {$0$} (0.5,0.866*3) -- node {$1$} (1,0.866*2) -- node {$0$} (1.5,0.866) -- node {$1$} (2,0);
\draw (-0.5, 0.866*3) -- node {$10$} (0.5, 0.866*3);    
\draw (-1, 0.866*2) -- node {$0$} (0, 0.866*2) -- node {$1$} (1, 0.866*2);
\draw (-0.5, 0.866*3) -- node {$0$} (0, 0.866*2) -- node {$1$} (0.5, 0.866*3);
\draw (-1.25, 0.866*1.5) node[left] {$\rightarrow$};
\draw (-2, 0) -- (2, 0);
\end{tikzpicture}
\end{center}
At this point there are two possible choices for the tile incident to the $1$ marked with $\rightarrow$, but making a choice uniquely determines the rest of the puzzle, giving the two possibilities
\begin{center}
\begin{tikzpicture}
\draw (-2, 0) -- node {$0$} (-1.5, 0.866) -- node {$1$} (-1, 0.866*2) -- node {$0$} (-0.5, 0.866*3) -- node {$1$} (0, 0.866*4);
\draw (0,0.866*4) -- node {$0$} (0.5,0.866*3) -- node {$1$} (1,0.866*2) -- node {$0$} (1.5,0.866) -- node {$1$} (2,0);
\draw (-0.5, 0.866*3) -- node {$10$} (0.5, 0.866*3);
\draw (-1, 0.866*2) -- node {$0$} (0, 0.866*2) -- node {$1$} (1, 0.866*2);
\draw (-0.5, 0.866*3) -- node {$0$} (0, 0.866*2) -- node {$1$} (0.5, 0.866*3);

\draw (-1.5, 0.866*1) -- node {$1$} (-0.5, 0.866*1) -- node {$0$} (0.5,0.866) -- node {$10$} (1.5,0.866);
\draw (-1, 0.866*2) -- node {$1$} (-0.5, 0.866) -- node {$10$} (0, 0.866*2) -- node {$1$} (0.5, 0.866) -- node {$1$} (1, 0.866*2);
\draw (1.5,0.866) -- node {$1$} (1,0) -- node {$0$} (0.5,0.866) -- node {$0$} (0,0) -- node {$0$} (-0.5,0.866) -- node {$0$} (-1,0) -- node {$10$} (-1.5,0.866);
\draw (-2, 0) -- node {$1$} (-1, 0) -- node {$0$} (0, 0) -- node {$0$} (1, 0)  -- node {$1$} (2, 0);
\end{tikzpicture} \qquad \qquad 
\begin{tikzpicture}
\draw (-2, 0) -- node {$0$} (-1.5, 0.866) -- node {$1$} (-1, 0.866*2) -- node {$0$} (-0.5, 0.866*3) -- node {$1$} (0, 0.866*4);
\draw (0,0.866*4) -- node {$0$} (0.5,0.866*3) -- node {$1$} (1,0.866*2) -- node {$0$} (1.5,0.866) -- node {$1$} (2,0);
\draw (-0.5, 0.866*3) -- node {$10$} (0.5, 0.866*3);
\draw (-1, 0.866*2) -- node {$0$} (0, 0.866*2) -- node {$1$} (1, 0.866*2);
\draw (-0.5, 0.866*3) -- node {$0$} (0, 0.866*2) -- node {$1$} (0.5, 0.866*3);
\draw (-1.5, 0.866*1) -- node {$10$} (-0.5, 0.866) -- node {$1$} (0.5, 0.866) --  node {$0$} (1.5, 0.866);
\draw (-1, 0.866*2) -- node {$0$} (-0.5, 0.866) -- node {$0$} (0, 0.866*2) -- node {$10$} (0.5, 0.866) -- node {$0$} (1, 0.866*2);
\draw (-1.5, 0.866) -- node {$0$} (-1, 0) -- node {$1$} (-0.5, 0.866) -- node {$1$} (0, 0) -- node {$1$} (0.5, 0.866) -- node {$1$} (1, 0) -- node {$10$} (1.5, 0.866);
\draw (-2,0) -- node {$0$} (-1,0) -- node {$1$} (0,0) -- node {$1$} (1,0) -- node {$0$} (2,0);
\end{tikzpicture}
\end{center}
Reading off the words on the bottom boundary, we see that $
c_{\ydiagram{1},\ydiagram{1}}^{\ydiagram{1,1}} =
c_{\ydiagram{1},\ydiagram{1}}^{\ydiagram{2}} = 1$.
\end{Example}

\bigskip

\subsubsection{2,3-step flag varieties}\label{subsub:2-3} 
\begin{figure}[h]
\begin{center}
\begin{tikzpicture}
\draw (-0.5,0) -- node {$0$} (0.5,0) -- node {$0$} (0,0.866) -- node {$0$} (-0.5,0);
\end{tikzpicture} \quad
\begin{tikzpicture}
\draw (-0.5,0) -- node {$1$} (0.5,0) -- node {$1$} (0,0.866) -- node {$1$} (-0.5,0);
\end{tikzpicture} \quad
\begin{tikzpicture}
\draw (-0.5,0) -- node {$2$} (0.5,0) -- node {$2$} (0,0.866) -- node {$2$} (-0.5,0);
\end{tikzpicture} \quad
\begin{tikzpicture}
\draw (-0.5,0) -- node {$3$} (0.5,0) -- node {$0$} (0,0.866) -- node {$1$} (-0.5,0);
\end{tikzpicture} \quad
\begin{tikzpicture}
\draw (-0.5,0) -- node {$4$} (0.5,0) -- node {$1$} (0,0.866) -- node {$2$} (-0.5,0);
\end{tikzpicture} \quad
\begin{tikzpicture}
\draw (-0.5,0) -- node {$5$} (0.5,0) -- node {$0$} (0,0.866) -- node {$2$} (-0.5,0);
\end{tikzpicture} \quad
\begin{tikzpicture}
\draw (-0.5,0) -- node {$6$} (0.5,0) -- node {$3$} (0,0.866) -- node {$2$} (-0.5,0);
\end{tikzpicture} \quad
\begin{tikzpicture}
\draw (-0.5,0) -- node {$7$} (0.5,0) -- node {$0$} (0,0.866) -- node {$4$} (-0.5,0);
\end{tikzpicture}
\end{center}
\caption{The seven labeled triangular puzzle pieces used in the
Buch-Kresch-Purbhoo-Tamvakis proof \cite{Buch-Kresch-Purbhoo.2-step-rule} of the 2-step puzzle rule for
Schubert calculus.}  \label{fig:2-step.pieces}
\end{figure}
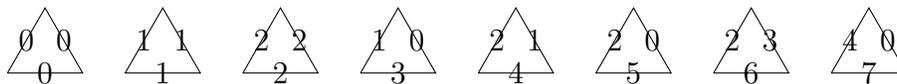

Can the
other Schubert structure constants be interpreted as counting
variations on puzzles?  For some cases, the answer is known to be
``yes''.  Recall, we think of the Grassmannian varieties as 1-step
flag varieties since $\Gr(k,n)= \Fl(n; \{k\})$.  The \textit{2-step}
flag variety is the partial flag variety $\Fl(n; {\bf d})$ with ${\bf
d} = \{d_1 < d_2\} \subseteq [n-1]$.  Knutson conjectured a puzzle
rule for the 2-step flag variety that uses more pieces with a larger
collection of labels on the triangles than we see in
\Cref{def:puzzle}.  Knutson's conjecture was proved by
Buch-Kresch-Purbhoo-Tamvakis in
\cite{Buch-Kresch-Purbhoo.2-step-rule}, which extends to equivariant
cohomology \cite{Buch.2015}.  In the 2-step puzzle rule, there are 7
puzzle pieces that can be used to fill the larger triangle and again
they must align on the boundaries. See \Cref{fig:2-step.pieces}.

\begin{figure}
\begin{tikzpicture}
\draw (-2, 0) -- node {$2$} (-1.5, 0.866) -- node {$0$} (-1, 0.866*2) -- node {$1$} (-0.5, 0.866*3) -- node {$3$} (0, 0.866*4);
\draw (0,0.866*4) -- node {$12$} (0.5,0.866*3) -- node {$0$} (1,0.866*2) -- node {$3$} (1.5,0.866) -- node {$12$} (2,0);
\draw (-0.5, 0.866*3) -- node {$G$} (0.5, 0.866*3);
\draw (-1, 0.866*2) -- node {$1$} (0, 0.866*2) -- node {$C$} (1, 0.866*2);
\draw (-0.5, 0.866*3) -- node {$12$} (0, 0.866*2) -- node {$3$} (0.5, 0.866*3);

\draw (-1.5, 0.866*1) -- node {$1$} (-0.5, 0.866*1) -- node {$0$} (0.5,0.866) -- node {$3$} (1.5,0.866);
\draw (-1, 0.866*2) -- node {$A$} (-0.5, 0.866) -- node {$0$} (0, 0.866*2) -- node {$0$} (0.5, 0.866) -- node {$3$} (1, 0.866*2);
\draw (1.5,0.866) -- node {$D$} (1,0) -- node {$R$} (0.5,0.866) -- node {$D$} (0,0) -- node {$12$} (-0.5,0.866) -- node {$1$} (-1,0) -- node {$12$} (-1.5,0.866);
\draw (-2, 0) -- node {$2$} (-1, 0) -- node {$1$} (0, 0) -- node {$3$} (1, 0)  -- node {$0$} (2, 0);
\end{tikzpicture} \qquad
\begin{tikzpicture}
\draw (-2, 0) -- node {$2$} (-1.5, 0.866) -- node {$0$} (-1, 0.866*2) -- node {$1$} (-0.5, 0.866*3) -- node {$3$} (0, 0.866*4);
\draw (0,0.866*4) -- node {$12$} (0.5,0.866*3) -- node {$0$} (1,0.866*2) -- node {$3$} (1.5,0.866) -- node {$12$} (2,0);
\draw (-0.5, 0.866*3) -- node {$G$} (0.5, 0.866*3);
\draw (-1, 0.866*2) -- node {$1$} (0, 0.866*2) -- node {$C$} (1, 0.866*2);
\draw (-0.5, 0.866*3) -- node {$12$} (0, 0.866*2) -- node {$3$} (0.5, 0.866*3);

\draw (-1.5, 0.866*1) -- node {$1$} (-0.5, 0.866*1) -- node {$0$} (0.5,0.866) -- node {$3$} (1.5,0.866);
\draw (-1, 0.866*2) -- node {$A$} (-0.5, 0.866) -- node {$0$} (0, 0.866*2) -- node {$0$} (0.5, 0.866) -- node {$3$} (1, 0.866*2);
\draw (1.5,0.866) -- node {$1$} (1,0) -- node {$E$} (0.5,0.866) -- node {$H$} (0,0) -- node {$E$} (-0.5,0.866) -- node {$1$} (-1,0) -- node {$12$} (-1.5,0.866);
\draw (-2, 0) -- node {$2$} (-1, 0) -- node {$3$} (0, 0) -- node {$0$} (1, 0)  -- node {$1$} (2, 0);
\end{tikzpicture} \qquad
\begin{tikzpicture}
\draw (-2, 0) -- node {$2$} (-1.5, 0.866) -- node {$0$} (-1, 0.866*2) -- node {$1$} (-0.5, 0.866*3) -- node {$3$} (0, 0.866*4);
\draw (0,0.866*4) -- node {$12$} (0.5,0.866*3) -- node {$0$} (1,0.866*2) -- node {$3$} (1.5,0.866) -- node {$12$} (2,0);
\draw (-0.5, 0.866*3) -- node {$G$} (0.5, 0.866*3);
\draw (-1, 0.866*2) -- node {$1$} (0, 0.866*2) -- node {$C$} (1, 0.866*2);
\draw (-0.5, 0.866*3) -- node {$12$} (0, 0.866*2) -- node {$3$} (0.5, 0.866*3);

\draw (-1.5, 0.866*1) -- node {$1$} (-0.5, 0.866*1) -- node {$0$} (0.5,0.866) -- node {$3$} (1.5,0.866);
\draw (-1, 0.866*2) -- node {$A$} (-0.5, 0.866) -- node {$0$} (0, 0.866*2) -- node {$0$} (0.5, 0.866) -- node {$3$} (1, 0.866*2);
\draw (1.5,0.866) -- node {$2$} (1,0) -- node {$F$} (0.5,0.866) -- node {$J$} (0,0) -- node {$F$} (-0.5,0.866) -- node {$L$} (-1,0) -- node {$F$} (-1.5,0.866);
\draw (-2, 0) -- node {$3$} (-1, 0) -- node {$1$} (0, 0) -- node {$0$} (1, 0)  -- node {$2$} (2, 0);
\end{tikzpicture}
\caption{All three completed puzzles used to compute $c_{u,v}^{w}$ using
the 3-step puzzle rule for Schubert calculus \cite{Buch.talk.2017} where $u=2314$, $v=2143$ with the resulting permutation $w=4213$, $3412$ and $3241$ from left to right.}  \label{fig:3-step}
\end{figure}
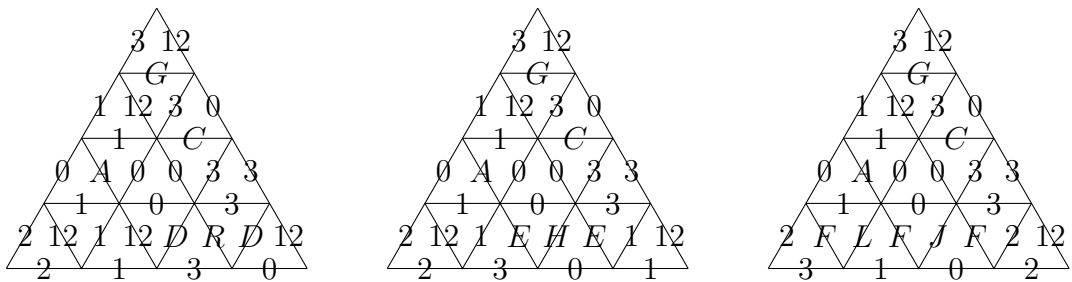

There is also a 3-step puzzle rule conjectured by Knutson-Buch and
proved by Knutson and Zinn-Justin
\cite{knutsonZin-Justin2020schubert}.  See also the special cases
proved earlier by Buch \cite{Buch.talk.2017}.  The 3-step rule uses
approximately 30 tiles with edges labeled by letters, numbers, or
pairs of numbers that have to match on their boundaries.  Some tiles
are triangles and some can be combined into diamonds. See
\Cref{fig:3-step} for an example. Here, we obtain a permutation from a
string by reading the positions of the $0$'s, then the positions of
the $1$'s and so on. For example, a string $1301220$ corresponds to
the permutation $3714562$ as $0$'s appear in positions $3$ and $7$,
$1$'s appear in positions $1$ and $4$, $2$'s appear in positions $5$
and $6$, and the unique $3$ appears in position $2$.

Alas, the consensus is there will never be an all-positive formula for
$d$-step flags for $d\geq 4$ using only a finite number of puzzle
pieces \cite{KnutsonICM2022,knutsonZin-Justin2020schubert}.  The
representation theory involved in the 1,2,3-step flags follows the
patterns for simply laced Dynkin diagrams.  However, for $d=4$, this
approach leads to some negative puzzle pieces and for $d\geq 5$
infinite-dimensional representations always occur.  For more details,
we recommend Allen Knutson's well-performed and insightful
video on the history and applications of puzzles in Schubert calculus
\cite{Knutson2022video}.

\subsection{Stanley Symmetric Functions}\label{sub:StanleySymmetrics}

Symmetric functions are power series in the formal power series ring
$\mathbb{Z}[[x_{1},x_{2},\ldots]]$ that are invariant under
interchanging any two variables and have finite degree.  The symmetric
functions form a subring of the power series ring, generally referred
to as $\mathrm{SYM}$.  Symmetric functions routinely appear in many
aspects of mathematics and theoretical physics, and they even have 
importance in quantum computation.  The symmetric functions are
contained in a bigger subring of the formal power series ring called
$\mathrm{QSYM}$ or the \textit{ring of quasisymmetric functions}. QSYM
is the subring of formal power series of finite degree $f(x_{1},x_{2},\dots )$,
which are invariant under a \textit{shift of indices: }
\begin{quote}
for every sequence of positive integers $a_1, a_2, \ldots, a_k$, the
coefficient of $x_{i_1}^{a_1} x_{i_2}^{a_2} \cdots x_{i_k}^{a_k}$ in
$f$ equals the coefficient of $x_{j_1}^{a_1} x_{j_2}^{a_2} \cdots
x_{j_k}^{a_k}$ in $f$ whenever $i_1 < i_2 < \cdots < i_k$ and $j_1 <
j_2 < \cdots < j_k$.
\end{quote}
For comparison, $f(x_{1},x_{2},\dots )$ is a symmetric function
provided the coefficient of $x_{i_1}^{a_1} x_{i_2}^{a_2} \cdots
x_{i_k}^{a_k}$ equals the coefficient of $x_{j_1}^{a_1} x_{j_2}^{a_2}
\cdots x_{j_k}^{a_k}$ in $f$ whenever $(i_1, i_2,\ldots, i_k)$ is any
permutation of $j_1 < j_2 < \cdots < j_k$.  Both $\mathrm{SYM}$ and $\mathrm{QSYM}$
have a wide range of applications related to Schubert calculus and
beyond \cite{Billey-McNamara,Gessel,S2}.

The constants in $\mathbb{Z}[[x_{1},x_{2},\ldots]]$ are of course both
symmetric and quasisymmetric. The rings $\mathrm{SYM}$ and $\mathrm{QSYM}$ also agree in
degrees 1 and 2, but for higher degrees the homogeneous symmetric
functions are a proper subset of the homogeneous quasisymmetric
functions of the same degree.  One important family of examples is the
\textit{fundamental quasisymmetric functions} given by 
\begin{equation}\label{eq:fundamentals}
F_{S}^{(d)}= \sum x_{i_{1} }x_{i_{2} }\cdots x_{i_{d} }
\end{equation}
summed over all weakly increasing sequences of positive integers $
(i_{1}\leq i_{2}\leq \cdots \leq i_{d})$ such that $i_{k}<i_{k+1} $
for all $k \in S \subset [d-1]$.  For $d=0$ and $S=\emptyset$, define
$F_{\emptyset}^{(0)}=1$.  The set $S$ is the \textit{required jump
sequence}, but additional jumps may occur or not.  For $d=3$ and
$S=\{1 \}$, we have $F_{S}^{(3)}= x_{1}x_{2}^{2}+x_{1}x_{2}x_{3} +
x_{1}x_{3}^{2}+x_{1}x_{2}x_{4}+ \cdots$, where the expansion has two
types of exponents $(1,2)$ and $(1,1,1)$, but not $(2,1)$, so
$F_{\{1 \}}^{(3)}$ is not in $\mathrm{SYM}$.

The homogeneous symmetric functions of degree $d>0$ have a natural
basis of monomial symmetric functions $m_{\lambda}$ indexed by
partitions $\lambda =(\lambda_{1}\geq \lambda_{2}\geq \dots \geq
\lambda_{k})$ such that each $\lambda_{i}$ is a positive integer and
$\sum \lambda_{i}=d$. Here $m_{\lambda}$ is the sum over all monomials
of the form $x_{i_1}^{\lambda_1} x_{i_2}^{\lambda_2} \cdots
x_{i_k}^{\lambda_k}$ with $i_{1},i_{2},\dots, i_{k}$ distinct positive
integers.  Similarly, the homogeneous quasisymmetric functions of
degree $d$ have a natural basis of monomial quasisymmetric functions
indexed by compositions $\alpha =(\alpha _{1}, \alpha
_{2},\dots,\alpha_{k})$ such that each $\alpha_{i}$ is a positive
integer and $\sum \alpha_{i}=d$ given by

\begin{equation}\label{eq:monomial.qsym}
M_{\alpha } = \sum_{i_{1}<i_{2}<\dots < i_{k}} x_{i_1}^{\alpha_1}
x_{i_2}^{\alpha_2} \cdots x_{i_k}^{\alpha_k}.
\end{equation}

As mentioned in the introduction, we expect the readers to have some
familiarity with symmetric function theory and the classical bases of
$\mathrm{SYM}$ such as monomial symmetric functions, elementary
symmetric functions $e_{\lambda}$ and Schur functions $s_{\lambda}$.
In particular, recall from \Cref{sub:123.StepFlags} that the Schur
function $s_{\lambda}$ equals the limit $\lim_{k \to \infty} s_{\lambda}(x_1, \ldots,
x_k)$ for a partition $\lambda$.  We don't expect any advanced
knowledge of $\mathrm{QSYM}$, but we recommend the reader try
expanding $m_\lambda$ in terms of the $M_{\alpha}'s$ in some cases for
practice.

The set of fundamental quasisymmetric functions $\{F_{S}^{(d)} \given
S \subset [d-1] \}$ forms another basis of the homogeneous elements of
$\mathrm{QSYM}$ of degree $d$.  Gessel first developed this basis and
named the quasisymmetric functions in his influential paper on
$P$-partitions \cite{Gessel}.  In that paper, he showed that the Schur
functions $s_{\lambda}$ have an elegant expansion into the basis of
fundamental quasisymmetric functions.  To state his theorem, we need some notation.  Let $\lambda=(\lambda_{1}\geq \dots \geq
\lambda_{k})$ be a partition of $d$ with its Young diagram drawn in
English notation with $\lambda_{1}$ boxes on the top row and
$\lambda_{k}$ boxes on the bottom, all left justified.  Recall,
$\mathrm{SYT}(\lambda)$ is the set of standard Young tableaux of shape
$\lambda$, which are bijective fillings of the boxes of $\lambda$ by
values $1,2,\dotsc , d$ increasing along rows and down the columns (see \Cref{sub:Balanced}).
The \textit{reading word} of $T \in \mathrm{SYT}(\lambda)$ is the permutation
obtained by reading the values in $T$ from left to right along rows
going bottom to top.   The \textit{ascent set} of $T$, denoted $A(T)$,
is the set of indices $i$ such that $i$ appears before $i+1$ in the
reading word of $T$.  The ascent set of $T$ is also the ascent set of
the inverse permutation for the reading word. For example,

\begin{equation}\label{eq:sty.ex}
T=\tableau{1 & 3 & 7\\
2 & 4 & 8\\
5 & 6
}
\end{equation}
\bigskip

\noindent is a standard Young tableau of shape $\lambda =(3,3,2)$ with
ascent set $A(T)=\{2,5,6 \}$ since the reading word is $56248137$.
The inverse permutation $63.741.2.85$ has ascents in exactly the
positions $A(T)$, as marked by periods. The following theorem has been
highly inspirational in the theory of quasisymmetric functions and for
everything that follows in this subsection.

\begin{Theorem}\label{thm:gessel}\cite{Gessel} For any partition
$\lambda$, the Schur function expansion into the basis of fundamental
quasisymmetric functions is given by
\begin{equation}\label{eq:schur}
s_{\lambda}=s_{\lambda}(x_{1},x_{2},\dots) =\sum_{T \in \mathrm{SYT}(\lambda)}
F^{(d)}_{A(T)}.
\end{equation}
\end{Theorem}

Compare the definition of the fundamental quasisymmetric functions
\eqref{eq:fundamentals} with the definition of compatible pairs for a
permutation $w$ in \Cref{ex:compatible.pair} and the BJS formula for
Schubert polynomials in \Cref{cor:bjs} given by

\begin{equation}\label{eq:BJS.redo}
\fS_{w}=\fS_{w}(x_1,x_2,\ldots,x_n) := \sum_{(r_{1},\dots , r_{p})
    \in R(w)} \sum_{\substack{(i_{1}\leq i_{2} \leq \cdots \leq i_{p}):\\i_{k}<i_{k+1}\text{ if $r_{k}<r_{k+1}$}\\i_{k}\leq  r_{k}}}
        x_{i_{1}}x_{i_{2}}\cdots x_{i_{p}}.
\end{equation}
By removing the upper bound constraints $i_{k}\leq r_{k}$, we get a
quasisymmetric variation of $\fS_w$, which happens to be a symmetric
function. These functions were  originally
defined by Stanley in order to count the number of reduced words for a
given permutation \cite{Sta84}, as below.

\begin{Definition}\label{def:stanleys}
Let $w \in S_{\infty}$ be a permutation with $p$ inversions.  Define
the \textit{Stanley symmetric function} by
\begin{equation}\label{eq:stanley.def}
    G_{w}= G_{w}(x_{1}, x_{2}, \dots ) :=\sum_{(r_{1},\dots , r_{p}) \in R(w)}
\sum_{\substack{(i_{1}\leq i_{2} \leq \cdots \leq
i_{p}):\\i_{k}<i_{k+1}\text{ if $r_{k}<r_{k+1}$}}}
        x_{i_{1}}x_{i_{2}}\cdots x_{i_{p}}
=
\sum_{\mathbf{r}\in R(w)}F_{A(\mathbf{r})}^{(p)}
\end{equation}
where $A(\mathbf{r})=\{k \in [p-1]\given r_{k}<r_{k+1} \}$ is the
\textit{ascent set} of $\mathbf{r}=(r_{1},\dots , r_{p}) \in R(w)$.
\end{Definition}

\begin{Example}
If $w = 1432$ then $R(w) = \{232,323\}$, so 
\begin{align*}
G_{1432} &= F^{(3)}_{\{1\}} + F^{(3)}_{\{2\}}\\
& = (x_{1}x_{2}^{2}+x_{1}x_{2}x_{3} + x_{1}x_{3}^{2}+x_{1}x_{2}x_{4}+ \cdots) + (x_1^2 x_2 + x_1 x_2 x_3 + x_1^2 x_3 + x_1 x_2 x_4 + \cdots)\\
&= s_{(2,1)}.
\end{align*}
\end{Example}

Note that the coefficient of $x_{1}x_{2}\cdots x_{p}$ in $G_{w}$ will
be exactly the number of reduced words for $w$, since $(1,2,3,\dots
,p)$ appears with coefficient 1 in every $F_{A(\mathbf{r})}$ with
$\mathbf{r}\in R(w)$.  This observation does not immediately help us
to compute the number of reduced words for $w$, but it does give some
indication that these quasisymmetric functions are relevant to the
problem.  The key observation made by Lascoux and Sch\"utzenberger is
that $G_{w}$ is the limit of a sequence of certain Schubert
polynomials related to $w$ \cite{LS.Stanley}. Let's spell out this
connection and use pipe dreams again to prove it.

For $w=[w_{1},\dots, w_{n}] \in S_{n}$, let $1 \times w$ be the
\textit{shift} of $w$: the permutation with one-line notation $[1,w_{1}+1,
w_{2}+1,\dots , w_{n}+1]$. More generally, let $1^k \times w = [1,2,\ldots,k,w_{1}+k, w_2+k, \ldots, w_n+k]$ be the result of applying this operation $k$ times.
\begin{Theorem}\label{thm:limit.of.Schubs}
For $w \in S_{\infty}$, $G_w$ is a symmetric function and 
\begin{equation*} G_w = \lim_{k \to \infty} \fS_{1^{k}\times w}.
\end{equation*}
\end{Theorem}

\begin{proof}
Note that $w$ and $1 \times w$ have the same
number of reduced words and
\begin{equation}\label{eq:stanley.shift.equality}
G_w = G_{1\times w}
\end{equation}
since there is an easy bijection between $R(w)$ and $R(1\times w)$,
which preserves the ascent set.  On the other hand, the reduced pipe
dreams for $w$ bijectively correspond with the subset of pipe dreams
for $1 \times w$ with empty first row by
\Cref{thm:chutes.and.ladders}, using only ladder moves preserving the top row.  So $|\rp(w)|\leq |\rp(1\times w)|$,
and the inequality is strict provided $w$ is not the identity. 

Recall, the Schubert polynomial $\fS_{w}$ is symmetric in $x_{i}$ and
$x_{i+1}$ if and only if $w_{i}<w_{i+1}$ by
\Cref{ex:divided-difference}.  So $\fS_{1\times w}$ is symmetric in
$x_{1}$ and $x_{2}$.  If we shift the permutation again, then
$\fS_{1^{2}\times w}$ will certainly be symmetric in $x_{1},
x_{2},x_{3}$.  Continuing this way, we see $\fS_{1^{k}\times w}$ is
symmetric in the first $k$ variables.  As $k$ gets larger, the upper
bound constraints $i_{k}\leq r_{k}$ from \eqref{eq:BJS.redo} get
weaker until they are no longer binding in the limit, exactly matching
Definition~\ref{def:stanleys}.
\end{proof}

We briefly review a few facts on Grassmannian permutations from \Cref{sub:Grassmannians}. Recall that a permutation $w \in
S_{\infty}$ is \textit{Grassmannian} if it has at most one descent.
The permutation $w$ is $k$-Grassmannian if $w_{1}<w_{2}<\dots <w_{k}$
and $w_{k+1}<w_{k+2}<\dots$.  The identity permutation is
$k$-Grassmannian for all positive integers $k$.  The Schur polynomials
are exactly the Schubert polynomials indexed by Grassmannian
permutations by \Cref{def:schur.poly}. Specifically, if $w$ is $k$-Grassmannian then $\fS_w = s_{\lambda(w)}(x_1, \ldots, x_k)$ where $\lambda(w)$ is the partition obtained by sorting the code $c(w)$ into decreasing order, discarding any $0$'s.
\begin{Corollary} \label{cor:grassmannian.stanleys}  If $w$ is Grassmannian, then $G_w = s_{\lambda(w)}$.
\end{Corollary}
\begin{proof}
Write $\lambda = \lambda(w)$, so $\fS_{w}=s_{\lambda}(x_1, \ldots, x_k)$.  Observe that the shapes $\lambda(w)$ and $\lambda(1 \times w)$ are the same, since the code $c(1\times w)$ is obtained from $c(w)$ by prepending a $0$ on the front.
Thus, $\fS_{1\times w}=s_{\lambda}(x_1, \ldots, x_k,x_{k+1})$. By Theorem~\ref{thm:limit.of.Schubs}, the Stanley
symmetric function of $w$ is then 
\begin{equation}\label{eq:grass.stanley}
G_{w} = \lim_{k \to \infty} \fS_{1^{k}\times w} = \lim_{k \to \infty}
s_{\lambda}(x_1, \ldots, x_k) =  s_{\lambda}, 
\end{equation}
where $s_{\lambda}$ is the Schur function defined above via the limit
of Schur polynomials.
\end{proof}

\begin{Example}
In \Cref{ex:322}, we considered the case where $w=346125$ is
$3$-Grassmannian.  Every $1^{k}\times w$ is $(k+3)$-Grassmannian and
$\fS_{1^{k}\times w}=s_{(3,2,2)}(x_{1},x_{2},\dots , x_{k+3}).$
Therefore, $G_{346125}$ is the Schur function $s_{(3,2,2)}(x_{1}, x_{2}, \dots )$.
\end{Example}

\Cref{thm:limit.of.Schubs} allows us to connect the Stanley symmetric
functions $G_w$ with the transition equations for Schubert polynomials
(\Cref{t:transitionA}). Specifically, we get a recurrence relation for Stanley symmetric functions whose base case is the family of Grassmannian permutations,
for which we know that $G_w$ is a single Schur function. The precise statement is given in Theorem~\ref{thm:Stanley.Transitions}, and Example~\ref{ex:stanley-transitions} shows the simplest nontrivial case. The reader may also wish to look
ahead to \Cref{ex:Schurproduct} for the resulting Schur expansion of a
more complicated Stanley symmetric function.

\begin{Definition}\label{defn:transition.set}
For $w\in S_{\infty}$ with $w \neq id$, let $(r,s)$ be the lex largest
inversion of $w$ and $v = wt_{rs}$.  Let 
\[
T(w)=\{w'=vt_{ir} \in S_{\infty}\given  1\leq i <
r \text{ and } \ell(vt_{ir})=\ell(w) \}
\]
denote the \textit{transition set} of $w$.  
\end{Definition}

\begin{Example}
Let $w = 4312$. Then the lex largest inversion is $(r,s) = (2,4)$,
$v=wt_{rs}=4213$, and the set $T(w)$ is empty. On the other hand, the
lex largest inversion of $1 \times w = 15423$ is $(3,5)$, and $T(1 \times w) = \{35124\} = \{(1 \times
w)t_{35}t_{13}\}$.
\end{Example}

\begin{Theorem}[\textbf{Transition Equation for Stanley Symmetric
Functions}]\label{thm:Stanley.Transitions} For $w \in S_{\infty}$, the
Stanley symmetric function $G_{w}$ is a Schur positive symmetric
function, and its expansion into Schur functions can be found
recursively via the recurrence

\begin{equation}\label{eq:stanley.transitions}
G_{w}  = \sum_{w' \in T(1\times w)} G_{w'} = \begin{cases}
\sum_{w' \in T(w)} G_{w'}& T(w) \neq \emptyset \\
\sum_{w' \in T(1\times w)} G_{w'}& T(w) = \emptyset 
\end{cases}
\end{equation}
along with the base case that $G_{w}=s_{\lambda(w)}$ whenever $w$ is
Grassmannian. 
\end{Theorem}
 
\begin{proof}
If $w$ is Grassmannian, then by \Cref{cor:grassmannian.stanleys} we
know $G_{w}=s_{\lambda (w)}(x_{1},x_{2},\dots)$.  This includes the
case $w$ is the identity permutation for which 
$G_{\id}=s_{\emptyset}=\fS_{\id}=1$.  Hence, the proof of the base
case is complete.

The transition set $T(w)$ may be empty, but if it is nonempty, then
there is a bijection from $T(w)$ to $T(1 \times w)$ simply by mapping
$w' \mapsto 1\times w'$.  Furthermore, the set $T(1\times w)$ is
nonempty for every permutation $w$ other than the identity.  It is
clear from the bijection between $R(w)$ and $R(1 \times w)$ that
$G_w=G_{1 \times w}$ for all $w$, hence the second equality in the
theorem follows. This equality is useful for computations, since it
sometimes saves us from having to move into a larger symmetric group.

To prove $G_{w}$ expands into a positive integral sum of Schur
functions, observe that the recursion in
\eqref{eq:stanley.transitions} leads to a \textit{transition tree}
with $w$ as the root, each $w' \in T(1\times w)$ as the
children of $w$, each of which have children in $T(1\times
w')$, and continuing recursively until the leaves are all Grassmannian
permutations.  Each branch of the tree will terminate.  This is not as
obvious as in the case of the proof of \Cref{lem:homogeneous} using
inversion order from \Cref{defn:lex.largest.inversion.order} since we move into progressively larger symmetric groups as the recursion proceeds.  One way to prove each branch is finite
is to observe that the number of southeast corners of the connected
components of the corresponding Rothe diagrams are compressed along
each branch until they end up only in the lowest occupied row which implies the Grassmannian property \cite[Remark 5.16]{Billey-Pawlowski}.  Thus, it remains to prove that the
recursion in \eqref{eq:stanley.transitions} holds for permutations
which are not Grassmannian.

There are two slick proofs of \eqref{eq:stanley.transitions}.  The
first follows from the transition equations for Schubert polynomials
(\Cref{t:transitionA}, \Cref{thm:limit.of.Schubs}) and the stability
properties of $T(1^{k}\times w)$ for large enough $k$.  

The second proof of \eqref{eq:stanley.transitions} uses the Little
Bump \Cref{algorithm:little bump original} to give an ascent
preserving bijective map
\[ \label{eq:stanley-transition-bijection}
R(w) \longrightarrow \bigcup_{w' \in T(w)} R(w'),
\]
so by \Cref{def:stanleys} one observes
that \eqref{eq:stanley.transitions} holds.  We outline the steps here
following Little's approach in \cite{little2003combinatorial}. 
Assume $w$ is not Grassmannian, and hence not the identity
permutation. For simplicity we assume $T(w)$ is nonempty, replacing $w$ with $1 \times w$ if necessary as described above.

Consider the wiring diagram of a reduced word $\mathbf{r} \in R(w)$,
labelling the wires in order on the right and by $w^{-1}$ on the left.
Find the crossing corresponding with the lex largest inversion $(r,s)$
of $w$.  Apply a decrement push to the $(r,s)$ crossing to initiate
the Little bump algorithm. In the notation of Algorithm~\ref{algorithm:little bump original}, compute $\mathbf{r}' = \LBump^{-}_{t_0}(\mathbf{r})$ where $t_0$ is the column in which the crossing corresponding to inversion $(r,s)$ occurs.
One can observe from the Little bump algorithm that the ascent set of the input word is
preserved under this map.  The details of proving that the resulting $\mathbf{r}'$
is a reduced word for some $w' \in T(w)$ and the map is a
bijection as in \eqref{eq:stanley-transition-bijection} is similar to the proof of Theorem~\ref{t:transitionA} and
left to the reader. That is where the miracles happen!
\end{proof}

\begin{Example} \label{ex:stanley-transitions}
    Let $w = 2143$. Then $T(w) = T(2143) = \{2314, 3124\}$, so $G_{2143} = G_{2314} + G_{3124}$ by \eqref{thm:Stanley.Transitions}. Both $2314$ and $3124$ are Grassmannian, with respective shapes $\lambda = (1,1)$ and $\lambda = (2)$, so $G_{2143} = s_{11} + s_2$ by \Cref{cor:grassmannian.stanleys}.
\end{Example}

\begin{Remark}\label{rem:EG-insertion}
The history of \Cref{thm:Stanley.Transitions} involved several steps
in addition to the citations given in the proof.  Stanley first proved
each $G_{w}$ is symmetric in \cite{Sta84}.  He conjectured that $G_w$
also expands into a positive integral sum of Schur functions.  Schur
positivity was originally proved by Edelman-Greene in their paper on
balanced tableaux \cite{EG}. See also Lascoux and Sch\"utzenberger's
work on the \emph{plactic monoid} for an alternate approach
\cite{LS7}. Edelman and Greene gave a variation on the
Robinson-Schensted-Knuth (RSK) correspondence that applies to all
reduced words.  Assuming familiarity with RSK, the Edelman-Greene
correspondence requires only one slight change: if when inserting an
$i$ into row $j$ and an $i$ already exists in row $j$, then skip row $j$
and insert $i+1$ into row $j+1$.  Under the Edelman-Greene
correspondence, every reduced word $\mathbf{r} \in R(w)$ maps to a
pair of tableaux $(P(\mathbf{r}), Q(\mathbf{r}))$ of the same
partition shape where $Q(\mathbf{r})$ is a standard Young tableau
recording the order in which boxes are added, and $P(\mathbf{r})$ is the
insertion tableau which will necessarily be increasing along rows and
columns but not necessarily standard. One can reverse the
Edelman-Greene insertion process to recover $\mathbf{r} \in R(w)$ as
with RSK.  Furthermore, for each tableau $P$ of shape $\lambda$ that
occurs as $P({\bf r})$ for some reduced word ${\bf r} \in R(w)$, and
for each $Q \in \mathrm{SYT}(\lambda)$, there exists some
$\mathbf{r'}\in R(w)$ such that $(P,Q)=(P(\mathbf{r'}),
Q(\mathbf{r'})).$ Thus, reading up the columns of $P$ from left to
right must yield a reduced word for $w$ as well.  Hamaker-Young
observed that the set $\{P({\bf r}) \given {\bf r} \in R(w)\}$ can be
recovered using the leaves of the transition tree for $w$ and the
inverse Little bump algorithm in \cite{hamaker-young}.  Therefore, the
combinatorics of reduced words and the Edelman-Green correspondence is
quite rich.  We highlight the following two important theorems.
\end{Remark}

\begin{Theorem}\cite{EG}\label{thm:counting.reduced.words}
For $w \in S_{\infty}$, if $G_{w}=\sum_{\lambda}
a_{\lambda}s_{\lambda}$, then each $a_{\lambda}$ is a nonnegative
integer and the number of reduced words for $w$ is given by
\begin{equation}\label{eq:counting.reduced.words}
|R(w)|= \sum_{\lambda} a_{\lambda} |\mathrm{SYT}(\lambda)|
\end{equation}
where $|\mathrm{SYT}(\lambda)|$ is easily determined by the
Frame-Robinson-Thrall Hook Length Formula.
\end{Theorem}

\begin{proof}
The first statement follows recursively from
\Cref{thm:Stanley.Transitions}.  Observe from
\Cref{def:stanleys} that if $\ell(w)=p$, then the coefficient of
$x_{1}x_{2}\cdots x_{p}$ in $G_{w}$ is exactly $|R(w)|$.  Similarly,
the coefficient of $x_{1}x_{2}\cdots x_{p}$ in $s_{\lambda}$ is
exactly $|\mathrm{SYT}(\lambda)|$ by Gessel's \Cref{thm:gessel}.  Therefore,
Equation \eqref{eq:counting.reduced.words} follows by comparing
coefficients on both sides of the equation $G_{w}=\sum_{\lambda}
a_{\lambda}s_{\lambda}$.
\end{proof}

In \cite{Greene-Nijenhuis-Wilf}, Greene, Nijenhuis, and Wilf gave a
beautiful algorithm to find a standard tableau of a given shape
$\lambda$ uniformly at random by what they call the ``Hook Walk Algorithm''.
Briefly, if $\lambda$ is a partition of $d$, then to place the letter
$d$ in the diagram of $\lambda$, choose any cell $c$ uniformly at
random.  If $c$ is a corner cell, place $d$ in that cell. If $c$ is
not a corner cell, uniformly at random choose a different cell $c'$ in
the hook of $c$, meaning either directly to the right of $c$ or directly
below it.  If $c'$ is in a corner cell, place $d$
there, and if not, again choose a different cell in the hook of $c'$, 
etc., until a corner cell is found and $d$ is placed.  Then continue
similarly to place $d-1$ in the remaining partition shape, etc.  The
Hook Walk Algorithm terminates when the values $1,2,\dots, d$ are all
placed.  The result is guaranteed to be standard. The next result
combines the Hook Walk Algorithm with the insertion algorithm of
Edelman-Greene mentioned in \Cref{rem:EG-insertion}.

\begin{Theorem}\label{thm:uniform.reduced.word}
For $w \in S_{\infty}$, let $\mathcal{P}=\{P(\mathbf{r})\given
\mathbf{r}\in R(w) \}$.  A random reduced word for $w$ can be found by
choosing a partition $\mu$ of $\ell(w)$ with probability
\[
\frac{a_{\mu }|\mathrm{SYT}(\mu)|}{\sum a_{\lambda}|\mathrm{SYT}(\lambda)|},
\]
choosing $P$ uniformly among all tableaux in $\mathcal{P}$ of shape
$\mu$, using the Hook Walk Algorithm to choose $Q \in \mathrm{SYT}(\mu)$
uniformly, and then applying inverse Edelman-Greene insertion to
$(P,Q)$ to find $\mathbf{r} \in R(w)$.
\end{Theorem}

The following class of permutations has many interesting properties relating to Stanley symmetric functions, some of
which are included in the next theorem.
\begin{Definition}\label{def:vexillary}
If $w \in S_{\infty}$ avoids the pattern 2143, then $w$ is a
\textit{vexillary permutation}.
\end{Definition}
In fact, we have already seen vexillary permutations in \Cref{ex:bpd-poset-vexillary}. They were first recognized and named
by Lascoux and Sch\"utzenberger in the 1980's.  The word vexillary comes from Roman times and could be translated as
``one of a special class of soldiers carrying a certain flag,'' hence
the connection with flag varieties.  Continuing the theme, one can
also find patriotic and heroic permutations in the literature
\cite{b-j-s}.\footnote{Hermann C\"asar Hannibal Schubert's first three names come from famous generals of antiquity. Coincidence?}

\begin{Theorem}\label{cor:vexillary.Stanleys}  For all $w \in
S_{\infty}$, the following are equivalent.
\begin{enumerate}
\item The permutation $w$ is vexillary.
\item The Stanley symmetric function is a single Schur function
$G_{w}=s_{\lambda(w)}$, where $\lambda (w)$ is the partition obtained
from the code $c(w)$ by sorting into decreasing order.
\item The number of reduced words for $w$ is $|\mathrm{SYT}(\lambda(w))|$. 
\item Every reduced word ${\bf r} \in R(w)$ inserts to the same tableau $P = P({\bf r})$ under
Edelman-Greene insertion. Furthermore, $P$ is easily found by inserting any
reduced word for $w$, so choosing a reduced word for $w$ uniformly at random
is efficiently done via the hook walk algorithm.
\item The Schubert polynomial $\fS_{w}$ is a flagged Schur function,
which is the sum over all semistandard tableaux satisfying the flagged
conditions in \Cref{def:balanced-flagged-column-injective}. The
flagged Schur functions can also be computed as a determinant
analogous to the Jacobi-Trudi formula. 
\item The rows and the columns of $D(w)$ can be ordered by containment.
\end{enumerate}
\end{Theorem}

\begin{proof}
The equivalence of (1), (2) and (6) is a good exercise for the reader.
See \Cref{exercise:vex}.  Statements (2) and (3) are equivalent by
\Cref{thm:counting.reduced.words}.  The equivalence with (4) follows
from \Cref{rem:EG-insertion}.  The proof of (5) that vexillary
Schubert polynomials are flagged Schur functions is due to Wachs
\cite{wachs}.
\end{proof}

Another application of Stanley symmetric functions and transitions is
to the computation of the Littlewood-Richardson coefficients for Schur
functions from \Cref{def:LR.coeffs}.  Say we want to expand the
product $s_{\lambda}(X_k) s_{\mu}(X_k)$ into a sum of Schur
polynomials.  We know we can use the puzzle rule as described in
\Cref{thm:KnutsonTao-puzzles}.  Alternatively, assuming $s_{\lambda}(X_k)
s_{\mu}(X_k)$ is not zero, one can find the unique $k$-Grassmannian
permutations $u,v \in S_{\infty}$ such that their shapes match the
input shapes $\lambda (u) =\lambda$ and $\lambda (v) =\mu$.  By \Cref{cor:grassmannian.stanleys} and 
\Cref{ex:product.perms}, we
know $G_{u \times v}=G_{u}G_{v}=s_{\lambda}s_{\mu}$.  On the
other hand, we can expand $G_{u \times v}$ into Schur functions by
using the transition equation \eqref{eq:stanley.transitions}.  As $k
\to \infty$, the expansion coefficients stabilize in the sense that if
\[
s_{\lambda}(X) s_{\mu}(X) = \sum c_{\lambda, \mu}^{\nu }
s_{\nu}(X) 
\]
where $X$ is the countably infinite set of variables
$x_{1},x_{2},\dots$, then 
\[
s_{\lambda}(X_{k}) s_{\mu}(X_k) = \sum c_{\lambda, \mu}^{\nu } s_{\nu}(X_k)
\]
where the second sum is restricted to those $\nu$ with at most $k$
parts.  This approach to computing Littlewood-Richardson coefficients
is due to Lascoux-Sch\"utzenberger \cite{LS2}.  It is reminiscent of
the Remmel-Whitney rule \cite{RW}.

\begin{Example}\label{ex:Schurproduct}
Consider the case $\lambda = (3,2)$, $\mu =(1,1)$.  To find the
expansion of $s_{(3,2)}s_{(1,1)}$ into Schur functions, we set
$u=35124$ and $v=231$ so $u \times v = 35124786$. The lex largest
inversion for $35124786$ is $(r,s)=(7,8)$ and the transition set is
$T(w)= \{35124768t_{27},\ 35124768t_{57}\}$.  Then we apply the
transition equation to both of these permutation in $T(w)$ recursively
creating the transition tree shown in
\Cref{fig:transition-tree-35124786}.  We chose to terminate the
branches as soon as the leaves are vexillary permutations since the
corresponding Stanley symmetric functions are then easily computed by
Part (2) of \Cref{cor:vexillary.Stanleys}.  Therefore, we obtain the
Schur expansion
\[
G_{35124786}=s_{(3,2)}s_{(1,1)} = s_{(3, 2, 1, 1)}+ s_{(3, 3, 1)}+s_{(4, 2, 1)}+s_{(4, 3)}.  
\]
\begin{figure}[h!]
\centering
\begin{tikzpicture}[scale=0.6]
\def\a{1.4};
\def\b{0.4};
\def\h{2.4};
\newcommand\Rec[3]{
\node at (#1,#2) {#3};
\draw(#1-\a,#2-\b)--(#1-\a,#2+\b)--(#1+\a,#2+\b)--(#1+\a,#2-\b)--(#1-\a,#2-\b);
}
\Rec{0}{0}{$35124786$}
\Rec{-2.5*\a}{-\h}{$35126748$}
\Rec{2.5*\a}{-\h}{$36124758$}
\Rec{-5*\a}{-2*\h}{$35146278$}
\Rec{0}{-2*\h}{$45126378$}
\Rec{5*\a}{-2*\h}{$36125478$}
\Rec{2.5*\a}{-3*\h}{$36142578$}
\Rec{7.5*\a}{-3*\h}{$46123578$}
\draw(0,-\b)--(-2.5*\a,-\h+\b);
\draw(0,-\b)--(2.5*\a,-\h+\b);
\draw(-2.5*\a,-\h-\b)--(-5*\a,-2*\h+\b);
\draw(-2.5*\a,-\h-\b)--(0*\a,-2*\h+\b);
\draw(2.5*\a,-\h-\b)--(5*\a,-2*\h+\b);
\draw(5*\a,-2*\h-\b)--(2.5*\a,-3*\h+\b);
\draw(5*\a,-2*\h-\b)--(7.5*\a,-3*\h+\b);
\node[below] at (-5*\a,-2*\h-\b) {$s_{(3,2,1,1)}$};
\node[below] at (0*\a,-2*\h-\b) {$s_{(3,3,1)}$};
\node[below] at (2.5*\a,-3*\h-\b) {$s_{(4,2,1)}$};
\node[below] at (7.5*\a,-3*\h-\b) {$s_{(4,3)}$};
\end{tikzpicture}
\caption{The transition tree for $35124786$ terminating with vexillary leaves.}\label{fig:transition-tree-35124786}
\end{figure}
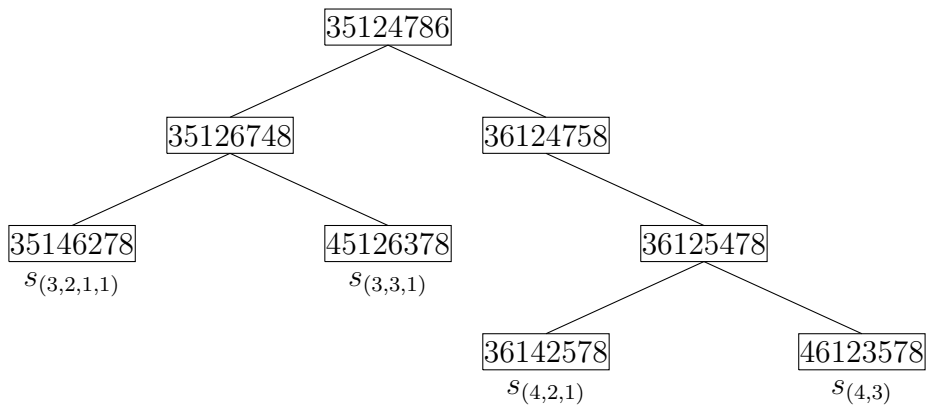

% CL-USER> (an-stanley (append '(3 5 1 2 4) (loop for i in '(2 3 1) collect (+ i 5))))
%   0: (STANLEY-REC (3 5 1 2 4 7 8 6) 7 A)
%     1: (ST-COMPUTE-WS (3 5 1 2 4 7 6 8) 7 7 A)
%     1: ST-COMPUTE-WS returned ((3 5 1 2 6 7 4 8) (3 6 1 2 4 7 5 8))
%     1: (STANLEY-REC (3 5 1 2 6 7 4 8) 7 A)
%       2: (ST-COMPUTE-WS (3 5 1 2 6 4 7 8) 6 7 A)
%       2: ST-COMPUTE-WS returned ((3 5 1 4 6 2 7 8) (4 5 1 2 6 3 7 8))
%       2: (STANLEY-REC (3 5 1 4 6 2 7 8) 7 A)
%       2: STANLEY-REC returned (<1.SHP:(3 2 1 1)>)
%       2: (STANLEY-REC (4 5 1 2 6 3 7 8) 7 A)
%       2: STANLEY-REC returned (<1.SHP:(3 3 1)>)
%     1: STANLEY-REC returned (<1.SHP:(3 2 1 1)> <1.SHP:(3 3 1)>)

%     1: (STANLEY-REC (3 6 1 2 4 7 5 8) 7 A)
%       2: (ST-COMPUTE-WS (3 6 1 2 4 5 7 8) 6 7 A)
%       2: ST-COMPUTE-WS returned ((3 6 1 2 5 4 7 8))
%       2: (STANLEY-REC (3 6 1 2 5 4 7 8) 7 A)
%         3: (ST-COMPUTE-WS (3 6 1 2 4 5 7 8) 5 7 A)
%         3: ST-COMPUTE-WS returned ((3 6 1 4 2 5 7 8) (4 6 1 2 3 5 7 8))
%         3: (STANLEY-REC (3 6 1 4 2 5 7 8) 7 A)
%         3: STANLEY-REC returned (<1.SHP:(4 2 1)>)
%         3: (STANLEY-REC (4 6 1 2 3 5 7 8) 7 A)
%         3: STANLEY-REC returned (<1.SHP:(4 3)>)
%       2: STANLEY-REC returned (<1.SHP:(4 2 1)> <1.SHP:(4 3)>)
%     1: STANLEY-REC returned (<1.SHP:(4 2 1)> <1.SHP:(4 3)>)
%   0: STANLEY-REC returned
%        (<1.SHP:(3 2 1 1)> <1.SHP:(3 3 1)> <1.SHP:(4 2 1)> <1.SHP:(4 3)>)
% (<1.SHP:(3 2 1 1)> <1.SHP:(3 3 1)> <1.SHP:(4 2 1)> <1.SHP:(4 3)>)
% CL-USER> 
The same expansion can also be obtained via the Edelman-Greene map (\Cref{fig:P-tableaux-35124786}).
\begin{figure}[h!]
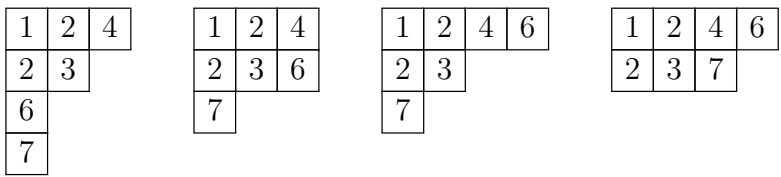

\centering
\[\ytableausetup{boxsize=1.3em}
\begin{ytableau}
1&2&4 \\ 2&3 \\ 6 \\ 7
\end{ytableau}\qquad
\begin{ytableau}
1&2&4 \\ 2&3&6 \\ 7
\end{ytableau}\qquad
\begin{ytableau}
1&2&4&6 \\ 2&3 \\ 7
\end{ytableau}\qquad
\begin{ytableau}
1&2&4&6 \\ 2&3&7
\end{ytableau}
\]
\caption{The $P$-tableaux using the   Edelman-Greene correspondence for the permutation $35124786$}
\label{fig:P-tableaux-35124786}
\end{figure}
\end{Example}

Theorem~\ref{thm:grassmannian-cohom} shows that the cohomology class
of any subvariety $Z \subseteq \Gr(k,n)$ can be represented as a
nonnegative linear combination of Schur functions. This suggests an
inverse question:

\begin{quote}
Given a Schur-positive symmetric function, does it
represent the class of a subvariety of $\Gr(k,n)$?
\end{quote}

If reducible subvarieties are allowed, the answer is trivially yes
(take an appropriate union of Schubert varieties), so let us insist
the subvariety be irreducible. Such questions can be very hard. For
instance, the \emph{permutohedral variety} $X_{A_{n}}$ is the closure of
the orbit of invertible diagonal matrices $T$ acting by left
multiplication on a generic flag $F_\bullet \in \Fl(n+1)$, and Huh
\cite{HuhThesis} showed that deciding whether a Chow class $\alpha \in
A^*(X_{A_{n}})$ comes from a subvariety of $X_{A_{n}}$ is at least as hard
as deciding representability of a matroid. However, there turns out to
be a nice answer for the Stanley symmetric functions $G_w$.

To describe a subvariety whose cohomology class is represented by $G_w$, we will generalize the notion of a Grassmannian Schubert variety. Fix a basis $e_1, \ldots, e_n$ of $\C^n$, choose a sequence of integers $m_0, m_1, \ldots, m_n \geq 0$, and consider the set
\begin{equation*}
    Y = \{V \in \Gr(k,n) \given \dim(V \cap \langle e_1, \ldots, e_j \rangle) \geq m_j \text{ for all $j \in [n]$}\}.
\end{equation*}
Recall from \eqref{eq:grassmannian-rank-conditions} that if $0 = m_0 \leq m_1 \leq \cdots \leq m_n = n$ and $m_i-m_{i-1} \in \{0,1\}$ for all $i$, then $Y$ is just the Schubert variety $X_I(E_\bullet)$ where $I = \{i \in [n] \given m_i > m_{i-1}\}$ and $E_j = \langle e_1, \ldots, e_j \rangle$. In general, such a $Y$ is Zariski-closed and is closed under the $B$-action on $\Gr(k,n)$. Therefore, it is a union of closures of $B$-orbits, which are exactly the Schubert varieties in $\Gr(k,n)$. 

Now generalize the definition further to allow \emph{arbitrary} intervals of basis vectors. That is, choose an array $(m_{ij})_{1 \leq i \leq j \leq n}$ of nonnegative integers and consider the varieties of the form
\begin{equation} \label{eq:interval-positroid}
    Y' = \{V \in \Gr(k,n) \given \dim(V \cap \langle e_i, \ldots, e_j \rangle) \geq m_{ij} \text{ for all $1 \leq i \leq j \leq n$}\}.
\end{equation}
One can give combinatorial conditions on the array $(m_{ij})$ guaranteeing that the associated $Y'$ is irreducible, in which case it is called an \emph{interval positroid variety}. In general, such a $Y'$ is a union of interval positroid varieties. In \cite{pawlowski-liu-conjecture}, Pawlowski showed that these varieties have the desired cohomology classes.
\begin{Theorem} \cite[Theorem 26]{pawlowski-liu-conjecture} \label{thm:interval-positroid} The cohomology class of an interval positroid variety is represented by a Stanley symmetric function.
\end{Theorem}
Theorem~\ref{thm:interval-positroid} gives a positive answer to our question in the case of Stanley symmetric functions. However, we would like to be more explicit and give a \emph{specific} variety with cohomology class $G_w$, for any given $w \in S_n$. For this, a dual perspective is useful. Define $\Sigma_w^\circ \subseteq \Gr(n,2n)$ as the set of $n$-planes $\Span \{v_1, \ldots, v_n\} \subseteq \C^{2n}$ with each $v_i \in \Span \{e_j \given j \in [w(i), n+i]\}$. For example,  $\Sigma_{231}^\circ$ is the set of points \begin{equation*}
    \operatorname{colspan} \begin{bmatrix}
0 & 0    & \ast \\
\ast & 0 & \ast \\
\ast & \ast & \ast \\
\ast & \ast & \ast \\
0   & \ast & \ast \\
0    & 0 & \ast \\
    \end{bmatrix}
\end{equation*}
in $\Gr(3,6)$. This is an irreducible set, so its closure $\Sigma_w = \overline{\Sigma_w^{\circ}}$ is an irreducible subvariety of $\Gr(n,2n)$.

\begin{Theorem}\cite[\S 6]{KLS}.\label{thm:KLS}
The cohomology class
$[\Sigma_{w}]$ is represented by $G_w$ in $H^{*}(\Gr(n,2n))$.
\end{Theorem}
It is natural to broaden this definition to $\Gr(k,n)$ and allow any set of $k$ intervals in $[n]$ in place of $\{[w(i), n+i] \given i \in [n]\}$. This results in the family of \emph{rank varieties} introduced by Billey and Coskun \cite{Billey-Coskun}, also called \emph{dual interval positroid varieties} \cite{interval-positroid}. No new cohomology classes are obtained this way, as can be shown by combining general results on matroid duality with Theorem~\ref{thm:interval-positroid}.

A very interesting further generalization is obtained by modifying \eqref{eq:interval-positroid} to allow for \emph{cyclic} intervals of vectors, i.e. ordinary intervals together with those of the form $\langle e_j, \ldots, e_n, e_1, \ldots, e_i \rangle$. The resulting irreducible varieties are called \emph{positroid varieties}. As the names suggest, this family of varieties includes both the interval positroid varieties and the dual interval positroid varieties. Knutson, Lam, and Speyer \cite{KLS} introduced positroid varieties building on work of Postnikov \cite{postnikov-2006}, and showed that their cohomology classes are represented by \emph{affine Stanley symmetric functions}. The nontrivial content of Theorem~\ref{thm:interval-positroid} is then that ordinary Stanley symmetric functions actually suffice in the case of interval positroid varieties. Positroid varieties also have a surprising connection to quantum field theory via their relationship to the \emph{amplituhedron} \cite{ScatteringAmplitudes}; see also \cite{BourjailyThomas} for a friendly introduction.

Positroid varieties can also be defined by collections of minors such
that some subspace $V \in \Gr (k,n)$ has all positive Pl\"ucker
coordinates on those minors and all others vanish.  Such subsets
determine a special type of matroid, called a \emph{positroid}
\cite{postnikov-2006}.  A famous theorem of Mn\"ev \cite{Mnev1988} asserts
that, for general subsets $\mathcal{P} \subseteq {[n] \choose k}$, the
matroid variety $\overline{X_{\mathcal P}}$ given by certain Pl\"ucker
coordinates vanishing can have arbitrarily complicated
singularities. By contrast, positroid varieties have very nice
geometric and combinatorial properties, as will be discussed further
in Chapter 4 \cite{speyer2024richardsonvarietiesprojectedrichardson}.

\begin{Exercise}\label{eq:qsym.play}
For a composition $\alpha$ of $d$, find the expansion of $M_{\alpha}$
in terms of the fundamental basis $F_{S}^{(d)}$.  Similarly, find the
inverse expansion giving  $F_{S}^{(d)}$ in terms of monomial
quasisymmetric functions $M_{\alpha}$.
\end{Exercise}

\begin{Exercise}\label{exercise:Grassmannian.Schur}
Let $w$ be a Grassmannian permutation with unique descent in position
$k$.  Prove $\fS_{w}$ is the Schur polynomial $s_{\lambda
(w)}(x_{1},\dots , x_{k}) $ by giving a bijection from pipe dreams in
$\rp(w)$ to semistandard Young tableaux of shape $\lambda (w)$ on the
alphabet $\{1,2,\dots , k \}$.   (Hint, consider the crossings along
one string at a time, or see \cite{Serrano.Stump.2012} for a
nice bijection along with many others.)
\end{Exercise}

\begin{Exercise}\label{exercise:vex} Prove the following are
equivalent for any permutation $w \in S_{\infty}$.
\begin{enumerate}
\item The permutation $w$ is vexillary, so it avoids the pattern
$2143$.
\item Expanding the transition set $T(w)$ recursively via
\Cref{thm:Stanley.Transitions} leads to a transition tree for $w$ that
is a path from $w$ to a Grassmannian permutation.
\item The rows of the diagram $D(w)$ are all comparable in the sense
that if you labeled each bubble in $D(w)$ by its column number, then
the subset of values appearing in row $i$ and the subset appearing in
row $j$ are comparable by subset inclusion for all $i,j$.
\item The columns of the diagram $D(w)$ are all comparable as subsets
of occupied rows 
\item The last two properties imply there is a permutation of the rows
and permutation of the columns that sort the bubbles in $D(w)$ into
the shape of a Young diagram.  Prove that every diagram for a
permutation in the transition tree of a vexillary permutation sorts to
the same partition shape.
\end{enumerate}

\end{Exercise}

\begin{Exercise}\label{ex:product.perms}
Say $u \in S_{m}$ and $v \in S_{n}$.  Given two permutations $u=(u_{1},u_{2},\dots , u_{m})$ and $v=(v_{1},\dots ,
v_{n})$, define
\[
u\times v = (u_{1},u_{2},\dots , u_{m}, v_{1}+m, v_{2}+m, \dots ,
v_{n}+m) \in S_{m+n}.  
\]
Show that $G_{u \times v} = G_{v \times u}=G_{u}G_{v}$.
\end{Exercise}

\begin{Exercise} \label{exer:graph-schubert}
    The \emph{graph Schubert variety} $\operatorname{GX}_w$ of $w \in S_n$ is the Zariski
closure of the set of points $\operatorname{colspan}
\left[ \begin{smallmatrix} A \\ I \end{smallmatrix} \right] \in
\Gr(n,2n)$ where $I$ is an $n \times n$ identity matrix and $A$
represents a flag in the \emph{opposite} Schubert cell
$C_w(\oppositeE_\bullet)$. Show that the interval positroid variety
$\Sigma_{w}$ is the graph Schubert variety $\operatorname{GX}_{w_0 w}$. (The terminology comes from the fact that
$\operatorname{colspan} \left[ \begin{smallmatrix} A \\
I \end{smallmatrix} \right]$ is the graph $\{(Av,A) \given v \in
\C^n\} \subseteq \C^n \times \C^n$ of the linear function $A$. If the reader consults our citation for Theorem~\ref{thm:KLS}, they will find graph Schubert varieties rather than our varieties $\Sigma_w$: the point of the exercise is that there is no difference.)
\end{Exercise}

Recall a Schubert polynomial $\fS_{w}$ expands with only zero-one
coefficients if and only if $w$ avoids 12 permutation patterns by
\Cref{thm:zero-one-schubert}.  Similar pattern avoiding
characterizations can be stated for Stanley symmetric functions
\cite{Billey-Pawlowski}.  Say $G_{w}$ is \emph{zero-one} if $G_{w} =
\sum c_{\alpha}x^{\alpha}$ and all $c_{\lambda}\in \{0,1 \}$, and say
$G_{w}$ is \emph{multiplicity-free} if $G_{w} = \sum
a_{\lambda}s_{\lambda}$ and all $a_{\lambda}\in \{0,1 \}$.  The two
notions of zero-one coefficients are studied further in
\cite{Hodges-Yong.2022}.

\begin{Exercise}
The Stanley symmetric function $G_w$ is zero-one if and only if
$|R(w)|=1$.
\end{Exercise}

\begin{Conjecture}\cite[Conj. 1]{Billey-Pawlowski}
The multiplicity-free Stanley symmetric functions $G_w$ are
characterized by avoiding 189 patterns in $S_{6}$ to $S_{11}$
\cite[P0045]{dppa}.  
\end{Conjecture}

%% file: section5.tex
\section{Combinatorial Algebraic Geometry of Schubert Varieties}\label{sec:CAGofSchubertVarieties}

In this chapter, we will present some of the algebraic geometry and
commutative algebra associated with Schubert varieties.  We first
describe the homogeneous equations defining Schubert varieties in all
partial flag varieties and their projective embeddings using Pl\"ucker
coordinates and Pl\"ucker relations.  Then we describe a family of
affine varieties with the same equations.  These are called matrix
Schubert varieties as first defined by Fulton.  Knutson and Miller
gave a very insightful theorem on matrix Schubert varieties showing
that the pipe dream formula for double Schubert polynomials gives rise
to the associated multidegree.  One might call this a ``geometric
reason'' that Schubert polynomials are indeed the right choice of
representatives for Schubert classes in the Borel presentation of the
cohomology ring.  Schubert varieties in the Grassmannians have also
been used in the study and proof of Horn's conjecture on eigenvalues
of Hermitian matrices.  We will survey this important topic here using
the Schubert structure constants and intersection theory. We also discuss the theory of degeneracy loci for vector bundles, which allows the theory we have developed for flag varieties be extended to varieties whose points have vector spaces associated with them, such as the tangent space at a point or a vector assigned by a vector field. We will    
give an overview of the study of the singular locus of Schubert
varieties.  The theory of pattern avoidance of permutations plays a
key role in this topic highlighting a fundamental tool in the
computational side of Schubert varieties and connections to
combinatorial machine learning.  Finally, we will present the complete
characterization of all isomorphism types of Schubert varieties in all
flag varieties due to Richmond-Slofstra, which is elegant and
surprisingly simple to state.  Furthermore, it was only discovered while this
chapter has been in preparation, proving this topic is very much alive
in modern research.

\subsection{Pl\"ucker Coordinates and Relations}\label{sub:Plucker}

In order to view each Schubert variety $X_{w}(E_{\bullet }) \subset
\flags$ as a projective variety, we embed it into a product of
projective spaces using Pl\"ucker coordinates and flag minors.  This
procedure is closely related to the Pl\"ucker coordinates from
\eqref{eq:plucker.coords.gnk} used to embed the the Grassmannian
variety $\Gr(k,n)$ in projective space in
\Cref{sub:Grassmannians.intro}.  Let $\binom{[n]}{k}$ denote the size
$k$ subsets of $[n]$.  For $I \in \binom{[n]}{k}$ and an $n\times n$
matrix $M$, let $\Delta_{I}(M)$ be the determinant of the submatrix of
$M$ in rows $I$ and columns $[k]=\{1,2,\ldots,k \}$.  The polynomial
functions $\{\Delta_{I}: I \in \binom{[n] }{k}\}$ are called
\textit{flag minors} because of their connection with flag varieties.
They also play a critical role in cluster algebras \cite[Ch. 1]{FWZ}.

\begin{Definition}\label{def:plucker}
For $F_{\bullet}=(F_{1} \subset F_{2} \subset \cdots \subset F_{n})
\in \flags$, let $P_{k}(F_{\bullet})=P_{k}(M) \in
\mathbb{P}^{\binom{n}{k}-1}$ be the list of all flag minors
$\Delta_{I}(M)$ in lex order on the $k$-subsets $I \subset [n]$ for
any $n \times n$ invertible matrix $M$ representing $F_{\bullet }$.
The list $P_{k}(F_{\bullet}) $ thought of as a homogeneous coordinate
constitutes the \textit{Pl\"ucker coordinates} of the $k$-dimensional
subspace $F_{k}$.  The concatenation
$P(F_{\bullet})=P_{1}(F_{\bullet})\times P_{2}(F_{\bullet}) \times
\cdots \times P_{n-1}(F_{\bullet})$ gives the \textit{Pl\"ucker coordinates}
of $F_{\bullet}$.  For partial flags in $\Fl(n;\mathbf{d})$, the
\textit{Pl\"ucker coordinates} are restricted to the subspaces indexed
by $\mathbf{d}$.   
\end{Definition}

Observe that if $M$ is invertible and represents $F_\bullet$, at least one of the $\binom{n}{k}$
flag minors of size $k \times k$ is nonzero for each $1\leq k\leq n$,
so each $P_{k}(F_{\bullet})$ is indeed a valid homogeneous coordinate.
Furthermore, for any $b \in B$, we have
$\Delta_{I}(Mb)=\Delta_{I}(M)\Delta_{I}(b)$, since columns $k+1,\dots
, n$ in $b$ will have no effect on the first $k$ columns of the
product $Mb$. Therefore, the homogeneous coordinates $P(F_{\bullet})$ are
the same no matter which $M$ in its coset $MB$ is chosen to represent
it, up to overall rescaling.  Since $P_{n}(F_{\bullet})=[1]$ always, we don't need to include
it in the Pl\"ucker coordinates of $F_{\bullet}$. 

\begin{Example}\label{example:flag.from.section3}
Recall the flag  $F_{\ci}=(2 e_{1}+ e_{2},\hspace{.1in} 2 e_{1} +e_{3},\hspace{.1in}
7e_{1}+e_{4},\hspace{.1in} e_{1})$ represented by matrices in
\eqref{eq:example.options}.  Then, 
\begin{equation}
P(F_{\ci}) = [2:1:0:0]\times [-2:2:0:1:0:0] \times [7:-2:2:1].
\end{equation}
\end{Example}

\begin{Exercise}
Prove that a basis for each subspace in a flag $F_{\bullet} \in
\Fl(n)$ can be recovered from its Pl\"ucker coordinates
$P(F_{\bullet})$. Hence, one can recover the canonical matrix
$M(F_{\bullet})$ from $P(F_{\bullet})$.  Consider the canonical
matrix representative of the flag from \Cref{def:canonical.mat}.
\end{Exercise}

Determinants of matrices have many interesting properties. Since the
flag minors are determinantal functions, some of these properties lead
to relations among the values in the Pl\"ucker coordinates of
a subspace or a flag.  For example, Sylvester's Lemma gives relations
on the product of determinants for two square matrices.  The proof
below is based on a sign reversing involution given by Konvalinka.
See also the proof in \cite[Lemma 2, \S 8.2]{Fulton-book} which cites
Sylvester (1851).

\begin{Lemma}[\textbf{Sylvester's Lemma}] \label{lem:sylvester}
For $n\times n$ matrices $A$ and $B$ and any fixed proper subset $T
\subset [n]$, we have
\begin{equation}\label{eq:sylvester}
\det(A) \det(B) = \sum_{S \in \binom{[n] }{r}} \det(A_{S}) \det(B_{S})
\end{equation}
where $r=|T|$,\ $A_{S}$ is obtained from $A$ by replacing the columns indexed by
$S=\{s_{1}<\dots <s_{r} \}$ with the columns of $B$ indexed by
$T=\{t_{1}<t_{2}<\dots <t_{r} \}$ in order, and $B_{S}$ is obtained
from $B$ by replacing the columns indexed by $T$ with the columns of
$A$ indexed by $S$, again maintaining the order.
\end{Lemma}

\begin{proof}
 By rearranging the columns of $B$ and each $B_{S}$ by a fixed
permutation if necessary, it suffices to prove the lemma in the case
$T=[r]$ by the alternating property of determinants.  The determinant
of $A=(a_{i,j})$ can be written as a signed sum over permutations as 
\[
\det(A) = \sum_{\pi \in S_{n}} (-1)^{\inv(\pi)}a_{\pi (1),1} a_{\pi
(2),2}\cdots a_{\pi (n),n}.
\]
Similarly, by construction  if $B=(b_{i,j}),$ we have 
\[
\det(A_{S}) = \sum_{\pi \in S_{n}} (-1)^{\inv(\pi)} \prod_{c \in [n]
\setminus S}a_{\pi (c),c} \cdot \prod_{d=1}^{r}b_{\pi (s_{d}),d} 
\]
and 
\[
\det(B_{S}) = \sum_{\sigma  \in S_{n}} (-1)^{\inv(\sigma )} \prod_{e
=1}^{r}a_{\sigma (e),s_{e}} \cdot \prod_{f=r+1}^{n}b_{\sigma(f),f}.
\]
Therefore, if we expand the right side
of \eqref{eq:sylvester}, each term corresponds to a triple
$(S,\pi,\sigma)$, where $S= \{s_{1}<\dots <s_{r} \} \in
\binom{[n]}{r}$ and $\pi,\sigma \in S_n$.  The \textit{weight} of the
$(S,\pi,\sigma)$ term is 
\begin{equation}\label{eq:weight.s.pi.sigma}
 \prod_{c \in [n] \setminus S}a_{\pi
(c),c} \cdot \prod_{d=1}^{r}b_{\pi (s_{d}),d} \cdot \prod_{e =1}^{r}a_{\sigma
(e),s_{e}} \cdot \prod_{f=r+1}^{n}b_{\sigma(f),f} 
\end{equation}
and the \textit{sign} is
\begin{equation}\label{eq:sign.s.pi.sigma}
(-1)^{\inv(\pi) + \inv(\sigma)}.
\end{equation}
  We want to find a weight-preserving, sign-reversing involution on these
terms with $(n!)^{2}$ fixed points that bijectively map to the terms
on the left side of \eqref{eq:sylvester}, and the final bijection must
also preserve sign and weight.

Define a map $\iota$ on the triple $(S,\pi,\sigma)$ as follows.  If
$\pi (S)=\sigma ([r])$
as sets, the triple is a fixed point under $\iota$.
%%%%%%%%%%%%%%%%
Otherwise, there exists a pair $(i,j)$ so that $i \notin S$, $j \leq
r$, and $\pi(i) = \sigma(j)$. If there are several such pairs, select
the one for which the value $\pi(i) = \sigma(j)$ is minimal. Assume
first that $i < s_j$, and let $s_h$ be the smallest element of $S$
larger than $i$.  Observe that $h \leq j$ since $i<s_{j}$.  Then, set
$\iota (S,\pi,\sigma)=(S',\pi',\sigma')$ where
\begin{itemize}
 \item $S' = S \setminus \{s_j\} \cup \{i\}$,
 \item $\pi' = \pi \cdot (i,s_h,s_{h+1},\ldots,s_j)$, 
 \item $\sigma' = \sigma \cdot (j,j-1,\ldots,h)$.
\end{itemize}
On the other hand, if $s_j < i$, assume that $s_h$ is the largest
element of $S$ smaller than $i$.  Note $h \geq j$ since $s_{j}<i$.
Then, set 
$\iota (S,\pi,\sigma)=(S',\pi',\sigma')$, where 
\begin{itemize}
 \item $S' = S \setminus \{s_j\} \cup \{i\}$,
 \item $\pi' = \pi \cdot (i,s_h,s_{h-1},\ldots,s_j)$,
 \item $\sigma' = \sigma \cdot (j,j+1,\ldots,h)$.
\end{itemize}

In the case $i<s_{j}$, we note $\sigma'(h)=\sigma (j)=\pi
(i)=\pi'(s_{j})$ is minimal among $\{\sigma'(1),\ldots,\sigma'(r)\}$
intersected with the complement of $\pi'(S')$, and furthermore $i$ is
the $h^{th}$ largest element of $S'$ by choice of $h$.  Thus when
computing $\iota(S',\pi', \sigma')$ the value $i$ and $s_{j}$ have
swapped roles.  Similarly, in the case $i>s_{j}$ we also observe that
$\sigma'(h)=\sigma (j)=\pi (i)=\pi'(s_{j})$ is minimal among
$\{\sigma'(1),\ldots,\sigma'(r)\}$ intersected with the complement of
$\pi'(S')$ and $i$ is the $h^{th}$ largest element of $S'$ by choice
of $h$, so again $i$ and $s_{j}$ have swapped roles.  By permutation
multiplication, one can verify $\pi''=\pi$,\ $\sigma''=\sigma$, and
the weight of the $(S,\pi,\sigma)$ and $(S',\pi',\sigma')$ terms are
exactly the same, though they appear in a permuted order.  
Furthermore, by construction of $\iota$ 
\[
(-1)^{\inv(\pi) + \inv(\sigma)}=-(-1)^{\inv(\pi') + \inv(\sigma')}.  
\]
Thus, $\iota$ is a weight-preserving, sign-reversing involution, and hence the right side of \eqref{eq:sylvester} reduces
to a signed sum of over the $(n!)^{2}$ fixed points of $\iota$.

A fixed point of the involution $(S, \pi ,\sigma)$ can be mapped
bijectively to a term from the left side \eqref{eq:sylvester} indexed
by of the same sign and weight by a pair $(\rho,\tau) \in S_n^2$, as
follows.  Define $\rho , \tau$ by
\begin{equation}\label{eq:rho}
\rho(i) = \begin{cases}
\pi(i)&  i \notin S,\\
\sigma(j) & i = s_j\in S 
\end{cases}
\end{equation}
and
\begin{equation}\label{eq:tao}
\tau(i) = \begin{cases}
\pi(s_i) &   i \leq r \\
\sigma(i)&  i > r
\end{cases}
\end{equation}
for all $1\leq i\leq n$.   The set $S=\{s_{1}<\dots <s_{r} \}$ can be
recovered from $(\rho,\tau)$ by considering the positions of the
values $\{\tau(1),\tau(2),\dotsc , \tau (r) \}$ in $\rho$, so $\pi$
and $\sigma$ can also be recovered.  Therefore,  the map is an injection
between two sets of the same size, hence a bijection.
The pair $(\rho,\tau) \in S_n^2$ is specifically constructed so that
the weight of $(S,\pi ,\sigma )$ given in \eqref{eq:weight.s.pi.sigma}
is the same as the weight of the $(\rho,\tau)$ term 
\[
 \prod_{i \in [n]}a_{\rho(i),i} b_{\tau (i),i}
\]
after commuting the factors.  It remains to prove this bijection
preserves the sign.  We leave this as an exercise for the reader, see
\Cref{ex:sign.reversing.map}.
\end{proof}

\begin{Example}\label{example:Sylvester}
Let $n=6, \ r=4$. Assume $A=(a_{i,j}),B=(b_{i,j})$ are $6 \times 6$
matrices. If $S=\{1,3,4,6 \}$,\ $\pi =513264$, \ and $\sigma =245136$,
then the triple $(S,\pi ,\sigma )$ indexes the term on the right side
of \eqref{eq:sylvester} given by the determinant of the product of
matrices
\begin{equation}\label{eq:syl.1}
M_{(S,\pi ,\sigma )}=\begin{pmatrix}
0 & a_{1,2} & 0 & 0 & 0 & 0\\
0 & 0 & 0 & b_{2,3} & 0 & 0\\
0 & 0 & b_{3,2} & 0 & 0 & 0\\
0 & 0 & 0 & 0 & 0 & b_{4,4}\\
b_{5,1} & 0 & 0 & 0 & 0 & 0\\
0 & 0 & 0 & 0 & a_{6,5} & 0
\end{pmatrix}
\begin{pmatrix}
0 & 0 & 0 & a_{1,6} & 0 & 0\\
a_{2,1} & 0 & 0 & 0 & 0 & 0\\
0 & 0 & 0 & 0 & b_{3,5} & 0\\
0 & a_{4,3} & 0 & 0 & 0 & 0\\
0 & 0 & a_{5,4} & 0 & 0 & 0\\
0 & 0 & 0 & 0 & 0 & b_{6,6}
\end{pmatrix}. 
\end{equation}
Note, the columns are determined by $A_{S},B_{S}$ and the nonzero
entries in each column are determined by $\pi,\sigma$ respectively.
To compute $\iota(S,\pi, \sigma)$, observe $\pi (S)=\{5,3,2,4 \}$ and
$\sigma([4])=\{2,4,5,1 \}$ are not equal.  The smallest value in
$\sigma([4])$ not in $\pi (S)$ is $1=\sigma (4)=\pi (2)$, so $i=2$ and
$j=4$.  Sorting $\pi (S)$ we get $s_{1}=2, s_{2}=3, s_{3}=4,s_{4}=6$,
so $s_{j}=6$.  Since $i<j$, we observe $s_{h}=3$ is the smallest value
of $S$ that is bigger than $i=2$.  Now we have all of the data to
compute $\iota (S,\pi ,\sigma )= (S',\pi' ,\sigma' )$ with
\begin{align}\label{eq:Syl.example}
S'&=S \setminus \{s_{j}\} \cup \{i \}=S \setminus \{6\} \cup \{2 \}=\{1,2,3,4\}\\
\pi' &= \pi (2,3,4,6) =532461\\
\sigma' &= \sigma (4,3,2) = 214536.
\end{align}
The term corresponding with $(S',\pi' ,\sigma ')$ is the determinant
of the product 
\begin{equation}\label{eq:syl.2}
M_{(S',\pi' ,\sigma' )}=\begin{pmatrix}
0 & 0 &  0 & 0 & 0 & a_{1,6}\\
0 & 0 &  b_{2,3} & 0 & 0 &0\\
0 & b_{3,2}   & 0 & 0 & 0 & 0\\
0 & 0 & 0  &  b_{4,4} & 0 &0\\
b_{5,1} & 0  & 0 & 0 & 0 & 0\\
0 & 0 & 0 & 0  & a_{6,5} & 0
\end{pmatrix}
\begin{pmatrix}
0 & a_{1,2}&  0 & 0 & 0 & 0\\
a_{2,1} & 0 & 0 & 0 & 0 & 0\\
0 & 0 & 0 & 0 &  b_{3,5} & 0\\
0 & 0 & a_{4,3}  & 0 & 0 & 0\\
0 & 0 & 0 & a_{5,4}  & 0 & 0\\
0 & 0 & 0 & 0 & 0 & b_{6,6}
\end{pmatrix}. 
\end{equation}
It is straightforward to observe $\det(M_{(S,\pi ,\sigma )}) =
-\det(M_{(S',\pi' ,\sigma' )})$, as expected since $\iota$ is a weight
preserving, sign reversing involution.  We leave it to the reader to
verify that $\iota (S',\pi', \sigma')=(S,\pi,\sigma)$ by confirming
$\pi'(6)=1=\sigma'(2)$ implies $2,6$ swap roles in the algorithm.  
\end{Example}

\bigskip

\begin{Exercise}\label{ex:sign.reversing.map}
Let $(S,\pi,\sigma)$ be a fixed point of $\iota$ with $S=\{s_{1}<\dots
<s_{r} \} \subset [n]$, and $\pi,\sigma \in S_{n}$.  If $(\rho,\tau)$
are given by \eqref{eq:rho} and \eqref{eq:tao}, then
\[
(-1)^{\inv(\pi) + \inv(\sigma)} = (-1)^{\inv(\rho) + \inv(\tau)} .
\]
\end{Exercise}

Similarly, the Pl\"ucker coordinates of two subspaces $V \in
\Gr (k,n)$ and $W \in \Gr (\ell,n)$ with $k<\ell$ satisfy certain
\textit{incidence relations} if and only if $V \subset W$.  These
incidence relations can be verified by considering Gauss elimination
on a matrix with columns consisting of a union of a basis for $V$ and
a basis for $W$.  Such a matrix has rank equal to $\dim(W)$ if and
only if $V \subset W$. These rank conditions in turn can be encoded by
checking that all size $\ell+1$ minors vanish.  Such minors will
always vanish if $V,W$ are two nested subspaces in a flag
$F_{\bullet}$, hence the Pl\"ucker coordinates $P(F_{\bullet})$ must
satisfy these additional incidence relations.  It is remarkable that
these relations can also be encoded as identities among flag minors in
a very similar way to Sylvester's Lemma, and that these identities
suffice to generate all relations on the Pl\"ucker coordinates of
flags as stated in the following theorem.  The full proof of this
theorem is established in \cite[\S 8.1-9.1]{Fulton-book} using some
essential ingredients from $GL_{n}$ representation theory and
invariant theory.  We encourage the reader to work through those
details.

\begin{Theorem}[\textbf{Pl\"ucker relations}]\label{thm:plucker}
For any $n\times n$ matrix $M,$ the flag minors satisfy the relations
\[
\Delta_{I}(M) \Delta_{J}(M) = \sum \Delta_{I'}(M) \Delta_{J'}(M)
\]
where the sum is over all pairs $(I',J')$ obtained from $(I,J)$ by
swapping a fixed subset of $J$ with every subset of $I'$ of the same
size.  Furthermore, for every point $p=(p_{I}: I \subsetneq [n])\in
\prod_{k=1}^{n-1} \mathbb{P}^{\binom{n}{k}-1}$ such that 
\[
p_{I} p_{J} = \sum p_{I'} p_{J'},
\]
where the sum is over all pairs $(I',J')$ as above,  there exists a flag $F_{\bullet} \in
\Fl(n)$ such that $p=P(F_{\bullet})$.  Therefore, these relations over
all subsets $I,J \subset [n]$ determine the image of $\Fl(n)$ in
$\prod_{k=1}^{n-1} \mathbb{P}^{\binom{n}{k}-1}$ as a projective
subvariety.
\end{Theorem}

\begin{Example}
Up to rearrangement of terms, there is just one nontrivial Pl\"ucker
relation using size $2$ subsets of $[4]$, namely
\begin{equation} \label{eq:plucker}
\Delta_{1,2}\Delta_{3,4} - \Delta_{1,3}\Delta_{2,4} +
\Delta_{1,4}\Delta_{2,3} = 0.
\end{equation}
That is, a point $[p_{12} : p_{13} : p_{14} : p_{23} : p_{24} :
p_{34}]$ in $\bP^5$ is the list of Pl\"ucker coordinates of a
$2$-plane in $\bC^4$ (or $F^4$ for any field $F$) if and only if
\eqref{eq:plucker} holds.
\end{Example}

The Pl\"ucker relations are quadratic relations in the freely
generated polynomial ring $\mathcal{P}_{n}$ on the symbols
$\Delta_{K}$'s for $K \subset [n]$.  They give rise to a family of
(nonzero) polynomials
\[
f_{I,J} =\Delta_{I} \Delta_{J} - \sum \Delta_{I'} \Delta_{J'} \in \mathcal{P}_{n}
\]
 for each pair of proper subsets $I,J \subset
[n]$.  Each $f_{I,J}$ is \textit{multihomogenous} in the sense that
for each $j$, every term of $f_{I,J}$ has the same degree in the
variables $\Delta_{K}$ as $K$ varies over the size $j$ subsets.
Therefore, $f_{I,J}$ is a well defined function on the product of
projective spaces $\prod_{k=1}^{n-1} \mathbb{P}^{\binom{n}{k}-1}$
where $\Delta_{K}$ is the coordinate function corresponding to $K$.
Multihomogeneity is an essential part of the following corollary of
\Cref{thm:plucker}.

\begin{Corollary}\label{cor:fln.projective.variety}\cite[Prop. 1, \S9.1]{Fulton-book}
The flag variety (and each partial flag variety) is an irreducible
projective variety. The equations derived from the Pl\"ucker relations
generate the prime ideal of all polynomial functions on
$\prod_{k=1}^{n-1} \mathbb{P}^{\binom{n}{k}-1}$ that vanish on the
image of $P(\Fl(n))$ in $\prod_{k=1}^{n-1} \mathbb{P}^{\binom{n}{k}-1}$.
\end{Corollary}

The ideal of quadratic polynomials coming from the Pl\"ucker relations
is central in both the representation theory of $GL_{n}$ and the study
of Schubert varieties.  Fulton's book on Young tableaux highlights
these relations in different contexts.  We state the following
important consequence and refer the interested readers to
\cite{Fulton-book} for more information.

Recall that $\rk(w)[i,j]$ is the rank of the submatrix of $M_w$ weakly
northwest of $(i,j)$.  The \emph{southwest rank table} of a matrix
$A$, whose $(i,j)$ entry is the rank of the submatrix of $A$ weakly
southwest of $(i,j)$, will be denoted as
$\rk_{\mathrm{SW}}(A)[i,j]$. We write $\rk_{\mathrm{SW}}(w)[i,j]$ for
$\rk_{\mathrm{SW}}(M_w)[i,j]$ as before.

\begin{Corollary}\label{cor:Schubert.variety.equations}
For $w \in S_{n}$, the vanishing minors of $C_{w}(E_{\bullet})$ are
determined by the southwest rank table for $w$. Therefore, the
Schubert variety $X_{w}(E_{\bullet})$ is isomorphic to the projective
variety determined by the Pl\"ucker relations, the incidence
relations, plus the vanishing minors $\Delta_{I,J}$ on rows $I$,
columns $J$ for all $1 \leq i,j \leq n$ and pairs of subsets $I
\subset \{i+1,i+2,\dots ,n \}$ and $J \subset [j]$ of size
$\rk_{\mathrm{SW}}(w)[i,j]+1$.
\end{Corollary}

\begin{Exercise}\label{ex:Schubert.}
Identify the Schubert variety $X_{2314}$ with a subvariety of a
product of projective spaces by specifying all of the necessary
multihomogenous equations.
\end{Exercise}

For the reader familiar with exterior algebras, recall that a $k$-plane $V = \langle v_1, \ldots, v_k \rangle$ uniquely determines a line $\bigwedge^k(V) = \C \cdot v_1 \wedge \cdots \wedge v_k \subseteq \bigwedge^k(\C^n)$, i.e. a point in the projective space $\bP(\bigwedge^k(\C^n)) \simeq \bP^{{n \choose k}-1}$. The Pl\"ucker coordinates of $V$ are  the projective coordinates of $\bigwedge^k(V)$ with respect to the basis $\{e_{i_1} \wedge \cdots \wedge e_{i_k} \given 1 \leq i_1 < \cdots < i_k \leq n\}$ of $\bigwedge^k(\C^n)$. However, a typical element of $\bigwedge^k(\C^n)$ cannot be written as a single pure tensor $v_1 \wedge \cdots \wedge v_k$, only a sum of several, in which case it does not represent any $k$-plane. The subset of pure tensors in $\bigwedge^k(\C^n)$ is exactly the subset where the Pl\"ucker relations hold. For instance, $\alpha = \sum_{1 \leq i < j \leq 4} \alpha_{ij} e_i \wedge e_j$ represents a $2$-plane in $\C^4$ if and only if the relation \eqref{eq:plucker} holds, i.e. 
\begin{equation} \label{eq:plucker2}
    \alpha_{12}\alpha_{34} - \alpha_{13}\alpha_{24} + \alpha_{14}\alpha_{23} = 0.
\end{equation}

\begin{Exercise} \label{exer:4-lines-3-lines}
    \hfill
    \begin{enumerate}[(a)]
        \item If $P$ is a line in $\C^3$, let $\lambda(P) \in
	\bigwedge^2(\C^4)$ denote any element representing the
	$2$-plane passing through $P$ and the origin in $\C^4$. Verify that $P \cap Q \neq \emptyset$ if and only if $\lambda(P) \wedge \lambda(Q) = 0$. Also show that $\alpha \in \bigwedge^2(\C^4)$ satisfies the Pl\"ucker relation \eqref{eq:plucker2} if and only if $\alpha \wedge \alpha = 0$. 
        \item Show that if $K \subseteq \bigwedge^2(\C^4)$ is a linear subspace of dimension at least $3$, then there are infinitely many distinct $2$-planes $V \in \Gr(2,4)$ represented by pure tensors in $K$.
        \item Deduce that if $P_1, P_2, P_3, P_4$ and $Q_1, Q_2, Q_3$ are distinct lines in $\C^3$ such that $Q_i \cap P_j \neq 0$ for all $i,j$, then there are infinitely many lines intersecting all 4 of $P_1, \ldots, P_4$.
        \item Show by example that if $\C$ is replaced by $\bR$, then (b) need not hold. Prove nevertheless that (c) still holds over $\bR$.
    \end{enumerate}
\end{Exercise}

\subsection{Matrix Schubert Varieties}\label{sub:MatrixSchubs}

Although we have seen some of the elegant combinatorics and algebra
surrounding pipe dreams, the cohomology class $[X_w]$ could be
represented by any polynomial congruent to $\fS_{w_{0}w}$ modulo the
ideal of positive degree invariants $I_{n}^{+}$.  Why choose this
representative?  One good reason is that it expands as a positive sum
of terms.  In this subsection, we describe a second good reason for
this choice. We outline a specific geometric interpretation for the
pipe dreams that clarifies their relationship to the geometry.  Kogan
\cite{kogan.phd} was the first to find such a connection; he showed
that certain toric varieties are determined by pipe dreams. Here we
describe the related approach of Knutson and Miller
\cite{knutson-miller-2005}. See also \cite{miller-sturmfels} for a
more detailed exposition of the material sketched in this section.

In order to describe the geometric interpretation for pipe dreams, we
consider affine varieties that are inspired by Schubert varieties, but
which contain slightly more information coming from including
non-invertible matrices. We will also replace the Chow ring with the
\emph{equivariant Chow ring}, which also records extra information
coming from a torus action. This machinery will produce a Schubert
polynomial exactly, rather than just a polynomial modulo an
ideal. Moreover, the terms in the monomial decomposition of the
Schubert polynomial will correspond to pieces in a deformation of the
variety into a union of linear subspaces.

Let $\Mat_n(S)$ be the space of $n \times n$ matrices with entries in
the set $S$. Given $A \in \Mat_n(\C)$, write $A_i$ for the span of the
leftmost $i$ columns of $A$.  Note, $A_\bullet = (A_{1}\subset \cdots
\subset A_{n}) \in \Fl_n(\C)$ if $A$ is invertible, and $A_\bullet$ is
still a nested list of subspaces even if $A$ is not invertible.  Let
$E_\bullet$ be the ``standard flag'' with $E_i = \mathrm{span} \langle
e_1, \ldots, e_i \rangle$.

Recall, we used the southwest rank table associated to $w$ to give
equations for $X_{w}$ in \Cref{cor:Schubert.variety.equations}.
Similar rank conditions are used to define an affine analog of
Schubert varieties.

\begin{Definition} \label{defn:matrix-schub}
    The \emph{matrix Schubert variety} of $w \in S_n$ is
    \begin{equation}\label{eq:matrix-schub}
        \MX_w = \{A \in \Mat_n(\C) \given \rk_{\mathrm{SW}}(A)[i,j] \leq \rk_{\mathrm{SW}}(w)[i,j] \text{ for all $1 \leq i,j \leq n$}\}.
    \end{equation}
\end{Definition}

\begin{Exercise}
Prove that
\[
\MX_w \cap \GL_n(\C)= \{A \in \Mat_n(\C) \given \dim(E_i \cap A_j)
\geq \rk(w)[i,j] \text{ for all $1 \leq i,j \leq n$}\}, 
\]
and show that $\MX_w \cap \GL_n(\C)/B= X_{w} \subset G/B$.  
\end{Exercise}

The rank conditions in \Cref{defn:matrix-schub} are equivalent to the
vanishing of all size $\rk_{\mathrm{SW}}(w)[i,j]+1$ minors in the
submatrix of $A$ southwest of $(i,j)$, so $\MX_w$ is indeed an affine
variety in $\Mat_n(\C) \simeq \C^{n \times n}$. In fact, something
stronger is true. Let $I(\MX_w)$ be the defining ideal in $\C[z_{11},
z_{12}, \ldots, z_{nn}]$ of all polynomials vanishing on $\MX_w$.  The
ideal of polynomials vanishing on an irreducible affine variety is always prime
\cite{Cox-Little-OShea}.  Let $[z_{pq}]$ denote the $n
\times n$ matrix of variables.

\begin{Theorem} \label{thm:MX-ideal} \cite[Prop. 3.3]{Fulton1} The
ideal $I(\MX_w)$ is generated by the $\rk_{\mathrm{SW}}(w)[i,j]+1$
minors of $[z_{pq}]_{i \leq p \leq n, 1 \leq q \leq j}$ for all $i,j$.
\end{Theorem}
Most of the rank conditions just described are redundant: for
instance, the conditions with $i = n$ or $j = n$ are completely
vacuous. Fulton identified the minimal set of irredundant rank
conditions using the essential set in the diagram of a permutation
\cite{Fulton1}.  Recall, the diagram $D(w)$ is based on the matrix for
$w^{-1}$ by \eqref{eq:diagram}.  Also let $D_{\mathrm{SW}}(w)$ be the
\emph{southwest Rothe diagram} of $w$: the subset of $[n]^2$ obtained
by removing all cells which are weakly right of or \emph{above} a 1 in
the matrix of $w^{-1}$.

\begin{Definition}\cite{Fulton1}\label{def:essential.set}
    The \emph{essential set}  of $w$ is the set of
    northeasternmost elements of the connected components of
    $D_{\mathrm{SW}}(w^{-1})$, denoted in  using matrix coordinates as 
\[
\Ess(w) = \{(i,j) \in  D_{\mathrm{SW}}(w^{-1}) \given (i-1,j),(i,j+1) \notin
    D_{\mathrm{SW}}(w^{-1}) \}.
 \]

\end{Definition}

\begin{Exercise} \label{exer:essential} Show that
the conditions defining $MX_{w}$ appearing in Equation
\eqref{eq:matrix-schub} still hold if we only let $(i,j)$ range over
$\Ess(w)$.
\end{Exercise}

\begin{Example}
Recall we saw the permutation diagrams $D(1437256)$ and $D(1437265)$
in \Cref{fig:v.w.diagrams}.  Compare those (northwest) diagrams with
the southwest diagrams shown in \Cref{fig:essential-sets-examples} for
the same two permutations.
The essential sets of $v=1532674=(1437256)^{-1}$ and
$w=1532764=(1437265)^{-1}$ are the cells shaded green in these
diagrams.  
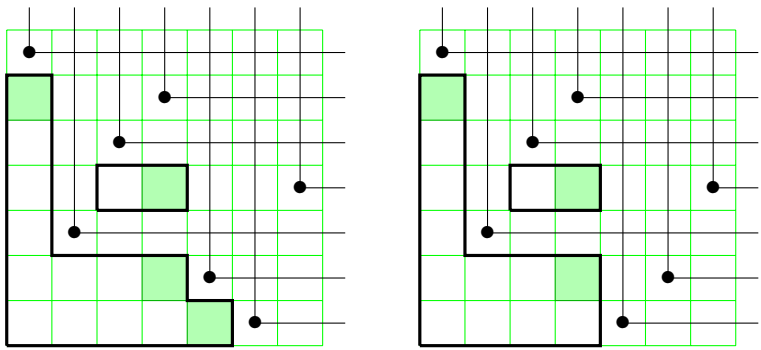
\begin{figure}[h!]
\centering
\begin{tikzpicture}[scale=0.6]
\draw[step=1.0,green,thin] (0,0) grid (7,7);
\draw[fill=green, opacity=0.3] (0,5)--(1,5)--(1,6)--(0,6)--(0,5);
\draw[fill=green, opacity=0.3] (3,3)--(4,3)--(4,4)--(3,4)--(3,3);
\draw[fill=green, opacity=0.3] (3,1)--(4,1)--(4,2)--(3,2)--(3,1);
\draw[fill=green, opacity=0.3] (4,0)--(5,0)--(5,1)--(4,1)--(4,0);
\draw[very thick] (0,0)--(0,6)--(1,6)--(1,2)--(4,2)--(4,1)--(5,1)--(5,0)--(0,0);
\draw[very thick] (2,3)--(4,3)--(4,4)--(2,4)--(2,3);
\node at (0.5,6.5) {$\bullet$};
\node at (1.5,2.5) {$\bullet$};
\node at (2.5,4.5) {$\bullet$};
\node at (3.5,5.5) {$\bullet$};
\node at (4.5,1.5) {$\bullet$};
\node at (5.5,0.5) {$\bullet$}; 
\node at (6.5,3.5) {$\bullet$};
\draw(0.5,7.5)--(0.5,6.5)--(7.5,6.5);
\draw(1.5,7.5)--(1.5,2.5)--(7.5,2.5);
\draw(2.5,7.5)--(2.5,4.5)--(7.5,4.5);
\draw(3.5,7.5)--(3.5,5.5)--(7.5,5.5);
\draw(4.5,7.5)--(4.5,1.5)--(7.5,1.5);
\draw(5.5,7.5)--(5.5,0.5)--(7.5,0.5);
\draw(6.5,7.5)--(6.5,3.5)--(7.5,3.5);
\end{tikzpicture}
\qquad
\begin{tikzpicture}[scale=0.6]
\draw[step=1.0,green,thin] (0,0) grid (7,7);
\draw[fill=green, opacity=0.3] (0,5)--(1,5)--(1,6)--(0,6)--(0,5);
\draw[fill=green, opacity=0.3] (3,3)--(4,3)--(4,4)--(3,4)--(3,3);
\draw[fill=green, opacity=0.3] (3,1)--(4,1)--(4,2)--(3,2)--(3,1);
\draw[very thick] (0,0)--(0,6)--(1,6)--(1,2)--(4,2)--(4,0)--(0,0);
\draw[very thick] (2,3)--(4,3)--(4,4)--(2,4)--(2,3);
\node at (0.5,6.5) {$\bullet$};
\node at (1.5,2.5) {$\bullet$};
\node at (2.5,4.5) {$\bullet$};
\node at (3.5,5.5) {$\bullet$};
\node at (4.5,0.5) {$\bullet$};
\node at (5.5,1.5) {$\bullet$};
\node at (6.5,3.5) {$\bullet$};
\draw(0.5,7.5)--(0.5,6.5)--(7.5,6.5);
\draw(1.5,7.5)--(1.5,2.5)--(7.5,2.5);
\draw(2.5,7.5)--(2.5,4.5)--(7.5,4.5);
\draw(3.5,7.5)--(3.5,5.5)--(7.5,5.5);
\draw(4.5,7.5)--(4.5,0.5)--(7.5,0.5);
\draw(5.5,7.5)--(5.5,1.5)--(7.5,1.5);
\draw(6.5,7.5)--(6.5,3.5)--(7.5,3.5);
\end{tikzpicture}
\caption{The essential sets (shaded) of $v=1532674=(1437256)^{-1}$ and $w=1532764=(1437265)^{-1}$}
\label{fig:essential-sets-examples}
\end{figure}
\end{Example}

\begin{Theorem}[\textbf{Fulton's Essential Set Theorem}]\label{thm:fulton.essential.set} The ideal
$I(\MX_w)$ is generated by the size $\rk_{\mathrm{SW}}(w)[i,j]+1$
minors of $[z_{pq}]_{i \leq p \leq n, 1 \leq q \leq j}$ for all $(i,j)
\in \Ess(w)$.  Furthermore, no subset of these essential set rank
equations determines $I(\MX_w)$.
\end{Theorem}

\Cref{exer:essential} asks the reader to prove the first statement in
this theorem.  The second statement is shown by an explicit
construction.  See \Cref{exer:minimal} below.  See \cite[Lemmas 3.10
and 3.14]{Fulton1} for more details. 

\begin{Exercise}\label{exer:minimal}
For $w \in S_{n}$ and $(p,q) \in \Ess(w)$, let $A=(a_{ij})\in
\Mat_{n}(\C)$ be defined by the partial permutation matrix with
$a_{ij}=1$ if $w(j)=i$, $i\geq p$, and $j\leq q$, and $a_{ij}=0$
otherwise. Let $A'$ be obtained from $A$ by additionally setting
$a_{p,q}=1$. Prove that $A'$ satisfies $\rk_{\mathrm{SW}}(A)[p,q] =
\rk_{\mathrm{SW}}(w)[p,q]+1$ and $\rk_{\mathrm{SW}}(A)[i,j] \leq
\rk_{\mathrm{SW}}(w)[i,j]$ for all $(i,j) \in \Ess(w) \setminus
\{(p,q) \}$.
\end{Exercise}

\begin{figure}
\begin{center}
\begin{tikzpicture}[scale=0.6]
\draw[step=1.0,black,thin] (0,0) grid (4,4);
\node at (0.5,0.5) {$0$};
\node at (0.5,1.5) {$1$};
\node at (0.5,2.5) {$1$};
\node at (0.5,3.5) {$1$};
\node at (1.5,0.5) {$0$};
\node at (1.5,1.5) {$1$};
\node at (1.5,2.5) {$1$};
\node at (1.5,3.5) {$2$};
\node at (2.5,0.5) {$0$};
\node at (2.5,1.5) {$1$};
\node at (2.5,2.5) {$2$};
\node at (2.5,3.5) {$3$};
\node at (3.5,0.5) {$1$};
\node at (3.5,1.5) {$2$};
\node at (3.5,2.5) {$3$};
\node at (3.5,3.5) {$4$};
\end{tikzpicture}
\qquad
\begin{tikzpicture}[scale=0.6]
\draw[step=1.0,green,thin] (0,0) grid (4,4);
\node at (0.5,1.5) {$\bullet$};
\node at (1.5,3.5) {$\bullet$};
\node at (2.5,2.5) {$\bullet$};
\node at (3.5,0.5) {$\bullet$};
\draw(0.5,4.5)--(0.5,1.5)--(4.5,1.5);
\draw(1.5,4.5)--(1.5,3.5)--(4.5,3.5);
\draw(2.5,4.5)--(2.5,2.5)--(4.5,2.5);
\draw(3.5,4.5)--(3.5,0.5)--(4.5,0.5);
\draw[fill=green, opacity=0.3] (1,2)--(2,2)--(2,3)--(1,3)--(1,2);
\draw[fill=green, opacity=0.3] (2,0)--(3,0)--(3,1)--(2,1)--(2,0);
\draw[very thick] (0,0)--(3,0)--(3,1)--(0,1)--(0,0);
\draw[very thick] (1,2)--(2,2)--(2,3)--(1,3)--(1,2);
\end{tikzpicture}
\end{center}
\caption{The southwest rank table and essential set for $3124$.} \label{fig:3124.sw}
\end{figure}
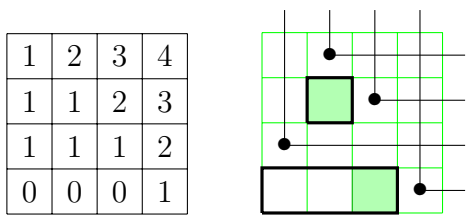

\begin{Example} \label{ex:matrix-schub} Consider a smaller example for
the purpose of explicit calculation. The southwest rank table of
$3124=(2314)^{-1}$ and $D_{SW}(2314)$ are shown in \Cref{fig:3124.sw}.
Therefore, the green shaded cells form $\Ess(3124)=\{(2,2),(4,3) \}$.
%\begin{equation*}
%\begin{array}{cccc}
%1 & 2 & 3 & 4 \\
%1 & 1 & 2 & 3 \\
%1 & 1 & 1 & 2\\
%0 & 0 & 0 & 1
%\end{array} \qquad 
%\begin{array}{cccc}
%    \cdot & 1 & \cdot & \cdot\\
%    \cdot & \blacksquare & 1 & \cdot\\
%    1 & \cdot & \cdot & \cdot\\
%    \square & \square & \blacksquare & 1 
%    \end{array}
%\end{equation*}
% The positions in $\Ess(3124)$ have been shaded.
Hence $\MX_{3124}$ is the set of $4\times 4$ matrices
% $\left[ \begin{smallmatrix} z_{11} & \cdots & z_{14} \\
%  \vdots & \ddots & \vdots\\
%  z_{41} & \cdots & z_{44} \end{smallmatrix}\right]$ with 
% \begin{equation*} \rank \begin{bmatrix} z_{41} & z_{42} & z_{43} \end{bmatrix} = 0 \quad \text{and} \quad \rank \begin{bmatrix} z_{21} & z_{22} \\ z_{31} & z_{32} \\ z_{41} & z_{42} \end{bmatrix} \leq 1,
% \end{equation*}
$(a_{ij})_{1\leq i,j\leq 4}$ with 
\begin{equation*} \rank \begin{bmatrix} a_{41} & a_{42} & a_{43} \end{bmatrix} = 0 \quad \text{and} \quad \rank \begin{bmatrix} a_{21} & a_{22} \\ a_{31} & a_{32} \\ a_{41} & a_{42} \end{bmatrix} \leq 1, 
\end{equation*}
or equivalently the vanishing set for the equations $z_{41} = z_{42} = z_{43} = z_{21}z_{32} - z_{22}z_{31} = 0$.
\end{Example}

Eriksson-Linusson showed the average size of the essential set for $w
\in S_{n}$ is asymptotically $(1/36)n^{2}$ \cite{EL.1996}, less than the
average size of the diagram of a permutation which is $\frac{1}{2}\binom{n
}{2}$.  However, \Cref{thm:fulton.essential.set} does \emph{not}
necessarily give a minimal generating set for $I(MX_w)$: indeed, two
of the three $2 \times 2$ minors in Example~\ref{ex:matrix-schub} are
redundant. S. Gao and Yong described a subset of the essential minors
which do form a minimal generating set \cite{gao-yong}.

Next we introduce some vocabulary from representation theory that
applies to all affine varieties with a torus action, including matrix
Schubert varieties.  See \cite{fulton-harris,Hum-LAG} for more
detailed background.  A nice overview of the representation theory
related to Grassmannians, Gr\"obner basis, and Schubert varieties
appears in \cite[Sect. 10]{Michalek-Sturmfels}.

Suppose $V$ is a vector space with an action of the complex torus $T
\simeq (\C^\times)^n$. Fixing an isomorphism $t \mapsto (t_1,…, t_n)$, the
\emph{weight space} of ${\bf a} \in \bZ^n$ in $V$ is
\begin{equation*}
    V_{\bf a} = \{v \in V \given tv = t_1^{a_1} \cdots t_n^{a_n}v \text{ for all $t \in T$}\}.
\end{equation*}
The \emph{character} of $V$ is then $\chr(V)=\sum_{{\bf a}
\in \bZ^n} \dim(V_{\bf a})x_1^{a_1} \cdots x_n^{a_n}$ viewed as a
formal power series in indeterminates $x_1^{\pm 1}, \ldots, x_n^{\pm
n}$, i.e. the formal trace of the action of $t = (t_1, \ldots, t_n)
\in T$ as a function on $V$.  We call an element $v \in V_{\bf a}$ a
\emph{weight vector} of \emph{weight} $\wt(v) = {\bf a}$.
We will only be concerned with
representations nice enough that all $\dim(V_{\bf a})$ are finite and
$V = \bigoplus_{{\bf a} \in \bZ^n} V_{\bf a}$, so the character is
well-defined and determines $V$. If we define the \emph{exponential
weight} $\expwt(v) = x_1^{\wt(v)_1} \cdots x_n^{\wt(v)_n}$, then the
character can be written as $\sum_{v} \expwt(v)$, where $v$ runs over
a basis of $V$ consisting of weight vectors.

Now assume $V$ is finite-dimensional, and let $I(Z)$ be the
\emph{ideal of an affine variety} $Z \subseteq V$, which is defined by 
\begin{equation*}
    I(Z) = \{f \in \C[V] \given f|_{Z} = 0\}
\end{equation*}
where $\C[V] \simeq \Sym^*(V^*)$ is the ring of polynomial functions
on $V$. The ring of polynomial functions on $Z$ is $\C[V]/I(Z)$.  If
$Z$ is $T$-stable, then the corresponding $T$-action on $\C[V]/I(Z)$
determines another representation of $T$. Its character is the
\emph{multigraded Hilbert series} of $I(Z)$ (or just of $Z$), denoted
\[
\mathcal{H}(I(Z)) = \chr(\C[V]/I(Z)).
\]

\begin{Definition}\label{def:k-poly}
The \emph{K-polynomial} of $I(Z)$ (or just $Z$) is $\mathcal{K}(I(Z)) = \mathcal{H}(I(Z)) / \mathcal{H}(0)$, and the \emph{multidegree} $\mathcal{M}(I(Z))$ is the lowest degree component of the power series resulting from making the substitutions $t_j \mapsto 1-t_j$ in $\mathcal{K}(I(Z))$.
\end{Definition}

\begin{Example} \label{ex:hilbert}
    Identify $(t_1, \ldots, t_n) \in T = (\C^\times)^n$ with the
$n\times n$ diagonal matrix with diagonal $t_1, \ldots, t_n$, and let
$T \times T$ act on $V=\Mat_n(\C)$ on both sides by the formula $(s,t)
\cdot A = t^{-1} As$. Write $\C[\Mat_n(\C)] \simeq \C[z_{11}, \ldots,
z_{nn}]=\C[z_{ij}]$. The induced $T \times T$-action on $\C[z_{ij}]$
is then given by $(s,t)\cdot z_{ij} = s_j t_i^{-1} z_{ij}$, so
monomials in the $z_{ij}$'s give a basis of weight vectors for
$\C[\Mat_n(\C)]$ with weights in $\bZ^{n\times n}$, which could also
be written as $\Mat_{n}(\bZ_{\geq 0})$. If $Z=V$, then 
$I(Z)=(0)$ and the multigraded Hilbert series is

    \begin{equation*}
        \mathcal{H}(0) = \sum_{{\bf e} \in \Mat_n(\bZ^{\geq 0})}  \expwt\left(\prod_{i,j=1}^n z_{ij}^{e_{ij}}\right) = \sum_{{\bf e}} \prod_{i,j=1}^n x_j^{e_{ij}}y_i^{-e_{ij}} = \prod_{i,j=1}^n \frac{1}{1-x_j y_i^{-1}},
    \end{equation*}
    where we have used $x_1, \ldots, x_n, y_1, \ldots, y_n$ as the indeterminates in our exponential weights.
    
    For a slightly more interesting example, let $Y = \{A \in \Mat_2(\C) \given \det A = 0\}$, so $I(Y) = (z_{11}z_{22} - z_{12}z_{21})$. As a basis of $\C[z_{ij}] / I(Y)$ we may take all monomials in $z_{11}, z_{12}, z_{21}, z_{22}$ not divisible by $z_{11}z_{22}$, so $\mathcal{H}(I(Y))$ is 
    \begin{equation*}
         \mathcal{H}(0) - \mathcal{H}(0)\expwt(z_{11}z_{22}) = \frac{1 - x_1 x_2 y_1^{-1} y_2^{-1}}{(1-x_1y_1^{-1})(1-x_1 y_2^{-1})(1-x_2 y_1^{-1})(1-x_2 y_2^{-1})}.
    \end{equation*}
% \label{ex:multidegree}
The K-polynomial of $Y$ is $\mathcal{K}(I(Y))=1-x_1 x_2 y_1^{-1}
y_2^{-1}$.  Also, the multidegree of $Y$ is the linear term
of \begin{equation*} 1 - (1-x_1)(1-x_2)(1-y_1)^{-1}
(1-y_2)^{-1}, \end{equation*} namely $\mathcal{M}(I(Y))=x_1 + x_2 -
y_1 - y_2$.  \end{Example}

Returning to the more general setting, we explain how both the
multidegree $\mathcal{M}(I(Z))$ and the K-polynomial
$\mathcal{K}(I(Z))$ represent analogs of the Chow class or cohomology
class for $Z$.  The K-polynomial is named this way because it
represents a class in K-theory, which is a more exotic cohomology
theory that we will not discuss here.  We will focus on the cohomology
theory that pertains to multidegrees of varieties, namely equivariant
cohomology.

One can define the \emph{$T$-equivariant Chow ring} $A_T^*(V)$
similarly to the Chow ring from \Cref{sub:Monk}.  We will briefly
introduce the topic here.  For a more thorough development of
equivariant cohomology with an eye towards Schubert calculus, we
recommend Anderson and Fulton's book \cite{Anderson-Fulton}.

In $A_T^*(V)$, a codimension $d$,  $T$-invariant subvariety $Z \subseteq
X$ defines a class $[Z] \in A_T^d(V)$ where $[Z] = [Z']$ if $Z$ and
$Z'$ are related by a suitable $T$-equivariant deformation, and
multiplication comes from intersection of varieties as before.  See
\cite[\S 17.1]{Anderson-Fulton} for example.  It turns out that
sending $[Z]$ to the multidegree of $Z$ defines an injective
homomorphism from $A_T^*(V)$ to the polynomial ring $\bZ[x_1, \ldots,
x_n]$. Note that the torus action is what makes any of this
interesting, as the ordinary Chow ring $A^*(V)=\bZ$ trivially since it
is one cell. Alternatively, one can use the \emph{equivariant
cohomology ring} $H_T^*(V)$, which in this case is isomorphic to
$A_T^*(V)$ but has the advantage of being defined for the action of
any topological group on a space.

Recall that $\rp(w)$ is the set of reduced pipe dreams of $w$ from
\Cref{def:pipedream}. In the next theorem, $T \times T$ acts on $\Mat_n(\C)$ as in \Cref{ex:hilbert}. Also recall from \eqref{eq:double.schubs} that 
    the \emph{double Schubert polynomial} of $w \in S_n$ is
    \begin{equation}
        \fS_w(X;Y)=\fS_{w}(x_1, \ldots, x_n; y_1,
\ldots, y_n) = \sum_{D \in \rp(w)} \prod_{(i,j) \in D} (x_i-y_j).
    \end{equation}

The following theorem is due to Knutson-Miller
\cite[Thm. A]{knutson-miller-2005} in this language.  See also the
earlier work by Feh\'{e}r-Rim\'{a}nyi \cite[Thm
4.2]{Feher-Rimanyi.2003}, where double Schubert polynomials were
interpreted at Thom polynomials. 

\begin{Theorem}\label{thm:KM-multidegree}
Under the $T\times T$-action on $\Mat_{n}(\C)$, the multidegree of the
matrix Schubert variety $\MX_w$ is the double Schubert polynomial,
$\mathcal{M}(I(\MX_{w}))=\fS_{w_0 w}(x_1, \ldots, x_n; y_n,
\ldots, y_1)$.  
If we forget the left action of $T$ and just use the right action,
then the multidegree of $\MX_w$ is the single Schubert polynomial
$\fS_{w_0 w} = \fS_{w_0 w}(X; 0)$.
\end{Theorem}

To connect all this back to flag varieties and ordinary Chow rings, we
make use of a ring homomorphism $A_T^*(\Mat_n(\C)) \to A^*(\flags)$,
which sends $[Z]$ to $[(Z \cap
\GL_n(\C))/B]$ whenever $Z \subseteq \Mat_n(\C)$ is a $B$-invariant
variety.  Furthermore, the following diagram commutes:
\begin{equation*}
    \begin{tikzcd}
        A_T^*(\Mat_n(\C)) \arrow[r] \arrow[d, "\text{multidegree}" left]   &   A^*(\flags) \arrow[d, "\text{Borel presentation}"]\\
        \bZ[x_1, \ldots, x_n] \arrow[r] & \bZ[x_1, \ldots, x_n] / I_n^+.
    \end{tikzcd}
\end{equation*}
Note that $\MX_w \cap \GL_n(\C)$ is exactly the set of matrices representing a flag in the Schubert variety $X_w$. The upshot is that the multidegree $\mathcal{M}(\MX_w)$ is necessarily a polynomial representing the \emph{ordinary} Chow ring class $[X_w]$, so Theorem~\ref{thm:KM-multidegree} provides a good reason to choose the Schubert polynomial $\fS_{w_0 w}$ in particular to represent this class.

We will point out some key steps in the proof of Theorem~\ref{thm:KM-multidegree} below.
The reader is encouraged to consider the full proof in
\cite{knutson-miller-2005}.

\begin{Example}
    The matrix Schubert variety $\MX_{312}$ is $\{A=(a_{ij}) \in \Mat_3(\C) \given
    a_{21}a_{32} - a_{31}a_{22} = 0\}$, which has multidegree $x_1+x_2 -
    y_2 - y_3$ as calculated (up to reindexing) in
    Example~\ref{ex:hilbert}. This agrees with
    $\fS_{132}(x_{1},x_{2};y_{1},y_{2}) = (x_1-y_2)+(x_2-y_1)$, calculated using the two pipe
    dreams $\begin{smallmatrix} \cdot & + \\ \cdot &
    \cdot \end{smallmatrix}$ and $\begin{smallmatrix} \cdot & \cdot \\ + &
    \cdot \end{smallmatrix}$ of $132=w_{0}312$.
    \end{Example}

Finally, we explain how pipe dreams naturally appear in the process of
calculating the multidegree $\mathcal{M}(\MX_w)$.  Fix a total order
$\prec$ on monomials in which (a) $1$ is minimal and (b) $m_1 \prec
m_2$ implies $mm_1 \prec mm_2$ for any monomial $m$: such a total
order is called a \emph{term order}.  Each nonzero polynomial $f$ has
a \emph{leading monomial} meaning the smallest monomial with nonzero
coefficient appearing in $f$ with respect to the fixed term order.
Let $I \subseteq \C[z_1, \ldots, z_m]$ be an ideal.  The \emph{initial
ideal} of $I$ is the ideal $\init(I)$ generated by the leading
monomials of all $f \in I$.  The operation $I \leadsto \init(I)$
preserves multigraded Hilbert series, and hence K-polynomials and
multidegrees, reducing the problem of computing these quantities to
the case of monomial ideals. In turn, generators for the monomial
ideal $\init(I)$ can be computed by finding a \emph{Gr\"obner basis}
of $I$, which is a generating set whose leading terms generate
$\init(I)$.  See \cite{Cox-Little-OShea} for more background on
initial ideals and Gr\"obner bases.

Call a term order on $\C[z_{ij}]$ \emph{antidiagonal} if the leading
term of any minor $\Delta_{I,J}:=\det [z_{ij}]_{i \in I, j \in J}$ is
the antidiagonal term. That is, if $I = \{i_1 < \cdots < i_r\}$ and $J
= \{j_1 < \cdots < j_r\}$, the leading term of $\Delta_{I,J}$ should
be $z_{i_1 j_r} z_{i_2 j_{r-1}} \cdots z_{i_r j_1}$. The next theorem
is due to Knutson and Miller \cite{knutson-miller-2005}.

\begin{Theorem} \label{thm:minors-gb}
The determinantal minors which generate $I(\MX_w)$ according to
Theorem~\ref{thm:MX-ideal} also form a Gr\"obner basis with respect to
any antidiagonal term order.
\end{Theorem}

% \sara{Martha asks if K-M use $B_{-}$- or $B_{+}$-orbits and if that
% means we need to replace antidiagonal order by diagonal order????
% Re-think.  Brendan says it doesn't matter if they use the lower
% triangulars because conjugation preserves antidiagonal  order.}

Theorem~\ref{thm:minors-gb} implies that $\init(I(\MX_w))$ is actually generated by $\emph{squarefree}$ monomials, which simplifies the process of computing multidegrees. 

\begin{Exercise} \label{exer:multidegree}
    Suppose $J \subseteq \C[z_1, \ldots, z_m]$ is an ideal generated
by squarefree monomials. Assume there is a $T$-action on
$\operatorname{span} \{z_1, \ldots, z_m\}$ for which the variables
$z_j$ are weight vectors, so that $J$ is $T$-stable for the induced
action on $\C[z_1, \ldots, z_m]$. Let $(z_i : i \in I)$ denote the
ideal generated by the $z_{i}$'s for $i \in I \subset [m]$.  Prove the
following statements hold.
\begin{enumerate}[(a)]
\item $J$ is an intersection $\bigcap_{I \in \mathcal{I}} (z_i : i \in I)$ over some collection $\mathcal{I}$ of subsets of $[m]$.
\item Each ideal $(z_i : i \in I)$ has multidegree $\prod_{i \in I} \langle \wt(z_i), (x_1, \ldots, x_m) \rangle$, where $\langle\, ,\rangle$ is dot product. 
\item $J$ has multidegree $\sum_{I \in \mathcal{I}} \prod_{i \in I} \langle \wt(z_i), (x_1, \ldots, x_m) \rangle$ where $\mathcal{I}$ is the same as in (a).
\end{enumerate}
\end{Exercise}

For instance, in Example~\ref{ex:hilbert} the torus element $(s_1, \ldots, s_n, t_1, \ldots, t_n)$ acts on the variable $z_{ij}$ 
with weight $(0,\ldots,1,\ldots,0,0,\ldots,-1,\ldots,0)$, where the
$1$ and $-1$ are in positions $j$ and $n+i$ respectively.  Thus
$\langle \wt(z_{ij}), (x_1, \ldots, x_n, y_1, \ldots, y_n) \rangle =
x_j - y_i$, so for any subset $D \subset [n]\times [n]$ the ideal
$(z_{ij} : (i,j) \in D)$ would have multidegree $\prod_{(i,j) \in D}
(x_j - y_i)$ according to Exercise~\ref{exer:multidegree}(b).
\bigskip

\begin{Theorem} \cite{knutson-miller-2005} \label{thm:PD-intersection}
For all $w \in S_n$,
\begin{equation*}
    \init(I(\MX_w)) = \bigcap_{D \in \rp(w_0 w)} (z_{n-j+1,i} : (i,j) \in D)
\end{equation*}
where $\init(I(\MX_w))$ is taken with respect to any antidiagonal term order.
\end{Theorem}
It follows by Exercise~\ref{exer:multidegree} that $\MX_w$ has multidegree $\sum_{D \in \rp(w_0 w)} \prod_{(i,j) \in D} (x_i - y_{n-j+1})$, proving Theorem~\ref{thm:KM-multidegree}.

\begin{Example}
    Let $w = 3124$. Using Example~\ref{ex:matrix-schub}, we have
    \begin{equation*}
        \init(I(\MX_w)) = (z_{41}, z_{42}, z_{43}, z_{21}z_{32}) = (z_{41}, z_{42}, z_{43}, z_{21}) \cap (z_{41}, z_{42}, z_{43}, z_{32})
    \end{equation*}
    corresponding to the two pipe dreams of $w_0 w = 2431$,  
    \begin{equation*}
\left\{
        \begin{smallmatrix} 
            + & \cdot & + & \cdot \\ 
            + & \cdot & \cdot & \cdot \\ 
            + & \cdot & \cdot & \cdot \\ 
            \cdot & \cdot & \cdot & \cdot \\ 
        \end{smallmatrix} \quad \text{,} \quad \begin{smallmatrix} 
            + & \cdot & \cdot & \cdot \\ 
            + & + & \cdot & \cdot \\ 
            + & \cdot & \cdot & \cdot \\ 
            \cdot & \cdot & \cdot & \cdot \\ 
        \end{smallmatrix}
\right\} = \rp(w).	
    \end{equation*}
    Hence $\mathcal{M}(I(\MX_{3124})) = (x_1-y_4)(x_2-y_4)(x_3-y_4)(x_1-y_2) + (x_1-y_4)(x_2-y_4)(x_3-y_4)(x_2-y_3)$.
\end{Example}

\begin{Remark}
    One cannot recover the weights of the $T$-action on $\C[V]/I$ uniquely from the multidegree $\mathcal{M}(I)$. By Exercise~\ref{exer:multidegree}, any $T \times T$-stable subvariety $Z \subseteq \Mat_n(\C)$ will have a multidegree which is a positive integer combination of products of the binomials $x_i-y_j$, but such expansions are not unique. For instance, $(x_1-y_2) + (x_2-y_1) = (x_1-y_1) + (x_2-y_2)$ is the multidegree of the ideal $(z_{12}z_{21})$ and also of $(z_{11}z_{22})$. These are the initial ideals of $I = (z_{11}z_{22}-z_{12}z_{21})$ with respect to two different term orders, so this ambiguity actually reflects some of the power and flexibility of the multidegree, namely the fact that $\mathcal{M}(I)$ can be computed as $\mathcal{M}(\init(I))$ with respect to any choice of term order.  See also \Cref{ex:PD-BPD-2143}.
\end{Remark}

The operation of taking initial ideals can be understood in terms of
the Chow ring. The ideals $I(Z)$ and $\init(I(Z))$ are part of a
continuous (more technically, \emph{flat}) family of ideals: intuitively, one can insert parameters in
front of the non-leading terms of polynomials in $I(Z)$ and then send them
to $0$. This corresponds to deforming $Z$ into a union of linear
subspaces (possibly with multiplicities), and this is exactly the kind of
deformation that preserves Chow ring or cohomology ring classes.

\begin{Corollary}
The reduced pipe dreams index the pieces in a deformation of a matrix
Schubert variety into linear subspaces, and each piece contributes a
single monomial term to the multidegree.  
\end{Corollary}

This corollary completes our geometric interpretation of pipe dreams.
Reflecting back on this section, it is natural to wonder what happens
if a diagonal term order is used rather than an antidiagonal one.  The
situation becomes more difficult: the initial ideal $\init(I(\MX_w))$
may depend on the particular choice of diagonal term order, it need
not be generated by squarefree monomials, and Gr\"obner bases are not
known for arbitrary term orders. Nevertheless, the multidegree can
still be defined by taking appropriate multiplicities into account,
and Klein and Weigandt \cite{klein-weigandt} derive the
\emph{bumpless} pipe dream formula in
Theorem~\ref{thm:bumpless-schubert} for $\fS_w$ by applying the tools
described here to appropriate diagonal term orders.

Another interesting open problem concerns symmetric matrix Schubert varieties, defined just like $\MX_w$ with the added constraint that the matrices be symmetric (so $w$ should satisfy $w = w^{-1}$). Theorems~\ref{thm:MX-ideal} and \ref{thm:minors-gb} describing Gr\"obner bases for their ideals seem to hold verbatim, but remain unproven. Marberg and Pawlowski \cite{marberg-pawlowski-ss} proved an analogue in the skew-symmetric case. As explained in \cite{marberg-pawlowski-k-theory}, these varieties bear the same relationship to the orbit closures described in Chapter 5 as $\MX_w$ does to the Schubert variety.

\subsection{Schubert Calculus and Hermitian Eigenvalue Problems}
\label{subsec:eigenvalues}

Hermitian matrices in $\Mat_{n}(\C)$ have the property that they are
equal to their complex conjugate transpose, denoted $H=H^{*}$.  By the
Spectral Theorem, they
have all real eigenvalues and an orthonormal basis of eigenvectors.
Let $\lambda(H) = (\lambda_1(H) \geq \cdots \geq \lambda_n(H))$ denote
the vector of eigenvalues of a Hermitian matrix $H$, sorted in
decreasing order.  In this subsection, we consider the set of all
possible triples of eigenvalue vectors
\begin{equation}
\mathcal{A}_n:=\{(\lambda(H_1), \lambda(H_2), \lambda(H_1+H_2)) \given \text{$H_1, H_2$ $n \times n$ Hermitian matrices}\} \subseteq \bR^{3n}.
\end{equation}
Starting with Weyl, many inequalities have been found relating
$\lambda(H_1+H_2)$ to $\lambda(H_1)$ and $\lambda(H_2)$; see
\cite{FultonEig} for a detailed survey of this and related topics. In
1962, Horn \cite{Horn} conjectured that $\mathcal{A}_n$ is a convex
polyhedron defined by a specific finite set of linear inequalities
defined recursively in terms of $n$. Horn's conjecture was eventually
resolved in full by work of Knutson and Tao \cite{KnutsonTao}; a
crucial step was a characterization of $\mathcal{A}_n$ in terms of
Schubert calculus on the Grassmannian by Klyachko \cite{KlyachkoEig},
which we describe in this section.

Some eigenvalue problems can be phrased as optimization problems. For
example, it is easy to prove using the Spectral Theorem that
\begin{equation}
 \lambda_n(H) = \min_{|v| = 1} v^* H v \quad \text{and} \quad \lambda_1(H) = \max_{|v| = 1} v^* H v
\end{equation}
when $H$ is an $n\times n$ Hermitian matrix. The next result, due to
Hersch and Zwahlen \cite{HerschZwahlen}, realizes \emph{any partial sum} of
eigenvalues 
\begin{equation}\label{eq:lambda.I}
\lambda_{I}(H)=\lambda_{i_1}(H) + \cdots + \lambda_{i_k}(H)
\end{equation}
 for $I = \{i_1 < \cdots < i_k\} \subseteq [n]$ as the solution to a
similar optimization problem. The formulation uses Schubert varieties
in the Grassmannian variety $\Gr(k,n)$ of $k$-dimensional subspaces in
$\C^n$ as defined in \Cref{sub:Grassmannians.intro} and further
developed in \Cref{sub:Grassmannians}.    Given a $k$-dimensional
subspace $L \in \Gr(k,n)$, let $U_L$ denote any $n \times k$ matrix
whose columns form an orthonormal basis for $L$.  Then,
$U_{L}^{*}HU_{L}$ is an $k \times k$ matrix with real values down the
diagonal, since $v^{*}Hv$ is real for all complex vectors $v$.  Let
$\tr(U_L^* H U_L)$ denote the trace.  Recall that 
\begin{equation}\label{eq:Icheck}
I^\vee = \{n-i+1 \given i \in I\}.
\end{equation}
Let $\Lambda^>_\bullet(H)$ be the flag $\langle v_1 \rangle \subseteq
\langle v_1, v_2 \rangle \subseteq \cdots \subseteq \C^n$ where $v_1,
\ldots, v_n$ is an ordered orthonormal basis of eigenvectors of $H$
whose eigenvalues come in weakly decreasing order.  Similarly, let
$\Lambda^<_\bullet(H)$ denote the flag defined by the basis $v_n,
\ldots, v_1$ in reverse order.  Recall, the Grassmannian Schubert
variety indexed by a $k$-subset $I$ is
\begin{equation}
X_I(F_{\bullet}) = \{L \in \Gr(k,n) \given \dim(L \cap F_{i_h}) \geq h
\text{ for $h = 1, \ldots, k$}\}
\end{equation}
for a fixed flag $F_{\bullet} \in \Fl(n)$.  Here we will consider
Schubert varieties with respect to the flags $\Lambda^>_\bullet(H)$
and $\Lambda^<_\bullet(H)$.

\begin{Lemma}\cite{HerschZwahlen} \label{lem:HZ} For an $n \times n$ Hermitian matrix $H$ and $k$-subset $I \subseteq [n]$,
    \begin{equation*}
    \lambda_I(H) = \min_{L \in X_I(\Lambda^>_\bullet(H))} \tr(U_L^* H U_L) \quad \text{and} \quad \lambda_{I^\vee}(H) = \max_{L \in X_I(\Lambda^<_\bullet(H))} \tr(U_L^* H U_L)
    \end{equation*}
\end{Lemma}
\bigskip

For a given pair of Hermitian matrices $(H_{1},H_{2})$, 
we have the identity 
\begin{equation} \label{eq:zero-trace}
0 = \tr(U_L^* H_1 U_L) + \tr(U_L^* H_2 U_L) - \tr(U_L^* (H_1+H_2) U_L).
\end{equation}
For all $k$-subsets $I,J,K \subset [n]$ such that
\begin{equation} \label{eq:eig-intersection}
    X_{I^\vee}(\Lambda_\bullet^<(H_1)) \cap X_{J^\vee}(\Lambda_\bullet^<(H_2)) \cap X_{K}(\Lambda_\bullet^>(H_1+H_2))
    \end{equation}
is non-empty,  we can choose $L$ from the intersection and deduce the inequality
\begin{equation} \label{eq:eig-ineq}
    \lambda_I(H_1) + \lambda_J(H_2) \geq \lambda_K(H_1 + H_2)
\end{equation}
by applying Lemma~\ref{lem:HZ} to each term of \eqref{eq:zero-trace}
simultaneously.  

Recall from \S\ref{sub:Grassmannians} that the Littlewood-Richardson coefficient $c_{IJ}^K$ is the number of points in the intersection $X_{I^\vee}(E_\bullet) \cap X_{J^\vee}(F_\bullet) \cap X_{K}(G_\bullet)$ with $E_\bullet, F_\bullet, G_\bullet$ being generic flags. Here, we are using the bijection between $k$-subsets of $[n]$ and partitions whose Young diagrams are contained in a $k \times (n-k)$ rectangle in order to write $c_{IJ}^K$ rather than the more common notation $c_{\lambda\mu}^{\nu}$ for partitions $\lambda, \mu, \nu$.  Although the flags in \eqref{eq:eig-intersection} may not be in general position, one can show that this can only increase the size of the intersection compared to the generic case. In particular, if $c_{IJ}^K \neq 0$, then the intersection \eqref{eq:eig-intersection} is always nonempty, and therefore the inequality \eqref{eq:eig-ineq} holds for that triple $(I,J,K)$ and all Hermitian matrices $H_1$ and $H_2$. That all these inequalities must hold for points $(\lambda(H_1), \lambda(H_2), \lambda(H_1+H_2)) \in \mathcal{A}_n$ was observed by Johnson \cite{JohnsonEig}. Klyachko proved the converse stated below, yielding a concrete description of $\mathcal{A}_n$ as a convex polyhedron with specified inequalities.

\begin{Theorem}\cite{KlyachkoEig}
Let $\alpha =(\alpha_{1},\dots ,\alpha_{n}), \beta =(\beta_{1},\dots
,\beta_{n}), \gamma = (\gamma_{1},\dots , \gamma_{n}) \in
\mathbb{R}^{n}$ be three decreasing sequences.  Then, there are $n
\times n$ Hermitian matrices $A,B,C$ with $A+B = C$ and $\lambda
(A)=\alpha$, $\lambda (B)=\beta$, $\lambda (C)=\gamma$ if and only if
$\sum_{i=1}^{n} \alpha_{i}+ \sum_{i=1}^{n} \beta_{i}=\sum_{i=1}^{n}
\gamma_{i}$ and $\alpha ,\beta , \gamma $ satisfy
\begin{equation} \label{eq:coordinate-ineq}
\sum_{i \in I} \alpha_{i}+ \sum_{j \in J}\beta_{j} \geq \sum_{k \in K} \gamma_{k}
\end{equation}
 for all subsets $I, J, K \subset [n]$ of the same size such that
$c_{IJ}^K \neq 0$.
\end{Theorem}

Independently, Belkale \cite{Belkale} and Knutson and Tao
\cite{KnutsonTao} later showed that it suffices to take exactly those
inequalities determined by triples $(I,J,K)$ for which $c_{IJ}^K = 1$.
That is, to define Horn's polyhedron, it suffices to consider the
inequalities in \eqref{eq:coordinate-ineq} for all $(I,J,K)$ with
$c_{IJ}^K = 1$, and this list is minimal.  
\begin{Example}
Consider the case $n = 2$. There is no real loss of generality in
assuming $A,B,C$ have trace zero, so say their eigenvalues are $\pm
\alpha, \pm \beta, \pm \gamma$ with $\alpha, \beta, \gamma \geq 0$.
The nonzero Littlewood-Richardson coefficients with $n=2$ are witnessed by the
three puzzles
\begin{equation*}
\begin{tikzpicture}
    \draw (0,0) -- node {$1$} (1,0) -- node {$0$} (2,0);
    \draw (0,0) -- node {$0$} (0.5,0.866) -- node {$1$} (1,2*0.866);
    \draw (1,2*0.866) -- node {$1$} (1.5,0.866) -- node {$0$} (2,0);
    \draw (0.5,0.866) -- node {$10$} (1,0) -- node {$0$} (1.5,0.866) -- node {$1$} (0.5,0.866);
\end{tikzpicture} \qquad \qquad 
\begin{tikzpicture}
    \draw (0,0) -- node {$1$} (1,0) -- node {$0$} (2,0);
    \draw (0,0) -- node {$1$} (0.5,0.866) -- node {$0$} (1,2*0.866);
    \draw (1,2*0.866) -- node {$0$} (1.5,0.866) -- node {$1$} (2,0);
    \draw (0.5,0.866) -- node {$1$} (1,0) -- node {$10$} (1.5,0.866) -- node {$0$} (0.5,0.866);
\end{tikzpicture} \qquad \qquad 
\begin{tikzpicture}
    \draw (0,0) -- node {$0$} (1,0) -- node {$1$} (2,0);
    \draw (0,0) -- node {$0$} (0.5,0.866) -- node {$1$} (1,2*0.866);
    \draw (1,2*0.866) -- node {$0$} (1.5,0.866) -- node {$1$} (2,0);
    \draw (0.5,0.866) -- node {$0$} (1,0) -- node {$1$} (1.5,0.866) -- node {$10$} (0.5,0.866);
\end{tikzpicture}
\end{equation*}
so the Littlewood-Richardson coefficients 
\begin{equation*}
    1 = c_{\{1\}\{2\}}^{\{2\}} = c_{\{2\}\{1\}}^{\{2\}} = c_{\{1\}\{1\}}^{\{1\}}
\end{equation*}
determine the facets of $\mathcal{A}_{2}$.  Alternatively, one can
compute the relevant Littlewood-Richardson coefficients via
intersection of Schubert varieties in $\Gr(1,2) = \C \bP^1$. Here
$X_{\{1\}}$ is a point and $X_{\{2\}} = \C\bP^1$, so the
corresponding nonempty zero dimensional intersections are easy to
identify.  These coefficients give the three inequalities
\begin{equation*}
    \lambda_{1}(A) + \lambda_{2}(B)\geq \lambda_2(C), \quad \lambda_{2}(A) + \lambda_{1}(B) \geq \lambda_2(C), \quad \text{and} \quad \lambda_1(A) + \lambda_1(B) \geq \lambda_1(C),
\end{equation*}
or equivalently
\[
\alpha-\beta \geq -\gamma, \quad -\alpha +\beta \geq - \gamma, \quad \text{and} \quad \alpha +\beta \geq \gamma.
\]
To visualize this polyhedron, divide through by $\alpha$ and plot
$\gamma/\alpha$ against $\beta/\alpha$ to get a rectangular polyhedron
as show in in \Cref{fig:eigenvalue.polyhedron}.  Then $\mathcal{A}_2$
is a cone over this region times a factor of $\bR \times \bR$
accounting for $\tr(A)$ and $\tr(B)$.

\begin{figure}
\begin{center}
\begin{tikzpicture}
    \draw[thick] (0,0) -- (3,0); \draw[thick] (0,0) -- (0,3);
\filldraw[fill=black!20!white, draw=blue] (3,2) -- (1,0) -- (0,1) --
(2,3); \draw (3,0) node[right] {$\scriptstyle \beta/\alpha$}; \draw
(0,3) node[above] {$\scriptstyle \gamma/\alpha$}; \node at (2.7,2.7)
[circle,fill,inner sep=.5pt]{}; \node at (2.8,2.8) [circle,fill,inner
sep=.5pt]{}; \node at (2.9,2.9) [circle,fill,inner sep=.5pt]{};
\end{tikzpicture}
\end{center}
\caption{A projection of the Horn polyhedron for $n=2$.}
\label{fig:eigenvalue.polyhedron}
\end{figure}
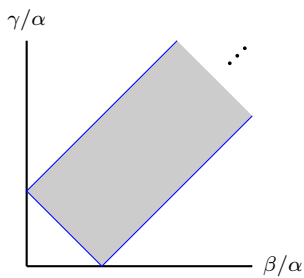
\end{Example}

In the context of this chapter, one might ask

\begin{quote}
Is there a natural generalization of this eigenvalue problem for which
Schubert structure constants in the flag varieties are relevant?
\end{quote}
See the work of Belkale \cite{Belkale.2006}, Belkale-Kumar
\cite{Belkale.Kumar.2006}, and Richmond \cite{Richmond.2009} for
interesting generalizations of these eigenvalue problems in a broader
context.

\subsection{Degeneracy Loci and Double Schubert Polynomials }\label{sub:DegeneracyLoci.DoubleSchubs}

In this section, we consider a profitable extension of the material of
\S\ref{sec:flags} to a more general setting where the matrices have
entries parameterized by a topological space. A \emph{vector bundle}
$E$ over a space $X$ is, roughly, an assignment of vector spaces $E_x$
to points $x \in X$ which fit together appropriately to make $E =
\bigcup_x E_x$ itself a space. Locally, we should be able to identify
$E$ and $X \times V$ for a fixed vector space $V$, in the same way
that a manifold looks locally like Euclidean space. The projection
$\pi : E \to X$ sending every $v \in E_x$ to $x$ should correspond
under this identification to the projection $X \times V \to
X$. Accordingly, the vector spaces $E_x$ are called the \emph{fibers}
of the bundle, since they are the sets $\pi^{-1}(x)$.  See
\cite{milnor-stasheff} for more details than we provide below.  

If $X$ is connected, then $\dim E_x$ is constant and is called the
\emph{rank} of $E$.  A vector bundle of rank $1$ is called a \emph{line
bundle}. We have not specified the field of scalars;
choosing $\bR$ or $\C$ defines real or complex vector bundles
respectively.

\begin{Example} Here are some commonly seen vector bundles.  
\begin{itemize}
    \item The rank $n$ \emph{trivial bundles} $\bR^n \times X$ and $\bC^n \times X$, including the $n = 0$ case $X$.
    \item The \emph{M\"obius bundle} $E$ over the circle $S^1 \subseteq \C$, where $E_{e^{i\theta}} = \bR e^{i\theta/2}$.
    \item If $M$ is a smooth manifold, the \emph{tangent bundle} $TM$, where $(TM)_x = T_x M$ is the tangent space at $x$.
    \item The rank $k$ \emph{tautological bundle} $E$ over $\Gr(k,n)$, where $E_V = V$. 
\end{itemize}
\end{Example}

A key special case of the last example is the tautological bundle over
$\C \bP^n$, often denoted $\cO(-1)$, but beware that this common
notation omits the base space $\C \bP^n$! Each $\cO(-1)$ is a complex \emph{line bundle}.  Also write $\cO(-1)_{\bR}$ for the tautological bundle over $\bR \bP^n$.

\begin{Example} \label{ex:mobius}
    Say $E$ is a real line bundle over $S^1$, which we view as the
    unit circle in the $xy$-plane. View each fiber $E_z \simeq \bR$ as
    a line perpendicular to $S^1$ at $z$, and then identify that line
    with the unit interval inside it. This gives a picture of $E$
    inside $\bR^3$. \Cref{fig:mobius} shows the trivial line bundle and the M\"obius
    bundle drawn this way.

\begin{figure}[h]
    \begin{center}
    \includegraphics[scale=0.3]{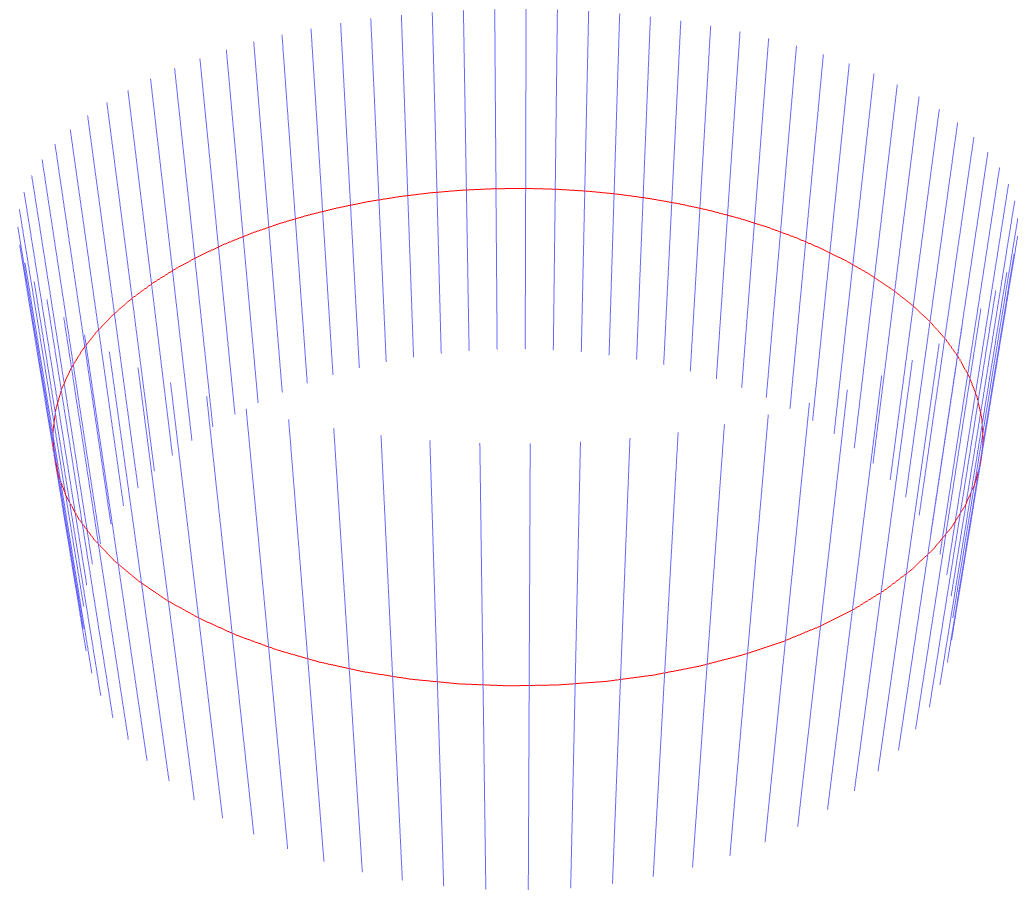} \qquad 
    \raisebox{4mm}{\includegraphics[scale=0.4]{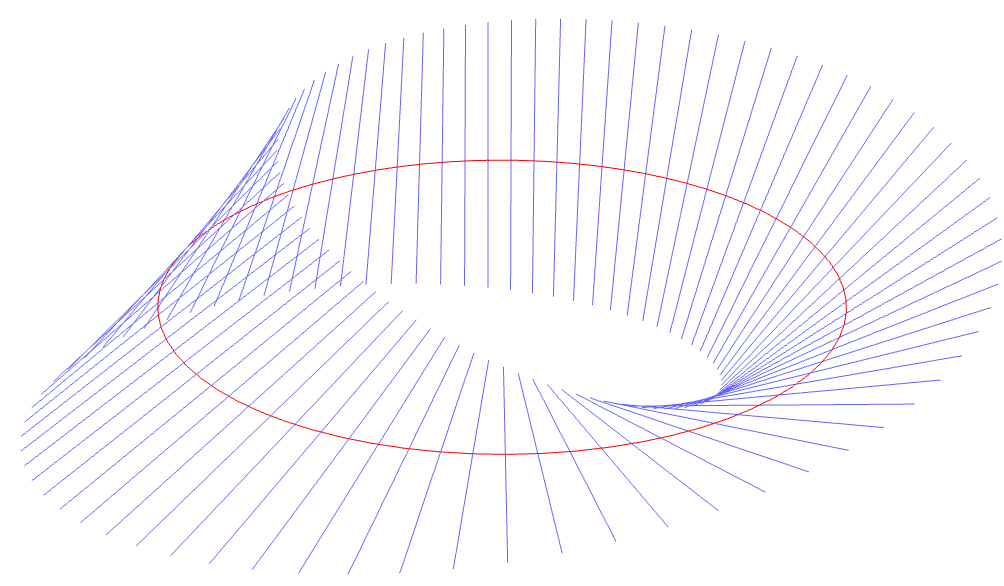}}
    %(s*cos(t/2)*cos(t), s*cos(t/2)*sin(t), s*sin(t/2))
    \end{center}
\caption{A trivial bundle and the M\"obius
    bundle over $S^1$.}
\label{fig:mobius}
\end{figure}

    To be precise, we have drawn some of the fibers $E_z$, and $E$ is the surface formed by all $E_z$. The projection $\pi$ sends each point in the fiber $E_z$ (in blue) to the point $z$ where that fiber intersects the base space $S^1$ (in red).

    Observe that $S^{1}$ is isomorphic to $\bR \bP^1$ by identifying a
point on the unit circle in the upper half-plane of $\bR^2$ with a
line through the origin in $\bR^2$ intersecting that point, and also
identifying the antipodes $(-1,0)$ and $(1,0)$.  Realizing $\bR^2$ as
the complex plane, this is the map $e^{i\theta} \mapsto \bR
e^{i\theta/2}$. Under that identification, the M\"obius bundle $E$ has
fibers $ E_{e^{i\theta}} = E_{\bR e^{i\theta/2}} = \bR e^{i\theta/2}$,
so the M\"obius bundle $E$ is isomorphic to the tautological bundle
$\cO(-1)_{\bR}$.
\end{Example}

Many basis-independent constructions from linear algebra can be
upgraded to vector bundles by using them on each fiber $E_x$
separately: linear maps, direct sums, quotients, duals, tensor
products, etc.  For instance, a morphism of vector bundles $f : E \to
F$ over $X$ is a suitably continuous collection of linear maps $f_x :
E_x \to F_x$ for $x \in X$. As another example, the double dual
$E^{**}$ is isomorphic to $E$, but the dual bundle $E^*$ is typically
\emph{not} isomorphic as a vector bundle to $E$ when $E$ is complex. A
construction that does not carry over well to vector bundles is the kernel of a morphism.  For example, the map $f :
(x,v) \mapsto (x,xv)$ is a morphism on the trivial bundle $\bR \times
\bR$ over $\bR$, but $\dim \ker f_x$ is not constant on $\bR$, so
$\ker f$ is not naturally a vector bundle over the first factor.

There is some special notation surrounding the tautological line
bundles $\cO(-1)$.  One writes $\cO(1) = \cO(-1)^*$ for the dual
tautological bundle, and $\cO(\pm r) = \cO(\pm 1)^{\otimes r}$ for $r
> 0$.

\begin{figure}
\begin{center}
\includegraphics[scale=0.4]{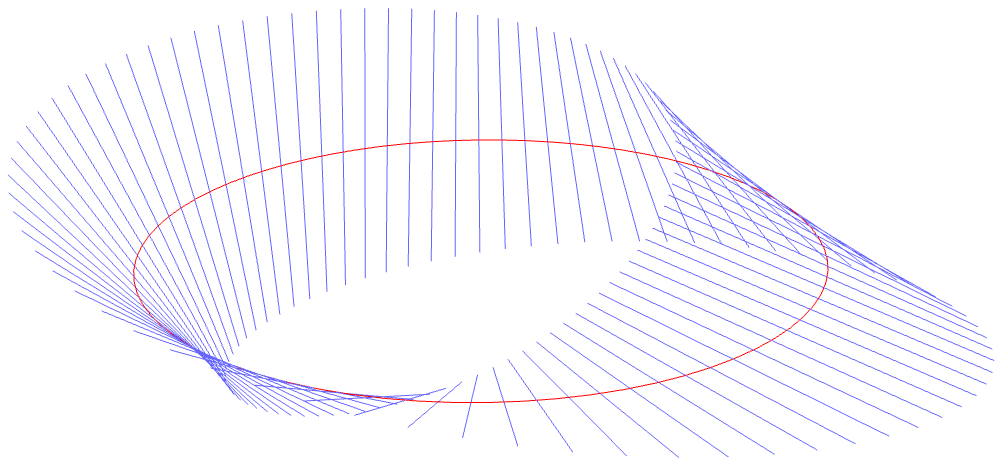} \qquad
\raisebox{-7mm}{\includegraphics[scale=0.4]{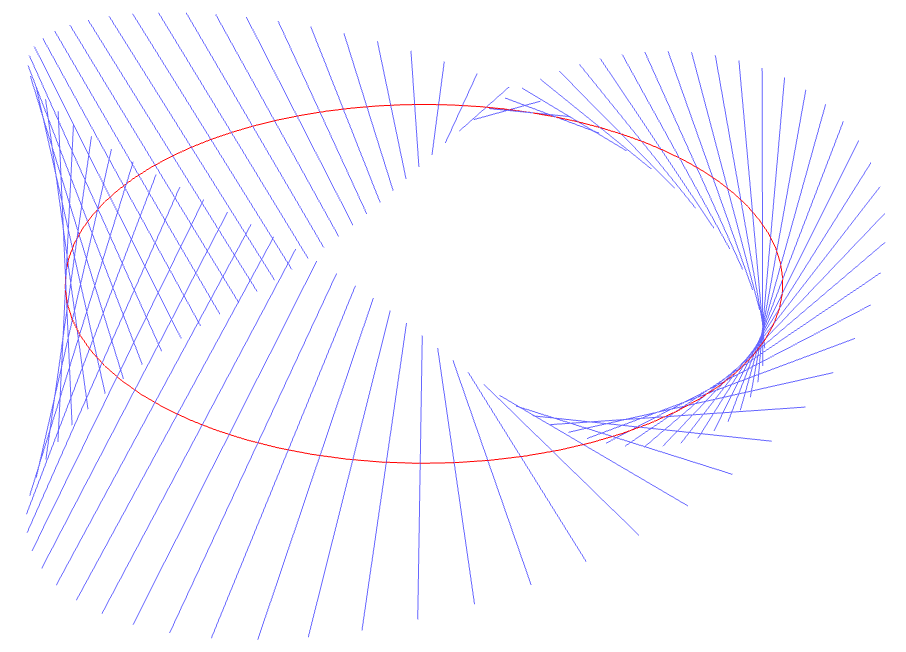}}
%(s*cos(t/2)*cos(t),s*cos(t/2)*sin(t), s*sin(t/2))
\end{center}
\caption{The bundles $\cO(1)_{\bR}$ and $\cO(-2)_{\bR}$ over $S^1$. }
\label{fig:mobius.2}
\end{figure}

\begin{Example} Let $E = \cO(-1)_{\bR}$ as in
Example~\ref{ex:mobius}. Then the dual $E^* = \cO(1)_{\bR}$ and $E
\otimes E = \cO(-2)_{\bR}$.  Their fibers over $e^{i\theta}$ are $\bR
e^{-i\theta/2}$ and $\bR e^{i\theta}$ respectively as shown in
\Cref{fig:mobius.2}.  There is an obvious isomorphism between
$\cO(1)_{\bR}$ and $\cO(-1)_{\bR}$: reflect across the $xy$-plane. In
turn, this means $\cO(-2)_{\bR} \simeq \cO(1)_{\bR} \otimes
\cO(-1)_{\bR}$, which is isomorphic to the trivial line bundle by the
map $\lambda \otimes v \mapsto \lambda(v)$. More concretely, the
morphism $\phi_z : (E \otimes E)_z = \bR z$ to $\bR$ defined by $v
\mapsto vz^{-1}$ is an isomorphism. We leave it to the reader to
picture the resulting homeomorphism between the cylinder $S^1 \times
[-1,1]$ and the doubly twisted M\"obius strip on the right above.
\end{Example}

A special case of a morphism of vector bundles deserving its own term
is a section of $E$.  A \emph{section} is specifically a morphism of
vector bundles $\phi:L \to E$ over the same space $X$ where $L$ is a
trivial line bundle. Equivalently, a section is determined by a
continuous map $s : X \to E$ such that $s(x) \in E_x$ for all
$x$.  

\begin{Example}\label{ex:section.line.bundle.Or}
Let us work out what it means to have a section $\phi$ of the line
bundle $\cO(r)$ over $\C \bP^n$ for a positive integer $r$. By definition, this means $\phi(\C v)$
is a linear functional on the line $(\C v)^{\otimes r} \subseteq
(\C^{n+1})^{\otimes r}$ for each $v \in \C^{n+1} \setminus
\{0\}$. Since $\{v^{\otimes r}\}$ is a basis of that line, $\phi$ is
determined by the formula $\phi(\C v)(cv^{\otimes r}) = cf(v)$ for all
$c \in \C$, where $f : \C^{n+1} \setminus \{0\} \to \C$ is some
function. However, not just any $f$ will do, because we must have
$\phi(\C v) = \phi(\C av)$ for all $a \in \C^{\times}$. Thus,
\begin{equation*}
    f(v) = \phi(\C v)(v^{\otimes r}) = \phi(\C av)(v^{\otimes r}) = \phi(\C av)(a^{-r}(av)^{\otimes r}) = a^{-r}f(av),
\end{equation*}
so $f(av) = a^r f(v)$ for all $a$. Conversely, if $f$ is homogeneous of degree $r$, then the formula for $\phi$ above is well-defined and we get a section of $\cO(r)$.
\end{Example}

\begin{Definition} \label{def:simple-degen-locus} Let $r$ be a
nonnegative integer.  The \emph{$r$\textsuperscript{th} degeneracy
locus} of a morphism $f : E \to F$ of vector bundles over $X$ is
$\Omega_r(f) = \{x \in X \given \rank f_x \leq r\}$. If $r = 0$, one
usually says $\Omega_{0}(f)$ is the \emph{zero locus}.
\end{Definition}

\begin{Example} \label{ex:degen} Here are some familiar examples of
degeneracy loci.

\begin{enumerate}[(a)]
\item The zero locus of the morphism $\bR \times \bR \to \bR \times \bR$, $(x,v) \mapsto (x,xv)$, both viewed as bundles over the first factor, is $\{0\}$.
\item A vector field on a smooth manifold is a section of the tangent bundle $TM$. The zero locus is its set of zeros.
\item A map $f : M \to N$ between smooth manifolds induces a morphism
$df : TM \to TN$, the differential of $f$. The set of critical points
of $f$ is the set in $M$ where $df$ has less than full rank, i.e. the
$r$\textsuperscript{th} degeneracy locus of $df$ where $r = \min(\dim
M, \dim N)-1$. A minor technical detail is that $df$ is not a morphism
according to our definition, since $df$ maps the fiber $(TM)_x$ not to
itself but to $(TN)_{f(x)}$. This can be fixed by replacing $TN$ with
the \emph{pullback bundle} $df^*(TN)$ over $M$, which has fibers
$(df^*(TN))_x = (TN)_{f(x)}$ for $x \in M$.
\item Recall from \Cref{ex:section.line.bundle.Or} that if $f$ is any
homogeneous function of degree $r$, then there is corresponding
section of $\cO(r)$.  If we restrict $f$ to be a polynomial, this
shows that degree $r$ hypersurfaces in $\C \bP^n$ are the same thing
as the zero loci of sections of $\cO(r)$.
\end{enumerate} 
\end{Example}

With appropriate assumptions on $X$ and $f: E \to  F$, each degeneracy locus
determines a cohomology class $[\Omega_r(f)] \in H^*(X)$. Thom
observed the remarkable fact that, for generic $f$, this class
$[\Omega_r(f)]$ equals the evaluation of a \emph{universal
polynomial}, depending only on $r$, $\rank E$, and $\rank F$, at
certain cohomological invariants of $E$ and $F$. An early example of
this phenomenon is the Poincar\'e-Hopf Theorem
\cite[p. 35]{milnor-diff} in the case of real manifolds.

\begin{Theorem}[\textbf{Poincar\'e-Hopf Theorem}] Let $V$ be a vector
field on a compact smooth manifold $M$ with isolated zeros. If $V_x =
0$, let $\operatorname{index}_V(x)$ be the degree of the map $y
\mapsto V_y/|V_y|$ from a sufficiently small sphere around $x$ to
$S^{\dim M - 1}$, having chosen local coordinates. Then, $\sum_x
\operatorname{index}_V(x)= \chi(M)$ where the sum is over all $x$ such
that $V_{x}=0$ and $\chi(M)$ is the Euler characteristic of $M$.  
\end{Theorem}
The \emph{hairy ball theorem} \cite[p. 30]{milnor-diff} is a
well-known corollary.  It says that any vector field on a real
even-dimensional sphere must have a zero.  This follows from the fact
that $\chi(S^{2m}) = 2$.

Our goal in this subsection is to state a formula for the cohomology
classes $[\Omega_r(f)]$ in terms of double Schubert polynomials. These
classes will be polynomials evaluated at invariants of a vector bundle
called \emph{Chern classes}. Chern classes are defined for any complex
vector bundle $E$ over a smooth manifold $X$, so from now on we assume
all vector bundles are complex bundles over a smooth base space. In general, if $E$ is a rank $r$
bundle over $X$, there are associated Chern classes $c_i(E) \in
H^{2i}(X;\bZ)$ for each $0 \leq i \leq r$.

We start with the simplest case: the Chern class $c_1(L)$ of a line
bundle $L$ over a compact smooth manifold.  
Recall the following facts from \S\ref{sub:Monk}.
\begin{itemize}
    \item $H^*(\C \bP^n; \bZ) = \bZ[\alpha]/(\alpha^{n+1})$ where $\alpha \in H^2(\C \bP^n)$ is the cohomology class of a linear hyperplane (cf. Exercise~\ref{ex:projective.space}).
    \item A map $f : X \to Y$ induces a ring homomorphism $f^* : H^*(Y) \to H^*(X)$, and if $g$ is homotopic to $f$ then $g^* = f^*$.
\end{itemize}
Like cohomology classes, vector bundles can be pulled back along
maps. If $E$ is a bundle over $Y$ and $f:X \to Y$ is a smooth map, the associated \emph{pullback
bundle} $f^*E$ over $X$ has fibers $(f^*E)_x = E_{f(x)}$.

\begin{Example} \label{ex:veronese} Take $\ell = [x_1 : x_2 : x_3] \in
\C \bP^2$. If $e_1, e_2, e_3$ is a basis of $\C^3$ then $\{e_i \otimes
e_j : i,j \in [3]\}$ is a basis of $\C^3 \otimes \C^3 \simeq \C^9$,
and in this basis the line $\ell \otimes \ell$ has coordinates $[x_i
x_j : 1 \leq i, j \leq 3]$. Alternatively, this is the image of
$(\ell,\ell)$ under the Segre embedding (cf. Exercise~\ref{ex:segre}).
We compute
\begin{equation*}
    (f^* \cO(-1))_{\ell} = \cO(-1)_{f(\ell)} = \cO(-1)_{\ell \otimes \ell} = \ell \otimes \ell,     
\end{equation*}
hence $f^* \cO(-1) = \cO(-2)$.  
\end{Example}

\begin{Definition} \label{def:chern}
    A \emph{classifying map} for a complex line bundle $L$ over $X$ is a smooth map $g : X \to \C \bP^n$, for some $n$, such that $L$ is the pullback $g^* \cO(1)$. The \emph{first Chern class} of $L$ is $c_1(L) = g^*(\alpha) \in H^2(X; \bZ)$ where $\alpha \in H^2(\C \bP^n; \bZ)$ is the class of a hyperplane.
\end{Definition}

When $X$ is a compact smooth manifold, a classifying map always exists
and is uniquely determined up to homotopy by $L$, making $c_1(L)$
well-defined. Intuitively, compactness allows all (dual) fibers
$L_x^*$ to be identified as lines in one sufficiently large $\C^{n+1}$, and then $g(x) = L_x^* \in \C \bP^{n}$. It is more common to define a
classifying map using $\cO(-1)$ instead of $\cO(1)$, but this way
seems more convenient for defining $c_1$. For a more detailed development of vector bundles over smooth manifolds, see \cite{milnor-stasheff}, especially \S 5 and \S 14 on classifying maps and Chern classes.

\begin{Example} \label{ex:chern}
    \hfill \begin{enumerate}[(a)] \item If $L = X \times \C$ is
trivial, the constant map $X \to \C \bP^0$ is a classifying map for
$L$, so $c_1(L) = 0$. The reader should verify this statement.  \item
The identity map $X=\C\bP^{1} \to \C\bP^{1}$ is a classifying map for
$\cO(1)$, so $c_1(\cO(1)) = \alpha$.  \item Taking duals in
Example~\ref{ex:veronese} shows that $f$ is a classifying map for
$\cO(2)$ over $\C\bP^2$.  \item The complex conjugation map $\theta :
\C\bP^n \to \C\bP^n, \ell \mapsto \overline{\ell}$ is a classifying
map for the tautological bundle $\cO(-1)$. Our Chow ring approach to
cohomology is not well-suited to computing the induced map $\theta^*$
on $H^*(\C\bP^n)$ since $\theta$ is not a morphism of complex
varieties, but one can show using more topological tools that
$\theta^*(\alpha) = -\alpha$. More generally, if $f$ is a classifying
map for a line bundle $L$, then $\theta \circ f$ is a classifying map
for $L^*$, which implies the general formula $c_1(L^*) =
-c_1(L)$.  \item Consider the external tensor product bundle $\cO(-1)
\boxtimes \cO(-1)$ over $\C\bP^m \times \C\bP^n$, whose fiber at
$(\ell_1, \ell_2)$ is $\ell_1 \otimes \ell_2$. The Segre embedding $f
: \C\bP^m \times \C\bP^n \to \C\bP^{(m+1)(n+1)-1}$ is a classifying
map for $\cO(-1) \boxtimes \cO(-1)$. One can show $f^*(\alpha) =
\alpha_1 + \alpha_2$ where $\alpha_i = \pi_i^*(\alpha)$ with $\pi_1,
\pi_2$ the projections from $\C\bP^m \times \C\bP^n$ onto its two
factors. In turn, this implies the general formula $c_1(L_1 \otimes
L_2) = c_1(L_1) + c_1(L_2)$.
    
    \item Let $T_i$ be the rank $i$ \emph{tautological bundle} over
$\Fl(n)$, with $(T_i)_{V_{\bullet}} = V_i$.  Then the quotient bundle
$T_i/T_{i-1}$ is a line bundle over $\Fl(n)$ and $x_i =
c_1((T_i/T_{i-1})^*)$ are the ``standard'' generators of $H^*(\Fl(n);
\bZ)$ with respect to which the Schubert variety $X_w$ has cohomology
class $\fS_{w_0 w}(x_1, \ldots, x_n)$.  \end{enumerate} \end{Example}

Definition~\ref{def:chern} can be understood by thinking of the line
bundle $L$ as providing a particular notion of ``functions on $X$''
and $c_1(L)$ as a ``hyperplane class'' with respect to such a
function. Indeed, a hyperplane in $\C \bP^n$ is the zero set of a
homogeneous linear polynomial, which is not quite a valid function on
$\C \bP^n$, but rather a section of $\cO(1)$ as per
Example~\ref{ex:degen}(d). More generally, we have the following
theorem from \cite[Appendix E]{Fulton-book}.  

\begin{Theorem} \label{thm:chern-zeroes}
    If $f$ is a generic section of a line bundle $L$ over a smooth
    projective variety $X$, then $c_1(L) \in H^2(X;\bZ)$ is the
    cohomology class of the zero locus of $f$. 
\end{Theorem}

\begin{Example}
    As discussed in Example~\ref{ex:degen}(d), the zero locus of a generic section of $\cO(r)$ is a degree $r$ hypersurface $Z \subseteq \C \bP^n$, so $c_1(\cO(r)) = [Z] \in H^2(\C \bP^n)$. If $H_1, \ldots, H_{n}$ are generic linear hyperplanes in $\C \bP^n$, then they intersect in a point by basic linear algebra, so $\alpha^{n} = [H_1] \cdots [H_n] = [H_1 \cap \cdots \cap H_n]$ is the class of a point. On the other hand, $Z \cap H_1 \cap \cdots \cap H_{n-1}$ is generically the union of $r$ points by B\'ezout's theorem, so $c_1(\cO(r))\alpha^{n-1} = r\alpha^{n}$. Since $c_1(\cO(r)) \in H^2(\C \bP^n) = \bZ \alpha$, this forces $c_1(\cO(r)) = r \alpha$.
\end{Example}
One might be tempted to use Theorem~\ref{thm:chern-zeroes} as a
definition of $c_1(L)$. This does not work in our algebraic setting
because there are line bundles that have no nonzero algebraic
sections. For instance, for $1 < i < n$ the Chern classes $x_i =
c_1((T_i/T_{i-1})^*)$ in Example~\ref{ex:chern}(f) are not positive
linear combinations of Schubert polynomials in $H^*(\flags)$, which by
Inherited Positivity \Cref{cor:Schub.basis} means that they are not
the cohomology classes of any subvarieties.

Theorem~\ref{thm:chern-zeroes} is a simple example of a nice
degeneracy locus formula in cohomology. We can now state a much more
general formula. Suppose $E_{\bullet}, F_{\bullet}$ are complete flags of vector bundles
over $X$ with a morphism $f:E_\bullet \to F_{\bullet}$, meaning that there are bundles
$E_1, \ldots, E_n, F_1, \ldots, F_n$ with $\rank E_j = \rank F_j = j$
and morphisms
\begin{equation*}
    0 = E_0 \hookrightarrow E_1 \hookrightarrow \cdots \hookrightarrow E_n  \xrightarrow{f} F_n \twoheadrightarrow \cdots \twoheadrightarrow F_1 \twoheadrightarrow F_0 = 0.
\end{equation*}
Given $w \in S_n$, define the \emph{flagged degeneracy locus}
\begin{equation} \label{eq:flagged-degen}
    \Omega_w(f) = \{x \in X \given \rank(f_x : (E_i)_x \to (F_j)_x) \leq \rk(w)[i,j]\},
\end{equation}
where by $f : E_i \to F_j$ we really mean the composition $E_i
\hookrightarrow E_{n} \xrightarrow{f} F_{n} \twoheadrightarrow F_j$. The
reader should be careful to note that the rank conditions defining
$\Omega_w(f)$ are those defining the \emph{opposite} Schubert variety
$X_{w_0 w}(\oppositeE_{\bullet})$, which is analyzed further in  \Cref{ex:schubert-locus}.

\begin{Theorem}[\textbf{Fulton's flagged degeneracy locus formula} \cite{Fulton1}]\label{thm:fulton-degeneracy-loci}
    Say $X$ is a smooth complex projective variety and $f : E_{\bullet} \to F_{\bullet}$
    is a morphism of flagged vector bundles over $X$. If all
    irreducible components of $\Omega_w(f)$ have codimension
    $\ell(w)$, then
\[
    [\Omega_w(f)] = \fS_w(x_1, \ldots, x_n; y_1, \ldots, y_n) \in H_{T}^{*}(X)
\]
    where $x_i = c_1(\ker(F_i \twoheadrightarrow F_{i-1}))$ and $y_i = c_1(E_i/E_{i-1})$.
\end{Theorem}

Note that no part of the expression $\fS_w(x_1, \ldots, x_n; y_1,
\ldots, y_n)$ depends on $f$. In this sense, the double Schubert
polynomials are the \emph{universal polynomials} which yield the
classes of all suitably generic flagged degeneracy loci upon
substituting appropriate Chern classes for the variables.

\begin{Example} \label{ex:section}
    The simplest nontrivial example of
Theorem~\ref{thm:fulton-degeneracy-loci} has $w = s_{1}=21345\cdots n
\in S_{n}$. It is convenient to use the essential set to reduce the
number of rank conditions under consideration
(cf. \Cref{exer:essential}), but we must remember that $\Omega_w$ is
defined using opposite Schubert conditions, namely the Schubert
conditions for $w_0 w$. Our pictures therefore get flipped across a
horizontal axis of symmetry compared to \S\ref{sub:MatrixSchubs}, and
the appropriate set is the \emph{northwest} essential set
$\Ess_{\mathrm{NW}}(w)$: the set of \emph{southeasternmost} elements
of the connected components of the \emph{northwest} Rothe diagram
$D(w^{-1})$. Figure~\ref{fig:2134.nw} gives an example of the
northwest essential set for 4123, which
we will come back to later.  When $w=s_{1}$, the diagram $D(w^{-1})$
is a single box in the northwest corner, so $\Ess_{\mathrm{NW}}(w) =
\{(1,1)\}$ and $\rank(w)[1,1] = 0$ regardless of $n$. Hence,
\begin{equation*} \Omega_{w}(f) = \{x \in X : \rank(f_x : (E_1)_x \to
(F_1)_x) = 0\}.
\end{equation*}
Since $\fS_{w}(x_{1},x_{2},\dots ;y_{1},y_{2},\dots ) = \fS_{21}(x_{1};y_{1}) = x_1-y_1$,
Theorem~\ref{thm:fulton-degeneracy-loci} says that \begin{equation}
\label{eq:21-example} [\Omega_{w}(f)] = x_1-y_1 = c_1(\ker(F_1
\twoheadrightarrow F_0)) - c_1(E_1) =
c_1(F_1)-c_1(E_1).  \end{equation} if every component of
$\Omega_{w}(f)$ has codimension $\ell(21345\cdots n) = 1$.
    
    Although relatively simple, this example illustrates a few important aspects of Fulton's degeneracy locus formula.
    \begin{itemize}
    \item If $E_1$ is a trivial bundle, then \eqref{eq:21-example} becomes $[\Omega_w(f)] = c_1(F_1)$. In this case, $f$ is a section of the line bundle $F_1$ and $\Omega_{w}(f)$ is its zero locus, and the formula $[\Omega_w(f)] = c_1(F_1)$ agrees with Theorem~\ref{thm:chern-zeroes} interpreting Chern classes in terms of zero loci. Fulton's formula can therefore be regarded as a substantial generalization of Theorem~\ref{thm:chern-zeroes}.
    \item Although Fulton's flagged degeneracy locus formula assumes flags
    of bundles $E_\bullet$ and $F_\bullet$, only $E_1$ and $F_1$ were
    actually relevant. The example applies to any pair of line bundles
    $E_1$ and $F_1$, since we could extend them to flags trivially by
    setting $E_i = E_1 \oplus \C^{i-1}$ and $F_i = F_1 \oplus \C^{i-1}$.
    \item The degeneracy loci $\Omega_{21}(f), \Omega_{213}(f),
    \Omega_{2134}(f), \ldots$ over a fixed base $X$ are all exactly the
    same, so we can simplify matters by just considering $\Omega_{21}(f)$.  More generally, $\Omega_{w \times 1^m}(f)$ is independent of $m$ since $\Ess_{\mathrm{NW}}(w \times 1^m)$ is. This stability property provides some justification for our switch to opposite Schubert
    conventions in this subsection, as it would be rather less apparent if the indexing permutations were $12, 231, 3421,
    \ldots$.
    \end{itemize}

\begin{figure}
    \begin{center}
    \begin{tikzpicture}[scale=0.6]
    \draw[step=1.0,black,thin] (0,0) grid (4,4);
    \node at (0.5,0.5) {$1$};
    \node at (0.5,1.5) {$0$};
    \node at (0.5,2.5) {$0$};
    \node at (0.5,3.5) {$0$};
    \node at (1.5,0.5) {$2$};
    \node at (1.5,1.5) {$1$};
    \node at (1.5,2.5) {$1$};
    \node at (1.5,3.5) {$1$};
    \node at (2.5,0.5) {$3$};
    \node at (2.5,1.5) {$2$};
    \node at (2.5,2.5) {$2$};
    \node at (2.5,3.5) {$1$};
    \node at (3.5,0.5) {$4$};
    \node at (3.5,1.5) {$3$};
    \node at (3.5,2.5) {$2$};
    \node at (3.5,3.5) {$1$};
    \end{tikzpicture}
    \qquad
    \begin{tikzpicture}[scale=0.6]
    \draw[step=1.0,green,thin] (0,0) grid (4,4);
    \node at (0.5,0.5) {$\bullet$};
    \node at (1.5,3.5) {$\bullet$};
    \node at (2.5,2.5) {$\bullet$};
    \node at (3.5,1.5) {$\bullet$};
    \draw(4.5,0.5)--(0.5,0.5)--(0.5,-0.5);
    \draw(4.5,3.5)--(1.5,3.5)--(1.5,-0.5);
    \draw(4.5,2.5)--(2.5,2.5)--(2.5,-0.5);
    \draw(4.5,1.5)--(3.5,1.5)--(3.5,-0.5);
    \draw[fill=green, opacity=0.3] (0,1) -- (1,1) -- (1,4) -- (0,4) -- (0,1);
    \draw[very thick] (0,1) -- (1,1) -- (1,2) -- (0,2) -- (0,1);
    \end{tikzpicture}
    \end{center}
    \caption{The northwest rank table and northwest essential set for $4123$.} \label{fig:2134.nw}
    \end{figure}
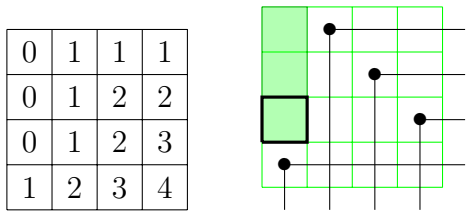

\end{Example}

\begin{Example} \label{ex:almost-isotropic}
    Fix a nondegenerate symmetric bilinear form $\langle \,\, , \,
\rangle$ on $\C^4$, like the standard dot product $\langle u, v
\rangle = \sum_i u_i v_i$. We use Fulton's formula to calculate the
cohomology class of the subvariety $Z = \{F_\bullet \in \Fl(4) : F_3 =
F_1^{\perp}\}$. Note that we are \emph{not} using a Hermitian inner
product where $F_1 \subseteq F_1^\perp$ would be impossible and our
algebraic tools would not apply. The form $\langle\,\, ,\, \rangle$
induces an linear isomorphism $\alpha : \C^4 \to (\C^*)^4$ by the
formula $\alpha(v)(w) = \langle v, w \rangle$. View $\alpha$ as a
morphism between the trivial bundles $\C^4 \times \Fl(4)$ and
$(\C^*)^4 \times \Fl(4)$, each flagged by tautological
bundles where
\begin{equation*} T_1 \subseteq T_2 \subseteq T_3 \subseteq T_4 = \C^4
\times \Fl(4) \xrightarrow{\alpha} (\C^*)^4 \times \Fl(4) = T_4^*
\twoheadrightarrow T_3^* \twoheadrightarrow T_2^* \twoheadrightarrow
T_1^*.
\end{equation*}
Each map $T_i^* \twoheadrightarrow T_{i-1}^*$ restricts linear functionals on $T_i$ to
the subbundle $T_{i-1}$. In terms of $\alpha$ and these maps, $Z$ is
the set of flags $F_\bullet$ such that $\rank(\alpha : F_3 \to F_1^*)
= 0$. We can write $Z$ as $\bigcup_w \Omega_w(\alpha)$ with the union
running over all $w \in S_4$ with $\rk(w)[3,1] = 0$, i.e. with $w_1 =
4$. Since $\Omega_w$ is defined using \emph{opposite} Schubert
conditions, we have $\Omega_v(f) \supseteq \Omega_w(f)$ whenever $v
\leq w$ in Bruhat order. The set $\{w \in S_4 : w_1 = 4\} =
\{4123,4132,4213,4231,4312,4321\}$ has $4123$ as its unique minimal
element in Bruhat order, so $Z = \Omega_{4123}(\alpha)$. The fact that
$\Ess_{\mathrm{NW}}(4123) = \{(3,1)\}$ and $\rk(4123)[3,1] = 0$ as
shown in \Cref{fig:2134.nw} confirms this calculation.

    Next we compute the appropriate Chern classes. First, $c_1(T_i/T_{i-1}) = -c_1((T_i/T_{i-1})^*) = -x_i$ is the negative of the ``standard'' generator $x_i \in H^*(\flags)$, by \Cref{ex:chern}(d,f). As for $c_1(\ker(T_i^* \twoheadrightarrow T_{i-1}^*))$, note that $\ker(T_i^* \twoheadrightarrow T_{i-1}^*)$ is the subspace of linear functionals on $T_i$ vanishing on $T_{i-1}$, which is naturally isomorphic to $(T_i/T_{i-1})^*$. Hence $c_1(\ker(T_i^* \twoheadrightarrow T_{i-1}^*))$ is $c_1((T_i/T_{i-1})^*) = x_i$. \Cref{thm:fulton-degeneracy-loci} then says 
    \begin{equation*}
[Z]=  \fS_{4123}(x_1, \ldots, x_4; -x_1, \ldots, -x_4) = 2x_1(x_1+x_2)(x_1+x_3),
    \end{equation*}
    at least provided that every irreducible component of $Z$ has codimension $\ell(4123) = 3$, or equivalently dimension $\dim \Fl(4) - \ell(4123) = {6 \choose 2} - 3 = 3$.

We give a parameterization of $Z$ from which we can compute $\dim Z$.
First, note that the subspaces $F_1$ and $F_3$ of any $F_\bullet \in
Z$ are completely determined by the choice of a line $F_1 \in \C\bP^3$
satisfying $\langle F_1, F_1 \rangle = 0$, since $F_3$ is determined
as $F_1^\perp$. The space of appropriate lines $F_1$ has dimension 2
because it is defined by one equation on $\C\bP^3$. For example, the
line spanned by $(1,i,0,0)$ is perpendicular to itself.  Once $F_1$
has been fixed, a full flag $F_\bullet \in Z$ is determined by the
choice of a 2-plane $F_2$ containing $F_1$ and contained in
$F_1^\perp$. This is equivalent to choosing a line in the
$2$-dimensional space $F_1^\perp/F_1$, i.e. a point of $\C\bP^1$. In
total, this shows $\dim Z = 2+1 = 3$.  We leave it to the reader to
verify $Z$ is irreducible.
\end{Example}

\begin{Example} \label{ex:isotropic}
    We give a non-example continuing on from
    \Cref{ex:almost-isotropic} which illustrates the importance of the
    expected codimension condition in
    \Cref{thm:fulton-degeneracy-loci}. A flag $F_\bullet \in \flags$
    is \emph{isotropic} if $F_i^\perp = F_{n-i   }$ for all $i \leq n/2$. Let $\IFl(n)$ be the subvariety of isotropic flags in $\flags$. Using the notation of \Cref{ex:almost-isotropic}, $\IFl(4)$ is the set of flags $F_\bullet$ with $\rank(\alpha : F_3 \to F_1^*) = \rank(\alpha : F_2 \to F_2^*) = 0$, and arguing as in \Cref{ex:almost-isotropic} gives $\IFl(4) = \Omega_{4312}(\alpha)$. The expected codimension is therefore $\ell(4312) = 5$.
    
    We claim $\IFl(4)$ is 2-dimensional. Indeed, the family of pairs
$(F_1 \subseteq F_3 = F_1^\perp)$ is 2-dimensional as explained in
\Cref{ex:almost-isotropic}, and completing such a pair to an isotropic
flag is equivalent to choosing a line $\ell$ in $F_3/F_1$ with
$\langle \ell, \ell \rangle = 0$. The latter is a nontrivial quadratic
equation on $\C\bP^1$ and its solution set is zero-dimensional so
$\IFl(4)$ is also 2-dimensional as claimed. Hence, $\codim \IFl(4) =
\dim \Fl(4) - 2 = 4$, so $[\IFl(4)]$ has degree 4 and cannot be
represented as any evaluation of the degree 5 double Schubert
polynomial $\fS_{4312}$.
\end{Example}

\begin{Example} \label{ex:schubert-locus}
    A Schubert variety is a basic example of a flagged degeneracy locus, although the rank conditions given in \eqref{eq:flagged-degen} are those of an opposite Schubert variety, so a little translation of conventions is required. Let $G_\bullet$ be a fixed flag, but viewed as a trivial flagged bundle over $\flags$, and let $T_\bullet$ be the flag of tautological bundles over $\flags$ as in Example~\ref{ex:chern}. Set $E_i = (\C^n/G_{n-i})^*$ and $F_j = T_j^*$, and let $f : E_n \to F_n$ be the identity. Then $f : E_i \to F_j$ is the map which takes a linear functional vanishing on $G_{n-i}$ and restricts it to $T_j$, so $f$ has image $(T_j / (G_{n-i} \cap T_j))^*$.

    Now consider the flagged degeneracy locus $\Omega_w(f)$. By the previous paragraph, this is the set of $H_\bullet \in \flags$ with 
    \begin{equation*}
        \dim (H_j / (G_{n-i} \cap H_j))^* = j - \dim(G_{n-i} \cap H_j)
	\leq \rk(w)[i,j], \ \forall i,j \in [n]
    \end{equation*}
    or equivalently, $\dim(G_{i} \cap H_j) \geq j - \rk(w)[n-i,j]$ for
    all $i,j \in [n]$. The upper left $(n-i) \times j$ corner of the
    permutation matrix $M_w$ and the lower left $i \times j$ corner
    together contain exactly $j$ many $1$'s, which shows $j - \rk(w)[n-i,j] = \rk(w_0 w)[i,j]$. Thus $\Omega_w(f)$ is the Schubert variety $X_{w_0 w}(G_\bullet)$. Finally,
    \begin{itemize}
        \item $c_1(\ker(F_i \twoheadrightarrow F_{i-1})) = c_1((T_i/T_{i-1})^*)$ is the standard generator $x_i \in H^*(\flags)$
        \item $c_1(E_i/E_{i-1}) = 0$ since the bundles $E_i$ are
	trivial bundles.
    \end{itemize}
    Theorem~\ref{thm:fulton-degeneracy-loci} therefore says $[X_{w_0 w}(G_{\bullet})] = \fS_w(x_1, \ldots, x_n; 0, \ldots, 0) = \fS_w(x_1, \ldots, x_n)$.
\end{Example}

\begin{Exercise}
Consider a dual version of Example~\ref{ex:schubert-locus}.  Set $E_i
= T_i$ and $F_j = \C^n / G_{n-j}$ for a fixed flag $G_\bullet$, all
over $\flags$ as before. Letting $f : E_n \to F_n$ be the identity,
show that the flagged degeneracy locus $\Omega_{w^{-1}}(f)$ is again
the Schubert variety $X_{w_0 w}(G_\bullet)$, hence 
\begin{equation}\label{eq:flag.degen.double}
[\Omega_{w^{-1}}(f)] =
\fS_w(x_1, \ldots, x_n).
\end{equation}
Also for comparison, prove \eqref{eq:flag.degen.double} using Fulton's
flagged degeneracy locus formula and the identity $c_1(L^*) = -c_1(L)$.
\end{Exercise}

\begin{Example}
    Suppose $A$ is an $n \times n$ matrix where each $A_{ij} = A_{ij}(x_0, \ldots, x_d)$ is a homogeneous polynomial of degree $a_i+b_j > 0$ where the $a_i, b_i$ are nonnegative integers. This condition guarantees that all minors of $A$ are homogeneous. For $w \in S_n$, let $\Omega_w$ be the set
    \begin{equation*}
        \{[x_0 : \cdots : x_d] \in \C\bP^d : \rk [A_{pq}(x_0, \ldots, x_d)]_{p \leq i, q \leq j} \leq \rk(w)[i,j]\, \forall i,j\}.
    \end{equation*}
    One can reason, as in Example~\ref{ex:degen}(d),  that specifying a matrix $A$ as described is equivalent to specifying a vector bundle morphism $f : \bigoplus_{i=1}^n \cO(-a_i) \to \bigoplus_{i=1}^n \cO(b_i)$ over $\C \bP^d$. Then $\Omega_w$ is the flagged degeneracy locus $\Omega_w(f)$, and the Schubert polynomial evaluation of Theorem~\ref{thm:fulton-degeneracy-loci} is 
\begin{equation}\label{eq:master.formula.doubleschub.degree}
 \fS_w(b_1, \ldots, b_n; -a_1, \ldots, -a_n)\alpha^{\ell(w)}
\end{equation}
where $\alpha \in H^2(\C \bP^n)$ is the cohomology class of a linear
hyperplane from Exercise~\ref{ex:projective.space}.
    This means that for generic $f$ (i.e.\ generic polynomials
    $A_{ij}$), the subvariety $\Omega_w(f) \subseteq \C\bP^d$ has
    degree $\fS_w(b_1, \ldots, b_n; -a_1, \ldots, -a_n)$ and
    codimension $\ell(w)$.

Here is an application of this degree calculation in
\eqref{eq:master.formula.doubleschub.degree} from
\cite{ED-degree}. Suppose $X \subseteq \bR^N$ is a real algebraic
variety of codimension $c$ and we wish to determine the closest points
on $X$ to some fixed $u = (u_1, \ldots, u_N) \in \bR^N$. These are the
minimizers of $\delta(x) = \sum_{i=1}^N (x_i-u_i)^2$ restricted to $X$
and are always critical points of $\delta$ on $X$. For general $u$ it
turns out there are a finite number of critical points to check---but
how many?
    
    The set $C$ of critical points can be characterized as the locus
where an appropriate Jacobian matrix of polynomials has prescribed
rank. Thus, $C$ is a degeneracy locus of the type considered in this
example, and if $C$ is 0-dimensional then its degree as a variety is
just the cardinality $|C|$. It can be shown using B\'ezout's theorem
and the calculation above that, if $X$ is defined by polynomial
equations $f_1 = \cdots = f_k = 0$ of degrees $d_1 \geq \cdots \geq
d_k$ and $c = \codim X$, then $|C| \leq d_1 \cdots d_c \sum_{i=0}^{N-c}
h_i(d_1-1,\ldots,d_c-1)$. Here, $h_i(x_1, \ldots, x_c)$ is the degree
$i$ homogeneous symmetric function in $c$ variables.
\end{Example}

\bigskip

We close this section by returning to the $r$\textsuperscript{th} 
degeneracy locus
\[ \Omega_r(f) = \{x \in X \mid \rank f_x \leq r\} \]
from Definition~\ref{def:simple-degen-locus} and giving a formula for
the cohomology class $[\Omega_r(f)]$ for any generic morphism of
vector bundles $f:E \to F$. Set $p = \rank
E, q = \rank F$, and suppose for the moment that $E = E_p$ and $F =
F_q$ fit into flags of bundles as in
Theorem~\ref{thm:fulton-degeneracy-loci}. Let $w(p,q,r) \in S_n$ be
the unique permutation in $S_{n}$ such that the (transposed) Rothe
diagram $D(w(p,q,r)^{-1})$ is the rectangle $[r+1,p] \times
[r+1,q]$. Explicitly, $w(p,q,r)$ is the $q$-Grassmannian permutation
\begin{equation*}
    [1,2,\ldots,r,p+1,p+2,\ldots,p+q-r,r+1,r+2,\dots , p, p+q-r+1,\ldots,n].
\end{equation*}
Since $\Ess_{\mathrm{NW}}(w(p,q,r)) = \{(p,q)\}$ and $\rk(w(p,q,r))[p,q] = r$, we have $\Omega_{r}(f) = \Omega_{w(p,q,r)}(f)$. If every component of $\Omega_r(f)$ has the expected codimension $\ell(w(p,q,r)) = (p-r)(q-r)$, then Fulton's formula gives 
\begin{equation} \label{eq:thom-porteous-1}
    [\Omega_r(f)] = \fS_{w(p,q,r)}(X,Y),
\end{equation}
    where $x_i = c_1(\ker(F_i \twoheadrightarrow F_{i-1}))$ and $y_i = c_1(E_i/E_{i-1})$.

This is not completely satisfying for several reasons. First, we might hope for a more explicit expression for $\fS_{w(p,q,r)}(X,Y)$. A more serious problem is that the bundles $E$ and $F$ may \emph{not} fit into a complete flag of bundles, making the analysis above impossible. The solution is to work with more general Chern classes and classifying maps.
\begin{Definition} \label{defn:higher-chern}
A \emph{classifying map} of a rank $k$ complex vector bundle $F$ over $X$ is a map $g : X \to \Gr(k,n)$ such that $g^*(T_k^*) = F$, where $T_k$ is the tautological bundle over $\Gr(k,n)$. The \emph{$d$\textsuperscript{th} Chern class} $c_d(F)$ is defined as $g^*(e_d)$ where the elementary symmetric polynomial $e_d(x_1,\ldots,x_k)$ corresponds to an element of $H^*(\Gr(k,n))$ as in Theorem~\ref{thm:grassmannian-cohom}.
\end{Definition}

This is completely analogous to Definition~\ref{def:chern} for the
first Chern class of a line bundle. If $X$ is a compact smooth
manifold, then classifying maps always exist and are unique up to
homotopy \cite[\S 5]{milnor-stasheff}.  So, the class $c_d(F) \in
H^d(X)$ is uniquely determined by Definition~\ref{defn:higher-chern}.

\begin{Proposition} \label{prop:elementary-chern}
    If $L_1, \ldots, L_n$ are line bundles, then
    \begin{equation*}
    c_d(L_1 \oplus \cdots \oplus L_n) = e_d(c_1(L_1), \ldots, c_1(L_n)).
\end{equation*}
\end{Proposition}
The reader might find it instructive to prove this proposition by constructing a
classifying map for $L_1 \oplus \cdots \oplus L_n$ from classifying
maps for the individual $L_i$, then using
Theorem~\ref{thm:grassmannian-cohom}. This statement is even more
useful than it appears because every short exact sequence of vector
bundles over a smooth manifold splits: if $E \subseteq F$ are vector
bundles over $X$,
then there is a subbundle $E' \subseteq F$ such that $F = E \oplus
E'$, so $E' \simeq F/E$. One way to prove this is to use the fact that
every smooth manifold has a Riemannian metric and take $E'$ to be the
orthogonal complement $E^{\perp}$. Now if $0 = E_0 \subseteq E_1
\subseteq \cdots \subseteq E_k$ is a flag of bundles, repeated use of this splitting technique shows that $E_k \simeq \bigoplus_{i=1}^k E_i/E_{i-1}$. Therefore
\Cref{prop:elementary-chern} gives $c_d(E_k) = e_d(X_k)$
where $x_i = c_1(E_i/E_{i-1})$ and $X_k$ denotes the variable set $\{x_1, \ldots, x_k\}$.

This framework allows for a completely general version of \eqref{eq:thom-porteous-1}.  Since $w(p,q,r)$ is $q$-Grassmannian and its inverse is $p$-Grassmannian, \Cref{prop:symmetry} and the identity $\fS_w(X;Y) = (-1)^{\ell(w)}\fS_{w^{-1}}(Y;X)$ imply that $\fS_{w(p,q,r)}(X;Y)$ is symmetric in $x_1, \ldots, x_q$ and in $y_1, \ldots, y_p$ separately. Therefore, $\fS_{w(p,q,r)}(X;Y)$ can be written uniquely as a polynomial in the elementary symmetric polynomials
\[ \{e_d(X_q) : 0 \leq d \leq q\} \qquad \text{and} \qquad \{e_d(Y_p):
0 \leq d \leq p\}. \]

\begin{Theorem} \label{thm:thom-porteous-bad}
Let $f:E \to F$
be a morphism of complex vector bundles, and set $p = \rank E$, $q = \rank F$.  
    If every component of $\Omega_r(f)$ has codimension $(p-r)(q-r) = \ell(w(p,q,r))$, then $[\Omega_r(f)]$ is obtained by making the substitutions $e_d(X_q) \mapsto c_d(F)$ and $e_d(Y_p) \mapsto c_d(E)$ in $\fS_{w(p,q,r)}$.
\end{Theorem}
If the vector bundles $E$ and $F$ fit into complete flags of bundles,
then this is equivalent to \eqref{eq:thom-porteous-1} by
Proposition~\ref{prop:elementary-chern}. However, the point of the
formulation in Theorem~\ref{thm:thom-porteous-bad} is that it holds for \emph{any} $E$ and $F$, regardless of whether they fit into complete flags.

\begin{Example} \label{ex:thom-porteous-1}
    Say $f : E \to F$ where $\rank(E) = 3$ and $\rank(F) = 2$. Then,
    to compute $[\Omega_1(f)]$, observe
    \begin{align*}
        \fS_{w(3,2,1)}(X;Y) &= \fS_{1423}(X;Y) = x_1^2 + x_1 x_2 + x_2^2 - x_1 y_1 - x_2 y_1 - x_1 y_2 - x_2 y_2 - x_1 y_3 - x_2 y_3\\
        & \qquad \qquad \qquad \qquad + y_1 y_2 + y_1 y_3 + y_2 y_3\\
        &= e_{1}(x_1,x_2)^2 - e_{2}(x_1,x_2) - e_1(x_1,x_2)e_1(y_1,y_2,y_3) + e_2(y_1,y_2,y_3).
    \end{align*}
    Thus, if $\Omega_1(f)$ has codimension $2 = \deg \fS_{1423} = \ell(1423)$, then
    \begin{equation*}
        [\Omega_1(f)] = c_1(F)^2 - c_2(F) - c_1(F)c_1(E) + c_2(E).
    \end{equation*}
\end{Example}

We can make Theorem~\ref{thm:thom-porteous-bad} more explicit
yet. Since $w(p,q,r)$ is vexillary, $\fS_{w(p,q,r)}(X;Y)$ is a
\emph{flagged double Schur function}
(cf. \Cref{cor:vexillary.Stanleys}), which satisfy a determinantal
formula in terms of elementary symmetric functions. To be precise, the
\emph{double elementary symmetric function} $e_d(X;Y)$ is
$\sum_{i=0}^d (-1)^{d-i} e_i(X)h_{d-i}(Y)$. Then, 
one can derive the following formula from \cite[(6.15)]{M2}:
\begin{equation} \label{eq:flagged-schubert}
    \fS_{w(p,q,r)}(X;Y) = \det [e_{q-r+j-i}(X_q;Y_p)]_{i,j \in [p-r]}.
\end{equation}
Let $c_k(F-E)$ denote the result of expanding $e_k(X_q; Y_p)$ as a polynomial in the elementary symmetric functions $\{e_d(X_q) : 0 \leq d \leq q\}$ and $\{e_d(Y_p) : 0 \leq d \leq p\}$, then making all substitutions $e_d(X_q) \mapsto c_d(F)$ and $e_d(Y_p) \mapsto c_d(E)$. Here $F-E$ should be regarded as just a formal symbol. Combining \eqref{eq:flagged-schubert} with Theorem~\ref{thm:thom-porteous-bad} gives, finally, the \emph{Thom-Porteous formula} expressing $[\Omega_r(f)]$ as a polynomial in the Chern classes of $E$ and $F$, all independent of $f$ assuming $f$ generic. 
\begin{Theorem}[\textbf{Thom-Porteous formula} \cite{Porteous}] \label{thm:thom-porteous} Suppose $f : E \to F$ is a morphism of complex vector bundles over a smooth complex projective variety $X$, and set $p = \rank E$, $q = \rank F$. If every component of $\Omega_r(f)$ has codimension $(p-r)(q-r) = \ell(w(p,q,r))$, then
    \begin{equation*}
    [\Omega_r(f)] = \det [c_{q-r+j-i}(F-E)]_{i,j \in [p-r]}.
    \end{equation*}
\end{Theorem}

\begin{Example}
Continue with the setup of Example~\ref{ex:thom-porteous-1}. We have 
\begin{equation*}
e_2(X_2;Y_3) = e_2(X_2) - e_1(X_2)h_1(Y_3) + h_2(Y_3) = e_2(X_2) - e_1(X_2)e_1(Y_3) + e_1(Y_3)^2 - e_2(Y_3),
\end{equation*}
which means $c_2(F-E) = c_2(F) - c_1(F)c_1(E) + c_1(E)^2 - c_2(E)$. Similarly, the expansions $e_1(X_2;Y_3) = e_1(X_2)-e_1(Y_3)$ and $e_0(X_2;Y_3)=1$ give $c_1(F-E) = c_1(F)-c_1(E)$ and $c_0(F-E) = 1$. The Thom-Porteous formula (Theorem~\ref{thm:thom-porteous}) now says that for generic $f$, the class $[\Omega_1(f)]$ equals
\begin{align*}
\det \begin{bmatrix} c_1(F-E) & c_2(F-E)\\ c_0(F-E) & c_1(F-E) \end{bmatrix} &= \det \begin{bmatrix} c_1(F) - c_1(E) & c_2(F) - c_1(F)c_1(E) + c_1(E)^2 - c_2(E) \\ 1 & c_1(F) - c_1(E) \end{bmatrix}\\
&= (c_1(F) - c_1(E))^2 - (c_2(F) - c_1(F)c_1(E) + c_1(E)^2 - c_2(E))\\
&= c_1(F)^2 - c_1(E)c_1(F) - c_2(F) + c_2(E),
\end{align*}
agreeing with Example~\ref{ex:thom-porteous-1}.
\end{Example}

Here is a useful special case of the Thom-Porteous formula (Theorem~\ref{thm:thom-porteous}), derived by letting $E$ be trivial of rank $1$ and $r = 0$. It provides a simple interpretation of the top Chern class of a vector bundle, generalizing Theorem~\ref{thm:chern-zeroes}.
\begin{Proposition} \label{prop:top-chern}
Suppose $F$ is a complex vector bundle of rank $q$ over a smooth projective variety $X$, with $f$ an algebraic section of $F$ and $Z$ its zero locus. If every component of $Z$ has codimension $q$, then $[X] = c_q(F)$.
\end{Proposition}

\begin{Example}\label{ex:27.lines}
A \emph{cubic surface} in $\C\bP^3$ is the vanishing set of a homogeneous cubic polynomial $s(x_0,x_1,x_2,x_3)$. A famous fact from classical algebraic geometry is that a smooth cubic surface contains exactly 27 lines. Here is a proof using Chern classes, modulo some technical details.

Recall our trick from \S\ref{sub:SchubertProblems2000.reprise} of identifying lines in $\C\bP^3$ with planes in $\C^4$, i.e.\ points in the Grassmannian $\Gr(2,4)$. Fix a smooth cubic surface $S = \{[x_0 : x_1 : x_2 : x_3] \in \C\bP^3 \mid s(x_0,\ldots,x_3) = 0\}$, and let $Z \subseteq \Gr(2,4)$ be the set of planes whose corresponding lines lie in $S$. We view $s$ in a coordinate-free way as an element of $\Sym^3((\C^4)^*)$. The restriction of $s$ to a plane $P \subseteq \C^4$ is an element $s|_P \in \Sym^3(P^*)$. This vector space has dimension 4: if $P^* = \langle y,z \rangle$, then $\Sym^3(P^*) = \langle y^3, y^2z, yz^2, z^3 \rangle$. With this notation,
\begin{equation}  \label{eq:section-27}
Z = \{P \in \Gr(2,4) : s|_{P} = 0\}.
\end{equation}

We can rephrase the above in terms of vector bundles. Let $T$ be the tautological bundle over $\Gr(2,4)$, and set $F = \Sym^3(T^*)$, a rank 4 vector bundle. Each $s|_{P}$ lies in the fiber $F_P$ by definition, so the map $P \mapsto s|_P$ is a section of $F$, and its zero locus is $Z$ by \eqref{eq:section-27}. By Proposition~\ref{prop:top-chern}, $[Z] = c_4(F)$ if every component of $Z$ has codimension $4$. Since $\dim \Gr(2,4) = 4$, this is the same as requiring that $Z$ be finite. Note that this condition is not true for every cubic surface, e.g.\ $x_0 x_1 x_2 = 0$ contains infinitely many lines. The fact that it holds for \emph{smooth} cubic surfaces is one technical detail we are eliding.

Assuming $Z$ is finite and none of its points occur with multiplicity (another missing technical detail!), the class $[Z] \in H^*(\Gr(2,4))$ is just $|Z|$ times the class of a point. Recall from \S\ref{sub:Grassmannians} that the Schubert variety $X_\lambda \subseteq \Gr(k,n)$ has dimension $|\lambda|$ and cohomology class represented by the Schur polynomial $s_{\lambda^{\vee}}(x_1,\ldots,x_k)$. In our case, this means the class of a point must be $[X_{\emptyset}] = s_{(2,2)}(x_1,x_2)$. Thus, we wish to show $[Z] = c_4(F) = 27[X_\emptyset]$.

By Theorem~\ref{thm:grassmannian-cohom}, the projection $\pi : \Fl(n) \to \Gr(k,n)$ induces an inclusion of rings $\pi^* : H^*(\Gr(k,n)) \to H^*(\Fl(n))$, which is in fact how we identified elements of $H^*(\Gr(k,n))$ with symmetric polynomials. This means nothing is lost by considering $\pi^* c_4(F) = c_4(\Sym^3(T_2^*))$ instead, where $T_j$ is the tautological rank $j$ bundle over $\Fl(4)$. The advantage of this perspective is that we can decompose $T_2^*$ into line bundles, specifically $T_2^* \simeq L_1 \oplus L_2$ where $L_1 = T_1^*$ and $L_2 = (T_2/T_1)^*$. 

Consider the general formula  
\begin{equation}  \label{eq:sym-decomp}
    \Sym^n(L_1 \oplus L_2) \simeq \bigoplus_{i=0}^n L_1^{\otimes i} \otimes L_2^{\otimes (n-i)}.
\end{equation}
This is nothing more than the fact that $\Sym^n(\langle y,z \rangle)$ has basis $\{y^i z^{n-i} : 0 \leq i \leq n\}$, applied fiberwise. Now by \eqref{eq:sym-decomp} and Proposition~\ref{prop:elementary-chern},
\begin{align*}
c_4(\Sym^3(T_2^*)) &= e_4(c_1(L_1^{\otimes 3}), c_1(L_1^{\otimes 2} \otimes L_2), c_1(L_1 \otimes L_2^{\otimes 2}), c_1(L_2^{\otimes 3}))\\
&= e_4(3x_1, 2x_1+x_2, x_1+2x_2, 3x_2) \qquad \text{(by Example~\ref{ex:chern}(e))}\\
&= 9x_1x_2(2x_1+x_2)(x_1+2x_2) = 18x_1^3 x_2 + 18x_1 x_2^3 + 54 x_1^2 x_2^2\\
&= 18s_{31}(x_1,x_2) + 27s_{22}(x_1,x_2).
\end{align*}
The first term is zero in $H^*(\Gr(2,4))$ since the partition $(3,1)$ does not fit in a $2\times 2$ rectangle  (Theorem~\ref{thm:grassmannian-cohom}), so $[Z] = c_4(F) = 27[X_{\emptyset}]$ as desired.
\end{Example}

\subsection{Singular Loci  and Permutation Patterns}\label{sub:Singularloci}

Say we wish to determine which Schubert varieties in $\flags$ or any
partial flag variety are smooth, and which are not.  Questions along
these lines were first addressed around the 1970s by Chevalley,
Demazure, Lakshmibai, Musili, and Seshadri, and were a major interest at
the Tata Institute around that time.  There are several combinatorial
and geometrical observations which make this determination easier to
characterize than for a typical variety.  In this section, we will give an
overview of the beautiful interplay between the combinatorics, algebra
and geometry involved in the study of singular loci of Schubert
varieties.  This subsection is based on the expository article by Abe
and Billey called ``Consequences of the Lakshmibai-Sandhya Theorem:
the ubiquity of permutation patterns in Schubert calculus and related
geometry'' \cite{Abe-Billey}.  For more details, see also the books
\cite{BLak} and \cite{kumar-book}.

For background, an affine variety is \textit{smooth at a point} if the
dimension of its tangent space equals the dimension of the variety
near that point \cite[Ch.9.6]{Cox-Little-OShea}.  If the variety is
given in terms of the vanishing set of certain polynomials, then one
can check the dimension of the tangent space by computing the rank of
the Jacobian matrix for those polynomials evaluated at the point.  The
rank is smaller than expected if and only if all minors of a certain
size vanish.  A point that is not smooth is \textit{singular}.  Thus,
the set of points where the variety is not smooth is itself a variety
called the \textit{singular locus}. Smooth/singular points in
projective varieties are defined similarly since locally near a point
every projective variety looks like an affine variety.

A priori, to determine if a variety is smooth at every point, one must
check the dimension of the tangent space at every point.  For Schubert
varieties, we make some easy observations that turn this problem into
a finite verification.  Recall from \eqref{prelim 100} that the
Schubert variety $X_{w}$ is equal to the union $\bigcup_{v\leq w}
C_{v}$ as a subset of $\flags$.  As in \Cref{sub:Schubert.Varieties},
we will assume the Schubert cells and Schubert varieties are defined
with respect to the standard flag $E_{\bullet}$.  Each Schubert cell
$C_{v}$ is a $B$-orbit.  Thus, every point in the same $B$-orbit looks
locally the same as every other point in its orbit as viewed in
$X_{w}$ or in the complete flag variety.  Hence the dimension of the
tangent space to $X_{w}$ at any point in a $B$-orbit $C_{v} \subset
X_{w}$ has the same dimension.  Therefore, a point $p\in C_{v}\subset
X_{w}$ is singular in $X_{w}$ if and only if every point in $C_{v}$ is
singular in $X_{w}$.  Furthermore, if $X_{w}$ is singular at any $p\in
C_{v}$, then it is also singular at every point in
$X_{v}=\bigcup_{u\leq v} C_{u}$ since the singular locus is a closed
subvariety of $X_{w}$.  This reasoning implies two important facts.
\begin{enumerate}
\item The Schubert variety $X_{w}$ is smooth if and only if $X_{w}$ is
smooth at the standard flag $E_{\bullet}$, which is represented by the
identity matrix.

\item The singular locus of a Schubert variety $X_{w}$ is the union of
certain Schubert varieties $X_{v}$ for $v<w$.  
\end{enumerate}
One can check if $X_{w}$ is singular at the identity matrix $I$, or
equivalently the standard flag $E_{\bullet}$, by writing down the
determinantal equations defining $X_{w}$ for an affine neighborhood of
$X_{w}$ around $I$ and then checking the rank of the Jacobian
matrix of these defining polynomials evaluated at $I$.  One example of
an affine neighborhood of $I$ is given by $X_{w} \cap w_{0}
C_{w_{0}}$, and the equations of $X_{w}$ were given in
\Cref{cor:Schubert.variety.equations}.  However, there are many
equations to consider, so this is not as efficient as possible.  There
is another way which provides a more unified tool for the study of
singularities of Schubert varieties using Lie algebras.  This line of
thought leads to a complete description of the singular locus in
\Cref{thm:sing.locus} given toward the end of the section.

Recall from Section~\ref{sub:flags} that the flag variety can be
identified with the quotient
\begin{align*}
\flags = GL_{n}(\mathbb{C})/B = SL_{n} / (B\cap SL_{n}).
\end{align*}
The tangent space of $SL_n$ at the identity is isomorphic as a complex vector
space to its Lie algebra, which is known to be the $n \times n$ trace
zero matrices over $\mathbb{C}$ \cite{Hum2}.  The Lie algebra of
$B\cap SL_{n}$ is the subalgebra in $\mathrm{Lie}(SL_{n})$ of upper
triangular matrices with trace zero.  Let $\mathfrak{g}
=\mathrm{Lie}(SL_{n})$ and $\mathfrak{b} =\mathrm{Lie}(B\cap
SL_{n})$. Then the tangent space of $G/B$ at the identity matrix is
isomorphic to $\mathfrak{g}/\mathfrak{b}$ as a vector space. Denoting
by $E_{i,j}$ the $n\times n$ matrix with 1 in the $(i,j)$-entry and
0's elsewhere, we obtain a basis for $\mathfrak{g}/\mathfrak{b}$ by
\[
\mathfrak{g}/\mathfrak{b} = \mathrm{span}\{E_{j,i}\given 1\leq i<j \leq n \}.
\]
Observe that there is a natural bijection between the basis elements
$\{E_{j,i}\given i<j \}$ and $R:=\{t_{i,j}\given i<j \}$ the set of
transpositions in $S_{n}$, also known as \textit{reflections}.

More generally, for any $v \in S_n$, we will think of
$v_{\bullet}=(e_{v(1)},\dots , e_{v(n)})$ as a point in $G/B$
represented by its permutation matrix $M_{v}$.  The coset $M_{v}B \in
G/B$ is fixed by the left action of the subgroup $M_{v}BM_{v}^{-1}
\subset G$, so $M_{v}( G/B) M_{v}^{-1} = G/(M_{v} B M_{v}^{-1})$ is an
isomorphic copy of the flag variety $G/B$ but with respect to the base
flag $v_\bullet$.  Therefore, the tangent space to $G/B$ at
$v_{\bullet}$ is
\begin{align}\label{geom of Schubert 100}
T_{v}(G/B) \simeq M_{v} \left(\mathfrak{g}/\mathfrak{b} \right)  M_{v}^{-1}
\end{align}
with basis $\{ E_{v(j),v(i)} \given 1\leq i<j\leq n \}$.

\begin{Exercise}
Let $v \in S_{n}$.  Prove that for any $1\leq i<j\leq n$, we have the identities
\begin{align*}
&M_{v} \ E_{ij}\ M_{v}^{-1}= E_{v(i),v(j)}, \\
&t_{v(i),v(j)}\ v = v\ t_{ij} 
\end{align*}
for all $v \in S_{n}$.  
\end{Exercise}

\begin{Exercise}\label{exercise:T-fixed points}
Let $T$ be the $n\times n$ invertible diagonal matrices.  Prove that the
permutation matrices are in bijection with the flags in $\Fl(n)$ that
are fixed under left multiplication by $T$.  
\end{Exercise}

The next theorem gives us an explicit description of a basis of the
tangent space of each Schubert variety at the permutation matrices, or
equivalently the $T$-fixed points by \Cref{exercise:T-fixed points}.
Thus, it can be used to identify the singular locus of $X_{w}$.

\begin{Theorem} [Lakshmibai-Seshadri \cite{Lak-Sesh.1984}]
\label{geom of Schubert 150}
For $v\leq w \in S_n$, the tangent space of $X_{w}$ at $v$ is given by
\begin{align*}
T_{v}(X_{w}) & \cong \mathrm{span}\{E_{v(j),v(i)} \given i<j,\ \ vt_{ij} \leq w \}, 
\end{align*}
and hence we obtain
\begin{align*}
\mathrm{dim} \ T_{v}(X_{w}) &= \# \{(i<j)\given vt_{ij} \leq w \}.
\end{align*}
\end{Theorem}
\begin{proof}
Assume $1\leq i<j \leq n$. By definition of the tangent space, $E_{v(j),v(i)} \in T_{v}(X_{w})$ if and only if $(I +
\varepsilon E_{v(j),v(i)} ) v_{\bullet} \in X_{w}$ for infinitesimal
$\varepsilon>0$ where we can assume $\varepsilon^2 =0$.  Think of $(I
+ \varepsilon E_{v(j),v(i)} )$ as a matrix in $G$ acting on the left
of the flag $v_\bullet$ by moving the flag a little bit in the
direction of $E_{v(j),v(i)}$.  Hence,  $E_{v(j),v(i)} \in
T_{v}(X_{w})$ if and only if $M_{v} + \varepsilon E_{v(j),i} \in X_{w}$
because as matrices 
\begin{equation}\label{eq:v.computation}
(I + \varepsilon E_{v(j),v(i)} )M_{v} \ = \ M_{v} + \varepsilon E_{v(j),v(i)}M_{v}
\ = \ M_v + \varepsilon E_{v(j),i}.
\end{equation}

Recall, $v(i)>v(j)$ if and only if $vt_{ij}<v$.  If $v(i)>v(j)$, then
$M_v + \varepsilon E_{v(j),i} \in X_{v}$ by the rank conditions defining
a Schubert variety given in \Cref{def:Schubert.variety}.  Since $v
\leq w$ implies $T_{v}(X_{v})\subset T_{v}(X_{w})$, we see that
$E_{v(j),v(i)}$ is in $T_{v}(X_{w})$ whenever $vt_{ij}<v$ by the
claim above.

On the other hand, if $v(i)<v(j)$ then $M_v + \varepsilon E_{v(j),i} \in
C_{vt_{ij}}$.  So $E_{v(j),v(i)} \in T_{v}(X_{w}) $ if and only if
$C_{vt_{ij}} \subset X_{w}$, which happens if and only if $vt_{ij}
\leq w$.  Thus, in either case $E_{v(j),v(i)} \in T_{v}(X_{w}) $ if
and only if $vt_{ij} \leq w$, which implies
\begin{equation}\label{eq:v.tang.space}
 \mathrm{span}\{E_{v(j),v(i)} \given i<j,\ \ vt_{ij} \leq w \} \subset
T_{v}(X_{w}).
\end{equation}

To prove $T_{v}(X_{w}) \subset \mathrm{span}\{E_{v(j),v(i)} \given
i<j,\ \ vt_{ij} \leq w \}$, assume there exist coefficients $a_{i,j}$
for $1\leq i<j\leq n$ such that $\sum a_{i,j} E_{v(j),i} \in
T_{v}(X_{w})$.  By definition of the tangent space, this implies $M_v + \varepsilon \sum a_{i,j}
E_{v(j),i} \in X_{w}$.  Say $M_v + \varepsilon \sum a_{i,j} E_{v(j),i}
\in C_{v'}$ for some $v'\leq w $.  Since $\varepsilon\ll 1$, none of the
minors in $M_{v}$ which are nonzero will vanish in $M_{v} +
\varepsilon \sum a_{i,j} E_{v(j),i} $, so the rank table for $M_{v} +
\varepsilon \sum a_{i,j} E_{v(j),i} $ dominates the rank table for
$M_{v}$ in every position.  Hence, $v\leq v'\leq w$.  Thus, for each
$a_{i,j}\neq 0$, we have $M_v + \varepsilon E_{v(j),i} \in X_{w}$ so
$E_{v(j),v(i)} \in T_{v}(X_{w})$ by the claim above, which in turn
occurs if and only if $vt_{ij} \leq w$.  Therefore, $\sum a_{i,j}
E_{v(j),i} \in \mathrm{span}\{E_{v(j),v(i)} \given i<j,\ \ vt_{ij}
\leq w \}$.
\end{proof}

\begin{Corollary}\label{cor:comb.smoothness}
The Schubert variety $X_{w}$ is smooth at $v\in S_n$ if and only if
\begin{align*}
 \# \{(i<j)\given vt_{ij} \leq w \} =\ell(w)
\end{align*}
or equivalently if and only if 
\begin{align*}
\# \{(i<j)\given  v<vt_{ij} \leq w \} =\ell(w)-\ell(v).
\end{align*}
Therefore, the singular locus of $X_{w}$ is the union of Schubert
varieties 
\[
\mathrm{Sing}(X_{w})=\bigcup X_{v}
\]
over all $v\leq w$ in $S_{n}$ such that $ \# \{(i<j)\given vt_{ij}
\leq w \} > \ell(w)$.
\end{Corollary}

\begin{Remark}
We will return to the question of finding the singular locus of a
Schubert variety in \Cref{thm:sing.locus} using permutation patterns
and in \Cref{thm:BK-smooth-test} using double Schubert polynomials.
\end{Remark}

\begin{Example}\label{ex:sing.locus.n4}
  Consider the case $n=4$ and $w=4231$. \Cref{fig.id.4231.interval}
shows the Hasse diagram of the interval $[\id,w]$ in Bruhat order.
From the Hasse diagram, we can determine that the Schubert variety
$X_{4231}$ is not smooth at the point $v=2143$ because for all 6
transpositions in $t_{ij} \in S_{4}$, we have $vt_{ij}\leq w$, but
$\ell(w)=5$.  Also, we see $6=\#\{t_{ij}\leq 4231\} = \mathrm{dim}\
T_{\id}(X_{4231})>\ell(4231)=5$.  Similarly, by referring back to the
Hasse diagram in \Cref{ex:3412.Hasse.diagram}, one can check
$X_{3412}$ is not smooth at $v=1324$ and is smooth at all $v'\leq w$
such that $v' \not \leq v$.  It follows that
\begin{align*}
&\mathrm{Sing}(X_{4231}) = X_{2143}\\
&\mathrm{Sing}(X_{3412}) = X_{1324}.
\end{align*}
All other Schubert varieties $X_{w}$ for $w$ in $S_{4}$ are smooth.
Note that $3412=w_{0} 2143=2143 w_{0}$ and $4231=w_{0}1324=1324
w_{0}$; however, this is not a pattern that extends beyond $S_{4}$.

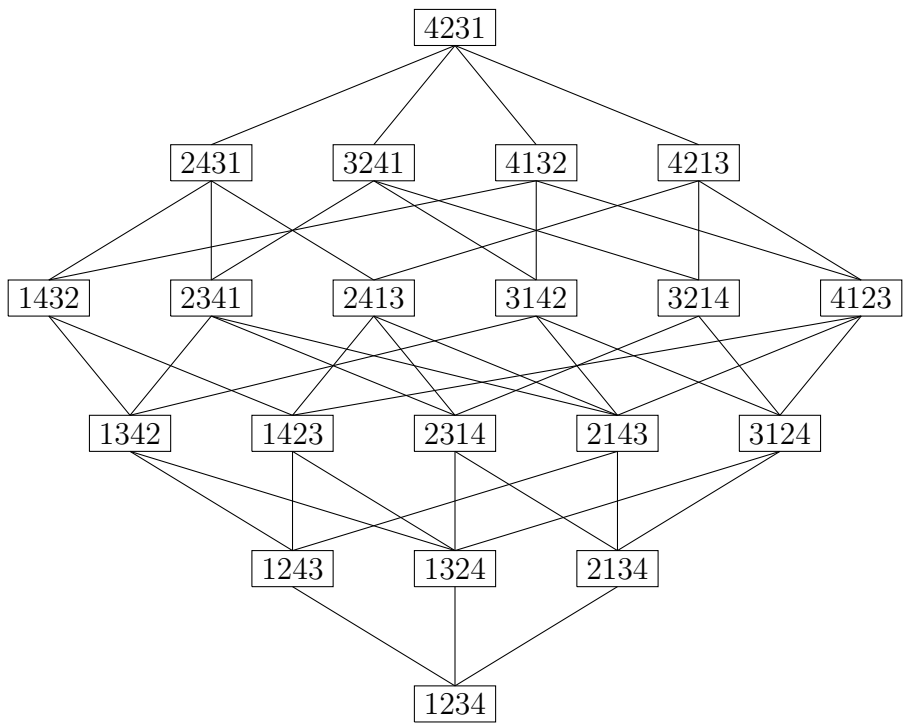
\begin{figure}
\centering
\begin{tikzpicture}[scale=0.6]
\def\a{0.9};
\def\b{0.4};
\def\h{3.0};
\newcommand\Rec[3]{
\node at (#1,#2) {#3};
\draw(#1-\a,#2-\b)--(#1-\a,#2+\b)--(#1+\a,#2+\b)--(#1+\a,#2-\b)--(#1-\a,#2-\b);
}
\Rec{0}{0}{$1234$}
\Rec{0}{\h}{$1324$}
\Rec{-4*\a}{\h}{$1243$}
\Rec{4*\a}{\h}{$2134$}
\Rec{-8*\a}{2*\h}{$1342$}
\Rec{-4*\a}{2*\h}{$1423$}
\Rec{0*\a}{2*\h}{$2314$}
\Rec{4*\a}{2*\h}{$2143$}
\Rec{8*\a}{2*\h}{$3124$}
\Rec{-10*\a}{3*\h}{$1432$}
\Rec{-6*\a}{3*\h}{$2341$}
\Rec{-2*\a}{3*\h}{$2413$}
\Rec{2*\a}{3*\h}{$3142$}
\Rec{6*\a}{3*\h}{$3214$}
\Rec{10*\a}{3*\h}{$4123$}
\Rec{-6*\a}{4*\h}{$2431$}
\Rec{-2*\a}{4*\h}{$3241$}
\Rec{2*\a}{4*\h}{$4132$}
\Rec{6*\a}{4*\h}{$4213$}
\Rec{0}{5*\h}{$4231$}
\draw(0,\b)--(-4*\a,\h-\b);
\draw(0,\b)--(0*\a,\h-\b);
\draw(0,\b)--(4*\a,\h-\b);
\draw(-4*\a,\h+\b)--(-8*\a,2*\h-\b);
\draw(-4*\a,\h+\b)--(-4*\a,2*\h-\b);
\draw(-4*\a,\h+\b)--(4*\a,2*\h-\b);
\draw(0*\a,\h+\b)--(-8*\a,2*\h-\b);
\draw(0*\a,\h+\b)--(-4*\a,2*\h-\b);
\draw(0*\a,\h+\b)--(0*\a,2*\h-\b);
\draw(0*\a,\h+\b)--(8*\a,2*\h-\b);
\draw(4*\a,\h+\b)--(0*\a,2*\h-\b);
\draw(4*\a,\h+\b)--(4*\a,2*\h-\b);
\draw(4*\a,\h+\b)--(8*\a,2*\h-\b);
\draw(-8*\a,2*\h+\b)--(-10*\a,3*\h-\b);
\draw(-8*\a,2*\h+\b)--(-6*\a,3*\h-\b);
\draw(-8*\a,2*\h+\b)--(2*\a,3*\h-\b);
\draw(-4*\a,2*\h+\b)--(-10*\a,3*\h-\b);
\draw(-4*\a,2*\h+\b)--(-2*\a,3*\h-\b);
\draw(-4*\a,2*\h+\b)--(10*\a,3*\h-\b);
\draw(0*\a,2*\h+\b)--(-6*\a,3*\h-\b);
\draw(0*\a,2*\h+\b)--(-2*\a,3*\h-\b);
\draw(0*\a,2*\h+\b)--(6*\a,3*\h-\b);
\draw(4*\a,2*\h+\b)--(-6*\a,3*\h-\b);
\draw(4*\a,2*\h+\b)--(-2*\a,3*\h-\b);
\draw(4*\a,2*\h+\b)--(2*\a,3*\h-\b);
\draw(4*\a,2*\h+\b)--(10*\a,3*\h-\b);
\draw(8*\a,2*\h+\b)--(2*\a,3*\h-\b);
\draw(8*\a,2*\h+\b)--(6*\a,3*\h-\b);
\draw(8*\a,2*\h+\b)--(10*\a,3*\h-\b);
\draw(-10*\a,3*\h+\b)--(-6*\a,4*\h-\b);
\draw(-10*\a,3*\h+\b)--(2*\a,4*\h-\b);
\draw(-6*\a,3*\h+\b)--(-6*\a,4*\h-\b);
\draw(-6*\a,3*\h+\b)--(-2*\a,4*\h-\b);
\draw(-2*\a,3*\h+\b)--(-6*\a,4*\h-\b);
\draw(-2*\a,3*\h+\b)--(6*\a,4*\h-\b);
\draw(2*\a,3*\h+\b)--(-2*\a,4*\h-\b);
\draw(2*\a,3*\h+\b)--(2*\a,4*\h-\b);
\draw(6*\a,3*\h+\b)--(-2*\a,4*\h-\b);
\draw(6*\a,3*\h+\b)--(6*\a,4*\h-\b);
\draw(10*\a,3*\h+\b)--(2*\a,4*\h-\b);
\draw(10*\a,3*\h+\b)--(6*\a,4*\h-\b);
\draw(-6*\a,4*\h+\b)--(0,5*\h-\b);
\draw(-2*\a,4*\h+\b)--(0,5*\h-\b);
\draw(2*\a,4*\h+\b)--(0,5*\h-\b);
\draw(6*\a,4*\h+\b)--(0,5*\h-\b);
\end{tikzpicture}
\caption{The Hasse diagram of the interval
$[\mathrm{id},4231]$}\label{fig.id.4231.interval}
%% formerly labeled: Interval [id,4231] in Bruhat order.
\end{figure}
\end{Example}

\begin{Exercise}
Find the singular locus of $X_{45312}$.  
\end{Exercise}

\subsubsection{Bruhat graphs}

The combinatorial data required to test smoothness of Schubert
varieties in \Cref{cor:comb.smoothness} gives rise to a graph which
contains the Hasse diagram of a Bruhat interval as a subgraph.  This
graph plays a pivotal role in the geometry of the Schubert varieties
and higher cohomology theories.

\begin{Definition}
For a permutation $w$, the \textit{Bruhat graph} for $w$ is a graph
whose vertex set is the \textit{identity interval} for $w$, denoted  $[\mathrm{id},w]=\{v\in
S_n\given v\leq w \}$, with an edge between $v$ and
$vt_{ij}$ if and only if both $v,vt_{ij} \leq w$.
\end{Definition}

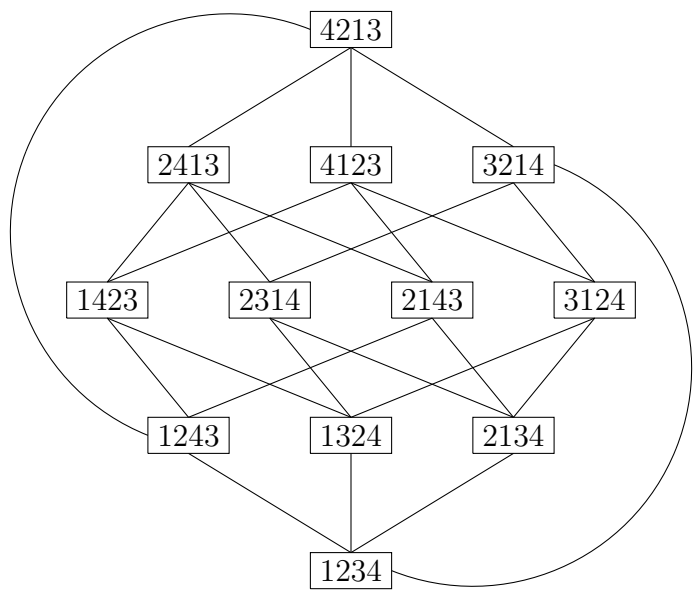
\begin{figure}[h]
\begin{tikzpicture}[scale=0.6]
\def\a{0.9};
\def\b{0.4};
\def\h{3.0};
\newcommand\Rec[3]{
\node at (#1,#2) {#3};
\draw(#1-\a,#2-\b)--(#1-\a,#2+\b)--(#1+\a,#2+\b)--(#1+\a,#2-\b)--(#1-\a,#2-\b);
}
\Rec{0}{0}{$1234$}
\Rec{0}{\h}{$1324$}
\Rec{-4*\a}{\h}{$1243$}
\Rec{4*\a}{\h}{$2134$}
\Rec{-6*\a}{2*\h}{$1423$}
\Rec{-2*\a}{2*\h}{$2314$}
\Rec{2*\a}{2*\h}{$2143$}
\Rec{6*\a}{2*\h}{$3124$}
\Rec{-4*\a}{3*\h}{$2413$}
\Rec{0*\a}{3*\h}{$4123$}
\Rec{4*\a}{3*\h}{$3214$}
\Rec{0*\a}{4*\h}{$4213$}

\draw(0,\b)--(-4*\a,\h-\b);
\draw(0,\b)--(0*\a,\h-\b);
\draw(0,\b)--(4*\a,\h-\b);
\draw(-4*\a,\h+\b)--(-6*\a,2*\h-\b);
\draw(-4*\a,\h+\b)--(2*\a,2*\h-\b);
\draw(0*\a,\h+\b)--(-6*\a,2*\h-\b);
\draw(0*\a,\h+\b)--(-2*\a,2*\h-\b);
\draw(0*\a,\h+\b)--(6*\a,2*\h-\b);
\draw(4*\a,\h+\b)--(-2*\a,2*\h-\b);
\draw(4*\a,\h+\b)--(2*\a,2*\h-\b);
\draw(4*\a,\h+\b)--(6*\a,2*\h-\b);
\draw(-6*\a,2*\h+\b)--(-4*\a,3*\h-\b);
\draw(-6*\a,2*\h+\b)--(0*\a,3*\h-\b);
\draw(-2*\a,2*\h+\b)--(-4*\a,3*\h-\b);
\draw(-2*\a,2*\h+\b)--(4*\a,3*\h-\b);
\draw(2*\a,2*\h+\b)--(-4*\a,3*\h-\b);
\draw(2*\a,2*\h+\b)--(0*\a,3*\h-\b);
\draw(6*\a,2*\h+\b)--(0*\a,3*\h-\b);
\draw(6*\a,2*\h+\b)--(4*\a,3*\h-\b);
\draw(-4*\a,3*\h+\b)--(0,4*\h-\b);
\draw(0*\a,3*\h+\b)--(0,4*\h-\b);
\draw(4*\a,3*\h+\b)--(0,4*\h-\b);
\draw[bend right=45](-\a,4*\h) to (-3*\a-1.5*\h,2.5*\h+2*\a);
\draw[bend right=45](-3*\a-1.5*\h,2.5*\h+2*\a) to (-5*\a,\h);
\draw[bend left=45](5*\a,3*\h) to (3*\a+1.5*\h,1.5*\h-2*\a);
\draw[bend left=45](3*\a+1.5*\h,1.5*\h-2*\a) to (\a,0);
\end{tikzpicture}
\caption{The Bruhat graph of $w=4213$}\label{Bruhat graph w=432}
\end{figure}

For example, the Bruhat graph of $w=4213$ is drawn in Figure
\ref{Bruhat graph w=432}.  Observe that the degree of $v$ (i.e. the
number of edges connected to $v$) in the Bruhat graph for $w$ is
$\mathrm{dim}\ T_{v}(X_{w})$.  Hence, this Bruhat graph is regular!

To describe the geometric interpretation of a Bruhat graph, let
$T\subset GL_{n}(\mathbb{C})$ be the set of invertible diagonal
matrices as above.  One can verify that the permutation matrices
exactly represent the $T$-fixed points in $GL_{n}/B$ and have the
following properties.
\begin{itemize}
\item[(i)] The permutations in $[\mathrm{id},w]$ are in bijection with
the $T$-fixed points of $X_{w}$ by \Cref{exercise:T-fixed points}.
\item[(ii)] If $v,vt_{ij}\leq w$, then the edge between $v$ and
$vt_{ij}$ in the Bruhat graph for $w$, corresponds to a 1-parameter
curve in $\flags$ passing through the flags $v_{\bullet}$ and
$(vt_{ij})_{\bullet}$ represented by matrices 
\[
L_{v} = \{M_{v}+ z E_{v(j),i}\given z \in \mathbb{C} \} \cup \{vt_{ij} \}   \simeq \mathbb{P}^{1}.
\]
This curve is $T$-invariant, and pointwise fixed by a torus $T'\subset T$ of codimension $1$.
\end{itemize}
Schubert varieties are examples of \textit{GKM-spaces} studied by
Goresky-Kottwitz-MacPherson \cite{GKM} and others.  It turns out that
much of the $T$-equivariant topology or geometry of GKM spaces can be
described in terms of their \emph{moment graph}.  The moment graph for a
Schubert variety $X_{w}$ is the Bruhat graph for $w$.   

\subsubsection{The Lakshmibai-Sandhya Theorem}\label{Lakshmibai-Sandhya Theorem}

There is a simple criterion for characterizing smooth Schubert
varieties using permutation pattern avoidance.  Recall from
\Cref{def:pattern-avoidance} that $w \in S_{n}$ \textit{avoids} a
pattern $v \in S_{k}$ for $k\leq n$ provided no subsequence $1\leq
i_{1}< i_{2}< \dots < i_{k}\leq n$ exists such that the inversions of
$[w(i_{1}), w(i_{2}),\dots , w(i_{k})]$ in terms of positions are the
same as for $v$.  We have seen some pattern avoidance
characterizations before for zero-one Schubert polynomials
(\Cref{thm:zero-one-schubert}) and vexillary permutations
(\Cref{cor:vexillary.Stanleys}).  Today, many families of permutations
are characterized by pattern avoidance or variations on that idea.

Lakshmibai-Sandhya proved the following criterion for the singularity
of Schubert varieties in 1990. See also the mutually independent work
by Haiman (unpublished), Ryan \cite{ryan}, and Wolper \cite{Wolper}.
This result has been highly influential on the field in terms of
bringing this technique into the study of Schubert varieties.

\begin{Theorem}[\textbf{Lakshmibai-Sandhya Theorem} \cite{Lak-San}]\label{Lak-San}
The Schubert variety $X_{w}$ is smooth if and only if $w$ avoids 3412
or 4231.
\end{Theorem}

\begin{proof}
Let us sketch one approach to proving \Cref{Lak-San} by applying
\Cref{geom of Schubert 150}.  Say $w$ contains a $3412$ or $4231$
pattern in positions $i_1<i_2<i_3<i_4$.  Let $v$ be the permutation
obtained from $w$ by rearranging the numbers $w_{i_1} w_{i_2} w_{i_3}
w_{i_4} $ according to the pattern for the corresponding singular
locus in $S_4$ as determined in \Cref{ex:sing.locus.n4}.
Specifically, if $w_{i_1} w_{i_2} w_{i_3} w_{i_4} $ is a
$4231$-pattern then replace $w_{i_1} w_{i_2} w_{i_3} w_{i_4} $ by the
$2143$-pattern $w_{i_2} w_{i_4} w_{i_1} w_{i_3} $ in the same
positions.  If $w_{i_1} w_{i_2} w_{i_3} w_{i_4} $ is a $3412$-pattern
then replace $w_{i_1} w_{i_2} w_{i_3} w_{i_4} $ by the $1324$-pattern
$w_{i_3} w_{i_1} w_{i_4} w_{i_2} $ in the same positions.  For
example, if $w=625431$ and we use the $6241$ instance of the pattern
$4231$, then $v=215634$ which contains a $2143$ pattern among the
values $1,2,4,6$.

We claim that $X_w$ is singular at the point $v$ by construction.  The
proof proceeds by comparing $\ell(w) - \ell(v)$ with the number of
$t_{ij}$ such that $v<vt_{ij}\leq w$.  For $i,j \in
\{i_1,i_2,i_3,i_4\}$, we know there will be strictly more such
transpositions than the length difference in these positions.  Since
$vt_{ij}$ and $w$ differ in at most 6 positions, we can apply
\Cref{ex:flatten.bruhat} to reduce the problem to a finite computer
verification on permutations in $S_{6}$, where we need to identify all
pairs $v,w$ such that
$$
\#\{t_{ij} \given  v < vt_{ij} \leq w\} > \ell(w) -\ell(v).
$$
It is a good exercise for the reader to do this computer verification
on your own, especially if computer assisted proofs are not part of
your standard repertoire of proof techniques.

In the other direction, assume that $w \in S_{n}$ avoids the patterns
$3412$ and $4231$. We use induction on permutation length to show that
$\#\{t_{ij}\given t_{ij}\leq w\}=\ell(w)$, which implies the standard
flag $E_{\bullet}$ is a smooth point in $X_{w}= X_{w}(E_{\bullet})$
and hence every point is smooth. The base case is when $w$ is the
identity permutation in which case the formula holds trivially.  If
$w(n)=n$, then we reduce to the case of $w \in S_{n-1}$, since all
permutations below such $w$ must have $n$ as a fixed point in Bruhat
order. So assume $w(n)<n$, and write $w(n)=n-k$ with $k\geq1$. We will
do a detailed analysis of what the permutation matrix $M_{w}$ looks
like when $w$ avoids 3412 and 4231.

\begin{figure}[h]
\centering
\begin{tikzpicture}[scale=0.6]
\draw[green, fill=green, fill opacity=0.4](0,3)--(9,3)--(9,4)--(0,4)--(0,3);
\draw[green, fill=green, fill opacity=0.4](2,0)--(3,0)--(3,9)--(2,9)--(2,0);
\draw[step=1.0,green,thin] (0,0) grid (9,9);
\node at (0.5,6.5) {$\bullet$};
\node at (1.5,8.5) {$\bullet$};
\node at (2.5,0.5) {$\bullet$};
\node at (3.5,1.5) {$\bullet$};
\node at (4.5,7.5) {$\bullet$};
\node at (5.5,2.5) {$\bullet$};
\node at (6.5,4.5) {$\bullet$};
\node at (7.5,5.5) {$\bullet$};
\node at (8.5,3.5) {$\bullet$};
\node at (1.5,4.5) {$A$};
\node at (3.5,4.5) {$B$};
\node at (1.5,2.5) {$C$};
\node at (3.5,2.5) {$D$};
\end{tikzpicture}
\caption{The matrix for the permutation $w=319827546$ with
$w'=31872654$ is shown with a dot in each position $(w(i),i)$ and
regions $A,B,C,D$ noted.  For example, $D=\{(8,4),(7,6) \}$.}
\label{fig:smooth-permutation-319827546}
\end{figure}
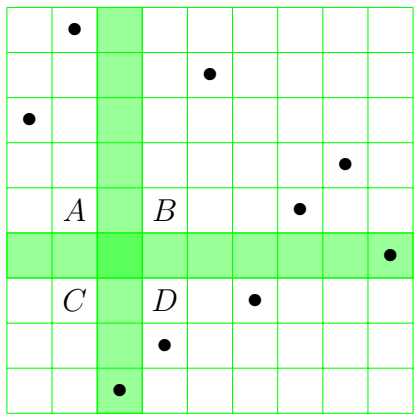

Consider the four regions of the matrix created by removing row $w(n)$
and column $w^{-1}(n)$, and define sets
\begin{align*}
A=&\{(w(i),i)\given i<w^{-1}(n),\ w(i)<w(n)\}, \\
B=&\{(w(i),i)\given i>w^{-1}(n),\ w(i)<w(n)\}, \\
C=&\{(w(i),i)\given i<w^{-1}(n),\ w(i)>w(n)\}, \\
D=&\{(w(i),i)\given i>w^{-1}(n),\ w(i)>w(n)\},
\end{align*}
as in \Cref{fig:smooth-permutation-319827546}.  In the figure, we put
a dot in each position $(w(i),i)$ for $1\leq i\leq 9$ for
visualization.  Since $w$ avoids $3412$, one of $B$ and $C$ must be
empty. By the symmetry of $w$ and $w^{-1}$, we can assume that
$C=\emptyset$. Moreover, as $w$ avoids $4231$, entries in $D$ must
correspond with a decreasing subsequence in the one-line notation for
$w$. In particular, when $C=\emptyset$, we know $D$ consists precisely
of all dots in rows $n-1,\ldots,n-k+1$.

Let $w'=\mathrm{fl}(w(1),\ldots,w(n{-}1))\in S_{n-1}$ be the
flattening of $w$ in the first $n-1$ indices. Since $w$ avoids $3412$
and $4231$, $w'$ also avoids $3412$ and $4231$. We see that
$\ell(w)-\ell(w')=k\geq 1$ since $w(n)=n-k$ creates precisely $k$
inversions with earlier indices.  By induction on permutation length,
we can assume $\#\{t_{ij}\given t_{ij}\leq w\}=\ell(w')$. It is also
straightforward to check that $\{t_{in}\given t_{in}\leq
w\}=\{t_{in}\given n-k\leq i<n\}$, and this set has cardinality $k$.
Hence, in order to prove $\#\{t_{ij}\given t_{ij}\leq w\}=\ell(w)$, it
suffices to show
\begin{equation}\label{eq:trans.below.w}
\{t_{ij}\given t_{ij}\leq w'\}=\{t_{ij}\given
i<j<n,t_{ij}\leq w\}.
\end{equation}

By \Cref{thm:Ehresmann}, the Ehresmann Tableau Criterion for the
Bruhat order, the ``$\subseteq$" direction of \eqref{eq:trans.below.w}
is immediate.  For the ``$\supseteq$" direction, consider any
$t_{ij}\leq w$ with $i<j<n$. If $t_{ij}\nleq w'$, then
\begin{eqnarray}\label{eq:double.gale.condition}
\{1,2,\ldots,i{-}1,j,i{+}1,\ldots,i{+}m\}&\trianglelefteq&\{w(1),\ldots,w(i{+}m)\}, \text{ and } \\
\{1,2,\ldots,i{-}1,j,i{+}1,\ldots,i{+}m\}&\ntrianglelefteq&\{w'(1),\ldots,w'(i{+}m)\} \label{eq:double.gale.condition.2}
\end{eqnarray}
in Gale order for some minimal $0\leq m<j-i$.  When $m=0$, the above
set on the left-hand side means $\{1,2,\ldots,i{-}1,j\}$ by
convention.  These conditions imply $i+m\geq w^{-1}(n)$, since
otherwise the two sets on the right-hand side are the same by
construction. Hence, $\max\{w'(1),\ldots,w'(i{+}m)\}=n-1$.  Since
$i+m<j\leq n-1=\max\{w'(1),\ldots,w'(i{+}m)\}$, in order for
\eqref{eq:double.gale.condition} to hold, we must have that
$\{w'(1),\ldots,w'(i{+}m)\}$ contains $\{1,\ldots,i\}$. At the same
time, $\{w(1),\ldots,w(i{+}m)\}$ does not contain all of
$\{1,\ldots,i\}$ by \eqref{eq:double.gale.condition.2}. Because $w'$
was obtained by flattening $w(1),\ldots,w(n{-}1)$, the only possible
explanation is that
\[
1,2,\ldots,i{-}1,i{+}1\in \{w(1),\ldots,w(i{+}m)\},
\]
$i\notin \{w(1),\ldots,w(i{+}m)\}$, and the dot in row $i+1$ of the
permutation matrix $M_{w}$ must be in the region $D$. Since entries in
the region $D$ are decreasing, if $i+1\in \{w(1),\ldots,w(i{+}m)\}$,
then any $i'\geq i+1$ is also in
$\{w(1),\ldots,w(i{+}m)\}$. Therefore, $[n]\setminus\{i\}\subset
\{w(1),\ldots,w(i{+}m)\}$, so $i+m\geq n-1$. However, $i+m\geq n-1$
contradicts the hypotheses $i+m<j<n$. Thus, $t_{ij} \leq w'$, and the
argument is now complete.
\end{proof}

Haiman's proof also contained the following enumerative formula as a
corollary.  Since his paper was never published due to the overlapping
content with \cite{Lak-San}, it wasn't until 2007 that this result had
a proof in the literature, due to Bousquet-M\'elou and Butler
\cite{BousquetMelou-Butler}.  

\begin{Corollary}\cite{BousquetMelou-Butler}\label{cor:smooth.perms}
  There is a closed form for the generating function $V(t)$ for the sequence counting the number of permutations $w \in S_{n}$ such that 
the Schubert variety $X_{w}$ is smooth:
\begin{align}\label{gf:haiman}
V(t) &= \ \frac{1-5t+3t^{2}+t^{2}\sqrt{1-4t}}{1-6t +8t^{2} -4 t^{3}}\\ \notag
 &= \ 1+ t+2t^{2}+6t^{3}+22t^{4}+88t^{5} \\ \notag
 &\qquad\qquad\qquad\qquad+366t^{6}+1552t^{7}+6652t^{8}+
 O(t^{9}).
\end{align}
\end{Corollary}

Note that by the Lakshmibai-Sandhya theorem, testing for smoothness of
Schubert varieties can be done naively in polynomial time, $O(n^{4})$,
based on the characterization of avoiding $3412$ and $4231$.
Remarkably, Guillemot and Marx \cite{Guillemot.Marx} showed that for
every permutation $v \in S_k$ there exists an algorithm to test if $w
\in S_n$ contains $v$ which runs in \textit{linear time}, $O(n)$!
This is a major improvement over brute force verification.  It is
often far from obvious that an $O(n)$ time algorithm exists for the
geometric or algebraic properties characterized by pattern avoidance.
However, combining the Lakshmibai-Sandhya and
Guillemot-Marx theorems, we know the problem of determining if a Schubert
varieties $X_{w}$ for $w \in S_{n}$ is smooth, denoted
(SmoothSchubert?), has linear time complexity, $O(n)$.

Another major advantage of permutation pattern characterizations is
that they provide efficient fingerprints for theorems
\cite{billey.tenner}.  Tenner's Database of Permutation Pattern
Avoidance (DPPA) provides a growing collection of known properties
characterized by patterns with references to the literature
\cite{dppa}.  This allows researchers to connect new theorems and
conjectures with known results in a format free of language or
notational differences.

%%%%%%%%%%%%%%%%%%%%%%%%%%%%%%%%%%
%%\section{10 Pattern Avoidance Properties}\label{s:properties}
%%%%%%%%%%%%%%%%%%%%%%%%%%%%%%%%%%

The Lakshmibai-Sandhya \Cref{Lak-San} from 1990 had a major impact on
the field.  Many researchers initiated a form of machine learning for
Schubert geometry to find permutation patterns characterizing
geometrical and algebraic concepts.  In \cite{Abe-Billey}, many
properties from Schubert geometry and related areas characterized by
pattern avoidance are spelled out with references and examples.  Here
we will just focus on the primary example of \textit{smooth
permutations} and all of the amazing properties related to them and
the Lakshmibai-Sandhya theorem.  We will give 10 such properties
first, and then given the notation and references after the theorem.

\begin{Theorem} \label{thm:smoothness}
The following are equivalent for $w \in S_n$.
\begin{enumerate}
\item The permutation $w$ avoids 3412 and 4231.
\item $X_{w}$ is smooth.
\item The tangent space dimension test holds:  $\ell(w) =
\#\{t_{ij}\leq w \}$.
\item \label{i:graph} The Bruhat graph for $w$ is regular and every vertex has degree $\ell(w)$.
\item \label{i:palindromic} The Poincar\'{e} polynomial for $w$, $P_{w}(t)= \sum_{v \leq w} t^{\ell(v)}$ is palindromic.
\item \label{i:factor} The Poincar\'{e} polynomial for $w$ factors as
\[
P_{w}(t)=\prod_{i=1}^{k}(1+t+t^{2}+\dots +t^{e_{i}})
\]
for positive integers $\{e_1, e_2, ..., e_k\}$ called
\textit{generalized exponents}, with $\ell(w)
= \sum e_i$.  
\item \label{i:hyper} 
  The Poincar\'{e} polynomial $P_w(t)$ is equal to the generating function
  $R_w(t)$ for the number of regions $r$ in the complement of the
  inversion hyperplane arrangement weighted by the distance of each
  region to the fundamental region.  In symbols,
\[
R_{w}(t)=\sum_r t^{d(r)} = \sum_{v \leq w} t^{\ell(v)} = P_{w}(t).
\]
Here, $d(r)$ is the number of hyperplanes crossed in a walk starting
at the fundamental region and going to the region $r$.
\item \label{i:hyper.free} 
The inversion arrangement for $w$ is free and the number of
  chambers of the arrangement is equal to the size of the Bruhat
  interval $[\id,w]$.
\item \label{i:KL} 
   The Kazhdan-Lusztig polynomial $P_{v,w}(t)= 1$ for all $v\leq w$.
\item \label{i:KL.id} 
   The Kazhdan-Lusztig polynomial $P_{\id,w}(t)= 1$.  
\end{enumerate}
\end{Theorem}

We have already discussed the equivalence of the first three items in
\Cref{geom of Schubert 150} and \Cref{Lak-San}.
Items~(\ref{i:graph}), (\ref{i:palindromic}), (\ref{i:KL}),
and~(\ref{i:KL.id}) are due to Carrell and Peterson \cite{carrell94}.
We will define the Kazhdan-Lusztig polynomials $P_{v,w}(t)$ below and
expand on some of their key properties.  Note, Carrell is the sole
author on the paper cited, but he always acknowledges Peterson as a
collaborator on this work so we give them both credit.  The term
\textit{palindromic} refers to the sequence of coefficients of the
polynomial, so the coefficient of $t^i$ equals the coefficient of
$t^{d-i}$ in a palindromic polynomial of degree $d$.

Item~(\ref{i:factor}) about factoring Poincar\'{e} polynomials is due
to Gasharov \cite{gasharov97}.  This factorization implies that the
geometry of smooth Schubert varieties has particularly nice structure
in terms of iterated fiber bundles over Grassmannians
\cite{GR2000,richmond.slofstra.2014,richmond.slofstra.2016,ryan,Wolper}.

\begin{Example} The permutation $w=4321$ avoids the patterns $3412$ and
  $4231$.   It has a palindromic Poincar\'e polynomial that also
  factors nicely, 
\begin{align*}
P_{4321}(t) 
&= (1+t)(1+t+t^{2})(1+t+t^{2}+t^{3})\\
&= 1+3t+5t^2+6t^3+5t^4+3t^5+t^6.
\end{align*}
\end{Example}

\begin{Example}
Recall from \Cref{ex:sing.locus.n4} or \Cref{Lak-San} that $3412$ and
$4231$ index the only two singular Schubert varieties in $S_4$.  Therefore, we
expect in each case the Poincar\'e polynomial does not have the nice
factorization, nor the palindromic property.  Here
$$P_{3412}(t) = 1+3t+5t^2+4t^3+t^4.
$$
What is $P_{4231}(t)$?
\end{Example}
\vspace{0.1in}

Item~(\ref{i:hyper}) about the inversion hyperplane arrangement is due
to Oh-Postnikov-Yoo \cite{OPY}.  This arrangement is the
collection of hyperplanes in $\mathbb{R}^{n}$ defined by
$x_{i}-x_{j}=0$ for all $1\leq i<j\leq n$ such that $w(i)>w(j)$.  This
generalizes the notion of the Coxeter arrangement of type $A_{n-1}$
given by all the hyperplanes $x_{i}-x_{j}=0$ for all $i<j$, so it is
the inversion arrangement for $w_0 \in S_{n}$.  The Coxeter
arrangement has $n!$ regions corresponding to all permutations.
In this case, the statistic $d(w)$ equals $\ell(w)$.  No explicit
bijective proof of Item~(\ref{i:hyper}) is known.  See Chapter 3 of
this book by Oh and Richmond for more \cite{oh.richmond.2024}.

Item~(\ref{i:hyper.free}) is due to Slofstra \cite{slofstra.2015}.
Here a central hyperplane arrangement in a Euclidean space $V$ is said
to be \textit{free} if the module of derivations of the complexified
arrangement is free as a module over the polynomial ring
$\mathbb{C}[V_{\mathbb{C}}]$.  Note he also gives an algebraic
interpretation for the generalized exponents in terms of degrees of a
homogeneous basis for the module of derivations.

As discussed below, the Kazhdan-Lusztig polynomials $P_{v,w}(q)$ are
closely related to the singularities in Schubert varieties.  They also
play an important role in representation theory and the study of Verma
modules.  Vazirani has a list of applications on her website
\cite{Vazirani}. We recall the definitions here, highlight some
important developments, and refer the reader to the textbooks by
Humphreys \cite{Hum} and Bj\"orner-Brenti \cite{b-b} for more details.

The \textit{Hecke algebra} $\mathcal{H}$ associated with $S_n$ is the
algebra over $\Z[q^{\frac{1}{2}},q^{-\frac{1}{2}}]$ generated by
$\{T_i\given 1\leq i\leq n-1\}$ with
the relations 
\begin{enumerate}
\item $(T_{i})^{2} = (q-1) T_{i} + q,$
\item $T_{i}T_{j} = T_{j} T_{i}$ \hspace{.1in} if $|i-j|>1,$ 
\item $T_{i}T_{i+1}T_{i} = T_{i+1} T_{i} T_{i+1}$   for all $1\leq i\leq n-2$.  
\end{enumerate}
This definition is a variation on the relations for the symmetric
group $S_n$ written in terms of its generating set of simple
transpositions and their relations in \eqref{eq:simple.relations}.  If
we take the specialization $q=1$, then the resulting algebra is the
group algebra of $S_n$.  The relations (2) and (3) are called the
\textit{braid relations}.  The braid relations imply that $T_{w} =
T_{i_{1}} T_{i_{2}}\cdots T_{i_{p}}$ is well-defined for any reduced
expression $w=s_{i_{1}} s_{i_{2}} \dots s_{i_{p}}$.  We will use the
notation $T_{\id} =1 \in \mathcal{H}$ for the empty product of
generators.  Note, $T_{\id}$ and $T_{1}$ are different elements of
$\mathcal{H}$.  

An easy observation is that $\{T_{w} \given w \in S_n \}$ is a vector
space basis for $\mathcal{H}$ over
$\Z[q^{\frac{1}{2}},q^{-\frac{1}{2}}]$.  One can also observe that the
$T_{w}$'s are invertible over $\Z[q, q^{-1}]$, which can be seen as
follows.  First check that $\left(T_{i} \right)^{-1} = q^{-1} T_{i} -
\left(1- q^{-1} \right)$ by multiplying by $T_{i}$ and using the
stated relations.  Then, we have $(T_{w})^{-1} = (T_{i_{p}})^{-1}
\cdots (T_{i_{1}})^{-1}$ for any reduced expression $w=s_{i_{1}}
s_{i_{2}} \dots s_{i_{p}}$.  Kazhdan and Lusztig observed that
$\{(T_{w})^{-1} \given w \in S_n \}$ is also a basis for
$\mathcal{H}$.  The change of basis matrix is determined by a family
of polynomials $R_{v,w}(q)$ which appear prominently in the
literature.  The $R$-polynomials can be computed recursively from the
expansion 
\[
(T_{w})^{-1}  = \sum_{v}q^{-\ell(v)}R_{v,w}(q) T_{v}.
\]
The $R$-polynomials contain geometric information that is closely
related to intersections of Schubert cells.  Deodhar proved the
following theorem, see \cite[Thm 1.3]{Deodhar85.BruhatOrderings} and
the subsequent paragraph for the proof.

\begin{Theorem}\label{thm:deodhar.finite.field}
If the flag variety $\Fl(n)$ over a finite field with $q$ elements,
then $R_{v,w}(q)$ counts the number of flags in the \emph{open
Richardson variety} $C_{w}(E_{\bullet}) \cap C_{w_{0}v}(\oppositeE_{\bullet})$.
\end{Theorem}

Kazhdan-Lusztig also considered an involution on the Hecke algebra
given by the $\Z$-linear transformation $i:\mathcal{H}\rightarrow
\mathcal{H}$ sending $T_{w}$ to $(T_{w^{-1}})^{-1}$ and $q$ to
$q^{-1}$, respectively.

\begin{Theorem} [\textbf{Kazhdan-Lusztig} \cite{k-l}]\label{thm:KL.polynomials}
There exists a unique basis $\{C'_{w} \given w \in S_n\}$ for the Hecke algebra $\mathcal{H}$ over $\mathbb{Z}[q^{\frac{1}{2}},q^{-\frac{1}{2}}]$  such that 
\begin{enumerate}
\item[(i)] $i( C'_{w} ) = C'_{w}$.
\item[(ii)] The change of basis matrix from $\{C'_w\}$ to $\{T_w\}$
is upper triangular when the elements of $S_n$ are listed in any total
order respecting Bruhat order, and the expansion coefficients $P_{v,w}(q)$ in
\begin{equation}
C_w' = q^{-\frac{1}{2} \ell(w)} \sum_{v\leq w} P_{v,w}(q) \ T_v 
\end{equation}
have the properties $P_{w,w}=1$ and for all $v<w$, \ $P_{v,w}(q) \in \mathbb{Z}[q]$
with degree at most
\begin{equation}\label{eq:KL.degree.bound}
\frac{\ell(w)-\ell(v)-1}{2}.
\end{equation}
\end{enumerate}
\end{Theorem}

The basis $\{C'_{w} \given w \in S_n\}$ is called the
\textit{Kazhdan-Lusztig basis} for $\mathcal{H}$, and $P_{v,w}(q)$ is
the \textit{Kazhdan-Lusztig polynomial} for $v,w\in S_n$.  This
theorem easily generalizes to all Coxeter groups for the reader
familiar with that topic.  Note, there is a related construction for a
basis $\{C_{w} \given w \in S_n\}$ using $P_{v,w}(q^{-1})$, but in
this context it is more natural to use the polynomial form
$P_{v,w}(q)$, instead of the Laurent polynomial form.

\begin{Example}
  We exhibit some computations with the Kazhdan-Lusztig basis indexed by
  permutations with the aid of \Cref{thm:KL.polynomials}.
  First, it is easy to verify
\begin{align*}
C'_{s_{i}} &= q^{-\frac{1}{2}} (1+ T_{i}) = q^{\frac{1}{2}} (1+ T_{i}^{-1})
\end{align*}
because of the degree bound \eqref{eq:KL.degree.bound}.  Then, for $i
\neq j$, the computation
\begin{align*}
C'_{s_{i}}C'_{s_{j}} &= q^{-1} (1+ T_{i})(1+ T_{j})= q^{-1} (1 + T_{i} + T_{j} + T_{i}T_{j})
\end{align*}
shows that $C'_{s_{i} s_{j}}=C'_{s_{i}}C'_{s_{j}}$ for $i \neq j$
since the coefficients of $C'_{s_{i}}C'_{s_{j}}$ satisfy the
requirements from \Cref{thm:KL.polynomials}.  Also, in the computation
\begin{align*}
C'_{s_{1}}C'_{s_{2}}C'_{s_{1}} &= q^{-\frac{3}{2}} (1+ T_{1})(1+ T_{2}) (1+ T_{1}) \\
                &= q^{-\frac{3}{2}} (1 + 2T_{1} + T_{2} + T_{1}T_{2} + T_{2}T_{1} + 
T_{1}^{2} + T_{1}T_{2}T_{1})\\
	        &= q^{-\frac{3}{2}} (1 + 2T_{1} + T_{2} + T_{1}T_{2} + T_{2}T_{1} \\
               &\hspace{1.4in}     +((q-1) T_{1} + q ) + T_{1}T_{2}T_{1}),
\end{align*}
one notices that $qT_{1} + q$, which comes from $(T_1)^2$, should not
appear for $C'_{s_{1}s_{2}s_{1}}$ because the degree of the polynomial
coefficient of $T_1$ and $T_{\id}$ are too large.  We need a correction
term.  Since $C'_{s_{i}} = q^{-\frac{1}{2}} (1+ T_{i})$, one can check
that $C'_{s_{1}s_{2}s_{1}} = C'_{s_{1}}C'_{s_{2}}C'_{s_{1}} -
C'_{s_{1}}$ by \Cref{thm:KL.polynomials}.
\end{Example}

\begin{Exercise}
 If $i_1,\cdots,i_k\in [n-1]$ are distinct, prove that
\begin{align*}
C'_{s_{i_1}\cdots s_{i_k}} = C'_{s_{i_1}} \cdots C'_{s_{i_k}}.
\end{align*}
More generally, a permutation $w \in S_{n}$ has the \textit{Deodhar
property} if $C'_{w}=C'_{s_{i_{1}}} C'_{s_{i_{2}}}\cdots
C'_{s_{i_{p}}} $ for some reduced expression
$w=s_{i_{1}}s_{i_{2}}\cdots s_{i_{p}}$.  Prove these permutations are
characterized by avoiding the five patterns 321, 56781234, 56718234,
46781235, 46718235.  These 5 patterns are known as the
\textit{321-hexagon avoiding permutations} \cite{b-w,Billey-Jones}.
\end{Exercise}

\begin{Example}
  The Kazhdan-Lusztig polynomials $P_{\id, w}$ for $w\in S_5$
  are completely determined by the following table and the fact that
  $P_{\id, w}=1$ if and only if $w$ is 3412- and 4231-avoiding.  
\begin{center}
\begin{tabular}{|c|c|}\hline
$w$ &  $P_{\id,w}$\\
\hline
$\displaystyle
\begin{array}{lll}
1 4 5 2 3&
1 5 3 4 2&
2 4 5 1 3\\
2 5 3 4 1&
3 4 1 2 5&
3 4 1 5 2\\
3 5 1 2 4&
3 5 1 4 2&
3 5 2 4 1\\
3 5 4 1 2&
4 1 5 2 3&
4 2 3 1 5\\
4 2 3 5 1&
4 2 5 1 3&
4 2 5 3 1\\
4 3 5 1 2&
4 5 1 3 2&
4 5 2 1 3\\
5 1 3 4 2&
5 2 3 1 4&
5 2 4 1 3\\
5 2 4 3 1&
5 3 1 4 2&
5 3 2 4 1\\
5 3 4 2 1&
5 4 2 3 1&
\end{array}$ &  $q+1 $\\
\hline
$ \displaystyle \begin{array}{ll}
3 4 5 1 2&
4 5 1 2 3 \\
4 5 2 3 1&
5 3 4 1 2
\end{array} $&  $2q + 1$\\
\hline
$ \displaystyle  \begin{array}{l}
5 2 3 4 1
\end{array}$ &  $ q^2 + 2q + 1$\\
\hline
$\begin{array}{l}
4 5 3 1 2
\end{array}$ &  $q^2 + 1$\\
\hline
\end{tabular}
\end{center}
\end{Example}

\bigskip

The reader might notice that all coefficients of Kazhdan-Lusztig
polynomials shown so far are non-negative integers.  In their 1979
paper, this property was stated as a conjecture for all
Kazhdan-Lusztig polynomials. In 1980, Kazhdan and Lusztig proved their
own conjecture using \textit{intersection homology} as introduced by
Goresky and MacPherson in 1974, see \cite{g-m} as a good starting
point for that theory.

\begin{Theorem} \cite{K-L2}
If $W$ is a Weyl group or affine Weyl group then 
\[
P_{v,w}(q) = \sum \operatorname{dim} \mathcal{IH}_{v}^{i}(X_{w}) \ q^{i}.
\]
\end{Theorem}

\begin{Corollary}
The coefficients of $P_{v,w}(q)$ are non-negative integers with constant term 1.
\end{Corollary}

In 1990, Soergel gave another proof that the coefficients of
$P_{v,w}(q)$ are non-negative integers using certain modules over the
cohomology ring of the flag variety in the case $W$ is a Weyl group.
However, it was an open problem to prove the same property holds more
generally for all Kazhdan-Lusztig polynomials for all Coxeter groups.
This conjecture was proved by Elias and
Williamson~\cite{Elias-Williamson} in 2014.  Their proof uses Soergel's
bimodules.

As stated above (\Cref{thm:smoothness}), Kazhdan-Lusztig polynomials can be used to
determine smoothness of Schubert varieties (in type A).  There are
several other interesting properties of Kazhdan-Lusztig polynomials
that have emerged since they were defined in 1979.  We cover some of
them here and recommend the Wikipedia page \cite{wiki:kl} 
for a very nice survey.  
\begin{enumerate}
\item 
  In 1981, Beilinson--Bernstein, and independently
  Brylinski--Kashiwara, proved another important conjecture due to
  Kazhdan and Lusztig.  They showed that the multiplicities which
  appear when expressing the formal character of a Verma module in
  terms of the formal character for the corresponding simple highest
  weight module are determined by evaluating Kazhdan-Lusztig
  polynomials at $q=1$ \cite{BeilBern, BryKash}.

\item 

  The coefficients of Kazhdan-Lusztig polynomials are increasing as
  one goes down in Bruhat order, while keeping the second index fixed.
  Specifically, if $u\leq v \leq w$, then $\displaystyle
  \mathrm{coef}_{q^{k}} P_{u,w}(q) \geq \mathrm{coef}_{q^{k}}
  P_{v,w}(q)$.  This monotonicity property was first published in 1988
  by Ron Irving \cite{Irving}, using the socle
  filtration of a Verma module.  In 2001, Braden and MacPherson gave a
  different proof using intersection homology \cite[Cor. 3.7]{B-M}.

\item Every polynomial with constant term 1 and nonnegative integer
coefficients is the Kazhdan-Lusztig polynomial of some pair of
permutations.  This is due to Patrick Polo, published in 1999
\cite{polo}.  He gives an explicit construction of a pair of
permutations for every such polynomial.  This was a surprising result
because from the small data that we can compute, say for $n\leq 9$,
the polynomials seem quite special.  They must get increasingly
complex as $n$ grows.

\item Let $\mu (v,w) $ be the coefficient of
$q^{\frac{\ell(w)-\ell(v)-1}{2}}$ in $P_{v,w}$.  Note, $\mu(v,w)$ can
be 0, so $q^{\frac{\ell(w)-\ell(v)-1}{2}}$ is not necessarily the
leading term of $P_{v,w}(q)$.  For $v,w \in S_{9}$,\ $\mu (v,w) \in
\{0,1 \}$.  MacLarnen and Warrington found an example in $S_{10}$
where $\mu(v,w) = 5$ \cite{MW02}.  Prior to their publication in 2003,
this was referred to as the ``0-1 Conjecture for Kazhdan-Lusztig
polynomials.''  This again demonstrates the increasing complexity as
$n$ grows.  The reader might be wondering how anyone could have
believed the 0-1 Conjecture after seeing Polo's theorem in (3).
However, Polo's theorem does not contradict the 0-1 Conjecture because
in his construction the length difference between $w$ and $v$ is large
enough that the leading term in $P_{v,w}(q)$ is typically not the
$\mu$-coefficient.

In special cases, Jones showed the 0-1 conjecture holds and classified
which pairs lead to $\mu(v,w)=1$ in those cases \cite{Jones.2009}.
Furthermore, Lusztig considered the statistic on permutations $w \in
S_{n}$ given by $a(w)=\frac{1}{2}\sum \lambda_{i}(\lambda_{i}-1)$
where $\lambda =(\lambda_{1},\dots , \lambda_{p})$ is the transposed
shape of the RSK tableaux for $w$, see \cite{lusztig.cells} and
\cite[Ch.6, Ex. 10]{b-b}.  In
\cite{xi2004leadingcoefficientcertainkazhdanlusztig}, Xi proved that
$\mu(v,w) \in \{0,1 \}$ whenever $a(v)<a(w)$.

\item \label{item:abstract.interval} There exists a formula for
$P_{v,w}(q)$ which only depends on the abstract interval
$[\mathrm{id},w]$ in Bruhat order.  See the work of du Cloux (2003)
\cite{ducloux.2004}, Brenti (2004) \cite{brenti04} and
Brenti-Caselli-Marietti (2006) \cite{BCM.2006}, and Open Problem~\ref{prob:comb-invariance} below.
\end{enumerate}

\bigskip

The next pattern property connects the 4231 and 3412 patterns to the
determination of the singular locus of a Schubert variety.  Recall
from Fact (2) from the start of \Cref{sub:Singularloci} that the
singular locus of a Schubert variety $X_w$ is a union of certain
Schubert varieties $X_v$ with $v<w$.  Because Bruhat order
characterizes when $X_{u} \subset X_{v}$, one just needs to identify
the maximal permutations $v <w$ such that $v_{\bullet}$ is a singular
point in $X_w$.  The theorem below shows that all such \textit{maximal
singular permutations} are obtained from $w$ by right multiplication
by a cycle determined by a subset of the 4231, 3412, and 45312
patterns in $w$.  The permutation 45312 contains a 3412 pattern, but
for the sake of identifying the maximal singular permutations for $w$
it is useful to consider it separately.

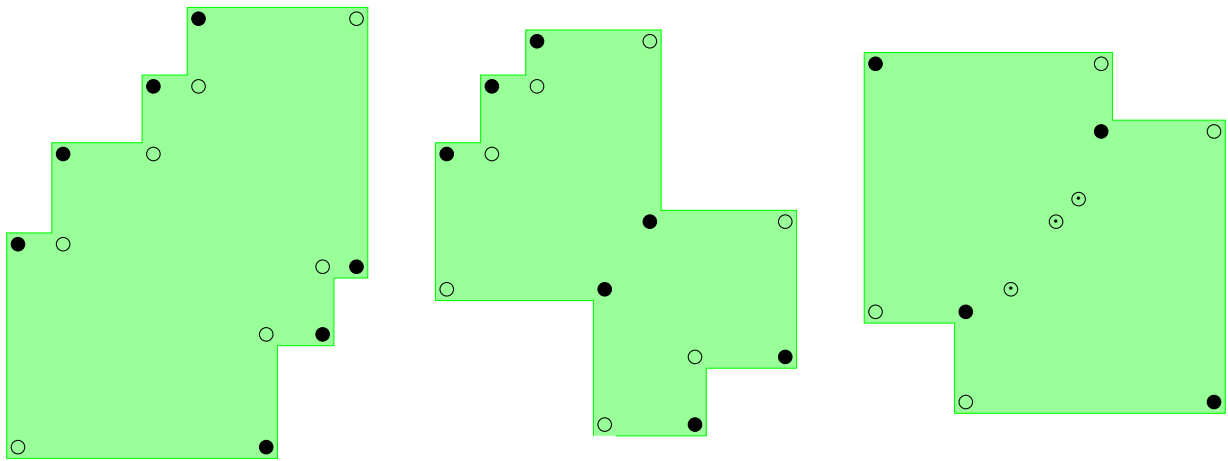
\begin{figure}[h] \centering
\begin{tikzpicture}[scale=0.3] \def\r{0.3}; \draw[green, fill=green,
fill
opacity=0.4](0,0)--(12,0)--(12,5)--(14.5,5)--(14.5,8)--(16,8)--(16,20)--(8,20)--(8,17)--(6,17)--(6,14)--(2,14)--(2,10)--(0,10)--(0,0);
\draw (0.5,0.5) circle (\r); \draw (11.5,5.5) circle (\r); \draw
(14,8.5) circle (\r); \draw (15.5,19.5) circle (\r); \draw (8.5,16.5)
circle (\r); \draw (6.5,13.5) circle (\r); \draw (2.5,9.5) circle
(\r); \filldraw (11.5,0.5) circle (\r); \filldraw (14,5.5) circle
(\r); \filldraw (15.5,8.5) circle (\r); \filldraw (8.5,19.5) circle
(\r); \filldraw (6.5,16.5) circle (\r); \filldraw (2.5,13.5) circle
(\r); \filldraw (0.5,9.5) circle (\r);

\def\b{19}; \draw[green, fill=green, fill
opacity=0.4](8+\b,1)--(12+\b,1)--(12+\b,4)--(16+\b,4)--(16+\b,11)--(10+\b,11)--(10+\b,19)--(4+\b,19)--(4+\b,17)--(2+\b,17)--(2+\b,14)--(0+\b,14)--(0+\b,7)--(7+\b,7)--(7+\b,1);
\draw (7.5+\b,1.5) circle (\r); \draw (11.5+\b,4.5) circle (\r); \draw
(15.5+\b,10.5) circle (\r); \draw (9.5+\b,18.5) circle (\r); \draw
(4.5+\b,16.5) circle (\r); \draw (2.5+\b,13.5) circle (\r); \draw
(0.5+\b,7.5) circle (\r); \filldraw (11.5+\b,1.5) circle (\r);
\filldraw (15.5+\b,4.5) circle (\r); \filldraw (9.5+\b,10.5) circle
(\r); \filldraw (4.5+\b,18.5) circle (\r); \filldraw (2.5+\b,16.5)
circle (\r); \filldraw (0.5+\b,13.5) circle (\r); \filldraw
(7.5+\b,7.5) circle (\r);

\def\c{38}; \draw[green, fill=green, fill
opacity=0.4](4+\c,2)--(16+\c,2)--(16+\c,15)--(11+\c,15)--(11+\c,18)--(0+\c,18)--(0+\c,6)--(4+\c,6)--(4+\c,2);
\draw (4.5+\c,2.5) circle (\r); \draw (15.5+\c,14.5) circle (\r);
\draw (10.5+\c,17.5) circle (\r); \draw (0.5+\c,6.5) circle (\r);
\filldraw (15.5+\c,2.5) circle (\r); \filldraw (10.5+\c,14.5) circle
(\r); \filldraw (0.5+\c,17.5) circle (\r); \filldraw (4.5+\c,6.5)
circle (\r); \draw (6.5+\c,7.5) circle (\r); \node at (6.5+\c,7.5)
{$\cdot$}; \draw (8.5+\c,10.5) circle (\r); \node at (8.5+\c,10.5)
{$\cdot$}; \draw (9.5+\c,11.5) circle (\r); \node at (9.5+\c,11.5)
{$\cdot$};
\end{tikzpicture}
\caption{Patterns for the singular locus of $X_w$ in the 4231, 3412,
and 45312 cases respectively.  Here $\circ$'s denote 1's in $w$ only,
$\bullet$'s denote 1's in $v$ only, and the circle around a dot
denotes a 1 in both $v,w$.  Multiplication of $w$ by the cycle
$(\alpha_{1},\dots , \alpha_{m}, \beta_{1},\dots , \beta_{m})$ rotates
the identified 1's in $w$ clockwise around the boundary of the shaded
region to obtain $v$.  }\label{f:sing.locus}
\end{figure}

\begin{Theorem}\cite{BW-sing,cortez,klr,manivel}\label{thm:sing.locus}
The Schubert variety $X_{v}$ is an irreducible component of the
singular locus of $X_{w}$ if and only if the following conditions
hold.
\begin{enumerate}
\item $v=w\cdot (\alpha_{1},\dots , \alpha_{m},
\beta_{k},\dots , \beta_{1})$ for a pair of disjoint sequences
\begin{align*}
    1&\leq \alpha_1 < \cdots < \alpha_m \leq n\text{, with } w(\alpha_1) > \cdots > w(\alpha_m),\text{ and }\\
    1&\leq \beta_1 < \cdots < \beta_k \leq n\text{, with } w(\beta_1)
> \cdots > w(\beta_k).
\end{align*}
  
\item The permutation matrices for $v$ and $w$ differ in one of three
ways shown in \Cref{f:sing.locus} depending on 4231, 3412, or 45312
patterns.
\item The interiors of the shaded regions determined by $w$ and
$(\alpha_{1},\dots , \alpha_{m}, \beta_{k},\dots , \beta_{1})$ shown
in \Cref{f:sing.locus} do not contain any other 1's in the permutation
matrix of $w$, except in the third case where the 1's in the central
region are allowed provided they form a decreasing sequence.
\end{enumerate}
\end{Theorem}

\Cref{thm:sing.locus} could be called the \textit{Billey-Warrington,
Kassel-Lascoux-Reutenauer, Manivel, Cortez Theorem}.  The result was
found around 2000 by 7 authors in 4 papers, plus Gasharov proved one
direction of the conjecture \cite{gasharov.2001} around the same
time. It must have been ripe for discovery.  It refined and proved a
conjecture due to Lakshmibai and Sandhya \cite{Lak-San}.  For the sake
of history, we note that the authors of \cite{BW-sing} were the first
to report this result to Lakshmibai.

\begin{Example}\label{ex:3.types.singular}
The permutation $w=87432651$ has a $4231$ pattern given by the
subsequence $7251$, so $X_{w}$ is singular. To determine if $7251$
gives rise to a maximal singular permutation for $w$, we use the
$4231$ pattern shown in \Cref{f:sing.locus} on the left.  One must
extend the subsequence $7251$ to include as many 1's as possible in
$M_{w}$ along the northwest and southeast boundaries of the shaded
region.  In this case, the subsequence $7251$ in $w$ extends to
include $7432651$ and there are indeed no additional 1's in the
interior of the corresponding region.  To find the $\alpha$ and $\beta$
sequence in this case, divide the extended subsequence at its unique
descent; then the first decreasing sequence is indexed by the
positions in $\alpha$ and the second decreasing sequence is in
positions indexed by $\beta$.  Here the extended subsequence is
$7432|651$, which appears in positions 2 to 8 in $w$.  The $\alpha$
sequence is $2<3<4<5$, and the $\beta$ sequence is $6<7<8$.  The cycle
$(2 \mapsto 3 \mapsto 4 \mapsto 5 \mapsto 8 \mapsto 7 \mapsto
6)=134582567$ in one-line notation.  By \Cref{thm:sing.locus}, $v=w
\cdot 134582567= 84321765$ indexes a maximal singular permutation for
$w$.  Note, all 4231 patterns in $7432651$ give rise to the same
extended subsequence in this case.  All of the $4231$ patterns in $w$
starting with 8 will give rise to a shaded region containing the 1
in position $(7,2)$ in $M_{w}$, so none of them will give rise to a
maximal singular permutation, and $w$ contains no $3412$ patterns.
Thus, $\mathrm{Sing}(X_{87432651})=X_{84321765}$ by
\Cref{thm:sing.locus}.
\end{Example}

\begin{Corollary}
  The codimension of the singular locus of a Schubert variety $X_w$ is
at least 3 for any $w \in S_n$.  
\end{Corollary}

\begin{proof}
From \Cref{thm:sing.locus}, we observe that each maximal singular
permutation for $w$ has at least 3 fewer inversions than $w$.  Hence,
the statement follows from \eqref{eq:dim.Schubert.cell}.
\end{proof}

The corollary is in fact true for all simply laced types, meaning
types $A,D,E$ in the classification of simple Lie groups and their
Weyl groups (cf. \Cref{subsec:type-B}). However, it is not true in type $B_n$ where the Weyl
group is the signed permutations.  The codimension of the singular
locus of a Schubert variety in that case can be 2.

Inspired by the Lakshmibai-Sandhya \Cref{Lak-San} and the construction
of the singular locus of a Schubert variety in \Cref{thm:sing.locus},
Woo and Yong \cite{WooYong} defined the notion of interval pattern
avoidance.  Given permutations $u <v \in S_m$ and $x<y \in S_n$ for
$m<n$, say $[u,v]$ \textit{interval pattern embeds} into $[x,y]$
provided
\begin{enumerate}
\item there exist indices $1\leq i_1< i_2 < \ldots < i_m \leq n$ such
that $\fl(x_{i_1}, \ldots, x_{i_m}) = u$ and $\fl(y_{i_1}, \ldots,
y_{i_m}) = v$ respectively,
\item the permutations $x,y$ agree in all positions other than $1\leq
i_1< i_2 < \ldots < i_m \leq n$, and
\item the Bruhat intervals $[u,v]$ and $[x,y]$ are isomorphic as
posets.
\end{enumerate}
In fact, if $x,y$ agree everywhere outside of the indices $1\leq i_1<
i_2 < \ldots < i_m \leq n$ and $u = \fl(x_{i_1}, \ldots, x_{i_m}),\ v=
\fl(y_{i_1}, \ldots, y_{i_m})$ then $[u,v]$ interval embeds in $[x,y]$
if and only if $\ell(v)-\ell(u) = \ell(y)-\ell(x)$ \cite[Lemma
2.1]{WooYong}.  Furthermore, if $[u,v]$ interval pattern embeds into
$[x,y]$, then for all $w \in S_n$ such that $x<w<y$,\ $w$ also agrees
everywhere with $y$ outside of the sequence and $[\fl(w_{i_1}, \ldots,
w_{i_m}), v]$ also interval embeds in $[x,y]$ \cite[Lemma
2.4]{WooYong}.  Combining these results with \Cref{thm:sing.locus},
one obtains the following result.  The verification is a useful
exercise for the reader.

\begin{Corollary}
The maximal singular locus of a Schubert variety is determined by
interval pattern conditions.
\end{Corollary}

Woo and Yong gave the following far-reaching generalization of this
observation.  This demonstrates the interplay between the
combinatorial data of permutations and the geometrical properties of
Schubert varieties.

\begin{Theorem}\cite[Thm. 4.2]{WooYong}\label{thm:wy-interval-embed}
Suppose $[u,v] \subset S_m$ interval pattern embeds into $[x,y] \in
S_n$. Then there exist affine neighborhoods of $X_{v} \subset \Fl(m)$
and $X_{y} \subset \flags$ respectively of $u_{\bullet}$ and
$x_{\bullet}$ which are isomorphic up to Cartesian products with
affine space.
\end{Theorem}

\begin{Corollary}\cite[Cor. 6.3]{WooYong}\label{thm:woo.yong} Suppose
$[u,v] \subset S_m$ interval pattern embeds into $[x,y] \in S_n$. Then
the Kazhdan-Lusztig polynomials $P_{u,v}(q) $ and $P_{x,y}(q) $ are
equal.
\end{Corollary}

For more details on interval pattern avoidance, see the chapter titled
``Schubert Geometry and Combinatorics'' by Woo and Yong in this book
\cite{WooYong2023schubertgeometrycombinatorics}.  They use interval
pattern avoidance to characterize many different properties of
Schubert varieties where simple pattern avoidance was insufficient
including factoriality, multiplicity, Gorensteinness, and
Cohen-Macaulay type.  Even more variations on this theme appear
in \cite{ulfarsson.2011}.

\begin{Remark}
A reader more familiar with symplectic geometry might notice that a
similar result to Theorem~\ref{thm:woo.yong} could be obtained using
the tools for intersection homology as explained by Braden-MacPherson
\cite{B-M}.
\end{Remark}

\bigskip

We conclude this section by uniting the theory of double Schubert
polynomials $\fS(X,Y)$ as defined in
\eqref{eq:divided.difference.for.doubleSchubs} with the singular locus
of a Schubert variety in the flag variety $\flags$.  To our knowledge,
this has not been stated explicitly in the literature previously. This
theorem is derived from Kumar's test for smoothness in the nil-Hecke
algebra \cite{kumar} and Billey's work on Kostant polynomials
\cite{B4}.  Kumar's test is the most general test for smoothness of
Schubert varieties that extends uniformly to all Kac-Moody groups, but
we are only using the type $A$ special case below.  Our presentation of the proof builds on \cite[Ch. 7]{BLak}.

\begin{Theorem}\label{thm:BK-smooth-test} Let $v\leq w$ in $S_{n}$.
Then $X_{w}$ is smooth at the flag $v_{\bullet}$ if
and only if 
%the double Schubert polynomial $\fS_{v}(X;Y)$ evaluated at permuted variables $y$-variables in place of $x$-variables is exactly the product of terms
\[
 \fS_{ww_0}( y_{vw_0(1)}, y_{vw_0(2)},\dots ,
y_{vw_0(n)};y_{1},y_{2},\dots , y_{n}) = \prod_{\substack{1\leq i<j\leq n\\
t_{ij}v\not \leq w}} (y_{j}-y_{i}).  
\]
\end{Theorem}

\begin{proof}
Fix a reduced word $\mathbf{b}=(b_{1}, b_{2},\ldots, b_{p}) \in R(w)$,
define a sequence of degree one polynomials in $\mathbb{Z}[y_{1},\dots
, y_{n}]$ by
\[
r_{\mathbf{b}}(j) = s_{b_{1}}s_{b_{2}}\cdots s_{b_{j}}(y_{b_{j}+1}-y_{b_{j}})
\]
for all $1\leq j\leq p$, where the action is defined by
$w(y_{k})=y_{w(k)}$.  The ordered list of polynomials
$r_{\mathbf{b}}(1),\dots , r_{\mathbf{b}}(p)$ is closely related to
the reflection order defined by $\mathbf{b}$ in \Cref{sub:Balanced}.
For all $v,w \in S_{n}$, let
\begin{equation}\label{eq:ajs.billey.formula}
d_{w,v}(Y) = \sum_{b_{i_{1}} b_{i_{2}} \dots b_{i_{k}} \in R(v)}
r_{\mathbf{b}}(i_{1})r_{\mathbf{b}}(i_{2})\cdots r_{\mathbf{b}}(i_{k})
\in \mathbb{Z}[y_{1},\dots , y_{n}]
\end{equation}
where the sum is over all subsequences  $1\leq i_{1}<i_{2}<\dots <
i_{k}\leq p$ such that  $b_{i_{1}} b_{i_{2}} \dots
b_{i_{k}}$ is a reduced word for $v$.  By definition, $d_{w,v}(Y)=0$
if $v \not \leq w$.  As suggested in \cite[Section 8, Remark 1]{B4},
the polynomial $d_{w,v}(Y)$ is exactly the double Schubert polynomial
$\fS_{w}(X,Y)$ with the $X$ variables evaluated at the permuted $Y$
variables given by $v.Y=(y_{v_{1}}, y_{v_{2}},\dots , y_{v_{n}})$.
Therefore, 
\begin{equation}\label{eq:unnormalized.kostant.poly.to.double.schub}
d_{w,v}(Y) = \fS_{w}(v.Y,Y) \ \text{ for all } v,w \in S_{n}. 
\end{equation}
In particular, the polynomial $d_{w,v}(Y)$ is well-defined and does not depend on
the choice of $\mathbf{b} \in R(w)$.

By \cite[Cor. 7.2.8]{BLak}, $X_{vw_{0}}$ is smooth at the flag
corresponding with $ww_0$ if and only if
\[
d_{w,v} = \prod_{\substack{1\leq i<j\leq n\\
v\not \leq t_{ij} w}} (y_{j}-y_{i}).  
\]
Note, $v \leq t_{ij} w$ if and only if $t_{ij} ww_{0} \leq vw_{0}$
since multiplication on the right by $w_{0}$ is order-reversing by
\Cref{ex:order-reversing}. Therefore, the result now follows by
swapping indices $w\mapsto vw_{0}$ and $v \mapsto ww_{0}$.
\end{proof}

\begin{Example}\label{ex:double.schub.smooth.test}
Consider $w=3241$. Since $w$ avoids both $3412$ and $4231$, $X_w$ is smooth, and in particular, it is smooth at $v=w=3241$. Now the double Schubert polynomial of $ww_0=1423$ is
\[
\fS_{1423}(X,Y)=(x_{2}-y_{1})(x_{2}-y_{2})+(x_{2}-y_{1})(x_{1}-y_{3})+(x_{1}-y_{2})(x_{1}-y_{3}).
\]
Since $vw_{0}=1423$, replace $x_{1}$ with $y_{1}$ and $x_{2}$ with
$y_{4}$ (only since there are no $x_{3}$'s or $x_{4}$'s) to obtain

\begin{align*}
\fS_{ww_{0}}(vw_0Y,Y) &= (y_{4}-y_{1})(y_{4}-y_{2})+(y_{4}-y_{1})(y_{1}-y_{3})+(y_{1}-y_{2})(y_{1}-y_{3})\\
&= (y_4 - y_2) (y_4 - y_3).
\end{align*}

The factored form of $\fS_{ww_{0}}(vw_0Y,Y)$ has roots indexed by
$\{1\leq i<j\leq 4 \given t_{ij}v \not \leq w \}$ since $t_{24}3241 =
3421 >3241,\ t_{34}3241 = 4231 >3241$, and all other transpositions
multiplied on the right of $3241$ reduce the number of
inversions. Therefore, $\fS_{ww_{0}}(vw_0Y,Y)$ has the factored form
given in \Cref{thm:BK-smooth-test}.

\bigskip

On the other hand, if $w=4231$, then $w w_{0}=1324$.  The double
Schubert polynomial of $1324$ is $\fS_{1324}(X,Y)=(x_{2}-y_{1})+(x_{1}-y_{2})$.
Replacing $x_{2}$ with $y_{3}$ and $x_{1}$ by $y_{4}$, we have
$\fS_{1324}(4321.Y,Y)=(y_{3}-y_{1})+(y_{4}-y_{2})$, which cannot be
written a as product of binomials.  Hence, $X_{4231}$ is singular at
$1234$ as expected.
\end{Example}

There are several interesting but difficult open problems in
Kazhdan-Lusztig theory and its connections with singular loci of
Schubert varieties.  There are many partial answers to these questions
in the literature, but we don't know of a complete source at this
time.  Perhaps there is a need for someone to start a wiki page.  See
Vazirani's webpage for a snapshot from 2005 \cite{Vazirani}.  The
following exercises lead into some conjectures.

\begin{Exercise}
Let $v\leq w$ in $S_{n}$.  What is the minimum value of
$\ell(w)-\ell(v)$ such that $P_{v,w}(1)>1$, $2$, and $3$?
\end{Exercise}

\begin{Exercise}\cite{Woo-Billey-Weed}
Characterize the permutations $w$ such that $P_{w}(1)=2$.  
\end{Exercise}

The following conjecture has been tested through $S_{8}$. Let $\ms(w)$
denote the set of permutations indexing the irreducible components of 
the singular locus of $X_{w}$.  

\begin{Conjecture}\cite[Conj. A.5]{Woo-Billey-Weed}  Let $w \in S_{n}$. 
\begin{enumerate}
\item If $P_{\mathrm{id}, w}(1) \leq 3$ then $|\ms(w)| \leq 3$.
\item If $P_{\mathrm{id}, w}(1) =3$ and  $|\ms(w)| =1$ then
$P_{\mathrm{id},w}=1+q^{a}+q^{b}$ for some  $0<a<b$. 
\item If $P_{\mathrm{id}, w}(1) =3$ and  $|\ms(w)| =2$ then
$P_{\mathrm{id},w}=1+q^{a}+q^{b}$ for some  $0<a<b$. 
\item If $P_{\mathrm{id}, w}(1) =3$ and  $|\ms(w)| =3$ then
$P_{\mathrm{id},w}=1+2q^{a}$ for some  $0<a$.   
\end{enumerate}
\end{Conjecture}

\begin{Problem}[\textbf{Combinatorial Invariance Conjecture} \label{prob:comb-invariance} due to Lusztig and Dyer
independently]\label{q:lusztig.interval} $P_{v,w}(q) = P_{x,y}(q)$ whenever $[v,w]$ and
$[x,y]$ are isomorphic as posets.
\end{Problem}

This conjecture has inspired a great deal of research.  For example,
recall the special case proved by Woo and Yong using interval pattern
avoidance stated above in \Cref{thm:woo.yong}.  Also,
\Cref{item:abstract.interval} on page \pageref{item:abstract.interval}
provides a weaker version of this conjecture.  Note that knowing the
abstract poset $[v,w]$ is not equivalent to knowing the abstract poset
$[\mathrm{id},w]$ since $[\mathrm{id},w]$ contains much more information in general.
For further background and special cases, see for example
\cite{Dyer.thesis,BCM.2006,Patimo.2021,BBDVW.2022}.

\begin{Problem}\label{q:kl.poly}
What is the computational complexity of computing the coefficients of
$P_{v,w}(q)$?  Is there a combinatorial interpretation for the
coefficients that would give a positivity test in polynomial time?
\end{Problem}

For example, Deodhar proposed that $P_{v,w}(q)$ can be expressed as a
generating function over subexpressions of a reduced expression for
$w$ using the ``defect statistic'', but the indexing set for such
a sum has been elusive in the general case \cite{Deodhar90}.  It is
known in the case $w$ is a 321-hexagon avoiding permutation
\cite{b-w}.  See also Lascoux's algorithm for the case $w$ is
$3412$-avoiding (co-vexillary) \cite{Lascoux95} or \cite[\S
6.3.28]{BLak}.

\begin{Problem}\label{q:kl.poly.deg}
  What is the degree of $P_{v,w}(q)$ in the case $\mu(v,w)=0$?  Here
$\mu(v,w)$ is the coefficient of $q^{\frac{\ell(w)-\ell(v)-1}{2}}$ in
$P_{v,w}$.
\end{Problem}

\begin{Problem}\label{q:interval.lattice}
Let $w\in S_{n}$ be a permutation such that $X_{w}$ is singular.
Assume the permutations $v^{(1)},\dots , v^{(s)}$ index the
irreducible components of $\mathrm{Sing}(X_{w})$, as 
determined by \Cref{thm:sing.locus}.  How does the lattice on the
intersections $X_{v^{(i_{1})}}\cap \cdots X_{v^{(i_{k})}}$ ordered by
containment relate to the computation of the Kazhdan-Lusztig
polynomials $P_{u,w}(q)$ for all $u\leq w$?
\end{Problem}

\subsection{Isomorphism Classes of Schubert Varieties}\label{sub:isomorphisms}

How can we determine if two Schubert varieties are isomorphic either
in the same flag variety or from different flag varieties $\Fl(m)$ and
$\Fl(n)$?  Having the same dimension is certainly a necessary
condition, but it is not sufficient.  For example, $X_{4231}$ is
singular and $X_{3421}$ is smooth, so they cannot be isomorphic even
though both have dimension 5.  For every simple transposition $s_{i}$,
we know $X_{s_{i}}$ is isomorphic to $\mathbb{P}^{1}$ so the
1-dimensional Schubert varieties in all flag varieties are all
isomorphic. What about the 2-dimensional Schubert varieties?
Hopefully, you solved the 2-dimensional isomorphism problem in
\Cref{ex:rank2.isom.types}.  For example, one can see via rank
conditions that $X_{2314}$ is isomorphic to $\mathbb{P}^{2}$, while
$X_{2143}$ is isomorphic to $\mathbb{P}^{1} \times \mathbb{P}^{1}$. In
fact, every 2-dimensional Schubert variety in $\Fl(4)$ other than
$X_{2143}$ is isomorphic to $\mathbb{P}^{2}$.  Observe that the Bruhat
intervals for $[\mathrm{id},2143]$ and $[\mathrm{id},2314]$ are
``diamonds'' so their intervals are isomorphic as posets.  Hence,
isomorphism of Bruhat intervals also does not suffice to distinguish
isomorphism classes.  What is different about $2143$ and $2314$?  For
one thing, the group structure on the supporting generators is 
different: $s_{1}s_{3}=s_{3}s_{1}$ and $s_{1}s_{2} \neq s_{2}s_{1}$.
Does that matter?

Can we use the cohomology rings for $X_{v}$ and $X_{w}$ to distinguish
the isomorphism classes? Each Schubert variety is the union of cells,
whose closures give rise to Schubert classes in cohomology. They
aren't all smooth however, so we cannot rely on the Chow ring
machinery from \Cref{sub:Monk} where smoothness is a hypothesis.  Recall,
$X_{w_{0}}=\Fl(n)$ and $H^{*}(\Fl(n),\mathbb{Z})$ is isomorphic to the
coinvariant algebra $R_{n}=\mathbb{Z}[x_{1},x_{2}, \dots ,
x_{n}]/I_{n}^{+}$ via the Borel presentation in \Cref{thm;borel}.
Thus, a presentation of $H^{*}(X_{w},\mathbb{Z})$ should be of the
form $R_{n}/I_{w}$.  The following theorem makes this precise.  See
\cite[Cor. 4.4]{carrell92} or \cite{ALP} for the proof.

\begin{Theorem}\label{thm:cohomlogy.Xw}
For $w \in S_{n}$, 
\[
H^{*}(X_{w},\mathbb{Z}) \simeq R_{n}/I_{w} \simeq \mathbb{Z}[x_{1},x_{2}, \dots]/J_{w}
\]
where $I_{w}$ is the ideal generated by the Schubert polynomials
$\fS_{y}$ for $y \not \leq w$ in Bruhat order on $S_{n}$, and $J_{w}$
is the ideal generated by the Schubert polynomials $\fS_{y}$ for $y
\not \leq w$ in Bruhat order on $S_{\infty}$.  The polynomials
$\{\fS_{v} \given v \leq w\}$ form the \textit{Schubert basis} for
$H^{*}(X_{w},\mathbb{Z})$ in this presentation, and multiplication is
determined by $\fS_{u} \fS_{v}$ modulo $J_{w}$.
\end{Theorem}

Akyildiz–Lascoux–Pragacz \cite[Thm 2.2]{ALP} improved on
\Cref{thm:cohomlogy.Xw} by showing both $I_{w}$ and $J_{w}$ can be
generated by those $\fS_{u}$ such that $u$ is Grassmannian and $u\not
\leq w$.  Reiner-Woo-Yong \cite{RWY} further reduced the number of
Schubert polynomials needed to generate the ideals, and they
conjectured a minimal set of generators.  This was proved by
St. Dizier and Yong in \cite[Thm. 1.2]{dizier.yong2023presenting}.

To determine the structure constants in the Schubert basis for
$H^{*}(X_{w},\mathbb{Z})$ by \Cref{thm:cohomlogy.Xw}, consider
\[
\fS_{u} \fS_{v} = \sum_{w' \in S_{n}} c_{u,v}^{w'} \fS_{w'}
\]
as polynomials.  Then, the expansion of $\fS_{u} \fS_{v}$ modulo
$J_{w}$ is given by
\[
\fS_{u} \fS_{v} = \sum_{w'\leq w} c_{u,v}^{w'} \fS_{w'}
\]
as elements in $H^{*}(X_{w},\mathbb{Z})$.  Hence, the $c_{u,v}^{w'}$
are either the same as in $H^{*}(\Fl (n),\mathbb{Z})$ or 0 depending
on whether $w'\leq w$ or not.  

\begin{Example}\label{ex:2314.2143}
If $w=2314 \in S_{4}$, then $w=s_{1}s_{2}$ so $[\mathrm{id},w]$
contains the 4 elements $\{\mathrm{id},s_{1}, s_{2},w \}$.  By Monk's
formula for Schubert polynomials (\Cref{thm:Schubert.Monk}), we can
compute all of the nontrivial products of Schubert classes in
$H^{*}(X_{w})$ to be 
\[
\begin{array}{rccclr}
(\fS_{2134})^{2} & = & \fS_{3124} & \equiv & 0 & \text{ modulo } J_{w}\\
(\fS_{1324})^{2} & = & \fS_{1423} + \fS_{2314}  & \equiv  & \fS_{w}  &  \text{ modulo } J_{w}\\
\fS_{1324}\fS_{2134}& = &\fS_{3124} + \fS_{2314} & \equiv & \fS_{w}  & \text{ modulo } J_{w}.
\end{array}
\]
On the other hand, if $y=2143 \in S_{5}$, then $[\mathrm{id},y]$
contains the 4 elements $\{\mathrm{id},s_{1}, s_{3},y \}$, and the
nontrivial products of Schubert classes in $H^{*}(X_{y})$ are
\[
\begin{array}{rccclr}
(\fS_{2134})^{2} & = & \fS_{3124} & \equiv  & 0  & \text{ modulo } J_{w}\\
\fS_{1243}^{2} & =  & \fS_{1342} + \fS_{12534}  &  \equiv  &  0  & \text{ modulo } J_{w}
\\
\fS_{1243}\fS_{2134} & =  & \fS_{2143}  &  =   & \fS_{y}. & 
\end{array}
\]
Thus, the Schubert classes have different structure constants in the
rings $H^{*}(X_{w}, \mathbb{Z})$ and $H^{*}(X_{y}, \mathbb{Z})$.
Indeed, we know $X_{w}$ and $X_{y}$ are not isomorphic varieties.
There is a graded ring isomorphism of $H^{*}(X_{w}, \mathbb{Q})$ and
$H^{*}(X_{y}, \mathbb{Q})$ given by mapping $\fS_{2134}$ to itself and
$\fS_{1324}$ to $\frac{1}{2}(\fS_{2134}+\fS_{1243})$, so the
coefficient ring does matter for differentiating isomorphism classes.
\end{Example}

Due to a remarkable theorem of Richmond-Slofstra
\cite{RichmondSlofstra2022isomorphismproblemschubertvarieties}
building on the intuition from \Cref{ex:2314.2143}, there is a simple
criterion to detect Schubert variety isomorphism classes in complete
flag varieties in terms of the group structure generated by the simple
transpositions below a permutation.  Their characterization can be
stated much more generally in terms of the Cartan matrix for all
Kac-Moody groups $G$ and their generalized flag varieties $G/B$.  We
will summarize their results in the case of Schubert varieties indexed
by permutations.

Let $w \in S_{\infty}$.  The \emph{support} of $w$ is the subset of
simple transpositions 
\[
S(w) =\{i \given s_{i} \leq w \}
\]
under Bruhat order.  Define $G(w)$ to be the simple graph on vertex
set $S(w)$ with an edge $(i, j)$ if $s_{i}s_{j} \neq
s_{j}s_{i}$. Because $w$ is in some $S_{n}$, we know $G(w)$ is a union
of path graphs, and its edges encode the noncommuting pairs of
generators in $S(w)$.  Let $A(w)$ be the \emph{adjacency matrix} of
$G(w)$, where $A(w)_{i,j}=1$ if $s_{i}s_{j} \neq s_{j}s_{i}$, and
$A(w)_{i,j}=0$ otherwise.  Technically, the \textit{Cartan matrix}
associated to the group generated by $\{s_{i} \given i \in S(w) \}$ is
$2I-A(w)$, but we won't need this here.

\begin{Definition}\cite[Def. 1.2]{RichmondSlofstra2022isomorphismproblemschubertvarieties}
\label{def:cartan.equivalent} Let $v, w \in S_{\infty}$.  The
permutation $v$ is \textit{Cartan equivalent} to $w$ if there exists a
bijection on the supports $\sigma : S(v) \longrightarrow S(w)$ such
that
\begin{enumerate}
\item there exists a reduced word $a_{1}a_{2}\cdots a_{p} \in R(v)$
with $\sigma (a_{1})\sigma (a_{2}) \cdots \sigma (a_{p}) \in R(w)
$, and 
\item if $s_{i}s_{j} \leq v$ for any $i \neq j \in S(v)$, then
$A(v)_{i,j}=A(w)_{\sigma (i),\sigma (j)}$.  
\end{enumerate}
\end{Definition}

Note, the definition of Cartan equivalence is not symmetric in $v$ and
$w$.  However, it is an instructive exercise to show Cartan
equivalence is indeed an equivalence relation.

\begin{Exercise}
Prove that if $v$ is Cartan equivalent to $w$ via the bijection
$\sigma$, then $w$ is Cartan equivalent to $v$ via the bijection
$\sigma^{-1}$.  Furthermore, prove $\sigma$ maps every reduced word
for $v$ to a reduced word for $w$.  Thus, $\sigma$ extends to an poset
isomorphism from $[\mathrm{id},v]$ to $[\mathrm{id},w]$.  
\end{Exercise}

The classification of all isomorphism types of Schubert varieties in
complete flag varieties is determined by Cartan equivalence.  As
mentioned above, this characterization extends to Schubert varieties
beyond $\GL_{n}/B$ as well.  We refer the reader to
\cite[Thm. 1.3]{RichmondSlofstra2022isomorphismproblemschubertvarieties}
for complete details.  See also \cite{Richmond.etal.2024} for the
isomorphism types in partial flag varieties of cominuscule type.

\begin{Theorem}[\textbf{Richmond-Slofstra Isomorphism Theorem} \cite{RichmondSlofstra2022isomorphismproblemschubertvarieties}]\label{thm:Richmond-Slofstra}
Let $v \in S_{m}$ and $w \in S_{n}$.  Then, the following are
equivalent.
\begin{enumerate}
\item The permutations $v$ and $w$ are Cartan equivalent. 
\item The Schubert varieties $X_{v} \subset \Fl(m)$ and $X_{w} \subset
\Fl(n)$ are isomorphic as projective varieties.  
\item There is a graded ring isomorphism $\phi : H^{*}(X_{v},
\mathbb{Z}) \longrightarrow H^{*}(X_{w}, \mathbb{Z})$ mapping the
Schubert basis to the Schubert basis.
\end{enumerate}
\end{Theorem}

Following
\cite{RichmondSlofstra2022isomorphismproblemschubertvarieties}, 
we briefly sketch the proof of \Cref{thm:Richmond-Slofstra} which
follows the (1) implies (2) implies (3) implies (1) approach.
Assuming $v$ and $w$ are Cartan equivalent via the bijection $\sigma:
S(v) \to S(w)$, then there is also a bijection between the inversion
sets of $v$ and $w$. Inversion sets have two nice properties.  For any
$i<j<k$, 
\begin{enumerate}
\item if $(i,j)$ and $(j,k)$ are in the inversion set of $w$, then
$(i,k)$ is also, and 
\item if $(i,k)$ is in the inversion set, then either $(i,j)$ or
$(j,k)$ is also.  
\end{enumerate}
The converse also holds:  any set of size 2 subsets of positive
integers such that (1) and (2) hold are the inversion set of some
permutation \cite[Prop. 2.2]{YO}
%% citation from Ziegler 1992
These nice properties imply that $\sigma$ also induces an isomorphism
between the subgroups generated by the supports determined by the
generators in $S(v) \simeq S(w)$, the corresponding Lie algebras of
parabolic subgroups of $\GL_{n}$, and the finite dimensional nilpotent
Lie subalgebras associated with the span of the root spaces indexed by
the inversion sets of $v$ and $w$, denoted $\mathfrak{n}^{+}_{v}$ and
$\mathfrak{n}^{+}_{w}$.  Furthermore, there exist suitably chosen
maximal, integrable, highest weight modules $V,W$ for the Lie algebra
of the parabolic subgroup of $\GL_{n}$ associated with the group
generated by $S(v)$ and $S(w)$ respectively with the same highest
weight $\lambda$, so $\pi: V \to W$ is an isomorphism.  By
construction, $\mathfrak{n}^{+}_{v}$ and $\mathfrak{n}^{+}_{w}$ act on
$V$ and $W$ respectively and choosing the dominant integral weights in
$V,W$ there exists a natural surjection $\pi: V \to W$.  They observe
that each point $bu \in X(v)$ can be expressed uniquely as a product
with $b=\mathrm{exp}(x)$, $x \in \mathfrak{n}^{+}_{v}$, and $u\leq v$
in Bruhat order.  Applying $bu$ to any fixed highest weight vector
$\omega$ of weight $\lambda$ in $V$ gives an embedding of $X(v)$ into
$\mathbb{P}(V)$.  Similarly, there is an embedding of $X(w)$ into
$\mathbb{P}(W)$ using the highest weight vector $\omega'=\pi(\omega)$.
Using the fact that Schubert varieties are normal over algebraically
closed fields \cite{Ramanan-Ramanathan.1985}, they conclude $\pi$
induces an isomorphism of $X(v)$ and $X(w)$ of Zariski closed subsets
of $\mathbb{P}(W)$ coming directly from $\sigma:S(w) \to S(w)$ and a
suitable choice of highest weight modules, completing (1) implies (2).

Given the explicit construction of an isomorphism $X(v)$ and $X(w)$ of
Zariski closed sets embedded in the same projective space, the
corresponding isomorphism of cohomology rings maps Schubert basis to
Schubert basis as graded rings.  So (2) implies (3).

To complete the proof sketch, observe that the adjacency matrices
$A(v)$ and $A(w)$ can be determined from  the Schubert basis
structure constants and restricting Monk's formula to permutations in
the intervals $[\mathrm{id},v]$ and $[\mathrm{id},w]$, as we have
already seen in \Cref{ex:2314.2143}.  Thus, a ring isomorphism $\phi :
H^{*}(X_{v}, \mathbb{Z}) \longrightarrow H^{*}(X_{w}, \mathbb{Z})$
mapping the Schubert basis to the Schubert basis induces a bijection
$\sigma : S(v) \to S(w)$ proving that $v,w$ are Cartan equivalent,
completing (3) implies (1) and our summary of Richmond-Slofstra's
proof of \Cref{thm:Richmond-Slofstra}.

\bigskip

The Richmond-Slofstra Isomorphism Theorem gives us an easy way to
compute the first few terms of the sequence $I(d)$ counting the number of 
distinct isomorphism classes of Schubert varieties in complete flag
varieties of dimension $d$.  For $d=0,\dots , 8$, we have
\[
1,1,2,6,14,39,106,298,839.
\]
These numbers were calculated as follows.  A permutation $w$
is said to be \textit{connected} if $S(w)=[n-1]$ for some $n$.  For
example, $2341=s_{1}s_{2}s_{3}$ is connected, while $2143=s_{1}s_{3}$
and $1342=s_{2}s_{3}$ are not.  The connected permutations in $S_{n}$
are enumerated in the sequence \cite[A003319]{oeis}. The involution
$\sigma:[n-1] \to [n-1]$ mapping $i$ to $n-i$ for all $1\leq i< n$
induces a group homomorphism on $S_{n}$ sending $s_{i} \to
s_{n-i}=w_{0}s_{i}w_{0}$.  Note, $\sigma$ is the only automorphism on
the path graph with $n-1$ vertices.  Therefore, if $w$ is connected,
its Cartan equivalence class restricted to $S_{n}$ is exactly $\{w,
w_{0}w w_{0} \}$, which could have one or two elements.  Let $CI(d)$
be the number of \textit{connected-isomorphism-types of dimension}
$d$, meaning the number of $w \in S_{\infty}$ such that $\ell(w)=d$,
$S(w)=[n-1]$ for some $n\geq 1$, and $w$ is maximal in lexicographic
order in its Cartan equivalence class.  The sequence $CI(d)$ starts out
for $d=0,1,2,\dots ,8$ with
\[
1, 1, 1, 4, 7, 21, 49, 139, 362.
\]
For example, the 7 connected permutations of length 4 which are maximal
in their Cartan equivalence classes are
\[
4213,4132,3412,51234,25134,31524,41253 
\]
written in one-line notation.  They correspond with the reduced words
\[(3212),(3213),(2132),(4321),(1432),(2143),(3214).
\]
Note, we consider the identity element in $S_{1}$ to be the unique
connected permutation with $\ell(w)=0$.  The Schubert varieties
indexed by identity elements in any $\Fl(n)$ are all isomorphic to a
point.

\begin{Corollary}\label{cor:isom.counts} For any positive integer $d$,
the number of distinct isomorphism classes of Schubert varieties of
dimension $d$ in all complete flag varieties is
\begin{equation}\label{eq:count.isom.types}
I(d)=\sum_{\lambda=(1^{m_{1}}2^{m_{2}}\cdots) \vdash
d} \  \prod_{j\geq 1} \binom{CI(j)+m_{j}-1 }{m_{j}}. 
\end{equation}
\end{Corollary}

\begin{Remark}
Does the formula above have a simpler closed form?
\end{Remark}

\begin{proof}
For $\mathrm{id} \neq w \in S_{n}$, assume $G(w)$ has $k$ connected
components, which we know are each isomorphic to a nonempty path.  By
\Cref{thm:Richmond-Slofstra}, the isomorphism type of the Schubert
variety $X_{w}$ as a subvariety of the complete flag manifold $\Fl(n)$
is determined by the factorization of $w=v^{(1)}\cdots v^{(k)}$ as a
product of commuting permutations indexed by the connected component
of $G(w)$. Each factor $v^{(i)}$ is Cartan equivalent to a unique
connected permutation with the same number of inversions and the same
size support which is also lexicographically maximal in its
equivalence class representing its isomorphism type, call it
$C(v^{(i)})$.  The isomorphism type of $X_{w}$ is
determined by the multiset $m(w)=\{C(v^{(1)}),\ldots, C(v^{(k)}) \}$.
Hence, the number $I(d)$ of Schubert variety isomorphism types
for dimension $d$ in all complete flag varieties $\Fl(n)$ for $n\geq 1$
can be determined by enumerating the number of multisets of
connected-isomorphism-types whose dimensions sum to $d$.

Say $\lambda=(\lambda_{1}\geq \lambda_{2}\geq \dots \geq
\lambda_{k}>0)$ is a partition of $d>0$.  If $\lambda \vdash d$ has
$m_{j}$ copies of the value $j$, then the number of multisets of size
$m_{j}$ of connected-isomorphism-types of dimension $j$ is given by
$\binom{CI(j)+m_{j}-1}{m_{j}}$.  The formula for $I(d)$ in
\eqref{eq:count.isom.types} now follows by enumerating over all
partitions of $d$ since each submultiset of $m(w)$ of dimension $j$
connected-isomorphism-types can each be chosen independently.
\end{proof}

\begin{Exercise}\label{ex:inverse.nonisomorphism}
Prove that $X_{w}$ and $X_{w^{-1}}$ need not be isomorphic by
identifying such a $w \in S_{4}$.
%% reduced word 231 and 132 have no Cartan equivalence bijection.
\end{Exercise}

\begin{Problem} Find a complete characterization of all isomorphism types
for Schubert varieties in all partial flag varieties.
\end{Problem}

%% file: section6.tex
\section{And Beyond}\label{sub:Beyond.Section6}

Many of the constructions in this introductory chapter can be
fruitfully generalized, often simultaneously, leading to a dizzying
array of further advances in the field and potential research
problems.  We end this chapter with a survey of some of these further
directions.  The rest of book will further develop many of the
beautiful structures inspired by the study of Schubert varieties.

\subsection{From $\GL_{n}$ to more general Lie groups: or, ``what about type B?''} \label{subsec:type-B}

A \emph{generalized flag variety} (over $\C$) is a smooth projective
variety of the form $G/P$ where $G$ is a complex connected linear
algebraic group and $P$ is a \emph{parabolic} subgroup, meaning it
contains a maximal connected solvable subgroup or \emph{Borel
subgroup} $B$ of $G$.  There are finitely many parabolic subgroups
containing each fixed $B$, but every conjugate of a Borel subgroup is
again a Borel subgroup.  Every Borel subgroup $B$ contains a maximal
torus $T$, generalizing the role of the invertible diagonal matrices.
One natural reason to study general parabolic subgroups $P$ is that they are those 
for which the quotient $G/P$ is a compact manifold.

The permutation matrices are special in $\GL_{n}$ because they form a
finite subgroup representing the entire set of flags in
$\GL_{n}/B$ that are fixed by the action of left multiplication by
$T$.  Furthermore, their left $B$-orbits are the Schubert cells in
$\GL_{n}/B$.  Note, the key property for a flag $wB$ to be a $T$-fixed
point is that for every $t \in T$, $tw=wt'$ for some $t'\in T $.  To generalize the
symmetric group, we consider the normalizer of $T$ in $G$, denoted
$N_{G}(T)$, modulo $T$.  The \emph{Weyl group} associated with $G$,
defined as $W=N_{G}(T)/T$, is a finite group much like the symmetric
group.  It is also isomorphic to a group generated by reflections in
some Euclidean space $\bR^d$, and hence is a Coxeter group with a
presentation similar to the simple transpositions with their
commutation and braid relations.  By choosing coset representatives,
one can view $W=N(T)/T$ as a subset of $G/B$, which necessarily
consists of the $T$-fixed points under left multiplication.
Therefore, much of the basic algebraic and geometric theory we have
developed surrounding flag varieties can be adapted to the generalized
flag variety case.

Familiar examples of complex connected linear algebraic groups include
the special linear groups $\SL_{n}$, as well as the special orthogonal groups
$\SO_{n}$ and symplectic groups $\Sp_{2n}$ over $\C$.  These are
called the \emph{classical groups}.  Their Weyl groups are
generalizations of permutations where one allows a possible sign flip
on each entry in the one-line notation, much like a little kid
shuffling a deck of cards and not caring if the cards all face up or
down.  The simply connected complex Lie groups are linear algebraic groups
that are completely classified by their Dynkin diagram types given by
the finite list
\[
A_{n}, B_{n}, C_{n},
D_{n}, E_{6}, E_{7},E_{8}, F_{4}, G_{2},
\]
where the subscript is called the \textit{rank} of the Lie group and
of its Weyl group.  The case we have studied throughout this chapter is
technically $A_{n-1}$, where $G$ is equivalently given by $\GL_n(\C)$
or $\SL_{n}(\C)$ and the Weyl group is $S_n$.  When people ask ``What
about type $B_{n}$?''  they mean the case $G=\SO_{2n+1}$ with its Weyl
group given by the set of all signed permutations on $[n]$.  Type
$C_{n}$ is the case $G=\Sp_{2n}$, and it just so happens to have an
isomorphic Weyl group, but its Dynkin diagram, which is a labeled
multigraph, is a slightly different encoding some of the algebraic and
geometric differences.  Type $D_{n}$ is the case $G=\SO_{2n}$, and its
Weyl group is the subgroup of signed permutations on $[n]$ with an even number of negative values.

A remarkable classification reduces the study of generalized flag
varieties to a manageable list of cases, as follows.  Specifically, any $G/P$
is a product of the form $G_1/P_1 \times \cdots \times G_m/P_m$ where
each $G_i$ is simple and $P_{i}$ is a parabolic subgroup. It can be shown that every Borel subgroup $B$
contains the center $Z(G)$, which implies that if $G/Z(G) \simeq
G'/Z(G')$, then $G$ and $G'$ have the same generalized flag
varieties.  In turn, the classification of the simple Lie groups by their Dynkin diagrams as listed above
 gives rise to the same classification for simple
complex linear algebraic groups modulo the equivalence relation $G
\sim G'$ if $G/Z(G) \simeq G'/Z(G')$.  Hence, by studying the flag
varieties for simple Lie groups, one essentially understands all
generalized flag varieties.   

What are the Schubert cells and Schubert varieties in generalized flag
varieties?  The reader may go back and review how the definitions from
\Cref{sec:flags} generalize to other matrix groups, even though the
generalization of a flag might not necessarily be so clear.  Let's go
through the basic setup.  Let $G$ be a complex linear algebraic group.
 The $B$-orbits on $G/B$ under left multiplication are exactly the
$B$-orbits of the elements $w \in W$, once again called \emph{Schubert
cells} in $G/B$.  The closure of the Schubert cells are the
\emph{Schubert varieties}. The cohomology ring (or Chow ring)
$H^*(G/B; \bZ)$ is isomorphic to a polynomial ring over $\mathbb{Z}$
modulo the ideal generated by positive degree $W$-invariants. Again,
the cohomology ring has an integral basis given by the classes of its
Schubert varieties.  Similar statements can be derived for any
generalized flag variety $G/P$ by projecting from $G/B$, as in
\S\ref{sub:123.StepFlags}.

Much work has been done to adapt the combinatorial theory surrounding
Schubert polynomials to other Lie groups, but $G = \GL_n(\C)$ remains
by far the best-understood case.  Historically,
Bernstein-Gelfand-Gelfand (BGG) \cite{BGG} and Demazure \cite{Dem}
gave a uniform description of polynomial representatives for Schubert
classes for all simple Lie group types.  For example, the BGG formula
for the analog of $[X_{\mathrm{id}}]$ is the product of all positive
roots as a polynomial function on $\mathrm{Lie}(T)$, where the roots are a special
finite set of vectors in $\mathbb{R}^{n}$ invariant under the action
of the Weyl group, and the positive ones lie on the positive side of a
hyperplane.  The other Schubert classes are all obtained by the analog
of the divided difference recurrence \eqref{eq:divided.difference}.
The BGG and Demazure representatives do not all have positive
expansions into monomials, and there are no known combinatorial or
geometrical interpretations for the terms.  Furthermore, using these
representatives to compute Schubert structure constants requires
multiplying and expanding modulo the ideal generated by the positive
degree invariants.  Lascoux and Sch\"uztenberger made a huge
contribution to the development of Schubert calculus when they
discovered the Schubert polynomials, which can be multiplied and
expanded in the basis of Schubert polynomials without the need to work
modulo the ideal!

During the 1980s, a major open problem was to find polynomial
representatives for the Schubert classes corresponding with the
classical Lie groups that had the nice properties of the type $A_{n}$
Schubert polynomials: they are stable under inclusion, they are a
basis for the space they span, and they have structure constants equal
to the intersection numbers of the corresponding Schubert classes.
The first solution to this problem was given by Billey-Haiman (1995)
\cite{BH}. Their generalized Schubert polynomials are defined via a
formula summing terms over all reduced expressions for the
corresponding Weyl group element in analogy with the BJS formula
\Cref{cor:bjs}, resulting in polynomials in the ring
$\mathbb{Q}[x_{1},x_{2},\dots; p_{1}, p_{3}, p_{5}, \ldots]$ where
$p_{k}=x_{1}^{k}+x_{2}^{k}+ \cdots$ is the $k^{th}$-power sum in the
countably infinite set of variables $x_{1},x_{2},\dots$.  Chevalley
had shown that the homogeneous finite degree power series in
$x_{1},x_{2},\dots$ that are invariant under the action of all signed
permutations are exactly the even power sums $p_{2n}$
\cite{Chevalley.invariants}.  Hence, $\mathbb{Q}[x_{1},x_{2},\dots;
p_{1}, p_{3}, p_{5}, \ldots]$ is the limiting case of the analog of
the coinvariant algebra (\Cref{def:coinvariant.algebra}) for types
$B,C,D$.  Note that the odd power sums and variables $x_{i}$ are
algebraically independent in the power series ring
$\mathbb{Q}[[x_{1},x_{2},\dots]]$, so one can treat
$\mathbb{Q}[x_{1},x_{2},\dots; p_{1}, p_{3}, p_{5}, \ldots]$ as a
freely generated polynomial ring in two sets of variables.  Taking
$W_{n}$ to be the Weyl group of one of the classical
groups $G$ of type $A_n$, $B_n$, $C_n$, or $D_n$, let $W_{\infty}=\bigcup_n W_{n}$.
The flag variety of rank $n$ for that classical group will be of the form
$G/B$.  For example, $\Sp_{10}/B$ is the symplectic flag variety of
rank $5$, with Weyl group $W_{5}$ given by the signed permutations on
$\{1,2,3,4,5\}$.  The generalization of Schubert polynomials given by the
\textit{Billey-Haiman polynomials} have the following key properties
proved in \cite{BH} unless otherwise cited.

\begin{enumerate}
\item Within each classical type, the flag variety of rank $n$
embeds naturally into the one of order $n+1$ as a Schubert variety,
inducing maps on cohomology that send Schubert classes to
Schubert classes.  Using this embedding, the Billey-Haiman polynomials
are stable under the inclusion of the corresponding Weyl groups $W_{n}
\hookrightarrow W_{n+1}$.

\item They are the unique solutions to the infinite system of BGG
divided difference equations together with the condition that the
constant term is 1 if and only if the indexing element is the
identity.  In particular, they include the type $A$ Schubert
polynomials indexed by permutations.

\item They generalize the Schubert classes for isotropic Grassmannians
in terms of Schur $P$-functions, denoted $P_{\lambda}$, when
$G=\SO_{n}$ and Schur $Q$-functions, denoted $Q_{\lambda}$, when
$G=\Sp_{2n}$ in analogy with the relationship between Grassmannians
and Schur functions in type $A_{n-1}$.  Both the Schur $P$-functions
and Schur $Q$-functions have combinatorial interpretations as sums
over \textit{shifted tableaux}.  The shifted analog of the
Littlewood-Richardson coefficients was given by Sagan \cite{Sag2} and
Worley \cite{Worley}.  The isotropic Grassmannian Schubert
calculus was developed by Pragacz \cite{P}.

\item They form a basis for $\mathbb{Q}[x_{1},x_{2},\dots; p_{1},
p_{3}, p_{5}, \ldots]$, and their structure constants are exactly the
Schubert structure constants for each classical type.

\item They expand positively as a sum of monomials of the form
$x^{\alpha }P_{\lambda}$ or $x^{\alpha }Q_{\lambda}$ depending on the
Lie group type.  

\item They satisfy an analog of the transition equation following from
Monk's formula, which is more generally referred to as Chevalley's
formula \cite{B1,chevalley}.

\item They have a pipe dream formula \cite{Smirnov-Tutubalina.2023}.

\item They can be obtained as the specialization of a generalization
of double Schubert polynomials as defined by Ikeda-Mihalcea-Naruse in
\cite{IKEDA.Mihalcea.Naruse.2011}.

\item The limits of the Schubert polynomials indexed by $1^{k}\times w$ 
as $k$ goes to infinity give rise to types $B_{n}$, $C_{n}$, $D_{n}$-Stanley symmetric
functions, which also have an expression as the sum over reduced words
$\mathbf{r}$ for $w$ of terms given by the quasisymmetric functions
\[
\Theta_{P(\mathbf{r})}^{n} = \sum_{i_{1}\leq  i_{2}\leq \cdots \leq
i_{p}} 2^{|\mathbf{i}|} x_{i_{1}} x_{i_{2}} \cdots x_{i_{p}}
\]
summed over all weakly increasing sequences $\mathbf{i}=(i_{1}\leq
i_{2}\leq \cdots \leq i_{p})$ such that if $r_{j}$ is a \textit{peak},
meaning $r_{j-1}<r_{j}>r_{j+1}$, 
then $i_{j-1}, i_{j}, i_{j+1}$ are not all equal.  These
quasisymmetric functions are the shifted analog of the fundamental
quasisymmetric functions indexed by the \emph{peak sets} of
$\mathbf{r}$ instead of the descent set \cite[\S 3]{BH}.  Gessel's \Cref{thm:gessel} for expanding Schur functions into
fundamental quasisymmetric functions has the analog
\[
Q_{\mu} = \sum \Theta_{P(T)}^{n}
\]
summed over all shifted standard tableaux of shape $\mu$
\cite[Prop. 3.2]{BH}.  The proof of this result uses Haiman's theory
of dual equivalence for shifted shapes \cite{haiman1992dual}.  Further
variations on $B,C,D$-Stanley symmetric functions also appeared around
the same time in \cite{FK2,kraskiewicz1989reduced,Lam1,Lam2,stembridge.1997}.

\item Peak sets have inspired a rich area of combinatorics.  The
\emph{shifted quasisymmetric functions} SQSYM spanned by the
$\Theta_{P(\mathbf{a})}^{n}$'s over the rational numbers form a
subspace of QSYM containing both the Schur $P$-functions and Schur
$Q$-functions.  It also contains the closely related \emph{peak subalgebra}
introduced by Stembridge in connection with enriched $P$-partitions
\cite{stembridge1997enriched}.  The peak subalgebra of QSYM was shown
to be a Hopf algebra, but SQSYM is only known to be a coalgebra
\cite{BMSvW.2002}.  There is also a peak algebra of the symmetric
group algebra, given by the linear span of sums of permutations in
$S_{n}$ with a given peak set as defined by Nyman \cite{Nyman.2003}. The
number of permutations in $S_{n}$ with a given peak set has a simple
enumerative formula in terms of peak polynomials
\cite{Billey-Burdzy-Sagan}, a positive expansion in the binomial basis
\cite{Diaz-Lopez.Harris.Insko.Perez-Lavin.Omar}, and extensions to
other Coxeter group types
\cite{Diaz-Lopez.Harris.Insko.Perez-Lavin,Orellana.et.al}
\end{enumerate}

There are some properties of type $A_{n}$ Schubert polynomials that
the Billey-Haiman polynomials do not satisfy.  One such key property
is a finite expansion into $x^{\alpha }$ monomials.  After the
Billey-Haiman polynomials were defined, Fomin-Kirillov examined
several different ways of defining Schubert polynomials for type
$B_{n}$ in \cite{FK2} using different subsets of nice properties that
occur for type $A_{n}$ Schubert polynomials.  They proved that no
family of representatives for Schubert classes in types $B,C,D$ satisfy
all such properties.

Around the same time, Fulton (1996) gave another approach to computing
the classes of degeneracy loci in all classical types generalizing his
formula \Cref{thm:fulton-degeneracy-loci}.  In his approach, the class
corresponding with the identity in each type is given by the
determinant of a matrix with polynomial entries similar to a
Vandermonde matrix, but in two sets of variables.  Therefore, his
approach was the first analog of double Schubert polynomials for the
classical groups.  More recently (2018), Anderson-Fulton
\cite{Anderson-Fulton.2018} returned to the topic of Chern class
formulas for degeneracy loci for classical groups using the double
Schubert polynomials of Ikeda-Mihalcea-Naruse along with the
multi-theta and multi-eta polynomials introduced by
Buch-Kresch-Tamvakis \cite{BKT.1,BKT.2} and Tamvakis-Wilson
\cite{TW}. Their degeneracy loci formulas cover all isotropic
Grassmannians of types $B_{n}$, $C_{n}$, $D_{n}$, they generalize the
previously known Pfaffian and determinantal formulas, and they have
(what they claim are) simple direct geometric proofs.

The Lakshmibai-Sandhya \Cref{Lak-San} characterizing the smooth
Schubert varieties $X_{w}$ by avoiding the patterns 3412 and 4231 has
been generalized in many ways.  Billey-Postnikov generalized the
notion of pattern avoidance using root subsystems in order to
completely characterize all Schubert varieties in $G/B$ for all
complex semisimple Lie groups $G$ independent of the type
\cite{BP-smooth}.  Their generalization of pattern avoidance extends
to all Coxeter groups and has applications in terms of Kazhdan-Lusztig
polynomials \cite{BiBr}.  Slofstra-Richmond have found enumerative
formulas to count all smooth Schubert varieties in all simple Lie
group types \cite{RS15}.  See Chapters 2, 3, and 5 of this book for more
pattern avoidance connections beyond type $A$ 
\cite{mcgovern2022patternavoidancekorbitclosures,oh2024coxetergroupsbilleypostnikovdecompositions,wooyong2023schubertgeometrycombinatoricsCh}.  There is also a complete characterization for smooth positroid
varieties using a circular analog of permutation patterns
\cite{Billey.Weaver.2025}, patterns characterizing smoothness for
Hessenberg Schubert varieties \cite{Cho.Huh.Park}, and patterns
characterizing spherical varieties \cite{Gaetz.2022}.

\subsection{Generalized Cohomology Theories}\label{sub:exotic}  

After over a century of development, the basic notion of a cohomology
ring for a topological space has been generalized in many directions,
such as quantum, equivariant, K-theory, and cobordism.  In fact, there
are an infinite number of cohomology theories to consider.  See
\cite{Conner-Floyd.1964} for some historical background.  Which ones
will be most beneficial?  That's hard to say at this time.  It depends
on the applications and how finely one needs to detect differences in
the topological objects of study.  However, studying flag varieties,
Grassmannians, Schubert varieties, and Richardson varieties in any
cohomology theory will lead to motivational examples.  Here, we
briefly mention some of the many known results about Schubert
varieties in these more ``exotic'' cohomology theories.

At the end of the inspirational paper by
Bernstein-Gelfand-Gelfand \cite{BGG}, they also identified another
uniform description of polynomial representatives for Schubert classes
for all simple Lie types, which they attribute to Bertram Kostant.
These representatives of Schubert classes are defined by certain
vanishing conditions instead of an explicit formula.  The
\emph{Kostant polynomials} were studied further by Kostant-Kumar
\cite{KK} and Billey \cite{B4}.  In particular, Billey \cite{B4}
described a uniform algorithm to obtain Schubert structure constants
using the Kostant polynomial extending to Kac-Moody groups.  Another
uniform approach to computing these constants is due to Duan
\cite{duan}.  Kostant-Kumar developed the nil-Hecke approach to
cohomology that will be further examined in Chapter 12 of this book.

The Kostant polynomials are closely related to evaluations of the
double Schubert polynomials, see \cite[Remark 1 in Sec.8]{B4} in
type $A$.  Kostant's approach and the so-called AJS/Billey formula,
used in \eqref{eq:ajs.billey.formula}, has since led to explicit
representatives for the equivariant cohomology classes of Schubert
varieties extending even beyond semisimple Lie groups to all Kac-Moody
groups.  The proof relies on the theory of GKM spaces due to
Goresky-Kottwitz-MacPherson \cite{GKM} who proved that any family of
polynomial functions satisfying the right vanishing properties on the
Bruhat graph must be the Schubert classes in equivariant
cohomology.  For a survey on applications to combinatorics, geometry
and topology, see \cite{Tymoczko.2016}. In \cite{Graham.2001}, Graham
proved that products of the polynomial functions giving rise to
equivariant cohomology classes of Schubert varieties expand in the
equivariant Schubert class basis with structure constants in the
positive span of monomials in the simple roots, proving conjectures of
Billey and Peterson.  This notion of \textit{Graham positivity} plays
the role of the nonnegativity of structure constants in ordinary cohomology 
for the flag variety.

For a comprehensive treatment of equivariant cohomology rings for
Schubert varieties and degeneracy loci see the book by Anderson-Fulton
\cite{Anderson-Fulton}.  Another helpful resource on equivariant
cohomology and equivariant intersection theory was written by Brion
\cite{Brion:equivariant}.  He spells out some of the connections
between the study of singular points in projective varieties with a
group action and equivariant cohomology, as hinted at in
\Cref{thm:BK-smooth-test}.

There is another generalization of the cohomology ring of a
nonsingular projective variety $X$ called the \textit{Grothendieck
ring} of $X$, denoted $K(X)$, related to the coherent sheaves and
vector bundles on $X$ modulo certain relations.  The product structure
comes from the tensor product of vector bundles on $X$. Like with the
Chow groups, the addition rule is formal addition of symbols
$[\mathcal{F}]+[\mathcal{G}]$ modulo commutation relations coming from
short exact sequences.  We briefly mentioned K-theory in connection
with the K-polynomial of an ideal (\Cref{def:k-poly}).  The analogs of
Schubert polynomials in this setting are the Grothendieck polynomials
defined by Lascoux and Sch\"utzenberger \cite{LS6}.  Fomin-Kirillov
proved that the Grothendieck polynomials can be realized as solutions
to Yang-Baxter equations \cite{fomin1994grothendieck}.  Weigandt gave
a formula for double Grothendieck polynomials as a sum over marked
bumpless pipe dreams \cite{weigandt-BPD-ASM}.  Fulton-Lascoux
\cite{Fulton-Lascoux.1994} gave an analog of Monk's formula for the
K-theory of generalized flag varieties, and these results were
extended to generalized flag varieties by Pittie-Ram
\cite{Pittie-Ram.1999}.  See also Lenart's variation on Monk's formula
for K-theory \cite{Lenart.2003}.  Buch first studied the analog of the
Littlewood-Richardson coefficients for the Grothendieck ring of a
Grassmannian using \emph{set-valued} tableaux \cite{buch2002littlewood}.
Ikeda-Naruse considered the set-valued tableau analogs of Schur
$P$-functions and Schur $Q$-functions to obtain Schubert classes in
the K-theory of torus equivariant coherent sheaves on the maximal
isotropic Grassmannians of symplectic and orthogonal types
\cite{Ikeda-Naruse.2013}.  Brion proved a conjecture of Buch's for the
structure constants in the Grothendieck ring of a flag variety in
\cite{Brion.2002}.  More generally, these classes were studied by
Demazure \cite{Dem} and Kostant-Kumar \cite{Kostant-Kumar.1990}.  A
nice survey of K-theory related to flag varieties was given in
``Lectures on the geometry of flag varieties'' by Brion
\cite{Brion:flag}.
% \url{https://arxiv.org/pdf/math/0410240.pdf}
As mentioned above, it can be productive to consider combining
cohomology theories.  See
\cite[Sect. 8]{CAMERON.DINU.MICHALEK.SEYNNAEVE} for a brief
introduction to the equivariant K-theory of Grassmannians and flag
varieties in connection with matroids.

Quantum cohomology is a variation incorporating ideas from
intersection theory with string theory in physics related to the work
of Vafa \cite{Vafa.1992} and Witten \cite{Witten.1991}.  The quantum
cohomology ring of the flag variety is a deformation of the classical
cohomology of a space depending on additional parameters incorporated
to measure the deformation. Fomin-Gelfand-Postnikov gave the quantum analog of Schubert
polynomials and the quantum variation of Monk's formula
\cite{Fomin.Gelfand.Postnikov.Quantum}.  This paper inspired great
interest and many extensions to equivariant quantum cohomology.  See
for example \cite{Anderson-Chen,Buch-Kresh-Tamvakis,
Ciocan-Fontanine.1999,Chen.2003,Givental-Kim.1995}.

The additive structure in quantum
cohomology is similar to that of ordinary cohomology as described via
the Chow ring in \Cref{sub:Monk}.  As an abelian group, it
is isomorphic to $H^{*}(\Fl(n),\mathbb{Z}) \otimes \mathbb{Z}[q_{1},
\dots , q_{n}]$ with an analogous basis of Schubert classes.  Using
this basis, the structure constants for multiplication in quantum
cohomology are known as the \emph{3-point Gromov-Witten invariants of
genus 0}.  Instead of counting points in an intersection, they count
equivalence classes of certain rational curves in $\Fl(n)$. Recall from  \Cref{subsec:eigenvalues} that the Littlewood-Richardson coefficients control all possible relations between eigenvalues of Hermitian matrices $H_1, H_2, H_1+H_2$. Agnihotri and Woodward \cite{agnihotri-woodward} and Belkale \cite{belkale-quantum} showed that 3-point Gromov-Witten invariants for Grassmannians solve the analogous multiplicative problem, relating the possible eigenvalues of unitary matrices $U_1, U_2, U_1 U_2$. For an interesting application of this result to quantum computing, see \cite{monodromy-polytope}---beware that the two appearances of the word ``quantum'' are unrelated!

New insights into the geometry and combinatorics of quantum Schubert
calculus came from a surprising connection with the $k$-Schur
functions originally defined by Lapointe-Lascoux-Morse \cite{LLM2003}
in the context of studying Macdonald polynomials.  Shimozono
conjectured, and Lam proved, that the $k$-Schur functions evaluated at
$q=1$ represent the Schubert basis for the homology ring of the affine
Grassmannians \cite{Lam.2008} inside the isomorphic ring of symmetric
functions.  Lapointe-Morse showed that the quantum Schubert structure
constants, and therefore also the ordinary Schubert structure
constants, are given by certain structure constants for $k$-Schur
functions \cite{Lapointe-Morse.2008}.  Many further developments in
this ``Quantum = Affine'' direction combining combinatorics, geometry, and
representation theory are further described in
\cite{LLMSSZ,LLMS,Lam-Shimo.2010,Lam.Shimozono.2012}.

\subsection{Connections to Representation Theory}\label{}

The importance of the Schur polynomials in representation theory stems from the fact that they are the irreducible characters of general linear groups. To elaborate, suppose $W$ is a polynomial representation of $\GL_n(\C)$, meaning that there is a group homomorphism $\phi : \GL_n(\C) \to \GL(W)$ for which the matrix entries of $\phi(A)$ are polynomials in the entries of $A$. Its \emph{character} $\chr(W)$ is defined as in \Cref{sub:MatrixSchubs} with respect to the action of any $n$-dimensional torus $T \subseteq \GL_n(\C)$. It can be shown that $\chr(W)$ is independent of the choice of $T$, and that it completely determines $W$ up to isomorphism. Concretely, $\chr(W)$ is the trace of $\rho(\operatorname{diag}(x_1, \ldots, x_n))$. This is a symmetric polynomial in $x_1,\ldots,x_n$ (exercise!) and it turns out that a symmetric polynomial $f(x_1,\ldots,x_n)$ is a character if and only if $f$ is a nonnegative linear combination of the Schur polynomials $s_\lambda(x_1,\ldots,x_n)$.

Given that Schubert polynomials generalize Schur polynomials, one
might wonder if they too have representation-theoretic
meaning. Kra\'skiewicz and Pragacz gave one answer, realizing each
Schubert polynomial $\fS_w$ as the character of a $B$-representation
\cite{kraskiewicz-pragacz}. Magyar later explained their construction
in terms of a very classical method of obtaining representations of
Lie groups, which we outline now.  Suppose $G$ and $B$ are as in
\S\ref{subsec:type-B}, and $\lambda : B \to \C^\times$ is a group
homomorphism. Define a topological space $L_\lambda$ as $G \times \C$
modulo the equivalence relation $(g,x) \sim (gb^{-1},\lambda(b)x)$ for
all $b \in B$.  The projection $\pi : L_\lambda \to G/B$ given by
$(g,x) \mapsto gB$ makes $L_{\lambda}$ into a line bundle on $G/B$.
Let $H^0(G/B, L_\lambda)$ denote\footnote{Observe that $H^{0}(X, L)$ is not the degree 0 graded
component of singular cohomology.  However, it is related to sheaf
cohomology in homological algebra.} the vector space of holomorphic
sections of the vector bundle $L_\lambda$ with base space
$G/B$.  Explicitly, sections of
$L_\lambda$ correspond to functions $s : G \to \C$ such that $s(gb) =
\lambda(b)s(g)$ for all $g \in G, b \in B$. If $h \in G$ and $s$ is
such a function, then so is the function $h \cdot s : g \mapsto
s(h^{-1}g)$. This $G$-action makes $H^0(G/B, L_\lambda)$ into a
representation of $G$.

\begin{Theorem}[\textbf{Borel-Weil Theorem}] \cite[\S 6.1.11]{chriss-ginzburg}
Every irreducible holomorphic representation of $G$ has the form $H^0(G/B, L_\lambda)$ for some $\lambda$.
\end{Theorem}

Magyar realized Schubert polynomials in the style of the Borel-Weil
Theorem, as characters of $B$-representations $\chr H^0(Z, L)$ for
some line bundles $L$ on \emph{Bott-Samelson varieties} $Z$
\cite{MagyarSchubert}. The Bott-Samelson varieties are resolutions
 of singularities for Schubert varieties based on
reduced expressions \cite{BottSamelson}.  In fact, Magyar obtained
characters of more general \emph{flagged Schur modules}
\cite{MagyarBorelWeil}, whose combinatorial properties were explored
by Reiner and Shimozono \cite{RSpercent}.  Magyar's approach to
defining Bott-Samelson varieties is elegantly explained in the context
of the combinatorics of commutation classes of reduced expressions and
Elnitsky tilings in \cite{Escobar.etal.2018}.

Alternatively, consider the line bundle $L_\lambda$ restricted to the
Schubert variety $X_w \subseteq \flags$. Then $H^0(X_w, L_\lambda)$ is
also a $B$-representation, and its character is known as a \emph{key
polynomial} $\kappa_{w,\lambda}(x_1,\ldots,x_n)$. These polynomials
are monomial-positive, and Schubert polynomials expand positively into
key polynomials with explicit combinatorial formulas \cite{LS5}, among
other interesting properties.  The key polynomials can be computed
using Kohnert's bubble rule on a key diagram, which was the
inspiration for chutes and ladder moves on pipe dreams \cite{Kohnert.1990}.

There has been some work leveraging this point of view to understand
Schubert polynomials and Stanley symmetric functions in the context of
representation theory
\cite{zero-one-schubert,Billey-Pawlowski,watanabe}. However, compared to the
substantial body of results connecting Schubert polynomials to the
geometry of Schubert varieties, cohomology of flag varieties, and so
on, this representation-theoretic perspective on Schubert polynomials
as characters seems underused. That said, the Borel-Weil theorem is foundational in the broader field of geometric representation theory. For more insight on representation theory in connection with the geometry of the flag variety, see
\cite{McGovern.2023}.

Many other combinatorial and representation-theoretic interpretations
of Schubert polynomials can be found in the literature.  See for
example \cite{Assaf-Searles,KP,RS2,Weigandt-Yong.2015}.

\subsection{Flag varieties over other fields}\label{sub:fields}

From the beginning of this chapter we have worked with flags in
$\C^n$. Why not replace $\C$ with some other field $\bK$? All of the
basic definitions of flag varieties, Schubert varieties,
Grassmannians, etc. still work, and some properties such as dimension
counts, the Bruhat decomposition, and Bruhat order are unchanged.  For
example, recall Deodhar's \Cref{thm:deodhar.finite.field} for counting
the number of flags in an open Richardson variety over a finite field.
The trouble begins when we start taking intersections. The problem can
be seen easily in the Pl\"ucker embedding, where intersecting Schubert
varieties amounts to solving some linear equations in the Pl\"ucker
coordinates together with the quadratic Pl\"ucker relations. These
quadratic equations may not have as many solutions in $\bK$ as they do
in its algebraic closure.

Consider the Grassmannian case first. Fix partitions $\lambda^1,
\ldots, \lambda^m$ fitting into a $k \times (n-k)$ rectangle such that
$\sum_j |\lambda^j| = k(n-k)$. We would at least like to know that
there exist \emph{some} flags $E^1_\bullet, \ldots, E^m_\bullet$ such
that $X_{(\lambda^1)^\vee}(E^1_\bullet) \cap \cdots \cap
X_{(\lambda^m)^\vee}(E^m_\bullet)$ is 0-dimensional with cardinality
matching the number computed by Chow ring techniques. Vakil
\cite{Vakil-B} showed that this is indeed the case for various fields
$\bK$, including the real numbers, any algebraically closed field, and
sufficiently large finite fields (compared to the fixed data $k,n,\lambda^{1},
\ldots, \lambda^{m}$).

However, Vakil's result is abstract and does not explain how the flags
$E^i_\bullet$ may be chosen. An earlier conjecture of Shapiro and
Shapiro suggested an explicit choice of flags for the important case
$\bK = \bR$. This conjecture was eventually proven in a more general
form by Mukhin, Tarasov, and Varchenko \cite{Mukhin-Tarasov-Varchenko,
Mukhin-Tarasov-Varchenko-2}, providing a way to choose flags in
$\bR^n$ making any desired intersection of Grassmannian Schubert cells
transverse and with expected cardinality. Their remarkable proof
passes through representation theory, integrable systems, and
differential operators, and along the way provides a direct
explanation for why the Littlewood-Richardson coefficients decompose
both products of Grassmannian Schubert classes as well as tensor
products of $\GL_n$ representations. A related result which is closer
to the language and ideas used in \Cref{sec:flags} is due to
Karp-Purbhoo \cite{karp-purbhoo}, which interprets the eigenvalues of
a commuting family of Hermitian operators in the group algebra
$\C[S_d]$ as Pl\"ucker coordinates for real points in $\Gr(k,n)$.  See
\cite{Sottile.book.2011} for an introduction to real algebraic
geometry with an emphasis on Schubert calculus and the
Mukhin-Tarasov-Varchenko theorem. 

The method of \emph{Frobenius splitting} developed by Mehta and
Ramanathan \cite{MehtaRamanathan1985} uses the Frobenius automorphism
of a field of characteristic $p$ to prove strong geometric results
about families of varieties. As one example, this method applies to
the family of Schubert varieties in $\flags$ with respect to a fixed
flag $E_\bullet$, and implies that any union or intersection of
Schubert varieties is reduced. In algebraic terms this means that any
intersection or sum of their defining ideals is automatically a
radical ideal. The truth of this statement in characteristic $p$ for
every $p$ even implies it in characteristic $0$! For an introduction
to Frobenius splitting methods, see \cite{BrionKumar2005}. More
generally, flag varieties and Schubert varieties over arbitrary fields
play an important role in representation theory for algebraic groups,
as in Jantzen's book on the topic \cite{jantzen}.

\subsection{Further Topics}

There are of course an uncountable number of subvarieties of a flag
variety, and almost as large a number have been studied. A few have
already appeared in this chapter, with more to come in later chapters
such as Kazhdan-Lusztig varieties (Chapter 2), generalized smooth
Schubert varieties (Chapter 3), Richardson varieties and positroid
varieties (Chapter 4), K-orbit closures (Chapter 5), torus orbit
closures (Chapter 6), spanning line configurations (Chapter 7),
different flavors of Hessenberg varieties (Chapters 8-11), and
generalizations to all Kac-Moody flag varieties (Chapter 12).  In this
final subsection, we attempt to list some topics that seem deserving of mention but have not fit in anywhere else. 

\begin{itemize}
\item Many aspects of Schubert calculus can be generalized to $K$-orbits on a flag variety $G/B$ where $K$ is a \emph{spherical subgroup}---which just means there are a finite number of orbits. There are analogues of weak and strong Bruhat order \cite{incitti.bruhat,richardson-springer}, Schubert and Grothendieck polynomials and associated pipe dream and transition combinatorics \cite{can-joyce-wyser, HamakerMarbergPawlowski2018, HamakerMarbergPawlowski2018Transition, HamakerMarbergPawlowski2022, wyser-yong}, reduced words and Stanley symmetric functions \cite{burks2018reducedwordsclans,HamakerMarbergPawlowski2019,hansson-hultman,marberg-pawlowski-signed-involutions}, and matrix Schubert varieties \cite{marberg-pawlowski-k-theory}, among other topics \cite{Gaetz.2022,gao-hodges-yong}. Smoothness for spherical orbits will be discussed in Chapter 5. Most of these results are specific to the $G = \GL_n$ case, however, so much work remains.
\item Generalizing Schubert varieties often results in interesting degeneracy locus formulas. For instance, one might consider vector bundles equipped with a bilinear form and impose conditions based on it \cite{Anderson-Fulton.2018, wyser-degen}; these loci are related to type B, C, and D Schubert varieties and to the spherical orbits mentioned above. Results on \emph{quiver loci} generalize Fulton's flagged degeneracy locus theorem by placing a vector bundle at every vertex of a directed multigraph and a morphism along each arrow \cite{buch-fulton-quivers, feher-rimanyi-quivers, KinserRajchgot2015, KnutsonMillerShimozono2006}. 

\item Insertion algorithms such as the Robinson-Schensted-Knuth
correspondence and Edelman-Greene insertion have inspired other types
of insertion algorithms.  For example, on the combinatorial side we
have
\cite{huang2023growthdiagramsschubertrsk,huang-RSK,LLMS,LENART.Growth.Diagrams.2010}.
There are also geometric variations such as \cite{BFPS.GeometricRSK}
% Crystal invariant theory I: Geometric RSK Benjamin Brubaker, Gabriel
% Frieden, Pavlo Pylyavskyy, Travis Scrimshaw
and representation theory variations such as
\cite{price2024representationsmatrixvarietiesfiltered}.

\item How can we further our understanding of the computational
complexity issues in the area of Schubert calculus along the lines of
the recent results mentioned at the end of \Cref{sec:flags} such as
\cite{Ikenmeyer-Pak-Panova,Narayanan,pak2022combinatorialinterpretation,pak.robichaux.2024vanishingschubertcoefficients}?
This area is a rich source of problems and results at the intersection
of algebra, combinatorics, geometry, and complexity theory.

\item As mentioned at the end of \Cref{sub:Fomin-Stanley} and in
\Cref{sub:StanleySymmetrics}, the \emph{permutohedral variety}
$X_{A_{n}}$ is the closure of the orbit of invertible diagonal
matrices $T$ acting by left multiplication on a generic flag
$F_\bullet \in \Fl(n+1)$.  It is also the toric variety associated to
the braid fan.  It was studied by Huh \cite{HuhThesis,Huh.2018} and
Klyachko \cite{Klyachko.1985} among many others.  A gentle
introduction to the Chow ring of the permutohedral variety can be
found in the survey article by Ardila-Mantilla \cite{ardila2024},
which covers  many exciting recent developments related to
intersection theory and matroids.

\item ``Schubert varieties and generalizations'' by Tonny Springer
is a survey with more directions you could go, including Kac-Moody
groups, spherical varieties, and the Steinberg variety \cite{Springer1998}.
\end{itemize}

\vspace{.7in}

Dear Readers, We hope you have enjoyed this introduction the
cohomology of the flag variety and the chapters to follow.  It is a
rich subject that ties together so many different aspects of
mathematics.  We have tried to give you an overview of many gems in
the field.  However, the list is endless!  And, we hope you find many
interesting ways to push the frontier even further.

%Sara: Add a page of the degeneration book here as an appendix.

%% file: notation.tex
\section{Notation}\label{sec:notation}

List of notation for use in \textbf{ Handbook of Combinatorial Algebraic Geometry: Subvarieties of the Flag Variety}.\\

\begin{itemize}
\item $[n]:= \{1,2,\ldots, n\}$
\item $S_n$ for the symmetric group on $[n]$; $u, v, w$ permutations
\item $\ell(w)$ and $\inv (w)$ both denote the (Coxeter) length of $w$ or
equivalently the number of inversions of $w$
\item Partitions of a number $n$ are denoted by $\lambda=(\lambda_{1},\lambda_{2}, \dots , \lambda_{k}) \vdash n$
\item We denote the one-line notation of a permutation $w\in S_n$ by
$w(1) w(2) \ldots w(n)$ or $[w_1, w_2, \ldots, w_n]$ if the commas are
helpful.
\item The flag variety is $\Fl(n):=\{F_\bullet:=(F_1\subset F_2\subset \cdots \subset F_n=\C^n)\ |\ \dim V_i=i\}.$
\item Generic flags $G_\bullet, H_\bullet \in \Fl(n)$.
\item  The standard flag is $E_\bullet\in \Fl(n)$ and its opposite flag is
$\oppositeE_{\bullet}$.
\item $C_w=C_{w}(E_{\bullet} = BwB/B$ denotes the Schubert cell with
respect to $E_\bullet$
\item Schubert variety with respect to a flag $F_{\bullet}$:
$$X_{w}=X_w(F_\bullet):=\{V_\bullet\in \Fl(n) \ |\ \dim (F_i\cap V_j)\geq
r_w(i,j)\}$$ where $r_{w}(i,j)$ is the rank of the submatrix of the
permutation matrix for $w$ weakly northwest of $(i,j)$.

\item Grassmannian: $\Gr(k,n):=\{V\subset \C^n\ |\ \dim(V)=k\}$.
\item Cohomology ring $H^*(X, \mathbb{Z})$ of a topological space $X$
over the integers, and its equivariant Cohomology is denoted
$H_T^*(X,\mathbb{Z})$,
\item $R_{n}=\mathbb{Z}[x_{1},x_{2}, \dots , x_{n}]/I_{n}^{+}$ is the
coinvariant algebra, which is isomorphic to $H^*(\Fl(n), \mathbb{Z})$ via
the Borel presentation
%\textcolor{red}{Grade by complex dimension so $H^i(X;\mathbb{C})=H^{2i}(X;\mathbb{R})$.}
\item The Schubert class $[X_{w}]=\fS_{w_{0}w}$ where Schubert
polynomial for $w$ is $\mathfrak{S}_w$;  a general cohomology class for a variety $X$ is $[X]$.
\item $G$ denotes an algebraic group; $\mathfrak{g}$ its Lie algebra
\item $GL_n(\mathbb{C})$ denotes the algebraic group of $n\times n$ invertible matrices over $\mathbb{C}$ with Lie algebra $\mathfrak{gl}_n(\mathbb{C})$
\item $B$ denotes a (fixed) Borel subgroup, and $P, Q$ typically denote parabolic subgroups, with corresponding Lie algebras $\mathfrak{b}, \mathfrak{p}, \mathfrak{q}$, respectively
\item $U$ denotes the unipotent subgroup of $B$, $\mathfrak{u}$ the Lie algebra of $U$
\item $T$ is the maximal torus in $B$
\item Denote the (generalized) flag variety of $G$ by $G/B$.
\item $W$ denotes the Weyl group of $G$, and $S\subseteq W$ the subset of simple reflections
\item Given $J\subseteq S$ let $W^J$ denote the set of shortest coset representatives for the left cosets $W/W_J$.  Similarly, one can use $W^P$ for a standard parabolic subgroup $P.$
\item Given a parabolic subgroup $P_J$ with Weyl group $W_J$, denote by $X_w^J$ the Schubert variety in the partial flag variety $G/P_J$ for each $w\in W^J$.  Similarly, the notation $X_w^P$ may be used.
\end{itemize}

%% file: acknowledgments.tex
\section{Acknowledgments}\label{sec:acknowledgments}

We would like to thank many people for discussions, insights, and
contributions to this manuscript. First, we would like to warmly thank
Erik Insko, Martha Precup, and Edward Richmond for their initiative
and inspiration to collect this material in a book format.  We would
also like to thank the other authors of chapters for their
contributions.  We believe this will strengthen the field and make the
subject more accessible.  Next, we thank William Fulton both for his
contributions to the field and for adding some historical context
about Schubert and the Italian school of algebraic geometry.  Many
people encouraged us along the way to write this chapter and
contributed to the final outcome. We thank Lakshmibai, George Lusztig,
Brendan Rhoades, Richard Stanley, and Greg Warrington for helpful
suggestions that have been incorporated into this work. We thank Dave
Anderson, Agnes Beaudry, Mackenzie Bookamer, Anders Buch, Herman Chau,
Elena Hafner, Clare Minnerath, Jasper Moxuan Liu, Bryan Lu, Connor
McCausland, Igor Pak, Rich\'{a}rd Rim\'{a}nyi, Colleen Robichaux,
Garcia Sun, Sheila Sundaram, Joshua Swanson, Vasu Tewari, Rachel Wu,
and Michael Ruofan Zeng for feedback on an early draft.  We thank
Matja\v z Konvalinka for explaining the proof we included for
Sylvester's Lemma.  We are still hoping to get more feedback and
comments so your name could go here \dots